\documentclass[reqno,a4paper,10pt]{amsbook}
\usepackage[pdftex,unicode,implicit]{hyperref}
\hypersetup{%
  pdftitle    = {Differential spinors on Lorentzian four-manifolds and supersymmetric evolution flows},
  pdfkeywords = {Geometric flows, Spin geometry,Lorentzian four manifolds, Cauchy hypersurface, supergravity, T duality}, 
  pdfsubject  = {math.dg,hep-th, gr-qc},
  pdfauthor   = {C.S.Shahbazi},
  pdfcreator  = {pdf\LaTeXe\ with package \flqq hyperref\frqq},
  pdffitwindow= false,  
  bookmarksopen   = false,
  bookmarksopenlevel = 1,
  bookmarksnumbered = true,
  unicode     = true,
  plainpages  = false,
  colorlinks  = true,  
  citecolor   = black,
  urlcolor    = black,
  linkcolor   = black,
}

\usepackage{amssymb,amstext,amsthm,amsfonts,mathrsfs,amsbsy,amsmath,array,bbm,multirow}

\usepackage{slashed,fancyhdr,hojadeestilotesisdoctoral,bbding}

\usepackage{a4wide} 
\usepackage[english,greek]{babel}
\usepackage[usenames,dvipsnames,svgnames,table]{xcolor}
\usepackage[pdftex]{graphicx}
\usepackage{xparse}
\usepackage{collref} 
\usepackage{tikz}
\usetikzlibrary{matrix,arrows,decorations.pathmorphing}
\usepackage{tikz-cd}

\usepackage{float}
\restylefloat{table}
\usepackage[utf8]{inputenc}
\usepackage{indentfirst}
\usepackage{pb-diagram}
\usepackage{fancyhdr}
\usepackage{fancybox}

\usepackage{concmath}
\usepackage[T1]{fontenc}

\usepackage{enumitem}
\setlist[itemize]{leftmargin=*}
\setlist[enumerate]{leftmargin=*}

\usepackage{lipsum}
\makeatletter
\def\subsection{\@startsection{subsection}{2}%
	\z@{.5\linespacing\@plus.7\linespacing}{-.5em}%
	{\normalfont\bfseries}}
\makeatother

\makeatletter
\let\@afterindenttrue\@afterindentfalse
\@afterindentfalse
\makeatother

\usepackage{parskip}
\usepackage{etoolbox} 
\makeatletter
\patchcmd\@thm
{\let\thm@indent\indent}{\let\thm@indent\noindent}%
{}{}
\makeatother

\theoremstyle{plain}
\newtheorem{thm}{Theorem}[chapter]

\newtheorem{prop}[thm]{Proposition}
\newtheorem{lemma}[thm]{Lemma}
\newtheorem{cor}[thm]{Corollary}

\theoremstyle{definition}
\newtheorem{definition}[thm]{Definition}

\theoremstyle{remark}
\newtheorem{remark}[thm]{Remark}
\newtheorem{ep}[thm]{Example}

\def\mA{\mathrm{A}}
\def\mT{\mathrm{T}}

\newcommand{\End}{\mathrm{End}}
\newcommand{\Aut}{\mathrm{Aut}}

\newcommand{\ev}{ev}

\newcommand{\frt}{\mathfrak{t}}
\newcommand{\frF}{\mathfrak{F}}
\newcommand{\frK}{\mathfrak{K}}
\newcommand{\frG}{\mathfrak{G}}
\newcommand{\frg}{\mathfrak{g}}
\newcommand{\frm}{\mathfrak{m}}

\DeclareFontFamily{U}{rsf}{}
\DeclareFontShape{U}{rsf}{m}{n}{<5> <6> rsfs5 <7> <8> <9> rsfs7 <10-> rsfs10}{}
\DeclareMathAlphabet\Scr{U}{rsf}{m}{n}
\newcommand{\KA}{K\"{a}hler-Atiyah~}

\def\Z{\mathbb{Z}}

\def\R{\mathbb{R}}

\def\bS{\mathbb{S}}
\def\s{\mathbb{s}}

\def\s{\mathbb{s}}
\def\rk{{\rm rk}}

\def\Der{{\rm Der}}

\def\GL{\mathrm{GL}}
\def\Gl{\mathrm{Gl}}
\def\dd{\mathrm{d}}

\def\supp{\mathrm{supp}}

\def\l\Xi{\overrightarrow{\Xi}}
\def\r\Xi{\overleftarrow{\Xi}}

\def\l{\partial^l}

\def\Ad{\mathrm{Ad}}

\def\cQ{\mathcal{Q}}

\def\pr{\mathrm{pr}}

\def\Diff{\mathrm{Diff}}
\def\Conf{\mathrm{Conf}}
\def\mConf{\mathfrak{Conf}}
\def\Sol{\mathrm{Sol}}
\def\mSol{\mathfrak{Sol}}

\def\Met{\mathrm{Met}}

\newcommand{\id}{\mathrm{id}}
\newcommand{\Tr}{\mathrm{Tr}}
\newcommand{\tr}{\mathrm{tr}}

\newcommand{\diag}{\mathrm{diag}}
\newcommand{\sign}{\mathrm{sign}}

\def\cR{{\mathcal R}}

\def\cC{{\mathcal C}}

\def\cB{\Scr B}

\def\cH{\mathcal{H}}

\def\cZ{{cal Z}}

\def\Cl{\mathrm{Cl}}

\def\cK{\mathcal{K}}
\def\odd{\mathrm{odd}}
\def\Spin{\mathrm{Spin}}

\def\Pin{\mathrm{Pin}}
\def\Spin{\mathrm{Spin}}

\def\SO{\mathrm{SO}}
\def\O{\mathrm{O}}
\def\U{\mathrm{U}}

\def\cD{\mathcal{D}}

\def\cA{\mathcal{A}}
\def\cE{\mathcal{E}}
\def\cI{\mathcal{I}}
\def\cP{\mathcal{P}}

\def\cG{\mathcal{G}}
\def\cT{\mathcal{T}}
\def\cF{\mathcal{F}}
\def\cC{\mathcal{C}}

\def\Sp{\mathrm{Sp}}
\def\G_2{\mathrm{G_2}}

\def\cO{\mathcal{O}}
\def\cL{\mathcal{L}}

\def\cZ{\mathcal{Z}}

\def\s{\mathfrak{s}}

\def\P{\mathbb{P}}

\def\frD{\mathfrak{D}}
\def\cX{\mathcal{X}}
\def\cY{\mathcal{Y}}
\def\frc{\mathfrak{c}}
\def\frw{\mathfrak{w}}

\newcommand{\Hom}{{\rm Hom}}

\def\Aut{\mathrm{Aut}}
\def\Ob{\mathrm{Ob}}

\def\Re{\mathrm{Re}}
\def\Im{\mathrm{Im}}

\def\G{\mathrm{G}}

\def\R{\mathbb{R}}

\def\Pic{\mathrm{Pic}}

\def\cW{\mathcal{W}}
\def\cL{\mathcal{L}}
\def\cH{\mathcal{H}}
\def\dd{\mathrm{d}}

\def\fra{\mathfrak{a}}
\def\frq{\mathfrak{q}}

\def\fX{\mathfrak{X}}

\DeclareSymbolFont{bbold}{U}{bbold}{m}{n}
\DeclareSymbolFontAlphabet{\mathbbold}{bbold}

\def\fre{\mathfrak{e}}
\def\frv{\mathfrak{v}}

\newcommand{\Id}{\operatorname{Id}}

\newcommand{\Ker}{\operatorname{Ker}}

\newcommand{\tll}{\Theta_{ll}}
\newcommand{\tnn}{\Theta_{nn}}
\newcommand{\tln}{\Theta_{ln}}
\newcommand{\tun}{\Theta_{un}}
\newcommand{\tul}{\Theta_{ul}}
\newcommand{\tuu}{\Theta_{uu}}
\newcommand{\Bt}{\mathcal{B}_t}

\newcommand{\Ric}{\mathrm{Ric}}

\newcommand{\surj}{\to\kern-1.8ex\to}

\newcommand{\frf}{\mathfrak{f}}
\newcommand{\frk}{\mathfrak{k}}

\newcommand{\ric}{\mathrm{Ric}}

\newcommand{\arxiv}[1]{{\tt
		\href{http://www.arXiv.org/abs/#1}{arXiv:#1}}}

\newcommand{\frS}{\mathfrak{S}}

\NewDocumentCommand\mycite{mgggg}{\IfNoValueTF{#5}{\IfNoValueTF{#3}{\IfNoValueTF{#2}{\singlecite{#1}}{\singlecitedetail{#1}{#2}}}{\multicite{#1}{#2}{#3}}}{\multimulticite{#1}{#2}{#3}{#4}{#5}}} 
\usepackage[authoryear,round,square,colon,sort&compress,comma]{natbib}
\bibpunct{[}{]}{,}{n}{}{,}

\begin{document}
\setcounter{tocdepth}{1}
\selectlanguage{english}
\thispagestyle{empty}
\phantom{a}
\vspace{2cm}

\vspace*{-1.3cm}
\newcommand{\HRule}{\rule{\linewidth}{0.3mm}}
\setlength{\parindent}{1cm}
\setlength{\parskip}{1mm}
\noindent
\HRule
\begin{center}

\setlength{\baselineskip}{1\baselineskip}
\textbf{Torsion parallel spinors on Lorentzian four-manifolds \\ and supersymmetric evolution flows on bundle gerbes}
\HRule
\vspace{2cm}

\begin{minipage}{12cm}
\begin{center}
\bf{C. S. Shahbazi} \\
\end{center}
\end{minipage}

\vspace{0.5cm}

\begin{minipage}{12cm}
\begin{center}
\emph{Dissertation zur Erlangung des Doktorgrades des Fachbereichs Mathematik der Universit\"at Hamburg} \\ 
\end{center}
\end{minipage}

\vspace{2cm}


\includegraphics[scale=0.3]{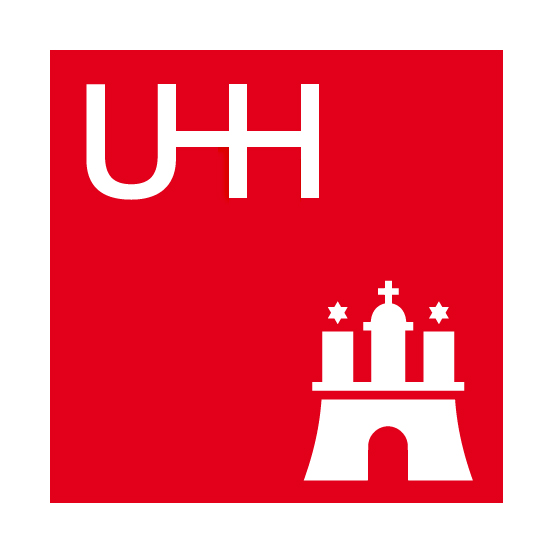}

\end{center}

\clearpage

\section*{Abstract}



This dissertation is concerned with the geometric investigation of differential spinors on oriented and spin Lorentzian four-manifolds via the theory of spinorial polyforms. The main results and applications are directed towards the study of torsion parallel spinors and the globally hyperbolic evolution flow determined by the globally hyperbolic solutions of the four-dimensional supersymmetric NS-NS system. This differential system, which originates in supergravity and string theory, involves skew-torsion parallel spinors subject to a curvature condition and provides a natural gauge-theoretic interpretation of skew-symmetric torsion as the curvature of a connection on an abelian bundle gerbe - a natural categorification of the notion of principal circle bundle.

\

\

\

\

\section*{Zusammenfassung}

Diese Arbeit befasst sich mit reelle und irreduziblen Spinoren auf einer orientierten und spin-Lorentzschen Viermannigfaltigkeit mit Blick auf Anwendungen in der global hyperbolischen Supersymmetrie und dem damit verbundenen Cauchy-Evolutionsfluss in vier Dimensionen.

\

\

\

\

\section*{Eigenanteilserklärung}

Hiermit erkläre ich an Eides statt, dass ich die vorliegende Dissertationsschrift selbst
verfasst und keine anderen als die angegebenen Quellen und Hilfsmittel benutzt
habe.

\begin{itemize}
	\item Kapitel \ref{chapter:spingeometryClifford} ist eine Adaption von \cite{Cortes:2019xmk} und \cite{SLSpin7}, mit gleichen Beiträgen aller Autoren.
	
	\item Der Unterabschnitt \ref{subsec:RKS} im Kapitel \ref{chapter:IrreducibleSpinors4d} wurde von \cite{RKS} übernommen, mit gleichen Beiträgen aller Autoren.
	
	\item Der Abschnitt \ref{sec:leftinvariant} im Kapitel \ref{chapter:Globallyhyperbolicsusy} wurde von \cite{Murcia:2021dur} übernommen, mit gleichen Beiträgen aller Autoren.
\end{itemize}

	
	
 
	
	

\vspace*{\fill}
\begin{itemize}
	\item Gutachter: Prof. Dr. Vicente Cortés
	\item Gutachter: Prof. Dr. Anton Galaev
	\item Gutachter: Prof. Dr. Thomas Leistner
	\item Vorsitz: Prof. Dr. Armin Iske 
	\item Stellvertretung: Prof. Dr. Jörg Teschner
	\item Schriftleitung: Prof. Dr. Ralf Holtkamp
\end{itemize}

\clearpage

\section*{Acknowledgements}

I would like to express my deepest gratitude to my supervisor, Professor Vicente Cortés, for his invaluable guidance and support ever since I joined the Mathematics Department of Hamburg University as a Humboldt Fellow back in 2017. His insightful advice, constructive feedback, and constant encouragement have been instrumental in my research. I am particularly grateful for his patience and dedication to mentoring me.

I am also deeply thankful for all the valuable interactions and collaborations that I have enjoyed over the years with my colleagues and excellent mathematicians Severin Bunk, Mario García-Fernández, Alejandro Gil, Calin Lazaroiu, Patrick Meessen, Andrei Moroianu, Vicente Muñoz, Tomás Ortín, and Marco Zambon, among others. Their insights and feedback have been invaluable to my research and career as a mathematician. A special thanks is reserved to Calin Lazaroiu for introducing me to the marvelous world of spinorial geometry, and to Ángel Murcia for embarking with me on the trip that led us to the world of \emph{parabolic pairs} and evolution flows.

Finally, I am profoundly grateful to my wonderful wife for her unwavering support and constant encouragement throughout this journey.

\thispagestyle{empty}
\phantom{1}


\addtocontents{toc}{\protect\thispagestyle{empty}}
\tableofcontents
\thispagestyle{empty}


\renewcommand{\thepage}{\arabic{page}}
\pagestyle{fancy}
\renewcommand{\leftmark}{Chapter \thechapter. Introduction}
\renewcommand{\headrulewidth}{0.0pt}
\fancyhead[LE]{\thepage}
\fancyhead[RE]{\slshape \nouppercase{\leftmark}}
\fancyhead[LO]{\slshape \nouppercase{\rightmark}}
\fancyhead[RO]{\thepage}
\fancyfoot{}
\renewcommand{\headrulewidth}{1pt}
\renewcommand{\footrulewidth}{1pt}


\renewcommand{\leftmark}{Introduction}

\chapter{Introduction}


This thesis is devoted to the study of real and irreducible spinors on an oriented and spin Lorentzian four-manifold with a view toward applications in globally hyperbolic supersymmetry and its associated Cauchy evolution flow. More specifically, let $\cC$ be an abelian bundle gerbe on an oriented and spin four-dimensional manifold $M$. The main purpose of this dissertation is to investigate the following coupled spinorial differential system on $(\cC,M)$:
\begin{eqnarray}
	\label{eq:introsystem}
	\nabla^{g,b}\varepsilon = 0 \, , \qquad H_b \cdot_g \varepsilon = \varphi_{\phi}\cdot_g \varepsilon
\end{eqnarray}

\noindent
with variables $(g,b,\phi,\varepsilon)$ consisting of a Lorentzian metric $g$ on $M$, a curving $b$ on $\cC$, a function $\phi\in C^{\infty}(Y)$ on $Y\to M$ satisfying an integrality condition on a certain submersion $Y\to M$ over $M$, and a section $\varepsilon\in \Gamma(S)$ of a bundle of irreducible real Clifford modules on $(M,g)$. Here $\varphi_{\phi}\in \Omega^1(M)$ and $H_b\in\Omega^3(M)$ respectively denote the \emph{curvatures} of $\phi$ and $b$, and $\nabla^{g,b}$ denotes the unique metric connection with completely skew-symmetric torsion $H_b$. This spinorial differential system exhibits multiple compelling features, of which we remark the following:
\begin{itemize}
	\item It defines a natural geometric system in both spin geometry and Lorentzian geometry with torsion. If both $\phi$ and $b$ are flat, then the system reduces to the parallelity of an irreducible real spinor under the Levi-Civita connection of the underlying Lorentzian manifold, which defines a classical problem in Lorentzian geometry and mathematical general relativity \cite{BergeryIkemakhen,EhlersKundt,Ikemakhen0,Leistner}.
	
	\item It involves choosing an underlying abelian bundle gerbe \cite{Murray}. Abelian bundle gerbes provide a natural categorification of the notion of principal U(1) bundle and define a \emph{geometric realization} of $H^3(M,\mathbb{Z})$, similarly to principal U(1) bundles determining a geometric realization of $H^2(M,\mathbb{Z})$. Consequently, Equation \eqref{eq:introsystem} yields an example of a geometric differential system in \emph{higher gauge theory} \cite{BaezSchreiber,BorstenEtAl24}.
	
	\item It gives a natural \emph{gauge theoretic} interpretation of pseudo-Riemannian skew-symmetric torsion as the curvature of a curving on a bundle gerbe. That is, in the system \eqref{eq:introsystem} the three-form torsion of $\nabla^{g,b}$ is not a variable in itself. Instead, it is a given as the curvature of a gauge-theoretic object that admits an interesting groupoid of automorphisms \cite{DjounvounaKrepski,MurrayStevenson,Waldorf}. 
	
	\item It couples higher gauge fields to spinors naturally through the canonical action of the curvature of the former on the latter via Clifford multiplication \cite{LawsonMichelsohn,Meinrenken}, as well as through the \emph{parallelicity} condition contained in \eqref{eq:introsystem}.
\end{itemize}

\noindent
The spinorial differential system \eqref{eq:introsystem} has its origins in the mathematical physics literature of supergravity and string theory \cite{BeckerBecker,FreedmanProeyen,GreenSchwarzWittenV1,GreenSchwarzWittenV2,Ortin,Tomasiello}: it encodes the \emph{supersymmetry conditions} of \emph{NS-NS supergravity} evaluated on a four-dimensional bosonic configuration. These supersymmetry conditions are typically called, as we will explain in more detail below, the \emph{Killing spinor equations} of the supergravity theory under consideration. NS-NS supergravity is the \emph{common sector} of ten-dimensional supergravity, and therefore plays a prominent role in the relationship between string theory and its low-energy limit, being a fundamental building block of the latter. As with any supergravity theory, the bosonic sector of NS-NS supergravity can be defined through its equations of motion, to which we will refer as the four-dimensional NS-NS system. This system is defined by the following second-order partial differential equations \cite[Equation (8.40)]{BeckerBecker}, \cite[Equations (2.2)]{Kawano:2003af}, \cite[Equation (21.1)]{Ortin}:
\begin{eqnarray*}
\label{eq:introNSNSsystem}
\Ric^{g,b} + \nabla^{g,b}\varphi_{\phi} = 0\, , \qquad \nabla^{g\ast}\varphi_{\phi} + \vert\varphi_{\phi}\vert^2_g = \vert H_b\vert^2_g
\end{eqnarray*}

\noindent
for triples $(g,b,\phi)$, where $\Ric^{g,b}$ is the Ricci tensor of $\nabla^{g,b}$ and $\nabla^g$ is the Levi-Civita connection of $g$.  A solution $(g,b,\phi)$ of the NS-NS system is \emph{supersymmetric} if and only if there exists an irreducible spinor bundle $S$ on $(M,g)$ and a spinor section $\varepsilon\in \Gamma(S)$ such that $(g,b,\phi,\varepsilon)$ is a solution of the NS-NS Killing spinor equations \eqref{eq:introsystem}. One of the most remarkable properties that makes the interaction between the NS-NS system and the NS-NS Killing spinor equations so mathematically rich is that the latter provides almost \emph{complete integrability} to the former and therefore defines a natural subclass of solutions which, based on extensive previous experience, can be expected to have remarkable mathematical properties regarding their stability, topology, moduli and applications to differential geometry and topology, see for instance \cite{Morgan} for an outstanding example of this phenomena. This dissertation initiates the systematic mathematical study of the supersymmetric NS-NS system. We begin by developing first a general geometric framework for investigating spinors that are parallel with respect to a general connection on the spinor bundle. Subsequently, we apply this framework to the study of Equation \eqref{eq:introsystem}. In more detail:

\begin{itemize}
	\item We begin in Chapter \ref{chapter:spingeometryClifford} by establishing a geometric setup to study spinorial differential equations of the form:
	\begin{eqnarray}
		\label{eq:diffspinorintro}
		\cD \varepsilon = 0 \, , \qquad \cQ(\varepsilon) = 0
	\end{eqnarray}
	
	\noindent
	where $\varepsilon\in \Gamma(S)$ is a section of the bundle of irreducible Clifford modules defined on a pseudo-Riemannian manifold $(M,g)$ of signature $(p,q)$ such that $(p-q) \equiv 0, 2 \mod(8)$. For lack of a better term, we will refer to such spinors as \emph{differential spinors}. Another possibility would be to refer to them simply as \emph{parallel spinors}, but that is usually reserved for spinors parallel for the Levi-Civita connection, or $\cD$-parallel spinors, which requires using a meaningless symbol, in this case $\cD$, as part of the definition.
	
	\item In Chapter \ref{chapter:IrreducibleSpinors4d} we apply the general framework developed in Chapter \ref{chapter:spingeometryClifford} to the case of differential spinors on four-dimensional Lorentzian four-manifolds. This covers in particular the case of every real and irreducible spinor parallel under the lift to the spinor bundle of any metric connection with torsion.
	
	\item We continue in Chapter \ref{chapter:parallelspinorstorsion} by applying the geometric framework developed in Chapter \ref{chapter:IrreducibleSpinors4d} to the particular case of real and irreducible spinors on Lorentzian four-manifolds equipped with a spinor parallel under a general metric connection with torsion, to which we refer as \emph{torsion parallel spinors}. Solutions to Equation \eqref{eq:introsystem} are particular instances of torsion parallel spinors.
	
	\item We proceed then with Chapter \ref{chapter:susyKundt4d}, in which we initiate the proper study of four-dimensional supersymmetric configurations and solutions, namely solutions of \eqref{eq:introsystem}. For this, we use the results and framework developed in the previous chapters.
	
	\item We culminate this dissertation in Chapter \ref{chapter:Globallyhyperbolicsusy}, which is devoted to the study of globally hyperbolic supersymmetric configurations and solutions building on the results obtained in previous chapters. 
\end{itemize}

\noindent
A detailed account of the contents and original results of each chapter can be found below in Section \ref{sec:mainresults}. The notion of differential spinor introduced in \eqref{eq:diffspinorintro} contains as particular cases most of the various types of special parallel spinors that have been considered in the literature, as well as some that are yet to appear or be discovered. This includes parallel spinors, Killing spinors, Codazzi spinors, Cauchy spinors, skew-torsion parallel spinors, skew-Killing spinors or generalized Killing spinors, just to name a few, see for instance \cite{Bar,BarGM,BaumIII,Bohle:2003abk,ContiDalmasso,ContiSalamon,FlamencourtMoroianu,FriedrichHyper,Friedrich:2001nh,FriedrichK,FriedrichKim,FriedrichKimII,Galaev,Kath,Ikemakhen,IkemakhenII,Leitner,MoroianuSemm,MoroianuSemmI,MoroianuSemmII} as well as their references and citations. NS-NS supergravity is a fundamental ingredient of ten-dimensional supergravity and superstring theory, and as such, it has been extensively studied in the physics and mathematical physics literature, specially coupled within more involved supergravities, such as Type IIA and Type IIB supergravity or Heterotic supergravity, see \cite{BeckerBecker,FreedmanProeyen,Figueroa-OFarrill:2003fkz,Figueroa-OFarrill:2003gow,Figueroa-OFarrill:2008lbd,Figueroa-OFarrill:2002ecq,Gauntlett:2002nw,Gauntlett:2002fz,Gauntlett:2003fk,Gauntlett:2003wb,Gran:2005wf,Gran:2018ijr,Grana:2004bg,Grana:2005sn,Legramandi:2018qkr,Ortin,Rosa:2013jja,Tomasiello:2007zq,Tomasiello:2011eb,Tomasiello} and their many references and citations. In particular, the supersymmetric solutions of the various supergravities that include NS-NS supergravity as a subsector have been widely studied in the theoretical physics and mathematical physics literature from diverse points of view, ranging from the very phenomenological \cite{Ellis:2013nka,Kallosh:2013yoa} to the very geometric \cite{Beckett24,Cecotti,Figueroa-OFarrill:2015rfh,Figueroa-OFarrill:2017tcy,Figueroa-OFarrill:2020gpr,Ortin,Tomasiello}. However, and interestingly enough, the study of supersymmetric solutions in Lorentzian signature has not \emph{transpired} into the mathematics literature and, to the best of my knowledge, this dissertation is the first systematic study of the class of Lorentzian supersymmetric configurations defined by Equation \eqref{eq:introsystem} from a pure mathematical point of view. Still more, spinorial geometry with torsion in Lorentzian signature is yet to be developed systematically in the mathematics literature, with some pioneering exceptions \cite{ErnstGalaev}. In striking contrast, equations \eqref{eq:introsystem} in Riemannian signature constitute one of the building blocks of the celebrated \emph{Hull-Strominger system} \cite{Garcia-Fernandez:2016azr,Hull,Strominger}, which was originally proposed by S. T. Yau \cite{FuTsengYau,LiYauI,LiYauII,Yau} as a natural generalization of the Calabi problem to non-Kähler complex geometry and has ever since evolved into a rapidly expanding area within complex geometry and geometric analysis. The main purpose of this thesis is to initiate the systematic mathematical study of this system in Lorentzian signature, which is the signature in which it was originally conceived as a supergravity theory, focusing on the geometry of its supersymmetric solutions and the evolution problem that they define in their globally hyperbolic regime. In the upcoming sections of this introductory chapter we will proceed to discuss the context and motivation of the mathematical study of \eqref{eq:introsystem} as well as more general supersymmetric systems that we plan to explore in the future. Just as \emph{mathematical} general relativity has become a well-established branch of mathematics, we hope the same will naturally happen for Lorentzian supergravity. Indeed, the application of modern geometric analysis to these theories is already proving fruitful, with seminal results on the stability of supersymmetric compactifications and Kaluza-Klein reductions \cite{Andersson:2020fuz,Huneau:2023cwa} paving the way for future work.


\section{The differential geometry of supergravity}


Supergravity theories are, by definition, gravitational theories invariant under a remarkable type of field-theoretic symmetry called \emph{supersymmetry}, which is characterized by the fact that it relates or \emph{maps} fields with integer spin to fields with fractional spin, and vice-versa. Supergravity theories were discovered in the mid-seventies \cite{DeserZumino,FreedmanNieuwenhuizen} and were consequently intensively investigated and developed in the physics literature, first as candidates for theories of everything and afterwards because of the prominent role they were discovered to play as the low-energy limits of superstring theory. This gave rise to a flurry of scientific activity that expanded and ramified into many different scientific areas in an ever-growing effort to understand and develop supergravity together with its applications to string theory and mathematics \cite{Cecotti,FreedmanProeyen,Ortin,Tomasiello}. 

Supergravity theories have been traditionally studied in the physics literature as local Lagrangian theories that are defined on a contractible open subset of $\mathbb{R}^n$ in terms of an explicit local Lagrangian function and its associated local equations of motion. Denote by $\cL[\Phi]$ a generic local supergravity Lagrangian function evaluated on a generic configuration element $\Phi$. That is, the symbol $\Phi$ represents the set of all variables of the given supergravity theory. In supergravity, and in contrast to other supersymmetric \emph{field} theories, this set of variables always includes a Lorentzian metric. Other common types of variables in $\Phi$ typically include connections on a principal bundle or sections of specific fibrations. A general phenomenon in supergravity is that each set of variables $\Phi$ can be divided into two disjoint subsets $\Phi = \left\{\Phi_b,\Phi_f\right\}$, defined as follows:
\begin{itemize}
	\item $\Phi_f$ consists of those variables in $\Phi$ that take values on a non-trivial spinor bundle defined on $(M,g)$ with $g\in\Phi$.
	
	\item $\Phi_b$ is the complement of $\Phi_f$ in $\Phi$, namely $\Phi_b$ consists of all those elements in $\Phi$ that do \emph{not} take values in any non-trivial spinor bundle.
\end{itemize} 

\noindent
Elements in $\Phi_f$ are called \emph{fermions}, \emph{fermionic fields} or \emph{variables}, whereas elements in $\Phi_b$ are called \emph{bosons}, \emph{bosonic fields} or \emph{bosonic variables}. Of course, this description is far away from the definition that a physicist would give for bosonic and fermionic fields in a given physical theory, but it is very convenient from a mathematical geometric perspective. An example of a bosonic variable that every supergravity contains is a Lorentzian metric, and an example of a fermionic variable that again every supergravity contains is a \emph{gravitino}, namely a one-form taking values on an appropriately chosen non-trivial spinor bundle. In fact, these are the variables of the simplest supergravity theory in four dimensions, sometimes called \emph{minimal supergravity} \cite[Chapter 5]{Ortin}. The defining property of a supergravity theory is that there exists a spinor bundle $S$ whose sections $\varepsilon\in \Gamma(S)$ generate \emph{infinitesimal} transformations $\delta_{\varepsilon} \Phi_b$ and $\delta_{\varepsilon} \Phi_f$ of $\left\{\Phi_b,\Phi_f\right\}$ that map the subsets $\Phi_b$ and $\Phi_f$ into each other (in an appropriate 
infinitesimal sense) and map solutions of the equations of motion to (infinitesimal) solutions. Schematically, we can write:
\begin{eqnarray*}
	\delta_{\varepsilon} \Phi_b = F_b(\varepsilon,\Phi_f)\, , \qquad \delta_{\varepsilon} \Phi_f = F_f(\varepsilon,\Phi_b)
\end{eqnarray*}

\noindent
where $F_b(-,-)$ and $F_f(-,-)$ are adequately chosen functions of the given entries that may contain derivatives of $\Phi$ and are always linear in the first entry. These are the \emph{infinitesimal supersymmetry transformations} of the given supergravity theory. Since they are generated by a spinorial parameter, they can be understood as the spinorial analog of the classical infinitesimal symmetries generated by vector fields, such as the infinitesimal transformation of a metric generated by a vector field, which is given by the Lie derivative of the former along the latter. As it happens for every theory defined via a differential system of equations that admits non-trivial symmetry transformations, particular solutions of the theory need not be preserved by the symmetries of the equations of motion of the theory. When this is the case for the aforementioned supersymmetry transformations then we talk about \emph{supersymmetric configurations} and \emph{supersymmetric solutions}. That is, a supersymmetric configuration is a configuration element $\Phi$ for which there exists a spinor $\varepsilon$ such that:
\begin{eqnarray*}
	\delta_{\varepsilon} \Phi_b = F_b(\varepsilon,\Phi_f) = 0\, , \qquad \delta_{\varepsilon} \Phi_f = F_f(\varepsilon,\Phi_b) = 0
\end{eqnarray*} 

\noindent
If in addition, $\Phi$ is a solution to the equations of motion of the theory then it is by definition a supersymmetric solution of the given supergravity. The supersymmetric solutions of a supergravity theory define a remarkable class of solutions that are expected to enjoy special properties regarding their stability, regularity, and topological behavior, and the extensive expertise gained by physicists in the last decades certainly provides strong evidence for this. There is however no proper and well-established mathematical theory of supergravity that can be used to investigate supersymmetric solutions in supergravity and their potential applications in mathematics. In fact, such a complete mathematical theory of supergravity is probably currently out of reach. A mathematical theory of supergravity would require at least implementing the following two points:

\begin{itemize}
	\item It should be developed in the appropriate framework of superspace or supergeometry \cite{Varadarajan}, in terms of which all the mathematical structures of a given supergravity, in particular its supersymmetry transformations, should occur naturally.
	
	\item It needs to be formulated on manifold/supermanifolds of general enough topology, certainly not necessarily contractible or parallelizable, equipped with appropriate bundles and their categorifications that allow describing supergravity in terms of globally defined differential operators acting on intrinsically defined infinite-dimensional spaces. This is the way mathematicians have successfully formalized other physical theories such as Yang-Mills theory or the Seiberg-Witten equations, leading to spectacular mathematical new results through a careful study of their moduli spaces \cite{KronheimerDonaldson,Morgan}.
\end{itemize}

\noindent
Both these requirements seem to be very far from being completely established in supergravity. Whereas supergeometry is being actively developed as a mathematical discipline with a very algebraic-geometric flavor \cite{Carmeli,Rogers}, its applications to the mathematical theory of Lorentzian supergravity beyond two dimensions seem to be for the moment limited, with some pioneering exceptions \cite{Liu:2018xtm,Liu:2019tet,Liu:2020kds}. This does not prevent physicists from efficiently using \emph{superspace methods} to study various aspects of supergravity or construct new supersymmetric actions \cite{FreedmanProeyen,Howe:1981gz}. In fact, these superspace methods are of the utmost importance to develop supergravity and understand it properly as full theory coupling bosons to fermions through supersymmetry. On the other hand, as a result of supergravity theories having been traditionally studied only locally, there is little emphasis on understanding the global topological and geometric structure of supergravity, and many supergravity theories are typically constructed in terms of mathematical structures that are in general only well-defined locally, such as Kähler potentials, moduli-space coframes or coset representatives. When necessary, physicists typically perform case-by-case gluing procedures, some of them remarkably sophisticated, to understand the global geometric and topological structure of specific classes of supergravity solutions, see for instance the seminal references \cite{Ferrero:2020laf,Ferrero:2020twa}. 

These remarks may point to supergravity being not mature enough to be studied as a mathematical discipline in the same sense as General Relativity, Yang-Mills theory or Seiberg-Witten theory have all become exceptionally successful mathematical disciplines and areas of study by themselves. However, we shall not despair: the physical purpose and motivation of supergravity point to a proper way of studying it mathematically. For this, we have to reflect on what is the ultimate purpose of supergravity and how it is used in the theoretical physics literature. In this context, supergravity is considered a fundamental theory of gravity that describes the \emph{macroscopic} gravitational interaction of fundamental string theoretic objects \cite{Ortin}. That is, supergravity is understood as a gravitational theory that describes the gravitational interaction of fundamental objects, such as \emph{black branes}, \emph{black holes} or \emph{gravitational waves}. In other words, supergravity describes, as a general relativistic theory, the macroscopic gravitational interaction in string theory. This has two immediate consequences: 
\begin{itemize}
	\item Supergravity should be considered as a very specific class of general relativistic models, determined by supersymmetry, and therefore its mathematical study should be approached analogously to the mathematical study of general relativity, that is, from the point of view of global Lorentzian geometry and geometric analysis.
	
	\item Since no \emph{fermionic degrees of freedom} are observed macroscopically, we can safely truncate the fermionic sector of supergravity and study instead its bosonic sector together with the supersymmetry conditions that remain after this truncation.
\end{itemize}

\noindent
From this perspective, and since the fermionic sector is truncated from the onset, we would only need to develop the mathematical theory of bosonic supergravity together with the remnant supersymmetry conditions. This is a task that does not require supergeometry and can be accomplished within \emph{standard} global differential geometry and geometric analysis in terms of manifolds, bundles, metrics, connections, and possibly their categorifications. Following standard usage in the literature, we will refer to the bosonic truncation of a supergravity theory simply as a \emph{bosonic supergravity}. This is the perspective adopted in \cite{Lazaroiu:2016iav,Lazaroiu:2017qyr,LSProc,Lazaroiu:2021iwd,Lazaroiu:2021vmb}. Although bosonic supergravity can have energy momentum tensors involving remarkable mathematical structures, ultimately they are nothing but very specific matter models in general relativity and are therefore formulated in terms of highly coupled second-order differential equations. The key point and the source of many of the mathematical wonders of supersymmetry and supergravity is that, even after truncating the fermionic sector, the supersymmetry transformations of a bosonic configuration remain non-trivial and define gauge theoretic first-order systems of equations of remarkable mathematical depth and applications. To see this more explicitly, let $(\Phi_b, \Phi_f = 0)$ be a bosonic configuration, namely a configuration $\Phi$ with its fermionic sector set to zero. We can still apply to it a supersymmetry transformation, obtaining:
\begin{eqnarray*}
	\delta_{\varepsilon} \Phi_b = F_b(\varepsilon,0)\, , \qquad \delta_{\varepsilon} 0 = F_f(\varepsilon,\Phi_b)	 
\end{eqnarray*} 

\noindent
Since $ F_b(\varepsilon,0) = 0$ identically, it follows that $\delta_{\varepsilon} \Phi_b = 0$ automatically. However, even though the fermionic variables have been set to zero, their infinitesimal variation by a supersymmetry transformation may be non-vanishing. Hence, since the \emph{value zero} of the fermions needs to be preserved by supersymmetry in order for $(\Phi_b,0)$ to be a supersymmetric solution, $(\Phi_b,0)$ is supersymmetric if and only if:
\begin{equation*}
	F_f(\varepsilon,\Phi_b)	 = 0
\end{equation*}

\noindent
This equation becomes a first-order differential equation for $\Phi_b$ and $\varepsilon$ which is linear in $\varepsilon$ but typically highly non-linear in $\Phi$. These are the \emph{supersymmetry Killing spinor equations}, or Killing spinor equations for short, of supergravity and constitute the source of a myriad of applications in physics and mathematics. Indeed, particular cases of this equation include the pseudo-holomorphicity equations of Gromov \cite{Gromov}, the $\G_2$ and $\Spin(7)$ instanton equations \cite{DonaldsonThomas}, of several Seiberg-Witten-like differential systems \cite{FineGhosh}, just to name a few. This is arguably one of the deepest contributions of supersymmetry and supergravity to mathematics and has led to a flurry of mathematical activity since it provides a plethora of techniques to deal with classes of solutions to coupled second-order differential equations whose moduli spaces have had extensive applications in differential geometry and topology. At the source of these many applications lies the fact that the Killing spinor equations of a supergravity theory provide \emph{partial first-order integrability} of the full second-order equations of the theory. Indeed, solving the Killing spinor equations on general configuration elements solves automatically \emph{most} of the second-order equations of motion of the theory, and in some cases all of them. The framework of supersymmetric configurations that we have just described has been exploited masterfully in \cite{TsengYauI,TsengYauII,TsengYauIII} and \cite{FeiGuo,FeiI,FeiII,Fei:2021yjr,GFJS23,GFRT24}, see also \cite{Phong,Picard:2024tqg,PhongPicardZhang}, to develop, among other applications, new cohomology theories and evolution flows in Hermitian complex geometry, especially in complex non-Kähler geometry. In this dissertation we will discuss the mathematical structure of supersymmetric solutions and configurations in the context of the NS-NS supergravity, that is, the NS-NS system as we have defined it in Equation \eqref{eq:introNSNSsystem}. By the previous discussion, we will only need to present the mathematical formulation of the bosonic sector of NS-NS supergravity together with its Killing spinor equations.


\section{Categorified geometry and supersymmetry}
\label{sec:categorified}


Mathematical gauge theory \cite{KronheimerDonaldson} is a mathematical discipline that grew out of the careful investigation of the Yang-Mills equations and instanton equations on a principal bundle, and has ever since evolved into a mathematical area that is concerned with the mathematical applications of the study of differential equations of \emph{gauge theoretic type}, namely equations that are formulated in terms of globally defined differential operators on principal bundles, spinor bundles and various fibrations and have as variables, sections, connections and spinors that enjoy a gauge principle. In mathematical gauge theory, there is typically an underlying rich \emph{infinite dimensional} automorphism group that maps solutions to solutions and plays a fundamental role in the study of the corresponding moduli space of solutions, since it is used to identify \emph{gauge-equivalent solutions}. The differential geometric study of bosonic supergravity together with its Killing spinor equations, which we have outlined and motivated in the previous section, fits perfectly within this framework with only two caveats:

\begin{enumerate}
	\item Some supergravity theories, especially those defined in \emph{higher dimensions}, require the use of \emph{categorified} notions of principal bundles and connections. The occurrence of higher category theory in mathematical gauge theory is not new, see for instance \cite{BehrendI,BehrendII}, although the way it specifically occurs as a consequence of supersymmetry is not fully understood yet.
	
	\item Whereas mathematical gauge theory usually deals with systems of equations of elliptic type, supergravity theories are Lorentzian gravitational theories and therefore are governed by differential systems of hyperbolic type. As we discuss in Section \ref{sec:susyevolutioninitial} below,  there is however an elegant canonical procedure to reduce supergravity to a system of \emph{Riemannian} equations that we conjecture is of elliptic type. The celebrated Hull-Strominger \cite{Hull,Strominger} system as well as other flavors of \emph{Heterotic systems} \cite{Moroianu:2021kit,Murcia:2022iba,Picard:2024tqg} would be particular instances of this procedure.  
\end{enumerate}

\noindent
The need for \emph{higher} or \emph{categorified geometry} in supergravity has an arguably very innocent origin. As mentioned in the previous sections, physicists typically study supergravity theories locally in terms of locally defined fields. The type of fields that can occur as part of a given supergravity theory is highly constrained by supersymmetry, sometimes to the point of being uniquely determined by supersymmetry at the local level, as it happens for instance with the maximally supersymmetric supergravities. Physicists quickly realized that some supergravity theories required for their consistency \emph{gauge fields} that did not correspond to local one-forms taking values on a Lie algebra, which were the type of gauge fields that physicists had studied up to that point in the context of Yang-Mills theories. Instead, imposing supersymmetry required the existence of fields that were locally represented by higher-order differential forms. These forms enjoyed a gauge principle and admitted a notion of curvature, similar to that of the traditional one-form gauge fields, albeit being of \emph{higher order}. At the local level, it is possible to work with these higher gauge fields by keeping track of the hierarchy of gauge transformations that they involve and the higher curvatures they define, which are sometimes coupled to each other, giving rise to deep algebraic structures \cite{Bonezzi,Borsten:2021ljb,BorstenEtAl24,Cagnacci:2018buk,deWit:2008ta,Giotopoulos:2024jcr,Giotopoulos:2024ynd,Hartong:2009az,Saemann,Trigiante:2016mnt}. A completely different story is trying to elucidate the global mathematical structure behind these higher gauge fields. By this, we refer to constructing the type of geometric objects on which these higher gauge fields can be defined as \emph{connections} involving the correct notion of curvature and gauge transformations. Doing this turns out to be a surprisingly sophisticated mathematical problem that requires some of the latest developments in the categorification of principal bundles, that is, in \emph{higher geometry} and is currently being intensively investigated \cite{Aschieri:2003mw,Bunk:2020wla,Bunk:2021quu,Bunk:2023jsj,BunkShahbazi,KristelLudewigWaldorf,NikolausSchreiberStevenson,NikolausWaldorf,Tellez,WaldorfII,SatiSchreiberEquivariant,Waldorf2018,WaldorfTdualCat,WaldorfTdualGeo,WaldorfLoops}. Therefore, from this simple minded point of view we can summarize the situation as follows:

\begin{itemize}
\item Certain supergravity theories require to be formulated on the categorification of traditional principal bundles and connections because, as a consequence of supersymmetry, they involve gauge fields that are locally represented as a $k$-forms with $k>1$.
\end{itemize}

\noindent
The study of gauge theories on various versions of categorified principal bundles has a rich history in the mathematical physics literature in various developments a priori unrelated to supergravity, see for instance \cite{BaezSchreiber,Schreiber} and its references and citations. In this dissertation we will be dealing with the simplest type of higher gauge field that occurs in supergravity, namely the so-called \emph{b-field}. This is a gauge field that is locally represented by a local two-form on the underlying manifold and consequently has a three-form as curvature. The geometric objects on which the $b$-field \emph{lives} are bundle gerbes with a connective structure \cite{Murray}, which define a natural geometric model for the categorification of a principal bundle. Consequently, we can think of the $b$-field as the higher analog of the \emph{Maxwell field}, which is mathematically modelled by a connection one-form on a principal U(1) bundle. Bundle gerbes were originally conceived in 1971 by Giraud \cite{Giraud} as certain sheaves of groupoids that, modulo the appropriate notion of isomorphism, provide a \emph{geometric realization} of $H^3(M,\mathbb{Z})$ akin to the geometric realization of $H^2(M,\mathbb{Z})$ given by principal $\U(1)$ bundles on $M$ modulo isomorphism. The original notion of gerbe as introduced by Giraud was relatively abstract and not particularly apt to do differential Riemannian or Lorentzian geometry \emph{on it}. Subsequently, more geometrically palatable models for the notion of abelian gerbes, all essentially equivalent, were proposed in the literature \cite{Brylinski}, culminating in the very geometric notion of \emph{abelian bundle gerbe} \cite{Murray}, which is a model constructed out of standard U(1) bundles and smooth fibrations. This makes Abelian bundle gerbes especially convenient for the supergravity applications we have in mind. One of the key properties of gerbes and other categorified geometric structures is that their \emph{natural symmetries} do not assemble a group but rather a categorified version of a group \cite{Murray,MurrayStevenson,Waldorf}, which reflects in turn the \emph{gauge symmetries} of the underlying supergravity. This is particularly important to study moduli spaces of supergravity solutions since such spaces are defined by identifying solutions via the higher symmetries of the theory. 

In this dissertation, we will propose a completely rigorous formulation of bosonic NS-NS supergravity on a bundle gerbe. This is a very natural formulation that offers no difficulties, in contrast to more complicated supergravity theories that involve \emph{non-abelian} higher geometry \cite{NikolausWaldorf,WaldorfII}. Nevertheless, it is important to keep track of the fact that the underlying variable of the theory is a curving $b$ on an abelian bundle gerbe and not a generic closed three-form $H$. Omitting this fact would be as erroneous as taking a Lie-algebra valued two-form as the variable in Yang-Mills theory, instead of the corresponding connection on a principal bundle. This is even \emph{experimentally incorrect} \cite{Aharonov}. Although the equations of motion of NS-NS supergravity only contain the $b$-field through its curvature, the importance of $H$ being a derived object becomes apparent as soon as one studies more subtle problems, such as the Cauchy and evolution problems or the supersymmetric moduli problem. Indeed, it is in this context where the fact that the variable is a gauge theoretic object with a rich groupoid of symmetries is apparent. We will encounter such a scenario in Chapter \ref{chapter:Globallyhyperbolicsusy}, where we will \emph{decompose} a curving in order to obtain the evolution flow defined by the NS-NS system.


\section{Spinorial exterior forms and the bilinear method}


The most important geometric tool in this dissertation is the theory of \emph{spinorial exterior forms}, developed originally in \cite{Cortes:2019xmk} for irreducible real spinors of simple type and extended to the non-simple case in \cite{Shahbazi3d}. We will abundantly apply the results of \cite{Cortes:2019xmk} to study irreducible differential spinors in four dimensions. This formalism covers as a particular case the spinorial differential system presented in \eqref{eq:introsystem}. The theory of spinorial exterior forms originated in the early supergravity literature, as physicists tried to understand the geometric structure of supersymmetric configurations and solutions, and hence were led to study highly coupled spinorial systems of the form \eqref{eq:diffspinorintro}. The systematic investigation of supersymmetric solutions in supergravity was initiated in the seminal work of Gibbons, Hull, and Tod \cite{Gibbons:1982fy,Tod:1983pm,Tod:1995jf}.  These pioneering efforts culminated in the early two thousand in what is nowadays called the \emph{spinorial bilinear method} \cite{Gauntlett:2002nw,Gauntlett:2003fk,Gauntlett:2003wb,Bellorin:2005zc}, which was consequently broadly applied to the study of the local structure of the supersymmetric solutions of many supergravity theories of increasing complexity, see \cite{Gran:2018ijr,Ortin} for general reviews on supersymmetric supergravity solutions from a physics but very geometric perspective. This bilinear method gives a systematic procedure to study the local geometric structure of supersymmetric supergravity solutions and more generally of differential spinors, by means of an equivalent associated differential system for a certain number of algebraically constrained differential forms. In essence, the bilinear method constructs a series of local differential forms by \emph{sandwiching} spinors with \emph{gamma} matrices, and then derives a number of differential and algebraic relations satisfied by these forms by applying the \emph{Fierz identities} \cite[Appendix D]{Ortin} on the one hand and the differential equations satisfied originally by the spinors on the other. These sandwiches of spinors and gamma matrices are called the \emph{bilinears} of the spinors, whence the name of the formalism. Hence, given a differential spinor, the bilinear method produces a set of forms that satisfies a system of algebraic relations derived by means of the Fierz identities and a system of differential conditions derived from the differential conditions satisfied by the differential spinor. The main advantage of this method is that it provides an algebraic procedure that does not depend on any representation theory, which for spinors on a pseudo-Riemannian signature can be remarkably non-trivial, and that it associates to every differential spinor a number of equivalent special differential forms that are in principle much more transparent to interpret and understand. This method directly relates to one of the two main challenges we encountered when trying to develop the mathematical theory of supersymmetric solutions in supergravity, namely:
\begin{enumerate}
	\item Identify mathematically the type of spinors that occur as supersymmetry generators of supersymmetric solutions in supergravity.
	
	\item Elucidate a geometric framework general enough to study systematically the type of spinorial equations that occur as supersymmetry conditions in \emph{Lorentzian supergravity}.
\end{enumerate}

\noindent
Point (1) above leads to the notion of \emph{differential spinor} introduced earlier together with the notion of \emph{spinorial Lipschitz structure} \cite{FriedrichTrautman,Lazaroiu:2016iav,LS2018}, which we shall not discuss here. Point (2) above led us to develop a general geometric framework to study equations of the type \eqref{eq:diffspinorintro}. The bilinear method stands as the perfect starting point for this since it has extensively proved its power in the literature and, as mentioned below, it addressed a fundamental problem in spin geometry, namely the relation between spinors and exterior forms \footnote{\emph{Spinors as the Square Root of Geometry}, conference given by Michael Atiyah at the IHES in 2013.}, in a manner that could be implemented globally with no reference to preferred choices of coordinates or coframes. In order to promote the remarkable method of bilinears into a proper mathematical theory, two challenges need to be addressed, that is:

\begin{itemize}
	\item We need to identify how the individual \emph{bilinears} occur or assemble into a natural geometric object associated to a section of a spinor bundle equipped with a bilinear pairing \cite{AC,ACDP}. This natural geometric object is the \emph{spinorial exterior form} associated to the given spinor.
	
	\item We need to obtain a precise correspondence between the differential system of type \eqref{eq:diffspinorintro} satisfied by a given spinor, and an equivalent differential system satisfied by its associated \emph{bilinears}.
\end{itemize}

\noindent
The aforementioned theory of spinorial exterior forms precisely solves these two points. First, it identifies the set of bilinears associated with a given spinor with its \emph{algebraic square}, appropriately defined via the choice of an \emph{admissible bilinear pairing} \cite{AC,ACDP} on the spinor space. Second, it establishes a correspondence between differential spinors, namely solutions to systems of the form \eqref{eq:diffspinorintro}, and their algebraic squares as solutions to a prescribed exterior differential system. In other words, this formalism reduces the original spinorial system to an exterior differential system for the square of a spinor, appropriately defined. Hence, in order to obtain a proper equivalence, we need to characterize which exterior forms can be constructed as the square of a spinor, that is, we need to characterize the image of the \emph{spinor square map}, which associates to every spinor its \emph{square} with respect to a chosen admissible bilinear pairing. In signature $(p-q) \equiv_8 0,1,2$, which crucially includes the case of signature $(3,1)$, this characterization is one of the main results of \cite{Cortes:2019xmk,Shahbazi3d}. In Op. Cit. the image of the spinor square map is characterized as a certain semi-algebraic body given as the solution set to a natural system of algebraic equations and inequalities elegantly written in terms of the \emph{geometric product}. The latter is the metric-dependent deformation of the wedge product that makes the exterior algebra bundle canonically isomorphic with the bundle of Clifford algebras. We note that the theory of spinorial exterior forms for irreducible real spinors in signatures different than $(p-q) \equiv_8 0,1,2$ is yet to be developed. In Lorentzian signature, the square of a spinor always contains the \emph{Dirac current} of the spinor, which is the rank-one component of its square and is therefore a one-form that is canonically associated with every irreducible spinor. In dimension larger than three the square of the spinor also contains forms of higher degree and therefore the Dirac current of the spinor does not encode all the information contained in the spinor. 

In this dissertation we will apply the theory of spinorial exterior forms to study differential spinors in four dimensions. For this particular case, the square of an irreducible real spinor can be solved explicitly, giving rise to the notions of \emph{parabolic pair} and \emph{isotropic parallelism}, see Section \ref{sec:4dAlgebraic} in Chapter \ref{chapter:IrreducibleSpinors4d} for more details. These fundamental objects give a very efficient way to deal with irreducible spinors in four Lorentzian dimensions. We will use them to study the four-dimensional supersymmetric solutions of NS-NS supergravity. The most basic property of nowhere vanishing irreducible spinors in four Lorentzian dimensions is that their Dirac current is necessarily isotropic, that is, of vanishing pointwise pseudo-norm. Hence, Lorentzian four-manifolds equipped with a differential spinor are necessarily equipped with a nowhere-vanishing isotropic vector satisfying a given differential equation. For irreducible real spinors parallel with respect the Levi-Civita connection this isotropic vector field is also parallel, and therefore the underlying Lorentzian manifold is \emph{Brinkmann} \cite{Brinkmann,Stephani}. The Brinkmann class of Lorentzian manifolds has generated ample interest in the mathematics literature since it encompasses a wide variety of interesting mathematical problems of topological and dynamical type \cite{Hanounah,Leistner,LeistnerSchliebner,MehidiZeghib}. Parallel spinors occurred relatively early in the mathematical general relativity literature \cite{EhlersKundt}. In modern terminology, in Op. Cit. it was shown that a Lorentzian four-manifold equipped with an irreducible real spinor parallel with respect to Levi-Civita is necessarily a particular type of Brinkmann space-time called \emph{pp-wave}, which stands for \emph{plane-fronted wave} with \emph{parallel propagation}. This name stems from the fact that they carry an isotropic vector field parallel with respect to the Levi-Civita connection and there is a local transverse Riemannian submanifold, encoding the local geometry of the front-wave of the wave, which is necessarily flat and thus \emph{plane-fronted}. This vector field is precisely the Dirac current of the parallel spinor and its dual belongs to the square of the latter. Such pp-waves define, from the physical point of view, idealized mathematical models for gravitational waves. This fact together with the spectacular observation of gravitational waves in the LIGO experiment makes this class of space-times of the utmost relevance in modern theoretical physics. From the mathematical point of view, generic Brinkmann space-times define Lorentzian manifolds of \emph{reducible} but \emph{non-decomposable} holonomy \cite{Leistner0I,Leistner0II,Leistner,LeistnerSchliebner}. That is, the distribution spanned by the parallel vector field is preserved by the Levi-Civita connection, yet the underlying Lorentzian metric does not split as a local product. This is a fascinating possibility that cannot occur in Riemannian signature by the de Rahm theorem. From the spinorial point of view, as soon as a differential spinor is isotropic with respect to the natural pseudo-norm on the spinor space then its Dirac current is automatically isotropic, nowhere vanishing, and satisfies a given differential equation. Hence, the relation between isotropic spinors in Lorentzian signature and isotropic vector fields runs deep in Lorentzian spinorial geometry. Beyond spinors parallel with respect to the Levi-Civita connection, previous experience \cite{AlonsoAlberca:2002gh,Brannlund,Shahbazi3d}, shows that supersymmetric solutions are in particular Kundt space-times \cite{Boucetta:2022vny,Coley,MehidiZeghibII}, at least in dimensions three and four, giving rise to various notions of \emph{supersymmetric} Kundt manifolds. The Kundt class of Lorentzian manifolds is a remarkable extension of the class of Brinkmann manifolds and is being actively investigated in the mathematics community \cite{MehidiZeghibII}. We hope the methods and techniques developed in the aforementioned references can be applied in the future to study in detail the mathematical structure of the supersymmetric classes of Kundt manifolds occurring in supergravity.


\section{Supersymmetric evolution flows and initial data sets}
\label{sec:susyevolutioninitial}


The mathematical study of supergravity theories in the mathematical literature has been thus far mostly limited to Euclidean signature, with some pioneering exceptions \cite{Andersson:2020fuz}. In the most celebrated supergravity system in the mathematics literature, namely, the supersymmetric Killing spinor equations of Heterotic supergravity, also known as the \emph{Hull-Strominger system} \cite{Garcia-Fernandez:2016azr}, this mathematical study is almost entirely done not only in Euclidean signature but also on compact manifolds. This is justified by the fact that solutions to the Hull-Strominger correspond to \emph{compactification backgrounds} of Heterotic supergravity \cite{BeckerBecker}. A more general but equally natural class of manifolds on which studying the Hull-Strominger system as well as other \emph{gauge-theoretic equations} in Euclidean signature are complete Riemannian manifolds. When addressing the mathematical study of supergravity in Lorentzian signature, it is convenient to consider the theory as defined on a particular class of manifolds akin to the complete or compact manifolds in Euclidean signature. Studying compact Lorentzian manifolds is usually avoided in the physics literature since these seem to be necessarily non-realistic, in the sense that they contain closed time-like curves giving rise to all sorts of physical paradoxes. On the other hand, imposing Lorentzian manifolds to be \emph{complete}, for instance, geodesically complete, is in general too strict, since many classes of remarkable gravitational space-times are known to be non-complete. In fact, there are \emph{singularity theorems} that, under certain conditions, predict the generic occurrence of singularities in the form of geodesic incompleteness with respect to the underlying Lorentzian metric \cite{SenovillaI,SenovillaII}. There is however, a natural class of Lorentzian manifolds that constitute the analog of complete Riemannian manifolds in Euclidean signature: these are the globally hyperbolic Lorentzian manifolds with complete Cauchy hypersurface. Similarly, the analog of compact Riemannian manifolds in Euclidean signature would correspond from this point of view with the globally hyperbolic Lorentzian manifolds with \emph{compact} Cauchy hypersurface. 

Globally hyperbolic Lorentzian manifolds can be defined in a number of different ways. From a practical point of view, globally hyperbolic Lorentzian manifolds constitute the class of Lorentzian manifolds on which the Cauchy problem for hyperbolic systems of equations can be well posed. In particular, a globally hyperbolic Lorentzian manifold admits a space-like submanifold on which we can define the initial data associated to any given hyperbolic system on $(M,g)$, and copies of this space-like submanifold foliate the entire space-time. Hence, globally hyperbolic Lorentzian manifolds are especially well-adapted to study analytic aspects of hyperbolic systems as well as to understand hyperbolic systems as evolution systems in Riemannian signature. The celebrated results of Choquet-Bruhat \cite{Choquet-BruhatI,Choquet-BruhatII} show that the vacuum equations of General Relativity define a well-posed Cauchy problem for globally hyperbolic Lorentzian metrics. This allowed, for the first time in the literature, for a precise understanding of the \emph{evolution problem} defined by General Relativity: in contrast to elliptic systems in Riemannian signature, hyperbolic-type systems in Lorentzian signatures admit a \emph{time-direction} with respect to which they can be understood as dynamical systems for certain quantities, typically tensors but not necessarily, that evolve in time. This brought the understanding of General Relativity on par with the common understanding of all the other theories in physics that are determined in terms of evolution equations for an explicit time parameter. In particular, Choquet-Bruhat established the modern paradigm to understand the evolution and Cauchy problem of certain second-order hyperbolic systems that define \emph{gravitational theories}. These are systems of equations that are obtained from diffeomorphism-invariant local Lagrangians depending on a Lorentzian metric and possibly other fields and giving rise to second-order equations of motion. Within this paradigm, these gravitational theories can be understood, in their globally hyperbolic regime, in terms of a system of evolution equations in \emph{Euclidean signature} for certain \emph{constrained} initial data that satisfies a Riemannian system of equations. In a sense, gravitational theories on globally hyperbolic Lorentzian metrics effectively reduce to a system of evolution flow equations in Riemannian signature with constrained initial data. Given these remarks and the fact that supergravity is a \emph{landmark} of a gravitational theory, it is natural to study the supergravity initial value or Cauchy problem on globally hyperbolic Lorentzian configurations, namely configurations whose associated metric is globally hyperbolic. However, it seems that this problem has not been considered systematically in the supergravity literature, with the remarkable exception of \cite{Choquet-BruhatIII}, in which Choquet-Bruhat herself proved well-posedness of the complete, second-order, eleven-dimensional supergravity. It should be noted that reference \cite{Choquet-BruhatIII} deals with the local Cauchy problem, meaning local both in \emph{space} and \emph{time} at the second-order derivative level. This, which may seem a minor restriction, is remarkably non-trivial for eleven-dimensional supergravity for two main reasons: 

\begin{itemize}
	\item The proper formulation of eleven-dimensional supergravity unavoidably involves higher-order curvature corrections that are expected to completely spoil the nice analytic properties of second-order equations of motion of eleven-dimensional supergravity \cite{BeckerBecker}.
	
	\item The global geometric and topological structure of the eleven-dimensional supergravity requires a quite sophisticated notion of bundle 2-gerbe \cite{Stevenson} whose higher groupoid of automorphisms needs to be properly understood in order to \emph{glue} the local well-posedness result into a global one. This is related to a subtle gauge quantization condition in M-theory, whose mathematical foundations \cite{Giotopoulos:2024xcg,Sati:2021uhj,SatiSchreiber} are still under development and require advanced tools in algebraic topology.
\end{itemize}

\noindent
The previous two bullet points are of common occurrence in many other supergravity theories, in the sense that they involve, as mentioned in Section \ref{sec:categorified} non-trivial topological and geometric structures in delicate equilibrium as dictated by supersymmetry together with categorified notions of automorphisms that need to be properly taken into account in order to obtain global results. Notably, it seems that the \emph{constraint equations} of many supergravity theories have not even been written explicitly in the literature. We fill this gap in the particular case of NS-NS supergravity in Chapter \ref{chapter:Globallyhyperbolicsusy}, as we explain in more detail in Section \ref{sec:mainresults}. It should be noted that the Cauchy problem in supergravity is richer than in other more \emph{standard} gravitational theories since a given supergravity involves two different systems of equations. As we have discussed earlier, a supergravity theory, by which we mean a \emph{bosonic} supergravity theory, consists of its equations of motion, which follow from a local Lagrangian principle, together with its Killing spinor equations, which consist of a spinorial differential system, typically of first order. The equations of motion determine their own Cauchy problem, defining the corresponding supergravity evolution flow and constraint equations, whereas the supersymmetric Killing spinor equations define in turn their own Cauchy problem, giving rise to their first-order supersymmetric evolution flow and constraint equations. Hence, we obtain two different systems of evolution and constraint equations that are intimately related since the Killing spinor equations provide partial first-order integrability of the supergravity equations of motion. Consequently, the well-posedness of the Cauchy problem in supergravity has several facets: we first have the standard Cauchy problem for the second-order equations of motion, which for the case of NS-NS supergravity has been proven to be well-posed in the upcoming doctoral dissertation of Oskar Schiller at Hamburg University. On the other hand, we have the Cauchy problem for the supersymmetric Killing spinor equations, which is yet to be studied in the literature. It is important to remark that the initial data for the supersymmetric Killing spinor equations is in general \emph{different} from the initial data for the supergravity equations of motion. However, experience shows \cite{Murcia:2020zig,Murcia:2021dur,Murcia:2022iba} that there is a natural map from the former to the latter that is at the basis of the rich interaction between these two Cauchy problems. Using this map it makes sense to talk about data that satisfies the constraint equations to \emph{both} Cauchy problems, and hence it is natural to compare the supersymmetric and supergravity, not necessarily supersymmetric, evolution flows starting on initial data admissible to both problems. In Theorem \ref{thm:compatibilityflows} we prove a compatibility result for these flows in NS-NS supergravity, which to the best of our knowledge is the first of its type in the supergravity literature. Many natural mathematical questions arise in this interacting framework of evolution flows. For instance, it is known that the solution space modulo isomorphism of the initial data equations of vacuum General Relativity is in general an \emph{infinite-dimensional space}. Similarly, we can expect that the initial data space modulo isomorphism of the supersymmetric Cauchy problem is also an infinite-dimensional space. However, based on previous experience of moduli spaces in supergravity, it is tempting to conjecture the following:
\begin{itemize}
	\item The moduli space of admissible initial data modulo isomorphism to both the supergravity and supersymmetric Cauchy problems is locally \emph{finite-dimensional}.
\end{itemize}

\noindent
If true, the previous statement should already hold for the case of parallel irreducible spinors on Ricci flat space-times, which we are currently studying in the specific case of four Lorentzian dimensions using the formalism of parabolic pairs and isotropic parallelism introduced in this dissertation. Contrary to more general supergravity systems, the Cauchy problem for parallel spinors in Lorentzian has been already studied in the mathematical literature. In the seminal work of Baum, Leistner and Lischewski \cite{BaumLeistnerLischewski,LeistnerLischewski,Lischewski}, the authors proved the well-posedness of the initial value problem of an irreducible complex spinor. The case of Killing spinor equations, of interest also in supergravity, was studied in arbitrary signature in \cite{ContiDalmasso}. On the other hand, the initial data set of the Cauchy problem defined by a parallel spinor has been considered in \cite{Ammann,Glockle:2019ovs,Glockle:2023ozz} in connection to the topology of the space of initial data conditions satisfying the dominant energy condition in General Relativity. For the specific case of four Lorentzian dimensions, the interaction between the evolution flow defined by a parallel spinor and the evolution flow defined by the Ricci-flatness equations has been investigated in \cite{Murcia:2020zig,Murcia:2021dur,RKS}, where it was proven that the former preserves the momentum and Hamiltonian constraints of the latter and therefore both flows \emph{coincide} on common initial data. It would be very interesting to continue studying these problems in the parallel case and also extend them to more general supergravity systems. It should be noted that, even if a given supergravity satisfies the strict dominant energy condition, the results of \cite{Ammann,Glockle:2019ovs,Glockle:2023ozz} may not directly apply, since the initial data is larger in the supergravity case and, as the NS-NS supergravity case considered in Chapter \ref{chapter:Globallyhyperbolicsusy} shows, initial data sets contain more variables than a Riemannian metric and a symmetric tensor. To the best of our knowledge, the topology of initial data sets in supergravity has not been considered in the literature, and it would be very interesting to generalize the results and techniques of \cite{Ammann,Glockle:2019ovs,Glockle:2023ozz} to supergravity. It seems that the Cauchy problem of bosonic supergravity and its supersymmetric Killing spinor equations has a very rich mathematical content that remains to be explored. We believe that their study could lead to many interesting mathematical problems and conjectures that will require the development of new tools in spin geometry, Riemannian geometry and analysis. 


\section{Main results}
\label{sec:mainresults}


In this section, we describe the main results of this dissertation. Generally speaking, all the results contained in the manuscript, except for Chapter \ref{chapter:spingeometryClifford}, subsection \ref{subsec:RKS} of Chapter \ref{chapter:IrreducibleSpinors4d}, Section \ref{sec:leftinvariant} of Chapter \ref{chapter:Globallyhyperbolicsusy}, and Appendix \ref{chapter:BundleGerbes}, are novel and have not been published elsewhere. 

\begin{itemize}
	\item In Chapter \ref{chapter:spingeometryClifford} we introduce in detail the theory of spinorial exterior forms associated to irreducible and real differential spinors in signature $(p-q) = 0,2\mod(8)$ closely following \cite{Cortes:2019xmk}. The main algebraic result of this chapter is Theorem \ref{thm:reconstruction}, which establishes a two-to-one correspondence between irreducible real spinors and those exterior forms satisfying the system of algebraic equations and inequalities given in the statement of the theorem. In Section \ref{sec:AlgebraicSpin(7)} of this chapter, which is extracted from \cite{SLSpin7}, we apply Theorem \ref{thm:reconstruction} to the case of irreducible chiral spinors in eight Euclidean dimensions. Despite the fact that this is a very classical case that has been amply studied in the literature \cite{LawsonMichelsohn}, Theorem \ref{thm:reconstruction} can be successfully applied to obtain a new result, namely an algebraic function whose self-dual critical points are precisely the Spin(7) structures of the underlying vector space. This result is presented in Theorems \ref{thm:Spin7algebraic} and \ref{thm:MetricPotentialVh}. In Theorem \ref{thm:Spin7algebraic} we compute the explicit algebraic square of an irreducible and chiral real spinor in eight Euclidean dimensions, whereas in Theorem \ref{thm:MetricPotentialVh} we elaborate on this result to construct the aforementioned potential. Although the main focus of the dissertation is on spinors in Lorentzian signature, we thought it was instructive to present an application in a very classical Riemannian geometry context, in order to illustrate the versatility of the methods employed. We hope this new potential can be of applicability in the study of Spin(7) structures via evolution flows and variational techniques. We continue in Section \ref{sec:differentialspinors} with the main \emph{differential result} of this chapter, which is given in Theorem \ref{thm:GCKS} and provides an equivalence between real and irreducible differential spinors in signature $(p-q) = 0,2\mod(8)$ and spinorial exterior forms satisfying the differential system given in the statement of the theorem. This result will be used extensively in the remaining of the dissertation. For applications to other signatures, see \cite{SAG2}.

	\item In Chapter \ref{chapter:IrreducibleSpinors4d} we apply the general framework developed in Chapter \ref{chapter:spingeometryClifford} to the case of differential spinors on Lorentzian four-manifolds. In Section \ref{sec:4dAlgebraic} we compute the square of an irreducible real spinor in four Lorentzian dimensions, and we associate to this construction a \emph{parabolic pair} and an \emph{isotropic parallelism}. These are the fundamental variables in terms of which we can rephrase every problem involving the aforementioned type of spinors into a problem involving exclusively equivalence classes of global coframes defined on the underlying manifold. The notion of isotropic parallelism seems to be new in the literature, and as exemplified in this dissertation, seems to be particularly well adapted to study the global geometric and topological structure of manifolds equipped with differential spinors. Choosing a representative in a given isotropic parallelism, we can naturally construct a \emph{complex tetrad} used in the Newman-Penrose formalism \cite{NewmanPenrose}, making thus contact with a classical object in four-dimensional General Relativity. Using these tools, we obtain in Theorem \ref{thm:differentialspinors4d} a characterization of the most general differential spinor on a strongly spin Lorentzian four-manifold $M$ in terms of a global coframe satisfying a system of partial differential equations that we give in its full generality in the statement of the theorem. We then use this result to study natural classes of differential spinors in four Lorentzian dimensions, focusing on the case of \emph{real Killing spinors}, which define a very natural class of differential spinors that occurs as the supersymmetric Killing spinor equations in minimal four-dimensional supergravity \cite{Ortin}. In particular, in Theorem \ref{thm:conformallyBrinkmann} we characterize the global geometry and topology of real Killing spinors on \emph{standard} conformally Brinkmann space-times. We end this chapter with Section \ref{sec:LorentzianInstantons}, in which we briefly discuss the notion of \emph{Lorentzian instanton} in four dimensions. This section illustrates the applicability of the theory of spinorial exterior forms to study spinors that are algebraically constrained to belong to the kernel of a given endomorphism. In this case the endomorphism is the curvature of a connection on a principal bundle, acting by Clifford multiplication on spinors, and the condition that it preserves a real and irreducible nowhere vanishing spinor reduces to a first-order differential equation for the connection which is the Lorentzian analog of the self-duality condition for connections in four Euclidean dimensions. These types of conditions appear extensively in supergravity but are yet to be considered in the mathematics literature. It would be very interesting to study them in more detail, developing their moduli space theory analogously to the Riemannian case \cite{KronheimerDonaldson}.
	
	\item In Chapter \ref{chapter:parallelspinorstorsion} we apply the geometric framework developed in Chapter \ref{chapter:IrreducibleSpinors4d} to \emph{torsion parallel spinors}, namely to real and irreducible spinors parallel under a general metric connection with torsion. This is a very natural condition to impose on a spinor and, whereas its Riemannian counterpart has been studied in the mathematics literature \cite{Agricola,AgricolaFriedrich,AgricolaHoll,Friedrich:2001nh}, the Lorentzian case seems to not have been systematically investigated, with the notable exceptions of \cite{ContiDalmasso,ErnstGalaev,Galaev}. In Theorem \ref{thm:existencenullcoframeII} we obtain a correspondence between torsion parallel spinors and solutions to an explicit exterior differential system \cite{BryantII,BryantBook} for isotropic parallelisms. This exterior differential system does not depend on any metric, hence realizing, at least to some extent, the motivation and ideology explained in \cite{Tomasiello:2011eb}. This opens up the possibility of applying the well-developed machinery of exterior differential systems to study torsion parallel spinors. It would certainly be very interesting to study aspects like Cartan's involutive test, possible prolongations, and associated Spencer cohomology \cite{BryantII,BryantBook}. Elaborating on this result, we extract, from the exterior differential system satisfied by the isotropic parallelism associated to every torsion parallel spinor, a pair of invariants. These invariants descend to cohomological invariants in de-Rham cohomology for Lorentzian four-manifolds equipped with an irreducible real spinor parallel with respect to a flat metric connection with torsion. It would be very interesting to elucidate a geometric interpretation of these invariants and clarify their role in the study of torsion parallel spinors with respect to a flat connection on a non-simply-connected Lorentzian four-manifold.
	
	\item In Chapter \ref{chapter:susyKundt4d} we initiate the proper study of four-dimensional NS-NS supersymmetric configurations and solutions. For this, we will extensively use the framework and associated results developed in the previous chapters. We proceed by first introducing in Section \ref{sec:NSNSsystem} the mathematical model of bosonic supergravity together with its supersymmetric Killing spinor equations, which is based on a choice of bundle gerbe with connective structure \cite{Murray} and principal $\mathbb{Z}$-bundle. In Section \ref{sec:sktorsionparallelspinors} of this chapter, as preparation for the study of Equation \eqref{eq:introsystem}, we investigate general skew-torsion parallel spinors, namely spinors parallel with respect to a connection with totally skew-symmetric torsion. In Section \ref{sec:susyconf} we consider the supersymmetric configurations of NS-NS supergravity, obtaining a general characterization in terms of a special class of isotropic parallelisms. In Theorem \ref{thm:susyfluxsolutions} of Section \ref{sec:susysolutions} we elaborate on the structure of supersymmetric configurations to obtain a characterization of supersymmetric NS-NS solution in the particular case in which the pseudo-norm of the curvature of the $b$-field of the solutions is nowhere vanishing. The general case is remarkably more complicated since the solution degenerates as the aforementioned norm tends to zero. We plan to study the general case in the future. In general terms, there is a lot to be learned about supersymmetric NS-NS solutions in four dimensions, and a general classification result seems out of reach, at least in the near future. We end this chapter with Section \ref{sec:dilatonflux}, in which we apply Theorem \ref{thm:susyfluxsolutions} to study the local structure of a canonical foliation that every NS-NS solution carries naturally. This foliation is generally singular, but for supersymmetric solutions with nowhere vanishing $b$-field curvature norm, it turns out to be regular. In particular, we obtain the general local form of the dual curvature of the $b$-field around a point in its pseudo-norm does not vanish.  
	
	\item Chapter \ref{chapter:Globallyhyperbolicsusy} is devoted to the study of globally hyperbolic four-dimensional supersymmetric configurations and solutions in NS-NS supergravity. For this, we first consider the \emph{globally hyperbolic reduction} of an isotropic parallelism, namely the decomposition and canonical representative that exists for isotropic parallelisms whose associated Lorentzian metric is globally hyperbolic. This framework provides a very efficient and transparent framework to deal with spinors on globally hyperbolic Lorentzian four-manifolds. We do this in Section \ref{sec:ReductionGloballyHyperbolic} of this chapter. In Section \ref{sec:evolutionskewtorsion} we study the evolution problem posed by a general skew-torsion parallel spinor on a globally hyperbolic Lorentzian four-manifold using the aforementioned notion of globally hyperbolic reduction of an isotropic parallelism. This allows to study the problem exclusively in terms of globally defined coframes on the underlying Cauchy hypersurface. Our first main result of this chapter is Theorem \ref{thm:Cauchytorsion}, in which we give simple criteria for a globally hyperbolic Lorentzian four-manifold to admit an irreducible real spinor parallel under a flat metric connection with skew-symmetric torsion. In Section \ref{sec:leftinvariant}, we demonstrate the efficacy of the isotropic parallelism formalism by explicitly solving the evolution flow determined by a parallel spinor on a globally hyperbolic Lorentzian four-manifold. This section is extracted from \cite{Murcia:2021dur}, to which the reader is referred for more details. In Section \ref{sec:reductionbundlegerbe} we embark on the globally hyperbolic reduction of abelian bundle gerbes with a connective structure of globally hyperbolic Lorentzian four-manifolds which, to the best of our knowledge, has not been considered elsewhere. This is a crucial step to study globally hyperbolic NS-NS solutions, as these are based on bundle gerbes. Finally, in Section \ref{sec:NSNSevolutionflow} we consider the supersymmetric NS-NS evolution flow, namely the evolution flow defined by the globally hyperbolic supersymmetric configurations of four-dimensional NS-NS supergravity. Along the way we obtain in detail the evolution and constraint equations of the bosonic NS-NS system together with the evolution and constraint equations of its Killing spinor equations. As discussed above, these define different but intimately related constrained evolution flows. Our main final result is Theorem \ref{thm:compatibilityflows}, in which we obtain the first compatibility criteria between these flows. Interestingly enough, this implies that both flows are compatible only if the Hamiltonian constraint of a globally hyperbolic Lorentzian four-manifold equipped with a spinor parallel with respect to the Levi-Civita connection is satisfied, illustrating the rigidity of globally hyperbolic supersymmetric NS-NS solutions.  
\end{itemize}
 
\begin{table}[h!]
	\centering
	\caption{Table of Theorems}
	\label{tab:original_results}
	\begin{tabular}{|l|m{4.5cm}|m{5.5cm}|}
		\hline
		\textbf{Theorem} & \textbf{Description} & \textbf{Potential significance} \\
		\hline
		Theorems \ref{thm:Spin7algebraic} \& \ref{thm:MetricPotentialVh} & Constructs an algebraic function on the space of self-dual four-forms whose critical points are precisely the Spin(7) structures of the underlying eight-dimensional vector space. & This provides a new variational framework to study Spin(7) structures, potentially opening them up to a new analysis via evolution flows. \\
		\hline
		Theorem \ref{thm:differentialspinors4d} & Characterizes the most general differential spinor on a Lorentzian four-manifold in terms of a global coframe satisfying a specific differential system. & This theorem translates a general spinorial problem into the more concrete language of the differential geometry on the frame bundle, using the newly introduced concept of an \emph{isotropic parallelism}. \\
		\hline
		Theorem \ref{thm:conformallyBrinkmann} & Obtains the global geometry and topology of real Killing spinors on standard conformally Brinkmann spacetimes. & This result gives a detailed picture for a specific, physically relevant class of solutions that appear in minimal supergravity. \\
		\hline
		Theorem \ref{thm:existencenullcoframeII} & Establishes a correspondence between torsion parallel spinors and solutions to an explicit exterior differential system for isotropic parallelisms, which does not depend on a metric. & This reframes the study of torsion parallel spinors in a metric-independent way, opening the door to applying the powerful machinery of exterior differential systems and Cartan's methods. In particular, it leads to the discovery of new invariants valued in $H^{1}(M,\mathbb{R})$. \\
		\hline
		Theorem \ref{thm:susyfluxsolutions} & Gives a characterization of supersymmetric NS-NS solutions in four dimensions for which the curvature of the $b$-field has a non-vanishing pseudo-norm. & This is a key step towards a classification of supersymmetric NS-NS solutions in four dimensions, handling the non-degenerate case and showing that the canonical foliation of the spacetime becomes regular. \\
		\hline
		Theorem \ref{thm:Cauchytorsion} & Provides a simple criteria for a globally hyperbolic Lorentzian four-manifold to admit an irreducible real spinor parallel under a flat metric connection with skew-symmetric torsion. & This result uses the framework of isotropic parallelisms to address the existence problem for a specific class of torsion parallel spinors in the globally hyperbolic setting. \\
		\hline
		Theorem \ref{thm:compatibilityflows} & Establishes the first compatibility criteria between the evolution flow of the NS-NS supergravity system and the flow of its first-order supersymmetric subsystem. & This result demonstrates the rigidity of the system, showing that the two flows are compatible only if a very restrictive condition is satisfied. \\
		\hline
	\end{tabular}
\end{table}

\renewcommand{\chaptername}{Chapter}

\renewcommand{\leftmark}{Chapter \thechapter. Spin geometry and bundles of Clifford modules}

\chapter{Spinorial exterior forms and bundles of irreducible Clifford modules}
\label{chapter:spingeometryClifford}


In this chapter, we follow \cite{Cortes:2019xmk} to develop the theory of spinorial exterior forms associated with sections of bundles of real irreducible Clifford modules in signature $(p-q)\equiv_8 0,2$, which includes the main case of interest in this dissertation, namely that of Lorentzian signature $(3,1)$ in four dimensions. The goal is to give a precise characterization of the \emph{square} of a spinor as an element in a given semi-algebraic real set and apply this description to the study of differential spinors, to be introduced in Section \ref{sec:differentialspinors}, see Definition \ref{def:generalizedKS}. We first introduce the algebraic theory, which we then proceed to extend to a differentiable theory on pseudo-Riemannian manifolds of signature $(p-q)\equiv_8 0,2$ equipped with a bundle of irreducible real Clifford modules. We begin with some algebraic preliminaries that will be of use in later sections.


\section{Real vectors as endomorphisms}
\label{sec:vectorasendo}


Let $\Sigma$ be a real vector space of positive even dimension $N$ equipped with a non-degenerate bilinear pairing $\cB \colon \Sigma\times \Sigma \to \R$, which we assume to be either symmetric or skew-symmetric. In this situation, the pair $(\Sigma,\cB)$ is called a {\em paired vector space}. We say that $\cB$ has {\em symmetry type} $s\in \mathbb{Z}_2$ if:
\begin{equation*}
\cB(\xi_1,\xi_2) = s \cB(\xi_2,\xi_1) \qquad \forall\, \xi_1, \xi_2\,\, \in \Sigma\, . 
\end{equation*}

\noindent
Hence $\cB$ is symmetric if it has symmetry type $s = +1$ and skew-symmetric if it has symmetry type $s = -1$.  Let $(\End(\Sigma),\circ)$ be the unital associative real algebra of linear endomorphisms of $\Sigma$, where $\circ$ denotes composition of linear maps. Given $E\in \End(\Sigma)$, we denote by $E^t\in \End(\Sigma)$ the adjoint of $E$ taken with respect to $\cB$, which is uniquely determined by the condition:
\begin{equation*}
\cB(\xi_1,E(\xi_2))=\cB(E^t(\xi_1),\xi_2)  \qquad \forall \,\, \xi_1,\xi_2\in \Sigma\, . 
\end{equation*}

\noindent
The map $E\rightarrow E^t$ is a unital anti-automorphism of the real algebra $(\End(\Sigma),\circ)$.

\begin{definition}
An endomorphism $E\in \End(\Sigma)$ is called {\em tame} if its rank satisfies $\rk(E)\leq 1$.
\end{definition}

\noindent 
Thus $E$ is tame if and only if it vanishes or is of rank one. Let:
\begin{equation*}
\cT :=  \{E\in \End(\Sigma) \,|\, \rk(E)\leq 1\} \subset \End(\Sigma)
\end{equation*}

\noindent
be the real determinantal variety of tame endomorphisms of $\Sigma$ and:
\begin{equation*}
\dot{\cT} :=  \cT\backslash\left\{ 0\right\}= \{E\in \End(\Sigma) \,|\, \rk(E)=1\}
\end{equation*}

\noindent
be its open subset consisting of endomorphisms of rank one.  We understand $\cT$ as a real affine variety of dimension $2N-1$ in the vector space $\End(\Sigma)\simeq \R^{N^{2}}$ and $\dot{\cT}$ as a semi-algebraic variety. Elements of $\cT$ can be written as:
\begin{equation*}
E = \xi\otimes \beta
\end{equation*}

\noindent
for some $\xi\in\Sigma$ and $\beta\in \Sigma^\ast$, where $\Sigma^\ast=\Hom(\Sigma,\R)$ denotes the real vector space dual to $\Sigma$. Notice that $\tr(E) = \beta(\xi)$. When $E\in \cT$ is non-zero, the vector $\xi$ and the linear functional $\beta$ appearing in the relation above are non-zero and determined by $E$ up to transformations of the form:
\begin{equation*}
(\xi,\beta)\rightarrow (\lambda \xi, \lambda^{-1}\beta)
\end{equation*}

\noindent
with $\lambda\in \R^\times$. In particular, $\dot{\cT}$ is a manifold diffeomorphic to the quotient $(\R^N\setminus\{0\})\times (\R^N\setminus\{0\})/\R^\times$, where $\R^\times$ acts with weights $+1$ and $-1$ on the two copies of $\R^N\setminus\{0\}$.

\begin{definition}
\label{def:squarevectorspace}
Let $\kappa \in \mathbb{Z}_2$ be a sign factor. The $\kappa$ {\em square map} of a paired vector space $(\Sigma,\cB)$ is the following quadratic map:
\begin{equation*}
\cE_{\kappa}\colon \Sigma \to \cT\, , \quad \xi \mapsto \cE_{\kappa}(\xi)=\kappa\, \xi \otimes \xi^{\ast} 
\end{equation*}

\noindent
where $\xi^\ast  :=  \cB(-, \xi)\in \Sigma^\ast$ is the linear map dual to $\xi$ relative to $\cB$.  
\end{definition}

\noindent 
Let $\kappa \in \mathbb{Z}_2$ be a sign factor and consider the open set $\dot{\Sigma} :=  \Sigma\backslash\left\{ 0\right\}$. Recall that $\cE_{\kappa}(\xi)=0$ if and only if $\xi=0$, hence $\cE_\pm(\dot{\Sigma})\subset \dot{\cT}$. Let $\dot{\cE}_{\kappa}\colon \dot{\Sigma}\to \dot{\cT}$ be the restriction of $\cE_{\kappa}$ to $\dot{\Sigma}$.

\begin{lemma}
\label{lemma:2to1E}
For each $\kappa\in \mathbb{Z}_2$, the restricted quadratic map $\dot{\cE}_{\kappa}\colon \dot{\Sigma}\to \dot{\cT}$ is two-to-one, namely:
\begin{equation*}
\dot{\cE}_{\kappa}^{-1}(\{\kappa\,\xi\otimes \xi^\ast\}) =  \left\{-\xi,\xi\right\}  \quad \forall \xi \in \dot{\Sigma}\, .
\end{equation*}

\noindent
Moreover, $\cE_\kappa$ is a real branched double cover of its image, which is ramified at the origin.
\end{lemma}

\begin{proof}
Follows from the presentation of $\dot{\cT}$ as a manifold diffeomorphic to the quotient $(\R^N\setminus\{0\})\times (\R^N\setminus\{0\})/\R^\times$, where $\R^\times$ acts with weights $+1$ and $-1$ on the two copies of $\R^N\setminus\{0\}$.
\end{proof}

\noindent
The map $\cE_{\kappa}$ need not be surjective. To characterize its image, we introduce the notion of \emph{admissible endomorphism}. Let $(\Sigma,\cB)$ be a paired vector space of symmetry type $s$.

\begin{definition}
An endomorphism $E$ of $\Sigma$ is called {\em $\cB$-admissible} if it satisfies the conditions:
\begin{equation*}
E\circ E = \tr(E) E~~\mathrm{and}~~ E^{t} = s E\, . 
\end{equation*}
\end{definition}

\noindent 
Let:
\begin{equation*}
\cC :=  \left\{ E\in\End(\Sigma)\,\, |\,\, E\circ E = \tr(E) E\, , \,\, E^{t} = s E \right\}
\end{equation*}

\noindent
denote the real cone of $\cB$-admissible endomorphisms of $\Sigma$.

\begin{remark} 
Tame endomorphisms are not related to admissible endomorphisms in any simple way. A tame endomorphism need not be admissible, since for instance it need not have a symmetry type with respect to $\cB$. On the other hand, an admissible endomorphism need not be tame, since it can have rank larger than one, as a quick inspection of explicit examples in four dimensions shows.
\end{remark}

\noindent
Let $\cZ := \cT \cap \cC$ denote the real cone of those endomorphisms of $\Sigma$ which are both tame and admissible with respect to $\cB$ and consider the open set $\dot{\cZ}  :=  \cZ\backslash\left\{ 0\right\}$.

\begin{lemma}
\label{lemma:EcE}
We have:
\begin{equation*}
\cZ = \Im(\cE_{+}) \cup \Im(\cE_{-})~~\mathrm{and}~~\Im(\cE_{+}) \cap \Im(\cE_{-})=\{0\}\, .
\end{equation*}

\noindent
Hence an endomorphism $E\in\End(\Sigma)$ belongs to $\dot{\cZ}$ if and only if there exists a non-zero vector $\xi \in \dot{\Sigma}$ and a sign factor $\kappa\in \mathbb{Z}_2$ such that:
\begin{equation*}
E = \cE_{\kappa}(\xi)\, .
\end{equation*}

\noindent
Moreover, $\kappa$ is uniquely determined by $E$ through this equation while $\xi$ is determined up to sign.
\end{lemma}

\begin{proof}
Let $E\in \dot{\cZ}$. Since $E$ has unit rank, there exists a non-zero vector $\xi\in\Sigma$ and a non-zero linear functional $\beta\in \Sigma^\ast$ such that $E = \xi\otimes \beta$.  Since $\cB$ is non-degenerate, there exists a unique non-zero $\xi_0\in \Sigma$ such that $\beta = \cB(-,\xi_0)=\xi_0^\ast$. The condition $E^t = s E$
amounts to:
\begin{equation*}
\cB(-,\xi_0) \xi = \cB(-,\xi) \xi_0\, .
\end{equation*}

\noindent
Since $\cB$ is non-degenerate, there exists an element $\chi\in \Sigma$ such that $\cB(\chi,\xi)\neq 0$, which by the previous equations also satisfies $\cB(\chi,\xi_0)\neq 0$. Hence:
\begin{equation*}
\xi_0 = \frac{\cB(\chi,\xi_0)}{\cB(\chi,\xi)} \xi=\frac{\cB(\xi_0, \chi)}{\cB(\xi, \chi)} \xi
\end{equation*}

\noindent
and:
\begin{equation*}
E = \frac{\cB(\xi_0, \chi)}{\cB(\xi, \chi)} \xi\otimes\xi^\ast\, .
\end{equation*}
	
\noindent	
Using the rescaling:
\begin{equation*}
\xi\mapsto \xi^{\prime} :=  \bigg\vert\frac{\cB(\xi_0, \chi)}{\cB(\xi, \chi)}\bigg\vert^{\frac{1}{2}} \xi
\end{equation*}

\noindent
the previous relation gives $E=\kappa\, \xi^{\prime}\otimes (\xi^{\prime})^\ast \in \Im(\cE_\kappa)$, where:
\begin{equation*}
\kappa := \sign\left(\frac{\cB(\xi_0, \chi)}{\cB(\xi, \chi)}\right)\, .
\end{equation*}

\noindent
This implies the inclusion $\cZ \subseteq \Im(\cE_{+}) \cup \Im(\cE_{-}) $. Lemma \ref{lemma:2to1E} now shows that $\xi^{\prime}$ is unique up to sign. The inclusion $\Im(\cE_{+})\cup \Im(\cE_{-})\subseteq \cZ$ follows by direct computation using the explicit form $E = \kappa\, \xi \otimes \xi^{\ast}$ of an endomorphism $E\in \Im(\cE_{\kappa})$, which implies:
\begin{equation*}
E\circ E = \cB(\xi,\xi) E\, ,\quad  E^t = s E\, , \quad \tr(E) = \cB(\xi,\xi)\, .
\end{equation*}

\noindent
Combining the two inclusions above gives $\cZ = \Im(\cE_{+}) \, \cup\, \Im(\cE_{-})$. Relation $\Im(\cE_{+}) \, \cap\, \Im(\cE_{-})=\{0\}$ follows immediately from Lemma \ref{lemma:2to1E}.
\end{proof}

\begin{definition}
The {\em signature} $\kappa\in \mathbb{Z}_2$ of an element $E\in \dot{\cZ}$ with respect to $\cB$ is the sign factor $\kappa$ determined as in Lemma \ref{lemma:EcE}.  
\end{definition}

\noindent 
In view of the above, given $\kappa\in \mathbb{Z}_2$ we define $\cZ_{\kappa} :=  \Im(\cE_{\kappa})$. Then:
\begin{equation*}
\cZ_-=-\cZ_+\, , \quad \cZ=\cZ_+\cup \cZ_-\, , \quad \mathrm{and}\quad \cZ_+\cap \cZ_-=\{0\}\, .
\end{equation*}
 
\begin{remark}
Let $\Z_2$ act on $\Sigma$ and on $\cZ\subset \End(E)$ by sign multiplication. Then $\cE_+$ and $\cE_-$ induce the same map between the quotients $\Sigma/\mathbb{Z}_2$ and $\cZ/\mathbb{Z}_2$ and this map is a bijection by virtue of Lemma \ref{lemma:EcE}.  
\end{remark}

\noindent
Given any endomorphism $A\in \End(\Sigma)$, define a possibly degenerate bilinear pairing $\cB_A$ on $\Sigma$ as follows:
\begin{equation}
\label{eq:cBA}
\cB_A(\xi_1,\xi_2)  :=  \cB(\xi_1, A(\xi_2)) \qquad  \forall \xi_1 , \xi_2 \in \Sigma\, .
\end{equation}

\noindent
Notice that $\cB_A$ is symmetric if and only if $A^t = s A$ and skew-symmetric if and only if $A^t = -s A$.

\begin{prop}
The open set $\dot{\cZ}$ has two connected components, which are given by:
\begin{equation*}
\dot{\cZ}_{+}  :=  \Im(\dot{\cE}_{+})=\Im(\cE_+)\setminus\{0\}\, , \qquad \dot{\cZ}_{-}  :=  \Im(\dot{\cE}_{-})=\Im(\cE_-)\setminus\{0\}
\end{equation*}

\noindent
and satisfy:
\begin{equation*}
\dot{\cZ}_{+}=\{E\in \cZ|\kappa_E=+1\} \subset \left\{ E\in \cZ\,\, \vert\,\, \cB_E\geq 0 \right\}\, , \quad  \dot{\cZ}_{-} =\{E\in \cZ|\kappa_E=-1\}\subset \left\{ E\in \cZ\,\, \vert\,\, \cB_E\leq 0 \right\}\, .
\end{equation*}

\noindent
Moreover, the map $\dot{\cE}_\kappa  \colon \dot{\Sigma}\to \dot{\cZ}_\kappa$ defines a principal $\Z_2$-bundle over $\dot{\cZ}_\kappa$.
\end{prop}

\begin{proof}
By Lemma \ref{lemma:EcE}, we know that $\dot{\cZ} = \dot{\cZ}_{+} \cup \dot{\cZ}_{-}$ and $\dot{\cZ}_{+}\cap \dot{\cZ}_{-} = \emptyset$. The open set $\dot{\Sigma}$ is connected because $N=\dim\Sigma\geq 2$. Fix $\kappa\in \mathbb{Z}_2$. Since the continuous map $\cE_{\kappa}$ surjects onto $\dot{\cZ}_{\kappa}$, it follows that $\dot{\cZ}_{\kappa}$ is connected. Let $E\in \dot{\cZ}_{\kappa}$. The pairing $\cB_E$ is symmetric since $E^t= s E$. By Lemma \ref{lemma:EcE}, we have $E = \kappa\, \cB(-,\xi_0) \xi_0$ for some non-zero $\xi_0\in \Sigma$ and hence:
\begin{equation*}
\cB_E(\xi,\xi) = \cB(\xi , E(\xi)) = \kappa \vert\cB(\xi,\xi_0)\vert^2 \qquad \forall\,\, \xi \in \Sigma\, .
\end{equation*}

\noindent
Since $\xi_0\neq 0$ and $\cB$ is non-degenerate, this shows that $\cB_E$ is nontrivial and that it is positive semidefinite when restricted to $\dot{\cZ}_{+}$ and negative semidefinite when restricted to $\dot{\cZ}_{-}$. The remaining statement follows from Lemma \ref{lemma:2to1E} and Lemma \ref{lemma:EcE}.
\end{proof}

\begin{prop}
\label{prop:admissiblendo}
$\dot{\cZ}$ is a manifold diffeomorphic to $\R^\times\times \R\P^{N-1}$, where $N=\dim \Sigma$. 
\end{prop}

\begin{proof}
Let $\vert\cdot\vert^2_0$ denotes the norm induced by any scalar product on $\Sigma$. The map:
\begin{equation*}
\dot{\cZ} \xrightarrow{\sim}\R^\times \times \R\P^{N-1}\, , \quad  \kappa\, \xi\otimes\xi^{\ast} \mapsto (\kappa\,\vert\xi\vert_0^2,[\xi])
\end{equation*}

\noindent
is a diffeomorphism.
\end{proof}

\noindent
The maps $\dot{\cE}_{+}\colon \dot{\Sigma}\to \End(\Sigma)\setminus \{0\}$ and $\dot{\cE}_{-}\colon \dot{\Sigma}\to \End(\Sigma)\setminus \{0\}$ induce the same map:
\begin{equation*}
\P\cE\colon \P(\Sigma) \to \P(\End(\Sigma))\, , \quad  [\xi]\mapsto  [\xi\otimes\xi^{\ast}]
\end{equation*}

\noindent
between the projectivizations $\P(\Sigma)$ and $\P(\End(\Sigma))$ of the real vector spaces $\Sigma$ and $\End(\Sigma)$. Setting $\P\cZ  :=  \dot{\cZ} /\R^\times\subset \P\End(\Sigma)$, Proposition \ref{prop:admissiblendo} gives the following result.

\begin{prop}
\label{prop:projectivediff}
The map $\P\cE:\P(\Sigma) \rightarrow \P\cZ$ is a diffeomorphism.
\end{prop}

\noindent
We will refer to $\P\cE\colon \P(\Sigma) \stackrel{\sim}{\to} \P\cZ$ as the {\em projective square map} of $(\Sigma,\cB)$. The following result characterizes an open subset of the cone $\cC$ of admissible endomorphisms which consists of rank one elements.

\begin{prop}
\label{prop:tametracenon0}
Let $E\in \cC(\Sigma,\cB)$ be a $\cB$-admissible endomorphism of $\Sigma$. If $\tr(E)\neq 0$, then $E$ is of rank one.
\end{prop}

\begin{proof}
Define $P :=  \frac{E}{\tr(E)}$. Then $P^2 = P$, which implies $\rk(P)=\tr(P)$ and $\tr(P) = 1$. Hence $\rk(E) = \rk(P) = \tr(P)=1$.
\end{proof}

\noindent 
Define:
\begin{equation*}
\cK_0 :=  \left\{ \xi \in \Sigma\,\, \vert\,\, \cB(\xi,\xi) = 0 \right\}\, ,\quad  \cK_{\mu} :=  \left\{ \xi \in \Sigma\,\, \vert\,\, \cB(\xi,\xi) = \mu \right\}
\end{equation*}

\noindent
where $\mu\in \mathbb{Z}_2$. When $\cB$ is symmetric, the set $\cK_0\subset \Sigma$ is the isotropic cone of $\cB$ and $\cK_{\mu}$ are the positive and negative unit pseudo-spheres defined by $\cB$. When $\cB$ is skew-symmetric, we have $\cK_0 = \Sigma$ and $\cK_{\mu} = \emptyset$. Lemma \ref{lemma:EcE} and Proposition \ref{prop:tametracenon0} imply:

\begin{cor}
Assume that $\cB$ is symmetric, that is, $s=+1$. For any $\mu\in \mathbb{Z}_2$, the set $\cE_+(\cK_{\mu})\cup \cE_-(\cK_{\mu})$ is the real algebraic submanifold of $\End(\Sigma)$ given by:
\begin{equation*}
\cE_+(\cK_{\mu})\cup \cE_-(\cK_{\mu}) = \left\{ E\in \End(\Sigma) \,\, \vert\,\, E\circ E = \mu\,E\, , \,\, E^{t} = E\, , \,\, \tr(E) = \mu \right\}
\end{equation*}
\end{cor}

\begin{prop}
\label{prop:tamescalarprod}
If $\cB$ is definite, then every non-zero $\cB$-admissible endomorphism $E\in \cC\setminus \{0\}$ is tame, whence $\cZ = \cC$. In this case, the signature of $E$ with respect to $\cB$ is given by $\kappa=\sign(\tr(E))$.
\end{prop}

\begin{proof}
Let $E\in\cC$. By Proposition \ref{prop:tametracenon0}, the first statement follows if we can show that $\tr(E)\neq 0$ when $E\neq 0$. Since $E$ is admissible, it is symmetric with respect to the scalar product $\cB$ and hence diagonalizable with eigenvalues $\lambda_1,\hdots , \lambda_N\in \R$. Taking the trace of equation $E^2 = \tr(E)\, E$ gives:
\begin{equation*}
\tr(E)^2 = \sum_{i=1}^{N} \lambda_i^2\, .
\end{equation*}

\noindent
Since the right-hand side is a sum of squares, it vanishes if and only if $\lambda_1 = \hdots = \lambda_N = 0$, that is, if and only if $E=0$. This proves the first statement.  To prove the second statement, recall from Lemma \ref{lemma:EcE} that any non-zero tame admissible endomorphism $E$ has the form $E=\kappa \, \xi\otimes \xi^\ast$  for some $\xi\in \dot{\Sigma}$ and $\kappa\in \mathbb{Z}_2$. Taking the trace of this relation gives:
\begin{equation*}
\tr(E)=\kappa\, \xi^\ast(\xi)=\kappa \cB(\xi,\xi)
\end{equation*}

\noindent
which implies $\kappa=\sign(\tr(E)) \sign(\cB)$ since either $\cB(\xi,\xi)>0$ or $\cB(\xi,\xi)<0$ for every $\xi\in \dot{\Sigma}$. 
\end{proof}

\noindent 
A quick inspection of examples shows that there exist non-trivial admissible endomorphisms which are not tame, and thus satisfy $\tr(E)= 0$, as soon as there exists a totally isotropic subspace of $\Sigma$ of dimension at least two. In these cases we need to impose further conditions on the elements of $\cC$ in order to guarantee tameness. To describe such conditions, we consider the more general equation:
\begin{equation*}
E\circ A\circ E = \tr(A\circ E) E  \qquad \forall\,\, A\in \End(\Sigma)
\end{equation*}

\noindent
which is automatically satisfied by every $E\in \Im(\cE) = \Im(\cE_+)\cup \Im(\cE_-)$.

\begin{prop}
\label{prop:characterizationtamecone}
The following statements are equivalent:
\begin{enumerate}
\item $E$ is $\cB$-admissible and $\rk(E) = 1$, that is, $E\in \Im(\dot{\cE})$.
\item $E$ is $\cB$-admissible and there exists an endomorphism $A\in \End(\Sigma)$ satisfying:
\begin{equation}
\label{eq:conditionstame0}
E\circ A\circ E = \tr(E\circ A) E~~\mathrm{and}~~ \tr(E\circ A)\neq 0\, .
\end{equation}
\item $E\neq 0$ and relations:
\begin{equation}
\label{eq:conditionstame1}    
E\circ A\circ E = \tr(E\circ A) E\, , \qquad E^t = s E
\end{equation}

\noindent
hold for every endomorphism $A\in \End(\Sigma)$.
\end{enumerate}
\end{prop}

\begin{proof}
We first prove the implication $(2)\Rightarrow (1)$. By Proposition \ref{prop:tametracenon0}, it suffices to consider the case $\tr(E) = 0$. Assume $A\in\End(\Sigma)$ satisfies \eqref{eq:conditionstame0}. Define:
\begin{equation*}
A_\epsilon = \Id +\frac{\epsilon}{\tr(E\circ A)} A
\end{equation*}

\noindent
where $\epsilon \in \R_{>0}$ is a positive constant. For $\epsilon>0$ small enough, $A_{\epsilon}$ is invertible and the endomorphism $E_{\epsilon}  :=  E\circ A_{\epsilon}$ has non-vanishing trace given by $\tr(E_{\epsilon}) = \epsilon$. The first relation in \eqref{eq:conditionstame0} gives:
\begin{equation*}
E_{\epsilon}\circ E_{\epsilon} = \epsilon E_{\epsilon}\, .
\end{equation*}

\noindent
Hence $P :=  \frac{1}{\epsilon} E_{\epsilon}$ satisfies $P^2 = P$ and $\tr(P) = 1$, whence $\rk(E_\epsilon)=\rk(P)=1$. Since $A_{\epsilon}$ is invertible, this implies $\rk(E) = 1$ and hence $(1)$ holds.
	
The implication $(1)\Rightarrow (3)$ follows directly from Lemma \ref{lemma:EcE}, which shows that $E\in \Im(\cE_\kappa)$ for some sign factor $\kappa$. For the implication $(3)\Rightarrow (2)$, notice first that setting $A=\Id$ in \eqref{eq:conditionstame1} gives $E^2=\tr(E)$.  Non-degeneracy of the bilinear form induced by the trace on the space $\End(\Sigma)$ now shows that we can choose $A$ in equation \eqref{eq:conditionstame1} such that $\tr(E\circ A)\neq 0$.
\end{proof}


\subsection*{Two-dimensional examples}


Let $\Sigma$ be a two-dimensional real vector space with basis $\left\{e_1,e_2\right\}$. Any vector $\xi\in
\Sigma$ expands as:
\begin{equation*}
\xi = \xi_1 e_1 + \xi_2 e_2~~\mathrm{with}~~\xi_1, \xi_2 \in \R\, .
\end{equation*}

\noindent
Let $E_\xi :=  \cE_\kappa(\xi) := \kappa \xi\otimes \xi^{\ast}\in \End(\Sigma)$, where $\kappa\in \mathbb{Z}_2$. For any $S\in \End(\Sigma)$, we denote by ${\hat S}$ the matrix of $S$ in the basis $\left\{e_1,e_2\right\}$.

\begin{ep}
\label{ep:2dEuclidean}
Let $\cB$ be a scalar product on $\Sigma$ having $\left\{e_1,e_2\right\}$ as an orthonormal basis. Then:
\begin{equation*}
{\hat E}_\xi = \kappa \begin{pmatrix} \xi_1^2 & \xi_1 \xi_2 \\ \xi_1 \xi_2 & \xi_2^2 \end{pmatrix}
\end{equation*}

\noindent
and the relations $E^2_{\xi} = \tr(E_{\xi}) E_{\xi}$ and $E_{\xi}^t = E_{\xi}$ follow from this form. Conversely, let $E\in \End(\Sigma)$ satisfy $E^2 = \tr(E) E$ and $E^t = E$. The second of these conditions implies:
\begin{equation*}
{\hat E} = \left( \begin{array}{ccc} k_1 & b \\ 
b & k_2 \end{array}\right)~~(\mathrm{with}~b,k_1,k_2\in \R)\, .
\end{equation*}

\noindent
Condition $E^2 = \tr(E) E$ amounts to $b^2 = k_1 k_2$, implying that $k_1$ and $k_2$ have the same sign unless at least one of them vanishes (in which case $b$ must also vanish). Since $E$ is $\cB$-symmetric (and hence diagonalizable), its trace $\tr(E)=k_1+k_2$ vanishes if and only if $E=0$. Assume $E\neq 0$ and set:
\begin{equation*}
\kappa :=  \sign(\tr(E))=\sign(k_1+k_2)\, , \quad \xi_1 :=  \sqrt{|k_1|}\, , \quad \xi_2 :=  \sign(b) \kappa \sqrt{|k_2|}\, .
\end{equation*}

\noindent
Then $k_1 = \kappa \xi_1^2$, $k_2 = \kappa \xi_2^2$ and $b=\kappa \xi_1\xi_2$, showing that $E=E_\xi$ for some $\xi\in \Sigma\setminus \{0\}$. Hence conditions $E^2 = \tr(E) E$ and $E = E^{t}$ characterize endomorphisms of the form $E_\xi$.
\end{ep}

\begin{ep}
\label{ep:2dsplitsig}
Let $\cB$ be a split signature inner product on $\Sigma$ having $\left\{e_1,e_2\right\}$ as an orthonormal basis:
\begin{equation*}
\cB(e_1,e_1) = 1\, , \quad  \cB(e_2,e_2) = - 1\, , \quad \cB(e_1,e_2) = \cB(e_2,e_1) = 0\, .
\end{equation*}

\noindent
We have:
\begin{equation*}
E_\xi = \kappa \left(\begin{array}{ccc} \xi_1^2 & 
-\xi_1 \xi_2 \\ \xi_1 \xi_2 & -\xi_2^2 \end{array} \right)
\end{equation*}

\noindent
and the relations $E^2_{\xi} = \tr(E_{\xi}) E_{\xi}$ and $E_{\xi}^t = E_{\xi}$ follow directly from this form, where $~^t$ denotes the adjoint taken with respect to $\cB$.  Conversely, let $E\in \End(\Sigma)$ satisfy $E^t = E$. Then:
\begin{equation*}
E = \left( \begin{array}{ccc} k_1 & -b \\ b & k_{2} \end{array}\right)\,~~\mathrm{and}
~~ E A = \left( \begin{array}{ccc} k_1 & b \\ b & -k_{2} \end{array}\right)\,~~\mathrm{with}~b,k_1,k_2\in \R
\end{equation*}

\noindent
where $A =\diag(+1,-1)$. A direct computation shows that the conditions $E^2 = \tr(E) E$ and $E\circ A\circ E=\tr(E\circ A) E$ are equivalent to each other in this two-dimensional example and amount to the relation $b^2 = -k_1 k_2$, which implies that $E$ vanishes if and only if $k_1=k_2$. Let us assume that $E\neq 0$ and set:
\begin{equation*}
\kappa :=  \tr(E\circ A)=\sign(k_1-k_2)\, , \quad \xi_1 :=   \sqrt{|k_1|}\, , \quad \xi_2 :=  \sign(b) \kappa \sqrt{|k_2|}
\end{equation*}

\noindent
where $\sign(b) :=  0$ if $b=0$. Then it is easy to see that $k_1 = \kappa \xi_1^2$, $k_2 =-\kappa \xi_2^2$ and $b= \kappa \xi_1\xi_2$, which implies $E=E_\xi$. In this example endomorphisms $E$ that can be written in the form $E_\xi$ are characterized by the condition $E^t=E$, together with either of the two equivalent conditions $E^2=\tr(E) E$ or $E\circ A\circ E=\tr(E\circ A) E$.  Notice that $\tr(E)=k_1+k_2$ can vanish in this case. However, and in contrast to higher dimensional cases, in this two-dimensional example, the conditions $E\circ E = \tr(E) E$ and $E^{t} = E$ suffice to characterize the endomorphisms of the form $E = E_\xi$, including those which satisfy $\tr(E) = 0$.
\end{ep}

\begin{ep}
\label{ep:2dskew}
Let $\cB$ a symplectic pairing on $\Sigma$ having $\left\{e_1,e_2\right\}$ as a Darboux basis:
\begin{equation*}
\cB(e_1,e_1) = \cB(e_2,e_2) =  0\, , \quad \cB(e_1,e_2) =- \cB(e_2,e_1) = 1\, .
\end{equation*}

\noindent
The complex structure $A$ of $\Sigma$ with matrix given by:
\begin{equation*}
A=\left( \begin{array}{ccc} 0 & 1 \\ -1
& 0 \end{array} \right)
\end{equation*}

\noindent
tames $\cB$ to the scalar product $(-,-)$ defined through:
\begin{equation*}
(e_1,e_1) = (e_2,e_2) =  1\, , \quad  (e_1,e_2) =(e_2,e_1) = 0\, .
\end{equation*}

\noindent
We have:
\begin{equation}
\label{eq:matrix22skew}
E_{\xi} = \kappa \left( \begin{array}{ccc} \xi_1 \xi_2 & -\xi^2_1 \\ \xi_2^2
& - \xi_1 \xi_2 \end{array} \right)
\end{equation}

\noindent
which implies $E^2_{\xi} = 0$ and $E_{\xi}^t = - E_{\xi}$, where $^t$ denotes transposition with respect to $\cB$. Conversely, let $E\in \End(\Sigma)$ be an endomorphism satisfying $E^t = - E$. This condition implies:
\begin{equation*}
E  = \left( \begin{array}{ccc} k & -b \\ c & - k \end{array}
\right)\, , \quad  E A = \left( \begin{array}{ccc} b & k \\ k &
		c \end{array} \right)\quad (\mathrm{with}~k,b,c\in \R)\, .
\end{equation*}

\noindent
Notice that $\tr(E)=0$. Direct computation shows that the conditions $E^2 = 0$ and $E\circ A \circ E=\tr(E\circ A) E$
are equivalent to each other in this two-dimensional example and amount to the relation $k^2 = b c$, which in particular shows that $E$ vanishes if and only if $b=-c$. Assume that $E\neq 0$ and set:
\begin{equation*}
\kappa :=  \tr(E\circ A)=\sign(b+c)\, , \quad \xi_1 :=  \sqrt{|b|}\, , \quad \xi_2 :=  \sign(k)\kappa\sqrt{|c|}
\end{equation*}

\noindent
where $\sign(k) :=  0$ if $k=0$. Then it is easy to see that $b=\kappa \xi^2_1$, $c = \kappa \xi^2_2$ and $k = \kappa \xi_1 \xi_2$, which shows that $E=E_\xi$. Hence endomorphism which can be written in this form are characterized by the condition $E^t=-E$ together with either of the conditions $E^2=0$ or $E\circ A \circ E=\tr(E\circ A) E$, which, in this low-dimensional example, are equivalent to each other.
\end{ep}


\subsection*{Linear constraints}


\noindent 
The following result will be used in later sections.

\begin{prop}
\label{prop:constraintendo} 
Let $Q\in \End(\Sigma)$ and $\kappa\in \mathbb{Z}_2$ be a fixed sign factor. A real spinor $\xi\in \Sigma$ satisfies $Q(\xi) = 0$ if and only if $Q\circ\cE_{\kappa}(\xi) = 0$ or, equivalently, $\cE_{\kappa}(\xi)\circ Q^{t} = 0$, where $Q^{t}$ is the adjoint of $Q$ with respect to $\cB$.
\end{prop}

\begin{proof}
Take $\xi\in \Sigma$ and assume $Q(\xi) = 0$. Then:
\begin{equation*}
(Q\circ\cE_{\kappa}(\xi))(\chi) = \kappa\, Q(\xi)\,\xi^{\ast}(\chi) = 0 \quad  \forall \chi\in \Sigma
\end{equation*}

\noindent
and hence $Q\circ\cE_{\kappa}(\xi) = 0$. Conversely, assume that $Q\circ\cE_{\kappa}(\xi) = 0$ and pick $\chi\in\Sigma$ such that $\xi^{\ast}(\chi) \neq 0$ (which is possible since $\cB$ is non-degenerate). Then the same calculation as before gives:
\begin{equation*}
Q(\xi)\,\xi^{\ast}(\chi) = 0 
\end{equation*}

\noindent
implying $Q(\xi)=0$. The statement for $Q^{t}$ follows from the fact that $\cB$-transposition is an anti-automorphism of the real algebra $(\End(\Sigma),\circ)$, upon noticing that the relation $\cE_\kappa(\xi)^t=\sigma \cE_\kappa(\xi)$ implies $(Q\circ \cE_{\kappa}(\xi))^{t}=\sigma \cE_{\kappa}(\xi)\circ Q^{t}$.
\end{proof}

\begin{ep}
Let $(\Sigma,\cB)$ be a two-dimensional Euclidean vector space with orthonormal basis $\left\{e_1,e_2\right\}$ as in Example \ref{ep:2dEuclidean}. Let $Q\in\End(\Sigma)$ have the following matrix expression:
\begin{equation*}
Q = \left( \begin{array}{ccc} q & 0 \\ 0 & 0 \end{array}\right)~~\mathrm{with}~q\in \R^\times
\end{equation*}

\noindent
in this basis. Given $\xi\in \Sigma$,  Example \ref{ep:2dEuclidean} gives:
\begin{equation*}
E_{\xi} = \kappa \left( \begin{array}{ccc} \xi_1^2 & \xi_1 \xi_2\\ 
\xi_1 \xi_2 & \xi_2^2 \end{array} \right)\, , \quad 
Q E_{\xi} = \kappa \left( \begin{array}{ccc} \xi_1^2 q & q\xi_1\xi_2 \\ 0
& 0 \end{array} \right)\, .
\end{equation*}	

\noindent
Thus $Q\circ E_{\xi}$ vanishes if and only if $\xi_1 = 0$, that is if and only if $Q(\xi) = 0$.
\end{ep}


\section{Irreducible Clifford modules and spinorial exterior forms}
\label{sec:SpinorsAsPolyforms}


In this section we develop the theory of spinorial exterior forms associated to irreducible spinors in signature $(p-q) \equiv_8 0,2$. As we explain below, these are the polyforms that can be constructed as \emph{squares} of spinors.


\subsection{Irreducible real Clifford modules}


Let $V$ be an oriented $d$-dimensional real vector space equipped with a non-degenerate metric $h$ of signature $p-q\equiv_8 0, 2$, which implies that the dimension $d=p+q$ of $V$ is even, and let $(V^{\ast},h^{\ast})$ be the quadratic space dual to $(V,h)$, where $h^{\ast}$ denotes the metric dual to $h$. Let $\Cl(V^{\ast},h^{\ast})$ be the real Clifford algebra of this dual quadratic space, viewed as a $\Z_2$-graded associative algebra with decomposition:
\begin{equation*}
\Cl(V^{\ast},h^{\ast}) = \Cl^{\ev}(V^{\ast},h^{\ast}) \oplus \Cl^{\odd}(V^{\ast},h^{\ast})\, .
\end{equation*}

\noindent
In our conventions, the Clifford algebra satisfies:
\begin{equation}
\label{Crel}
\theta^2 = h^\ast(\theta,\theta) \quad \forall\,\, \theta\in V^\ast\, .
\end{equation}

\noindent
Let $\pi$ denote the standard automorphism of $\Cl(V^{\ast},h^{\ast})$, which acts as minus the identity on $V^{\ast}\subset \Cl(V^{\ast},h^{\ast})$, and let $\tau$ denote its standard anti-automorphism, which acts as the identity on $V^{\ast}\subset \Cl(V^{\ast},h^{\ast})$. These two commute and their composition is an anti-automorphism denoted by $\hat{\tau} = \pi\circ\tau=\tau\circ \pi$. Let $\Cl^\times(V^{\ast},h^{\ast})$ denote the group of units in $\Cl(V^{\ast},h^{\ast})$. Its \emph{twisted adjoint representation} is the morphism of groups $\widehat{\Ad}\colon \Cl^{\times}(V^{\ast},h^{\ast}) \to \Aut(\Cl(V^{\ast},h^{\ast}))$ defined through:
\begin{equation*}
\widehat{\Ad}_x(y) = \pi(x)\, y\, x^{-1} \quad \forall\,\, x , y  \in \Cl^{\times}(V^{\ast},h^{\ast})\, .
\end{equation*}

\noindent
We denote by $\Gamma(V^{\ast},h^{\ast})\subset \Cl(V^{\ast},h^\ast)$ the Clifford group of $\Cl(V^{\ast},h^{\ast})$, which is defined as the subgroup of $\Cl^{\times}(V^{\ast},h^{\ast})$ that preserves $V^{\ast}$ via $\widehat{\Ad}$, that is:
\begin{equation*}
\Gamma(V^{\ast},h^{\ast})  :=  \left\{ x \in 
\Cl^{\times}(V^{\ast},h^{\ast})\quad \vert\quad \widehat{\Ad}_x(V^{\ast}) = V^{\ast} \right\}
\end{equation*}

\noindent
The Clifford group fits into the short exact sequence:
\begin{equation}
\label{eq:GammaSeq}
1 \to \R^\times \hookrightarrow  \Gamma(V^{\ast},h^{\ast}) \xrightarrow{\widehat{\Ad}} \O(V^{\ast},h^{\ast})\to 1
\end{equation}

\noindent
where $\O(V^{\ast},h^{\ast})$ is the orthogonal group defined by the the quadratic space $(V^{\ast},h^{\ast})$. The special orthogonal group and its connected component of the identity will be denoted by $\SO(V^{\ast},h^{\ast})$ and $\SO_o(V^{\ast},h^{\ast})$, respectively.  Recall that the pin and spin groups of $(V^{\ast},h^{\ast})$ are the subgroups of $\Gamma(V^{\ast},h^{\ast})$ defined as follows:
\begin{equation*}
\Pin(V^{\ast},h^{\ast})  :=  \left\{x\in \Gamma(V^{\ast},h^{\ast})\,\, \vert\,\, N(x)^2 = 1 \right\}\, , \quad 
\Spin(V^{\ast},h^{\ast})  :=  \Pin(V^{\ast}, h^{\ast})\cap \Cl^{\ev}(V^{\ast},h^{\ast})
\end{equation*}

\noindent
where $N\colon \Gamma(V^{\ast},h^{\ast}) \to \R^\times$ is the \emph{Clifford norm morphism}, given by:
\begin{equation*}
N(x)  :=  \hat{\tau}(x)\, x \quad \forall\,\, x\in \Gamma(V^{\ast},h^{\ast})\, .
\end{equation*}

\noindent
We have $N(x)^2 = N(\pi(x))^2$ for all $x\in\Gamma(V^{\ast},h^{\ast})$. For $p q \neq 0$, the groups $\SO(V^\ast,h^\ast)$, $\Spin(V^{\ast},h^{\ast})$ and $\Pin(V^{\ast},h^{\ast})$ are disconnected; the first have two connected components while the latter has four. The connected components of the identity in $\Spin(V^{\ast},h^{\ast})$ and $\Pin(V^\ast,h^\ast)$ coincide, being given by:
\begin{equation*}
\Spin_o(V^{\ast},h^{\ast}) = \left\{x\in \Gamma(V^{\ast},h^{\ast})\quad \vert\,\, N(x) = 1 \right\}
\end{equation*}

\noindent
and we have $\Spin(V^\ast,h^\ast)/\Spin_o(V^\ast,h^\ast)\simeq \Z_2$ and $\Pin(V^{\ast},h^{\ast})/\Spin_o(V^\ast,h^\ast)\simeq \Z_2\times \Z_2$.

Let $\Sigma$ be a finite-dimensional $\R$-vector space and $\gamma\colon \Cl(V^{\ast},h^{\ast})\to \End(\Sigma)$ a Clifford
representation. Then $\Spin(V^{\ast},h^{\ast})$ acts on $\Sigma$ through the restriction of $\gamma$ and \eqref{eq:GammaSeq} induces the following short exact sequence:
\begin{equation}
\label{eq:SpinSeq}
1\to \Z_2 \to \Spin(V^{\ast},h^{\ast}) \xrightarrow{\widehat{\Ad}} \SO(V^{\ast},h^{\ast}) \to 1
\end{equation}

\noindent
which in turn gives the exact sequence:
\begin{equation*}
1\to \Z_2 \to \Spin_o(V^{\ast},h^{\ast})  \xrightarrow{\widehat{\Ad}} \SO_o(V^{\ast},h^{\ast}) \to 1 
\end{equation*}

\noindent
where $\SO_o(V^{\ast},h^{\ast})$ denotes the identity component of the special orthogonal group $\SO(V^{\ast},h^{\ast})$. In signature $p-q\equiv_8 0, 2$, the Clifford algebra $\Cl(V^{\ast},h^{\ast})$ is simple and admits a unique irreducible representation $\gamma\colon\Cl(V^{\ast},h^{\ast})\xrightarrow{\simeq} \End(\Sigma)$ on the endomorphisms $\End(\Sigma)$ of a real vector space $\Sigma$ of dimension $2^{\frac{d}{2}}$. In particular, $\gamma$ defines an isomorphism of unital and associative real algebras.  

We would like to equip $\Sigma$ with a non-degenerate bilinear pairing \emph{compatible} with the representation $\gamma\colon\Cl(V^{\ast},h^{\ast})\xrightarrow{\simeq} \End(\Sigma)$. Ideally, such compatibility condition should translate into the invariance of the bilinear pairing at least under the natural action of the pin group induced by $\gamma$. However, when $p q\neq 0$ it may not be possible to satisfy this condition. Instead, we consider the weaker notion of \emph{admissible bilinear pairing} introduced in \cite{AC,ACDP}, which encodes the \emph{best} compatibility condition with Clifford multiplication that can be generally imposed on a bilinear pairing on $\Sigma$ in arbitrary dimension and signature. The following result, which can also be found in {\rm \cite[Theorem 13.17]{HarveyBook}}, summarizes the main properties of admissible bilinear pairings.  

\begin{thm} 
\label{thm:admissiblepairings}
Suppose that $h$ has signature $p-q\equiv_8 0,2$. Then the irreducible real Clifford module $\Sigma$ admits two non-degenerate bilinear pairings $\cB_{+}\colon \Sigma\times\Sigma\to \R$ and $\cB_{-}\colon \Sigma\times\Sigma\to \R$ (each determined up to multiplication by a non-zero real number) such that:
\begin{equation}
\label{eq:admissiblepairins}
\cB_{+}(\gamma(x)(\xi_1),\xi_2) = \cB_{+}(\xi_1, \gamma(\tau(x))(\xi_2))\, , \quad \cB_{-}(\gamma(x)(\xi_1),\xi_2) = \cB_{-}(\xi_1, \gamma(\hat{\tau}(x))(\xi_2))
\end{equation}

\noindent
for all $x\in \Cl(V^{\ast},h^{\ast})$ and $\xi_1, \xi_2 \in \Sigma$. The symmetry properties of $\cB_{+}$ and $\cB_{-}$ are as follows in terms of the modulo $4$ reduction of $k  :=  \frac{d}{2}$:
\begin{center}
\begin{tabular}{ | l | p{3cm} | p{3cm} | p{3cm} | p{3cm} |}
\hline
$k\,mod$ 4 & 0 & 1 & 2 & 3 \\ \hline
$\cB_{+}$ & Symmetric & Symmetric & Skew-symmetric & Skew-symmetric   \\ \hline
$\cB_{-}$ & Symmetric & Skew-symmetric & Skew-symmetric & Symmetric  \\ \hline
\hline
\end{tabular}
\end{center}

\noindent
In addition, if $\cB_{\sigma}$, with $\sigma\in \mathbb{Z}_2$, is symmetric, then it is of split signature unless $pq=0$, in which case $\cB_{\sigma}$ is definite. 
\end{thm}

\begin{proof}
We pick an $h^\ast$-orthonormal basis $\left\{ e^{i}\right\}_{i=1,\ldots, d}$ of $V^\ast$ and let:
\begin{equation*}
	\mathrm{K}(\left\{ e^{i}\right\})  :=  \{1\}\cup \left\{ \pm e^{i_1}\cdot \ldots \cdot e^{i_k}\, | \,\, 1\leq i_1<\ldots <i_k\leq d\, , \, 1\leq k\leq d\right\}
\end{equation*}

\noindent 
be the finite multiplicative subgroup of $\Cl(V^{\ast},h^{\ast})$ generated by the elements $\pm e^i$. Averaging over
$\mathrm{K}(\left\{ e^{i}\right\})$, we construct an auxiliary positive-definite inner product $(-,-)$ on $\Sigma$ which is invariant under the action of this group. This product satisfies:
\begin{equation*}
	(\gamma(x)(\xi_1), \gamma(x)(\xi_2)) = (\xi_1, \xi_2) \quad \forall\,\, x\in \mathrm{K}(\left\{e^{i}\right\})\quad \forall \,\, \xi_1, \xi_2 \in \Sigma\, .
\end{equation*}

\noindent
Write $V^{\ast} = V^{\ast}_{+} \oplus V^{\ast}_{-}$, where $V^{\ast}_{+}$ is a $p$-dimensional subspace of $V^{\ast}$ on which $h^{\ast}$ is positive definite and $V^{\ast}_{-}$ is a $q$-dimensional subspace of $V^{\ast}$ on which $h^{\ast}$ is negative-definite. Fix an orientation on $V^{\ast}_{+}$, which induces a unique orientation on $V^{\ast}_{-}$ compatible with the orientation of $V^{\ast}$ induced from that of $V$, and denote by $\nu_{+}$ and $\nu_{-}$ the corresponding pseudo-Riemannian volume forms. We have $\nu = \nu_{+} \wedge \nu_{-}$.  For $p$ (and hence $q$) odd, define:
\begin{equation}
	\label{eq:Bpmodd}
	\cB_{\pm}(\xi_1 , \xi_2) = (\gamma(\nu_{\pm})(\xi_1),\xi_2) \quad \forall\,\, \xi_1, \xi_2 \in \Sigma
\end{equation}

\noindent
whereas for $p$ (and hence $q$) even, set:
\begin{equation}
	\label{eq:Bpmeven}
	\cB_{\pm}(\xi_1 , \xi_2) = (\gamma(\nu_{\mp})(\xi_1),\xi_2) \quad\forall \,\, \xi_1, \xi_2 \in \Sigma\, .
\end{equation}

\noindent
Then $\cB_{\pm}$ are the desired admissible pairings.
\end{proof}

\begin{definition}
The sign factor $\sigma$ appearing in the previous theorem is called the {\em adjoint type} of $\cB_{\sigma}$, hence $\cB_{+}$ is of positive adjoint type whereas $\cB_{-}$ is of negative adjoint type. 
\end{definition}

\noindent 
The relations appearing in \eqref{eq:admissiblepairins} can be equivalently written as:
\begin{equation}
\label{gammat}
\gamma(x)^t=\gamma((\pi^{\frac{1-\sigma}{2}}\circ \tau)(x)) \quad \forall x\in \Cl(V^{\ast},h^{\ast})\, , 
\end{equation}

\noindent
where the superscript $(-)^t$ denotes the adjoint with respect to $\cB_{\sigma}$. The symmetry type of an admissible bilinear form $\cB_{\sigma}$ will be denoted by $s\in \mathbb{Z}_2$. If $s = +1$ then $\cB$ is symmetric whereas if $s = -1$ then $\cB$ is skew-symmetric. Notice that $s$ depends both on $\sigma$ and on the modulo $4$ reduction of $\frac{d}{2}$.

\begin{definition}
A real {\em paired irreducible Clifford module} for $(V^\ast,h^\ast)$ is a triplet $(\Sigma,\gamma,\cB)$, where $(\Sigma,\gamma)$ is a irreducible $\Cl(V^\ast, h^\ast)$-module and $\cB$ is an admissible pairing on $(\Sigma,\gamma)$. We say that $(\Sigma,\gamma,\cB)$ has adjoint type $\sigma\in \mathbb{Z}_2$ and symmetry type $s\in \mathbb{Z}_2$ if $\cB$ has these adjoint and symmetry types, respectively.
\end{definition}

\begin{remark}
\label{rem:cBrelation} 
Admissible bilinear pairings of positive and negative adjoint types are related through the pseudo-Riemannian volume form $\nu$ of $(V^{\ast},h^{\ast})$:
\begin{equation}
\label{eq:cBpm}
\cB_{+} = C\,\cB_{-}\circ (\gamma(\nu)\otimes \Id)\, , 
\end{equation}

\noindent
for an appropriate non-zero real constant $C$. In applications, we will choose to work with $\cB_{+}$ or with $\cB_{-}$ depending on which admissible pairing yields the computationally simplest polyform associated to a given spinor $\xi\in \Sigma$. When $pq = 0$, we will take $\cB_{\sigma}$ to be positive-definite, which can always be achieved by rescaling it with a non-zero constant of appropriate sign. We refer to \cite{LazaroiuBC} for a useful discussion of properties of admissible pairings in various dimensions and signatures.
\end{remark}

\begin{remark}
Directly from their definition, the pairings $\cB_{+}$ and $\cB_{-}$ satisfy:
\begin{equation*}
\cB_{\sigma}(\gamma(\pi^{\frac{1+\sigma}{2}}(x))(\xi_1),\gamma(x)(\xi_2)) = N(x) \cB_{\sigma}( \xi_1, \xi_2) \quad \forall\,\, x\in\Cl(V^{\ast},h^{\ast})\, , \quad \forall\,\, \xi_1, \xi_2 \in \Sigma\, .
\end{equation*} 

\noindent
This relation yields:
\begin{equation*}
\cB_{\sigma}(\gamma(x)(\xi_1),\xi_2) + \cB_{\sigma}(\xi_1,\gamma(x)(\xi_2)) = 0 \quad \forall\,\,\xi_1 , \xi_2 \in \Sigma
\end{equation*}

\noindent
for all $x = \theta_1\cdot\theta_2$ for orthonormal $\theta_1, \theta_2 \in V^{\ast}$. This implies that $\cB_{\sigma}$ is invariant under the action of $\Spin_o(V^{\ast},h^{\ast})$. If $h$ is positive-definite, then $\cB_{+}$ is $\Pin(V^{\ast},h^{\ast})$-invariant, since it satisfies:
\begin{equation*}
\cB_{+}(\gamma(\theta)(\xi_1),\gamma(\theta)(\xi_2)) = \cB_{+}(\xi_1,\xi_2) \quad \forall \xi_1, \xi_2 \in \Sigma
\end{equation*}

\noindent
for all $\theta\in V^{\ast}$ of unit norm. On the other hand, if $h$ is negative-definite, then $\cB_{-}$ is $\Pin(V^{\ast},h^{\ast})$-invariant.
\end{remark}

\noindent
A direct computation using equations \eqref{eq:Bpmodd} and \eqref{eq:Bpmeven} gives the following result, which fixes the constant $C$ appearing in Remark \ref{rem:cBrelation}.

\begin{prop}
The admissible pairings $\cB_{+}$ and $\cB_{-}$ constructed above are related as follows:
\begin{equation}
\label{eq:cB_pm}
\cB_{+} = (-1)^{[\frac{q}{2}]} \cB_{-}(\gamma(\nu)\otimes \Id)\, .
\end{equation}

\noindent
Thus we can normalize $\cB_\pm$ such that the constant in \eqref{eq:cBpm} is given by $C= (-1)^{[\frac{q}{2}]}$.
\end{prop}


\subsection{The K\"ahler-Atiyah model of $\Cl(V^\ast,h^\ast)$}


To identify spinors with polyforms, we will use an explicit realization of $\Cl(V^{\ast},h^{\ast})$ as a deformation of the exterior algebra $\wedge V^{\ast}$. This model, which can be traced back to the work of Chevalley and Riesz \cite{Chevalley1,Chevalley2,Riesz}, has an interpretation as a deformation quantization of the odd symplectic vector space obtained by parity change from the quadratic space $(V,h)$ \cite{Berezin,Meinrenken,Voronov}. It can be constructed using the \emph{symbol map} and its inverse, the \emph{quantization map}. Consider first the linear map $\bar{f}\colon V^{\ast}
\to \End(\wedge V^{\ast})$ given by:
\begin{equation*}
\bar{f}(\theta)(\alpha) = \theta\wedge \alpha + \iota_{\theta^{\sharp}} \alpha \quad \forall\,\, \theta\in V^{\ast} \quad \forall \,\, \alpha \in \wedge V^{\ast}\, .
\end{equation*}

\noindent
We have:
\begin{equation*}
\bar{f}(\theta)\circ \bar{f}(\theta) = h^{\ast}(\theta,\theta) \quad  \forall\,\, \theta \in V^{\ast}\, .
\end{equation*}

\noindent
By the universal property of Clifford algebras, it follows that $\bar{f}$ extends to a unique morphism $f \colon \Cl(V^{\ast},h^{\ast}) \to \End(\wedge V^{\ast})$ of unital associative algebras such that $f\circ \iota = \bar{f}$, where $\iota\colon V^{\ast}\hookrightarrow \Cl(V^{\ast},h^{\ast})$ is the canonical inclusion of $V^{\ast}$ in $\Cl(V^{\ast},h^{\ast})$.

\begin{definition}
The {\em symbol map} is the linear map $\mathfrak{l} \colon \Cl(V^{\ast},h^{\ast}) \to \wedge V^{\ast}$ defined by:
\begin{equation*}
\mathfrak{l}(x) = f(x)(1)\quad \forall\,\, x\in \Cl(V^\ast,h^\ast)\, , 
\end{equation*}

\noindent
where $1\in \R$ is viewed as an element of $\wedge^0(V^{\ast}) = \R$.
\end{definition}

\noindent
The symbol map is an isomorphism of filtered vector spaces. We have:
\begin{equation*}
\mathfrak{l}(1) = 1\, , \quad  \mathfrak{l}(\theta) = \theta\, , \quad  \mathfrak{l}(\theta_1\cdot \theta_2) = \theta_1\wedge \theta_2 + h^{\ast}(\theta_1,\theta_2) \quad \forall \,\,\theta, \theta_1, \theta_2  \in V^{\ast}\, . 
\end{equation*}

\noindent
As expected, $\mathfrak{l}$ is not a morphism of algebras. The inverse:
\begin{equation*}
\Psi :=  \mathfrak{l}^{-1}\colon \wedge V^{\ast} \to \Cl(V^{\ast},h^{\ast})\, . 
\end{equation*}

\noindent
of $\mathfrak{l}$, called the {\em quantization map} of $\wedge V^{\ast}$, allows one to view  $\Cl(V^{\ast},h^{\ast})$ as a deformation of the exterior algebra $(\wedge V^\ast,\wedge)$, see \cite{Berezin,Voronov} for more details. Using $\mathfrak{l}$ and $\Psi$, we transport the algebra product of $\Cl(V^{\ast},h^{\ast})$ to an $h$-dependent unital associative product defined on $\wedge V^{\ast}$, which deforms the wedge product.

\begin{definition}
The {\em geometric product} $\diamond_h :\wedge V^{\ast}\times \wedge V^{\ast}\rightarrow \wedge V^{\ast}$ is defined through:
\begin{equation*}
\alpha_1 \diamond_h \alpha_2  :=  \mathfrak{l}(\Psi(\alpha_1)\cdot \Psi(\alpha_2)) \quad \forall\,\, \alpha_1,\alpha_2\in \wedge V^\ast 
\end{equation*}

\noindent
where $\cdot$ denotes multiplication in $\Cl(V^\ast,h^\ast)$.
\end{definition}

\noindent
By definition, the map $\Psi$ is an isomorphism of unital associative real algebras from $(\wedge V^{\ast},\diamond)$ to $\Cl(V^{\ast},h^{\ast})$. Through this isomorphism, the inclusion $V^{\ast}\hookrightarrow \Cl(V^{\ast},h^{\ast})$ corresponds to the natural inclusion
$V^{\ast}\hookrightarrow \wedge V^{\ast}$. We shall refer to $(\wedge V^{\ast}, \diamond)$ as the \emph{K\"ahler-Atiyah algebra} of the quadratic space $(V,h)$ \cite{Kahler,Graf}. It is easy to see that the geometric product satisfies:
\begin{equation*}
\theta\diamond \alpha = \theta \wedge \alpha + \iota_{\theta^{\sharp}} \alpha \quad \forall\,\, \theta\in
V^{\ast}\quad \forall\,\, \alpha\in \wedge V^{\ast}\, . 
\end{equation*}

\noindent
Also notice the relation:
\begin{equation*}
\theta\diamond\theta = h^{\ast}(\theta,\theta) \quad \forall\,\, \theta\in V^\ast\, . 
\end{equation*}

\noindent
The maps $\pi$ and $\tau$ transfer through $\Psi$ to the K\"ahler-Atiyah algebra, producing unital (anti)-automorphisms of the latter which we denote by the same symbols. With this notation, we have:
\begin{equation}
\label{pitauPsi}
\pi\circ \Psi=\Psi\circ \pi\, , \quad \tau\circ \Psi=\Psi\circ \tau\, . 
\end{equation}

\noindent
For any orthonormal basis $\left\{ e^i \right\}_{i=1,\ldots,d}$ of $V^{\ast}$ and any $k\in \{1,\ldots,d\}$, we have $e^{1}\diamond \dots\diamond e^{k}=e^{1}\wedge\dots\wedge e^{k}$ and:
\begin{equation*}
\pi(e^{1}\wedge\dots\wedge e^{k}) = (-1)^{k} e^{1}\wedge\dots\wedge e^{k}\, ,\quad \tau(e^{1}\wedge\dots\wedge e^{k}) = e^{k}\wedge\dots\wedge e^{1}\, .  
\end{equation*}

\noindent
Let $\cT(V^{\ast})$ denote the tensor algebra of the parity change of $V^{\ast}$, viewed as a $\Z$-graded associative super-algebra whose $\Z_2$-grading is the reduction of the natural $\Z$-grading; thus elements of $V$ have integer degree one and they are odd. Let:
\begin{equation*}
\Der(\cT(V^{\ast}))  :=  \bigoplus_{k\in \Z} \Der^k(\cT(V^{\ast}))
\end{equation*}

\noindent
denote the $\Z$-graded Lie superalgebra of all superderivations. The minus one integer degree component $\Der^{-1}(\cT(V^{\ast}))$ is linearly isomorphic with the space $\Hom(V^{\ast},\R) = V$ acting by contractions:
\begin{equation*}
\iota_v (\theta_1\otimes \dots \otimes \theta_k) = \sum_{i=1}^{k} (-1)^{i-1} \theta_1\otimes \dots \otimes \iota_v\theta_i\otimes \dots \otimes \theta_k \quad \forall\,\, v\in V\quad \forall \,\,\theta_1,\hdots , \theta_k\in V^{\ast}
\end{equation*}

\noindent
while the zero integer degree component $\Der^0(\cT(V^{\ast})) = \End(V^{\ast}) = \mathfrak{gl}(V^{\ast})$ acts through:
\begin{equation*}
L_A(\theta_1\otimes \dots \otimes \theta_k) = \sum_{i=1}^{k} \theta_1\otimes \dots \otimes A(\theta_i)\otimes \dots \otimes \theta_k \quad \forall A\in \mathfrak{gl}(V^{\ast})\, . 
\end{equation*}

\noindent
We have an isomorphism of super-Lie algebras:
\begin{equation*}
\Der^{-1}(\cT(V^{\ast})) \oplus \Der^0(\cT(V^{\ast})) \simeq V\rtimes \mathfrak{gl}(V^{\ast})\, . 
\end{equation*}

\noindent
The action of this super Lie algebra preserves the ideal used to define the exterior algebra as a quotient of $\cT(V^{\ast})$ and hence descends to a morphism of super Lie algebras $\cL_{\Lambda}\colon V\rtimes \mathfrak{gl}(V^{\ast}) \to \Der(\wedge V^{\ast},\wedge)$. Contractions also preserve the ideal used to define the Clifford algebra as a quotient of $\cT(V^{\ast})$. On the other hand, endomorphisms $A$ of $V^{\ast}$ preserve that ideal if and only if $A\in\mathfrak{so}(V^{\ast},h^{\ast})$. Together with contractions, they induce a morphism of super Lie algebras $\cL_{\Cl} \colon V\rtimes \mathfrak{so}(V^{\ast},h^{\ast}) \to \Der(\Cl(V^{\ast},h^{\ast}))$. The following result states that $\cL_{\Lambda}$ and $\cL_{\Cl}$ are compatible with $\mathfrak{l}$ and $\Psi$.

\begin{prop}{\rm \cite[Proposition 2.11]{Meinrenken}}
\label{prop:equivariancePsi}
The quantization and symbol maps intertwine the natural action of $V\rtimes \mathfrak{so}(V^{\ast},h^{\ast})$ on $\Cl(V^{\ast},h^{\ast})$ and $\wedge V^{\ast}$, that is, the following holds:
\begin{equation*}
\Psi(\cL_{\Lambda}(\varphi)(\alpha)) = \cL_{\Cl}(\varphi)(\Psi(\alpha))\, , \quad \mathfrak{l}(\cL_{\Cl}(\varphi)(x)) = \cL_{\Lambda}(\varphi)(\mathfrak{l}(x))\, , 
\end{equation*}

\noindent
for all $\varphi \in V\rtimes \mathfrak{so}(V^{\ast},h^{\ast})$, $\alpha  \in \wedge V^{\ast}$ and $x\in \Cl(V^{\ast},h^{\ast})$.
\end{prop}

\noindent
This proposition shows that quantization is equivariant with respect to affine orthogonal transformations of $(V^\ast, h^\ast)$. In signatures $p-q\equiv_8 0, 2$, composing $\Psi$ with the irreducible representation $\gamma\colon \Cl(V^{\ast},h^{\ast}) \to \End(\Sigma)$ (which in such signatures is a unital isomorphism of algebras) gives an isomorphism of unital associative real algebras:
\begin{equation}
\label{Psigamma}
\Psi_{\gamma} :=  \gamma\circ\Psi\colon (\wedge V^{\ast}, \diamond_h )\ \stackrel{\sim}{\to} (\End(\Sigma),\circ)\, . 
\end{equation}

\noindent
fitting into the following commutative diagram:
\begin{center}
\begin{tikzcd}
\Cl(V^{\ast},h^{\ast})  \arrow[r,"\gamma"]   & (\End(\Sigma),\circ)  \\
(\wedge V^{\ast}, \diamond_h) \arrow[u,"\Psi"] \arrow[ur,"\Psi_{\gamma}"]	& 
\end{tikzcd}
\end{center}

\noindent
Since $\Psi_{\gamma}$ is an isomorphism of algebras and $(\wedge V^{\ast}, \diamond_h)$ is generated by $V^{\ast}$, the identity together with the elements $\Psi_\gamma(e^{i_1}\wedge \ldots \wedge e^{i_k})=\gamma(e^{i_1})\circ \hdots \circ \gamma(e^{i_k})$ for $1\leq i_1<\hdots <i_k\leq d$ and $k=1,\ldots, d$ form a basis of $\End(\Sigma)$. In accordance with the terminology that refers to $\Psi \colon \wedge V^{\ast} \to \Cl(V^{\ast},h^{\ast})$ as the \emph{quantization map}, we introduce the following definition.

\begin{definition}
The \emph{dequantization} of an endomorphism $A\in \End(\Sigma)$ is the polyform determined by $\fra := \Psi_{\gamma}^{-1}(A)\in \wedge V^{\ast}$.
\end{definition}

\begin{remark}
For ease of notation we will sometimes denote the action of a polyform $\alpha \in \wedge V^{\ast}$ as an endomorphism on $\Sigma$ thorugh $\Psi_\gamma$ simply by a dot, namely:
\begin{equation*}\alpha\cdot\xi  :=  \Psi_{\gamma}(\alpha)(\xi) \quad \forall\,\, \alpha\in \wedge V^{\ast}\quad \forall\,\, \xi\in \Sigma\, . 
\end{equation*} 
\end{remark}

\noindent 
The trace on $\End(\Sigma)$ transfers to the K\"ahler-Atiyah algebra through the isomorphism $\Psi_\gamma$:

\begin{definition}
The {\em K\"ahler-Atiyah trace} is the linear functional:
\begin{equation*}
\cT r\colon \wedge V^{\ast}  \to \R\, , \quad  \alpha\mapsto \mathrm{Tr}(\Psi_{\gamma}(\alpha)) 
\end{equation*}

\noindent
where $\mathrm{Tr}(-)$ denotes the trace in $\End(\Sigma)$.
\end{definition}

\noindent
We will see in a moment that $\cT r$ does not depend on $\gamma$ or $h$. Since $\Psi_{\gamma}$ is a unital morphism of algebras, we have:
\begin{equation*}
\cT r(1) = \dim(\Sigma)= 2^{\frac{d}{2}}~~\mathrm{and}~~\cT r(\alpha_1\diamond \alpha_2) =  \cT r(\alpha_2\diamond \alpha_1) \quad \forall \,\, \alpha_1,\alpha_2\in \wedge V^\ast
\end{equation*}

\noindent
where $1\in \R=\wedge^0 V^\ast$ is the unit element of the field $\R$ of real numbers.

\begin{lemma}
\label{lemma:tracev}
For any $0< k \leq d$, we have $\cT r\vert_{\wedge^{k} V^{\ast}} = 0$.
\end{lemma}

\begin{proof}
Let $( e^1, \hdots , e^n )$ be an orthonormal basis of $(V^{\ast},h^{\ast})$. For $i\neq j$ we have $e^i\diamond e^j = - e^j\diamond e^i$ and hence $(e^i)^{-1}\diamond e^j\diamond e^i = -e^j$. Let $0 < k \leq d$ and $1\leq i_1 < \cdots < i_k \leq d$. If $k$ is even, then:
\begin{equation*}
\cT r(e^{i_1}\diamond \dots \diamond e^{i_{k}}) = \cT r(e^{i_k}\diamond e^{i_1}\diamond \dots \diamond e^{i_{k-1}}) = (-1)^{k-1} \cT r(e^{i_1}\diamond \dots \diamond e^{i_{k}})\, , 
\end{equation*}

\noindent
and hence $\cT r(e^{i_1}\diamond \dots \diamond e^{i_{k}}) = 0$. Here we used cyclicity of the K\"ahler-Atiyah trace and the fact that $e^{i_{k}}$ anticommutes with $e^{i_1}, \dots, e^{i_{k-1}}$. If $k$ is odd, let $j\in \left\{1,\dots,d\right\}$ be such that $j\not \in \left\{i_1,\dots,i_k\right\}$ (such a $j$ exists since $k<d$). We have:
\begin{equation*}
\cT r(e^{i_1}\diamond \dots \diamond e^{i_{k}}) = \cT r(( e^j)^{-1} \diamond e^{i_1}\diamond \dots 
\diamond e^{i_{k}}\diamond e^j) = -  \cT r( e^{i_1}\diamond \dots \diamond e^{i_{k}}) = 0
\end{equation*}

\noindent
and we conclude.
\end{proof}

\noindent 
Let $\alpha^{(k)} \in \wedge^k V^{\ast}$ denote the degree $k$ component of $\alpha\in \wedge V^\ast$. Lemma \ref{lemma:tracev} immediately implies the following result.
\begin{prop}
\label{prop:cS}
The K\"ahler-Atiyah trace is given by:
\begin{equation*}
\cT r(\alpha) = \dim (\Sigma)\, \alpha^{(0)} =2^{\frac{d}{2}} \alpha^{(0)} \quad \forall\,\, \alpha\in \wedge V^\ast\, . 
\end{equation*}

\noindent
In particular, $\cT r$ does not depend on the irreducible representation $\gamma$ of $\Cl(V^\ast,h^\ast)$ or on $h$.
\end{prop}

\begin{lemma}
\label{lemma:adjointpoly}
Let $\alpha\in \wedge V^{\ast}$ and $\cB$ be an admissible bilinear pairing of $(\Sigma,\gamma)$ of adjoint type $\sigma \in \mathbb{Z}_2$. Then the following equation holds:
\begin{equation}
\label{eq:traceD}
\Psi_{\gamma}(\alpha)^{t} = \Psi_{\gamma}((\pi^{\frac{1-\sigma}{2}}\circ\tau)(\alpha))\, , 
\end{equation}

\noindent
where $\Psi_{\gamma}(\alpha)^{t}$ is the $\cB$-adjoint of $\Psi_{\gamma}(\alpha)$.
\end{lemma}

\begin{proof}
Follows directly from \eqref{gammat} and relations \eqref{pitauPsi}.
\end{proof}


\subsection{Spinor square maps}


In this subsection we introduce one of the most important concepts of this dissertation, namely that of the algebraic \emph{square of a spinor}, which is obtained in terms of the \emph{spinor square maps}, which we proceed to introduce. 

\begin{definition}
\label{def:squarespinor}
Let $(\Sigma,\gamma,\cB)$ be a paired irreducible Clifford module for $(V^\ast, h^\ast)$. The {\em signed spinor square maps} of $\Sigma$ are the quadratic maps:
\begin{equation*}
\label{eq:spinorsquaremap}
\cE^{\gamma}_{\kappa}  :=  \Psi_{\gamma}^{-1} \circ \cE_{\kappa}: \Sigma \to \wedge V^{\ast}\, , \quad \kappa\in \mathbb{Z}_2
\end{equation*}

\noindent
where $\cE_{\kappa}\colon \Sigma\to \End(\Sigma)$ is the signed square maps of the paired vector space $(\Sigma,\cB)$, which was introduced in Definition \ref{def:squarevectorspace}. 
\end{definition}

\noindent
Given a spinor $\xi\in \Sigma$, the polyforms $\cE^{\gamma}_{+}(\xi)$ and $\cE^{\gamma}_{-}(\xi)$ are respectively called the positive and negative {\em squares} of $\xi$ relative to the admissible pairing $\cB$. A polyform $\alpha\in \wedge V^\ast$ is called a {\em signed square} of $\xi\in \Sigma$ if either $\alpha=\cE^{\gamma}_{+}(\xi)$ or $\alpha=\cE_{-}^{\gamma}(\xi)$. The spinor square maps fit into the following commutative diagram.

\begin{center}
\begin{tikzcd}
\Cl(V^{\ast},h^{\ast})  \arrow[r,"\gamma"]   & (\End(\Sigma),\circ)  \arrow[dl,"\Psi^{-1}_{\gamma}"] & \arrow[l,"\cE_{\kappa}"] \Sigma   \\
(\wedge V^{\ast}, \diamond) \arrow[u,"\Psi"]  	& 
\end{tikzcd}
\end{center}

\noindent 
Since $\Psi_{\gamma}$ is a linear isomorphism, the results of Section \ref{sec:vectorasendo} imply that $\cE_{\kappa}$ is two-to-one except at $0\in \Sigma$. Furthermore:
\begin{equation*}
\Im (\cE^{\gamma}_{-})-=-\Im (\cE^{\gamma}_{+})\, , \qquad \Im (\cE^{\gamma}_{+})\cap \Im (\cE^{\gamma}_{-})=\{0\}\, . 
\end{equation*}

\noindent
We define:
\begin{eqnarray*}
\Im (\cE^{\gamma}) := \Im (\cE^{\gamma}_{+})\cup \Im (\cE^{\gamma}_{-})\, .
\end{eqnarray*}

\noindent
Moreover, both $\cE^{\gamma}_{-}$ and $\cE^{\gamma}_{+}$ induce the same bijective map on the $\mathbb{Z}_2$ quotients:
\begin{equation}
\label{eq:hcEspinor}
\cE^{\gamma}\colon \Sigma/\Z_2 \stackrel{\sim}{\rightarrow} \Im (\cE^{\gamma})/\Z_2\, . 
\end{equation}

\noindent
Notice that $\Im (\cE^{\gamma})$ is a cone in $\wedge V^\ast$, which is the union of the opposite half cones $\Im (\cE^{\gamma}_{+})$ and $\Im (\cE^{\gamma}_{-})$. For simplicity in the exposition we will sometimes denote by $\alpha_{\xi} :=  \cE_+^{\gamma}(\xi)$ the positive polyform square of $\xi\in \Sigma$. Polyforms in the image of of the spinor square map $\cE^{\gamma}$ will be generically called \emph{spinorial exterior forms}.

\begin{remark}
The representation map $\gamma$ is an isomorphism when $p-q\equiv_8 0, 2$. This does not hold in other signatures, for which the construction of spinor square maps is more delicate, see \cite{LazaroiuBC} for more details.
\end{remark}

\noindent
The following result is a direct consequence of Proposition \ref{prop:equivariancePsi}.

\begin{prop}
\label{prop:equivariancesquaremap}
The quadratic map $\cE_{\kappa}^{\gamma}\colon \Sigma \to \wedge V^{\ast}$ is $\Spin_o(V^{\ast},h^{\ast})$ equivariant, that is:
\begin{equation*}
\cE_{\kappa}^{\gamma}(u\,\xi) = \Ad_u(\cE_{\kappa}^{\gamma}(\xi)) \quad \forall\,\, u\in \Spin_o(V^{\ast},h^{\ast}) \quad \forall\,\, \xi\in\Sigma 
\end{equation*}

\noindent
where the right hand side denotes the natural action of $\Ad_u\in\SO(V^{\ast},h^{\ast})$ on $\wedge V^{\ast}$.
\end{prop}

\noindent 
We are ready to give the algebraic characterization of spinors in terms of polyforms.

\begin{thm}
\label{thm:reconstruction} 
Let $(V^\ast, h^\ast)$ be a quadratic vector space of signature $(p-q) \equiv_8 0,2$, and let $(\Sigma,\gamma,\cB)$ be an associated paired irreducible Clifford module of symmetry type $s$ and adjoint type $\sigma$. Then the following statements are equivalent for a polyform $\alpha\in \wedge V^{\ast}$:
\begin{enumerate}
\itemsep 0.0em

\item $\alpha$ is a signed square of some spinor $\xi\in \Sigma$, that is, it lies in the set $\Im(\cE^{\gamma})$.

\item $\alpha$ satisfies the following relations:
\begin{equation}
\label{eq:thmdefequationsequiv}
\alpha\diamond\alpha =2^{\frac{d}{2}} \alpha^{(0)} \alpha \, , \quad (\pi^{\frac{1-\sigma}{2}}\circ\tau)(\alpha) = s\,\alpha\, , \quad  \alpha\diamond \beta\diamond\alpha = 2^{\frac{d}{2}} (\alpha\diamond\beta)^{(0)} \alpha
\end{equation}
	
\noindent	
for a fixed polyform $\beta \in \wedge V^{\ast}$ satisfying $\cT r(\alpha\diamond\beta) \neq 0$.

\item The following relations hold:
\begin{equation}
\label{eq:thmdefequations}
(\pi^{\frac{1-\sigma}{2}}\circ\tau)(\alpha) = s\,\alpha \, , \quad \alpha\diamond \beta\diamond\alpha = \cT r(\alpha\diamond\beta)\, \alpha
\end{equation}

\noindent
for every polyform $\beta\in \wedge V^\ast$.
\end{enumerate}

\noindent
In particular, the set $\Im(\cE^{\gamma})$ depends only on $\sigma$, $s$ and $(V^\ast, h^\ast)$. 
\end{thm}
 
\begin{remark}
In view of this result, we will occasionally denote $\Im(\cE^{\gamma})$ by $Z_{\sigma,s}(V^\ast,h^\ast)$, and $\Im(\cE^{\gamma}_{\kappa})$ by $Z^{\kappa}_{\sigma,s}(V^\ast,h^\ast)$.
\end{remark}

\begin{proof}
Since $\Psi\colon \Cl(V^{\ast},h^{\ast}) \to \End(\Sigma)$ is a unital isomorphism of associative algebras, $\alpha$ satisfies \eqref{eq:thmdefequations} if and only if:
\begin{equation}
\label{eq:Esquarethm}
E^t = \sigma\,E\, , \quad E\circ A\circ E = \tr(E\circ A) E \quad \forall\, A\in \End(\Sigma) 
\end{equation}

\noindent
where $E  :=  \Psi_\gamma(\alpha)$, $A :=  \Psi_\gamma(\beta)$ and we have used Lemma \ref{lemma:adjointpoly} together with the definition and properties of the \KA trace. The conclusion now follows from Proposition \ref{prop:characterizationtamecone}.
\end{proof}

\noindent
The second equation in \eqref{eq:thmdefequations} implies the following result.

\begin{cor}
\label{cor:signcriteria}
Let $\alpha \in Z_{\sigma,s}(V^\ast,h^\ast)$. If $k\in \{1,\ldots, d\}$ satisfies:
\begin{equation*}
(-1)^{k \frac{1-s}{2}} (-1)^{\frac{k(k-1)}{2}} = -\sigma\, , 
\end{equation*}

\noindent
then $\alpha^{(k)} = 0$.
\end{cor}

\noindent
Polyforms $\alpha\in Z_{\sigma,s}(V^\ast,h^\ast)$ admit an explicit presentation, see for instance \cite{LazaroiuB, LazaroiuBII, LazaroiuBC}.

\begin{prop}
Let $(e^1,\hdots , e^n)$ be an orthonormal basis of $(V^{\ast},h^{\ast})$ and let $\kappa\in \mathbb{Z}_2$. Then every polyform $\alpha\in Z_{\sigma,s}(V^\ast,h^\ast)$ can be written as:
\begin{equation}
\label{eq:bilinears}
\alpha = \frac{\kappa}{2^{\frac{d}{2}}} \sum_{k=0}^{d}  \,\sum_{i_1 < \dots < i_k} \cB((\gamma(e^{i_k})^{-1} \circ \dots  \circ \gamma(e^{i_1})^{-1})(\xi),\xi)\, e^{i_1}\wedge \hdots \wedge e^{i_k}\, , 
\end{equation}

\noindent
where the spinor $\xi\in \Sigma$ is determined by $\alpha$ up to sign.
\end{prop}

\begin{remark}
\label{rem:Dirac}
We have:
\begin{equation*}
\gamma(e^i)^{-1} = h^{\ast}(e^i,e^i) \gamma(e^i) = h(e_{i},e_{i}) \gamma(e^i)\, , 
\end{equation*}

\noindent
where $( e_1, \hdots, e_n)$ is the contragradient orthonormal basis of $(V,h)$. For simplicity, set:
\begin{equation*}
\gamma^i  :=  \gamma(e^i)~~ \mathrm{and} ~~ \gamma_i  :=  h(e_i,e_i)\gamma(e^i)\, , 
\end{equation*}

\noindent
so that $(\gamma^i)^{-1}=\gamma_i$. Then the degree one component in \eqref{eq:bilinears} reads:
\begin{equation*}
\alpha^{(1)}=\frac{\kappa}{2^{\frac{d}{2}}} \cB(\gamma_i(\xi),\xi)e^i
\end{equation*}

\noindent
and its dual vector field $(\alpha^{(1)})^\sharp=\frac{\kappa}{2^{\frac{d}{2}}} \cB(\gamma_i(\xi),\xi)e_i$ is called the signed {\em Dirac vector} of $\xi$ relative to $\cB$. 
\end{remark}

\begin{proof}
It is easy to see that the set:
\begin{equation*}
P  :=  \left\{\Id\right\}\cup \left\{ \gamma^1 \circ \cdots \circ \gamma^{i_1} \circ \cdots \gamma^{i_k}\circ \cdots \circ \gamma^{d} \,|\,1\leq i_1 < \cdots < i_k \leq d, k=1,\ldots, d \right\} 
\end{equation*}
is an orthogonal basis of $\End(\Sigma)$ with respect to the nondegenerate and symmetric bilinear pairing induced by the trace:
\begin{equation*}
\End(\Sigma)\times \End(\Sigma) \to \R\, , \quad  (A_1,A_2) \mapsto \tr(A_1 A_2)\, . 
\end{equation*}

\noindent
In particular, the endomorphism $E  =  \Psi_{\gamma}(\alpha)\in \Im(\cE)$ expands as:
\begin{eqnarray*}
E = \frac{1}{2^{\frac{d}{2}}}\sum_{k=0}^{d} \sum_{i_1 < \dots < i_k} \tr((\gamma^{i_1}\circ \cdots \circ \gamma^{i_k})^{-1}\circ E)\, \gamma^{i_1}\circ \cdots \circ \gamma^{i_k} \\ 
= \frac{\kappa}{2^{\frac{d}{2}}}\sum_{k=0}^{d} \sum_{i_1 < \dots < i_k} \cB((\gamma^{i_1}\circ \cdots \circ \gamma^{i_k})^{-1} (\xi), \xi)\, \gamma^{i_1}\circ \cdots \circ \gamma^{i_k} \, , 
\end{eqnarray*}

\noindent
where $\xi\in \Sigma$ is a spinor such that $E=\cE_\kappa(\xi)$ and we have used that:
\begin{equation*}
\tr(B\circ \cE_\kappa(\xi))=\kappa\, \tr(B(\xi)\otimes \xi^\ast)=\kappa\, \xi^\ast(B(\xi))=\kappa\, \cB(B\xi,\xi)
\end{equation*}

\noindent
for all $B\in \End(\Sigma)$. The conclusion follows now by applying the isomorphism algebras $\Psi_{\gamma}^{-1}:(\End(\Sigma),\circ)\rightarrow (\wedge V^\ast, \diamond)$ to the previous equation.
\end{proof}

\begin{lemma}
\label{lemma:actionnu}
The following identities hold for all $\alpha \in \wedge V^{\ast}$:
\begin{equation}
\label{eq:nuaction}
\alpha \diamond\nu  = \ast\, \tau(\alpha)\, , \qquad \nu \diamond \alpha =  \ast\, (\pi\circ\tau) (\alpha)\, . 
\end{equation}
\end{lemma}

\begin{proof}
Since multiplication by $\nu$ is $\R$-linear, it suffices to prove the statement for homogeneous elements $\alpha=e^{i_1}\wedge \cdots \wedge e^{i_k}$ with $1\leq i_1 < \cdots < i_k \leq d$, where $(e^1, \hdots, e^n)$ is an orthonormal basis of $(V^{\ast},h^{\ast})$. We have:
\begin{eqnarray*}
& e^{i_1}\wedge \cdots \wedge e^{i_k} \diamond \nu =  e^{i_1}\diamond \cdots \diamond e^{i_k} \diamond e^1\diamond\cdots \diamond e^d \\ 
& = (-1)^{i_1 + \cdots + i_k} (-1)^k e^1\diamond\cdots \diamond (e^{i_1})^2 \diamond e^{i_1+1}\diamond \cdots \diamond (e^{i_k})^2\diamond e^{i_k + 1} \diamond \cdots \diamond e^d\\
& = h^{\ast}(e^{i_1},e^{i_1})\cdots h^{\ast}(e^{i_k},e^{i_k})\,(-1)^{i_1 + \cdots + i_k} (-1)^k e^1\diamond \cdots \diamond e^{i_1-1}\diamond e^{i_1+1}\cdots \diamond e^{i_k-1} \diamond e^{i_k+1}\diamond \cdots\diamond e^d\\
& = (-1)^{\frac{k(k - 1)}{2}}(-1)^{2(i_1 + \cdots + i_k)} (-1)^{2k} \ast(e^{i_1}\wedge  \cdots \wedge e^{i_k} ) = \ast \tau(\alpha)\, , 
\end{eqnarray*}

\noindent
which implies $ \alpha \diamond \nu = \ast\, \tau (\alpha)$. Using the relation $\alpha \diamond \nu = (\nu\diamond\pi)(\alpha)$, we conclude.
\end{proof}

\noindent 
The following shows that the choice of admissible pairing used to construct the spinor square map is a matter of taste, see also Remark \ref{rem:cBrelation}.

\begin{prop} 
Let $\xi\in \Sigma$ and denote by $\alpha_{\xi}^\pm \in Z_+$ the {\em positive} polyform squares of $\xi$ relative to the admissible pairings $\cB_+$ and $\cB_-$ of $(\Sigma,\gamma)$, which we assume to be normalized such that they are related through \eqref{eq:cB_pm}. Then the following relation holds:
\begin{equation*}
\ast \,\alpha_{\xi}^{+} = (-1)^{[\frac{q+1}{2}]+ p(q+1) } (-1)^{d} c(\alpha_{\xi}^{-})\, . 
\end{equation*}

\noindent	
where $c\colon \wedge V^{\ast} \to \wedge V^{\ast}$ is the linear map which acts as multiplication by $\frac{k!}{ (d - k)!}$ in each degree $k$.
\end{prop}

\begin{proof}
We compute:
\begin{eqnarray*}
& \ast (\alpha_{\xi}^{+})^{(k)} = \frac{1}{2^{\frac{d}{2}}} \cB_{+}((\gamma_{i_k}\circ \hdots \circ \gamma_{i_1})(\xi), \xi)   \ast (e^{i_1}\wedge \hdots \wedge e^{i_k}) \\ 
& =  (-1)^{[\frac{q+1}{2}]+ pq }  (-1)^{\frac{k (k-1)}{2}}\frac{\sqrt{\vert h\vert }}{2^{\frac{d}{2}} (d - k)!} \cB_{-}((\gamma(\nu)\circ\gamma_{i_1}\circ \hdots \circ \gamma_{i_k})(\xi), \xi) \epsilon^{i_{i_1}\hdots i_k a_{k+1} \hdots a_d} e_{a_{k+1}}\wedge \hdots \wedge e_{a_d} \\ 
& =  (-1)^{[\frac{q+1}{2}]+ pq } (-1)^{\frac{k (k-1)}{2}} (-1)^{k(d-k)}  \frac{k!}{2^{\frac{d}{2}} (d - k)!} \cB_{-}(\gamma(\nu) \gamma(\ast (e^{a_{k+1}}\wedge \hdots \wedge e^{d}))(\xi), \xi)  e_{a_{k+1}}
\wedge \hdots \wedge e_{d} \\ 
& =  (-1)^{[\frac{q+1}{2}]+ pq } (-1)^{\frac{k (k-1)}{2}} (-1)^{\frac{(d-k)(d+k+1)}{2}}  \frac{k!}{2^{\frac{d}{2}} (d - k)!} \cB_{-}((\gamma(\nu)^2 \circ \gamma^{a_{k+1}} \circ \hdots  \circ \gamma^{a_d})(\xi) , \xi)  e_{a_{k+1}}\wedge \hdots \wedge e_{d} \\ 
& = (-1)^{[\frac{q+1}{2}]+ p(q+1) } (-1)^{k} \frac{k!}{ (d - k)!}  (\alpha_{\xi}^{-})^{(d-k)} = (-1)^{[\frac{q+1}{2}]+ p(q+1) } (-1)^{d}\frac{k!}{ (d - k)!}  \pi(\alpha_{\xi}^{-})^{(d-k)}\, , 
\end{eqnarray*}

\noindent
where we used the identity $\nu\diamond \alpha = \ast (\pi \circ\tau)(\alpha)$ proved in Lemma \ref{lemma:actionnu}.
\end{proof}


\subsection{Linear constraints}


Let $(\Sigma,\gamma,\cB)$ be a paired irreducible Clifford module for $(V^{\ast},h^{\ast})$. Given any endomorphism $Q\in \End(\Sigma)$ we will refer to:
\begin{eqnarray*}
\frq \in \Psi^{-1}_{\gamma}(Q) \in \wedge V^{\ast}
\end{eqnarray*}

\noindent
as the \emph{symbol} of $Q\in \End(\Sigma)$.

\begin{prop}
\label{prop:constraintendopoly} 
A spinor $\xi\in \Sigma$ lies in the kernel of an endomorphism $Q\in \End(\Sigma)$ if and only if:
\begin{equation*}
\frq \diamond \alpha_{\xi} = 0\, , 
\end{equation*}

\noindent
where $\alpha_{\xi} :=  \cE^{\gamma}_{+}(\xi)$ is the positive polyform square of $\xi$.
\end{prop}

\begin{remark}
Equation $\frq \diamond \alpha_{\xi} = 0$ is equivalent to 
\begin{equation*}
\alpha_{\xi}\diamond (\pi^{\frac{1-\sigma}{2}}\circ\tau)(\frq) = 0  
\end{equation*}

\noindent
by taking the transpose of equation $Q(\xi)$ with respect to $\cB$.
\end{remark}

\begin{proof} 
Follows as a consequence of Proposition \ref{prop:constraintendo}, using the fact that $\Psi_{\gamma}\colon (\wedge V^{\ast},\diamond)\to \End(\Sigma)$ is an isomorphism of unital associative algebras.
\end{proof}


\subsection{Real chiral spinors}


Theorem \ref{thm:reconstruction} can be refined for chiral spinors of real type, which exist in signature $p-q\equiv_8 0$. In this case, the Clifford volume form $\nu_h\in \Cl(V^{\ast},h^{\ast})$ squares to the identity in $\Cl(V^{\ast},h^{\ast})$ and lies in the center of $\Cl^{\ev}(V^{\ast},h^{\ast})$, giving the following decomposition of the latter as a direct sum of simple associative algebras:
\begin{equation*}
\Cl^{\ev}(V^{\ast},h^{\ast}) = \Cl^{\ev}_{+}(V^{\ast},h^{\ast}) \oplus \Cl^{\ev}_{-}(V^{\ast},h^{\ast}) 
\end{equation*}

\noindent
where we have defined:
\begin{equation*}
\Cl^{\ev}_{\pm}(V^{\ast},h^{\ast})  := \frac{1}{2}(1 \pm \nu_h )\, \Cl(V^{\ast},h^{\ast})\, . 
\end{equation*}

\noindent
We decompose $\Sigma$ accordingly:
\begin{equation*}
\Sigma = \Sigma^{(+)} \oplus \Sigma^{(-)}\, , \quad \mathrm{where}~~ \Sigma^{(\pm)}  :=  \frac{1}{2} (\Id \pm \gamma(\nu_h))(\Sigma)\, . 
\end{equation*}

\noindent
The subspaces $\Sigma^{(\pm)}\subset\Sigma$ are preserved by the restriction of $\gamma$ to $\Cl^{\ev}(V^{\ast},h^{\ast})$, which therefore decomposes as a sum of two irreducible representations:
\begin{equation*}
\gamma^{(+)} \colon \Cl^{\ev}(V,h)\to \End(\Sigma^{(+)})~~ \mathrm{and} ~~ \gamma^{(-)} 
\colon \Cl^{\ev}(V,h)\to \End(\Sigma^{(-)})
\end{equation*}

\noindent
distinguished by the value which they take on the volume form $\nu_h\in \Cl^{\ev}(V^{\ast},h^{\ast})$:
\begin{equation*}
\gamma^{(+)} (\nu_h) = \Id\, , \quad  \gamma^{(-)} (\nu_h) = -\Id\, . 
\end{equation*}

\noindent
A spinor $\xi \in \Sigma$ is called chiral of chirality $\mu\in\mathbb{Z}_2$ if it belongs to $\Sigma^{(\mu)}$. Setting $\alpha_{\xi} :=  \cE_{+}^{\gamma}(\xi)$, Proposition \ref{prop:constraintendopoly} shows that this amounts to the condition:
\begin{equation*}
\nu_h \diamond \alpha_{\xi} = \mu\, \alpha_{\xi}\, . 
\end{equation*}

\noindent
Given $\mu, \kappa\in \mathbb{Z}_2$, we define:
\begin{equation*}
\cE_{\kappa}^{\gamma \mu} \colon \Sigma^{\mu}\rightarrow \wedge V^{\ast}
\end{equation*}

\noindent
as the restriction of $\cE_{\kappa}^{\gamma} \colon \Sigma \to \wedge V^{\ast}$ to $\Sigma^{\mu}\subset \Sigma$.  We have:
\begin{eqnarray*}
\Im(\cE_{+}^{\gamma \mu}) = - \Im(\cE_{-}^{\gamma \mu})\, , \qquad \Im(\cE_{+}^{\gamma \mu})\cap \Im(\cE_{-}^{\gamma \mu}) = \left\{ 0\right\}
\end{eqnarray*}

\noindent
Theorem \ref{thm:reconstruction}, Proposition \ref{prop:constraintendopoly} and Lemma \ref{lemma:actionnu} all together give:

\begin{cor}
\label{cor:reconstructionchiral}
Let $\Sigma$ be a paired irreducible $\Cl(V^\ast,h^\ast)$-module of symmetry type $s$ and adjoint type $\sigma$. The following statements are equivalent for $\alpha\in \wedge V^{\ast}$, where $\mu\in \mathbb{Z}_2$ is a fixed chirality type:
\begin{enumerate}
\itemsep 0.0em
\item $\alpha$ belongs to $\Im(\cE_{\kappa}^{\gamma \mu})$, that is, it is a signed square of a chiral spinor of chirality $\mu$.
\item The following conditions are satisfied:
\begin{equation}
\label{eq:reconstructionchiral0II}
(\pi^{\frac{1-\sigma}{2}}\circ\tau)(\alpha) = s\, \alpha \, , \quad \ast\, (\pi\circ \tau)(\alpha) = \mu\, \alpha \, , \quad \alpha\diamond \alpha =\cT r(\alpha) \, \alpha\, , \quad \alpha\diamond \beta\diamond \alpha = \cT r(\alpha\diamond\beta)\, \alpha
\end{equation}

\noindent		
for a fixed polyform $\beta \in \wedge V^{\ast}$ which satisfies $\cT r(\alpha\diamond\beta) \neq 0$.
\item The following conditions are satisfied:
\begin{equation}
\label{eq:reconstructionchiral0}
(\pi^{\frac{1-\sigma}{2}}\circ\tau)(\alpha) = s \alpha\, , \quad \ast\, (\pi\circ \tau)(\alpha) = \mu\, \alpha\, , \quad \alpha\diamond \beta\diamond \alpha =\cT r(\alpha\diamond\beta)\, \alpha
\end{equation}

\noindent
for every polyform $\beta \in \wedge V^{\ast}$.
\end{enumerate}

\noindent
In this case, the real chiral spinor of chirality $\mu$ which corresponds to $\alpha$ through the either of the maps $\cE_{+}^{\gamma \mu}$ or $\cE_{-}^{\gamma \mu}$ is unique up to sign and vanishes if and only if $\alpha = 0$.
\end{cor}

\noindent 
In particular, $Z^{\mu}_{\sigma,s}(V^\ast, h^\ast) := \Im(\cE_{+}^{\gamma \mu}) \cup \Im(\cE_{-}^{\gamma \mu})$ depends only on $\mu,\sigma,s \in \mathbb{Z}_2$ and $(V^\ast, h^\ast)$. The following corollary is useful to simplify computations.

\begin{cor}
Let $\alpha\in Z^{+}_{\sigma,s}(V^\ast, h^\ast)\cup Z^{-}_{\sigma,s}(V^\ast, h^\ast)$. If $k\in \{1,\ldots, d\}$ satisfies:
\begin{equation*}
s = - (-1)^{k \frac{\sigma-1}{2}} (-1)^{\frac{k(k-1)}{2}} 
\end{equation*}

\noindent
then we have $\alpha^{(k)} = 0$ and $\alpha^{(d-k)} = 0$. 
\end{cor}

\begin{proof}
Follows from Corollary \ref{cor:signcriteria} and the second relation in \eqref{eq:reconstructionchiral0II}.
\end{proof}


\subsection{Low-dimensional examples}


\noindent 
We describe now $Z_{\sigma,s}(V^\ast, h^\ast)$ and $Z^{\mu}_{\sigma,s}(V^\ast, h^\ast)$ for some low-dimensional cases. 


\subsubsection{Signature $(2,0)$}
\label{sec:Riemanexample}


Let $(V^{\ast},h^{\ast})$ be a two-dimensional real vector space with a scalar product $h^{\ast}$. Its irreducible Clifford module $(\Sigma,\gamma)$ is two-dimensional with an admissible pairing $\cB$ which is symmetric and positive definite. Theorem \ref{thm:reconstruction} with $\beta =1$ shows that $\alpha \in \wedge V^{\ast}$ is a signed square of $\xi\in\Sigma$ if and only if:
\begin{equation}
\label{eq:2dEuclideanpolyform}
\alpha \diamond \alpha = 2\,\alpha^{(0)}\,\alpha\, , \quad \tau(\alpha) = \alpha\, . 
\end{equation}

\noindent
Writing $\alpha = \alpha^{(0)}\oplus \alpha^{(1)} \oplus \alpha^{(2)}$, the second of these relations reads:
\begin{equation*}
\alpha^{(0)} + \alpha^{(1)} - \alpha^{(2)} = \alpha^{(0)} + \alpha^{(1)} + \alpha^{(2)}\, . 
\end{equation*}

\noindent
This gives $\alpha^{(2)} = 0$, whence the first equation in \eqref{eq:2dEuclideanpolyform} becomes $(\alpha^{(0)})^2 = h^{\ast}(\alpha^{(1)},\alpha^{(1)})$ and we conclude that $\alpha$ is a signed square of a spinor if and only if:
\begin{equation*}
\alpha = \pm h^{\ast}(\alpha^{(1)},\alpha^{(1)})^{\frac{1}{2}} \oplus \alpha^{(1)} ~~ \mathrm{with} ~~ \alpha^{(1)}\in V^{\ast}\, . 
\end{equation*}

\noindent
Let $(e^1,e^2)$ be an orthonormal basis of $(V^{\ast},h^{\ast})$ and $\alpha = \cE_\Sigma^+(\xi)$ for some $\xi\in\Sigma$. Then:
\begin{equation*}
2\,\alpha = \cB(\xi,\xi) + \cB(\gamma_i(\xi),\xi)\, e^i\, . 
\end{equation*}

\noindent
Therefore:
\begin{equation*}
4\,h^{\ast}(\alpha^{(1)},\alpha^{(1)}) = \cB(\xi,\xi)^2
\end{equation*}

\noindent
and hence the norm of $\xi$ determines the norm of one-form $\alpha^{(1)}\in V^{\ast}$.


\subsubsection{Signature $(1,1)$}
\label{sec:LorentzExample}


Let $(V^{\ast},h^{\ast})$ be a two-dimensional vector space $V^{\ast}$ equipped with a Lorentzian metric $h^{\ast}$. Its irreducible Clifford module $(\Sigma,\gamma)$ is two-dimensional and equipped with a symmetric admissible bilinear pairing $\cB$ of split signature and positive adjoint type (see Theorem \ref{thm:admissiblepairings}). To guarantee that $\alpha\in \wedge V^\ast$ belongs to $Z_{+,+}(V^\ast, h^\ast)$, we should in principle consider the first equation in \eqref{eq:thmdefequations} of Theorem \ref{thm:reconstruction} for all $\beta \in \wedge V^{\ast}$. However, $V^{\ast}$ is two-dimensional and Example \ref{ep:2dsplitsig} shows that it suffices to take $\beta=1$. Thus $\alpha$ belongs to the set $Z_{+,+}(V^\ast,h^\ast)$ if and only if:
\begin{equation}
\label{eq:2dLorentzpolyform}
\alpha \diamond \alpha = 2 \, \alpha^{(0)}\,\alpha\, , \quad \tau(\alpha) = \alpha\, . 
\end{equation}

\noindent
Writing $\alpha = \alpha^{(0)} + \alpha^{(1)} + \alpha^{(2)}$, the second condition gives $\alpha^{(2)} = 0$, while the first condition becomes:
\begin{equation*}
(\alpha^{(0)})^2 = h^{\ast}(\alpha^{(1)},\alpha^{(1)})\, . 
\end{equation*}

\noindent
In particular, $\alpha^{(1)}$ is space-like or null. Hence $\alpha$ is the signed square of a spinor if and only if:
\begin{equation}
\label{eq:alpha11}
\alpha = \pm h^{\ast}(\alpha^{(1)},\alpha^{(1)})^{\frac{1}{2}} + \alpha^{(1)}
\end{equation}

\noindent
for a one-form $\alpha^{(1)}\in V^{\ast}$.  As in the Euclidean case, we have:
\begin{equation*}
2\,\alpha = \cB(\xi,\xi) + \cB(\gamma_i(\xi),\xi) \, e^i 
\end{equation*}

\noindent
whence:
\begin{equation*}
4\,h^{\ast}(\alpha^{(1)},\alpha^{(1)}) = \cB(\xi,\xi)^2\, . 
\end{equation*}

\noindent
Thus $\alpha^{(1)}$ is null if and only if $\cB(\xi,\xi) = 0$. In this signature the volume form squares to the identity and therefore we have a well-defined notion of chirality. Fix $\mu\in \mathbb{Z}_2$. By Corollary \ref{cor:reconstructionchiral}, $\alpha$ lies in the set $Z^{\mu}_{+,+}(V^\ast,h^\ast)$ if and only if it has the form \eqref{eq:alpha11} and satisfies the supplementary condition:
\begin{equation*}
\ast\,(\pi\circ\tau)(\alpha) = \mu \, \alpha\, . 
\end{equation*}

\noindent
This amounts to the following system, where $\nu_h$ is the volume form of $(V^{\ast},h^{\ast})$:
\begin{equation*}
\pm h^{\ast}(\alpha^{(1)},\alpha^{(1)})^{\frac{1}{2}} \nu_h - \ast \alpha^{(1)} = 
\pm \mu\, h^{\ast}(\alpha^{(1)},\alpha^{(1)})^{\frac{1}{2}} + \mu\, \alpha^{(1)}\, . 
\end{equation*} 

\noindent
Thus $h^{\ast}(\alpha^{(1)},\alpha^{(1)}) = 0$ and $\ast \alpha^{(1)} = -\mu\, \alpha^{(1)}$. Hence a signed polyform square of a chiral spinor of chirality $\mu$ is a null one-form which is anti-self-dual when $\mu = +1$ and self-dual when $\mu = - 1$. Notice that the nullity condition on $\alpha^{(1)}$ is 
equivalent with (anti-)selfduality.


\subsubsection{Signature $(2,2)$}
\label{sec:(2,2)}


Let $(V^{\ast},h^{\ast})$ be four-dimensional and equipped with a metric $h^{\ast}$ of split signature. Its irreducible real Clifford module $(\Sigma,\gamma)$ is four-dimensional and is equipped with a skew-symmetric admissible pairing $\cB$ of positive adjoint type (see Theorem \ref{thm:admissiblepairings}). This dimension and signature admits chiral spinors. Let:
\begin{equation*}
\alpha = \sum_{k=0}^{4} \alpha^{(k)}\in \wedge V^{\ast}~~\mathrm{with}~~\alpha^{(k)}\in \wedge^k V^\ast \quad \forall\,\, k=1,\ldots 4\, . 
\end{equation*}

\noindent
Fixing an orthonormal basis $(e^1,e^2,e^3,e^4)$ of $(V^{\ast},h^{\ast})$ with $e^1, e^2$ timelike, define \emph{timelike} and \emph{spacelike} volume forms through $\nu_{-} = e^1\wedge e^2$ and $\nu_{+} = e^3\wedge e^4$. By Corollary \ref{cor:reconstructionchiral}, we have $\alpha\in Z^{(\mu)}_{-,+}(V^\ast,h^\ast)$ if and only if:
\begin{equation}
\label{eq:eqs2,2}
\alpha\diamond_h \alpha = 0\, , \quad \tau(\alpha) = -\alpha\, , \quad  \ast \pi(\tau(\alpha)) = \mu\, \alpha\, , \quad     \alpha\diamond_h \beta \diamond \alpha = 4\, (\beta\diamond\alpha)^{(0)}\, \alpha
\end{equation}

\noindent
for a polyform $\beta\in \wedge V^{\ast}$ such that $(\beta\diamond\alpha)^{(0)}\neq 0$. Here we used skew-symmetry of $\cB$, which implies $ \alpha^{(0)} = 0$. The condition $\tau(\alpha) = -\alpha$ amounts to:
\begin{equation*}
\alpha^{(0)} = \alpha^{(1)} = \alpha^{(4)} = 0 
\end{equation*}

\noindent
whereas the condition $\ast \pi(\tau(\alpha)) = \mu\, \alpha$ is equivalent with:
\begin{equation*}
\ast\alpha^{(2)} = - \mu\, \alpha^{(2)}\, , \quad  \alpha^{(3)} = 0\, . 
\end{equation*}
Thus it suffices to consider $\alpha = \omega$, where $\omega$ is selfdual if $\mu = -1$ and anti-selfdual if $\mu = 1$. In signature $(2,2)$, the Hodge star operator squares to the identity and yields a decomposition:
\begin{equation*}
\wedge^2 V^{\ast} = \wedge^2_{+} V^{\ast} \oplus \wedge^2_{-} V^{\ast}
\end{equation*}

\noindent
into self-dual and anti-selfdual two-forms. This corresponds to the decomposition $\mathfrak{so}(2,2) = \mathfrak{sl}(2) \oplus \mathfrak{sl}(2)$ of the Lie algebra $\mathfrak{so}(2,2) = \wedge^2 V^{\ast}$.
Expanding the geometric product shows that the first equation in \eqref{eq:eqs2,2} reduces to the following condition for a self-dual or anti-selfdual two-form $\alpha = \omega$:
\begin{equation*}
\langle \omega,\omega\rangle_h = 0\, . 
\end{equation*}

\noindent
For simplicity of exposition we set $\mu = -1$ in what follows, in which case $\omega$ is self-dual (analogous results hold for $\mu =1$). Consider the basis $\left\{u_a\right\}_{a = 1, 2, 3}$ of $\wedge^2_+V^\ast$ given by:
\begin{equation*}
u_1  :=  e^1\wedge e^2 + e^3\wedge e^4\, , \quad   u_2  :=  e^1\wedge e^3 +  e^2\wedge e^4\, , \quad  u_3  :=  e^1\wedge e^4 - e^2\wedge e^3 
\end{equation*}

\noindent
and expand:
\begin{equation*}
\omega = \sum k^a  u_a\, . 
\end{equation*}

\noindent
We have:
\begin{equation*}
\nu_{-}\diamond u_1 = u_1\diamond \nu_{-} = -1 + \nu_h\, , \quad 
\nu_{-}\diamond u_2 = - u_2\diamond\nu_{-} = -u_3\, , \quad 
\nu_{-}\diamond u_3 = -u_3\diamond\nu_{-} = u_2
\end{equation*}

\noindent
which gives:
\begin{equation*}
(\nu_{-}\diamond\omega)^{(0)} = - k^1\, . 
\end{equation*}

\noindent
Furthermore, we compute:
\begin{eqnarray*}
& u_1\diamond u_3 = - u_3 \diamond u_1 = 2\, u_2 \, , \quad  u_1\diamond u_2 = -u_2\diamond u_1 = -2\, u_3\, , \quad  u_2\diamond u_3 = - u_3\diamond u_2 = 2\, u_1\, , \\ & u_1\diamond u_1 = - u_2 \diamond u_2 = - u_3 \diamond u_2 = -2 + 2 \, \nu_h\, . 
\end{eqnarray*}

\noindent
These products realize the Lie algebra $\mathfrak{sl}(2,\R)$ upon defining a Lie bracket by the commutator:
\begin{equation*}
[u_1 , u_2] = u_1 \diamond u_2 - u_2 \diamond u_1= - 4\, u_3\, , \quad 
[u_1 , u_3] = u_1 \diamond u_3 - u_3 \diamond u_1= 4\, u_2\, , \quad 
[u_2 , u_3] = u_2 \diamond u_3 - u_3 \diamond u_2= 4\, u_1\, . 
\end{equation*}

\noindent
Since $\wedge^2_{+} V^{\ast} = \mathfrak{sl}(2,\R)$, the Killing form $\mathfrak{B}$ of $\mathfrak{sl}(2,\R)$ gives a symmetric non-degenerate pairing of signature $(1,2)$ on $\wedge^2_{+}
V^{\ast}$, which can be rescaled to coincide with that induced by $h$. Then:
\begin{equation*}
\mathfrak{B}(\omega,\omega) = \langle \omega,\omega\rangle_h^2 = 2 \left[(k^ 1)^2 - (k^2)^2 - (k^3)^2\right] \quad \forall \omega \in \wedge^2_{+} V^{\ast}\, . 
\end{equation*}

\begin{prop}
\label{prop:(2,2)}
A polyform $\alpha \in \wedge V^{\ast}$ is a signed square of a real chiral spinor $\xi \in \Sigma^{(-)}$ of negative chirality if and only if $\alpha$ is a self-dual two-form of zero norm.
\end{prop}

\begin{proof} 
It suffices to consider the case $\alpha\neq 0$. By the discussion above, a non-zero polyform $\alpha\neq 0$ belongs to the set $Z^{-}_{-,+}(V^\ast,h^\ast)$ only if $\alpha=\omega$ is self-dual and of zero norm (which is equivalent to the first three equations in \eqref{eq:eqs2,2}). Once these conditions are satisfied, the only equation that remains to be solved is the fourth equation in \eqref{eq:eqs2,2}. To solve it, we take $\beta = \nu_{-}$. Since $(\nu_{-}\diamond\omega)^{(0)} = -4\, k^1$ (as remarked above), we conclude that $(\nu_{-}\diamond\omega)^{(0)} \neq 0$ if and only if $\omega \neq 0$, whence taking $\beta = \nu_{-}$ is a valid choice. A computation shows that this equation is automatically satisfied and thus we conclude.
\end{proof}

\begin{remark} 
Subsection \ref{sec:LorentzExample} together with Proposition \ref{prop:(2,2)} show that the square of a chiral spinor in signatures $(1,1)$ and $(2,2)$ is given by an (anti-)self-dual form of zero norm in middle degree. The reader can verify, through a computation similar to the one presented in this subsection, that the same statement holds in signature $(3,3)$. It is tempting to conjecture that the square of a chiral spinor in general split signature $(p,p)$ corresponds to an (anti-)self-dual $p$-form of zero norm, the latter condition being automatically implied when $p$ is odd. Verifying this conjecture would be useful in the study of manifolds of split signature which admit parallel chiral spinors \cite{Dunajski}.
\end{remark}


\section{Algebraic Spin(7) structures}
\label{sec:AlgebraicSpin(7)}


The spinorial exterior forms associated to chiral irreducible real spinors in eight Euclidean dimensions deserve a separate study because of their applications to the theory of Spin(7) structures. In this section we apply Corollary \ref{cor:reconstructionchiral} in eight euclidean dimensions to construct an algebraic function whose critical points describe algebraic Spin(7) structures.


\subsection{$\Spin(7)$ structures on an eight-dimensional vector space}


By definition, a {\em $\Spin(7)$ structure} on $V$ is a $\Spin(7)$ subgroup of the group $\GL_+(V)$. A well-known way to describe such a structure is to give a non-zero four-form on $V$, called the {\em Cayley form} of the structure, with certain properties. We start by recalling this description. Consider first the real vector space $\mathbb{R}^8$ and denote by $(e_1,\hdots,e_8)$ its standard basis, with dual basis $(e^1,\hdots,e^8)$. Let $h_0$ be the usual Euclidean metric on $\mathbb{R}^8$ and endow this space with its canonical orientation, for which the Euclidean volume form reads:
\begin{equation*}
\nu_{o} = e^1\wedge e^2 \wedge e^3 \wedge e^4 \wedge e^5 \wedge e^6 \wedge e^7 \wedge e^8\, .
\end{equation*}
The standard $\Spin(7)_+$ structure on $\R^8$ is described by the {\em canonical Cayley form} $\Phi_0 \in \wedge^4(\mathbb{R}^8)^{\ast}$, which is defined as follows \cite{Bonan,HarveyBook,Joyce2007,SalamonWalpuski}:
\begin{eqnarray}
	& \Phi_0 := e^1\wedge e^2\wedge e^3 \wedge e^4+e^1\wedge e^2\wedge e^5 \wedge e^6 + e^1\wedge e^2\wedge e^7 \wedge e^8 +e^1\wedge e^3\wedge e^5 \wedge e^7 \label{Phi_0}\\
	& - e^1\wedge e^3\wedge e^6 \wedge e^8 - e^1\wedge e^4\wedge e^5 \wedge e^8 - e^1\wedge e^4\wedge e^6 \wedge e^7 + e^5\wedge e^6\wedge e^7 \wedge e^8 + e^3\wedge e^4\wedge e^7 \wedge e^8 \nonumber\\
	& + e^3\wedge e^4\wedge e^5 \wedge e^6 + e^2\wedge e^4\wedge e^6 \wedge e^8 - e^2\wedge e^4\wedge e^5 \wedge e^7 - e^2\wedge e^3\wedge e^5 \wedge e^8 - e^2\wedge e^3\wedge e^6 \wedge e^7\nonumber
\end{eqnarray}
This four-form is self-dual with respect to the Euclidean metric $h_0$ and the canonical orientation of $\R^8$. Furthermore, we have $\Phi_0\wedge \Phi_0 = 14\nu_0$ and hence the square norm of $\Phi_0$ with respect to $h_0$ is $\vert\Phi_0\vert_{h_0}^2 = 14$. The general linear group $\Gl(8,\mathbb{R})$ acts naturally on $\wedge^4(\mathbb{R}^8)^{\ast}$ and consequently on $\Phi_0$. The stabilizer of $\Phi_0$ under this action is isomorphic to the Lie group $\Spin(7)$ \cite{Bonan,HarveyBook, Joyce2007,SalamonWalpuski}. The stabilizer preserves $h_0$ and $\nu_0$ and hence is a subgroup of the special orthogonal group $\SO(8)\subset \GL(8,\mathbb{R})$ determined by $h_0$ and $\nu_0$. 

\begin{definition}
\label{def:Spin7forms} 
A {\em $\Spin(7)_+$ form} on $V$ is a four-form $\Phi \in \wedge^4 V^{\ast}$ for which there exists an orientation-preserving linear isomorphism $f\colon V\to \mathbb{R}^8$ satisfying $\Phi = f^{\ast}\Phi_0$, where $f^{\ast}\colon \wedge(\mathbb{R}^8)^{\ast} \to \wedge V^{\ast}$ denotes the pull-back of forms by $f$. A {\em $\Spin(7)_{-}$ form} is defined similarly but using an orientation-reversing linear isomorphism $f$. A \emph{$\Spin(7)$ form} on $V$ is either a $\Spin(7)_+$  or a $\Spin(7)_-$ form defined on $V$.
\end{definition}

\noindent 
In particular, pulling back the canonical Cayley form by any orientation-reversing linear automorphism of $\R^8$ produces a $\Spin(7)_-$ form. A four-form $\Phi\in \wedge^4 V$ is a $\Spin(7)_+$ form on $V$ if and only if there exists a positively-oriented basis $(\epsilon_1,\ldots, \epsilon_8)$ of $V$ in which $\Phi$ is given by the relation obtained from \eqref{Phi_0} by replacing $e_i$ with $\epsilon_i= f^{-1}(e_i)$ for all $i=1 ,\ldots,8$.  Every $\Spin(7)_+$ form $\Phi$ comes together with an Euclidean metric on $V$ given by:
\begin{equation*}
h_\Phi = f^{\ast} h_0=\sum_{i=1}^8 \epsilon^i\otimes \epsilon^i
\end{equation*}
which makes $\epsilon_1,\ldots, \epsilon_8$ into an orthonormal basis. Similar statements hold for $\Spin(7)_-$ forms, except that the relevant bases of $V$ are negatively oriented. If $\Phi\in \wedge^4 V$ is a $\Spin(7)$ form, then $\lambda \Phi$ is also a $\Spin(7)$ form for any {\em positive} $\lambda\in \R_{>0}$.

A $\Spin(7)$ form $\Phi$ determines its associated metric $h_\Phi$ algebraically as explained in \cite[Section 4.3]{KarigiannisDefs}; we say that $h_\Phi$ is {\em induced by $\Phi$}. Let $\nu_{h_\Phi}=f^\ast \nu_0$ be the corresponding volume form. A $\Spin(7)_+$ form $\Phi$ is self-dual with respect to $h_\Phi$ in the given orientation of $V$ and satisfies $\Phi\wedge \Phi=14 \nu_{h_\Phi}$ (as can be verified for $\Phi_0$ on $\R^8$ and pulling back by $f$), a condition which amounts to $|\Phi|_{h_\Phi}=\sqrt{14}$. On the other hand, a $\Spin(7)_-$ form $\Phi$ is anti-self-dual and satisfies $\Phi\wedge \Phi=-14 \nu_{h_\Phi}$, which again amounts to $|\Phi|_{h_\Phi}=\sqrt{14}$. Notice that a $\Spin(7)$ form can induce {\em any} Euclidean metric on $V$. The stabilizer of a four-form $\Phi\in \wedge^4 V^\ast$ inside $\GL_+(V)$ is isomorphic with the group $\Spin(7)$ if and only if there exists a sign factor $\pm$ such that $\pm \Phi$ is a $\Spin(7)$ form. There are exactly two conjugacy classes of $\Spin(7)$ subgroups in $\Gl_+(V)$ \cite{Varadarajan}. The subgroups belonging to one of these stabilize $\Spin(7)_+$ forms, while those belonging to the other stabilize $\Spin(7)_-$ forms defined on $(V,h)$. Moreover, two $\Spin(7)$ forms have the same stabilizer if and only if they differ by multiplication with a {\em positive} real number. This establishes a bijection between the conjugacy class of $\Spin(7)$ subgroups of $\GL_+(V)$ and the set of positive homothety classes of $\Spin(7)$ forms. The conjugacy classes of $\Spin(7)_+$ and $\Spin(7)_-$ forms inside $\GL(V)$ combine into a single conjugacy class within $\Gl(V)$; for example, the reflection of $V$ in any hyperplane contained in $V$ conjugates $\Spin(7)_+$ forms to $\Spin(7)_-$ forms inside $\GL(V)$. 

\begin{remark}
If $\Phi$ is a $\Spin(7)_{+}$ form on $V$, then $\Phi$ and $-\Phi$ have the same stabilizer, which is a $\Spin(7)_+$ subgroup of $\GL_+(V)$. However, $-\Phi$ is {\em not} a $\Spin(7)_{+}$ form although it can be written as $f^\ast(-\Phi_0)$ for some orientation-preserving isomorphism $f:V\rightarrow \R^8$. A simple continuity argument shows that the overall sign of the canonical Cayley form cannot be changed by acting with an element of $\GL_+(8,\mathbb{R})$ on $\R^8$. It is traditional to fix the overall sign of the canonical Cayley form in order to avoid double counting of $\Spin(7)$ subgroups of $\GL_+(V)$.
\end{remark}


\subsection{$\Spin(7)$ structures on an eight-dimensional Euclidean space}


Let $(V,h)$ be an oriented Euclidean vector space of dimension eight. Let $\O(V,h)\subset \GL(V)$ be the disconnected group of all orthogonal transformations of $(V,h)$ and let $\SO(V,h)\subset \O(V,h)$ be its identity component. These groups act naturally on $\wedge V^\ast$. The Hodge operator $\ast_h$ of $h$ squares to the identity on four-forms and hence gives a decomposition:
\begin{equation*}
\wedge^4 V^{\ast} = \wedge^4_{+} V^{\ast} \oplus \wedge^4_{-} V^{\ast} 
\end{equation*}
where $\wedge^4_{+} V^{\ast}$ and $\wedge^4_{-} V^{\ast}$ are the eigenspaces of self-dual and anti-self-dual four forms. We denote by $\langle \cdot ,\cdot \rangle_h$ the scalar product induced by $h$ on the exterior algebra: 
\begin{equation*}
\wedge V^\ast:=\oplus_{k=0}^8 \wedge^k V^\ast
\end{equation*}
and by $\vert \cdot \vert_h$ the corresponding norm. For later convenience, we will often work with the dual oriented Euclidean space $(V^\ast, h^\ast)$ instead of $(V,h)$. By definition, a {\em metric $\Spin(7)$ structure} on $(V,h)$ is a $\Spin(7)$ subgroup of $\SO(V,h) \subset \GL_+(V)$. There exist two conjugacy classes of such subgroups in $\SO(V,h)$, which correspond to the two conjugacy classes in $\GL_+(V)$. They combine into a single conjugacy class of $\O(V,h)$. 

\begin{definition}	
\label{def:Spin7formsII}
A {\em metric} $\Spin(7)$ form on $(V,h)$ is a $\Spin(7)$ form $\Phi$ on $V$ which satisfies $h_\Phi=h$. A \emph{conformal} $\Spin(7)$ form on $(V,h)$ is a $\Spin(7)$ form $\Phi$ on $V$ which satisfies $h_\Phi= c_{\Phi} h$ for some constant $c_{\Phi} > 0$.
\end{definition}

\noindent 
The positive number $c_{\Phi}$ is uniquely determined by the conformal $\Spin(7)$ structure $\Phi$; we call it the {\em conformal constant} of $\Phi$ relative to $h$. The stabilizer of a $\Spin(7)$ form $\Phi$ on $V$ is a metric $\Spin(7)$ structure on $(V,h)$ if and only if $\Phi$ is a conformal $\Spin(7)$ form on $(V,h)$. The conformal constant $c_{\Phi}$ of a conformal  $\Spin(7)$ form on $(V,h)$ can be expressed through the norm of $\Phi$ as follows. 

\begin{lemma}
\label{lemma:conformalconstant}
Let $\Phi\in \wedge^4 V^{\ast}$ be a conformal $\Spin(7)$ form on $(V,h)$. Then we have: 
\begin{equation*}
c_{\Phi} = 14^{-\frac{1}{4}} \vert\Phi\vert^{\frac{1}{2}}_h
\end{equation*}

\noindent
Thus $\Phi$ is a metric $\Spin(7)$ form on the oriented Euclidean space $(V,h_\Phi)$, where: 
\begin{equation} 	
\label{hPhi}
h_\Phi=14^{-\frac{1}{4}} \vert\Phi\vert^{\frac{1}{2}}_h h
\end{equation} 	

\noindent
is the metric induced by $\Phi$. In particular, a conformal $\Spin(7)$ form on $(V,h)$ is a metric $\Spin(7)$ form if and only if $\vert\Phi\vert_h = \sqrt{14}$.
\end{lemma}

\begin{proof}
Let  $\Phi\in \wedge^4 V^{\ast}$ be a conformal $\Spin(7)$ form on $(V,h)$ and let $f:V \rightarrow \R^8$ be a linear isomorphism such that $\Phi=f^\ast(\Phi_0)$. We have $h_{\Phi} = f^{\ast}(h_0) = c_{\Phi} h$, which implies:
\begin{equation*}
\sqrt{14} = |\Phi|_{h_\Phi}=|\Phi|_{c_\Phi h} = (c_\Phi)^{-2} |\Phi|_{h}  
\end{equation*}
This gives:
\begin{equation*}
c_\Phi = 14^{-\frac{1}{4}} \vert\Phi\vert^{\frac{1}{2}}_h
\end{equation*}
Hence $h_\Phi$ equals \eqref{hPhi}. 
\end{proof}

\begin{cor} 
Let $\Phi\in \wedge^4 V^\ast$ be a four-form on $V$. Then $\Phi$ is a conformal $\Spin(7)$ form on $(V,h)$ if and only if:
\begin{equation*}
\frac{\sqrt{14}}{|\Phi|_h} \Phi\in \wedge^4_{+} V^{\ast}
\end{equation*}
	
\noindent
is a metric $\Spin(7)$ form on $(V,h)$. 
\end{cor}
 
\begin{cor}
\label{cor:conformal}
The set of conformal $\Spin(7)_\pm$ forms on $(V,h)$ is in bijection with the set of positive homothety classes of metric $\Spin(7)_\pm$ structures on $(V,h)$. 
\end{cor}

\noindent
Each $\Spin(7)_\pm$ subgroup of $\SO(V,h)$ corresponds to the positive homothety class of a metric $\Spin(7)_\pm$ form, which consists of conformal  $\Spin(7)_\pm$ forms. In particular, two conformal $\Spin(7)_\pm$ forms on $(V,h)$ have the same stabilizer inside $\SO(V,h)$ if and only if they differ through multiplication by a constant. Two metric $\Spin(7)_\pm$ forms on $(V,h)$ have the same stabilizer inside $\SO(V,h)$ if and only if they coincide. The following endomorphism of $\wedge^4 V^\ast$ will be used later on.

\begin{definition}
\label{def:Lambda}
For any $\Phi\in \wedge^4 V^\ast$, define $\Lambda^h_\Phi\in \End(\wedge^4 V^\ast)$ by:
\begin{equation} 	
\label{Lambda}
\Lambda^h_\Phi(\omega):=2 \Phi \Delta_2^h \omega :=  (\iota_{e^{i_1}} \iota_{e^{i_2}}\Phi)\wedge(\iota_{e^{i_1}} \iota_{e^{i_2}} \omega)\quad\forall\,\, \omega\in \wedge^4 V^\ast 
\end{equation} 	

\noindent
where $(e^1, \hdots e^8)$ is any orthonormal basis of $(V^{\ast} , h^{\ast})$.
\end{definition}

\begin{remark}
\label{rem:Lambda}
In an arbitrary basis $(v^1,\ldots, v^8)$ of $V^\ast$, we have: 
\begin{equation} 	
\label{PhiDeltaomega}
\Phi\Delta_2^h \omega=\frac{1}{8}\Phi_{ijmn}\omega^{mn}_{\,\,\,\,\,\,\,\,\,\, kl} v^i\wedge v^j\wedge v^k\wedge v^l\quad\forall \Phi,\omega\in \wedge^4 V^\ast 
\end{equation} 	

\noindent
where: 
\begin{equation*}
\omega^{m n}_{\,\,\,\,\,\,\,\,\, kl}=h^{mp} h^{nq} \omega_{pqkl}
\end{equation*}

\noindent
and we expanded $\Phi$ and $\omega$ as: 
\begin{equation*}
\Phi=\frac{1}{4!}\Phi_{ijkl} v^i\wedge v^j\wedge v^k\wedge v^l~~,~~\omega=\frac{1}{4!}\omega_{ijkl} v^i\wedge v^j\wedge v^k\wedge v^l 
\end{equation*}

\noindent
Here we use the so-called \emph{determinant convention} for the wedge product of forms. Using \eqref{PhiDeltaomega}, it is easy to check that the operator $\Lambda^h_\Phi$ defined in \eqref{Lambda} coincides with the operator denoted by $\Lambda_\Phi$ in \cite[Definition 2.7]{KarigiannisFlows} and \cite[Definition 2.5]{DwivediGrad}. 
\end{remark}


\subsection{The algebraic characterization of conformal $\Spin(7)$ forms}


In this section, we present an algebraic characterization of conformal $\Spin(7)_\pm$ forms which follows from the description of the signed square of a real irreducible chiral spinor given in Corollary \ref{cor:reconstructionchiral}. Since in eight euclidean dimensions $\nu_h$ squares to the identity in the K\"ahler-Atiyah algebra, we have $\Psi_{\gamma}(\nu_h)^2=\id_\Sigma$ and the vector space $\Sigma$ splits as a $\cB$-orthogonal direct sum $\Sigma = \Sigma^{+}\oplus \Sigma^{-}$, where $\Sigma^\pm$ are the eigenspaces of $\Psi_{\gamma}(\nu_h)$ corresponding to the eigenvalues $\pm 1$. The subspaces define the \emph{chiral} irreducible representations of the even Clifford subalgebra $\Cl^e(V^\ast,h^\ast)$ of chirality $\pm 1$, respectively. The spin group $\Spin(V^\ast,h^\ast)\subset \Cl^e(V^\ast,h^\ast)$ naturally acts on $\Sigma^\pm$ through the representation induced by $\gamma$. The stabilizer of any nonzero chiral spinor is a $\Spin(7)$ of subgroup of $\Spin(V^{\ast},h^{\ast})\simeq \Spin(8)$. There exists two conjugacy orbits of such subgroups, which correspond respectively to the stabilizers of nonzero chiral spinors of positive and negative chirality and are the $\Spin(7)_\pm$ conjugacy orbits of $\Spin(V^\ast,h^\ast)$. The latter project onto the $\Spin(7)_\pm$ conjugacy classes of $\SO(V^\ast,h^\ast)\simeq \SO(V,h)$ through the double covering morphism $\lambda:\Spin(V,h)\rightarrow \SO(V,h)$. Since the stabilizer of a chiral spinor depends only on its homothety class, this gives a bijection between the $\Spin(7)_\pm$ subgroups of $\Spin(V^\ast,h^\ast)$ and the real projective space $\P(\Sigma^\pm)$. Let  $(e^1 , \hdots , e^8 )$ be an orthonormal basis  of $(V^{\ast},h^{\ast})$. It is well-known that a nonzero four-form $\Phi\in \wedge^4 V^\ast$ is a conformal $\Spin(7)_\pm$ form on $(V,h)$ if and only if there exists a nonzero chiral spinor $\xi\in\Sigma^\pm\setminus\{0\}$ such that:
\begin{equation} 	
	\label{PhiSpinor}
	\Phi=\sum_{1\leq i_1 < \cdots < i_4\leq 8} \cB(\gamma_h(e^{i_1})\cdots \gamma_h(e^{i_4})\xi \, , \, \xi)\, e^{i_1}\wedge \cdots \wedge e^{i_4} 
\end{equation} 	
The $\Spin(7)_\pm$ stabilizer of this form in $\SO(V,h)$ coincides with the image of the $\Spin(7)_\pm$ stabilizer of $\xi$ through the double covering morphism of $\Spin(8) \to \SO(8)$. The chiral spinor $\xi$ with this property is determined by $\Phi$ up to sign.  

\begin{lemma}
\label{lemma:squarespinor} 	
Let $\kappa\in \mathbb{Z}_2$ be a sign factor. A non-zero polyform $\alpha\in \wedge V^\ast$ is the $\kappa$-signed square of a non-zero chiral spinor $\xi\in \Sigma^\mu$ of chirality $\mu\in \{-1,1\}$ if and only if it has the form:
\begin{equation} 	
\label{eq:squareformalgebraic}
\alpha =\frac{\kappa}{16}\left[\frac{1}{\sqrt{14}}|\Phi|_h +\Phi + \frac{\mu}{\sqrt{14}}|\Phi|_h \nu_h\right]~,
\end{equation} 	
where $\Phi \in \wedge^4 V^{\ast}$ is a uniquely determined four-form which satisfies the condition $\ast_h \Phi=\mu \Phi$ as well as the algebraic equation: 
\begin{equation} 	
\label{eq:quadsys}
\sqrt{14} \Phi\Delta_2^h \Phi + 12 |\Phi|_h \Phi = 0 
\end{equation} 	
In this case, we have:
\begin{equation} 	
\label{xiNormPhi}
\cB(\xi,\xi)=\frac{|\Phi|_h}{\sqrt{14}}
\end{equation} 	
and hence $\xi$ has unit $\cB$-norm if and only if $|\Phi|_h = \sqrt{14}$. The nonzero chiral spinor $\xi\in \Sigma^\mu$ is determined by the polyform $\alpha$ up sign. 
\end{lemma}

\begin{remark}
Notice that we can replace $\Phi$ by $-\Phi$ in \eqref{eq:squareformalgebraic} provided that we also do so in \eqref{eq:quadsys}. This changes the sign of the middle term in  \eqref{eq:squareformalgebraic}  and of the second term in the left hand side of \eqref{eq:quadsys}. Which of the signs we choose is a matter of convention, provided that the signs of these terms in  \eqref{eq:squareformalgebraic} and \eqref{eq:quadsys} are changed simultaneously. 
\end{remark}

\begin{proof} 
By Corollary \ref{cor:reconstructionchiral}, a polyform $\alpha \in \wedge V^{\ast}$ is the $\kappa$-signed square of a chiral spinor $\xi$ of chirality $\mu$ if and only if $\kappa\alpha^{(0)}>0$ and equations \eqref{eq:reconstructionchiral0II} or \eqref{eq:reconstructionchiral0} are satisfied. For the present case they reduce to:
\begin{equation}
\label{eq:algebraiceqs1}
\alpha\diamond_h \alpha = 16\, \alpha^{(0)}\, \alpha \, , \qquad \tau(\alpha) = \alpha\, , \qquad \nu_h\diamond_h\alpha= \mu\, \alpha
\end{equation}
 
\noindent
The general solution of the second and third equations in \eqref{eq:algebraiceqs1} can be written as:
\begin{equation} 	
\label{alphasol}
\alpha = \frac{\kappa}{16}(c + \Phi + \mu c\nu_h)~,
\end{equation} 	

\noindent
where $\Phi$ is a four-form which satisfies $\ast_h \Phi=\mu \Phi$ and $c\in\mathbb{R}$ is a constant. This constant must be positive since $\alpha^{(0)}=\frac{\kappa c}{16}$. Plugging \eqref{alphasol} into the first equation of \eqref{eq:algebraiceqs1} gives:
\begin{equation} 	
\label{eq:QuadraticPhiLambda}
\Phi\diamond_h \Phi = 12\,c\,\Phi+ 14\, c^2 (1+\mu\nu_h)~,
\end{equation} 	

\noindent
where we used the last equation in \eqref{eq:algebraiceqs1} and the fact that $\nu_h$ squares to $1$ in the \KA algebra. Expanding the geometric product $\diamond_h$ gives:
\begin{equation} 	
\label{PhiSquareExp}
\Phi\diamond_h\Phi = \vert\Phi\vert_h^2- \Phi \Delta^h_2 \Phi +\Phi\wedge \Phi \in \R\oplus \wedge^4 V^\ast\oplus \wedge^8 V^\ast 
\end{equation} 	
Separating degrees, this shows that \eqref{eq:QuadraticPhiLambda} is equivalent with the conditions: 
\begin{equation} 	
\label{conds}
|\Phi|_h^2=14 c^2 \, , \qquad \Phi\Delta_2^h \Phi + 12 c \Phi = 0\, , \qquad \Phi\wedge \Phi = 14\, c^2 \mu \nu_h 
\end{equation} 	
The first and last of these conditions are equivalent since $\ast_h\Phi=\mu \Phi$ , and give:
\begin{equation} 	
c = \frac{\vert\Phi\vert_h}{\sqrt{14}}  
\end{equation} 	
Substituting this into the middle equation of \eqref{conds} gives \eqref{eq:quadsys}. Relation \eqref{alphasol} gives $\cT r(\alpha)=\kappa c$, which implies $\cB(\xi,\xi)=c=\frac{\vert\Phi\vert_h}{\sqrt{14}}$. Substituting the value of $c$ into \eqref{alphasol} gives \eqref{eq:squareformalgebraic}.
\end{proof}

\begin{lemma}
\label{lemma:ClassicalForm}
Let $\xi\in \Sigma^\mu$ be a nonzero chiral spinor of chirality $\mu\in \mathbb{Z}_2$ and let $\Phi$ be the self-dual four-form given in \eqref{PhiSpinor}, where $(e^1 , \hdots , e^8 )$ is any orthonormal basis  of $(V^{\ast},h^{\ast})$. Then the following statements hold: 
\begin{enumerate}
\item We have: 
\begin{equation} 	
\label{normxi}
\cB(\xi,\xi)=\frac{|\Phi|_h}{\sqrt{14}} 
\end{equation} 	

\item The $\kappa$-signed square of $\xi$ is given by: 
\begin{equation} 	
\label{sqxi}
\cE_\gamma^\kappa(\xi)=\frac{\kappa}{16}\left[\frac{1}{\sqrt{14}}\vert\Phi\vert_h +\Phi + \frac{\mu}{\sqrt{14}}|\Phi|_h\nu_h\right] 
\end{equation} 	
\end{enumerate}
\end{lemma}

\begin{proof}
The quantity $\cB(\xi, \gamma^{i_1}\cdots \gamma^{i_k}\xi)$ with $i_1<\ldots< i_4$ vanishes by $\cB$-symmetry of $\gamma^i$ unless $\frac{k(k-1)}{2}$ is even, which requires $k\in \{0,1, 4, 5, 8\}$. On the other had, this quantity vanishes when $k$ is odd since in that case $\gamma^{i_1}\cdots \gamma^{i_k}\xi$ has different chirality from $\xi$ and since the spaces $\Sigma^+$ and $\Sigma^-$ are $\cB$-orthogonal.  Hence we obtain:
\begin{equation*}
\cB(\xi,\xi)^2=\frac{1}{16} (2\cB(\xi,\xi)^2 +|\Phi|_h^2)~~,
\end{equation*}
which gives \eqref{normxi}. On the other hand, the expansion:
\begin{eqnarray}
\label{eq:alphaexp}
\cE_\gamma^\kappa(\xi) = \frac{\kappa}{16} \sum_{k=0}^8 \sum_{1\leq i_1 < \cdots < i_k\leq 8} \cB(\gamma_h(e^{i_k})\cdots \gamma_h(e^{i_1})\xi,\xi)\, e^{i_1}\wedge \cdots \wedge e^{i_k} 
\end{eqnarray}

\noindent
reduces to: 
\begin{equation*}
\cE_\gamma^\kappa(\xi)=\frac{\kappa}{16}\left[\cB(\xi,\xi)+\Phi+\mu\cB(\xi,\xi)\nu_h\right]~,
\end{equation*}
where we used the relation $\gamma_h(\nu_h)\xi=\mu \xi$. Combing this with \eqref{normxi} gives \eqref{sqxi}.
\end{proof}

\begin{thm}
	\label{thm:Spin7algebraic} 
	The following statements are equivalent for a nonzero self-dual four-form $\Phi\in \wedge^4_{+} V^{\ast}$: 
	\begin{enumerate}
		\item $\Phi$ is a conformal $\Spin(7)_+$ form on $(V,h)$.
		\item The polyform:
		\begin{equation} 	
		\label{eq:alphapolyformform}
		\alpha =\frac{1}{16}\left[\frac{1}{\sqrt{14}}\vert\Phi\vert_h +\Phi + \frac{1}{\sqrt{14}}|\Phi|_h\nu_h\right]
		\end{equation} 	
		is the positive square of a nonzero positive chirality spinor. 
		\item $\Phi$ satisfies the following algebraic equation:
		\begin{equation} 	
		\label{eq:algebraicconditionII}
		\sqrt{14} \Phi\Delta_2^h \Phi + 12 |\Phi|_h \Phi = 0 
		\end{equation} 	
	\end{enumerate}
	In particular, there is a one-to-one correspondence between the set of conformal $\Spin(7)_+$ forms on $(V,h)$ and the set of nonzero self-dual solutions of equation \eqref{eq:algebraicconditionII}.
\end{thm}

\begin{proof} 
	Suppose that $(1)$ holds. Then $\Phi$ is given by \eqref{PhiSpinor} for some non-zero positive chirality 
	spinor $\xi$. Lemma \ref{lemma:ClassicalForm} implies that the polyform $\alpha:=\cE_\gamma^+(\xi)$ is given by \eqref{eq:alphapolyformform} with $\kappa=1$ and hence $(b)$ holds. Thus $(1)$ implies $(2)$. Lemma \ref{lemma:squarespinor} shows that $(2)$ implies $(3)$. Finally, let us assume that $(3)$ holds. Then Lemma \ref{lemma:squarespinor} shows the polyform $\alpha$ defined in terms of $\Phi$ by formula \eqref{eq:alphapolyformform} is the positive square of a positive chirality spinor $\xi$, and hence is given by \eqref{eq:alphaexp} with $\kappa=1$, so that $(a)$ holds. This shows that $(c)$ implies $(a)$ and we conclude. 
\end{proof}

\noindent
Theorem \ref{thm:Spin7algebraic} and Lemma \ref{lemma:conformalconstant} imply the following algebraic characterization of metric $\Spin(7)_+$ forms. 

\begin{cor}
	\label{cor:Spin7metricalgebraic} 
	A self-dual four-form $\Phi$ is a metric $\Spin(7)_+$ form on $(V,h)$ if and only if $|\Phi|_h = \sqrt{14}$ and $\Phi$ satisfies the equation:
	\begin{equation} 	
	\label{metricSpin7criterion}
	\Phi\Delta^h_2 \Phi+12 \Phi=0 
	\end{equation} 	 
\end{cor}

\noindent
Verifying that a self-dual four-form is a conformal $\Spin(7)_{+}$ form by traditional methods involves dealing with the cumbersome task of checking if there exists a basis of $V^{\ast}$ in which $\Phi$ is proportional to the expression given in equation \eqref{Phi_0}. Theorem \ref{thm:Spin7algebraic} provides a different criterion which can be used to verify if such a four-form is a conformal $\Spin(7)_{+}$ form: we only need to check if equation \eqref{eq:algebraicconditionII} is satisfied. This gives an intrinsic characterization of conformal $\Spin(7)_{+}$ forms that do not require the use of any privileged basis of $V^{\ast}$ and that we hope it can be useful for applications, see for instance \cite{fol1,fol2} for early applications of this framework. 

\begin{remark}
	By Lemma \ref{lemma:squarespinor}, a conformal $\Spin(7)_+$ form on $(V,h)$ is metric if and only if the corresponding positive chirality spinor $\xi\in \Sigma^+$ (which is determined up to sign) has unit $\cB$-norm. Notice that equation \eqref{metricSpin7criterion} for metric $\Spin(7)_+$ forms is quadratic but inhomogeneous in $\Phi$, unlike equation \eqref{eq:algebraicconditionII} for conformal $\Spin(7)_+$ forms, which is not quadratic in $\Phi$ but is homogeneous of order two under rescaling  $\Phi$ by a positive constant. Also notice that $\Phi=0$ is a special solution of \eqref{eq:algebraicconditionII}, as expected.
\end{remark}


\subsection{An algebraic potential for conformal \texorpdfstring{$\Spin(7)$}{Spin(7)} forms} 
\label{sec:Wh}


To describe the set of conformal $\Spin(7)_+$ forms on $(V,h)$, consider the following cubic function defined on the vector space of self-dual four-forms:
\begin{equation} 	
\label{eq:cubicfunction}
W_h:\wedge^4_{+} V^{\ast} \rightarrow \mathbb{R}\, ,\quad \Phi \mapsto  W_h(\Phi) := -\frac{\sqrt{14}}{3}\langle \Phi\diamond_h \Phi , \Phi \rangle_h + 4 \langle \Phi , \Phi \rangle^{\frac{3}{2}}_h=\frac{\sqrt{14}}{3}\langle \Phi\Delta_2^h\Phi , \Phi \rangle_h + 4 \langle \Phi , \Phi \rangle_h^{\frac{3}{2}} 
\end{equation} 	
The second equality above follows from \eqref{PhiSquareExp}, which implies: 
\begin{equation} 	
\label{PhiDiamondDelta}
\langle \Phi\diamond_h \Phi , \Phi\rangle_h=-\langle \Phi\Delta_2^h\Phi , \Phi \rangle_h
\end{equation} 	
since $\langle ~ ,~ \rangle_h $ is block-diagonal with respect to the rank decomposition of $\wedge V^\ast$. Note that we can equivalently write $W_h$ as: 
\begin{equation} 	
\label{WS}
W_h(\Phi)=-\frac{\sqrt{14}}{48}\cT r(\Phi\diamond_h \Phi\diamond_h\Phi) + \frac{1}{16} \cT r(\Phi\diamond_h\Phi)^{\frac{3}{2}}~,
\end{equation} 	
where we noticed that $\tau(\Phi)=\Phi$ since $\Phi$ is a four-form. Using the rotational invariance of the K\"{a}hler-Atiyah trace and the rotational equivariance of the generalized products, relation \eqref{WS} implies that the potential $W_h:\wedge^4_+ V^\ast\rightarrow \R$ is invariant under the natural action of $\SO(V,h)$ on $\wedge^4_+ V^\ast$. Moreover, notice that $W_h$ is homogeneous of degree three under positive rescalings $\Phi\rightarrow \lambda \Phi$ ($\lambda>0$) of its argument and that we have $W_h(0)=0$. In particular, the restriction of $W_h$ to $\wedge^4_+ V^\ast$ descends to a section of the real line bundle $\cO(3)\rightarrow \P(\wedge^4_+ V^\ast)$.

\begin{remark}
	\label{rem:WBasis}
	Direct computation gives: 
	\begin{equation} 	
	\langle \Phi\Delta^h_2\Phi , \Phi \rangle_h=\frac{1}{8}\Phi_{i j k l}\Phi^{i j}_{\,\,\,\,\, m n} \Phi^{m n k l }~\mathrm{and}~\langle \Phi , \Phi \rangle_h=\frac{1}{24} \Phi_{i j k l}\Phi^{i j k l}~.
	\end{equation} 	
	Thus: 
	\begin{equation} 	
	W_h(\Phi)=\frac{\sqrt{14}}{24}\Phi_{i j k l}\Phi^{i j}_{\,\,\,\,\, m n} \Phi^{m n k l }+4\left(\frac{1}{24}\Phi_{i j k l}\Phi^{i j k l}\right)^{3/2}~,
	\end{equation} 	
	in any basis $(e^1,\ldots, e^8)$ of $(V,h)$. Notice that we use the determinant inner product of forms, with respect to which the volume form of $(V,h)$ has norm one. Also notice that $W_h$ is of class $C^2$ on $\wedge^4_+ V^\ast$, since the 3-rd power of the norm function $||\cdot||: \wedge^4_+ V^\ast\rightarrow \R$ is of class $C^2$ (and its first two differentials vanish at the origin). See \cite[Theorem 3.1]{RN}.
\end{remark}

\begin{prop}
	\label{prop:d2Phi}
	Let $\Phi\in \wedge^4_+ V^\ast$ be a self-dual four-form on $(V,h)$. For any $q\in \wedge^4_+ V^\ast$, we have: 
	\begin{eqnarray}
		& W_h(\Phi+q)=W_h(\Phi)-\sqrt{14}\left[\langle \Phi\diamond_h \Phi, q\rangle+
		\langle q\diamond_h q, \Phi\rangle_h+\frac{1}{3}\langle q\diamond_h q, q\rangle_h\right] \nonumber \\
		& +4\left[\left(|\Phi|_h^2+|q|_h^2+2\langle \Phi,q\rangle_h\right)^{\frac{3}{2}}-|\Phi|_h^3\right] \label{WhExpansion}
	\end{eqnarray}
	
	\noindent
	When $\Phi\neq 0$ and $|q|_h\ll 1$, the Taylor expansion of $W_h(\Phi)$ around $\Phi$ is: 
	\begin{equation} 	
	\label{WhExp}
	W_h(\Phi+q)=W_h(\Phi)+\langle -\sqrt{14}\Phi\diamond_h\Phi+12 |\Phi|_h\Phi,q\rangle_h +
	6|\Phi|_h|q|_h^2-\sqrt{14}\langle q\diamond_h q,\Phi\rangle_h+ 6\frac{\langle q,\Phi\rangle_h^2}{|\Phi|_h}+\O(\epsilon^3)
	\end{equation} 	
	for $\epsilon:=\frac{|q|_h}{|\Phi|_h}\ll 1$
\end{prop}

\noindent 
Consider a nonzero self-dual four-form $\Phi\in \wedge_{+}^4 V^{\ast}\setminus \{0\}$. We denote by:
\begin{equation*}
\dd_{\Phi} W_h \colon \wedge^4_{+} V^{\ast} \to \R\, , \qquad q \mapsto (\dd_{\Phi} W_h) (q)
\end{equation*}
the differential of $W_h$ computed at $\Phi$ and by:
\begin{equation*}
\dd^2_{\Phi}W_h \colon \wedge^4_{+} V^{\ast} \odot \wedge^4_{+} V^{\ast}\to \R
\end{equation*}
the second differential of $W_h$ at $\Phi$. The latter coincides with the Euclidean Hessian of $W_h$ at $\Phi$ computed relative to the scalar product $\langle~,~\rangle_h$. 

\begin{cor}
	\label{cor:dWh}
	For any nonzero self-dual four-form $\Phi\in \wedge^4_+ V^\ast \setminus \{0\}$, we have: 
	\begin{equation} 	
	\label{dWh}
	(\dd_{\Phi} W_h) (q)=\langle\sqrt{14}\Phi\Delta_2^h\Phi + 12 |\Phi|_h \Phi, q\rangle_h \quad \forall q\in \wedge^4_+ V^\ast
	\end{equation} 	
	and:
	\begin{equation} 	
	\label{ddWh}
	(\dd^2_{\Phi}W_h) (q_1 , q_2) =  \sqrt{14}\,\langle  q_1\Delta_2^h q_2+q_2\Delta_2^h q_1 , \Phi \rangle_h + 12 \vert \Phi\vert_h  \langle   q_1, q_2 \rangle_h + \frac{12}{\vert\Phi\vert_h}  \langle q_1 , \Phi\rangle_h \langle \Phi, q_2 \rangle_h  
	\end{equation} 	
	
	\noindent
	for every $q_1, q_2\in \wedge^4_{+} V^{\ast}$.
\end{cor}

\begin{proof}
	Follows immediately from \eqref{WhExp} by polarization.
\end{proof}

\begin{remark}
	Using \cite[Theorem 3.1]{RN}, one easily checks that $(\dd^2 W_h)(0)=0$. This also follows from the fact that $W_h$ is of class $C^2$. The Cauchy–Schwarz inequality gives:
\begin{equation*}
	\frac{1}{\vert\Phi\vert_h} \vert \langle q_1 , \Phi\rangle_h \langle \Phi, q_2 \rangle_h| \leq |\Phi|_h |q_1|_h |q_2|_h 
\end{equation*}
	and hence the $\lim_{\Phi \to 0} \dd^2_{\Phi}W_h$ exists and equals zero. 
\end{remark}

\begin{thm}
\label{thm:PotentialVh}
A self-dual four-form $\Phi\in \wedge^4_{+} V^{\ast}\setminus \{0\}$ is a conformal $\Spin(7)_+$ form on $(V,h)$ if and only if it is a critical point of the function $W_h$. In this case, we have $W_h(\Phi)=0$.  
\end{thm}

\begin{proof}
Since $\langle~,~\rangle_h$ restricts to a scalar product on $\wedge^4_+ V^\ast$, Corollary \eqref{cor:dWh} implies that $\dd_\Phi W_h$ vanishes if and only if $\Phi$ satisfies equation \eqref{eq:algebraicconditionII}, which by Theorem \ref{thm:Spin7algebraic} happens if and only if $\Phi$  is a conformal $\Spin(7)_+$ form on $(V,h)$. In this case, we have:
\begin{equation*}
W_h(\Phi)  = \frac{\sqrt{14}}{3}\langle \Phi\Delta_2^h \Phi , \Phi \rangle_h + 4 |\Phi|_h^3 =-4|\Phi|^3_h +4 |\Phi|_h^3 =0
\end{equation*}

\noindent
and hence we conclude.
\end{proof}

\noindent
Let $\Phi\in \wedge^4_+ V^\ast$ be a conformal $\Spin(7)_+$ form on $(V,h)$. We denote the induced metric of $\Phi$ by $h_{\Phi}$. Recall that $\Phi$ is a \emph{metric} $\Spin(7)_+$ form with respect to $h_{\Phi}$. Consider the orthogonal decomposition of $\wedge^4_{+} V^{\ast}$ into irreducible representations under the linear action of the $\Spin(7)$ stabilizer of $\Phi$ which is obtained by restricting the action of $\SO(V,h)$ on $\wedge^4_{+} V^{\ast}$ (see, for example, \cite{KarigiannisFlows}):
\begin{equation} 	
\label{eq:odeg}
\wedge^4_+ V^{\ast} = \wedge^4_1 V^{\ast} \oplus \wedge^4_7 V^{\ast} \oplus \wedge^4_{27} V^{\ast} 
\end{equation} 	
The subscript in this decomposition denotes the real dimension of the corresponding irreducible representation of $\Spin(7)$. Notice that $\Phi\in \wedge^4_1 V^\ast$.  Since this subspace is one-dimensional, we have: 
\begin{equation*}
\wedge^4_1 V^\ast=\R \Phi 
\end{equation*}
Thus any $q\in \wedge^4_1 V^\ast$ is of the form $q=\lambda(q) \Phi$ with $\lambda(q)\in \R$. Taking norms gives $|\lambda(q)|=\frac{|q|_h}{|\Phi|_h}$ and hence:
\begin{equation} 	
\label{q1}
q=\epsilon(q) \frac{|q|_h}{|\Phi|_h}\Phi \quad \forall q\in \wedge^4_1 V^\ast~,
\end{equation} 	
where $\epsilon(q):=\sign(\lambda(q))$. As shown in \cite[Proposition 2.8]{KarigiannisFlows} (see also \cite[Proposition 2.6]{DwivediGrad}), the invariant subspaces appearing in the right-hand side of \eqref{eq:odeg} can be written as eigenspaces of the operator $\Lambda^h_\Phi$ of  Definition \ref{def:Lambda} (see Remark \ref{rem:Lambda}):
\begin{eqnarray*}
& \wedge^4_1 V^\ast = \left\{q\in \wedge^4_{+} V^{\ast}\,\, \vert\,\, \Phi \Delta^{h_{\Phi}}_2 q = -12 q \right\}\, , \quad 
\wedge^4_7 V^\ast = \left\{ q \in \wedge^4_{+} V^{\ast}\,\, \vert\,\, \Phi \Delta^{h_{\Phi}}_2 q = - 6 q \right\}\\
& \wedge^4_{27} V^\ast = \left\{ q\in \wedge^4_{+} V^{\ast}\,\, \vert\,\, \Phi \Delta^{h_{\Phi}}_2 q= 2 q\right\} 
\end{eqnarray*}
We also have $\wedge^4_{-} V^\ast=\wedge^4_{35} V^\ast=\ker \Lambda_\Phi^h$. Using Lemma \ref{lemma:conformalconstant}, we write the relations above in terms of $h$:
\begin{eqnarray}
\label{InvSubspaces}
& \wedge^4_1 V^\ast = \left\{ q \in \wedge^4_{+} V^{\ast}\,\, \vert\,\, \Phi\Delta^{h}_2 q = -  \frac{12}{\sqrt{14}} \,\vert \Phi\vert_h q \right\}\, , \quad \wedge^4_7 V^\ast = \left\{ q \in \wedge^4_{+} V^{\ast}\,\, \vert\,\, \Phi\Delta^{h}_2 q = -  \frac{6}{\sqrt{14}} \,\vert \Phi\vert_h  q \right\}\nn\\
& \wedge^4_{27} V^\ast = \left\{q\in \wedge^4_{+} V^{\ast}\,\, \vert\,\, \Phi \Delta^{h}_2 q=  \frac{2}{\sqrt{14}} \,\vert \Phi\vert_h q \right\} 
\end{eqnarray}

\begin{prop}
\label{prop:ddWh}
Let $\Phi\in \wedge^4_+ V^\ast$ be a conformal $\Spin(7)_+$ form on $(V,h)$. Then $\dd^2_\Phi W_h$ is block-diagonal with respect to the decomposition \eqref{eq:odeg}. Moreover, the restrictions of $\dd^2_\Phi W_h$ to the subspaces appearing in this decomposition are given by: 
\begin{equation} 	
\label{HessRestrictions}
(\dd^2_{\Phi}W_h)(q_1,q_2) =
\begin{cases}
	0 & \text{if } q_1,q_2\in \wedge^4_1 V^\ast \\
	0 & \text{if } q_1,q_2\in \wedge^4_7 V^\ast \\
	16 |\Phi|_h \langle q_1,q_2\rangle_h & \text{if } q_1 , q_2\in \wedge^4_{27} V^\ast
\end{cases}
\end{equation} 	
\end{prop}


\subsection{Metric deformations of the potential}
\label{sec:algebraicmetricdef}


In this section, we study the expression in the right hand side of \eqref{eq:cubicfunction} as a function of pairs $(h,\Phi)$, where $h$ is an Euclidean metric on $V$ and $\Phi\in \wedge^4 V^{\ast}$ is a four-form that need not be self-dual with respect to $h$. We thus consider the function:
\begin{equation} 	
\label{eq:metriccubicfunction}
W\colon \Met(V)\times\wedge^4 V^{\ast} \rightarrow \mathbb{R}\, ,\quad (h,\Phi) \mapsto  W(h,\Phi) :=\frac{\sqrt{14}}{3}\langle \Phi\Delta_2^h\Phi , \Phi \rangle_h + 4 \langle \Phi , \Phi \rangle_h^{\frac{3}{2}}~, 
\end{equation} 	

\noindent
where $\Met(V)$ denotes the cone of Euclidean metrics on $V$. 

\begin{definition}
A pair $(h,\Phi)\in \Met(V)\times\wedge^4 V^{\ast}$ is called \emph{self-dual} if $\Phi$ is self-dual relative to $h$. A pair $(h,\Phi)\in \Met(V)\times\wedge^4 V^{\ast}$ is called a \emph{conformal $\Spin(7)_+$ pair} on the oriented eight-dimensional space $V$ if $\Phi$ is  conformal $\Spin(7)_+$ form on the oriented Euclidean vector space $(V,h)$.
\end{definition}

\noindent
We denote by:
\begin{equation*}
	\dd_{(h,\Phi)} W\colon V^{\ast}\odot V^{\ast} \oplus \wedge^4 V^{\ast} \to \mathbb{R}
\end{equation*}

\noindent
the differential of $W$ at the point $(h,\Phi)\in \Met(V)\times \wedge^4 V^{\ast}$. 

\begin{lemma}
\label{lemma:metricdeformations}
Let $(h,\Phi)\in \Met(V)\times \wedge^4 V^{\ast}$ be a conformal Spin(7) structure on $V$. Then:
\begin{equation*}
(\dd_{(h,\Phi)} W)(k,0) = 0
\end{equation*}
	
\noindent
for every $k\in V^{\ast}\odot V^{\ast}$.
\end{lemma}

\begin{proof}
First, recall that by Equation \eqref{PhiDeltaomega} we have:
\begin{eqnarray*}
\Phi\Delta_2^h \Phi=\frac{1}{8}\Phi_{ijmn}\Phi^{mn}_{\,\,\,\,\,\,\,\,\,\, kl} e^i\wedge e^j\wedge e^k\wedge e^l
\end{eqnarray*}
	
\noindent
for every $\Phi\in \wedge^4 V^\ast$ and any basis $(e^1,\hdots ,e^8)$. Hence:
\begin{equation*}
\langle \Phi\Delta_2^h\Phi , \Phi \rangle_h = \frac{1}{8}\Phi_{ijmn}\Phi^{mn}_{\,\,\,\,\,\,\,\,\,\, kl} \Phi^{ijkl} = \frac{1}{8} \Phi_{i_1 i_2 i_3 i_4}\Phi_{j_1 j_2 j_3 j_4}\Phi_{k_1 k_2 k_3 k_4} h^{i_1 k_1}h^{i_2 k_2}h^{i_3 j_1}h^{i_4 j_2}h^{j_3 k_3}h^{j_4 k_4}
\end{equation*}
	
\noindent
as well as:
\begin{eqnarray*}
\vert\Phi\vert^2_h = \frac{1}{4!} \Phi_{ijkl} \Phi^{ijkl} = \frac{1}{4!} \Phi_{i_1 i_2 i_3 i_4} \Phi_{j_1 j_2 j_3 j_4} h^{i_1 j_1}h^{i_2 j_2}h^{i_3 j_3}h^{i_4 j_4}
\end{eqnarray*}
	
\noindent
Therefore:
\begin{eqnarray*}
& W(h,\Phi) = \frac{\sqrt{14}}{24} \Phi_{i_1 i_2 i_3 i_4}\Phi_{j_1 j_2 j_3 j_4}\Phi_{k_1 k_2 k_3 k_4} h^{i_1 k_1}h^{i_2 k_2}h^{i_3 j_1}h^{i_4 j_2}h^{j_3 k_3}h^{j_4 k_4} \\ 
& + \frac{4}{(4!)^{\frac{3}{2}}} (\Phi_{i_1 i_2 i_3 i_4} \Phi_{j_1 j_2 j_3 j_4} h^{i_1 j_1}h^{i_2 j_2}h^{i_3 j_3}h^{i_4 j_4})^{\frac{3}{2}} 
\end{eqnarray*}
	
\noindent
Using this expression for $W(h,\Phi)$, which we have written explicitly in terms of all the metric contractions for clarity in the exposition, a direct computation gives:
\begin{eqnarray*}
& (\dd_{(h,\Phi)} W)(k,0) = - \frac{\sqrt{14}}{4} \Phi_{i_1 i_2 i_3 i_4}\Phi_{j_1 j_2 j_3 j_4}\Phi_{k_1 k_2 k_3 k_4} k^{i_1 k_1}h^{i_2 k_2}h^{i_3 j_1}h^{i_4 j_2}h^{j_3 k_3}h^{j_4 k_4} \\
& - \vert\Phi\vert_h  \Phi_{i_1 i_2 i_3 i_4} \Phi_{j_1 j_2 j_3 j_4} k^{i_1 j_1}h^{i_2 j_2}h^{i_3 j_3}h^{i_4 j_4} \\
& = - 4^{-1} 14^{-\frac{3}{4}} \vert\Phi\vert_h^{\frac{5}{2}}\Phi_{i_1 i_2 i_3 i_4}\Phi_{j_1 j_2 j_3 j_4}\Phi_{k_1 k_2 k_3 k_4} k^{i_1 k_1}h^{i_2 k_2}_{\Phi} h^{i_3 j_1}_{\Phi} h^{i_4 j_2}_{\Phi} h^{j_3 k_3}_{\Phi} h^{j_4 k_4}_{\Phi}\\
& - 14^{-\frac{3}{4}} \vert\Phi\vert^{\frac{5}{2}}_h  \Phi_{i_1 i_2 i_3 i_4} \Phi_{j_1 j_2 j_3 j_4} k^{i_1 j_1}h^{i_2 j_2}_{\Phi} h^{i_3 j_3}_{\Phi} h^{i_4 j_4}_{\Phi}
\end{eqnarray*}
	
\noindent
where we have used Equation \eqref{hPhi} to relate $h$ and the metric $h_{\Phi}$ induced by $\Phi$. Since $h_{\Phi}$ is induced by $\Phi$ and the latter is a Spin(7) structure, the following well-known identities hold \cite{KarigiannisDefs,KarigiannisFlows}:
\begin{eqnarray*}
& \Phi_{i_1 i_2 i_3 i_4} \Phi_{j_1 j_2 j_3 j_4} h^{i_3 j_3}_{\Phi} h^{i_4 j_4}_{\Phi} = 6 (h_{\Phi})_{i_1 j_1} (h_{\Phi})_{i_2 j_2} - 6 (h_{\Phi})_{i_1 j_2} (h_{\Phi})_{i_2 j_1} - 4 \Phi_{i_1 i_2 j_1 j_2}\\
& \Phi_{i_1 i_2 i_3 i_4} \Phi_{j_1 j_2 j_3 j_4} h^{i_2 j_2}_{\Phi} h^{i_3 j_3}_{\Phi} h^{i_4 j_4}_{\Phi} = 42 (h_{\Phi})_{i_1 j_1}
\end{eqnarray*}
	
\noindent
Using these identities, we obtain:
\begin{eqnarray*}
& (\dd_{(h,\Phi)} W)(k,0)  =  14^{-\frac{3}{4}} 42 \vert\Phi\vert_h^{\frac{5}{2}}  k^{i_1 k_1} (h_{\Phi})_{i_1 j_1} - 14^{-\frac{3}{4}} 42 \vert\Phi\vert^{\frac{5}{2}}_h  k^{i_1 k_1} (h_{\Phi})_{i_1 j_1} = 0
\end{eqnarray*}
	
\noindent
and thus we conclude.
\end{proof}

\noindent
Using this lemma we adapt Theorem \ref{thm:PotentialVh} to the case of the function $W\colon \Met(V)\times\wedge^4 V^{\ast} \rightarrow \mathbb{R}$.

\begin{thm}
	\label{thm:MetricPotentialVh}
	A pair $(h,\Phi)\in \Met(V)\times \wedge^4 V^{\ast}$ is a conformal Spin(7) structure on $V$ if and only if it is a self-dual critical point of the function $W\colon \Met(V)\times\wedge^4 V^{\ast} \rightarrow \mathbb{R}$. 
\end{thm}

\begin{proof}
	Let $(h,\Phi)\in \Met(V)\times \wedge^4 V^{\ast}$ be a self-dual pair. Then, by Corollary \ref{cor:dWh} we have: 
\begin{eqnarray*}
	(\dd_{(h,\Phi)} W) (0,q)=\langle\sqrt{14}\Phi\Delta_2^h\Phi + 12 |\Phi|_h \Phi, q\rangle_h \quad \forall q\in \wedge^4_+ V^\ast
\end{eqnarray*}
	
	\noindent
	Hence, condition $(\dd_{(h,\Phi)} W) (0,q) = 0$ holds for all $q\in \wedge^4_+ V^\ast$ if and only if $\langle\sqrt{14}\Phi\Delta_2^h\Phi + 12 |\Phi|_h \Phi, q\rangle_h = 0$, which by Theorem \ref{thm:Spin7algebraic} is equivalent to $\Phi$ being a conformal Spin(7) structure on $(V,h)$. Since $(g,\Phi)$ is a conformal Spin(7) structure, Lemma \ref{lemma:metricdeformations} implies that $(\dd_{(h,\Phi)} W) (k,0) = 0$ for every $k\in V^{\ast}\odot V^{\ast}$ and thus we conclude.
\end{proof}

\noindent
Hence by considering the variations of $W\colon \Met(V)\times \wedge^4 V^{\ast} \to \mathbb{R}$ with respect to pairs $(h,\Phi)$ we can describe all 
$\Spin(7)_+$ structures on $V$, as opposed to those that are conformal relative to a fixed Euclidean metric. It would be interesting to investigate the geometric significance of those critical points of $W$ which are not self-dual.


\section{Bundles of Clifford modules and differential spinors}
\label{sec:differentialspinors}


To study differential spinors of real type, we will extend the theory developed in Section \ref{sec:SpinorsAsPolyforms} to bundles of real irreducible Clifford modules equipped with an arbitrary connection. Throughout this section, let $(M,g)$ denote a  connected pseudo-Riemannian manifold of signature $(p,q)$, namely $p$ \emph{pluses} and $q$ \emph{minuses}, and even dimension $d=p+q\geq 2$, where $p-q\equiv_8 0,2$. Since $M$ is connected, the pseudo-Euclidean vector bundle $(TM,g)$ is modeled on a fixed quadratic vector space denoted by $(V,h)$. For any point $m\in M$, we thus have an isomorphism of quadratic spaces $(T_mM,g_m)\simeq (V,h)$. Accordingly, the cotangent bundle $T^\ast M$ (endowed with the dual metric $g^\ast$) is modeled on the dual quadratic space $(V^{\ast},h^{\ast})$. We denote by $\Cl(M,g)$ the bundle of real Clifford algebras of the {\em cotangent} bundle $(T^\ast M,g^\ast)$, which is modeled on the real Clifford algebra $\Cl(V^{\ast},h^{\ast})$. Let $\pi$ and $\tau$ be the canonical automorphism and anti-automorphism of the Clifford bundle, given by fiberwise extension of the corresponding objects defined in Section \ref{sec:SpinorsAsPolyforms} and set ${\hat \pi}=\pi\circ \tau$.  We denote by $(\wedge M,\diamond)$ the exterior bundle $\wedge M=\oplus_{j=0}^d \wedge^j T^\ast M$, equipped with the point-wise extension $\diamond$ of the geometric product of Section \ref{sec:SpinorsAsPolyforms}, which depends on the metric $g$. This bundle of unital associative algebras is called the \emph{K\"ahler-Atiyah bundle} of $(M,g)$ (see \cite{LazaroiuB,LazaroiuBC}). The map $\Psi$ of Section \ref{sec:SpinorsAsPolyforms} extends to a unital isomorphism of bundles of algebras:
\begin{equation*}
\Psi\colon (\wedge M,\diamond_g) \xrightarrow{\sim} \Cl(M,g)
\end{equation*}

\noindent
which allows us to view the K\"ahler-Atiyah bundle as a model for the Clifford bundle. We again denote by $\pi$, $\tau$ and ${\hat \pi}=\pi\circ \tau$ the (anti-)automorphisms of the \KA bundle obtained by transporting the corresponding objects from the Clifford bundle through $\Psi$. The \KA trace introduced in Section
\ref{sec:SpinorsAsPolyforms} pointwise extends to a morphism of vector bundles:
\begin{equation*}
\cT r:\wedge M\rightarrow M\times \mathbb{R} 
\end{equation*}

\noindent
whose induced map on smooth sections satisfies:
\begin{equation*}
\cT r(1_M) = 2^{\frac{d}{2}} 1_M ~~ \mathrm{and} ~~ \cT r(\omega_1\diamond_g\omega_2)=\cT r(\omega_2\diamond_g
\omega_1)\quad \forall\,\, \omega_1,\omega_2\in \Omega^\ast(M) 
\end{equation*}

\noindent
where $1_M\in C^{\infty}(M)$ is the unit function defined on $M$. By Proposition \ref{prop:cS}, we have:
\begin{equation*}
\cT r(\omega)=2^{\frac{d}{2}}\omega^{(0)}\quad \forall \,\,\omega \in \Omega(M)\, . 
\end{equation*}

\noindent
In particular, $\cT r$ does not depend on the metric $g$.  The following result encodes a well-known property of the Clifford bundle, which also follows from the definition of $\diamond_g$.

\begin{prop}
\label{prop:LCderClifford}
The canonical extension to $\wedge M$ of the Levi-Civita connection $\nabla^g$ of $(M,g)$, which we again denote by $\nabla^g$, acts by derivations of the geometric product, namely:
\begin{equation*}
\nabla^g(\alpha\diamond_g \beta) = (\nabla^g\alpha)\diamond_g \beta +  \alpha\diamond_g (\nabla^g\beta) 
\end{equation*}

\noindent
for every $\alpha, \beta \in \Omega(M)$.
\end{prop}


\subsection{Bundles of irreducible real Clifford modules}


In this section we introduce the type of \emph{bundle of spinors} that we will consider throughout this dissertation.

\begin{definition}
A {\em bundle of real Clifford modules}, or a \emph{real spinor bundle} for short, on $(M,g)$ is a pair $(S,\Gamma)$, where $S$ is a real vector bundle on $M$ and $\Gamma:\Cl(M,g)\rightarrow \End(S)$ is a unital morphism of bundles of algebras.
\end{definition}

\noindent 
Since $M$ is connected, any bundle of Clifford modules $(S,\Gamma)$ on $(M,g)$ is modelled on a Clifford module $(\Sigma,\gamma)$ called its \emph{type}. That is, for every point $m\in M$, the Clifford module $\Gamma_m:\Cl(T_m^\ast M, g_m^\ast) \to (\End(S_m),\circ)$ is isomorphic to the Clifford module $\gamma\colon \Cl(V^{\ast},h^{\ast})\to \End(\Sigma)$ via an \emph{unbased} isomorphism of Clifford modules, see \cite{Lazaroiu:2016vov} for more details.

\begin{definition}
A \emph{bundle of irreducible real Clifford modules} $(S,\Gamma)$, or an {\em irreducible real spinor bundle} for short, is a bundle of real Clifford modules whose type $(\Sigma,\gamma)$ is irreducible. In this case, global sections $\epsilon\in \Gamma(S)$ of $S$ are {\em irreducible spinors} on $(M,g)$.
\end{definition}

\noindent 
In the signatures $p-q\equiv_8 0,2$ considered in this dissertation, the rank of a bundle of irreducible real Clifford modules is $\rk\, S=\dim \Sigma = 2^{\frac{d}{2}}$, where $d$ is the dimension of $M$. In this situation, Reference \cite{Lazaroiu:2016vov} proves that $(M,g)$ admits a bundle of irreducible spinors if and only if it admits a {\em real Lipschitz structure} of type $\gamma$. In signatures $p - q\equiv_8 0,2$, the latter corresponds to an adjoint-equivariant, also known as \emph{untwisted}, $\Pin(V^{\ast},h^{\ast})$-structure $\cQ$ on $(M,g)$ and furthermore $(S,\Gamma)$ is isomorphic to the bundle of real Clifford modules associated to $\cQ$ through the natural representation of $\Pin(V^{\ast},h^{\ast})$ induced by $\gamma$ on $\Sigma$. The obstructions to existence of such structures were given in \cite{Lazaroiu:2016vov}; when $p-q\equiv_8 0,2$, they are a slight modification of those given in \cite{Karoubi} for ordinary twisted adjoint-equivariant $\Pin(V^\ast,h^\ast)$-structures.

\begin{prop}
\label{prop:Ltensor}
Let $(S,\Gamma)$ be a bundle of real Clifford modules on $(M,g)$, $L$ a real line bundle on $M$ and set $S_L :=  S\otimes L$. Then there exists a natural unital morphism of bundles of algebras $\Gamma_L:\Cl(M,g)\rightarrow \End(S\otimes L)$ such that $(S_L,\Gamma_L)$ is a bundle of real Clifford modules. In particular, the set of isomorphism classes of bundles of irreducible real Clifford modules defined over $(M,g)$ is a torsor over the real Picard group $\Pic(M)$ of $M$.
\end{prop}

\begin{proof}
Let $\psi_L: \End(L) \to M\times \R$ be the canonical trivialization of the real line bundle $\End(L)$. This induces a unital isomorphism of bundles of algebras $\varphi_L\colon \End(S\otimes L) \xrightarrow{\sim} \End(S)$ given by composing the natural isomorphism of bundles of real algebras $\End(S\otimes L) \xrightarrow{\sim} \End(S)\otimes \End(L)$ with $\Id_{\End(S)}\otimes \psi_L$. The conclusion follows by setting $\Gamma_L :=  \varphi_L^{-1}\circ \Gamma$.
\end{proof}

\noindent 
The map $\Psi_\gamma$ of Section \ref{sec:SpinorsAsPolyforms} extends to a unital isomorphism of bundles of algebras which we denote by:
\begin{equation*}
\Psi_{\Gamma} :=  \Gamma\circ \Psi \colon (\wedge M,\diamond) \xrightarrow{\sim} (\End(S),\circ)\, .
\end{equation*}

\noindent
This map allows us to identify bundles $(S,\Gamma)$ of modules over $\Cl(T^\ast M, g^\ast)$ with bundles of modules $(S,\Psi_\Gamma)$ over the \KA algebra. For ease of notation we again denote by a \emph{dot} the Clifford multiplication of $(S,\Psi_\Gamma)$, whose action on global sections is:
\begin{equation*}
\alpha\cdot\epsilon  :=  \Psi_{\Gamma}(\alpha)(\epsilon) \quad \forall\,\, \alpha \in \Omega(M) :=  \Gamma(\wedge M) \quad \forall \,\, \epsilon\in \Gamma(S)\, . 
\end{equation*}

\begin{definition}
\label{def:symbol}
Let $(S,\Gamma)$ be a real spinor bundle on $(M,g)$ and $W$ be any vector bundle on $M$. The {\em symbol} of a section $\cQ\in \Gamma(\End(S)\otimes W)$ is the section $\mathfrak{q}\in \Gamma(\wedge M \otimes W)$ defined through:
\begin{equation*}
\mathfrak{q} :=  (\Psi_\Gamma\otimes \Id_W)^{-1}(\cQ)\in \Gamma(\wedge T^\ast M \otimes W)
\end{equation*}

\noindent
where $\Id_W$ is the identity endomorphism of $W$. 
\end{definition}

\begin{remark}
In particular, the symbol of an endomorphism $\cQ\in \Gamma(\End(S))$ is a polyform $\mathfrak{q}\in \Omega(M)$, while the symbol of an $\End(S)$-valued one-form $\cA\in \Gamma(T^\ast M\otimes End(S))$ is an element $\fra \in \Gamma(M,T^\ast M\otimes \wedge T^\ast M)=\Omega^1(M,\wedge M)=\Omega^\ast(M,T^\ast M)$, which can be viewed as a $T^\ast M$-valued polyform or as a $\wedge M$-valued one-form.
\end{remark}


\subsection{Paired spinor bundles}


\begin{definition}
Let $(S,\Gamma)$ be a real spinor bundle on $(M,g)$. A fiberwise-bilinear pairing $\cB$ on $S$ is called {\bf admissible} if $\cB_m:S_m\times S_m\rightarrow \R$ is an admissible pairing on the simple Clifford module $(S_m,\Gamma_m)$ for all $m\in M$.  A real {\em paired spinor bundle} on $(M,g)$ is a triplet $(S,\Gamma,\cB)$, where $(S,\Gamma)$ is a real spinor bundle on $(M,g)$ and $\cB$ is an admissible pairing on $S$.
\end{definition}

\noindent 
Since $M$ is connected, the symmetry and adjoint type $s, \sigma\in \mathbb{Z}_2$ of the admissible pairings $\cB_m$, which are non-degenerate by definition, are constant on $M$; they are called the {\em symmetry type} and {\em adjoint type} of $\cB$ and of $(S,\Gamma,\cB)$ respectively. Since $M$ is paracompact, the defining algebraic properties of an admissible
pairing can be formulated equivalently as follows using global sections when viewing $(S,\Gamma)$ as a bundle $(S,\Psi_\Gamma)$ of modules over the \KA algebra of $(M,g)$:
\begin{enumerate}
\itemsep 0.0em
\item $\cB(\xi_1,\xi_2)= s \cB(\xi_2,\xi_2)$
\item
$\cB(\Psi_\Gamma(\omega)(\xi_1),\xi_2)=\cB(\xi_1,\Psi_\Gamma((\pi^{\frac{1-\sigma}{2}}\circ\tau)(\omega))(\xi_2))\quad \forall\,\, \omega\in \Omega(M) \quad \forall\,\, \xi_1,\xi_2\in \Gamma(S)$
\end{enumerate}

\noindent
for every $\xi_1,\xi_2\in \Gamma(S)$ and every $\omega\in \Omega(M)$.

\begin{definition}
Let $(M,g)$ be a pseudo-Riemannian manifold. We say that $(M,g)$ is \emph{strongly orientable} if its orthonormal coframe bundle admits a reduction to an $\SO_o(V^{\ast},h^{\ast})$-bundle. We say that $(M,g)$ is {\em strongly spin} if it admits a $\Spin_o(V^{\ast},h^{\ast})$-structure  --- which we call a {\em strong spin structure}.  
\end{definition}  

\begin{remark}
\label{rem:stronglyspin}
When $pq=0$, the special orthogonal and spin groups are connected while the pin group has two connected components. In this case, orientability and strong orientability are equivalent, as are the properties of being spin and strongly spin. When $pq\neq 0$, the groups $\SO(V^\ast,h^\ast)$ and $\Spin(V^\ast,h^\ast)$ have two connected components, while $\Pin(V^\ast,h^\ast)$ has four and we have $\Pin(V^\ast,h^\ast)/\Spin_o(V^\ast,h^\ast)\simeq \Z_2\times \Z_2$. In this case, $(M,g)$ is strongly orientable if and only if it is orientable and in addition the principal $\Z_2$-bundle associated to its bundle of oriented coframes through the group morphism $\SO(V^\ast,h^\ast)\rightarrow \SO(V^\ast,h^\ast)/\SO_o(V^\ast,h^\ast)$ is trivial, while an untwisted $\Pin(V^\ast,h^\ast)$-structure $\cQ$ reduces to a $\Spin_o(V^\ast, h^\ast)$-structure if and only if the principal $\Z_2\times \Z_2$-bundle associated to $\cQ$ through the group morphism $\Pin(V^\ast,h^\ast)\rightarrow
\Pin(V^\ast,h^\ast)/\Spin_o(V^\ast,h^\ast)$ is trivial. When $(M,g)$ is strongly spin, the short exact sequence:
\begin{equation*}
1 \to \Z_2 \hookrightarrow \Spin_o(V^{\ast},h^{\ast}) \rightarrow \SO_o(V^{\ast},h^{\ast})\to 1
\end{equation*}

\noindent
induces a sequence in \u{C}ech cohomology which implies that $\Spin_o(V^\ast,h^\ast)$-structures form a torsor over $H^1(M,\Z_2)$. A particularly simple case arises when $H^1(M,\Z_2)=0$ (for example, when $M$ is simply-connected). In this situation, $M$ is strongly orientable and any untwisted $\Pin(V^\ast,h^\ast)$-structure on $(M,g)$ reduces to a $\Spin_o(V^\ast,h^\ast)$-structure since $H^1(M,\Z_2\times \Z_2)=H^1(M,\Z_2\oplus \Z_2)=0$. Similarly, any $\Spin(V^\ast,h^\ast)$-structure on $(M,g)$ reduces to a $\Spin_o(V^\ast,h^\ast)$-structure. Up to isomorphism, in this special case there exists at most one $\Spin(V^\ast,h^\ast)$-structure, one $\Spin_o(V^\ast,h^\ast)$-structure and one real spinor bundle on $(M,g)$, which is automatically strong.
\end{remark}

\noindent 
By the results of \cite{Lazaroiu:2016vov}, a strongly orientable pseudo-Riemannian manifold $(M,g)$ of signature $(p,q)$ satisfying $(p-q) \equiv_8 0,2$ admits an irreducible real spinor bundle if and only if it is strongly spin, in which case every irreducible spinor bundle is associated to a strong spin structure via the tautological representation of $\Spin_o(V^\ast,h^\ast)$ induced by an irreducible real Clifford module.

\begin{prop}
\label{prop:SpinorialConnection}
Suppose that $(S,\Gamma)$ be an irreducible real spinor bundle of type $(\Sigma,\gamma)$ on a strongly orientable pseudo-Riemannian manifold $(M,g)$ of signature $(p-q)\equiv_8 0,2$. Then every admissible pairing on $(\Sigma,\gamma)$ extends to an admissible pairing $\cB$ on $(S,\Gamma)$. Moreover, the Levi-Civita connection $\nabla^g$ of $(M,g)$ lifts to a unique connection on $S$, denoted for simplicity by the same symbol, which acts on $\Gamma(S)$ via module derivations:
\begin{equation*}
\nabla^g_v(\alpha\cdot \epsilon) = (\nabla^g_v\alpha)\cdot \epsilon + \alpha\cdot (\nabla^g_v\epsilon) \quad \forall\,\, \alpha\in \Omega(M)\quad \forall \,\,\epsilon \in \Gamma(S)\quad \forall\,\, v\in \fX(M)
\end{equation*}

\noindent
and is compatible with $\cB$:
\begin{equation*}
v(\cB(\epsilon_1,\epsilon_2)) = \cB(\nabla^g_v\epsilon_1,\epsilon_2) + \cB(\epsilon_1,\nabla^g_v\epsilon_2) \quad \forall\,\, \epsilon_1, \epsilon_2 \in \Gamma(S) \quad \,\,\forall \,\, v\in \fX(M)\, . 
\end{equation*}
\end{prop}

\begin{proof}
The first statement follows from the associated bundle construction since admissible pairings are invariant under $\Spin_o(V^\ast,h^\ast)$ transformations. The second and third statements are standard, see for instance \cite[Chapter 3]{Friedrich}. The last statement holds since the holonomy of $\nabla^g$ as a connection on $S$ is contained in $\Spin_o(V^{\ast},h^{\ast})$, whose action on $\Sigma$ preserves $\cB$.
\end{proof}

\noindent 
The spinorial connection $\nabla^g$ induces a linear connection, which we denote by the same symbol for ease of notation, on the bundle of endomorphisms $\End(S)=S^\ast \otimes S$.  Given $v\in \fX(M)$, by definition we have:
\begin{equation*}
(\nabla^g_v A)(\epsilon) = \nabla^g_v (A(\epsilon)) - A(\nabla^g_v\epsilon) \quad \forall\,\, A\in \Gamma(\End(S)) \quad \forall\,\, \epsilon \in \Gamma(S) \quad \forall\,\, v\in \fX(M) 
\end{equation*}

\noindent
for every $A\in \Gamma(\End(S))$ and $\epsilon \in \Gamma(S)$.

\begin{prop}
Let $(S,\Gamma,\cB)$ be a paired irreducible real spinor bundle. Then $\nabla^g\colon \Gamma(\End(S))\to \Gamma(T^{\ast}M\otimes \End(S))$ acts by derivations:
\begin{equation*}
\nabla^g_v(A_1\circ A_2) = \nabla^g_v(A_1)\circ A_2 + A_1\circ \nabla^g_v(A_2) \quad \forall\,\, A_1, A_2 \in \Gamma(\End(S))\quad \forall\,\, v\in \fX(M)\, . 
\end{equation*}

\noindent
Moreover, $\Psi_{\Gamma}$ induces a unital isomorphism of algebras $(\Omega(M),\diamond_g)\simeq \Gamma(\End(S))$ which is compatible with $\nabla^g$. In other words, the following equation holds:
\begin{equation*}
\nabla^{g}_v(\Psi_{\Gamma}(\alpha)) = \Psi_{\Gamma}(\nabla^g_v\alpha)  
\end{equation*}

\noindent 
for every $\alpha \in \Omega(M)$ and $v\in \fX(M)$.
\end{prop}

\begin{proof}
The fact that $\nabla^g$ acts by algebra derivations of $\Gamma(\End(S))$ is standard. Proposition \ref{prop:SpinorialConnection} gives:
\begin{equation*}
(\nabla^g_v A)(\epsilon)= \nabla^g_v A(\epsilon) - A(\nabla^g_v\epsilon) = \nabla^g_v(\Psi_{\Gamma}(\alpha)(\epsilon)) - \Psi_{\Gamma}(\alpha)(\nabla^g_v\epsilon) = \Psi_{\Gamma}(\nabla^g_v\alpha)(\epsilon)
\end{equation*} 

\noindent
for all $A\in \Gamma(\End(S))$, $\epsilon \in \Gamma(S)$ and $v\in \fX(M)$, where $\alpha :=  \Psi_\Gamma^{-1}(A)\in \Omega(M)$.
\end{proof} 


\subsection{Constrained differential spinors}


In this section we introduce the notion of \emph{constrained differential spinor}, which constitute the main object of study in this dissertation.

\begin{definition}
\label{def:generalizedKS}
Let $(S,\Gamma)$ be a real spinor bundle on $(M,g)$ equipped with a connection $\cD$ and let $\cQ\in \Gamma(\End(S)\otimes \cW)$ be an endomorphism of $S$ taking values in an auxiliary vector bundle $\cW$ on $M$. A section $\epsilon \in \Gamma(S)$ is a {\em constrained differential spinor} with respect to $(\cD,\cQ)$ if:
\begin{equation}
\label{GKSE}
\cD\epsilon = 0\, , \qquad \cQ(\epsilon) = 0\, .
\end{equation}
\end{definition}

\begin{remark}
Supersymmetric solutions of supergravity theories can often be characterized in terms of manifolds admitting certain systems of constrained differential spinors, see for instance \cite{LazaroiuB,LazaroiuBII}. This extends the notion of \emph{generalized Killing spinors} considered \cite{BarGM,FriedrichKim,FriedrichKimII,MoroianuSemm}.
\end{remark}

\noindent
Suppose that $(S,\Gamma)$ is an irreducible real spinor bundle. Then we can write $\cD=\nabla^g-\cA$ with $\cA\in \Omega^1(End(S))$. In this case, the equations satisfied by a constrained differential spinor can be written as:
\begin{equation*}
\nabla^g\epsilon = \cA (\epsilon)\, , \quad \cQ(\epsilon) = 0
\end{equation*}

\noindent
and their solutions are called constrained differential spinors {\em relative to $(\cA,\cQ)$}. Using connectedness of $M$ and the parallel transport of $\cD$, equation \eqref{GKSE} implies that the space of constrained differential spinors relative to $(\cA,\cQ)$ is finite-dimensional and that a constrained differential spinor which is not zero at some point of $M$ is automatically nowhere-vanishing on $M$; in this case, we say that such $\epsilon$ is {\em nontrivial}.


\subsection{Spinor square maps}


Let $(S,\Gamma,\cB)$ be a paired irreducible real spinor bundle on $(M,g)$. The admissible pairing $\cB$ allows us to construct point-wise extensions to $M$ of the  spinor square maps $\cE_{\kappa}^{\gamma}\colon \Sigma \to \wedge V^{\ast}$. $\kappa \in \mathbb{Z}_2$, of Section \ref{sec:SpinorsAsPolyforms}. We denote these by:
\begin{equation*}
\cE_{\kappa}^{\Gamma} \colon S\rightarrow   \wedge M 
\end{equation*}

\noindent
which fit in the following commutative diagram:

\begin{center}
	\begin{tikzcd}
		\Cl(M,g)  \arrow[r,"\gamma"]   & (\End(S),\circ)   \arrow[dl,"\Psi^{-1}_{\Gamma}"] &  \arrow[l,"\cE_{\kappa}"] S  \\
		(\wedge M, \diamond_g) \arrow[u,"\Psi"]  	& 
	\end{tikzcd}
\end{center}

\noindent
which extends to maps of sections that we denote by the same symbol for ease of notation. Although $\cE_{\kappa}^{\Gamma}$ preserves fibers, it is not a morphism of vector bundles since it is fiberwise quadratic. By the results of Section \ref{sec:SpinorsAsPolyforms}, this map is two to one away from the zero section of $S$, where it branches, and its image is a subset of the total space of $\wedge M$ which fibers over $M$. We have $\Im(\cE_{-}^{\Gamma})=-\Im(\cE_{+}^{\Gamma})$ and $\Im(\cE_{+}^{\Gamma})\cap \Im(\cE_{-}^{\Gamma})=0$. The fiberwise sign action of $\Z_2$ on $S$ permutes the sheets of these covers (fixing the zero section), hence $\cE_{\kappa}^{\Gamma}$ gives a bijection from $S/\Z_2$ to its image as well as a single bijection:
\begin{equation*}
\cE^{\Gamma} \colon S/\Z_2 \xrightarrow{\sim}(\Im(\cE_{+}^{\Gamma})\cup \Im(\cE_{-}^{\Gamma}))/\Z_2
\end{equation*}

\noindent
The sets $\Im(\dot{\cE}_{\kappa}^{\Gamma}):= \Im(\cE_{\kappa}^{\Gamma})\setminus 0$, $\kappa \in \mathbb{Z}_2$, are connected submanifolds of the total space of $\wedge M$ and the restrictions:
\begin{equation}
\label{dcE}
\cE_{\kappa}^{\Gamma}:\dot{S}\rightarrow \Im(\dot{\cE}_{\kappa}^{\Gamma})
\end{equation}

\noindent
of $\cE_{\kappa}^{\Gamma}$ away from the zero section are surjective morphisms of fiber bundles which are two to one.

\begin{definition}
The {\em signed spinor square maps} of the paired irreducible real spinor bundle $(S,\Gamma,\cB)$ are the maps $\cE_{\kappa}^{\Gamma}\colon \Gamma(S)\to \Omega(M)$, $\kappa\in \mathbb{Z}_2$, induced by $\cE_{\kappa}^{\Gamma}$ on sections.
\end{definition}

\begin{definition}
A {\em spinorial polyform} is a a polyform $\alpha\in \Gamma(\wedge M)$ that belongs to the image of either $\cE_{+}^{\Gamma}$ or $\cE_{-}^{\Gamma}$.
\end{definition}

\noindent 
By the results of Section \ref{sec:SpinorsAsPolyforms}, $\cE_{\kappa}^{\Gamma}$ are quadratic maps of $\cC^\infty(M)$-modules and satisfy:
\begin{equation*}
\supp(\cE_{\kappa}^{\Gamma}(\epsilon))=\supp(\epsilon) \quad \forall\,\, \epsilon\in \Gamma(S)\, . 
\end{equation*}

\noindent
In the following denote by $\cE_{\kappa}^{\Gamma}(\Gamma(S))$ the image of $\cE_{\kappa}^{\Gamma} \colon \Gamma(S) \to \Omega(M)$ as a map of sections. Likewise we set $\cE^{\Gamma}(\Gamma(S)) := \cE_{+}^{\Gamma}(\Gamma(S))\cup \cE_{-}^{\Gamma}(\Gamma(S))$. 

\begin{prop}
\label{eq:obstructionliftpoly}
Let $(S,\Gamma,\cB)$ be a paired real irreducible spinor bundle associated to a $\Spin_o(V^{\ast},h^{\ast})$-structure $Q$ on $(M,g)$. Then every nowhere-vanishing polyform in $\Im(\cE^{\Gamma})$ determines a cohomology class $c_Q(\alpha)\in H^1(M,\Z_2)$ encoding the obstruction to existence of a globally-defined spinor $\epsilon\in \Gamma(S)$, necessarily nowhere-vanishing, such that $\alpha\in \{\cE_{+}^{\Gamma}(\epsilon) , \cE_{-}^{\Gamma}(\epsilon)\}$. In particular, such $\epsilon\in\Gamma(S)$ exists if and only if $c_\cQ(\alpha) = 0$.  
\end{prop}

\begin{proof}
We have $\alpha\in \Im(\cE^{\Gamma}_{\kappa})$ for some $\kappa\in \mathbb{Z}_2$. Let $L_\alpha$ be the real line sub-bundle of $\wedge M$ determined as the span of $\alpha$. Since the projective spinor square map $\P\cE_\bS:\P(S)\rightarrow\P(\wedge(M))$ is bijective when correstricted to its image, $L_\alpha$ determines a real line sub-bundle $L_\cQ(\alpha) :=  (\P\cE_\bS)^{-1}(L_\alpha)$ of $S$. A section $\epsilon$ of $S$ such that $\cE^\kappa_\bS(\epsilon)=\alpha$ is a section of $L_\cQ(\alpha)$. Since such $\epsilon$ must be nowhere-vanishing (because $\alpha$ is), it exists if and only if $L_\cQ(\alpha)$ is trivial, which happens if and only if its first Stiefel-Whitney class vanishes. The conclusion follows by setting $c_\cQ(\alpha)  :=  w_1(L_\cQ(\alpha))\in H^1(M,\Z_2)$. Notice that $c_\cQ(\alpha)$ depends only on $\alpha$ and $\cQ$, since the Clifford bundle $(S,\Gamma)$ is associated to $\cQ$ while all admissible pairings of $(S,\Gamma)$ are related to each other by automorphisms of $S$, see Remark \ref{rem:cBrelation}.
\end{proof}

\noindent
We will refer to the cohomology class $c_\cQ(\alpha)\in H^1(M,\Z_2)$ occurring the previous proposition as the {\em spinor class} of the nowhere-vanishing polyform $\alpha\in \Im(\cE^{\Gamma}_{\kappa})$. Note that $c_\cQ(\alpha)$ is not a characteristic class of $S$, since it depends on $\alpha$. 

\begin{lemma}
\label{lemma:L}
Let $(S,\Gamma,\cB)$ be a paired real spinor bundle on $(M,g)$, $(S_L,\Gamma_L)$ be the modification of $(S,\Gamma)$ by a real line bundle $L$ on $M$. For every vector bundle trivialization $q\colon L\otimes L \to M\times \mathbb{R}$, $\cB$ extends naturally an admissible bilinear pairing $\cB_L$ on $(S_L,\Gamma_L)$, making $(S_L,\Gamma_L,\cB_L)$ into a paired real spinor bundle. 
\end{lemma}

\begin{proof}
Recall from Proposition \ref{prop:Ltensor} that $\Gamma_L=\varphi_L^{-1}\circ \Gamma$. A simple computation gives:
\begin{equation*}
\cB_L(\xi_1\otimes l_1,\xi_2\otimes l_2)=q(l_1\otimes l_2)\cB(\xi_1,\xi_2)\quad \forall\,\, \xi_1,\xi_2\in \Gamma(S)\, \quad\forall\,\, l_1,l_2\in \Gamma(L)
\end{equation*}

\noindent
which immediately implies the conclusion.   
\end{proof}

\noindent 
The following proposition shows that $c_Q(\alpha)$ can be made to vanish by changing $Q$.

\begin{prop}
\label{eq:anotherspin}
Let $(S,\Gamma,\cB)$ be a paired real spinor bundle associated to a $\Spin_o(V^{\ast},h^{\ast})$-structure $Q$ on $(M,g)$. For every nowhere-vanishing polyform $\alpha\in \Im(\cE^{\Gamma})$, there exists a unique $\Spin_o(V^{\ast},h^{\ast})$-structure $Q^{\prime}$ on $(M,g)$, possibly distinct from $Q$, such that $c_{Q^{\prime}}(\alpha) = 0$.
\end{prop}

\begin{proof}
Suppose for definiteness that $\alpha\in \Im(\cE^{\Gamma}_{+})$. Let $(S,\Gamma)$ be the irreducible real spinor bundle associated to $Q$ and set $L :=  L_\cQ^+(\alpha)\subset S$. By Remark \ref{rem:stronglyspin}, isomorphism classes of $\Spin_o(V^{\ast},h^{\ast})$-structures on $(M,g)$ form a torsor over $H^1(M,\Z_2)$. Let $Q^{\prime} = c_Q(\alpha)\cdot Q$ be the spin structure obtained from $Q$ by acting in this torsor with $c_Q(\alpha)$. Then the irreducible real spinor bundle associated to $Q^{\prime}$ coincides with $(S_L,\Gamma_L)$.  Pick an isomorphism $q:L^{\otimes 2}\simeq M\times \mathbb{R}$ and equip $S_L$ with the admissible pairing $\cB_L$ constructed as in Lemma \ref{lemma:L}. Since $\Psi_{\Gamma_L}=\varphi_L^{-1}\circ \Psi_\Gamma$, it follows $\cE^{\Gamma}_{+}=\Psi_{\Gamma}^{-1}\circ (\ast\otimes \Id_S)$ and the positive spinor square maps of $(S_L,\Gamma_L,\cB_L)$ and $(S,\Gamma,\cB)$  are related through:
\begin{equation*}
\cE^{\Gamma_L}_{+} = \cE^{\Gamma}_{+} \circ \Id_{S\otimes S}\otimes q\, . 
\end{equation*}

\noindent
Since $(\Id_{S\otimes S}\otimes q)(L^{\otimes 2}\otimes L^{\otimes 2})=L^{\otimes 2}$ (where $L^{\otimes 2}$ is viewed as a sub-bundle of $S_L=S\otimes L$), this gives $\cE^{\Gamma_L}_{+}(L^{\otimes 2}\otimes L^{\otimes 2})=\cE^{\Gamma}_{+}(L\otimes L)$, which implies $\cE^{\Gamma_L}_{+}(L^{\otimes 2})=\cE^{\Gamma}_{+}(L)=L_\alpha$ Hence the line sub-bundle of $S_L$ determined by $\alpha$ is the trivializable real line bundle $L^{\otimes 2}\simeq \R_M$. Thus $c_{Q'}(\alpha) = 0$.
\end{proof}


\subsection{Differential spinors and spinorial exterior forms}


Let $(S,\Gamma,\cB)$ be a paired spinor bundle. Let $\mathfrak{q} :=  (\Psi_{\Gamma}\otimes \Id_\cW)\in \Omega^\ast(M,\cW)$ be the symbol of $\cQ\in \Gamma(\End(S)\otimes \cW)$ (see Definition \ref{def:symbol}). Proposition \ref{prop:constraintendopoly} implies:

\begin{lemma}
\label{lemma:constrainedspinor}
A spinor $\epsilon\in \Gamma(S)$ satisfies:
\begin{equation*}
\cQ(\epsilon) = 0
\end{equation*}

\noindent
if and only if one of the following mutually-equivalent relations holds:
\begin{equation*}
\mathfrak{q}\diamond_g \alpha = 0\, , \quad \alpha \diamond_g (\pi^{\frac{1-\sigma}{2}}\circ\tau)(\mathfrak{q}) = 0
\end{equation*}

\noindent
where $\alpha :=  \cE^{\Gamma}_{+}\in \Omega(M)$ is the positive polyform square of $\epsilon$.
\end{lemma}

\noindent 
Now assume that $(S,\Gamma,\cB)$ is an irreducible real spinor bundle on $(M,g)$. Set $\cA :=  \nabla^S-\cD\in \Omega^1(M,\End(S))$ and let $\fra  :=  (\Psi_{\Gamma}\otimes \Id_{T^\ast M})^{-1}(\cA)\in \Omega^1(M, \wedge M)$ be the symbol of $\cA$, viewed as a $\wedge M$-valued one-form. We have:

\begin{lemma}
\label{lemma:GKSiff}
A spinor $\epsilon \in \Gamma(S)$ satisfies $\cD\epsilon = 0$ if and only if:
\begin{equation}
\label{eq:GKSI}
\nabla^g\alpha = \fra\diamond \alpha +  \alpha \diamond (\pi^{\frac{1-\sigma}{2}}\circ\tau)(\fra)
\end{equation}

\noindent
where $\alpha  :=  \cE^{\Gamma}_{+}(\varepsilon)$ is the positive polyform square of $\epsilon$.
\end{lemma}

\begin{proof}
Assume that $\epsilon$ satisfies $\nabla^g\epsilon = \cA(\epsilon)$. We have $E_{\varepsilon} := \cE_{+}(\varepsilon) \in \Gamma(\End(S))$ and:
\begin{eqnarray*}
(\nabla^g E_{\varepsilon})(\chi) = \nabla^g(E_{\varepsilon}(\chi)) - E_{\varepsilon}(\nabla^g\chi) = \nabla^g(\cB(\chi,\epsilon)\,\epsilon) - \cB(\nabla^g\chi,\epsilon)\,\epsilon \\ 
= \cB(\chi,\nabla^g\epsilon)\,\epsilon + \cB(\chi,\epsilon)\,\nabla^g\epsilon = \cB(\chi,\cA(\epsilon))\,\epsilon + \cB(\chi,\epsilon)\,\cA(\epsilon) = E_{\varepsilon}(\cA^t (\chi)) + \cA(E_{\varepsilon}(\chi))
\end{eqnarray*}

\noindent
for all $\chi \in \Gamma(S)$, where $\cA^t$ denotes the adjoint of $\cA$ with respect to $\cB$. The previous equation implies:
\begin{equation}
\label{eq:DEepsilon}
\nabla^g E_{\varepsilon} =  \cA\circ E_{\varepsilon} + E_{\varepsilon}\circ\cA^t\, . 
\end{equation}

\noindent
Applying $\Psi^{-1}_{\Gamma}$ and using Lemma \ref{lemma:adjointpoly} and Proposition \ref{prop:LCderClifford} gives \eqref{eq:GKSI}. Conversely, assume that $\alpha$ satisfies \eqref{eq:GKSI}. Applying $\Psi_{\Gamma}$ gives equation \eqref{eq:DEepsilon}, which in turn can be rewritten as follows:
\begin{equation}
\label{eq:Depsilonrelation}
\cB(\chi,\cD_{v}\epsilon)\,\epsilon + \cB(\chi,\epsilon)\,\cD_{v}\epsilon = 0  
\end{equation}

\noindent
for every $\chi \in \Gamma(S)$ and every $v\in \fX(M)$. Hence $\cD_{v}\epsilon = \beta(v) \epsilon$ for some $\beta \in \Omega^{1}(M)$. Using this in \eqref{eq:Depsilonrelation} gives:
\begin{equation*}
\cB(\chi,\epsilon)\,\beta\otimes \epsilon = 0 \quad \forall\,\, \chi \in \Gamma(S)\, . 
\end{equation*}

\noindent
This implies $\beta=0$, since $\cB$ is non-degenerate and $\epsilon$ is nowhere-vanishing. Hence $\cD\epsilon = 0$.
\end{proof}

\begin{remark}
If $\cA$ is skew-symmetric with respect to $\cB$, then \eqref{eq:GKSI} simplifies to:
\begin{equation}
\label{eq:GKSII}
\nabla^g\alpha = \fra\diamond \alpha -  \alpha\diamond \fra\, . 
\end{equation}

\noindent
For the type of differential spinors occurring in the supersymmetric configurations of supergravity $\cA$ need {\em not} be skew-symmetric relative to $\cB$.
\end{remark}

\noindent
We arrive to the final characterization of differential spinors in terms of their associated spinorial polyform that we will use throughout this dissertation.

\begin{thm}
\label{thm:GCKS}
Let $(S,\Gamma,\cB)$ be a paired spinor bundle associated to a $\Spin_o(V^{\ast},h^{\ast})$-structure $Q$ and whose admissible form $\cB$ has adjoint type $\sigma\in\mathbb{Z}_2$ and symmetry type $s\in\mathbb{Z}_2$. Let $\cA\in \Omega^1(M,\End(S))$ and  $Q\in \Gamma(\End(S)\otimes \cW)$. Then the following statements are equivalent:
\begin{enumerate}
\itemsep 0.0em
\item There exists a nontrivial  constrained differential spinor $\epsilon\in \Gamma(S)$ with respect to $(\cA,\cQ)$.
\item There exists a nowhere-vanishing polyform $\alpha\in \Omega(M)$ with vanishing cohomology class $c_\cQ(\alpha)$ which satisfies the following algebraic and differential equations for every polyform $\beta \in\Omega(M)$:
\begin{equation}
\label{eq:Fierzglobal}
\alpha\diamond_g \beta\diamond_g \alpha = 2^{\frac{d}{2}} (\alpha \diamond_g \beta)^{(0)} \alpha\, , \quad (\pi^{\frac{1-\sigma}{2}}\circ\tau)(\alpha) = s\,\alpha  
\end{equation}

\begin{equation}
\label{eq:GKSeqiff}
\nabla^g\alpha = \fra\diamond_g \alpha + \alpha\diamond_g (\pi^{\frac{1-\sigma}{2}}\circ\tau)(\fra)\, , \quad   \mathfrak{q}\diamond_g \alpha  = 0
\end{equation}

\noindent
or, equivalently, satisfies the equations:
\begin{equation}
\alpha\diamond_g\alpha = 2^{\frac{d}{2}} \alpha^{(0)} \alpha\, ,  \quad (\pi^{\frac{1-\sigma}{2}}\circ\tau)(\alpha) = s\,\alpha\, ,  \quad \alpha\diamond_g \beta\diamond_g \alpha = 2^{\frac{d}{2}} (\alpha\diamond\beta)^{(0)} \alpha 
\end{equation}
\begin{equation}
\nabla^g\alpha = \fra\diamond_g \alpha + \alpha\diamond_g (\pi^{\frac{1-\sigma}{2}}\circ\tau)(\fra)\, \, , \quad  \mathfrak{q}\diamond_g \alpha = 0 
\end{equation}

\noindent
for some fixed polyform $\beta \in\Omega(M)$ such that $(\alpha\diamond\beta)^{(0)} \neq 0$.
\end{enumerate}

\noindent
If $\epsilon\in \Gamma(S)$ has chirality $\mu\in\mathbb{Z}_2$, then we have to add the condition:
\begin{equation*}
\ast_g (\pi\circ\tau)(\alpha) = \mu \, \alpha\, . 
\end{equation*}

\noindent
The polyform $\alpha$ as above is determined by $\epsilon$ through the relation:
\begin{equation*}
\alpha=\cE_{\kappa}^{\Gamma}(\epsilon)
\end{equation*}

\noindent
for some $\kappa\in \mathbb{Z}_2$. Moreover, $\alpha$ satisfying the conditions above determines a nowhere-vanishing real spinor $\epsilon$ satisfying this relation, which is unique up to sign.
\end{thm}

\begin{remark}
\label{remark:obstructionliftingspinor} 
Suppose that $\alpha\in \Omega(M)$ is nowhere-vanishing and satisfies \eqref{eq:Fierzglobal} and \eqref{eq:GKSeqiff} but we have $c_Q(\alpha)\neq 0$. Then Proposition \ref{eq:anotherspin} implies that there exists a unique $\Spin_o(V^\ast,h^\ast)$-structure $Q^{\prime}$ such that $c_{Q^{\prime}}(\alpha) = 0$. Thus $\alpha$ is the square of a global section of a paired spinor bundle $(S',\Gamma',\cB')$ associated to $\cQ'$. Hence a nowhere-vanishing polyform $\alpha$ satisfying \eqref{eq:Fierzglobal} and \eqref{eq:GKSeqiff} corresponds to the square of a differential spinor with respect to a uniquely-determined $\Spin_o(V^\ast,h^\ast)$-structure.
\end{remark}

\begin{proof} 	
The algebraic conditions in the Theorem follow from the pointwise extension of Theorem \ref{thm:reconstruction} and Corollary \ref{cor:reconstructionchiral}. The differential condition follows from Lemma \ref{lemma:GKSiff}, which implies that $\nabla^g\epsilon = \cA(\epsilon)$ holds if and only if \eqref{eq:GKSeqiff} does upon noticing that $\epsilon \in \Gamma(S)$ vanishes at a point $m\in M$ if and only if its positive polyform square $\alpha$ satisfies $\alpha\vert_m = 0$. Condition $c_Q(\alpha) = 0$ follows from Proposition \ref{eq:obstructionliftpoly}.
\end{proof}
 
\noindent
In practical applications, such as in the study of supersymmetric configurations of supergravity theories, the specific form of $\cA$ and $\cQ$ as vector-valued endomorphisms of $S$ is usually not relevant. Instead, only their \emph{symbols} are relevant. This is because usually the supersymmetry conditions are directly expressed in terms of polyforms acting via Clifford multiplication on locally defined spinors. If this type of situations we would be interested in studying differential spinors with respect to any pair $(\cA,\cQ)$ that has a given fixed symbol $(\fra,\frq)$ (which does not depend on the choice of paired spinor bundle but only on the underlying pseudo-Riemannian structure). In fact, if $H^1(M,\mathbb{Z}_2) = 0$ then there is a unique paired spinor bundle modulo isomorphism and in this case $(\cA,\cQ)$ and $(\fra,\frq)$ are equivalent via a canonical isomorphism. 

\begin{cor}
Let $(M,g)$ be a strongly spin pseudo-Riemmanian manifold of signature $(p-q) \equiv_8 0,2$ and let $\fra \in \Omega^1(M,\wedge M)$ and $\frq \in \Gamma(\wedge M\otimes \cW)$ given. Then, $(M,g)$ admits a constrained differential spinor with respect to any pair $(\cA,\cQ)$ whose symbol is $(\fra,\frq)$ if and only if there exists a polyform $\alpha \in \Omega(M)$ satisfying the following algebraic and differential equations:
\begin{eqnarray*}
& \alpha\diamond_g \alpha = 2^{\frac{d}{2}} \alpha^{(0)} \alpha\, ,  \quad (\pi^{\frac{1-\sigma}{2}}\circ\tau)(\alpha) = s\,\alpha\, ,  \quad \alpha\diamond_g \beta\diamond_g \alpha = 2^{\frac{d}{2}} (\alpha\diamond\beta)^{(0)} \alpha\\ 
& \nabla^g\alpha = \fra\diamond_g \alpha + \alpha\diamond_g (\pi^{\frac{1-\sigma}{2}}\circ\tau)(\fra)\, \, , \quad  \frq \diamond_g \alpha = 0 
\end{eqnarray*}
	
\noindent
for a fixed polyform $\beta \in\Omega(M)$ such that $(\alpha\diamond_g\beta)^{(0)} \neq 0$.
\end{cor}

\noindent
The previous corollary is especially well-adapted to study the geometric and topological consequences of the existence of constrained differential spinors when the specific expression for the spinors themselves is not in itself relevant. We will make extensive use of it throughout this dissertation. 


\renewcommand{\chaptername}{Chapter}

\renewcommand{\leftmark}{Chapter \thechapter. Irreducible spinors on Lorentzian four-manifolds}

\chapter{Irreducible spinor bundles and parabolic pairs}
\label{chapter:IrreducibleSpinors4d}


In this chapter we specialize the theory of differential spinors and spinorial polyforms developed in the previous chapter to the specific case of Lorentzian four-manifolds, which will be the case of interest for the remainder of this dissertation. We begin with the study of the algebraic theory of irreducible real spinors in four Lorentzian dimensions to then proceed to investigate the most general type of differential spinor in this dimension and signature.

 
 \section{Algebraic spinorial polyforms}
 \label{sec:4dAlgebraic}
 


\subsection{Parabolic pairs on Minkowski space} 

 
Let $(V,h)$ be a Minkowski space of \emph{mostly plus} signature $(-+++)$ which we consider to be oriented and time-oriented. That is, we fix a time-like element in $V$ and a volume form $\nu_h\in \wedge^4 V$, which is necessarily of negative unit norm and therefore also time-like. The dual quadratic vector space, endowed with the induced orientation, will be denoted by $(V^{\ast},h^{\ast})$. The unique, modulo isomorphism, irreducible real Clifford module $(\Sigma,\gamma)$ of $\Cl(V^{\ast},h^{\ast})$ is four-dimensional and its admissible pairings $\cB_{+}$ and $\cB_{-}$ are both skew-symmetric, see the table provided in Theorem \ref{thm:admissiblepairings}. In the following we choose to work with the admissible pairing $\cB=\cB_{-}$ of negative adjoint type since we find it yields a spinorial polyform that is more convenient for computations, although we could have equivalently chosen to work with $\cB_{+}$. For every non-zero isotropic vector $u\in V^{\ast}$ we define the following equivalence relation $\sim_u$ on $V^{\ast}$:
\begin{equation*}
l_1 \sim_u l_2 \quad \mathrm{if\,and\,only\,if}\quad l_1 = l_2 + c\, u\, , \quad c\in \mathbb{R}\, .
\end{equation*}
 
\noindent
We will refer to transformations of the form $V^{\ast}\ni l \mapsto l + c u\in V^{\ast}$ as \emph{gauge transformations} generated by $u$, and we will denote by $[l]_u$ the equivalence class determined by any element $l\in V^{\ast}$.  

\begin{definition}
\label{def:parabolicvector}
A {\em parabolic pair} is a pair $(u,[l]_u) \in V^{\ast}\times (V^{\ast}/\sim_u)$ satisfying:
\begin{equation}
\label{eq:ulcond}
u\neq 0 \, , \quad h^\ast (u,u) = 0\, ,\quad  h^\ast (l,l) = 1\, ,\quad  h^{\ast}(u,l) = 0
\end{equation}

\noindent
for any, and hence for all, representatives $l\in [l]_u$. 
\end{definition}
 
\noindent 
Non-zero elements in $V^{\ast}$ or $V$ of zero norm will be called \emph{isotropic} in the following. We denote by $\mathfrak{P}(V)$ the category whose objects $\Ob(\mathfrak{P}(V))$ are tuples $(h,u,[l]_u)$ consisting of a Lorentzian metric $h$ on $V$ and a parabolic pair $(u,[l]_u)$ on $(V,h)$, and whose morphisms are invertible isometries that preserve the corresponding parabolic pairs. That is, if $(h_1,u_1,[l_1]_{u_1})$ and  $(h_2,u_2,[l_2]_{u_2})$ are elements in $\Ob(\mathfrak{P}(V))$, then a morphism:
\begin{equation*}
T \colon (h_1,u_1,[l_1]_{u_1}) \to  (h_2,u_2,[l_2]_{u_2})
\end{equation*}

\noindent
in the category $\mathfrak{P}(V)$ is a linear map:
\begin{equation*}
T\colon V\to V
\end{equation*}

\noindent
such that:
\begin{equation*}
T^{\ast} h_1  = h_2 \, , \qquad T^{\ast} u_1 = u_2 \, , \qquad T^{\ast}([l_1]_{u_1}) := [T^{\ast}(l_1)]_{T^{\ast}(u_1)} = [l_2]_{u_2} 
\end{equation*}

\noindent
The last equation above is consistent since $T$ is an isometry that preserves $u$. Hence, for any $c\in \mathbb{R}$ we have:
\begin{equation*}
T^{\ast}(l_1 + c u_1) = T^{\ast}(l_1) + c u_2
\end{equation*}

\noindent
showing that the action of $T\colon V\to V$ on $[l_1]_{u_1}$ does not depend on the representative chosen. The main theorem of this section is the following characterization of irreducible real spinors on four-dimensional Minkowski space, which justifies the notion of parabolic pair introduced above.
 
\begin{thm}
\label{thm:squarespinorMink}
A polyform $\alpha\in \wedge V^{\ast}$ is the square of a nonzero spinor, that is, it belongs to the set $Z_{-,-}(V^\ast,h^\ast)$, if and only if there exists a parabolic pair $(u,[l]_u)\in \Ob(\mathfrak{P}(V))$ such that
\begin{equation}
\label{eq:alphaul}
\alpha = u + u\wedge l
\end{equation}

\noindent
where $l\in V^{\ast}$ is any representative of $[l]_u \in  V^{\ast}/\sim_u$ . If this is the case, the parabolic pair $(u,[l]_u)$ is uniquely determined by $\alpha$.  
\end{thm}
 
\begin{proof}
Let:
\begin{equation*}
\alpha = \sum_{k=0}^{4} \alpha^{(k)}\in \wedge V^\ast\, , \quad \alpha^{(k)}\in \wedge^k V^\ast\quad  \forall\,\, k=0,\ldots, 4\, .
\end{equation*}

\noindent
By Theorem \ref{thm:reconstruction}, $\alpha$ lies in $Z_{-,-}(V^\ast,h^\ast)$ if and only if the following relations hold for $\beta =1 $ and for a polyform $\beta$ such that $(\beta\diamond\alpha)^{(0)} \neq 0$:
\begin{equation}
\label{eq:eqs3,1}
\alpha\diamond\beta\diamond \alpha = 4\, (\beta\diamond\alpha)^{(0)}\, \alpha\, , \qquad (\pi\circ\tau)(\alpha) = -\alpha\, .
\end{equation}

\noindent
Condition $(\pi\circ\tau)(\alpha) = -\alpha$ implies $\alpha^{(0)} = \alpha^{(3)} =\alpha^{(4)} = 0$.  Thus $\alpha = u + \omega$, where $u :=  \alpha^{(1)}\in \wedge^1 V^\ast$ and $\omega  := \alpha^{(2)}\in \wedge^2 V^\ast$. For $\beta =1$, the first condition in \eqref{eq:eqs3,1} gives $(u + \omega)\diamond (u + \omega) = 0$, which reduces to the following relations upon expanding the geometric product:
\begin{equation}
\label{eq:eqseqiv3,1}
h^{\ast}(u,u) = \langle \omega,\omega\rangle_h\, , \qquad \omega\wedge u = 0\, .
\end{equation}

\noindent 
Here $\langle\cdot,\cdot\rangle_h$ is the determinant metric induced by $h$ on $\wedge V^\ast$. The second condition in \eqref{eq:eqseqiv3,1} amounts to $\omega = u\wedge l$ for some $l\in V^{\ast}$ determined up to gauge transformations generated by $u\in V^{\ast}$. Using this in \eqref{eq:eqseqiv3,1} gives the condition:
\begin{equation}
\label{eq:norms}
(h^\ast(l,l)-1)\, h^\ast (u,u) = h^{\ast}(u,l)^2
\end{equation}

\noindent
which is invariant under gauge transformations generated by $u$. For $\beta = u$, the first equation in \eqref{eq:eqs3,1} amounts to $h^\ast(u,u) = 0$, whence the right-hand side of \eqref{eq:norms} must vanish, implying $h^{\ast}(u,l)= 0$. It remains to show that $h^\ast(l,l)=1$. Since $u$ is non-zero and null, there exists a non-zero null one-form $v\in V^{\ast}$ such that $h^\ast(v,u)= 1$. We have $(v\diamond v)^{(0)}=(v\diamond (u + u\wedge l))^{(0)} = h^\ast(v,u)=1$ and therefore $(v\diamond \alpha)^{(0)}\neq 0$, as required in Theorem \ref{thm:reconstruction} for $v$ to be an appropriate choice of $\beta$. Taking $\beta=v$, the first condition in \eqref{eq:eqs3,1} reduces to:
\begin{equation*}
(u + u\wedge l)\diamond v \diamond (u + u\wedge l) =  4\,(u + u\wedge l)~~.
\end{equation*}

\noindent
A direct computation shows that this equation amounts to $h^\ast(l,l) = 1$ and thus we conclude.
\end{proof}
 
\begin{remark} 
Given the time-like orientation $\frt\in V^{\ast}$ of $V^{\ast}$, denote by $P_{\frt}\colon V^{\ast}\to \R\, \frt$ the orthogonal projection onto the line spanned by $\frt$. It follows that for every parabolic pair $(u,[l]_u)$ on $(V,h)$ a canonical choice of representative $l\in [l]_u$ is obtained by imposing the condition:
\begin{equation*}
P_{\frt}(l) = 0\, .
\end{equation*}

\noindent
Indeed, given $l\in V^\ast$ of unit norm and orthogonal to $u$, there exists a unique $c\in \R$ such that $P_{\frt} (l+c\, u) = 0$. As we will see in Chapter \ref{chapter:Globallyhyperbolicsusy} this \emph{gauge choice} is useful to study spinors on globally hyperbolic Lorentzian four-manifolds.
\end{remark}
  
\begin{remark}
As a consistency check we can informally count the degrees of freedom encoded in $\alpha = u + u\wedge l$. A priori, the null one-form $u$ has three degrees of freedom while the space-like one-form $l$ has four, which are reduced to two by the requirements that $l$ be of unit norm and orthogonal to $u$. Since $l$ is defined only up to gauge transformations $l\mapsto l + c\, u$, $c \in \R$, its number of degrees of freedom further reduces from two to one. This gives a total of four degrees of freedom, matching those of an irreducible real spinor in four-dimensional Lorentzian signature.
\end{remark}


\subsection{Null basis and parabolic pairs} 


Let $(u,[l]_u)$ be a parabolic pair on $(V,h)$. The orthogonal complement $u^{\perp_h}\subset V^{\ast}$ of $u\in V^{\ast}$ in $(V^{\ast},h^{\ast})$ is a three-dimensional vector subspace of $V^{\ast}$ that contains again $u\in u^{\perp_h}$. Hence, the quotient:
\begin{equation*}
\frS_u = \frac{u^{\perp_h}}{u}
\end{equation*}

\noindent
is well-defined and determines a two-dimensional vector space on which the metric $h^{\ast}$ restricts to a Riemannian metric that we denote by $q_u$. Since $V$ is oriented and time-oriented, $\frS_u$ inherits a canonical orientation which, together with the metric $q_u$ defines a Hodge operator that we denote by $\ast_{q_u} \colon \wedge \frG_u \to \wedge \frG_u$. Applying this Hodge operator to $[l]_u$ we obtain a canonical element:
\begin{equation*}
[n]_u = \ast_{q_u} [l]_u
\end{equation*}

\noindent
for a unit one-form $n\in V^{\ast}$ orthogonal to $u$ and unique modulo gauge transformations of the form $n\mapsto n + c u$ with $c\in \mathbb{R}$. In particular, we can think of $(u,[n]_u)$ as a parabolic pair canonically associated to $(u,[l]_u)$. This implies, in particular, that parabolic pairs occur always in canonical pairs. Given $(u,[l]_u)$, when necessary we will refer to $(u,[n]_u)$ as its \emph{associated parabolic pair}. We obtain a short exact sequence:
\begin{eqnarray}
\label{eq:shortexactu}
0 \to \langle \mathbb{R} u\rangle \to (u^{\perp_h},h^{\ast}\vert_{u^{\perp_h}}) \to (\frS_u,h^{\ast}_u)  \to 0
\end{eqnarray}

\noindent
where $\langle \mathbb{R} u\rangle$ denotes the line spanned by $u$ in $V^{\ast}$. Given an isotropic vector $u\in V^{\ast}$, we say that an element $v\in V^{\ast}$ is conjugate to $u$ if $v$ is also isotropic and $h^{\ast}(u,v) =1$. Similarly, given a parabolic pair $(u,[l]_u)$, we say that an element $v\in V^{\ast}$ is conjugate to $(u,[l]_u)$ if $v$ is conjugate to $u$. Conjugate vectors to a given isotropic vector $u\in V^{\ast}$ are not unique. Instead, a direct computation gives the following result.

\begin{prop}
\label{prop:conjugatevariation}
Let $v\in V^{\ast}$ be conjugate to an isotropic vector $u\in V^{\ast}$. Any other vector $v^{\prime}\in V^{\ast}$ conjugate to $u$ is given by:
\begin{equation}
\label{eq:conjugatevectors}
v^{\prime} = v -\frac{1}{2} \vert \frw \vert_h^2 u + \frw
\end{equation}

\noindent
for a unique element $\frw\in ( \langle \mathbb{R}\, u\rangle  \oplus \langle \mathbb{R}\, v \rangle)^{\perp_h}$, where $( \langle \mathbb{R}\, u\rangle  \oplus \langle \mathbb{R}\, v \rangle)^{\perp_h} \subset V^{\ast}$ denotes the orthogonal complement of the span on $u$ and $v$ in $V^{\ast}$.
\end{prop}

\begin{remark}
Vectors conjugate to a given isotropic vector $u\in V^{\ast}$ are in one-to-one correspondence with splittings of the short exact sequence \eqref{eq:shortexactu}. Such splittings are in turn an affine space of $\mathbb{R}^2$, as expected from the previous proposition. 
\end{remark}

\begin{prop}
\label{prop:adaptedbasis}
Let $(u,[l]_u)$ be a parabolic pair. Every choice of isotropic vector $v\in V^{\ast}$ conjugate to $(u,[l]_u)$ canonically determines a positively oriented basis $(u,v,l,n)$ of $V^{\ast}$, where $l\in [l]_u$ and $n\in [n]_u$ are unit vectors uniquely determined by:
\begin{equation*}
l,n\in ( \langle \mathbb{R}\, u\rangle  \oplus \langle \mathbb{R}\, v \rangle)^{\perp_h}
\end{equation*}

\noindent
that is, uniquely determined by the condition of being orthogonal to the span of $u$ and $v$.
\end{prop}

\begin{proof}
Let $(u,[l]_u)$ be a parabolic pair and let $v\in V^{\ast}$ be conjugate to $(u,[l]_u)$. Let $l\in [l]_u$ be the unique representative in $[l]_u$ that is orthogonal to $v$. Such $l$ is orthogonal to both $u$ and $v$, and thus the three-dimensional real span of $(u,v,l)$ is a non-degenerate vector subspace of $(V^{\ast},h^{\ast})$. Therefore, there exists a unique element $n\in V^{\ast}$ such that:
\begin{eqnarray*}
h(u,n) = h(v,n) = h(l,n) = 0\, , \qquad \vert n \vert_h^2 = 1 
\end{eqnarray*}

\noindent
and such that the basis $(u,v,l,n)$ of $V^{\ast}$ is positively oriented. Equivalently, $n$ can be explicitly defined by the equation $n = - \ast_h (u\wedge v \wedge l)$.
\end{proof}

\noindent
Therefore, a choice of isotropic vector $v\in V^{\ast}$ conjugate to a parabolic pair $(u,[l]_u)$ determines a canonical basis $(u,v,l,n)$ of $V^{\ast}$ to which we will refer as the \emph{null basis} determined by $(u,[l]_u)$ and $v$. Given $(u,v,l,n)$, we denote again by $q_u$ the restriction of $h$ to $( \langle \mathbb{R}\, u\rangle  \oplus \langle \mathbb{R}\, v \rangle)^{\perp_h}$ since the abstract vector bundle $\frG_u$ is isomorphic to $( \langle \mathbb{R}\, u\rangle  \oplus \langle \mathbb{R}\, v \rangle)^{\perp_h}$. Note that $q_u = l\otimes l + n\otimes n$.
\begin{prop}
\label{prop:primerparallelization}
Let $(u,[l]_u)$ be a parabolic pair on $(V,h)$. Given null basis $(u,v,l,n)$ and $(u,v^{\prime}, l^{\prime}, n^{\prime})$ associated to $(h,u,[l]_u)$, there exists a unique vector $\frw \in ( \langle \mathbb{R}\, u\rangle  \oplus \langle \mathbb{R}\, v \rangle)^{\perp_h}$ such that:
\begin{equation}
\label{eq:algebraictransform}
(u,v^{\prime},l^{\prime},n^{\prime}) = (u,v-\frac{1}{2} \vert \frw\vert_h^2 u + \frw,l - \frw(l)u,n - \frw(n) u)
\end{equation}

\noindent
In particular, the set of null basis associated to $(h,u,[l]_u)$ is a torsor over the vector space $( \langle \mathbb{R}\, u\rangle  \oplus \langle \mathbb{R}\, v \rangle)^{\perp_h}$.
\end{prop}
 
\begin{proof}
By Proposition \ref{prop:conjugatevariation} we have:
\begin{equation*}
v^{\prime} = v -\frac{1}{2} \vert\frw\vert_h^2 u + \frw
\end{equation*}

\noindent
for a unique element $\frw\in ( \langle \mathbb{R}\, u\rangle  \oplus \langle \mathbb{R}\, v \rangle)^{\perp_h}$. Since $(u,v^{\prime},l^{\prime},n^{\prime})$ is a positively oriented null basis of $V^{\ast}$ associated to the same parabolic pair $(u,[l]_u)$ and $(u,v,l,n)$ is in particular a basis of $V^{\ast}$, we can write:
\begin{equation*}
l^{\prime} = l + a_u u  \, , \qquad n^{\prime} = b_u u + b_v v +  b_l l + b_n n 
\end{equation*}

\noindent
for constants $a_u, b_u, b_v, b_l, b_n \in \mathbb{R}$. Imposing that:
\begin{equation*}
h^{\ast}(u,l^{\prime}) = h^{\ast}(u,n^{\prime}) = h^{\ast}(v^{\prime},l^{\prime}) = h^{\ast}(v^{\prime},n^{\prime}) = h^{\ast}(l^{\prime},n^{\prime}) = 0\, ,\quad h^{\ast}(l^{\prime},l^{\prime}) = h^{\ast}(n^{\prime},n^{\prime})  = 1
\end{equation*}

\noindent
we obtain:
\begin{equation*}
  b_v  = b_l = 0\, , \quad a_u = -h^{\ast}(\frw,l)\, , \quad b_u = - h^{\ast}(\frw,n)\, , \quad b_n = 1
\end{equation*}

\noindent
and hence we conclude. 
\end{proof}

\noindent
Identifying $( \langle \mathbb{R}\, u\rangle  \oplus \langle \mathbb{R}\, v \rangle)^{\perp_h} = \mathbb{R}^2$ we conclude that set of null basis associated to $(h,u,[l]_u)$ is a $\mathbb{R}^2$ torsor precisely with respect to the action defined by \eqref{eq:algebraictransform}.
 
\begin{remark}
\label{remark:NewmanPenrose}
By the previous discussion, a null basis associated to the parabolic pair corresponding to a real irreducible spinor in four Lorentzian dimensions determines a \emph{real} version of the notion of \emph{complex tetrad} used in the Newman-Penrose formalism \cite{NewmanPenrose}. More precisely, given a null basis $(u,v,l,n)\in \cP(u,[l]_u)$ as introduced above, we obtain a complex null tetrad as follows:
\begin{equation*}
(u,v, \Psi = l +i n, \bar{\Psi} = l - i n)
\end{equation*}

\noindent
When necessary, we will refer to the previous tuple as the \emph{complex null tetrad} associated to the null basis $(u,v,l,n)$. 
\end{remark}

\noindent
We introduce the category $\mathfrak{F}(V)$ whose objects are equivalence classes $[u,v,l,n]$ of oriented basis $(u,v,l,n)$ of $V^{\ast}$ with equivalence relation given by:
\begin{equation*}
(u^{\prime},v^{\prime},l^{\prime},n^{\prime}) \in [u,v,l,n]
\end{equation*}

\noindent
if and only if $u^{\prime} = u$ and there exists  an element $w\in ( \langle \mathbb{R}\, u\rangle  \oplus \langle \mathbb{R}\, v \rangle)^{\perp_h}$ such that Equation \eqref{eq:algebraictransform} holds. Here the metric $h$ is given by:
\begin{equation*}
h = u\odot v + l\otimes l + n\otimes n = u^{\prime}\odot v^{\prime} + l^{\prime}\otimes l^{\prime} + n^{\prime}\otimes n^{\prime}
\end{equation*}

\noindent
Morphisms in $\mathfrak{F}(V)$ are linear automorphisms of $V$ that map equivalence classes to equivalence classes in the natural way. That is, if $[u,v,l,n],[u^{\prime},v^{\prime},l^{\prime},n^{\prime}] \in \Ob(\mathfrak{F}(V))$,  then a morphism:
\begin{equation*}
	T : [u,v,l,n] \to  [u^{\prime},v^{\prime},l^{\prime},n^{\prime}]
\end{equation*}

\noindent
in the category $\mathfrak{F}(V)$ is a linear map $T\colon V\to V$ such that:
\begin{eqnarray*}
	T([u^{\prime},v^{\prime},l^{\prime},n^{\prime}] ) := [T^{\ast} u^{\prime},T^{\ast} v^{\prime},T^{\ast} l^{\prime},T^{\ast} n^{\prime} ] = [u ,v ,l ,n ]
\end{eqnarray*}

\noindent 
This is well defined, since, for any other representative in $[u^{\prime},v^{\prime},l^{\prime},n^{\prime}]$ we have:
\begin{eqnarray*}
T^{\ast}(v^{\prime}-\frac{1}{2} \vert w^{\prime}\vert_{h^{\prime}}^2 u^{\prime} + w^{\prime}) = v -\frac{1}{2} \vert w^{\prime}\vert_{h^{\prime}}^2 u  + T^{\ast} w^{\prime} = v -\frac{1}{2} \vert T^{\ast} w^{\prime}\vert_{h}^2 u  + T^{\ast} w^{\prime}
\end{eqnarray*}

\noindent
where  $w^{\prime}\in ( \langle \mathbb{R}\, u^{\prime}\rangle  \oplus \langle \mathbb{R}\, v^{\prime} \rangle)^{\perp_{h^{\prime}}}$. Note that we have $T^{\ast} w^{\prime}\in ( \langle \mathbb{R}\, u \rangle  \oplus \langle \mathbb{R}\, v\rangle)^{\perp_{h}}$. Similarly:
\begin{eqnarray*}
	T^{\ast}(l^{\prime} - \langle w^{\prime},l^{\prime}\rangle_{h^{\prime}} u^{\prime}) =   l  - \langle T^{\ast}w^{\prime},T^{\ast}l^{\prime}\rangle_{h} u\, , \quad T^{\ast}(l^{\prime} - \langle w^{\prime},l^{\prime}\rangle_{h^{\prime}} u^{\prime}) =   l  - \langle T^{\ast}w^{\prime},T^{\ast}l^{\prime}\rangle_{h} u
\end{eqnarray*}

\noindent
Hence, the equivalence class of $[u,v,l,n]$ is preserved. Altogether previous discussion yields the following result.

\begin{prop}
\label{prop:categoricalequivalence}
There exists a natural equivalence of categories: 
\begin{equation*}
\mathbb{E}\colon \mathfrak{F}(V) \to \mathfrak{P}(V)
\end{equation*}

\noindent
that maps $(u,v,n,l)\in \Ob(\mathfrak{F}(V))$ to $(h,u,[l]_u)\in \Ob(\mathfrak{P}(V))$, where $h = u\odot v + l\otimes l + n \otimes n$. 
\end{prop}

\begin{remark}
We recall that two categories are considered to be \emph{equivalent} if there exists a fully faithful and essentially surjective functor between them.  
\end{remark}

\noindent
Proposition \ref{prop:categoricalequivalence} allows us to equivalently use parabolic pairs or equivalence families of coframes to study geometric problems involving differential spinors in four Lorentzian dimensions. This will be especially crucial in Chapter \ref{chapter:parallelspinorstorsion} to describe skew-torsion parallel spinors in terms of an exterior differential system. To end this subsection we recall the identities:
\begin{eqnarray*}
& \ast_h u = -u\wedge l\wedge n \, , \qquad \ast_h v = v\wedge l\wedge n \, , \qquad \ast_h l = u\wedge v\wedge n \, , \qquad \ast_h n = - u\wedge v\wedge l\\
& \ast_h(u\wedge v) = - l\wedge n\, , \quad \ast_h(u\wedge l) = u\wedge n \, , \quad \ast_h(v\wedge l) = - v\wedge n\, , \\
& \ast_h \nu_g = -1\, , \quad \ast_h^2\vert_{\wedge^1 V^{\ast}} = 1\, , \quad \ast_h^2\vert_{\wedge^2 V^{\ast}} = - 1\, , \quad\iota_w \ast_h \beta = \ast_h(\beta\wedge w)  
\end{eqnarray*}

\noindent
for every $w\in V^{\ast}$ and every $\beta \in \wedge V^{\ast}$. These formulas will be extensively used throughout the dissertation. 


\subsection{The stabilizer of a parabolic pair} 
\label{subsec:algebraicstabilizer}


Since the spinor square map is equivariant with respect to the natural action of $\Spin_o (V^{\ast},h^{\ast})$ on both its source and target spaces, we can study the stabilizer $\cI_{\varepsilon} \subset \Spin_o(V^{\ast},h^{\ast})$ in $\Spin_o(V^{\ast},h^{\ast})$ of a spinor $\varepsilon\in \Sigma$ through the stabilizer $\cI_{\alpha_{\varepsilon}}\subset \SO_o(V^{\ast},h^{\ast})$ of its associated spinorial polyform $\alpha_{\varepsilon}\in \wedge V^{\ast}$. More precisely, if $\pi\colon \Spin_o(V^{\ast},h^{\ast})\to \SO_o(V^{\ast},h^{\ast})$ is the double cover of $\SO_o(V^{\ast},h^{\ast})$, then:
\begin{equation*}
\pi(\cI_{\varepsilon}) = \cI_{\alpha_{\varepsilon}}
\end{equation*}

\noindent
implying that either $\cI_{\varepsilon}$ is a double cover of $\cI_{\alpha_{\varepsilon}}$ or are otherwise isomorphic. The latter case occurs in particular if $\cI_{\alpha_{\varepsilon}}$ is simply connected. As proven in Theorem \ref{thm:squarespinorMink}, the square of an irreducible real spinor in four Lorentzian dimensions is a spinorial polyform of the following type:
\begin{eqnarray*}
\alpha_{\varepsilon} = u + u\wedge l
\end{eqnarray*}

\noindent
for a uniquely determined null one-form $u\in V^{\ast}$ and a unit one-form $l\in V^{\ast}$ that is uniquely determined modulo gauge transformations generated by $u$. The pair $(u,l)$ uniquely determines the parabolic pair $(u,[l]_u)$ associated to $\varepsilon$ that we use to characterize the stabilizer of the latter.

\begin{lemma}
Let $T\in\SO_o(V^{\ast},h^{\ast})$ and $\varepsilon\in \Sigma$ with associated spinorial polyform $\alpha_{\varepsilon}\in \wedge V^{\ast}$. Then, $T(\alpha_{\varepsilon}) = \alpha_{\varepsilon}$ if and only if $T(u) = u$ and $T(l) = l + c u$ for a certain constant $c\in\mathbb{R}$, that is, if and only if $T$ preserves the parabolic pair associated to $\varepsilon\in \Sigma$.
\end{lemma}

\begin{remark}
The natural action of the special orthogonal group $\SO_o(V^{\ast},h^{\ast})$ on any element $\beta\in \wedge V^{\ast}$ is denoted simply by the symbol $T(\beta)$ for every $T\in\SO_o(V^{\ast},h^{\ast})$
\end{remark}

\begin{proof}
The result follows from the standard formula:
\begin{eqnarray*}
T(\alpha_{\varepsilon}) = T(u + u\wedge l) = T(u) + T(u)\wedge T(l)
\end{eqnarray*}

\noindent
Imposing $T(\alpha_{\varepsilon}) = \alpha_{\varepsilon}$ gives $T(u) = u$ as well as $u\wedge (l-T(l)) = 0$ and hence we conclude.
\end{proof}

\noindent
Hence, we can study the stabilizer of a spinorial polyform in terms of the stabilizer of the corresponding parabolic pair. In particular:
\begin{equation*}
\cI_{\alpha_{\varepsilon}} = \Aut_{\mathfrak{P}(V)}(h,u,[l]_u)
\end{equation*}

\noindent
that is, $\cI_{\alpha_{\varepsilon}}$ is the automorphim group of $(h,u,[l]_u)$ in the category $\mathfrak{P}(V)$. This allows for an alternative and immediate characterization of $\cI_{\varepsilon}$ using the theory of parabolic pairs developed in the previous subsection. Since Proposition \ref{prop:categoricalequivalence} is an equivalence of categories, $T\in \Aut_{\mathfrak{P}(V)}(h,u,[l]_u)$ if and only if $T\in \Aut_{\mathfrak{F}(V)}([u,v,l,n])$, where $\mathbb{E}([u,v,l,n]) =(h,u,[l]_u) $. However:
\begin{eqnarray*}
T\in \Aut_{\mathfrak{F}(V)}([u,v,l,n]) \quad \Leftrightarrow \quad T([u,v,n,l]) = [u,v,l,n]
\end{eqnarray*}

\noindent
and therefore there exists a unique $w \in ( \langle \mathbb{R}\, u \rangle  \oplus \langle \mathbb{R}\, v\rangle)^{\perp_{h}}$ such that:
\begin{equation*}
\label{eq:stabilizerframe}
	(T^{\ast}u,T^{\ast} v ,T^{\ast} l, T^{\ast} n) = (u,v-\frac{1}{2} \vert w\vert_h^2 u + w,l - h^{\ast}(w,l)u,n - h^{\ast}(w,n) u)
\end{equation*}

\noindent
and therefore we obtain a natural group isomorphism $\mathbb{R}^2 = \Aut_{\mathfrak{P}(V)}(h,u,[l]_u)$. Elaborating on this result, we obtain the following characterization of the stabilizer of  $\cI_{\varepsilon}$.
 
\begin{prop}
The stabilizer $\cI_{\varepsilon}$ of an irreducible real spinor $\varepsilon \in \Sigma \backslash \left\{ 0 \right\}$ is isomorphic to $\mathbb{R}^2$ and projects via the spin double cover map to the following subgroup of $\SO_o(3,1)$:
\begin{equation*}
	\cI_o =   \left( \begin{array}{cccc} 
		1 & -\frac{1}{2} (c_1^2 + c_2^2)  & - c_1 & -c_2\\ 
		0  & 1 & 0 & 0 \\
		0  & c_1 & 1 & 0 \\
		0 & c_2 & 0 & 1\end{array} \right)  
\end{equation*}

\noindent
where $c_1 , c_2 \in \mathbb{R}$.
\end{prop}

\begin{proof}
By the previous discussion, the stabilizer of a parabolic pair is isomorphic to $\mathbb{R}^2$. Since $\mathbb{R}^2$ is connected and simply connected, $\Spin_o (V^{\ast},h^{\ast})$ is connected and the spinor square map is equivariant with respect to the double cover $\Spin_o (V^{\ast},h^{\ast}) \to \SO_o (V^{\ast},h^{\ast})$ we conclude that the connected component of the stabilizer of an irreducible real spinor in $\Spin_o (V^{\ast},h^{\ast})$ is also isomorphic to $\mathbb{R}^2$. Its image in $\SO_o (3,1)$ is computed from Equation \eqref{eq:stabilizerframe} after identifying $\SO_o(3,1) = \SO_o (V^{\ast},h^{\ast})$ using the given basis $(u,v,l,n)$.
\end{proof}
 
\noindent
This result is of course well-known in the literature. Here we have highlighted how this result follows easily and naturally from the theory of parabolic pairs.

 
\section{Differential spinors}
\label{sec:Differentialspinors4d}
 

Let $(M,g)$ be a strongly spin Lorentzian four-manifold of \emph{mostly plus} signature with time-like orientation $\frt\in\mathfrak{X}(M)$. By \emph{strongly spin} we mean that $(M,g)$ admits spin structures that reduce to the identity component $\Spin_o(3,1) \subset \Spin(3,1)$ of the spin group $\Spin(3,1)$. This condition guarantees the existence of admissible bilinear pairings on the irreducible real spinor bundles associated to such spin structures. In this section we extend the algebraic notion of parabolic pair introduced in the previous section to $(M,g)$ and we apply it to the study of differential spinors on $(M,g)$. For ease of notation we will always consider vector fields to be identified with one-forms by means of $g$. 

\begin{definition}
A one-form on $(M,g)$ is \emph{isotropic} if it is nowhere vanishing and $g(u,u)=0$.
\end{definition}

\noindent
For every isotropic one-form $u\in \Omega^1(M)$ we define the following equivalence relation $\sim_u$ on $\Omega^1(M)$:
\begin{equation*}
l_1 \sim_u l_2 \quad \mathrm{if\,and\,only\,if}\quad l_1 = l_2 + f u\, , \quad f\in C^{\infty}(M)\, .
\end{equation*}

\noindent
Transformations of the form $\Omega^1(M)\ni l \mapsto l + f u\in \Omega^1(M)$ will be again referred to as \emph{gauge transformations} generated by $u$. We will denote by $[l]_u$ the equivalence class determined by an element $l\in \Omega^1(M)$.

\begin{definition}
\label{def:parabolic}
A {\em parabolic pair} on $(M,g)$ is an element $(u,[l]_u) \in \Omega^1(M)\times (\Omega^1(M)/\sim_u)$ consisting of an isotropic one-form $u\in \Omega^1(M)$ and an equivalence class $[l]_u$ satisfying:
\begin{equation}  
\label{eq:ulcondM}
h^\ast (l,l) = 1\, ,\qquad  h^{\ast}(u,l) = 0
\end{equation}
	
\noindent
for any, and hence for all, representative $l\in [l]_u$. 
\end{definition}

\noindent
The point behind the previous definition is of course the following characterization of non-vanishing spinors on $(M,g)$, which follows directly from Theorem \ref{thm:squarespinorMink}.

\begin{prop}
A strongly spin Lorentzian four-manifold $(M,g)$ admits a nowhere vanishing irreducible real spinor if and only if it admits a parabolic pair.
\end{prop}

\noindent 
Note that a parabolic pair $(u,[l]_u)$ determines the corresponding spinor $\varepsilon$ uniquely modulo a global sign. Let $\varepsilon$ be a differential spinor on $(M,g)$ with associated parabolic pair $(u,[l]_u)$. Since $u$ is isotropic, we can consider its associated \emph{screen bundle} $\frS_u$, defined as the following vector bundle quotient: 
\begin{equation*}
\frS_u = u^{\perp_g}/ \langle \mathbb{R}\, u\rangle
\end{equation*}

\noindent
where $u^{\perp_g}\subset TM$ is the regular distribution orthogonal to the line bundle $\langle \mathbb{R} \, u\rangle\subset T^{\ast}M$ defined as the distribution spanned by $u$. Hence, $\frS_u$ is a rank-two vector bundle over $(M,g)$ fitting into the following short exact sequence of vector bundles: 
\begin{eqnarray}
0 \to \langle \mathbb{R} u\rangle \to u^{\perp_g}  \to  \frS_u \to 0
\end{eqnarray}

\noindent
which gives the global analog of \eqref{eq:shortexactu}. We will refer to $\frS_u$ as the screen bundle associated to either the differential spinor $\varepsilon$ or its parabolic pair $(u,[l]_u)$. Since $u$ is isotropic, the Lorentzian metric descends to $\frS_u$, where it defines a Riemannian metric $q_u$ as follows:
\begin{equation*}
q_u([v_1]_u,[v_2]_u) = g(v_1,v_2)
\end{equation*}

\noindent
where $[v_1]_u , [v_2]_u \in \frS_u$ denote the equivalence classes determined by $v_1, v_2 \in u^{\perp_g}$. In contrast to the screen bundle associated to a general isotropic vector, the screen bundle associated to a differential spinor is always trivial, as the following result shows.  

\begin{prop}
The screen bundle associated to the parabolic pair $(u,[l]_u)$ is canonically trivialized by $[l]_u$ and its Hodge dual $\ast_{q_u}[l]_u$ with respect to $q_u$. 
\end{prop}

\begin{proof}
Since $(M,g)$ is oriented and time-oriented by assumption, $\frS_u$ inherits a canonical orientation and hence a Hodge operator associated to the metric $q_u$. The equivalence class $[l]_u$ is a section of $\frS_u$ and its Hodge dual with respect to $q_u$ and the canonical orientation determines an orthogonal section $[n]_u \in \frS_u$. The pair $([l]_u,[n]_u)$ defines a canonical trivialization of $\frS_u$.
\end{proof}

\noindent
Given a isotropic one-form $u\in \Omega^1(M)$ on $(M,g)$, we say that a one-form $v\in \Omega^1(M)$ is conjugate to $u$ if $v$ is also isotropic and $g(u,v) = 1$. The fact that $(M,g)$ is oriented and time-oriented implies, via decomposition of $u$ and a direct construction of $v$, that every isotropic one-form on $(M,g)$ admits conjugates. Furthermore, Proposition \ref{prop:conjugatevariation} immediately implies that if $v, v^{\prime} \in \Omega^1(M)$ are both conjugate to $u$, then:
\begin{equation*}
v^{\prime} = v - \frac{1}{2} \vert \frw\vert_g^2 u + \frw
\end{equation*} 

\noindent
for a unique one-form $\frw\in\Omega^1(M)$ orthogonal to both $u$ and $v$. Proposition \ref{prop:adaptedbasis} gives the following analog result globally on $(M,g)$.

\begin{prop}
Let $(u,[l]_u)$ be a parabolic pair on $(M,g)$. Every choice of isotropic one-form $v\in \Omega(M)$ conjugate to $(u,[l]_u)$ canonically determines a positively oriented coframe $(u,v,l,n)$ on $M$, where $l\in [l]_u$ and $n\in [n]_u$ are unit vectors uniquely determined by the condition $l,n\in ( \langle \mathbb{R}\, u\rangle  \oplus \langle \mathbb{R}\, v \rangle)^{\perp_h}$.
\end{prop}

\noindent
We will refer to the global coframes $(u,v,l,n)$ occurring in the previous proposition as \emph{isotropic coframes} associated to $(u,[l]_u)$. Similarly to the previous result, Proposition \ref{prop:primerparallelization} gives the following analog result globally on $(M,g)$.

\begin{prop}
\label{prop:primerparallelizationglobal}
Let $(u,[l]_u)$ be a parabolic pair on $(M,g)$. Given isotropic coframes $(u,v,l,n)$ and $(u,v^{\prime},l^{\prime},n^{\prime})$ associated to $(u,[l]_u)$, there exists a unique vector field $\frw\in \Gamma( \langle \mathbb{R}\, u\rangle  \oplus \langle \mathbb{R}\, v \rangle)^{\perp_h}$ such that:
\begin{equation}
\label{eq:equivalencerelation}
(u,v^{\prime},l^{\prime},n^{\prime}) = (u,v-\frac{1}{2} \vert\frw\vert_g^2 u + \frw,l - \frw(l) u,n - \frw(n) u)
\end{equation}
	
\noindent
In particular, the set of isotropic coframes associated to $(u,[l]_u)$ is a torsor over the infinite-dimensional vector space $\Gamma( \langle \mathbb{R}\, u\rangle  \oplus \langle \mathbb{R}\, v \rangle)^{\perp_h}$.
\end{prop}

\begin{remark}
\label{remark:topologicaltriviality}
Any isotropic coframe $(u,v,l,n)$ defines a smooth trivialization of $TM$. In particular, every strongly spin Lorentzian four-manifold admitting a nowhere vanishing irreducible real spinor is parallelizable. Indeed, since the stabilizer of a nowhere vanishing spinor is the non-compact Lie group $\mathbb{R}^2\subset \Spin_o(3,1)$ the frame bundle reduces to the identity and consequently a strongly spin Lorentzian four-manifold admits a nowhere vanishing spinor field if and only if it is parallelizable. This will pose no topological obstruction for the study of differential spinors on globally hyperbolic Lorentzian four-manifolds in later chapters, since every orientable three-manifold is parallelizable. 
\end{remark}

\noindent
A choice of one-form $v\in\Omega^1(M)$ conjugate to a given parabolic pair $(u,[l]_u)$ defines a \emph{realization} of the abstract screen bundle $\mathfrak{G}_u$ as the codimension-two distribution in $T^{\ast}M$ spanned by $l$ and $n$. That is:
\begin{eqnarray*}
\mathfrak{G}_u \xrightarrow{\simeq} ( \langle \mathbb{R}\, u\rangle  \oplus \langle \mathbb{R}\, v \rangle)^{\perp_g} 
\end{eqnarray*}

\noindent
Given an isotropic coframe on $(u,v,l,n)$ on $M$ we will denote by $q$ the Euclidean metric defined by $g$ on the distribution $( \langle \mathbb{R}\, u\rangle  \oplus \langle \mathbb{R}\, v \rangle)^{\perp_g}$ and by $\ast_q$ the corresponding Hodge dual associated to $q$ and the induced natural orientation. 

We denote by $\mathfrak{P}(M)$ the category whose objects are tuples $(g,u,[l]_u)$ consisting of a Lorentzian metric $g$ on $M$ and a parabolic pair $(u,[l]_u)$ on $(M,g)$, and whose morphisms are invertible isometries that preserve the corresponding parabolic pairs, in complete analogy to the category $\mathfrak{P}(V)$ introduced in the previous section. In addition, and similarly to the category $\mathfrak{F}(V)$ defined in the previous section, we also introduce the category $\mathfrak{F}(M)$ whose objects are equivalence classes $[u,v,l,n]$ of oriented coframes $(u,v,l,n)$ on $M$ with respect to the same pointwise equivalence relation as in the definition of $\mathfrak{F}(V)$. Given their importance in this dissertation, we introduce the following terminology to refer to the objects of $\mathfrak{F}(M)$.

\begin{definition}
An \emph{isotropic parallelism} is an element $[u,v,l,n]\in \Ob(\mathfrak{F}(M))$, that is, it is an equivalence class of isotropic coframes on $M$ with respect to the equivalence relation defined by Equation \eqref{eq:equivalencerelation}.
\end{definition}

\noindent
The equivalence of categories of Proposition \ref{prop:categoricalequivalence} extends naturally to $\mathfrak{P}(M)$ and $\mathfrak{F}(M)$ and therefore we have a fully faithful essentially surjective functor:
\begin{equation*}
	\mathbb{E}\colon \mathfrak{F}(M) \to \mathfrak{P}(M)
\end{equation*}

\noindent
defined in complete analogy to Proposition \ref{prop:categoricalequivalence}. The main goal of this section is to characterize differential spinors relative to  $\cA\in \Omega^1(M,\End(S))$ and $\cQ\in \Gamma(\End(S)))$ with symbol given by a fixed pair $(\fra \in \Omega(M , \wedge M),\mathfrak{q}\in \Gamma(\wedge M))$ in terms of solutions to a differential system for the associated parabolic pairs or isotropic parallelisms. As remarked earlier, we are not interested \emph{per se} on the specific form of the differential spinor $\varepsilon \in \Gamma(S)$ but on the precise geometric and topological consequences of its existence. For simplicity in the exposition, we will sometimes refer to the parabolic pair associated to a differential spinor as a \emph{differential parabolic pair}.

\begin{prop}
\label{prop:differentialspinors4d}
A strongly spin Lorentzian four-manifold $(M,g)$ admits a differential spinor relative to $\fra \in \Omega^1(M,\wedge M)$ if and only if it admits a parabolic pair $(u,[l]_u)$ satisfying the following differential system:
\begin{eqnarray}
\label{eq:generaldiffpoly1}
& \frac{1}{2} \nabla^g_w u =  \fra_w^{(0)} u -  u \wedge l \lrcorner_g \fra_w^{(1)} +  l \wedge u \lrcorner_g \fra_w^{(1)} -  u \lrcorner_g \fra_w^{(2)}+  u\lrcorner_g l \lrcorner_g \fra_w^{(3)}\\
& \frac{1}{2} u\wedge (\nabla^g_w l + 2 \mathrm{P}^{\perp_g}_l(\fra_w^{(1)}) + 2 l\lrcorner_g \fra_w^{(2)}) = l\wedge u\lrcorner_g l \lrcorner_g \fra^{(3)}_w + u \lrcorner_g\fra^{(3)}_w + u\lrcorner_g l\lrcorner_g \fra^{(4)}_w\label{eq:generaldiffpoly2}
\end{eqnarray}
	
\noindent
for every vector field $w\in\mathfrak{X}(M)$ and for any, and hence for all, representatives $(u,l)\in (u,[l]_u)$, where $\mathrm{P}^{\perp_g}_l\colon T^{\ast}M \to T^{\ast}M$ denotes projection to the orthogonal complement of $l$ in the cotangent bundle $T^{\ast}M$.  
\end{prop}

\begin{proof}
By Theorem \ref{thm:GCKS}, the statement is equivalent to $(M,g)$ admitting a parabolic pair $(u,[l]_u)$ satisfying:
\begin{equation}
\label{eq:iffdifferentialspinor4d}
\nabla^g_w (u + u \wedge l) = \fra_w\diamond_g (u + u \wedge l) + (u + u \wedge l) \diamond_g (\pi\circ\tau)(\fra_w)
\end{equation}
	
\noindent
for every vector field $w\in\mathfrak{X}(M)$ and for any representative $(u,l)\in (u,[l]_u)$. Writing:
\begin{eqnarray}
\label{eq:expansioncA}
\fra = \sum_{k=0}^{4} \fra^{(k)}\, , \qquad \fra^{(k)}\in \Omega^1(M,\wedge^k M)
\end{eqnarray}
	
\noindent
we have:
\begin{equation*}
(\pi\circ\tau)(\fra) = \fra^{(0)} - \fra^{(1)} - \fra^{(2)} + \fra^{(3)} + \fra^{(4)}\, .
\end{equation*}

\noindent
We compute:
\begin{eqnarray*}
& \fra\diamond_g u = u\diamond_g \pi(\fra) - 2\, u \lrcorner_g \pi(\fra)\\
& \fra\diamond_g (u \wedge l) = \fra\diamond_g u \diamond_g l = u \diamond_g l \diamond_g \fra - 2\, u\wedge l \lrcorner_g \fra + 2\, l \wedge u \lrcorner_g \fra
\end{eqnarray*}
	
\noindent
Hence:
\begin{eqnarray*}
& \nabla^g_w (u + u \wedge l) = 2\, u \diamond_g (\fra^{(0)}_w - \fra^{(1)}_w + \fra^{(4)}_w) + 2\,(u \wedge l)\diamond_g (\fra^{(0)}_w + \fra^{(3)}_w + \fra^{(4)}_w)\\
& - 2 u \lrcorner_g \pi(\fra_w) - 2 u\wedge l \lrcorner_g \fra_w + 2 l \wedge u \lrcorner_g\fra_w
\end{eqnarray*}
	
\noindent
for every vector field $w\in \mathfrak{X}(M)$. Isolating terms in the previous equation by degree, we obtain non-trivial conditions only in first and second degree:
\begin{eqnarray*}
& \frac{1}{2} \nabla^g_w u =  u \,\fra_w^{(0)}  -  u \wedge l \lrcorner_g \fra_w^{(1)} +  l \wedge u \lrcorner_g \fra_w^{(1)} -  u \lrcorner_g \fra_w^{(2)}+  u\lrcorner_g l \lrcorner_g \fra_w^{(3)} \\
&  \frac{1}{2}\nabla^g_w (u\wedge l) =  u\wedge l \,\fra^{(0)}_w -  u \wedge \fra^{(1)}_w   +  l\wedge u \lrcorner_g \fra^{(2)}_w - u\wedge l \lrcorner_g \fra^{(2)}_w  +  u \lrcorner_g \fra^{(3)}_w +  u\lrcorner_g l \lrcorner_g \fra^{(4)}_w  
\end{eqnarray*}
	
\noindent
The first equation above gives the first equation in the statement of the proposition. Together with the second equation, we obtain:
\begin{equation*}
\frac{1}{2} u\wedge (\nabla^g_w l + 2 \mathrm{P}^{\perp_g}_l(\fra_w^{(1)}) + 2 l\lrcorner_g \fra_w^{(2)}) = l\wedge u\lrcorner_g l \lrcorner_g \fra^{(3)}_w + u \lrcorner_g\fra^{(3)}_w + u\lrcorner_g l\lrcorner_g \fra^{(4)}_w
\end{equation*}
	
\noindent
where $\mathrm{P}^{\perp_g}_l(\fra_w^{(1)})$ denotes the projection of $\fra_w^{(1)}$ to the orthogonal complement of $l$ in $T^{\ast}M$.
\end{proof}

\noindent
We can reformulate Proposition \ref{prop:differentialspinors4d} in simpler terms as follows. 

\begin{cor}
\label{cor:differentialspinors4d}
A strongly spin Lorentzian four-manifold $(M,g)$ admits a differential spinor if and only if it admits a parabolic pair $(u,[l]_u)$ satisfying Equations \eqref{eq:generaldiffpoly1} and \eqref{eq:generaldiffpoly2} for some one-form valued polyform  $\fra\in \Omega^1(M,\wedge M)$.
\end{cor}

\noindent
Let $(S,\Gamma,\cB)$ be a paired spinor bundle and let $\cA\in \Omega^1(M,\End(S))$ with symbol $\fra\in \Omega^1(M,\wedge M)$. In the previous proposition we have expanded $\cA$ as  prescribed in Equation \eqref{eq:expansioncA}. Exploiting the fact that $M$ is oriented and four-dimensional, it is sometimes convenient to instead expand $\fra$ as follows:
\begin{equation*}
\fra = \fra^{(0)} + \fra^{(1)} + \fra^{(2)} + \ast_g \mathfrak{c}^{(1)} +  \mathfrak{c}^{(0)}\otimes \nu_g
\end{equation*}

\noindent
where $\mathfrak{c}^{(0)} \in \Omega^1(M)$, $\mathfrak{c}^{(1)}\in \Omega^1(M,\wedge^1 M)$ are respectively uniquely determined by $\fra^{(4)}$ and $\fra^{(3)}$, and the Hodge dual acts only on the $\wedge^1 M$ factor of $\mathfrak{c}^{(1)}$.  
	
\begin{thm}
\label{thm:differentialspinors4d} 
A strongly spin Lorentzian four-manifold $(M,g)$ admits a differential spinor relative to $\fra\in \Omega^1(M,\wedge M)$ if and only if it admits an isotropic parallelism $[u,v,l,n]\in \Ob(\mathfrak{F}(M))$, satisfying the following differential system:
\begin{eqnarray}
& \frac{1}{2}\nabla^g_w u =  (\fra_w^{(0)}  + \mathfrak{c}^{(1)}_w(n) -   \fra_w^{(1)}(l))\, u +    \fra_w^{(1)}(u)\, l - \fra_w^{(2)}(u) - \mathfrak{c}^{(1)}_w(u)\, n  \label{eq:nablau}\\
&\frac{1}{2} \nabla^g_w v =   (  \fra_w^{(1)}(l) +   \fra_w^{(2)}(u,v) - \fra_w^{(0)} - \frc_w^{(1)}(n))\, v  \nonumber \\
& - (\fra^{(2)}_w(v,n) + \rho(w) ) \, n + (\fra^{(1)}_w(v) + \fra^{(2)}_w(l,v) - \kappa(w)) \, l  \label{eq:nablav}\\
& \frac{1}{2}\nabla^g_w l \label{eq:nablal} =  \kappa(w) u - \mathfrak{a}^{(1)}_w(u)v- \mathfrak{a}^{(1)}_w(v)u -   \fra_w^{(2)}(l)   - (\mathfrak{a}^{(1)}_w(n) + \mathfrak{c}^{(1)}_w(l) + \mathfrak{c}^{(0)}_w)\, n  \\
& \frac{1}{2}\nabla^g_w n = \rho(w) u + (\mathfrak{a}^{(1)}_w(n) + \mathfrak{c}^{(1)}_w(l) + \mathfrak{c}^{(0)}_w)\, l + \mathfrak{c}^{(1)}_w(u) \, v -  \fra^{(2)}_w(n) \label{eq:nablan} 
\end{eqnarray}

\noindent
for any, and hence all representatives $(u,v,l,n)\in [u,v,l,n]$, and for any given one-forms $\kappa , \rho \in \Omega^1(M)$.
\end{thm}

\begin{proof}
By Proposition \ref{prop:differentialspinors4d}, $(M,g)$ admits a differential spinor $\varepsilon$ relative to $\fra$ if and only the parabolic pair associated to $\varepsilon$ satisfies equations \eqref{eq:generaldiffpoly1} and \eqref{eq:generaldiffpoly2}. Assume then that $(M,g)$ admits a differential spinor whose associated parabolic pair satisfies \eqref{eq:generaldiffpoly1} and \eqref{eq:generaldiffpoly2}. Equation \eqref{eq:nablau} corresponds to Equation \eqref{eq:generaldiffpoly1} in Proposition \ref{prop:differentialspinors4d} after noticing that:
\begin{equation*}
u\lrcorner_g l \lrcorner_g \fra^{(3)}_w  = u\lrcorner_g l \lrcorner_g \ast_g\frc^{(1)}_w =   \ast_g (\frc^{(1)}_w\wedge l\wedge u)  =  \mathfrak{c}^{(1)}_w(n)\, u - \mathfrak{c}^{(1)}_w(u)\, n
\end{equation*} 

\noindent
where $w\in \mathfrak{X}(M)$. To obtain Equation \eqref{eq:nablal}, we consider Equation \eqref{eq:generaldiffpoly2} and we compute:
\begin{eqnarray*}
& l\wedge u\lrcorner_g l \lrcorner_g \fra^{(3)}_w  = \mathfrak{c}^{(1)}_w(n)\, l\wedge u - \mathfrak{c}^{(1)}_w(u)\, l\wedge n\\
& u \lrcorner_g\fra^{(3)}_w = \ast_g (\mathfrak{c}^{(1)}_w \wedge u) = \mathfrak{c}^{(1)}_w(u)\, l\wedge n - \mathfrak{c}^{(1)}_w(l)\, u\wedge n + \mathfrak{c}^{(1)}_w(n)\, u\wedge l\\
& u\lrcorner_g l\lrcorner_g \fra^{(4)}_w = \mathfrak{c}^{(0)}_w \ast_g (l\wedge u) = - \mathfrak{c}^{(0)}_w u \wedge n
\end{eqnarray*}

\noindent
Hence:
\begin{equation*}
l\wedge u\lrcorner_g l \lrcorner_g \fra^{(3)}_w + u \lrcorner_g\fra^{(3)}_w + u\lrcorner_g l\lrcorner_g \fra^{(4)}_w =  - \mathfrak{c}^{(1)}_w(l)\, u\wedge n - \mathfrak{c}^{(0)}_w u \wedge n
\end{equation*}

\noindent
Plugging the previous relations into Equation \eqref{eq:generaldiffpoly2} we obtain:
\begin{eqnarray*}
\frac{1}{2} u\wedge (\nabla^g_w l + 2 \mathrm{P}^{\perp_g}_l(\fra_w^{(1)}) + 2 l\lrcorner_g \fra_w^{(2)}   + 2( \mathfrak{c}^{(1)}_w(l) + \mathfrak{c}^{(0)}_w)\, n ) = 0     
\end{eqnarray*}

\noindent
whence:
\begin{equation*}
\nabla^g_w l + 2 \mathrm{P}^{\perp_g}_l(\fra_w^{(1)}) + 2 l\lrcorner_g \fra_w^{(2)}   + 2( \mathfrak{c}^{(1)}_w(l) + \mathfrak{c}^{(0)}_w)\, n = 2 \kappa(w) u
\end{equation*} 

\noindent
for a one-form $\kappa \in \Omega^1(M)$. This recovers Equation \eqref{eq:nablal}. To compute the covariant derivative of $n$ we apply $\nabla^g_w$ to $\ast_g (u\wedge l) = u\wedge n$, obtaining:
\begin{eqnarray*}
u\wedge   \nabla^g_w n =   \ast_g \nabla^g_w (u\wedge l) +    n \wedge \nabla^g_w u  = \ast_g ( (\nabla^g_w u)\wedge l) + \ast_g (  u \wedge  \nabla^g_wl) +    n \wedge \nabla^g_w u
\end{eqnarray*}

\noindent
Using equations \eqref{eq:nablau} and \eqref{eq:nablal}, we compute:
\begin{eqnarray*}
& \frac{1}{2} n \wedge \nabla^g_w u = (\fra_w^{(0)}      -   \fra_w^{(1)}(l)+ \mathfrak{c}^{(1)}_w(n))\, n \wedge u +    \fra_w^{(1)}(u)\,n \wedge l - n \wedge\fra_w^{(2)}(u)  \\
& \frac{1}{2} (\nabla^g_w u)\wedge l = (\fra_w^{(0)} - \fra_w^{(1)}(l) + \mathfrak{c}^{(1)}_w(n))\, u\wedge l - \fra_w^{(2)}(u)\wedge l - \mathfrak{c}^{(1)}_w(u)\, n\wedge l\\
& \frac{1}{2} u \wedge  \nabla^g_wl =   - u\wedge \mathrm{P}^{\perp_g}_l(\fra_w^{(1)}) -  u\wedge \fra_w^{(2)}(l)   - ( \mathfrak{c}^{(1)}_w(l) + \mathfrak{c}^{(0)}_w)\, u\wedge n 
\end{eqnarray*}

\noindent
From this expressions it follows that:
\begin{eqnarray*}
& \frac{1}{2} \ast_g((\nabla^g_w u)\wedge l) = (\fra_w^{(0)} -   \fra_w^{(1)}(l) + \mathfrak{c}^{(1)}_w(n))\, u\wedge n  - \ast_g(\fra_w^{(2)}(u)\wedge l) + \mathfrak{c}^{(1)}_w(u)\, u\wedge v   \\
& \frac{1}{2} \ast_g(u \wedge  \nabla^g_w l) =   - \ast_g (u\wedge \mathrm{P}^{\perp_g}_l(\fra_w^{(1)})) - \ast_g ( u\wedge \fra_w^{(2)}(l))   + ( \mathfrak{c}^{(1)}_w(l) + \mathfrak{c}^{(0)}_w)\, u\wedge l 
\end{eqnarray*}

\noindent
Adding the corresponding contributions and simplifying we obtain:
\begin{eqnarray*}
\frac{1}{2} u\wedge   \nabla^g_w n = (\mathfrak{a}^{(1)}_w(n) + \mathfrak{c}^{(1)}_w(l) + \mathfrak{c}^{(0)}_w)\, u\wedge l + \mathfrak{c}^{(1)}_w(u) u\wedge v - u\wedge \fra^{(2)}_w(n)
\end{eqnarray*}

\noindent
and therefore:
\begin{eqnarray*}
\frac{1}{2}\nabla^g_w n = \rho(w) u + (\mathfrak{a}^{(1)}_w(n) + \mathfrak{c}^{(1)}_w(l) + \mathfrak{c}^{(0)}_w)\, l + \mathfrak{c}^{(1)}_w(u) \, v -  \fra^{(2)}_w(n)
\end{eqnarray*}

\noindent
for a one-form $\rho\in \Omega^1(M)$. This gives Equation \eqref{eq:nablan}. Finally, to compute the covariant derivative of $v$ we apply $\nabla^g_w$ to $\ast_g (l\wedge n) = u\wedge v$, obtaining:
\begin{eqnarray*}
	u\wedge   \nabla^g_w v =   \ast_g \nabla^g_w (l\wedge n) +    v \wedge \nabla^g_w u  = \ast_g ( (\nabla^g_w l)\wedge n) + \ast_g (  l \wedge  \nabla^g_w n) +    v \wedge \nabla^g_w u
\end{eqnarray*}

\noindent
Using equations \eqref{eq:nablau}, \eqref{eq:nablal} and \eqref{eq:nablan}, we compute:
\begin{eqnarray*}
	& \frac{1}{2} v \wedge \nabla^g_w u = (\fra_w^{(0)} + \frc_w^{(1)}(n)      -   \fra_w^{(1)}(l) -   \fra_w^{(2)}(u,v))\, v \wedge u \\
	& + ( \fra_w^{(1)}(u) + \fra_w^{(2)}(l,u))\,v \wedge l - ( \frc_w^{(1)}(u) + \fra_w^{(2)}(u,n))\,v \wedge n  \\
	& \frac{1}{2} (\nabla^g_w l)\wedge n = (\kappa(w) - \fra^{(1)}_w(v) - \fra^{(2)}_w(l,v) )u\wedge n - (\fra^{(1)}_w(u) + \fra^{(2)}_w(l,u)) v\wedge n   \\
	& \frac{1}{2} l \wedge  \nabla^g_w n =  (\rho(w)   + \fra^{(2)}_w(v,n)) l \wedge u  + (\frc^{(1)}_w(u)  + \fra^{(2)}_w(u,n)) l\wedge v
 \end{eqnarray*}

\noindent
from which it follows that:
\begin{eqnarray*}
	& \frac{1}{2} \ast_g((\nabla^g_w u)\wedge l) =  (\fra^{(1)}_w(v) + \fra^{(2)}_w(l,v) - \kappa(w)) u\wedge l - ( \fra^{(1)}_w(u) + \fra^{(2)}_w(l,u) ) v\wedge l   \\
	& \frac{1}{2} \ast_g(u \wedge  \nabla^g_w l) =   - (\fra^{(2)}_w(v,n) + \rho(w) ) u\wedge n + (\frc^{(1)}_w(u)   + \fra^{(2)}_w(u,n)) v\wedge n
\end{eqnarray*}

\noindent
Summing the corresponding contributions and simplifying:
\begin{eqnarray*}
	& \frac{1}{2} u\wedge   \nabla^g_w v = (\fra_w^{(0)} + \frc_w^{(1)}(n)      -   \fra_w^{(1)}(l) -   \fra_w^{(2)}(u,v))\, v \wedge u \\
	&- (\fra^{(2)}_w(v,n) + \rho(w) ) u\wedge n + (\fra^{(1)}_w(v) + \fra^{(2)}_w(l,v) - \kappa(w)) u\wedge l 
\end{eqnarray*}

\noindent
and therefore there exists a one-form $\theta\in\Omega^1(M)$ such that:
\begin{eqnarray*}
	&   \frac{1}{2} \nabla^g_w v = \theta(w) u +  (  \fra_w^{(1)}(l) +   \fra_w^{(2)}(u,v) - \fra_w^{(0)} - \frc_w^{(1)}(n)) v\\
	&- (\fra^{(2)}_w(v,n) + \rho(w) )  n + (\fra^{(1)}_w(v) + \fra^{(2)}_w(l,v) - \kappa(w))  l 
\end{eqnarray*}
\noindent
for a one-form $\theta\in \Omega^1(M)$. On the other hand, taking the covariant derivative of $g(v,v) = 0$ we have:
\begin{eqnarray*}
0 = \frac{1}{2}\nabla^g_w g(v,v) =  g( \nabla^g_w v,v) = 2 \theta(w) = 0 
\end{eqnarray*}

\noindent
for every $w\in\mathfrak{X}(M)$. Hence, we obtain Equation \eqref{eq:nablav}. The converse follows from Proposition \ref{prop:differentialspinors4d} together with the observation that $\nabla^g$ as prescribed in Equations \eqref{eq:nablav} has vanishing torsion and preserves:
\begin{equation*}
g = u\odot v + l\otimes l + n\otimes n
\end{equation*}

\noindent
and thus is the Levi-Civita connection of this metric. 
\end{proof}

\begin{remark}
\label{remark:constraintequation4d}
If, in the situation of the previous theorem, we impose the condition $\cQ(\varepsilon) = 0$ for a given endomorphism $\cQ \in \Gamma(\End(S))$, then by Lemma \ref{lemma:constrainedspinor} we have to supplement the differential system \eqref{eq:nablau}, \eqref{eq:nablav}, \eqref{eq:nablal} and \eqref{eq:nablan} with the algebraic equation:
\begin{equation*}
(u+u\wedge l) \diamond_g \frq = u\wedge \frq + \frq(u) + u\wedge l \wedge \frq - l \wedge \frq(u) + u \wedge \frq(l) + \frq(l,u) =  0
\end{equation*}

\noindent
where $(u,[l]_u)$ is the parabolic pair associated to $\varepsilon$ and $\frq\in \Omega^{\bullet}(M)$ is the symbol of $\cQ$. 
\end{remark}

\begin{remark}
We will refer to $g = u\odot v + l\otimes l + n\otimes n$ as the \emph{associated metric} of $(u,v,l,n)$. Note that if we consider the isotropic parallelism $[u,v,l,n]\in\mathfrak{F}(M)$ determined by $(u,v,l,n)$, then the associated metric does not depend on the representative in $[u,v,l,n]$. Therefore, we can talk about $g$ as the metric associated to the isotropic parallelism $[u,v,l,n]$, and we obtain a natural functor:
\begin{equation*}
\mathfrak{F}(M) \to \mathrm{Lor}(M) \, ,\qquad [u,v,l,n] \mapsto g = u\odot v + l\otimes l + n\otimes n
\end{equation*}
	
\noindent
Here $\mathrm{Lor}(M)$ is understood as the groupoid of Lorentzian metrics and invertible isometries on $M$.
\end{remark}

\noindent
We can think of Theorem \ref{thm:differentialspinors4d} as giving the \emph{master formula} for the existence of differential spinors on Lorentzian four-manifolds. Despite its full generality, it can already be used to obtain several direct corollaries, of which we mention a pair which are particularly transparent.

\begin{cor}
Let $\varepsilon$ be a differential spinor relative to $\fra\in\Omega^1(M,\wedge M)$ and with associated Dirac current $u\in \Omega^{1}(M)$. Then, the Levi-Civita connection $\nabla^g$ preserves the distribution $u^{\perp_g}\subset TM$ if and only if:
\begin{eqnarray*}
\fra^{(2)}_w (u) \vert_{u^{\perp_g}}= 0  
\end{eqnarray*}

\noindent
for every $w\in\mathfrak{X}(M)$. 
\end{cor}

\begin{cor}
Let $\varepsilon$ be a differential spinor relative to $\fra\in\Omega^1(M,\wedge M)$ with associated Dirac current $u\in \Omega^{1}(M)$. Then, the Levi-Civita connection $\nabla^g$ descends to the rank-two vector bundle $\mathfrak{G}_u$ if and only if:
\begin{equation*}
\fra^{(1)}_w (u) = 0\, , \quad\frc^{(1)}_w (u) = 0\, , \quad \fra^{(2)}_w (u) \vert_{u^{\perp_g}}= 0	
\end{equation*}	

\noindent
for every $w\in\mathfrak{X}(M)$.
\end{cor}
 
\noindent
More generally, Theorem \ref{thm:differentialspinors4d}  is especially well-adapted to characterize geometric properties of Lorentzian four-manifolds equipped with differential spinors $\varepsilon$ in terms of conditions satisfied by the data $\fra\in \Omega^1(M,\wedge M)$ relative to which $\varepsilon$ is a differential spinor. 


\section{Natural classes of differential parabolic pairs}
\label{sec:KillingDiracCurrent} 


We have allowed for general endomorphism-valued one-forms $\cA \in \Omega^1(M,\End(S))$ to occur in the definition of differential spinors as well as in Theorem \ref{thm:differentialspinors4d}. Given a paired spinor bundle $(S,\Gamma,\cB)$, it is natural to consider differential spinors relative endomorphism-valued one-forms $\cA \in \Omega^1(M,\End(S))$ that are compatible with $\cB$ in the sense that:
\begin{equation*}
	\cB(\cA_w(\varepsilon_1),\varepsilon_2) + \cB(\varepsilon_1 ,\cA_w(\varepsilon_2)) = 0\, , \qquad \forall\,\, \varepsilon_1, \varepsilon_2 \in \Gamma(S)
\end{equation*}

\noindent
for every $w\in \mathfrak{X}(M)$ and $\varepsilon_1 , \varepsilon_2 \in \Gamma(S)$. Equivalently, $\cA$ is compatible with $\cB$ if and only if the connection $\cD= (\nabla^g - \cA)$ on $S$ preserves $\cB$, that is:
\begin{equation*}
	\cD\cB(\varepsilon_1 ,\varepsilon_2) = \cB(\cD\varepsilon_1,\varepsilon_2) + \cB(\varepsilon_1 ,\cD\varepsilon_2) \, , \qquad \forall\,\, \varepsilon_1, \varepsilon_2 \in \Gamma(S)
\end{equation*}

\noindent
Note that if $\cA$ is compatible with $\cB$ then $\cD$ is a \emph{symplectic connection} on $(S,\omega)$ and therefore has holonomy in the symplectic real group $\Sp(4,\mathbb{R})$.   

\begin{prop}
\label{prop:compatiblekilling}
	Let $(S,\Gamma,\cB)$ be a paired spinor bundle on $(M,g)$. An endomorphism-valued one-form $\cA\in \Omega^1(M,\End(S))$ is compatible with $\cB$ if and only if:
	\begin{equation*}
		\fra = \fra^{(1)}  + \fra^{(2)}   
	\end{equation*}
	
	\noindent
	where $\fra$ denotes the symbol of $\cA$.
\end{prop}

\begin{proof}
$\cA$ is compatible with $\cB$ if and only if $\cA^t_w = -\cA_w$ for every $w\in \mathfrak{X}(M)$. Since $\Psi_{\Gamma}\colon (\wedge M, \diamond_g) \to \End(S)$ is an isomorphism of bundles of unital associative algebras, this condition is equivalent to:
\begin{equation*}
\fra_w = -(\pi\circ\tau)(\fra_w)\, , \qquad \forall\,\, w\in \mathfrak{X}(M)
\end{equation*}
	
\noindent
where we have used that $\Psi^{-1}_{\Gamma}(\cA_v^{t}) = (\pi\circ\tau)(\fra_v)$. Hence, the previous equation is equivalent to:
\begin{equation*}
\fra^{(0)} + \fra^{(1)} + \fra^{(2)} + \fra^{(3)} + \fra^{(4)} = - \fra^{(0)} + \fra^{(1)} + \fra^{(2)} - \fra^{(3)} - \fra^{(4)}
\end{equation*}
	
\noindent
from which the statement follows.  
\end{proof}

\noindent
Interestingly enough, supersymmetric solutions in supergravity may require considering differential spinors relative to data $\cA$ that is not compatible with the given admissible pairing, although we will not encounter such cases in this dissertation. Another natural class of differential spinors $\varepsilon\in\Gamma(S)$, inspired by the type of spinors that occur as supersymmetry parameters in supergravity and string theory, is given by the following two conditions:

\begin{itemize}
\item The Dirac current $u\in \Omega^1(M)$ of $\varepsilon\in\Gamma(S)$ defines a Killing vector field for $g$.
	
\item The distribution $\Ker(u)\subset TM$ determined by the kernel of the Dirac current $u\in \Omega^1(M)$ of $\varepsilon\in\Gamma(S)$ is integrable.
\end{itemize}

\noindent
Indeed, experience shows that these two properties hold for all the known supersymmetric configurations in Lorentzian supergravity that have an isotropic supersymmetry parameter, or equivalently, an isotropic Dirac current \cite{Gran:2018ijr,Ortin}. The fact that the Dirac current associated to the supersymmetry parameters of supersymmetric configurations in supergravity are Killing is one of the reasons why the supersymmetry conditions are sometimes called the \emph{Killing spinor equations} of the corresponding supergravity theory.  
   
\begin{remark}
Differential spinors whose Dirac current is Killing are especially interesting since they give rise to \emph{exact idealized gravitational waves}. These are, by definition, Lorentzian manifolds equipped with an isotropic Killing field that is sometimes required to satisfy further curvature conditions. From this point of view, the study of differential spinors whose associated Dirac current is Killing can be understood as the study of a particular class of \emph{spinorial}
exact gravitational waves.
\end{remark}

\noindent
Let $(S,\Gamma,\cB)$ be a paired spinor bundle on $(M,g)$ and let $\varepsilon\in \Gamma(S)$ be a differential spinor relative to $\cA$ with associated parabolic pair $(u,[l]_u)$. The standard formula for the Lie derivative of a $g$ in terms of $\nabla^g$ immediately implies that $\varepsilon$ is Killing if and only if the symmetrization of Equation \eqref{eq:generaldiffpoly1} vanishes. This yields a complicated condition that explicitly involves the parabolic pair $(u,[l]_u)$ associated to the given differential spinor. Instead of working with this general solution to the Killing condition, we will extract natural geometric conditions that do not depend on the given differential spinor and that at the same time guarantee that the corresponding Dirac current is Killing. A simple verification using Equation \eqref{eq:generaldiffpoly1} proves the following result.

\begin{prop}
\label{prop:KillingcAif}
Let $(S,\Gamma,\cB)$ be a paired spinor bundle on $(M,g)$. The Dirac current of a differential spinor $\varepsilon\in \Gamma(S)$ relative to $\cA\in \Omega^1(M,\End(S))$ is Killing if:
\begin{equation}
\label{eq:KillingcA}
\fra^{(0)}= 0\, , \quad \fra^{(1)}=  \sigma g\, , \quad\fra^{(2)} \in \Omega^3(M)\subset \Omega^1(M,\wedge^2 M)\, , \quad \fra^{(3)} \in \Omega^4(M) \subset \Omega^1(M,\wedge^3 M)
\end{equation}

\noindent
where $\fra$ is the symbol of $\cA$ and $\sigma\in C^{\infty}(M)$ is a function.
\end{prop}

\noindent
By Proposition \ref{prop:KillingcAif}, if $\fra$ satisfies Equation \eqref{eq:KillingcA} then we can write:
\begin{equation*}
\fra =  \sigma g + \ast_g\alpha + (\tau + \beta) \otimes \nu_g
\end{equation*}

\noindent
for uniquely determined functions $\sigma, \tau \in C^{\infty}(M)$ and one-forms $\alpha , \beta \in \Omega^1(M)$. The previous remarks together with Theorem \ref{thm:differentialspinors4d} imply the following characterization of differential spinors relative to data satisfying Equation \eqref{eq:KillingcA}.

\begin{prop}
\label{prop:differentialspinors4dadapted}
A strongly spin Lorentzian four-manifold $(M,g)$ admits a differential spinor relative to $\fra\in \Omega^1(M,\wedge M)$ satisfying Equation \eqref{eq:KillingcA} if and only if it admits a parabolic pair $(u,[l]_u)$ satisfying the following differential system:
\begin{eqnarray}
\label{eq:generaldiffpoly1adapted}
& \frac{1}{2} \nabla^g  u = \sigma\, u\wedge l + \ast_g (\alpha\wedge u) + \tau \ast_{g} (l\wedge u)\\
& \frac{1}{2} u\wedge (\nabla^g_w l + 2 \sigma (w - w(l) l) + 2 \ast_g (\alpha\wedge w \wedge l)  ) = \nonumber\\
& = \tau\, l\wedge \ast_g (w \wedge l \wedge u)  + \tau \ast_g (w\wedge u) + \beta(w) \ast_g (l\wedge u) \label{eq:generaldiffpoly2adapted}
\end{eqnarray}
	
\noindent
for every vector field $w\in\mathfrak{X}(M)$ and for any, and hence for all, representatives $(u,l)\in (u,[l]_u)$.
\end{prop}

\noindent
Proposition \ref{prop:differentialspinors4dadapted} allows us to easily characterize a class of differential spinors whose Dirac current is guaranteed to be Killing and defines an integrable distribution.

\begin{prop}
\label{prop:adaptedKillingintegrable}
The distribution $u^{\sharp_g}\subset TM$ defined by the Dirac current $u\in\Omega^1(M)$ of a differential spinor $\varepsilon$ relative to data satisfying \eqref{eq:KillingcA} is integrable if and only if $\alpha(u) = 0$.
\end{prop}

\begin{proof}
Frobenius theorem implies that the kernel of $u$ is integrable if and only if $u\wedge\dd u = 0$. By equation \eqref{eq:generaldiffpoly1adapted} we have:
\begin{eqnarray*}
\frac{1}{4} \dd  u = \sigma\, u\wedge l + \ast_g (\alpha\wedge u) + \tau \ast_{g} (l\wedge u) 
\end{eqnarray*}

\noindent
Wedging the previous equation with $u$ we obtain:
\begin{eqnarray*}
u\wedge \dd u = u\wedge \ast_g (\alpha\wedge u) = \alpha(u) u\wedge l\wedge n
\end{eqnarray*}

\noindent
and hence we conclude since $u , l , n \in \Omega^1(M)$ are all nowhere vanishing.
\end{proof}

\noindent
 
\begin{remark}
\label{remark:compatiblekilling}
Proposition \ref{prop:KillingcAif} immediately implies that $\cA$ is compatible with $\cB$ and in addition satisfies Equation \eqref{eq:KillingcA}  if and only if $\fra^{(1)} $ is proportional to $g$ and  $\fra^{(2)}$ is totally skew-symmetric.
\end{remark}

\noindent
We collect in the following several well-known types of differential spinors on $(M,g)$, all fitting as particular cases of Theorem \ref{thm:differentialspinors4d}, most of which have been already considered in the literature, especially in the Riemannian case.

\begin{itemize}
\item If $\cA = 0$ vanishes identically then $\varepsilon$ is parallel with respect to the Levi-Civita connection of $g$ and we conclude that $(M,g)$ admits a parallel irreducible real spinor if and only if it admits an isotropic vector field parallel with respect to the Levi-Civita connection as well as a unit one-form $l\in\Omega^1(M)$ orthogonal to $u$ and satisfying $\nabla^g l = \kappa\otimes u$ for a one-form $\kappa \in \Omega^1(M)$.  
	
\item If  $\fra^{(k)} = 0$ for $k = 0,2,3,4$ and $2 \fra^{(1)} = c g$ for a non-zero real constant $c$ then $\varepsilon$ is a standard Killing spinor and from Proposition \ref{prop:adaptedKillingintegrable} we conclude that $(M,g)$ admits an irreducible real Killing spinor if and only if it admits a parabolic pair $(u,[l]_u)$ satisfying:
\begin{equation}
\label{eq:RKSeq0}
\nabla^g u = c\, u\wedge l\, , \qquad \nabla^g l = \kappa\otimes u + c\, (l\otimes l - g)
\end{equation}
	
\noindent
These real Killing spinors are precisely the supersymmetry parameters of the supersymmetric configurations of AdS minimal supergravity in four dimensions. 
	
\item If  $\fra^{(k)} = 0$ for $k = 0,2,3,4$ and $\fra^{(1)} \in \End(TM)$ is an endomorphism of $TM$ symmetric with respect to $g$ then $\varepsilon$ is formally an irreducible Lorentzian generalized Killing spinor, namely a type of Lorentzian analog of the notion of generalized Killing spinor introduced in \cite{MoroianuSemm,MoroianuSemmI,MoroianuSemmII} for Riemannian manifolds.

\item  If  $\fra^{(k)} = 0$ for $k = 0,1,3,4$ and $\fra^{(2)} \in \Omega^3(M)$ is totally skew-symmetric, then $\varepsilon$ is a \emph{skew-torsion parallel spinor}, namely a spinor parallel under a connection with totally skew-symmetric torsion $2\fra^{(2)}$. The remaining chapters of this dissertation will be devoted to the study of this type of parallel spinors.
	
\item If  $\fra^{(k)} = 0$ for $k = 1,2,3,4$ then $\cA\in \Omega^1(M,\End(S))$ takes values in the line spanned by the identity isomorphism, that is, a parallel $\varepsilon\in \Gamma(S)$ relative to such $\cA$ satisfies:
\begin{equation*}
\nabla^g\varepsilon = \theta\otimes\varepsilon 
\end{equation*}
	
\noindent
for a one-form $\theta\in\Omega^1(M)$. Hence, Theorem \ref{thm:differentialspinors4d} implies in this case that $(M,g)$ admits such a differential spinor if and only if it admits a parabolic pair $(u,[l]_u)$ satisfying:
\begin{equation*}
\frac{1}{2} \nabla^g_w u = \theta\otimes u\, , \qquad  \frac{1}{2}  \nabla^g_w l = \rho\otimes u
\end{equation*}
	
\noindent
for a one-form $\rho \in \Omega^1(M)$. In particular, $(M,g)$ admits a recurrent isotropic line, namely a bundle of isotropic lines preserved by the Levi-Civita connection of $(M,g)$. Manifolds equipped with such a recurrent isotropic line, called \emph{Walker manifolds} in \cite{WalkerManifolds} and \emph{weakly abelian} in \cite{MehidiZeghib}, have been already considered in the literature, see for instance \cite{Galaev:2009ie,Galaev:2010jg} and their references and citations.  
\end{itemize}

\noindent
We end this section by considering two examples of differential spinors. We proceed by first considering real Killing spinors, which have been extensively considered in the literature \cite{ContiDalmasso}, to then study parallel spinors relative to data of the form $\fra = \fra^{(4)} = \beta\otimes \nu_h$, which is a type of differential spinor that to the best of our knowledge has not been investigated before. Interestingly enough, these parallel spinors define a particular class of Brinkmann spacetimes. 


\subsection{A class of real Killing spinors} 
\label{subsec:RKS}


In this subsection we consider a class of differential spinors called real Killing spinors on a particular class of Lorentzian four-manifolds that is adapted to the geometric structure that they induce. We follow \cite{RKS}.

\begin{definition}
Let $(M,g)$ be a strongly spin Lorentzian four-manifold $(M,g)$ equipped with a paired irreducible spinor bundle $(S,\Gamma,\cB)$. A \emph{real Killing spinor} on $(M,g)$ is a spinor $\varepsilon\in \Gamma(S)$ satisfying:
\begin{eqnarray*}
\nabla^g_w\varepsilon = \frac{c}{2} \, w\cdot_g \varepsilon\, , \qquad \forall\,\, w\in \mathfrak{X}(M)
\end{eqnarray*}

\noindent
where $c\in \mathbb{R}^{\ast}$ is a non-zero real constant.
\end{definition}

\noindent
By Theorem \ref{thm:GCKS}, a spinor $\varepsilon$ on $(M,g)$ is Killing if and only if its associated parabolic pair $(u,[l]_u)$ satisfies Equation \eqref{eq:RKSeq0}, which is in turn equivalent to: 
\begin{equation}
\label{eq:RKSeq}
\cL_{u^{\sharp_g}}g= 0\, , \quad \dd u = 2 c\, u\wedge l\, , \quad \dd l = \kappa\wedge u\, , \quad \cL_{l^{\sharp_g}} g=\kappa \odot u+2 c\,( l\otimes l-g)
\end{equation}

\noindent
Instead of considering general real Killing spinors we consider in the following a particular class of Lorentzian four-manifolds equipped with a real Killing spinor whose associated parabolic pair is of a special form. 

\begin{definition}
A four-dimensional space-time $(M,g)$ is \emph{standard conformally Brinkmann} if $(M,g)$ has the following isometry type: 
\begin{equation}
\label{eq:isostanads}
(M,g) = (\mathbb{R}^2\times X, \cH_{x_u} \dd x_u\otimes \dd x_u + e^{\cF_{x_u}} \dd x_u\odot \dd x_v + \dd x_u \odot \alpha_{x_u} + q_{x_u})\, ,
\end{equation}
	
\noindent
where $X$ is an oriented two-dimensional manifold, $(x_u,x_v)$ are the Cartesian coordinates of $\mathbb{R}^2$, and where: 
\begin{equation*}
\left\{ \cH_{x_{u}} ,\cF_{x_{u}} \right\}_{x_{u}\in \mathbb{R}}\, , \quad \left\{\alpha_{x_u} \right\}_{x_{u}\in \mathbb{R}}\, , \quad \left\{ q_{x_u}\right\}_{x_{u}\in \mathbb{R}}\, ,
\end{equation*}
	
\noindent
respectively denote a family of pairs of functions $\cH_{x_{u}}$ and $\cF_{x_{u}}$, a family of one-forms $\alpha_{x_u}$ and a family of complete Riemannian metrics $q_{x_u}$ on $X$ parametrized by $x_{u}\in \mathbb{R}$.
\end{definition} 

\begin{definition}
A spinor $\varepsilon$ on a standard conformally Brinkmann space-time $(M,g)$ is \emph{adapted} if the parabolic pair $(u,[l]_u)$ associated to $\varepsilon$ satisfies $u^{\sharp_g} = \partial_{x_v}$, in which case we will refer to $(u,[l]_u)$ as an \emph{adapted parabolic pair} on $(M,g)$. 
\end{definition}

\noindent
By \cite[Lemma 5.10]{Cortes:2019xmk} every Lorentzian four-manifold equipped with an irreducible real Killing spinor is locally isomorphic to a standard conformally Brinkmann space-time equipped with an adapted Killing spinor. 

\begin{definition}
A \emph{standard real Killing spinor triple} is a triple $(M,g,\varepsilon)$ consisting of a standard conformally Brinkmann space-time $(M,g)$ and an adapted real Killing spinor $\varepsilon$ on $(M,g)$.	
\end{definition}

\noindent
In particular, the parabolic pair associated to a standard real Killing spinor is adapted. We will characterize standard real Killing spinors $(M,g,\varepsilon)$ in the following. Let $\varepsilon$ be an adapted spinor with associated Killing pair $(u,[l]_u)$ on a standard conformally Brinkmann space-time $(M,g)$. Since $u^{\sharp_g} = \partial_{x_v}$ by definition, the first equation in \eqref{eq:RKSeq} is automatically satisfied. Using that  $u = e^{\cF_{x_u}} \dd x_u$, the second equation in \eqref{eq:RKSeq} reduces to:
\begin{equation}
(c\, l + \frac{1}{2} \dd_X \cF_{x_u})\wedge u = 0 \, ,
\end{equation}

\noindent
for any representative $l\in [l]$. Here $\dd_X\colon \Omega^{\bullet}(X) \to \Omega^{\bullet}(X)$ denotes the exterior derivative on $X$. The general solution to this equation reads:
\begin{equation*}
	l=\ell_{x_u}+\sigma_l u\, , \qquad \ell_{x_u}=-\frac{1}{2 c} \dd_X \cF_{x_u}\, ,
\end{equation*}

\noindent
for a function $\sigma_l\in C^\infty(M)$ and a family of functions $\cF_{x_u}$ in $X$ parametrized by $x_u$. Since $l$ is defined modulo gauge transformations generated by $u$, we assume, without loss of generality, the following expression for $l$:
\begin{equation}
	\label{eq:expressionl}
	l = \ell_{x_u} = -\frac{1}{2  c} \dd_X \cF_{x_u}\, .
\end{equation}

\noindent
With this choice, the third equation in \eqref{eq:RKSeq} reduces to:
\begin{equation*}
\dd l=\dd x_u \wedge \partial_{x_u} \ell_{x_u} + \dd_X \ell_{x_u} = \dd x_u \wedge \partial_{x_u} \ell_{x_u} = e^{\cF_{x_{u}}}\kappa \wedge \dd x_u\,.
\end{equation*}

\noindent
where we have used that $\dd_X \ell_{x_u}=0$ by Equation \eqref{eq:expressionl}. The previous equation  is solved by isolating $\kappa$ as follows:
\begin{equation}
	\label{eq:algunaskapicas}
	\kappa = -e^{-\cF_{x_{u}}}\partial_{x_u} \ell_{x_u} + \sigma_\kappa \dd x_u\, ,
\end{equation}

\noindent
for a function $\sigma_\kappa\in C^\infty(M)$. Recalling that we must have $\vert l \vert_g^2=1$, the previous discussion implies the following result.

\begin{prop}
Let $(u,[l]_u)$ be an adapted Killing parabolic pair on a standard conformally Brinkmann space-time. Then, there exists a representative $l\in[l]$ such that $(u,l)$ is of the form:
\begin{equation*}
u = e^{\cF_{x_{u}}} \dd x_u\, , \qquad l = -\frac{1}{2 c} \dd_X \cF_{x_u}
\end{equation*}
	
\noindent
where the family of functions $\left\{\cF_{x_{u}}\right\}_{x_u\in \mathbb{R}}$ satisfies the differential equation:
\begin{equation*}
\vert \dd_X \cF_{x_{u}}\vert_{q_{x_u}}^2 = 4 c^2 
\end{equation*}
	
\noindent
Here $\vert \cdot \vert_{q_{x_u}}^2$ denotes the norm defined by $\left\{q_{x_{u}}\right\}_{x_u\in \mathbb{R}}$ on $X$.
\end{prop}

\noindent
It only remains to solve the fourth equation in \eqref{eq:RKSeq}. For this, we need first the following lemma.

\begin{lemma}
\label{lemma:Lieg}
The Lie derivative of $g$ with respect to  $l^{\sharp_g} = \ell_{x_u}^{\sharp_g}$ is given by:
\begin{eqnarray*}
& \cL_{l^{\sharp_g}} g = ( \dd_X \cH_{x_u}(\ell_{x_u}^{\sharp_q})+2\alpha_{x_u} ( \partial_{x_u} \ell_{x_u}^{\sharp_q}) + 2e^{\cF_{x_{u}}}\partial_{x_u} \frk_{x_u} )  \dd x_u \otimes \dd x_u+ \dd_X e^{\cF_{x_{u}}} (\ell_{x_u}^{\sharp_q}) \dd x_u \odot \dd x_v\\ 
& + e^{\cF_{x_{u}}} \dd_X \frk_{x_u} \odot \dd x_u+ \cL_{\ell_{x_u}^{\sharp_q}}^X \alpha_{x_u} \odot \dd x_u +\cL_{\ell_{x_u}^{\sharp_q}}^X q_{x_u}+\partial_{x_u} \ell_{x_u}  \odot \dd x_u- (\partial_{x_u} q_{x_u}) (\ell_{x_u}^{\sharp_q}) \odot \dd x_u\, , 
\end{eqnarray*}
	
\noindent
where the symbol $\cL^{X}$ denotes the Lie derivative operator on $X$ and $\left\{ \frk_{x_u}\right\}_{x_u\in \mathbb{R}}$ is the family of functions on $X$ determined by:
\begin{equation*}
\frk_{x_u} = \frac{e^{-\cF_{x_{u}}}}{2 c} \dd_X \cF_{x_u}(\alpha^{\sharp_{q}}_{x_u}) 
\end{equation*}
	
\noindent
for every $x_u\in \mathbb{R}$. The symbol $\sharp_q$ denotes the musical isomorphism on $X$ determined by $\left\{ q_{x_u}\right\}_{x_u\in \mathbb{R}}$.
\end{lemma}

\begin{proof}
We first compute the metric dual of $\ell_{x_u}$ with respect to $g$, which is given by:
\begin{equation}
\label{eq:defil}
l^{\sharp_g}=\ell_{x_u}^{\sharp_{q}} - e^{-\cF_{x_{u}}} \alpha_{x_u}(\ell_{x_u}^{\sharp_{q}}) \partial_{x_v} 
\end{equation} 
	
	\noindent 
	Using the previous expression for $l^{\sharp_g}$ compute:
	\begin{eqnarray*}
		& \cL_{l^{\sharp_g}} \dd x_u = 0\, , \qquad \cL_{l^{\sharp_g}} \dd x_v=\dd_X \frk_{x_u}+\partial_{x_u} \frk_{x_u} \dd x_u \, , \qquad \cL_{l^{\sharp_g}} \alpha_{x_u}=\cL_{\ell_{x_u}^{\sharp_q}}^X \alpha_{x_u}+\alpha_{x_u}\left ( \partial_{x_u} \ell_{x_u}^{\sharp_q} \right) \dd x_u \\ 
		&\cL_{l^{\sharp_g}} q_{x_u}  =\cL_{\ell_{x_u}^{\sharp_q}}^X q_{x_u}+q_{x_u}(\partial_{x_u} \ell^{\sharp_{q}}_{x_u})  \odot \dd x_u =\cL_{\ell_{x_u}^{\sharp_q}}^X q_{x_u}+\partial_{x_u} \ell_{x_u}  \odot \dd x_u- (\partial_{x_u} q_{x_u}) (\ell_{x_u}^{\sharp_q}) \odot \dd x_u
	\end{eqnarray*}
	
	\noindent
	from which the Lie derivative of $g$ as given in Equation \eqref{eq:isostanads} follows directly.
\end{proof}

\begin{lemma}
	\label{lemma:iffRKS}
	A standard conformally Brinkmann space-time admits an adapted Killing spinor if and only if the tuple $\left\{\cF_{x_{u}}, \alpha_{x_u}, q_{x_u}\right\}_{x_u\in\mathbb{R}}$ satisfies the following differential system on $X$:
	\begin{eqnarray}
		& \vert \dd_X \cF_{x_{u}}\vert_{q_{x_u}}^2 = 4 c^2\, , \quad \nabla^{q_{x_u}} \dd_X \cF_{x_{u}} + \frac{1}{2}  \dd_X \cF_{x_{u}}\otimes \dd_X \cF_{x_{u}} = 2 c^2 q_{x_u}\label{eq:constraint2dX}\\ 
		& 4 c^2  \dd_X \alpha_{x_u} = \langle \partial_{x_u} \ast_{q_{x_u}}\dd_X \cF_{x_{u}} + \ast_{q_{x_u}}\partial_{x_u}\dd_X\cF_{x_{u}} - 4 c^2 \ast_{q_{x_u}} \alpha_{x_u} , \dd_X \cF_{x_{u}} \rangle_{q_{x_u}} \nu_{q_{x_u}} \label{eq:evolution2dX}
	\end{eqnarray}
	
	\noindent
	where $\langle\cdot,\cdot\rangle_{q_{x_u}}$ denotes the norm defined by $q_{x_u}$ and $\nu_{q_{x_u}}$ denotes the Riemannian volume form of $(X,q_{x_u})$.
\end{lemma}

\begin{proof}
	By the previous discussion we only need to consider the fourth equation in \eqref{eq:RKSeq} evaluated on a parabolic pair $(u,[l])$ of the form:
	\begin{equation}
		u = e^{\cF_{x_{u}}} \dd x_u \, , \qquad l = \ell_{x_u}=  -\frac{1}{2 c} \dd_X \cF_{x_u}\, ,
	\end{equation}
	
	\noindent
	and with respect to a one-form $\kappa \in \Omega^1(\mathbb{R}^2\times X)$ given by:
	\begin{equation*}
		\kappa = -e^{-\cF_{x_{u}}}  \partial_{x_u} \ell_{x_u} + \sigma_\kappa \dd x_u\, ,
	\end{equation*}
	
	\noindent
	where $\left\{\cF_{x_{u}}\right\}_{x_u\in\mathbb{R}}$ satisfies $\vert \dd_X \cF_{x_{u}}\vert_{q_{x_u}}^2 = 4 c^2$. Plugging Lemma \ref{lemma:Lieg} into the fourth equation in \eqref{eq:RKSeq} and isolating we obtain the following system of equivalent equations: 
	\begin{eqnarray}
		\label{eq:kapicau} & 2\sigma_\kappa e^{\cF_{x_{u}}}  = 2 c \cH_{x_u} + \dd_X \cH_{x_{u}}(\ell_{x_u}^{\sharp_q})+ 2 e^{\cF_{x_{u}}}  \partial_{x_u} \frk_{x_u} +  2\alpha_{x_u} ( \partial_{x_u} \ell_{x_u}^{\sharp_q}) \\ 
		\label{eq:lieq} & \dd_X \cF_{x_u} ( \ell_{x_u}^{\sharp_q}) = - 2 c \, , \quad \cL_{\ell_{x_u}^{\sharp_q}}^X q_{x_u} =2 c ( \ell_{x_u} \otimes \ell_{x_u}- q_{x_u})\\
		\label{eq:killings3}& e^{\cF_{x_{u}}}  \dd_X \frk_{x_u}  = (\partial_{x_u} q_{x_u})(\ell_{x_u}^{\sharp_q}) - 2 c \alpha_{x_u} - \cL^X_{\ell_{x_u}^{\sharp_q}}\alpha_{x_u} - 2\partial_{x_u} \ell_{x_u} 
	\end{eqnarray}
	
	\noindent
	Equation \eqref{eq:kapicau} is solved by isolating $\sigma_\kappa$, which determines it unambiguously. The first equation in \eqref{eq:lieq} follows from the first equation in \eqref{eq:constraint2dX} together with equation \eqref{eq:expressionl} whereas the second equation in \eqref{eq:lieq} is equivalent to the second equation in \eqref{eq:constraint2dX} upon use of Equation \eqref{eq:expressionl}. On the other hand, Equation \eqref{eq:killings3} can be shown to be equivalent to:
	\begin{equation}
		2 c \, \alpha_{x_u}+ \iota_{\ell_{x_u}^{\sharp_q}} \dd_X \alpha_{x_u}-\frac{1}{2 c} \dd_X \cF_{x_u} (\alpha_{x_u}^{\sharp_q}) \dd_X \cF_{x_u} + 2 \partial_{x_u} \ell_{x_u}- (\partial_{x_u} q_{x_u})(\ell_{x_u}^{\sharp_q})=0\, .
	\end{equation}
	
	\noindent
	Projecting the previous equation along $\ell_{x_u}^{\sharp_q}$ we obtain an identity. On the other hand, projecting along $(\ast_{q_{x_u}}\ell_{x_u})^{\sharp_q}$ we obtain, after some manipulations, Equation in \eqref{eq:evolution2dX} and hence we conclude.
\end{proof}

\begin{remark}
Recall that equations \eqref{eq:constraint2dX} and \eqref{eq:evolution2dX} do not involve $\cH_{x_{u}}$, which hence can be chosen at will while preserving the existence of an adapted Killing spinor.
\end{remark}

\begin{prop}
\label{prop:standardqxu}
Let $(M,g,\varepsilon)$ be a standard real Killing spinor triple. Then, there exists a diffeomorphism identifying either $X= \mathbb{R}^2$ or $X= \mathbb{R}\times S^1$ and a smooth family of closed one-forms $\omega_{x_u}$ on $X$ such that:
\begin{equation}
\label{eq:standardqxu}
q_{x_u} = \frac{1}{4 c^2} \dd_X \cF_{x_{u}}\otimes \dd_X\cF_{x_{u}} + e^{\cF_{x_u}} \omega_{x_u}\otimes \omega_{x_u}\, .
\end{equation}
	
\noindent
In particular, for every $x_u\in \mathbb{R}$ the Riemann surface $(X , q_{x_u})$ is a hyperbolic Riemann surface of constant scalar curvature $s^{q_{x_u}} = -2 c^2$.
\end{prop}

\begin{proof}
For any fixed $x_u\in \mathbb{R}$, the function $\cF_{x_{u}}$ on $(X , q_{x_u})$ has unit-norm gradient. Hence, since $(X , q_{x_u})$ is complete, $\cF_{x_{u}}(X) = \mathbb{R}$ and we have a diffeomorphism:
\begin{equation*}
X = \mathbb{R}\times \cF_{x_{u}}^{-1}(0)
\end{equation*}
	
\noindent
Since $X$ is assumed to be connected, either $\cF_{x_{u}}^{-1}(0)$ is diffeomorphic to $\mathbb{R}$ or $S^1$. On the other hand, the family of one-forms $\left\{ \ell_{x_u}\right\}_{x_u\in \mathbb{R}}$ has unit norm with respect to $\left\{ q_{x_u}\right\}_{x_u\in \mathbb{R}}$. Therefore, defining:
\begin{equation*}
n_{x_u} = \ast_{q_{x_u}} \ell_{x_u}\, , \qquad x_u\in \mathbb{R}
\end{equation*}
	
\noindent
we obtain a family $\left\{ \ell_{x_u},n_{x_u}\right\}_{x_u\in \mathbb{R}}$ of orthonormal coframes on $X$. The second equation in \eqref{eq:constraint2dX} implies the following equations for $\left\{ \ell_{x_u},n_{x_u}\right\}_{x_u\in \mathbb{R}}$:
\begin{equation*}
\nabla^{q_{x_u}}\ell_{x_u} = - c\, n_{x_u}\otimes n_{x_u}\, , \qquad \nabla^{q_{x_u}} n_{x_u} = c n_{x_u}\otimes \ell_{x_u}
\end{equation*}
	
\noindent
which in turn implies, since $\left\{ \ell_{x_u},n_{x_u}\right\}_{x_u\in \mathbb{R}}$ is orthonormal, that for every fixed $x_u\in\mathbb{R}$, the scalar curvature of the Riemannian metric $q_{x_u}$ satisfies $s^{q_{x_u}} = -2 c^2$ and hence $(X,q_{x_u})$ is a hyperbolic Riemann surface. Furthermore, by the previous equation we have:
\begin{equation*}
\dd_X n_{x_u} = c\, n_{x_u}\wedge \ell_{x_u}
\end{equation*}
	
\noindent
Since $2 c \ell_{x_u} = - \dd_X\cF_{x_u}$, the previous equation is equivalent to:
\begin{equation*}
\dd_X (e^{-\frac{1}{2} \cF_{x_{u}}}n_{x_u}) = 0 
\end{equation*}
	
\noindent
whence there exists a family $\left\{ \omega_{x_u}\right\}_{x_u\in \mathbb{R}}$ of closed one-forms such that:
\begin{equation*}
\omega_{x_u} = e^{-\frac{1}{2} \cF_{x_{u}}}n_{x_u}\, , \qquad \forall\,\, x_u\in \mathbb{R}
\end{equation*}
	
\noindent
Therefore:
\begin{equation*}
q_{x_u} = \ell_{x_u}\otimes \ell_{x_u} + n_{x_u}\otimes n_{x_u} = q_{x_u} = \frac{1}{4c^2} \dd_X \cF_{x_{u}}\otimes \dd_X\cF_{x_{u}} + e^{\cF_{x_u}} \omega_{x_u}\otimes \omega_{x_u}
\end{equation*}
	
\noindent
and hence we conclude.
\end{proof}

\begin{remark}
If $X = \mathbb{R}^2$ then there exists a family of functions $\cG_{x_u}$ such that:
\begin{equation*}
e^{-\frac{1}{2} \cF_{x_{u}}}n_{x_u} = \dd_X \cG_{x_{u}}\, ,
\end{equation*}
	
\noindent
which immediately implies that:
\begin{equation*}
q_{x_u} = \frac{1}{4 c^2} \dd_X \cF_{x_{u}}\otimes \dd_X\cF_{x_{u}} + e^{\cF_{x_u}} \dd_X\cG_{x_u}\otimes \dd_X\cG_{x_u}\, .
\end{equation*}
	
\noindent
On the other hand, if $X=\mathbb{R}\times S^1$ then $H^1(X,\mathbb{Z}) = \mathbb{Z}$ and there exists a family of constants $c_{x_u}$ and a family of functions $\cG_{x_u}$ such that:
	\begin{equation*}
		e^{-\frac{1}{2} \cF_{x_{u}}}n_{x_u} = c_{x_u}\,\omega + \dd_X \cG_{x_{u}}\, , 
	\end{equation*}
	
	\noindent
	where $\omega$ is a volume form on $S^1$. This parametrizes $q_{x_u}$ in terms of families of functions $\cF_{x_{u}}$ and $\cG_{x_{u}}$, and a family of constants $c_{x_{u}}$.
\end{remark}

\noindent
Since $X$ is either diffeomorphic to $\mathbb{R}^2$ or $\mathbb{R}\times S^1$, the uniformization theorem for Riemann surfaces yields the following result.

\begin{cor}
Let $(M,g,\varepsilon)$ be a standard real Killing spinor triple. For every fixed $x_u\in \mathbb{R}$ the pair $(X,q_{x_u})$ is isometric to an elementary hyperbolic surface, namely it is isometric to either the Poincar\'e upper space, to a hyperbolic cylinder or to a parabolic cylinder, in all cases of scalar curvature $-2 c^2$.
\end{cor}

\noindent
More explicitly, for every $x_u\in \mathbb{R}$ the Riemann surface $(X,q_{x_u})$ is locally isometric to the model:
\begin{equation*}
	(\mathbb{R}^2,q = \frac{\dd\rho\otimes \dd\rho}{4c ^2} + e^{\rho} \dd w\otimes \dd w)\, ,
\end{equation*}

\noindent
which in turn is isometric to the Poincar\'e upper space of scalar curvature $-2 c^2$. The hyperbolic cylinder is obtained from $(\mathbb{R}^2,q)$ via quotient by the cyclic group generated by $(\rho,w) \mapsto (\rho - 2 l , e^l w)$, $l>0$. The parabolic cylinder is on the other hand obtained instead via taking the quotient by the cyclic group generated by $(\rho,w) \mapsto (\rho,w+1)$. In both cases, the globally defined one-form $\dd\rho \in \Omega^1(\mathbb{R}^2)$ is invariant and therefore descends to the quotient and defines the constant-norm gradient that Lemma \ref{lemma:iffRKS} requires for an adapted Killing spinor to exist on $\mathbb{R}^2\times X$.

\begin{lemma}
\label{lemma:alphabeta}
Let $(M,g,\varepsilon)$ be a standard real Killing spinor triple. Then the family of one-forms $\alpha_{x_{u}}$ can be written as follows:
\begin{eqnarray}
\alpha_{x_u} = e^{\cF_{x_{u}}} \beta_{x_{u}}\, , 
\end{eqnarray}
	
\noindent
where $\beta_{x_{u}}$ is a family of one-forms satisfying the following differential equation:
\begin{equation}
\label{eq:betaequation}
\dd_X \beta_{x_{u}} = \frac{e^{-\cF_{x_{u}}}}{4c ^2} \langle \partial_{x_u} \ast_{q_{x_u}}\dd_X \cF_{x_{u}} + \ast_{q_{x_u}}\partial_{x_u}\dd_X\cF_{x_{u}} , \dd_X \cF_{x_{u}} \rangle_{q_{x_u}} \nu_{q_{x_u}}\, .
\end{equation}
	
\noindent
Since $H^2(X,\mathbb{Z}) = 0$, this equation has always an infinite-dimensional vector space of solutions.
\end{lemma}

\begin{proof}
Let $(M,g,\varepsilon)$ be a standard real Killing spinor triple. Then, the family of one-forms $\alpha_{x_{u}}$ satisfies Equation \eqref{eq:evolution2dX}. Plugging $\alpha_{x_{u}} = e^{\cF_{x_{u}}} \beta_{x_{u}}$ into Equation \eqref{eq:evolution2dX} and isolating $\beta_{x_u}$ we obtain Equation \eqref{eq:betaequation}. The fact that $H^2(X,\mathbb{Z}) = 0$ follows from Proposition \ref{prop:standardqxu}, which proves the that either $X=\mathbb{R}^2$ or $X=\mathbb{R}\times S^1$.
\end{proof}

\noindent
We can explicitly solve Equation \eqref{eq:betaequation} in terms of local coordinates $(y,z)$, which can be taken to be global if $X=\mathbb{R}^2$. Using these coordinates, the most general solution can be found to be:
\begin{eqnarray}
	\label{eq:standardalphaxu}
	\beta_{x_u} =  \left (\gamma_{x_u}(y)+\int_0^z \left[ \partial_y  \kappa_{x_u}(y,z')-   \Upsilon_{x_u}(y,z') \right ]\, \dd z'  \right )\,\dd y+ \kappa_{x_u}(y,z) \, \dd z   \nonumber
\end{eqnarray}

\noindent
for families of functions $\left\{\kappa_{x_{u}}(y,z),\gamma_{x_{u}}(y)\right\}_{x_u\in \mathbb{R}}$ depending on the indicated variables, and where we have defined:
\begin{equation*}
	\Upsilon_{x_u} := \frac{e^{-\cF_{x_{u}}}}{4c ^2 } \langle \partial_{x_u} \ast_{q_{x_u}} \dd_X \cF_{x_u} + \ast_{q_{x_u}} \partial_{x_u} \dd_X \cF_{x_u}, \dd_X \cF_{x_u} \rangle_{q_{x_u}} q_{x_u}^{\frac{1}{2}}\, ,
\end{equation*}

\noindent
with the understanding that the square root $q^{\frac{1}{2}}_{x_u}$ of the determinant of $q_{x_u}$ is to be taken in the coordinates $(y,z)$. This provides an explicit parametrization of the infinite-dimensional space of local solutions.

\begin{thm}
\label{thm:conformallyBrinkmann}
A triple $(\mathbb{R}^2\times X,g, \varepsilon)$ is a standard real Killing spinor triple if and only if there exist a family $\omega_{x_u}$ of closed one-forms on $X$ in terms of which the metric $g$ reads:    
\begin{eqnarray}
\label{eq:thmcondition1}
g = \cH_{x_u} \dd x_u\otimes \dd x_u +   e^{\cF_{x_u}} \dd x_u\odot ( \dd x_v + \beta_{x_u}) + \frac{1}{4c ^2} \dd_X \cF_{x_{u}}\otimes \dd_X\cF_{x_{u}} + e^{\cF_{x_u}} \omega_{x_u}\otimes \omega_{x_u}
\end{eqnarray}
	
\noindent
where $\beta_{x_u}$ is a family of one-forms on $X$ satisfying the following equation:
\begin{equation}
\label{eq:thmcondition2}
\dd_X \beta_{x_{u}} = \frac{e^{-\cF_{x_{u}}}}{4c ^2} \langle \partial_{x_u} \ast_{q_{x_u}}\dd_X \cF_{x_{u}} + \ast_{q_{x_u}}\partial_{x_u}\dd_X\cF_{x_{u}} , \dd_X \cF_{x_{u}} \rangle_{q_{x_u}} \nu_{q_{x_u}}\, .
\end{equation}
	
\noindent
In particular, for every $x_u\in\mathbb{R}$ the pair:
\begin{equation*}
(X,q_{x_u} = \frac{1}{4c ^2} \dd_X \cF_{x_{u}}\otimes \dd_X\cF_{x_{u}} + e^{\cF_{x_u}} \omega_{x_u}\otimes \omega_{x_u})
\end{equation*}
	
\noindent
is an elementary hyperbolic surface of scalar curvature $-2 c^2$, and therefore diffeomorphic to either $\mathbb{R}^2$ or $\mathbb{R}\times S^1$.
\end{thm} 

\begin{proof}
The \emph{only if} direction follows from Proposition \ref{prop:standardqxu} and Lemma \ref{lemma:alphabeta}. For the \emph{if} direction, consider a standard conformally Brinkmann space-time $(M=\mathbb{R}^2\times X,g)$ whose metric $g$ is as prescribed in Equation \eqref{eq:thmcondition1} for a family of one-forms $\beta_{x_u}$ satisfying Equation \eqref{eq:thmcondition2}. To prove that such Lorentzian four-manifold admits an adapted Killing spinor it is enough to prove that the necessary and sufficient conditions of Lemma \ref{lemma:iffRKS} are satisfied. For each fixed $x_u\in \mathbb{R}$, the functions $\cF_{x_{u}}$ and $\cG_{x_{u}}$ have linearly independent differentials at every point in $X$ (otherwise $g$ would be degenerate) and therefore provide local coordinates $(\rho,w)$ on $X$ in terms of which the metric $q_{x_u}$ reads:
\begin{equation*}
q_{x_u} = \frac{1}{4c ^2} \dd\rho\otimes \dd\rho + e^{\rho} \dd w\otimes \dd w\, ,
\end{equation*} 
	
\noindent
which is isometric to the Poincar\'e metric of Ricci curvature $-c^2$. In these coordinates, we have $\dd_X \cF_{x_{u}} = \dd \rho$ and a quick computation shows that:
\begin{equation*}
\vert\dd\cF_{x_{u}}\vert^2_{q_{x_u}} = 4 c^2 \quad \forall\,\, x_u\in\mathbb{R}
\end{equation*}
	
\noindent
whence the first equation in \eqref{eq:constraint2dX} is satisfied. Furthermore, we have:
\begin{equation*}
\nabla^{q_{x_u}} \dd_X\cF_{x_{u}} = \nabla^{q_{x_u}}\dd\rho = 2c ^2 e^\rho \dd w \otimes \dd w\, ,
\end{equation*}
	
\noindent
whence the second equation in \eqref{eq:constraint2dX} also follows. Finally, Equation \eqref{eq:thmcondition2} is equivalent to Equation \eqref{eq:evolution2dX} by Lemma \ref{lemma:alphabeta} and hence we conclude. 
\end{proof}

\noindent
Using the fact that if $M$ is simply connected then $X=\mathbb{R}^2$ and hence contractible, we obtain the following refinement of Equation \eqref{eq:thmcondition1} in Theorem \ref{thm:conformallyBrinkmann}.

\begin{cor}
Every choice of families of functions $\cF_{x_u}$ and $\cG_{x_u}$ on $\mathbb{R}^2$ with everywhere linearly independent differentials determines a standard real Killing spinor triple on $\mathbb{R}^4$ with metric:
\begin{equation}
\label{eq:metricBspinorgeneral}
g = \cH_{x_u} \dd x_u\otimes \dd x_u +   e^{\cF_{x_u}} \dd x_u\odot ( \dd x_v + \beta_{x_u}) + \frac{1}{4c ^2} \dd_X \cF_{x_{u}}\otimes \dd_X\cF_{x_{u}} + e^{\cF_{x_u}} \dd_X\cG_{x_u}\otimes \dd_X\cG_{x_u} 
\end{equation} 
	
\noindent
for a choice of $\beta_{x_{u}}$ as prescribed in Equation \eqref{eq:thmcondition2}. Conversely, every simply connected  standard real Killing spinor triple can be constructed in this way.
\end{cor}

\noindent
Theorem \ref{thm:conformallyBrinkmann} characterizes as well the local isometry type of every standard real Killing spinor triple.

\begin{cor}
Every four-dimensional space-time $(M,g)$ admitting a real Killing spinor is locally isometric to an open set of $\mathbb{R}^4$ equipped with the metric:
\begin{equation*}
g = \cH_{x_u} \dd x_u\otimes \dd x_u +   e^{\cF_{x_u}} \dd x_u\odot ( \dd x_v + \beta_{x_u}) + \frac{1}{4c ^2} \dd_X \cF_{x_{u}}\otimes \dd_X\cF_{x_{u}} + e^{\cF_{x_u}} \dd_X\cG_{x_u}\otimes \dd_X\cG_{x_u}   
\end{equation*} 
	
\noindent
for a choice of families of functions $\cF_{x_u}$ and $\cG_{x_u}$, and of one-forms $\beta_{x_{u}}$ as prescribed in Equation \eqref{eq:thmcondition2}. 
\end{cor}

\begin{ep}
Assume that $X = \mathbb{R}^2$ with Cartesian coordinates $(y_1,y_2)$, take $c = \frac{1}{2}$ and write $q_{x_u}$ as follows:
\begin{equation*}
q_{x_u} =  \dd_X \cF_{x_{u}}\otimes \dd_X\cF_{x_{u}} + e^{\cF_{x_u}} \dd_X\cG_{x_u}\otimes \dd_X\cG_{x_u}
\end{equation*}
	
\noindent
in terms of families of functions $\cF_{x_{u}}$ and $\cG_{x_{u}}$ on $X$. Equation \eqref{eq:thmcondition2} can be equivalently written as follows:
\begin{equation}
\label{eq:thmcondition2alt}
\dd_X \beta_{x_{u}} =   (   \partial_{x_u}  \dd_X \cG_{x_{u}} + e^{-\cF_{x_{u}}/2}   \ast_{q_{x_u}}\partial_{x_u}\dd_X\cF_{x_{u}} ) \wedge \dd_X \cG_{x_{u}} 
\end{equation}
	
\noindent
where we have used that $\ast_{q_{x_u}}\dd_X\cF_{x_{u}} =  e^{\cF_{x_{u}}/2} \dd_X\cG_{x_{u}}$. Assume that $\cF_{x_{u}}$ and $\cG_{x_{u}}$ are such that:
\begin{equation*}
\dd_X \cF_{x_{u}} = a_{x_u} \dd y_1 \, , \qquad \dd_X \cG_{x_{u}} = f_{x_u} \dd y_1  + b_{x_u} \dd y_2\, , \qquad x_u\in \mathbb{R}
\end{equation*}
	
\noindent
where $a_{x_u} , b_{x_u}$ and $f_{x_u}$ are families of constant functions on $X = \mathbb{R}^2$. In particular:
\begin{equation*}
q_{x_u} =   (a_{x_u}^2 + f_{x_u}^2 e^{a_{x_u} y_1 + k_{x_u}}) \dd y_1 \otimes \dd y_1 +  e^{a_{x_u} y_1 + k_{x_u}} (b_{x_u} f_{x_u} \dd y_1\odot \dd y_2 + b_{x_u}^2 \dd y_2 \otimes \dd y_2)
\end{equation*} 
	
\noindent
where $k_{x_u}$ is a family of constants. A quick computation shows that:
\begin{equation*}
e^{-\cF_{x_{u}}/2}   \ast_{q_{x_u}}\partial_{x_u}\dd_X\cF_{x_{u}} = \partial_{x_u}\log(a_{x_u}) \dd_X\cG_{x_{u}}
\end{equation*}
	
\noindent
where we assume that the family of constants $a_{x_u}$ are strictly positive. Hence:
\begin{equation*}
\dd_X \beta_{x_{u}} =  \partial_{x_u}  \dd_X \cG_{x_{u}}  \wedge \dd_X \cG_{x_{u}} = (b_{x_u}\partial_{x_u}f_{x_u} - f_{x_u}\partial_{x_u}b_{x_u}) \dd y_1 \wedge \dd y_2\, .
\end{equation*}
	
\noindent
Solutions to the previous equation can be easily found by direct inspection. For instance:
\begin{equation*}
\beta_{x_{u}} = \frac{1}{2}(b_{x_u}\partial_{x_u}f_{x_u} - f_{x_u}\partial_{x_u}b_{x_u}) (y_1   \dd y_2	- y_2   \dd y_1)
\end{equation*}
	
\noindent
which yields the following four-dimensional metric $g$ on $\mathbb{R}^4$:
\begin{eqnarray*}
& g = \cH_{x_u} \dd x_u\otimes \dd x_u +   \frac{1}{2} e^{a_{x_u} y_1 + k_{x_u}} \dd x_u\odot ( \dd x_v + (b_{x_u}\partial_{x_u}f_{x_u} - f_{x_u}\partial_{x_u}b_{x_u}) (y_1   \dd y_2	- y_2   \dd y_1))\\
& + (a_{x_u}^2 + f_{x_u}^2 e^{a_{x_u} y_1 + k_{x_u}}) \dd y_1 \otimes \dd y_1 +  e^{a_{x_u} y_1 + k_{x_u}} (b_{x_u} f_{x_u} \dd y_1\odot \dd y_2 + b_{x_u}^2 \dd y_2 \otimes \dd y_2)
\end{eqnarray*} 
	
\noindent
This provides an example of standard real Killing spinor triple for which the crossed term $\beta_{x_{u}}$ is not trivial and, in particular, not closed. 
\end{ep}

\noindent
It is instructive to present the metric $g$ occurring in Equation \eqref{eq:thmcondition1} of Theorem \ref{thm:conformallyBrinkmann} in an alternative equivalent form. Using the notation of Theorem \ref{thm:conformallyBrinkmann}, define:
\begin{equation*}
\cY_{x_u}  = e^{-\cF_{x_{u}}/2} \, , \quad x_u \in \mathbb{R}\, .	
\end{equation*}

\noindent
This defines a family $\left\{\cY_{x_u} \right\}_{x_u\in\mathbb{R}}$ of strictly positive functions on $X$. Substituting this expression into Equation \eqref{eq:thmcondition1} and relabeling some of the symbols adequately we obtain the following equivalent expression for $g$:
\begin{equation}
\label{eq:gmetricalternative}
g = \frac{1}{c^2 \cY_{x_u}^2} (\cK_{x_u}  \dd x_u\otimes \dd x_u + \dd x_u\odot (\dd x_v +   c^2 \beta_{x_u}) + \dd_X \cY_{x_{u}}\otimes \dd_X\cY_{x_{u}} +  \omega_{x_u} \otimes \omega_{x_u})
\end{equation}

\noindent
In this form it becomes apparent that standard real Killing spinor triple provide a vast generalization of Siklos space-times, which in turn can be interpreted as a deformation of the AdS$_4$ space-time. Indeed, assume that $(M,g,\varepsilon)$ is a standard real Killing spinor triple for which $\left\{ \cF_{x_{u}} \right\}_{x_u\in \mathbb{R}}$ and  $\left\{ \cG_{x_{u}} \right\}_{x_u\in \mathbb{R}}$ are both independent of the coordinate $x_u$. Then, by Equation \eqref{eq:gmetricalternative} it is clear that there exists local coordinates in which the metric $g$ reads:
\begin{equation*}
	g = \frac{1}{c^2 y^2} (\cK_{x_u}  \dd x_u\otimes \dd x_u + \dd x_u\odot (\dd x_v +  c^2 \beta_{x_u}) + \dd y \otimes \dd y + \dd w\otimes \dd w)
\end{equation*}

\noindent
In addition Equation \eqref{eq:thmcondition2} implies in this case that $\beta_{x_{u}}$ is closed on $X$, whence locally exact, a fact that can be used to redefine $\cK_{x_{u}}$ as well as the coordinate $x_v$ in order to absorb the one-form $\beta_{x_{u}}$ in such a way that $g$ locally reads:
\begin{equation*}
g = \frac{1}{c^2 y^2} (\cK_{x_u}  \dd x_u\otimes \dd x_u + \dd x_u\odot \dd x_v + \dd y \otimes \dd y + \dd w\otimes \dd w)\, .
\end{equation*}

\noindent
This is precisely the local four-dimensional metric constructed by Siklos in \cite{Siklos}, which defines what are nowadays called \emph{Siklos space-times} or \emph{Siklos gravitational waves}. The latter describe exact idealized gravitational waves moving through anti-de Sitter space-time \cite{GibbonsRuback,Podolsky:1997ik}. The fact that every Siklos space-time admits real Killing spinors was explicitly noticed in \cite{GibbonsRuback}. Therefore, standard real Killing spinor triples  provide a broad generalization of the Siklos class of space-times that also admits Killing spinors and reduces to the latter in certain special cases.  


\subsection{Axionic parallel spinors} 


In this subsection we briefly consider a natural class of differential spinors that are parallel relative to data $\fra \in \Omega^1(M,\wedge M)$ satisfying:  
\begin{equation}
\label{eq:axionica}
\fra = \fra^{(4)} = \beta\otimes \nu_g \in \Omega^1(M,\wedge^4 M)
\end{equation}
 
\noindent
where $\beta\in \Omega^1(M)$ is a given one-form. In other words, we consider spinors $\epsilon\in \Gamma(S)$ on a strongly spin Lorentzian four-manifold $(M,g)$ that satisfy the following equation:
\begin{eqnarray*}
\nabla^g_w \epsilon = \beta(w)\, \nu_g\cdot_g\epsilon\, , \qquad w\in \mathfrak{X}(M)
\end{eqnarray*}

\noindent
We will refer to these spinors as \emph{axionic parallel spinors}, since an element of the form \eqref{eq:axionica} acts on the corresponding spinor as a volume form, or using physics jargon, as an \emph{axionic scalar field}. Note that these spinors are parallel with respect to a connection that does not preserve the admissible bilinear pairing $\cB$, which to the best of my knowledge is a case that has not been studied in the literature. Theorem \ref{thm:differentialspinors4d} immediately implies the following existence result.
\begin{prop}
A strongly spin Lorentzian four-manifold $(M,g)$ admits a axionic parallel spinor relative to $\fra = \beta\otimes \nu_g \in \Omega^1(M,\wedge^4 M)$ if and only if it admits an isotropic parallelism $[u,v,l,n]\in \Ob(\mathfrak{F}(M))$, satisfying the following differential system:
\begin{eqnarray}
\label{eq:axionicparallelism}
\nabla^g u = 0\, , \quad \frac{1}{2} \nabla^g v = - \rho\otimes  n  - \kappa \otimes l  \, , \quad \frac{1}{2}\nabla^g l   =  \kappa\otimes u  - \beta \otimes n \, ,\quad  \frac{1}{2}\nabla^g n = \rho\otimes u  + \beta \otimes l 
\end{eqnarray}

\noindent
for one-forms $\kappa , \rho \in \Omega^1(M)$.
\end{prop}

\noindent
In particular, every Lorentzian manifold equipped with an axionic parallel spinor is Brinkmann, a fact that leads us to define the notion of \emph{axionic Brinkmann four-manifolds} as Lorentzian four-manifolds equipped with an axionic parallel spinor. Consequently, we will refer to isotropic parallelisms satisfying Equation \eqref{eq:axionicparallelism} as \emph{axionic}.

\begin{prop}
An isotropic parallelism $[u,v,l,n]$ is axionic relative to $\fra = \fra^{(4)} = \beta\otimes \nu_h$ if and only if any of its representatives $(u,v,l,n)\in [u, v, l, n]$ satisfies the following differential system:
\begin{eqnarray}
\label{eq:axionicparallelismdifferential}
\dd u = 0\, , \frac{1}{2}\quad \dd v = - \rho\wedge  n  - \kappa \wedge l  \, , \quad \frac{1}{2}\dd l   =  \kappa\wedge u  - \beta \wedge n \, ,\quad  \frac{1}{2}\dd n = \rho\wedge u  + \beta \wedge l 
\end{eqnarray}

\noindent
for one-forms $\kappa , \rho \in \Omega^1(M)$.
\end{prop}
 
\noindent
The fact that if a representative $(u,v,l,n)\in [u, v, l, n]$ satisfies \eqref{eq:axionicparallelismdifferential} then any other representative in $(u,v^{\prime},l^{\prime},n^{\prime}) \in [u,v,l,n]$ also satisfies the same system \eqref{eq:axionicparallelismdifferential} for possibly different one-forms $\kappa^{\prime} , \rho^{\prime} \in \Omega^1(M)$ can be verified explicitly. Write:
\begin{equation*} 
(u,v^{\prime},l^{\prime},n^{\prime}) = (u,v-\frac{1}{2} \vert\frw\vert_g^2 u + \frw,l - \frw(l) u,n - \frw(n) u)
\end{equation*}

\noindent
for a unique element $\frw\in \Gamma( \langle \mathbb{R}\, u\rangle  \oplus \langle \mathbb{R}\, v \rangle)^{\perp_h}$. We compute:
\begin{eqnarray*}
& \dd l^{\prime} = \dd (l - \frw(l) u) = \kappa^{\prime} \wedge u - \beta \wedge n^{\prime} = \kappa^{\prime} \wedge u - \beta \wedge (n - \frw(n) u)\\
& \dd n^{\prime} = \dd (n - \frw(n) u) = \rho^{\prime} \wedge u + \beta \wedge l^{\prime} = \rho^{\prime} \wedge u + \beta \wedge (l - \frw(l) u)
\end{eqnarray*}

\noindent
and thus:
\begin{equation*}
\kappa = \kappa^{\prime} + \dd (\frw(l)) + \frw(n) \beta\, , \qquad \rho = \rho^{\prime} + \dd (\frw(n)) - \frw(l) \beta
\end{equation*}

\noindent
Using these relations a computation shows that $\dd v^{\prime} = - \rho^{\prime}\wedge  n^{\prime}  - \kappa^{\prime} \wedge l^{\prime}$ holds if and only if $\dd v = - \rho\wedge  n  - \kappa \wedge l$ holds and thus $(u,v^{\prime},l^{\prime},n^{\prime})$ is axionic if and only if $(u,v,l,n)$ is. Axionic spinors seem to define an interesting class of Brinkmann four-manifolds which we plan to study in more detail in the future.  


\section{Four-dimensional Lorentzian instantons}
\label{sec:LorentzianInstantons} 


Thus far we have focused on the study of parallel spinors with respect to a given, possibly completely general, connection. It is however important to also consider \emph{algebraic constraints} on such spinors. These constraints can become, in fact, important coupled differential equations when they are adequately constructed in terms of variables of a given differential system. We illustrate this possibility in the following by considering a particular \emph{constraint} that occurs in the supersymmetry conditions of several supergravity equations, and whose Euclidean analog in four dimensions corresponds to the celebrated self-duality condition for the curvature of a connection on a principal bundle. We plan to study these supersymmetric constraints in more detail in the future, in particular in higher dimensions, where they relate to the notion of higher-dimensional instantons \cite{DonaldsonThomas}. 

Let $P\to M$ be a principal bundle with structure group $\G$ defined over the oriented four-manifold $M$, which we consider to be equipped with a Riemannian metric that we denote momentarily by $g_r$. The affine space of connections on $P$ will be denoted by $\cA_P$. Since $M$ is four-dimensional, the Hodge dual operator associated to $g_r$ and the given orientation maps two-forms to two-forms and squares to the identity:
\begin{eqnarray*}
\ast_{g_r} \colon \Omega^2(M) \to \Omega^2(M)\, , \qquad \ast_{g_r}^2 = \mathrm{Id}
\end{eqnarray*}

\noindent
Hence, the bundle of two-forms $\wedge^2 M = \wedge^2_{+} M \oplus \wedge^2_{-} M$ on $M$ splits as a direct sum of the eigenbundles of $\ast_{g_r}$ with positive and negative eigenvalue, respectively, and consequently we obtain an analogous decomposition at the level of sections:
\begin{equation*}
\Omega^2(M) = \Omega^2_{+}(M) \oplus \Omega^2_{-}(M)
\end{equation*}  

\noindent
Elements in $\Omega^2_{+}(M)$ are called self-dual, whereas elements in $\Omega^2_{-}(M)$ are called anti-self-dual. Given a connection $A\in \Omega^1(P,\frg)$, where $\frg$ denotes the Lie algebra of $\G$, its curvature defines a two-form with values in the adjoint bundle $\frg_P$ of $P$, namely:
\begin{eqnarray*}
	F_A \in \Omega^2(M,\frg_P)
\end{eqnarray*}

\noindent
The Hodge dual operation can be applied to the two-form part of $\Omega^2(M,\frg_P)$, where it again induces a splitting in eigenspaces:
\begin{equation*}
\Omega^2(M,\frg_P) = \Omega^2_{+}(M,\frg_P) \oplus \Omega^2_{-}(M,\frg_P)
\end{equation*}

\noindent
Connections $A\in \cA_P$ whose curvature $F_A$ belongs to $\Omega^2_{+}(M,\frg_P)$ are called \emph{self-dual instantons}, whereas connections $A\in\cA_P$ whose curvature $F_A$ belongs to $\Omega^2_{-}(M,\frg_P)$ are called \emph{anti-self-dual instantons}. By virtue of the Bianchi identity, both self-dual and anti-self-dual instantons satisfy the Yang-Mills equation:
\begin{equation*}
\dd_A \ast_{g_r} F_A = 0
\end{equation*}

\noindent
which can be obtained via a variational problem of the Yang-Mills functional. If $M$ is compact then the instantons on $P$, if any, are absolute minima of the Yang-Mills functional. Their study has led to outstanding results in the topology of compact low-dimensional manifolds and has evolved into a mathematical area of its own, sometimes called \emph{mathematical gauge theory} \cite{DonaldsonThomas}. Now, if we try to reproduce the previous discussion on a \emph{Lorentzian} four-manifold $(M,g)$, we immediately run into a wall: in four Lorentzian dimensions the Hodge operator on two-forms squares to minus the identity, and hence the concept of self-dual or anti-self dual instanton cannot be defined in terms of the Hodge operator within the framework described above. Instead, we can consider the following equation:
\begin{equation}
\label{eq:LorentzianInstanton}
F_A \cdot \varepsilon = 0
\end{equation}

\noindent
which in four Euclidean dimensions is equivalent to self-duality condition for $A$ but which, in contrast to the self-duality condition, also makes sense in four-dimensional Lorentzian signature and thus can be interpreted as the Lorentzian analog of a Euclidean instanton. Note that the \emph{dot} symbol in Equation \eqref{eq:LorentzianInstanton} denotes the Clifford multiplication of the \emph{two-form} part of $F_A$ on $\varepsilon$.  

\begin{definition}
A \emph{Lorentzian instanton} on $(P,M,g,\varepsilon)$ is a connection $A\in \cA_P$ that satisfies Equation \eqref{eq:LorentzianInstanton}.
\end{definition}

\noindent
Here $\varepsilon \in \Gamma(S)$ is an irreducible real spinor on $(M,g)$, and therefore the notion of Lorentzian instanton in four dimensions depends not only on the choice of a Riemannian metric and orientation, as it happens in the Riemannian case, but also on the choice of such spinor. By Remark \ref{remark:constraintequation4d}, a connection $A\in \cA_P$ satisfies Equation \eqref{eq:LorentzianInstanton} if and only if:
\begin{equation*}
(u+u\wedge l) \diamond_g F_A= u\wedge F_A +  F_A(u) + u\wedge l \wedge  F_A - l \wedge  F_A(u) + u \wedge  F_A(l) +  F_A(l,u) =  0
\end{equation*}

\noindent
Isolating by degree, we conclude that $A$ is a Lorentzian instanton on $(P,M,g,\varepsilon)$ if and only if the following equations hold:
\begin{equation*}
F_A(u) = 0\, , \qquad u \wedge  F_A = 0 
\end{equation*}

\noindent
where $u$ is the Dirac current of $\varepsilon \in \Gamma(S)$. Hence, we could consistently define the notion of Lorentzian instanton on a tuple $(P,M,g,u)$, where $u$ is an isotropic one-form on $M$, instead of on $(P,M,g,\varepsilon)$. Using the notion of isotropic parallelism we can obtain the following equivalent characterization of a Lorentzian instanton.

\begin{prop}
Let $[u,v,l,n]$ be the isotropic parallelism canonically associated to $\varepsilon \in \Gamma(S)$. Then, $A\in\cA_P$ is a Lorentzian instanton on $(P,M,g,\varepsilon)$ if and only if for any representative $(u,v,l,n) \in [u,v,l,n]$ there exists a section:
\begin{eqnarray*}
\omega_A \in \Gamma(( \langle \mathbb{R}\, u\rangle  \oplus \langle \mathbb{R}\, v \rangle)^{\perp_g} \otimes \frg_P)
\end{eqnarray*}

\noindent
such that $F_A = u\wedge \omega_A$.
\end{prop}

\begin{proof}
The curvature $F_A$ of $A$ satisfies $u\wedge F_A = 0$ if and only if there exists a one-form $\omega_A \in \Omega^1(M,\frg_P)$ taking values on the adjoint bundle $\frg_P$ of $P$ such that $F_A = u \wedge \omega_A$. This bundle valued one-form is unique modulo modifications of the form
\begin{eqnarray}
\label{eq:gaugetransformationA}
\omega_A \mapsto \omega_A + u\otimes \tau
\end{eqnarray}

\noindent
where $\tau \in \Gamma(\frg_P)$. Let $(u,v,l,n) \in [u,v,l,n]$. We can choose $\omega_A$ such that $\omega_A(v) = 0$ and expand:
\begin{eqnarray*}
\omega_A = v\otimes \omega_A(u) + l\otimes \omega_A(l) + n\otimes\omega_A(n)
\end{eqnarray*}

\noindent
Using this expansion, it follows that $F_A(u) = 0$ if and only if $\omega_A (u) = 0$, that is, if and only if $\omega_A \in \Gamma(( \langle \mathbb{R}\, u\rangle  \oplus \langle \mathbb{R}\, v \rangle)^{\perp_g} \otimes \frg_P)$. In particular:
\begin{eqnarray*}
F_A = u\wedge( l\otimes \omega_A(l) + n\otimes\omega_A(n))
\end{eqnarray*}

\noindent
The converse follows by construction and hence we conclude. 
\end{proof}
   
\noindent
Changing the isotropic parallelism $(u,v,l,n)$ within its class $[u,v,l,n]$ changes $\omega_A$ by a transformation of the form \eqref{eq:gaugetransformationA}. Hence, every Lorentzian instanton $A\in \cA_P$ defines a canonical element:
\begin{eqnarray*}
[\omega_A] \in \Gamma(\frG_u \otimes \frg_P)
\end{eqnarray*}  

\noindent
where $\frG_u$ is the screen bundle defined by $u\in \Omega^1(M)$. Hence, we obtain the following characterization of Lorentzian instantons. 

\begin{cor}
A connection $A\in \cA_P$ is an instanton on $(P,M,g,\varepsilon)$ if and only if there exists a section $[\omega_A] \in \Gamma(\frG_u \otimes \frg_P)$ such that $F_A = u\wedge \omega_A$, where $\omega_A\in [\omega_A]$ is any representative and $u$ is the Dirac current of $\varepsilon$.
\end{cor}

\noindent
The previous corollary shows that the Lorentzian instanton condition depends only on the Dirac current associated to $\varepsilon$ and not on the full information encoded in the latter. Hence, we can speak of a Lorentzian instanton on a tuple $(P,M,g,u)$ instead of $(P,M,g,\varepsilon)$. Given a Lorentzian instanton $A$ on $(P,M,g,u)$, a direct computation gives:
\begin{eqnarray*}
\ast_g F_A = u\wedge \ast_{q_u} [\omega_A]
\end{eqnarray*}

\noindent
where $\ast_{q_u}$ is the induced Hodge dual on $\frG_u$. Hence, a priori we cannot expect that a Lorentzian instanton in four dimensions is automatically a solution to the Yang-Mills equations, in sharp contrast with the Euclidean case. Note that $A$ being a Lorentzian instanton is a condition that will depend in general not only on the underlying orientation and Lorentzian metric $g$, but also on the choice of isotropic one-form $u$. 


\renewcommand{\leftmark}{\MakeUppercase{Chapter \thechapter. Parallel spinors with skew-symmetric torsion}}

\chapter{Parallel spinors with torsion}
\label{chapter:parallelspinorstorsion}


In this chapter we apply the theory of parabolic pairs and isotropic parallelisms to study irreducible real spinors parallel with respect to a metric connection possibly with non-parallel torsion on a Lorentzian four-manifold $(M,g)$. These spinors, to which we refer as \emph{torsion parallel spinors}, define a natural class of differential spinors and constitute the main object of study in this dissertation, especially in the case in which the torsion is totally skew-symmetric. The latter is the case of interest for applications to supergravity and will be considered in Chapter \ref{chapter:susyKundt4d}.


\section{Torsion parallel spinors}


Let $(M,g)$ be an oriented and time-oriented Lorentzian four-manifold. Denote by $\nabla^g$ the Levi-Civita connection on $(M,g)$. Every other metric connection $\nabla$ on $(M,g)$ can be written as follows:
\begin{equation*}
	\nabla_{w_1}w_2 = \nabla^g_{w_1} w_2 + \mathrm{A}(w_1,w_2)\, , \qquad w_1 , w_2 \in \mathfrak{X}(M)
\end{equation*}

\noindent
in terms of a uniquely defined tensor $\mA\in \Gamma(T^{\ast}M\otimes T^{\ast}M\otimes TM)$ that satisfies:
\begin{equation*}
	g(\mA(w_1,w_2),w_3) +g(w_2, \mA(w_1,w_3)) = 0\, , \quad \forall\,\, w_1 , w_2 , w_3 \in \mathfrak{X}(M)
\end{equation*}

\noindent
for every $w_1,w_2,w_3\in\mathfrak{X}(M)$. Hence $\mA\in \Gamma(T^{\ast}M\otimes \End_g (TM))$, where $\End_g(TM)$ denotes the endomorphisms of $TM$ that preserve $g$. We will refer to $\mA$ as the \emph{contorsion tensor} of the given metric connection with torsion on $(M,g)$, which we will consequently denote by $\nabla^{g,\mA}$. The vector space of all metric contorsion tensors on $(M,g)$ identifies with the space of sections $\Gamma(T^{\ast}M\otimes \End_g (TM))$. We introduce the notation:
\begin{equation*}
	\mA_{w} := \mA(w) \in \Gamma(\End_g (TM))\, , \qquad w  \in \mathfrak{X}(M)\, .
\end{equation*}

\noindent
which we will use occasionally in the following. For simplicity in the exposition we will identify $\mA \in \Gamma(T^{\ast}M\otimes \End_g (TM))$ with $\mA\in \Omega^1(M,\wedge^2 M)$ by means of the musical isomorphism defined by the underlying Lorentzian metric. 

\begin{remark}
We have described the space of connections compatible with a given Lorentzian metric using the notion of \emph{contorsion}. Equivalently, we could have used the notion of \emph{torsion}, which for a metric connection contains the same information as the contorsion. Given a metric connection $\nabla^{g,\mA}$ with contorsion $\mA$, its torsion reads:
\begin{equation*}
	\mT(w_1,w_2) = \mA(w_1,w_2) - \mA(w_2,w_1)\, .
\end{equation*}

\noindent
Hence, $\mT\in \Gamma(\wedge^2 M \otimes TM)$ and the previous formula identifies the space of contorsion tensors with the space of torsion tensors. 
\end{remark}

\noindent
Let $(M,g)$ be strongly spin and let $(S,\Gamma,\cB)$ be a paired irreducible spinor bundle on $(M,g)$. Since $\nabla^{g,\mA}$ is metric and $(M,g)$ is spin, $\nabla^{g,\mA}$ lifts canonically to the spinor bundle $S$ and defines a connection on $S$, which for simplicity we denote by the same symbol. Note that this connection is compatible with both Clifford multiplication and the admissible pairing $\cB$, see Proposition \ref{prop:compatiblekilling}. More explicitly, we have:
\begin{equation*}
	\nabla^{g,\mA}_w \varepsilon = \nabla^g_w \varepsilon - \frac{1}{2} \Psi_{\Gamma}(\mA_w)(\varepsilon) = \nabla^g_w \varepsilon - \frac{1}{2} \mA_w \cdot \varepsilon
\end{equation*}

\noindent
where, using the Lorentzian metric $g$, we have identified $\mA_w \in \Gamma(\End_g(TM))$ with $\mA_w \in \Omega^2(M)$ for every $w\in\mathfrak{X}(M)$. 


\begin{definition}
Let $(M,g)$ be a strongly spin Lorentzian four-manifold. A \emph{torsion parallel spinor} on $(M,g)$ with contorsion $\mA \in \Gamma(T^{\ast}M\otimes \End_g (TM))$ is a section $\varepsilon\in \Gamma(S)$ of a paired spinor bundle $(S,\Gamma,\cB)$ over $(M,g)$ that is parallel with respect to $\nabla^{g,\mA}$.
\end{definition}

\noindent
Torsion parallel spinors clearly define a particular class of differential spinors. More precisely, a parallel spinor with contorsion $\mA\in\Omega^3(M)$ is a differential spinor relative to an endomorphism-valued one-form $\cA\in \Omega^1(M,\End(S))$ whose symbol is given by:
\begin{equation}
	\label{eq:frawgeneral}
	\fra_w = \frac{1}{2} \mA_w  \in \Omega^2(M)\, , \qquad w\in \mathfrak{X}(M)
\end{equation}

\noindent
Hence, a direct application of Theorem \ref{thm:differentialspinors4d} gives the following result.

\begin{prop}
	\label{prop:TorsionSpinorParabolicPairgeneral}
	A strongly spin Lorentzian four-manifold $(M,g)$ admits a torsion parallel spinor with contorsion $\mA\in \Gamma(T^{\ast}M\otimes \End_g (TM))$ if and only if it admits a parabolic pair $(u,[l]_u)$ satisfying:
	\begin{equation}
		\label{eq:TorsionSpinorParabolicPairgeneral}
		\nabla^{g , \mA} u = 0	 \, , \qquad   \nabla^{g,\mA} l = \kappa\otimes u 
	\end{equation}
	
	\noindent 
	for any, and hence for all, representatives $l\in [l]_u$.
\end{prop}

\noindent
In the previous proposition $\nabla^{g , \mA}$ denotes the metric connection with contorsion $\mA$ induced on the contangent bundle by $\nabla^{g,\mA}$, which for simplicity is denoted by the same symbol.

\begin{remark}
\label{remark:various(u,l)general}
Let $(u,l)$ be a pair satisfying equations \eqref{eq:TorsionSpinorParabolicPairgeneral}. Any other representative $l^{\prime}\in [l]_u$ can be written as $l^{\prime} = l + f u$ for a function $f\in C^{\infty}(M)$. Then:
\begin{equation*}
\nabla^{g,\mA} l^{\prime} = \nabla^{g , \mA}(l + f u) =  \kappa\otimes u + \dd f \otimes u  	= \kappa^{\prime}\otimes u  
\end{equation*}
	
\noindent
where $\kappa^{\prime} = \kappa + \dd f$. Hence, if a representative of $(u, l \in [l]_u)$ satisfies the differential system \eqref{eq:TorsionSpinorParabolicPairgeneral} then any other representative also satisfies the same differential system.
\end{remark}

\noindent
In the following will refer to the parabolic pair associated to a spinor parallel with respect to a metric connection with contorsion simply as a \emph{torsion parabolic pair}.  

\begin{definition}
A \emph{torsion parabolic structure} on $M$ is a tuple $(g,(u,[l]_u),\mA)$ consisting of a Lorentzian metric $g$ on $M$ and  torsion parabolic pair $(u,[l]_u)$ on $(M,g)$ with contorsion $\mA\in \Gamma(\wedge^1 M\otimes \wedge^2 M)$. A \emph{parabolic torsion four-manifold} is a four-manifold equipped with a torsion parabolic structure.
\end{definition}

\noindent
Equivalently, a parabolic torsion Lorentzian four-manifold  $(M,g,(u,[l]_u),\mA)$ can be considered as a tuple $(M,g,\mA,\varepsilon)$, where $\varepsilon\in \Gamma(S)$ is a torsion spinor relative to $\mA$. Three-dimensional torsion Lorentzian manifolds have already been considered in \cite{Shahbazi3d} and, differently to the four-dimensional case considered in this dissertation, they can be equivalently studied exclusively in terms of their associated Dirac current. We rephrase now the result of Proposition \ref{prop:TorsionSpinorParabolicPairgeneral} in the language of isotropic parallelisms. 
\begin{prop}
\label{prop:existencenullcoframegeneral}
An oriented and strongly spin four-manifold $M$ admits a torsion parallel spinor with respect to a Lorentzian metric $g$ and a contorsion tensor $\mA$ if and only if it there exists an isotropic parallelism $[u,v,l,n]$ on $M$ satisfying the following differential system:
\begin{eqnarray}
\label{eq:existencenullcoframegeneral}
\nabla^{g,\mA} u = 0\, , \qquad \nabla^{g,\mA} v = - \kappa\otimes l - \rho\otimes n  \, , \qquad \nabla^{g,\mA} l = \kappa\otimes u \, , \qquad \nabla^{g,\mA} n = \rho\otimes u  
\end{eqnarray}
	
\noindent
for a given pair of one-forms $\kappa , \rho \in \Omega^1(M)$, where $(u,v,l,n)\in [u,v,l,n]$ is any representative in $[u,v,l,n]$.
\end{prop}

\begin{proof}
By Proposition \ref{prop:TorsionSpinorParabolicPairgeneral}, an oriented Lorentzian four-manifold $(M,g)$ admits a torsion parallel spinor if and only if its associated parabolic pair $(u,[l]_u)$ satisfies equations \eqref{eq:TorsionSpinorParabolicPairgeneral}. Assume then that there exists a torsion parabolic structure $(g,u,[l]_u)$.  We proceed by taking the covariant derivative of the identity $u\wedge n = \ast_g (u\wedge l)$ and expanding the result, obtaining, for every $w\in\mathfrak{X}(M)$:
\begin{eqnarray*}
u\wedge \nabla^{g,\mA}_w n = n\wedge \nabla^{g,\mA}_w u + \ast_g (\nabla^{g,\mA}_w u\wedge l) + \ast_g (u\wedge \nabla^{g,\mA}_w l) = 0
\end{eqnarray*}
	
\noindent
where we have used the expressions for $\nabla^{g,\mA}_w u$ and $\nabla^{g,\mA}_w l$ given in \eqref{eq:TorsionSpinorParabolicPairgeneral}. Hence $\nabla^{g,\mA}_w n = \rho\otimes u$ for a uniquely determined one-form $\rho\in \Omega^1(M)$. This gives the last equation in \eqref{eq:existencenullcoframegeneral}. Then, we take the covariant derivative of the identity $u\wedge v = \ast_g (l\wedge n)$, obtaining, for every $w\in\mathfrak{X}(M)$:
\begin{equation*}
u\wedge \nabla^{g,\mA}_w v = \ast_g (\nabla^{g,\mA}_w l  \wedge n ) + \ast_g ( l\wedge \nabla^{g,\mA}_w n ) + v\wedge \nabla^{g,\mA}_w u = - u \wedge (\kappa(w) l + \rho(w) n  ) 
\end{equation*}
	
\noindent
and thus:
\begin{equation*}
\nabla^{g,\mA} v = \omega \otimes u - \kappa\otimes l - \rho\otimes n 
\end{equation*}
	
\noindent
for a certain one-form $\omega\in \Omega^1(M)$. However, since $g(v,v) = 0$, we must have:
\begin{equation*}
0 =	g(\nabla^{g,\mA}_w v , v) = \omega(w) = 0 \, , \qquad \forall\,\, w \in \mathfrak{X}(M)
\end{equation*}
	
\noindent
and thus we obtain the second equation in \eqref{eq:existencenullcoframegeneral}. For the converse, we note that the differential system \eqref{eq:existencenullcoframegeneral} implies equations \eqref{eq:TorsionSpinorParabolicPairgeneral} in Proposition \ref{prop:TorsionSpinorParabolicPairgeneral} for a torsion parabolic pair. Hence, we only need to verify that the differential system \eqref{eq:existencenullcoframegeneral} gives a consistent prescription for a metric connection with contorsion $\mA$ compatible with the Lorentzian metric:
\begin{equation*}
g = u \odot v + l\otimes l + n\otimes n
\end{equation*}
	
\noindent
when applied to $(u,v,l,n)$. We compute:
\begin{eqnarray*}
&\nabla^{g,\mA}_w g = u \odot \nabla^{g,\mA}_w v + l\odot \nabla^{g,\mA}_w l + n \odot \nabla^{g,\mA}_w n \\
& = - u \odot (\kappa(w) l + \rho(w) n) + \kappa(w) u\odot l + \rho(w) u \odot n = 0
\end{eqnarray*}
	
\noindent
and therefore the prescription given by \eqref{eq:existencenullcoframegeneral} preserves $g$. A direct computation shows that it also preserves the orthogonality properties of $(u,v,l,n)$. Thus, we only need to check that the prescription given by \eqref{eq:existencenullcoframegeneral} defines a connection with precisely contorsion tensor given by $\mA \in \Gamma(T^{\ast}M\otimes \End_g(TM))$. For this, we need to compute the Lie brackets of the frame $(u^{\sharp_g} , v^{\sharp_g} , l^{\sharp_g} , n^{\sharp_g})$. A tedious computation gives:
\begin{eqnarray*}
& \left[ u^{\sharp_g} , v^{\sharp_g} \right]= - \kappa(u) l^{\sharp_g} - \rho(u) n^{\sharp_g} - \mT(u^{\sharp_g} , v^{\sharp_g})\, , \qquad  \left[ u^{\sharp_g} , l^{\sharp_g}\right] =  \kappa(u)  u^{\sharp_g} - \mT(u^{\sharp_g} , l^{\sharp_g})  \\
& \left[ u^{\sharp_g} , n^{\sharp_g}\right] = \rho(u)  u^{\sharp_g} - \mT(u^{\sharp_g} , n^{\sharp_g})\, , \qquad  \left[ v^{\sharp_g} , l^{\sharp_g}\right] = \kappa(v) u^{\sharp_g} + \kappa(l) l^{\sharp_g} +  \rho(l)  n^{\sharp_g} - \mT(v^{\sharp_g},l^{\sharp_g})\\
& \left[ v^{\sharp_g} , n^{\sharp_g}\right] = \rho(v) u^{\sharp_g} + \kappa(n) l^{\sharp_g} + \rho(n) n^{\sharp_g} - \mT(v^{\sharp_g},n^{\sharp_g})\, , \,\, \left[ l^{\sharp_g} , n^{\sharp_g}\right] =(\rho(l) - \kappa(n)) u^{\sharp_g}  - \mT(l^{\sharp_g},n^{\sharp_g})
\end{eqnarray*}
	
\noindent
where $\mT(w_1,w_2) = \mA(w_1,w_2) - \mA(w_2,w_1)$. Using these Lie brackets it can be checked that the contorsion is indeed  $\mA$ and thus we conclude. 
\end{proof}

\noindent
Let $(u,[l]_u)$ be a torsion parabolic pair on $(M,g)$ and let $(u,v,l,n)\in \mathbb{E}^{-1}(g,u,[l]_u)$ be an associated isotropic coframe satisfying the differential system \eqref{eq:existencenullcoframegeneral}. Recall that $\mathbb{E}\colon \mathfrak{P}(M)\to \mathfrak{F}(M)$ denotes the natural equivalence between the category of parabolic structures on $M$ and the category of isotropic parallelisms on $M$ introduced in Proposition \ref{prop:categoricalequivalence}. By Proposition \ref{prop:primerparallelizationglobal}, for any other isotropic coframe $(u,v^{\prime},l^{\prime},n^{\prime})\in \mathbb{E}^{-1}(g,u,[l]_u)$ we can write:
\begin{equation*}
(u,v^{\prime},l^{\prime},n^{\prime}) = \mathfrak{w}\cdot (u,v,l,n)  = (u,v-\frac{1}{2} \vert \frw\vert^2_g u + \frw,l - \frw (l) u,n - \frw (n) u)
\end{equation*}

\noindent
for a unique vector field $\frw\in \Gamma( \langle \mathbb{R}\, u\rangle  \oplus \langle \mathbb{R}\, v \rangle)^{\perp_g}$ in the orthogonal complement of the distribution spanned by $u$ and $v$. In particular, $[u,v,l,n] = \mathbb{E}^{-1}(g,u,[l]_u)$ is a torsor over $\Gamma( \langle \mathbb{R}\, u\rangle  \oplus \langle \mathbb{R}\, v \rangle)^{\perp_g}$ with respect to the action given above. Then, by Remark \ref{remark:various(u,l)general} we know that:
\begin{equation*}
	\nabla^{g,\mA} l^{\prime} = \kappa^{\prime} \otimes u 
\end{equation*}

\noindent
for $\kappa^{\prime} = \kappa - \dd(\frw(l))$. Similarly, we obtain:
\begin{equation*}
	\nabla^{g,\mA} n^{\prime} = \rho^{\prime} \otimes u 
\end{equation*}

\noindent
for $\rho^{\prime} = \rho - \dd(\frw(n))$. A satisfying computation shows now that the remaining equation in \eqref{eq:existencenullcoframegeneral} holds automatically, that is, we have:
\begin{equation*}
	\nabla^{g,\mA} v^{\prime} = - \kappa^{\prime}\otimes l^{\prime} - \rho^{\prime}\otimes n^{\prime} 
\end{equation*}

\noindent
and hence we confirm the following crucial result.
\begin{lemma}
	\label{lemma:allnullparallelisms}
	If an isotropic coframe $(u,v,l,n)\in \mathbb{E}^{-1}(g,u,[l]_u)$ satisfies the differential system \eqref{eq:existencenullcoframegeneral} with respect to $\kappa,\rho\in \Omega^1(M)$, then any other global isotropic coframe:
	\begin{equation*}
		(u,v^{\prime},l^{\prime},n^{\prime}) = \mathfrak{w} \cdot (u,v,l,n) \in \mathbb{E}^{-1}(g,u,[l]_u) 
	\end{equation*}
	
	\noindent
	in the same isotropic parallelism $\mathbb{E}^{-1}(g,u,[l]_u)$ also satisfies it with respect to $\kappa^{\prime} = \kappa - \dd(\frw(l)),\rho^{\prime} = \rho - \dd(\frw(n))\in \Omega^1(M)$. 
\end{lemma}

\noindent
As mentioned earlier, the contorsion of a metric connection is a section of $T^{\ast}M\otimes \End_g(TM)$ or, alternatively, a section of $T^{\ast}M\otimes \wedge^2 M$. As a fiber-wise representation of the Lorentz algebra, $T^{\ast}M \otimes \wedge^2 M$ is reducible and splits in terms of three irreducible representations:
\begin{equation*}
	T^{\ast}M\otimes \wedge^2 M = T^{\ast}M \oplus \cT M \oplus \wedge^3 M
\end{equation*}

\noindent
where:
\begin{equation*}
	\cT M = \left\{ \tau \in T^{\ast}M\otimes \wedge^2 M \,\, \vert \,\, \tau_{w_1}(w_2,w_3) + \tau_{w_3}(w_1,w_2) + \tau_{w_2}(w_3,w_1) = 0\,\,\vert\,\, \epsilon^{\mu\nu} \tau_{e_{\mu}}(e_{\nu}) = 0 \,\, \right\}
\end{equation*}

\noindent
where $\epsilon^{\mu\nu} = g(e^{\mu},e^{\nu})$ in terms of any orthonormal coframe $(e^0,\hdots,e^3)$ with dual frame $(e_0 , \hdots , e_3)$ and summation over repeated indices is assumed. In particular, for every contorsion tensor $\mA\in \Gamma(T^{\ast}M\otimes \wedge^2 M)$ there exists a unique one-form $\xi\in\Omega^1(M)$, a unique three-form $H\in \Omega^3(M)$ and a tensor $\tau \in \cT M$ such that:
\begin{equation*}
\mA_{w_1}(w_2) := \mA (w_1,w_2) = g(w_1,w_2) \xi - \xi(w_2) w_1 + \tau(w_1 , w_2) + \frac{1}{2} H(w_1, w_2) 
\end{equation*}

\noindent
and consequently:
\begin{eqnarray*}
\nabla^{g,\mA}_{w} \beta = \nabla^g_{w} \beta + \beta(w) \xi - \beta(\xi) w + \tau(w , \beta) + \frac{1}{2} H(w,\beta) 
\end{eqnarray*}

\noindent
where as usual we are using the same symbol to denote a vector or one-form and its metric dual. Skew-symmetrization of the previous formula together with the cyclic property satisfied by $\tau \in \Gamma(\cT M)$ gives the following result.
\begin{lemma}
\label{lemma:skewsymconnection}
The following formulae holds:
\begin{eqnarray*}
(\nabla^{g,\mA}_{w_1} \beta) (w_2) - (\nabla^{g,\mA}_{w_2} \beta) (w_1) = \dd\beta (w_1,w_2) + (\beta \wedge \xi)(w_1 , w_2)  + \tau(\beta ,w_1,w_2)  +   H(w_1,\beta , w_2) 
\end{eqnarray*}
	
\noindent
for every $w_1 , w_2 \in \mathfrak{X}(M)$.
\end{lemma}

\noindent
Using this result together with Proposition \ref{prop:existencenullcoframegeneral}, we can describe torsion parallel spinors in terms of an exterior differential system for global coframes on $M$ that does not use explicitly any Lorentzian metric. As we will see in the following chapters, this reformulation is particularly convenient to study moduli spaces of solutions and initial data in applications of skew-torsion parallel spinors to supergravity. 

\begin{thm}
\label{thm:existencenullcoframeII}
A strongly spin four-manifold $M$ admits a torsion parallel spinor with contorsion $\mA= \xi \oplus \tau \oplus H$ if and only if it there exists a global coframe $(u,v,l,n)$ and a pair of one-forms $\kappa , \rho \in \Omega^1(M)$ satisfying the following differential system:
\begin{eqnarray}
& \dd u = H_{u} + \xi\wedge u - \tau_u  \, , \qquad \dd v = H_{v} + \xi\wedge v - \tau_v - \kappa\wedge l - \rho\wedge n   \label{eq:duvtorsiongeneral}\\
& \dd l = H_{l} + \xi\wedge l - \tau_l + \kappa\wedge u \, , \qquad \dd n = H_{n} + \xi\wedge n - \tau_n + \rho\wedge u   \label{eq:dlntorsiongeneral}
\end{eqnarray}
	
\noindent
where $g = u\odot v + l\otimes l + n\otimes n$ is the Lorentzian metric associated to $[u,v,l,n]$.
\end{thm}

\begin{remark}
The previous theorem should be understood as giving the necessary and sufficient conditions for a strongly spin four-manifold $M$ to admit a torsion parallel spinor for \emph{some} Lorentzian metric and bundle of irreducible spinors. Hence, we do not fix an \emph{a priori} Lorentzian metric; such metric can be determined \emph{a posteriori} by $(u,v,l,n)$ and depends on the isotropic-parallelism class of the latter.
\end{remark}

\begin{proof}
The \emph{only if} direction follows by skew-symmetrization of equations \eqref{eq:existencenullcoframegeneral} in Proposition \ref{prop:existencenullcoframegeneral} upon use of Lemma \ref{lemma:skewsymconnection}. This immediately gives equations \eqref{eq:duvtorsiongeneral} and \eqref{eq:dlntorsiongeneral}. For the converse, we observe that the symmetrization of \eqref{eq:existencenullcoframegeneral} gives the following system of equations:
\begin{eqnarray*}
& \cL_u g = 2 \xi(u) g	- u\odot \xi  - \tau^s_u\, , \qquad \cL_v g = 2 \xi(v) g	- v\odot \xi - \tau^s_v - \kappa\odot l - \rho\odot n \\
& \cL_l g = 2 \xi(l) g	- l\odot \xi - \tau^s_l + \kappa\odot u\, , \qquad \cL_n g = 2 \xi(n) g	- n\odot \xi - \tau^s_n  + \rho\odot u 
\end{eqnarray*}

\noindent
where $\tau^s_{\beta}(w_1 , w_2) = \tau(w_1,\beta,w_2) + \tau(w_2,\beta,w_1)$, $w_1 , w_2 \in \mathfrak{X}(M)$. To obtain this symmetrization we have used the following identity:
\begin{eqnarray*}
& (\nabla^{g,\mA}_{w_1} \beta) (w_2) + (\nabla^{g,\mA}_{w_2} \beta) (w_1) = \cL_{\beta} g (w_1,w_2) + (\beta \odot \xi)(w_1 , w_2) \\
& - 2 \beta(\xi) g(w_1,w_2)  + \tau (w_1,\beta,w_2) + \tau(w_2,\beta,w_1) 
\end{eqnarray*}

\noindent
applied to each of the elements in $(u,v,l,n)$. Explicitly computing using equations \eqref{eq:duvtorsiongeneral} and \eqref{eq:dlntorsiongeneral}, we obtain:
\begin{eqnarray*}
& \cL_u g = \dd u (u)\odot \dd v  + \dd u \odot \dd v(u) + l \odot \dd l (u) + n \odot \dd n (u)\\
& = (\xi(u) u - \tau(u,u)) \odot v + u\odot (H(v,u) + \xi(u) v - \xi -\tau(v,u) - \kappa(u) l - \rho(u) n) \\
& + l \odot (H(l,u) + \xi(u) l - \tau(l,u) + \kappa(u) u) + n \odot (H(n,u) + \xi(u) n - \tau(n,u) + \rho(u) u)\\
& =  2 \xi(u) g	- u\odot \xi  - \tau^s_u\, .
\end{eqnarray*}
	
\noindent
This gives the symmetrization of the first equation in \eqref{eq:existencenullcoframegeneral}. Computing similarly for $\cL_v g$, $\cL_l g$ and $\cL_n g$, we obtain the symmetrization of the remaining equations in \eqref{eq:existencenullcoframegeneral} and thus we conclude. 
\end{proof} 

\begin{ep}
As an immediate consequene of the previous theorem, it follows that a strongly spin four-manifold $M$ admits a vectorial-torsion parallel spinor with vectorial contorsion $\xi \in \Omega^1(M)$ if and only if it there exists a global coframe $(u,v,l,n)$ and a pair of one-forms $\kappa , \rho \in \Omega^1(M)$ satisfying the following differential system:
\begin{eqnarray*}
\dd u =  \xi\wedge u  \, , \qquad \dd v =  \xi\wedge v  - \kappa\wedge l - \rho\wedge n  \\
\dd l = \xi\wedge l + \kappa\wedge u \, , \qquad \dd n = \xi\wedge n + \rho\wedge u   
\end{eqnarray*}
	
\noindent
where $g = u\odot v + l\otimes l + n\otimes n$ is the Lorentzian metric associated to $[u,v,l,n]$. Taking $\rho = \kappa = 0$ as a particular case, we obtain: 
\begin{eqnarray*}
\dd u =  \xi\wedge u  \, , \qquad \dd v =  \xi\wedge v  \, , \qquad  \dd l = \xi\wedge l  \, , \qquad \dd n = \xi\wedge n  
\end{eqnarray*}

\noindent
The integrability conditions of this system reduce to $\dd\xi = 0$, and thus we conclude that locally conformally parallel coframes defined on a locally conformally flat, strongly spin, Lorentzian four-manifold define parallel spinors with respect to a metric connection with vectorial torsion. 
\end{ep}

\noindent
By Lemma \ref{lemma:allnullparallelisms} we have the following alternative version of Proposition \ref{prop:existencenullcoframeII} to characterize strongly spin Lorentzian four-manifolds that admit torsion parallel spinors.

\begin{cor}
A strongly spin four-manifold $M$ admits a torsion parallel spinor if and only if it there exists an isotropic parallelism satisfying \eqref{eq:duvtorsiongeneral} and \eqref{eq:dlntorsiongeneral} for some contorsion tensor $\mA= \xi \oplus \tau \oplus \frac{1}{2} H$.
\end{cor}

\noindent
We will refer to the isotropic parallelisms that satisfy the differential system \eqref{eq:duvtorsiongeneral} and \eqref{eq:dlntorsiongeneral} as \emph{torsion isotropic parallelisms} when needed. There is a canonical equivalence of categories between the category of torsion parabolic pairs and the category of torsion isotropic parallelisms given by the restriction of $\mathbb{E}\colon \mathfrak{P}(M) \to \mathfrak{F}(M)$ to the corresponding full subcategories. Hence, we can study torsion parallel spinors in terms of either torsion parabolic structures or torsion isotropic parallelisms. Theorem \ref{thm:existencenullcoframeII} establishes that a Lorentzian four-manifold equipped with a torsion parallel spinor is naturally equipped with a nowhere vanishing isotropic vector field satisfying a differential equations \eqref{eq:duvtorsiongeneral} and \eqref{eq:dlntorsiongeneral}. Nowhere vanishing isotropic vector fields are vital in Lorentzian geometry and mathematical general relativity, since when they are Killing they determine idealized models for gravitational waves such as the class of Brinkmann space-times or more generally the class of Kundt space-times. 

\begin{definition}
\label{def:Kundt}
A \emph{Kundt four-manifold} is a triple $(M,g,u)$ consisting of a Lorentzian four-manifold $(M,g)$ equipped with an isotropic one-form $u\in \Omega^1(M)$ such that:
\begin{eqnarray*}
\nabla^{g\ast}u = 0\, , \qquad \vert\dd u\vert_g^2 = 0\, , \qquad \vert \cL_{u^{\sharp^g}} g\vert_g^2 = 0
\end{eqnarray*}

\noindent
Kundt space-times define a remarkable class of space-times that has been intensively studied in the general relativity literature, see \cite{Boucetta:2022vny} for a geometric characterization of this class of Lorentzian manifolds. Equivalently, using general relativity jargon, $u\in \Omega^1(M)$ defines a non-expanding, non-twisting and non-shear geodesic null congruence on $(M,g)$.
\end{definition}

\noindent
The differential system \eqref{eq:duvtorsiongeneral} and \eqref{eq:dlntorsiongeneral} can be immediately applied to obtain conditions on a torsion parallel spinor that imply the underlying Lorentzian manifold is Kundt with respect to the Dirac current of the torsion parallel spinor. 
\begin{prop}
Let $\varepsilon$ be a torsion parallel spinor on a Lorentzian four-manifold $(M,g)$ with associated parabolic pair $(u,[l]_u)$. Then, $(M,g,u)$ is a Kundt four-manifold if and only if:
\begin{eqnarray*}
& \nabla^{g\ast} u = 3 u(\xi)\, , \qquad \cL_{u^{\sharp_g}} g + u\odot \xi - 2 u(\xi) g + \tau^s(u) = 0\\
& \vert H_u \vert_g^2 - 2 \langle H_u , \tau_u\rangle_g + 2\tau_u (u,\xi) - u(\xi)^2 + \vert \tau_u\vert^2_g = 0
\end{eqnarray*}

\noindent
where $\tau^s(u) \in \Gamma(T^{\ast}M\odot T^{\ast}M)$ is the symmetrization of $\tau(u)$.
\end{prop}

\begin{proof}
By Theorem \ref{thm:existencenullcoframeII} we have:
\begin{equation*}
\nabla^{g,\mA}_{w} u = \nabla^g_{w} u +  u(w) \xi - u(\xi) w + \tau_w (u) + \frac{1}{2} H(w,u) = 0
\end{equation*}

\noindent
Taking the trace of this equation we obtain:
\begin{eqnarray*}
\nabla^{g\ast} u = 3 u(\xi)
\end{eqnarray*}

\noindent
whereas taking the symmetrization we obtain:
\begin{eqnarray*}
\cL_{u^{\sharp_g}} g + u\odot \xi - 2 u(\xi) g + \tau^s(u) = 0
\end{eqnarray*}

\noindent
Again by Theorem \ref{thm:existencenullcoframeII}, we have:
\begin{eqnarray*}
\vert \dd u \vert_g^2 = \vert H_u \vert_g^2 - 2 \langle H_u , \tau_u\rangle_g + 2\tau_u (u,\xi) - u(\xi)^2 + \vert \tau_u\vert^2_g = 0
\end{eqnarray*}

\noindent
and hence we conclude.
\end{proof}

\noindent
By Lemma \ref{lemma:allnullparallelisms}, every torsion parallel spinor defines two \emph{invariants} that are determined in terms of the pair $(\kappa,\rho)$ associated to any torsion isotropic parallelization. More precisely, we have a natural map:
\begin{eqnarray*}
\mathfrak{M}(M) \to \frac{\Omega^1(M)}{\dd C^{\infty}(M)}\times \frac{\Omega^1(M)}{\dd C^{\infty}(M)}\, , \qquad (g,u,[l]_u) \mapsto ([\kappa],[\rho]) 
\end{eqnarray*}

\noindent
where $\mathfrak{M}(M)$ is the set of equivalence classes in $\mathfrak{P}(M)$, that is the set of torsion parabolic structures modulo diffeomorphisms isotopic to the identity. 

\begin{definition}
Let $(u,[l]_u)$ be a torsion parabolic pair. Using the notation introduced above, the pair $([\kappa],[\rho])$ are the \emph{rank-one invariants} of $(u,[l]_u)$. 
\end{definition}

\noindent
It would be interesting to further elucidate the basic properties of this map, which we plan to study in more detail in the future.  We end this section with a corollary that follows immediately by our earlier computation given in the proof of Proposition \ref{prop:existencenullcoframegeneral} of the Lie brackets of any representative $(u,v,l,n)\in [u,v,l,n]$ of a torsion isotropic parallelism. As explained in Chapter \ref{chapter:IrreducibleSpinors4d}, every torsion isotropic parallelism $[u,v,l,n]$ defines a screen bundle $\mathfrak{G}_u$ associated to its Dirac current $u\in\Omega^1(M)$. A choice of one-form conjugate to $u$ defines a representative $(u,v,l,n)\in [u,v,l,n]$ and a rank-two distribution in $TM$ that is isomorphic to $\mathfrak{G}_u$ and is spanned by $(l^{\sharp_g},n^{\sharp_g})$. It is then natural to ask if there exists a choice of conjugate vector field $v$ to $u$ such that the rank-two distribution that it defines is integrable. This leads us to the following result.

\begin{cor}
Let $[u,v,l,n]$ be a torsion isotropic parallelism on $M$ with contorsion $\mA$. There exists an integrable realization of the screen bundle determined by $u$ if and only if there exists an isotropic coframe $(u,v,l,n)\in [u,l,v,n]$ satisfying:
\begin{eqnarray*}
\xi\lrcorner_g (l\wedge n) + \tau(l , n) - \tau(n , l) +   H(l, n) = (\rho(l) - \kappa(n)) u  
\end{eqnarray*}
	
\noindent
where $\mA = \xi \oplus \tau \oplus \frac{1}{2} H$ is split in its irreducible components. If that is the case, then the corresponding foliation has flat leaves when endowed with the metric induced by $g$.
\end{cor}

\begin{proof}
By the proof of Proposition \ref{prop:existencenullcoframegeneral} we know that:
\begin{equation*}
\left[ l^{\sharp_g} , n^{\sharp_g}\right] =(\rho(l) - \kappa(n)) u^{\sharp_g}  - \mT(l^{\sharp_g},n^{\sharp_g})
\end{equation*}

\noindent
and thus the span of $l^{\sharp_g}$ and $n^{\sharp_g}$ is integrable if and only if:
\begin{equation*}
\mT(l^{\sharp_g},n^{\sharp_g}) = (\rho(l) - \kappa(n)) u^{\sharp_g} 
\end{equation*}

\noindent
On the other hand, we have:
\begin{equation*}
\mT(l^{\sharp_g},n^{\sharp_g}) = \mA(l^{\sharp_g},n^{\sharp_g}) - \mA(n^{\sharp_g},l^{\sharp_g}) = \xi(l) n - \xi(n) l + \tau(l , n) - \tau(n , l) +   H(l, n)
\end{equation*}

\noindent
and thus the span of $l^{\sharp_g}$ and $n^{\sharp_g}$ is integrable if and only if:
\begin{equation*}
\xi\lrcorner_g (l\wedge n) + \tau(l , n) - \tau(n , l) +   H(l, n) = (\rho(l) - \kappa(n)) u 
\end{equation*}

\noindent
Using the previous equation together with the differential system given in \eqref{eq:duvtorsiongeneral} and \eqref{eq:dlntorsiongeneral}, we compute:
\begin{equation*}
	\dd l (l,n) = - \xi(n) - \tau_l(l,n) = 0\, , \qquad \dd n (l,n) = \xi(l) - \tau_l(l,n) = 0\, , 
\end{equation*}

\noindent
and therefore every leaf of the distribution determined by the span of $l^{\sharp_g}$ and $n^{\sharp_g}$ admits an orthonormal flat frame given by the pull-back of $(l,n)$ and thus we conclude. 
\end{proof}

\noindent
The foliation determined by an integrable screen bundle corresponds to the \emph{wave front} of $(M,g)$ when the latter is interpreted as a \emph{gravitational wave}. Hence, and interestingly enough, the previous corollary implies that the gravitational waves defined by torsion parallel spinors have all flat wave fronts.


\section{Lorentzian curvature and cohomological invariants}
\label{sec:generalcurvature}


The goal of this section is to \emph{compute} the curvature of a Lorentzian four-manifold equipped with a torsion parallel spinor. We will find that it is highly constrained in a quite elegant way. In our conventions the curvature tensor of $\nabla^{g,\mA}$, which we denote by $\cR^{g,\mA}$, is given by:
\begin{equation*}
	\cR^{g,\mA}_{w_1 , w_2} w_3 = \nabla^{g,\mA}_{w_1}\nabla^{g,\mA}_{w_2}w_3 - \nabla^{g,\mA}_{w_2}\nabla^{g,\mA}_{w_1}w_3 - \nabla^{g,\mA}_{[w_1 , w_2]}w_3\, , \quad \forall\,\, w_1, w_2, w_3 \in \mathfrak{X}(M)
\end{equation*}

\noindent
which expands as follows in terms of the torsion $\mA$ and its covariant derivative:
\begin{eqnarray*}
	& \cR^{g,\mA}_{w_1 , w_2} w_3 = \cR^{g}_{w_1 , w_2} w_3 + (\nabla^g_{w_1}\mA)(w_2,w_3)  -  (\nabla^g_{w_2}\mA)(w_1,w_3) \\
	& +  \mA (w_1,\mA (w_2,w_3)) -  \mA (w_2,\mA (w_1,w_3)) \, , \quad \forall\,\, w_1, w_2, w_3 \in \mathfrak{X}(M)
\end{eqnarray*}

\noindent
where $\cR^{g}$ denotes the Riemann tensor of $g$. For any representative $(u,l,v,n)\in [u,v,l,n]$ of a torsion isotropic parallelism $[u,v,l,n]$ and every pair of vector fields $w_1, w_2\in\mathfrak{X}(M)$, a computation using \eqref{eq:existencenullcoframegeneral} gives:
\begin{eqnarray*}
& \cR^{g,\mA}_{w_1 w_2} u = 0\, , \qquad  \cR^{g,\mA}_{w_1 w_2} v = -\dd \kappa(w_1,w_2) l - \dd\rho(w_1,w_2) n\\
& \cR^{g,\mA}_{w_1 w_2} l = \dd\kappa(w_1,w_2) u\, , \qquad \cR^{g,\mA}_{w_1 w_2} n = \dd\rho(w_1,w_2) u  
 \end{eqnarray*}

\noindent
and therefore, we obtain the following result.
\begin{prop}
\label{prop:curvaturetorsion}
Let $(M,g)$ be a strongly spin Lorentzian four-manifold equipped with a torsion parallel spinor with associated torsion isotropic parallelism $[u,v,l,n]$. Then:
\begin{eqnarray}
\label{eq:torsioncurvaturegeneral}
\cR^{g,\mA}_{w_1 w_2}   =   - \dd \kappa(w_1,w_2) \,u\wedge l - \dd\rho(w_1,w_2)\, u\wedge n\, , \qquad \forall\,\, w_1 , w_2 \in \mathfrak{X}(M)
\end{eqnarray}

\noindent
where $(u,v,l, n) \in [u,v,l,n]$ satisfies the differential system \eqref{eq:existencenullcoframegeneral} relative to $\kappa , \rho \in \Omega^1(M)$.
\end{prop}

\noindent
Note that Equation \eqref{eq:torsioncurvaturegeneral} is clearly invariant with respect to the choice of isotropic coframe representative $(u,v,l,n)\in [u,v,l,n]$.

\begin{definition}
A torsion isotropic parallelism $[u,l,v,n]$ on $M$ is \emph{torsion-flat} if $\cR^{g , \mA} = 0$, where $g = u\odot v + l\otimes l + n\otimes n$. 
\end{definition}

\noindent
We define similarly the notion of torsion-flat torsion parabolic pairs. Remarkably enough, by Equation \eqref{eq:torsioncurvaturegeneral} it follows that the rank-one invariants of torsion-flat isotropic parallelisms descend to de Rham cohomology, that is:
\begin{eqnarray*}
\mathfrak{M}_o(M) \to H^1(M,\mathbb{R})\times H^1(M,\mathbb{R})\, , \qquad [u,v,l,n] \mapsto ([\kappa],[\rho]) 
\end{eqnarray*}

\noindent
where $\mathfrak{M}_o(M)$ denotes the moduli space of torsion-flat parabolic pairs modulo diffeomorphisms isotopic to the identity. In other words, we have the following result.
\begin{cor}
\label{cor:torsionflat}
A torsion parallel spinor is torsion-flat if and only if its rank-one invariants descend to de Rham cohomology.
\end{cor}

\noindent
Interestingly enough, the torsion-flatness condition on a torsion parallel spinor does not involve explicitly the torsion $\mA \in \Gamma(T^{\ast}M\otimes \wedge^2 M)$. In particular, the study of torsion-flat torsion parallel spinors reduces to the study of the differential system \eqref{eq:duvtorsion} and \eqref{eq:dlntorsion} with $\kappa , \rho \in \Omega^1(M)$ closed. Appropriately tracing Equation \eqref{eq:torsioncurvaturegeneral}, we obtain the Ricci and scalar curvatures of a Lorentzian manifold admitting a torsion parallel spinor as a corollary of Proposition \ref{prop:curvaturetorsion}.

\begin{cor}
\label{cor:curvaturetorsion}
Let $(M,g)$ be a strongly spin Lorentzian four-manifold equipped with a torsion parallel spinor with associated torsion isotropic parallelism $[u,v,l,n]$. Then:
\begin{eqnarray*}
&\Ric^{g,\mA}(w)   =   (\dd\kappa(w,l) + \dd\rho(w,n)) u - \dd\kappa(w,u) l - \dd\rho(w,u) n\\
& s^{g,\mA} = 2(\dd\kappa(u,l) + \dd\rho(u,n))  
\end{eqnarray*}

\noindent
where $w \in \mathfrak{X}(M)$ and $(u,v,l, n) \in [u,v,l,n]$ satisfies the differential system \eqref{eq:existencenullcoframegeneral} relative to $\kappa , \rho \in \Omega^1(M)$.
\end{cor} 

\noindent
Interestingly enough, by the previous corollary the scalar curvature $s^{g,\mA}$ does not depend \emph{explicitly} on the conjugate vector $v\in\Omega^1(M)$. In particular, we obtain the following obstruction for a Lorentzian four-manifold equipped with torsion parallel spinor to be torsion Ricci flat.
\begin{prop}
Let $(M,g)$ be a strongly spin Lorentzian four-manifold equipped with a torsion isotropic parallelism $[u,v,l,n]$ with contorsion $\mA$. Then, $\Ric^{g,\mA} = 0$ only if:
\begin{equation*}
\dd\kappa\vert_{\Ker(u)} = 0\, , \qquad \dd\rho\vert_{\Ker(u)} = 0
\end{equation*}

\noindent
where $(u,v,l, n) \in [u,v,l,n]$ satisfies the differential system \eqref{eq:existencenullcoframegeneral} relative to $\kappa , \rho \in \Omega^1(M)$.
\end{prop} 

\noindent
The previous result can be interpreted more clearly in terms of invariants in \emph{foliated cohomology} associated to every Lorentzian four-manifold equipped with torsion parallel spinors with respect to a Ricci-flat metric connection with torsion. All the results obtained so far hint at torsion parallel spinors on Lorentzian four-manifolds having a remarkably rich geometry which has not been systematically studied in the  mathematical literature, modulo some pioneering exceptions \cite{Galaev:2010jg,Galaev:2009ie}. In later chapters of this dissertation we will focus on supersymmetric NS-NS solutions, which involve a very special class of torsion parallel spinors on Lorentzian four-manifolds. It would be however very interesting to further develop the theory of torsion parallel spinors as an exterior differential system, which is the perspective given by Theorem \ref{thm:existencenullcoframeII}. Note that the system \eqref{eq:duvtorsiongeneral} and \eqref{eq:dlntorsiongeneral} is interpreted as an exterior differential system on the total space of the frame bundle of $M$, in which case both $\kappa$ and $\rho$ should be considered as depending on the coframe $(u,v,l,n)$ chosen and thus depending on \emph{all coordinates} of the frame bundle of $M$, which results in a rich exterior differential system. Alternatively, though not equivalently in general, we can consider tuples of the form $(u,v,l,n,\kappa,\rho)$ as variables of the differential system \eqref{eq:duvtorsiongeneral} and \eqref{eq:dlntorsiongeneral}. In this case, the system is interpreted as an exterior differential system on the bundle:
\begin{equation*}
F(M)\times_M (T^{\ast}M\oplus T^{\ast}M)
\end{equation*}

\noindent
given by the fibered product of the bundle of coframes $F(M)$ of $M$ with the direct sum of $T^{\ast}M$ with itself. Either way, it would be interesting to apply Cartan's involutive test to the exterior differential system defined by \eqref{eq:duvtorsiongeneral} and \eqref{eq:dlntorsiongeneral} and study its prolongations and associated Spencer cohomology \cite{BryantII,BryantBook}.

\renewcommand{\leftmark}{\MakeUppercase{Chapter \thechapter. Supersymmetric Kundt four-manifolds}}

\chapter{Supersymmetric Kundt four-manifolds}
\label{chapter:susyKundt4d}


In this chapter we consider the particular type of torsion parallel spinors that occurs in the supersymmetric configurations and solutions of four-dimensional NS-NS supergravity. This requires introducing the notion of a \emph{abelian gerbe}, which we review in Appendix \ref{chapter:BundleGerbes}, together with other geometric preliminaries needed to establish the mathematical foundations of the bosonic sector of NS-NS supergravity and its four-dimensional Killing spinor equations. The title of this chapter is justified by the fact that, as stated in Corollary \ref{cor:susykundt}, every supersymmetric configuration in NS-NS supergravity is in particular a Kundt Lorentzian four-manifold of special type.


\section{The NS-NS system on a bundle gerbe}
\label{sec:NSNSsystem}

 
In the following $\cC$ will denote a bundle gerbe $\cC := (\cP,\cA,Y)$ \cite{Murray,BehrendXu}, with underlying submersion $Y\to M$, equipped with a fixed connective structure $\cA$ defined on a fixed four-dimensional manifold $M$, and $\cX$ will denote a principal $\mathbb{Z}$ bundle defined on $M$. Given a curving $b\in \Omega^2(Y\times_M Y)$ on $\cC$ we denote its curvature by $H_b\in \Omega^3(M)$. We recall also that an equivariant function $\phi \in \cC^{\infty}(\cX)$ defined on $\cX$ does not necessarily descend to $M$ but its exterior derivative does descend and defines a closed one-form $\varphi_{\phi}\in \Omega^1(M)$ on $M$. Given an equivariant function $\phi \in \cC^{\infty}(\cX)$, we will refer to $\varphi_{\phi}\in \Omega^1(M)$ as its \emph{curvature}. More precisely, we can understand $\mathbb{Z}$-bundles $\cX$
 on $M$ as $(-1)$-gerbes, and we can consider equivariant functions on $\cX$ with equivariance rule:
\begin{eqnarray*}
\phi(x n) = \phi(x) + 2\pi n\, , \quad \forall\,\, x\in \cX\, , \quad \forall\,\, n\in \mathbb{Z}
\end{eqnarray*}

\noindent
as \emph{connections} on the $(-1)$-gerbe $\cX \to M$. Every such equivariant function $\phi$ descends to a $S^1$-valued function $\bar{\phi} \colon M \to S^1 = \mathbb{R}/\mathbb{Z}$. Pulling back the standard volume form on $\mathbb{R}/\mathbb{Z}$ to $M$ gives a closed one-form which coincides with $\varphi_{\phi}$ and which, when it is appropriately normalized, defines an integer class in de Rham cohomology. Hence, our notion of curvature of $\phi$ can be truly understood as the curvature of a \emph{connection} on a $(-1)$-gerbe. 


\subsection{The equations of motion}


A pair $(\cC,\cX)$ determines a unique bosonic NS-NS supergravity on $M$, as prescribed in the following definition. For clarity of exposition, we will omit the term \emph{bosonic} in the following.
 
\begin{definition}
\label{label:NSNSsystem}
The \emph{NS-NS supergravity system}, or \emph{NS-NS system} for short determined by $(\cC,\cX)$ on $M$ is the following differential system \cite{Ortin,Tomasiello}:
\begin{equation}
\label{eq:NSNSsystem}
\Ric^g + \nabla^g \varphi_{\phi} - \frac{1}{2} H_b \circ_g H_b = 0\, , \qquad \nabla^{g\ast}H_b + \varphi_{\phi}\lrcorner_g H_b = 0\, , \qquad    \nabla^{g\ast}\varphi_{\phi} +  \vert \varphi_{\phi}\vert_g^2 = \vert H_b\vert^2_g 
\end{equation}
 	
\noindent
for triples $(g,b,\phi)$ consisting of a Lorentzian metric $g$ on $M$, a curving $b\in \Omega^2(Y)$ on $\cC$, and an equivariant function $\phi\in C^{\infty}(\cX)$, where $H_b\in \Omega^3(M)$ and $\varphi_{\phi}\in \Omega^1(M)$ are the curvatures of $b$ and $\phi$, respectively.
\end{definition}
 
\noindent
By the previous definition we can consider the NS-NS system as a natural \emph{gauge theoretic} system that couples a Lorentzian metric to \emph{connections} $b$ and $\phi$ respectively on a $1$-gerbe, namely a bundle gerbe, and a $(-1)$-gerbe. 

\begin{remark}
Solutions $(g,b,\phi)$ to equations \eqref{eq:NSNSsystem} are \emph{NS-NS solutions}. The first equation in \eqref{eq:NSNSsystem} is the so-called \emph{Einstein equation}, whereas the second equation in \eqref{eq:NSNSsystem} is called the \emph{Maxwell equation} and the third equation in \eqref{eq:NSNSsystem} is referred to as the \emph{dilaton equation} in the literature. In this context, an equivariant function on $\cX$ corresponds to the \emph{dilaton} of the theory, and therefore we will refer to equivariant functions on $\cX$ as \emph{dilatons}, whereas curvings on $\cC$ are usually called \emph{b-fields}, a terminology that we will use occasionally. 
\end{remark}
 
\noindent
Given a NS-NS solution $(g,b,\phi)$, we will refer to the cohomology class $\sigma = [\varphi_{\phi}] \in H^1(M,\mathbb{R})$ determined by $\phi\in \cC^{\infty}(\cX)$ as the \emph{Lee class} of $(g,b,\phi)$. Clearly, different dilatons on $\cX$ define the same Lee class through their curvature. Given $(\cC,\cX)$, we will denote by $\Conf(\cC,\cX)$ the \emph{configuration space} of the NS-NS system on $(\cC,\cX)$, namely the set of triples $(g,b,\phi)$ consisting of a Lorentzian metric $g$ on $M$, a curving $b$ on $\cC$ and an equivariant function $\phi$ on $\cX$. Similarly, we will denote by $\Sol(\cC,\cX)\subset \Conf(\cC,\cX)$ the solution set of NS-NS supergravity on $(\cC,\cX)$. The notion of NS-NS system that we have introduced in Definition \ref{label:NSNSsystem} is valid in any dimension. Since in our case $M$ is four-dimensional, we can simplify the NS-NS system accordingly. Given a curving $b\in \Omega^2(Y\times_M Y)$ with curvature $H_b \in \Omega^3(M)$ we set $H_b = \ast_g \alpha_b$ for a uniquely determine one-form. A quick computation gives the following formulae:
 \begin{equation*}
 H_b\circ_g H_b = \alpha_b\otimes \alpha_b - \vert\alpha_b\vert_g^2 g\, , \quad \vert H_b\vert^2_g  = - \vert \alpha_b\vert^2_g\, , \quad \varphi_{\phi}\lrcorner_g H_b =\ast_g (\alpha_b \wedge \varphi_{\phi})\, , \quad \nabla^{g\ast} H_b = \ast_g \dd\alpha_b 
 \end{equation*}
 
 \noindent
 Hence, the four-dimensional NS-NS system is equivalent to the following system of equations:
 \begin{equation}
 	\label{eq:NSNSsystem4d}
 	\Ric^g + \nabla^g \varphi_{\phi} - \frac{1}{2} \alpha_b \otimes \alpha_b + \frac{1}{2} \vert\alpha_b\vert_g^2 g= 0\, , \qquad \dd \alpha_b = \varphi_{\phi}\wedge \alpha_b \, , \qquad    \nabla^{g\ast}\varphi_{\phi} +  \vert \varphi_{\phi}\vert_g^2 + \vert\alpha_b\vert^2_g = 0
 \end{equation}
 
\noindent
for triples $(g,b,\phi)\in \Conf(\cC,\cX)$. Note that, since $H_b = \ast_g\alpha$ is the curvature of a curving on an abelian gerbe it follows that we always have $\nabla^{g\ast}\alpha = 0$. In the following, when referring to the NS-NS system we will always refer to the differential system \eqref{eq:NSNSsystem4d}.
 
\begin{definition}
A triple $(g,b,\phi)\in \Conf(\cC,\cX)$ is \emph{flux-less} if $b$ is flat, namely if $\alpha_b=0$ identically on $M$, and is \emph{flux} otherwise. A triple $(g,b,\phi)\in \Conf(\cC,\cX)$ is \emph{trivial} if its flux-less and $g$ is flat. 
\end{definition}
 
\noindent
For a trivial triple $(g,b,\phi)\in \Sol(\cC,\cX)$ the NS-NS system reduces to:
\begin{equation*}
\nabla^g\varphi_{\phi} = 0\, , \qquad \vert\varphi_{\phi}\vert_g^2 = 0
\end{equation*}
 
\noindent
and thus it follows that $\varphi_{\phi}$ is a parallel light-like one-form on $M$ and consequently $(M,g,\varphi)$ defines a four-dimensional flat Brinkmann space-time \cite{MehidiZeghib}. Note that this class of space-times can be very \emph{non-trivial} as Lorentzian manifolds \cite{Carriere}. The study of NS-NS solutions can be separated into the study of flux-less NS-NS solutions, namely solutions with flat $b$-field, and \emph{flux} NS-NS solutions, namely solutions for which $\alpha_b \neq 0$ at a point. For flux-less triples $(g,b,\phi)\in \Conf(\cC,M)$ the NS-NS system reduces to:
 
\begin{equation*}
\Ric^g + \nabla^g \varphi_{\phi} = 0\, ,   \qquad    \nabla^{g\ast}\varphi_{\phi} +  \vert \varphi_{\phi}\vert_g^2  = 0
\end{equation*}
 
\noindent
and thus it effectively reduces to differential system for pairs $(g,\varphi)$, where $g$ is a Lorentzian metric on $M$ and $\varphi$ is a closed one-form.  

\begin{remark}
Solutions to the first equation above correspond to a mild generalization of steady Lorentzian Ricci solitons \cite{Gavino} for which the gradient of the function is taken to be the metric dual of a closed, not necessarily exact, one-form.
\end{remark}


\subsection{The supersymmetry conditions}


We have introduced the NS-NS system as the equations of motion of the bosonic sector of NS-NS supergravity in four dimensions. If that was the end of the story, then the NS-NS system would simply be a very particular general relativity model, with a very particular \emph{matter content} and energy momentum tensor. However, and thanks to underlying supersymmetry of the theory, the NS-NS system comes equipped with a first-order \emph{spinorial} differential system that provides \emph{partial integrability} to the full second order NS-NS system and defines the notion of \emph{supersymmetric solution} or \emph{BPS state}. 

Given an element $(g,b,\phi)\in \Conf(\cC,\cX)$ we denote by $\nabla^{g,b}$ the unique metric connection on $(M,g)$ with totally skew-symmetric torsion given by $H_b\in \Omega^3(M)$, namely the curvature of $b\in\Omega^2(Y\times_M Y)$. This is the natural way in which supersymmetry realizes geometrically the notion of torsion for a metric connection. For ease of notation we denote by the same symbol $\nabla^{g,b}$ its lift to any irreducible spinor bundle defined on $(M,g)$.

\begin{definition}
A \emph{supersymmetric configuration} on $(\cC,\cX)$ is triple $(g,b,\phi)$ consisting of a Lorentzian metric $g$ on $M$, an equivariant function $\phi\in \cC^{\infty}(\cX)$, and a curving $b\in \Omega^2(Y)$ satisfying the following differential system:
\begin{equation}
\label{eq:KSE4d}
\nabla^{g,b} \varepsilon = 0\, , \qquad  \varphi_{\phi} \cdot \varepsilon = H_b\cdot \varepsilon
\end{equation}
	
\noindent
for a nowhere vanishing section $\varepsilon\in\Gamma(S)$ of an irreducible paired spinor bundle $(S,\Gamma,\cB)$  on $(M,g)$. Such spinor $\varepsilon\in \Gamma(S)$ is a \emph{supersymmetry parameter} or \emph{supersymmetry generator} for the supersymmetric configuration $(g,\varphi,b)$.
\end{definition}

\begin{remark}
Note that, using the rigorous notation introduced in Chapter \ref{chapter:spingeometryClifford} the second equation in \eqref{eq:KSE4d} would read as follows:
\begin{equation*}
	\Psi_{\Gamma}(\varphi_{\phi}) (\varepsilon) = \Psi_{\Gamma}(H_b)(\varepsilon)
\end{equation*}

\noindent
It is however standard to denote Clifford multiplication with just a \emph{dot}.
\end{remark}

\noindent
Therefore, the supersymmetry parameter of a supersymmetric configuration is in particular a torsion parallel spinor with fully skew-symmetric torsion. We denote by $\Conf_s(\cC,\cX)$ the set of supersymmetric configurations on $(\cC,\cX)$, whose elements consist by definition of tuples $(g,b,\phi,\varepsilon)$ satisfying the differential system \eqref{eq:KSE4d}. Because of their physical origin in the supergravity literature, the first equation in \eqref{eq:KSE4d} is called the \emph{gravitino equation}, whereas the second equation in \eqref{eq:KSE4d} is called the \emph{dilatino equation}. We will occasionally use this terminology in the following. Since supersymmetric NS-NS configurations involve torsion parallel spinors with completely skew-symmetric torsion, we consider the latter separately in the following subsection.


\section{Skew-symmetric torsion parallel spinors}
\label{sec:sktorsionparallelspinors}


We say that a metric connection $\nabla = \nabla^g + \mA$ on $(M,g)$ has completely skew-symmetric torsion, or that is \emph{skew} for short, if $\mA \in \Omega^3(M)\subset \Gamma(T^{\ast}M\otimes \wedge^2 T^{\ast}M)$ is a three-form, namely, it is skew-symmetric in all of its entries. Introducing the three-form $H\in\Omega^3(M)$ as:
\begin{eqnarray*}
	H = 2 \mA \in \Omega^3(M)
\end{eqnarray*}

\noindent
we write $\nabla^{g,H} := \nabla$. It follows that $\nabla^{g,H}$ has torsion precisely $H$:
\begin{equation*}
	\nabla^{g,H}_{w_1}w_2 - \nabla^{g,H}_{w_2}w_1 - [w_1,w_2] = H(w_1,w_2)\, , \quad \forall\,\, w_1, w_2  \in \mathfrak{X}(M)\, .
\end{equation*}

\noindent
In our conventions the curvature tensor of $\nabla^{g,H}$, which we denote by $\cR^{g,H}$, is given by:
\begin{equation*}
	\cR^{g,H}_{w_1 , w_2} w_3 = \nabla^{g,H}_{w_1}\nabla^{g,H}_{w_2}w_3 - \nabla^{g,H}_{w_2}\nabla^{g,H}_{w_1}w_3 - \nabla^{g,H}_{[w_1 , w_2]}w_3\, , \quad \forall\,\, w_1, w_2, w_3 \in \mathfrak{X}(M)
\end{equation*}

\noindent
and expands as follows in terms of the torsion $H$ and its covariant derivative:
\begin{eqnarray*}
	& \cR^{g,H}_{w_1 , w_2} w_3 = \cR^{g}_{w_1 , w_2} w_3 + \frac{1}{2}(\nabla^g_{w_1}H)(w_2,w_3)  - \frac{1}{2}(\nabla^g_{w_2}H)(w_1,w_3) \\
	& + \frac{1}{4} H(w_1,H(w_2,w_3)) - \frac{1}{4} H(w_2,H(w_1,w_3)) \, , \quad \forall\,\, w_1, w_2, w_3 \in \mathfrak{X}(M)
\end{eqnarray*}

\noindent
where $\cR^{g}$ denotes the Riemann tensor of $g$. The associated Ricci and scalar curvatures, which we denote by $\mathrm{Ric}^{g,H}$ and $s^{g,H}$ respectively, are given by:
\begin{equation*}
	\mathrm{Ric}^{g,H} = \mathrm{Ric}^{g} - \frac{1}{2}\nabla^{g\ast}H - \frac{1}{2} H\circ_g H \, , \qquad  s^{g,H} = s^g - \frac{3}{2} \vert H\vert^2_g
\end{equation*}

\noindent
Since we consider $M$ to be oriented and four-dimensional, we can \emph{replace} $H$ by its Hodge dual $\alpha = \ast_g H$, in which case we will sometimes write $\nabla^{g,\alpha}$ for unique metric connection with totally skew torsion $H$. This is convenient for computations as well as for the geometric understanding of torsion in four Lorentzian dimensions, and allows for introducing the notion of the \emph{causal character} of $H$, which will be useful in the following. 

\begin{definition}
	Let $\nabla^{g,H}$ be a metric connection with torsion and $O\subset M$ a subset of $M$. The torsion of $\nabla^{g,H}$ is \emph{time-like}, \emph{null}, or \emph{space-like} on $O$ if $\vert\alpha\vert^2_g\vert_O <0$, $\vert\alpha\vert^2_g\vert_O =0$ and $\vert\alpha\vert^2_g\vert_O >0$, respectively.
\end{definition}

\begin{remark}
\label{remark:curvatureskewtorsion}
The curvature tensors $\cR^{g,H}$, $\mathrm{Ric}^{g,H}$ and $s^{g,H}$ can be written as follows in terms of the skew-torsion $\alpha = \ast_g H$:
\begin{eqnarray*}
	& \cR^{g,\alpha}_{w_1 , w_2} = \cR^{g,H}_{w_1 , w_2}   = \cR^{g}_{w_1 , w_2} + \frac{1}{2} \ast_g ( \nabla^g_{w_1}\alpha\wedge w_2 - \nabla^g_{w_2}\alpha\wedge w_1) \\
	& + \frac{1}{4} (\alpha(w_1) \alpha \wedge w_2 - \alpha(w_2) \alpha \wedge w_1 - \vert\alpha\vert_g^2 w_1\wedge w_2) \, , \quad \forall\,\, w_1, w_2  \in \mathfrak{X}(M)\\
	& \mathrm{Ric}^{g,H} = \mathrm{Ric}^{g} - \frac{1}{2}\ast_g \dd\alpha + \frac{1}{2} \vert\alpha\vert_g^2  g - \frac{1}{2}\alpha\otimes\alpha \, , \qquad  s^{g,\alpha} = s^g + \frac{3}{2} \vert \alpha\vert^2_g
\end{eqnarray*}

\noindent
We will use this formulae later in this section to compute the curvature of a supersymmetric NS-NS configuration. 
\end{remark}

\begin{definition}
\label{def:spinorswithtorsion}
Let $(M,g)$ be a strongly spin Lorentzian four-manifold. A \emph{skew-torsion parallel spinor} on $(M,g)$ with torsion $H\in\Omega^3(M)$ is a section $\varepsilon\in \Gamma(S)$ of a paired spinor bundle $(S,\Gamma,\cB)$ over $(M,g)$ that is parallel with respect to $\nabla^{g,H}$.
\end{definition}

\noindent
Proposition \ref{prop:differentialspinors4d} immediately implies the following characterization of skew-torsion parallel spinors on Lorentzian four-manifolds.
\begin{prop}
\label{prop:skewtorsionspinors4d}
A strongly spin Lorentzian four-manifold $(M,g)$ admits a skew-torsion parallel spinor with torsion $H\in \Omega^3(M)$ if and only if it admits a parabolic pair $(u,[l]_u)$ satisfying:
\begin{equation}
\label{eq:TorsionSpinorParabolicPair}
\nabla^g u = \frac{1}{2} \ast_g(\alpha\wedge u)	 \, , \qquad   \nabla^g l = \kappa\otimes u + \frac{1}{2} \ast_g( \alpha \wedge l)
\end{equation}
	
\noindent 
for any, and hence for all, representatives $l\in [l]_u$.
\end{prop}

\noindent
In the following will refer to the parabolic pair associated to a spinor parallel with respect to a metric connection with skew torsion simply as a \emph{skew-torsion parabolic pair}. The torsion $H$ of the connection that preserves the spinor will be referred to as the torsion of the associated skew-torsion parabolic pair.

\begin{definition}
A \emph{skew-torsion parabolic structure} on $M$ is a tuple $(g,(u,[l]_u),\alpha)$ consisting of a Lorentzian metric $g$ on $M$ and skew-torsion parabolic pair $(u,[l]_u)$ on $M$ with torsion $H = \ast_g\alpha$. A \emph{skew-torsion parabolic four-manifold} consists of a four-manifold equipped with a skew-torsion parabolic structure.
\end{definition}

\noindent
By Proposition \ref{prop:skewtorsionspinors4d} a skew-torsion Lorentzian four-manifold  $(M,g,(u,[l]_u),\alpha)$ can be equivalently considered as a tuple $(M,g,\varepsilon,\alpha)$, where $\varepsilon\in \Gamma(S)$ is a skew-torsion spinor with respect to $H = \ast_g \alpha$. Three-dimensional skew-torsion Lorentzian manifolds have been already considered in \cite{Shahbazi3d} and, in contrast to the four-dimensional case considered in this dissertation, in three dimensions they can be equivalently studied exclusively in terms of their associated Dirac current.

\begin{remark}
For further reference, we denote by $\mathfrak{P}_{sk}(M)$ the full subcategory of $\mathfrak{P}(M)$ whose objects are skew-torsion parabolic structures $(g,u,[l]_u,\alpha)$, and we denote by $\mathfrak{P}^{\alpha}_{sk}(M)$, where $\alpha\in \Omega^1(M)$ the full subcategory of $\mathfrak{P}_{sk}(M)$ whose objects all have torsion $H = \ast_g \alpha$. 
\end{remark}

\noindent
From Equations \eqref{eq:TorsionSpinorParabolicPair} it immediately follows that if $(u,[u]_l)\in \mathfrak{P}^{g}_{sk}(M)$ then $u^{\sharp_g}\in\mathfrak{X}(M)$ is Killing and its integral curves are geodesics, that is:
\begin{equation*}
	\nabla_u^{g} u = 0
\end{equation*}

\noindent
Hence, every Lorentzian four-manifold $(M,g)$ admitting a parallel spinor with torsion comes equipped with a geodesic null congruence generated by the Dirac current of the latter. In particular, the optical invariants associated to this congruence read as follows:
\begin{equation*}
	\theta := \frac{1}{2}\mathrm{Tr}_g(\nabla^g u) = 0\, , \qquad 4 \omega^2 := \vert \dd u \vert^2_g = -  \alpha(u)^2\, , \qquad \sigma^2 := \frac{1}{8} \vert \cL_u g\vert^2_g -\theta^2 = 0
\end{equation*}

\noindent
In the terminology of mathematical general relativity, $\theta$ is the \emph{expansion}, $\omega$ the \emph{twist} and $\sigma$ the \emph{shear} of the given null congruence. The previous expression for the expansion of the null congruence defined by $u$ together with Proposition \ref{prop:adaptedKillingintegrable} yields the following result.
\begin{cor}
	Let $(M,g,(u,[l]_u),\alpha)$ be a skew-torsion Lorentzian four-manifold. Then the twist of $u$ is zero if and only if $\alpha(u) = 0$, if and only if $u^{\perp_g}\subset TM$ is integrable, if and only if $(M,g,u)$ is Kundt.
\end{cor}

\noindent
The previous corollary suggests that skew-torsion Lorentzian four-manifolds are not necessarily Kundt, in contrast to the three-dimensional case considered in \cite{Shahbazi3d}. As a direct application of Proposition \ref{prop:existencenullcoframegeneral} we obtain the following characterization of skew-torsion parabolic structures.

\begin{prop}
\label{prop:existencenullcoframe}
An oriented Lorentzian four-manifold $(M,g)$ admits a skew-torsion parabolic pair $(u,[l]_u)$ if and only if it there exists null coframe $(u,v,l,n)\in \mathbb{E}^{-1}(g,u,[l]_u)$ satisfying the following differential system:
\begin{eqnarray}
& \nabla^g u =  \frac{1}{2} \ast_g(\alpha\wedge u) \, , \qquad \nabla^g v = - \kappa\otimes l - \rho\otimes n + \frac{1}{2}\ast_g (\alpha\wedge v) \label{eq:nablatorsion1}\\
& \nabla^g l = \kappa\otimes u + \frac{1}{2} \ast_g(\alpha\wedge l) \label{eq:nablaltorsion}\, , \qquad \nabla^g n = \rho\otimes u + \frac{1}{2} \ast_g(\alpha\wedge n) \label{eq:nablatorsion2}
\end{eqnarray}
	
	\noindent
	for a given pair of one-forms $\kappa , \rho \in \Omega^1(M)$.
\end{prop}

\noindent
Similarly, as a direct application of Theorem \ref{thm:existencenullcoframeII} we obtain the following equivalent characterization of Lorentzian four-manifolds equipped with skew-torsion parallel spinors.
\begin{prop}
	\label{prop:existencenullcoframeII}
	A strongly spin four-manifold $M$ admits a skew-torsion parallel spinor if and only if it there exists a global coframe $(u,v,l,n)$ and a pair of one-forms $\kappa , \rho \in \Omega^1(M)$ satisfying the following exterior differential system:
	\begin{eqnarray}
		& \dd u =    \ast_g(\alpha\wedge u)  \, , \qquad \dd v = - \kappa\wedge l - \rho\wedge n + \ast_g (\alpha\wedge v) \label{eq:duvtorsion}\\
		& \dd l = \kappa\wedge u + \ast_g(\alpha\wedge l) \, , \qquad \dd n = \rho\wedge u + \ast_g(\alpha\wedge n) \label{eq:dlntorsion}
	\end{eqnarray}
	
	\noindent
	where $g = u\odot v + l\otimes l + n\otimes n$ is the Lorentzian metric associated to $[u,v,l,n]$.
\end{prop}

\begin{cor}
A strongly spin four-manifold $M$ admits a skew-torsion parallel spinor if and only if it there exists a null coframe satisfying \eqref{eq:duvtorsion} and \eqref{eq:dlntorsion}.
\end{cor}

\noindent
We denote by $\mathfrak{F}_{sk}(M)$ the full subcategory of $\mathfrak{F}(M)$ whose objects $[u,v,l,n]$ satisfy the differential system \eqref{eq:duvtorsion} and \eqref{eq:dlntorsion} for any and hence all representatives $(u,v,l,n)\in [u,v,l,n]$. We will refer to such $[u,v,l,n]$ as skew-torsion isotropic parallelisms, and to any of the representatives $(u,v,l,n)\in [u,v,l,n]$ as \emph{skew-torsion isotropic coframes} when needed. The full subcategory $\mathfrak{F}^{\alpha}_{sk}(M)$ is defined analogously to $\mathfrak{P}^{\alpha}_{sk}(M)$.

\begin{remark}
\label{remark:Liebrackets}
We recall below the Lie brackets of the dual of a skew-torsion isotropic coframe. 
	
\begin{eqnarray*}
& \left[ u^{\sharp_g} , v^{\sharp_g} \right]= - (\alpha(n) + \kappa(u))\, l^{\sharp_g} + (\alpha(l) - \rho(u))\, n^{\sharp_g}\\ 
& \left[ u^{\sharp_g} , l^{\sharp_g}\right] = (\kappa(u) + \alpha(n))\, \mu^{\sharp_g} - \alpha(u)\, n^{\sharp_g} \\
& \left[ u^{\sharp_g} , n^{\sharp_g}\right] = (\rho(u) - \alpha(l))\, \mu^{\sharp_g} + \alpha(u)\, l^{\sharp_g}\\ 
& \left[ v^{\sharp_g} , l^{\sharp_g}\right] = \kappa(v) u - \alpha(n) v  + \kappa(l) l + (\rho(l) + \alpha(v)) n\\
& \left[ v^{\sharp_g} , n^{\sharp_g}\right] =  \alpha(l) v + \rho(v) u + (\kappa(n) - \alpha(v)) l + \rho(n) n\\
& \left[ l^{\sharp_g} , n^{\sharp_g}\right] = (\rho(l) - \kappa(n) + \alpha(v)) u  - \alpha(u) v    
\end{eqnarray*}
	
\noindent
We will use these expressions occasionally in the following.
\end{remark}

\noindent
Proposition \ref{prop:existencenullcoframeII} allows to give an equivalent point of view on skew-torsion parabolic pairs which is very convenient both conceptually and computationally.   

\begin{prop}
There is a canonical equivalence of categories between $\mathfrak{P}_{sk}(M)$ and $\mathfrak{F}_{sk}(M)$ given by the restriction of $\mathbb{E}\colon \mathfrak{P}(M) \to \mathfrak{F}(M)$. 
\end{prop}

\noindent
Another direct consequence of Proposition \ref{prop:existencenullcoframe} is the following uniqueness result with respect to the torsion of a skew-torsion parabolic pair.

\begin{cor}
	\label{cor:equaltorsion}
	Let $(u,[l]_u)$ be a skew-torsion parabolic pair such that $(g,u,[l]_u)\in \mathfrak{P}^{\alpha_1}_{sk}(M)$ and $(g,u,[l]_u)\in \mathfrak{P}^{\alpha_2}_{sk}(M)$, where $\alpha_1 , \alpha_2 \in \Omega^1(M)$. Then $\alpha_1 = \alpha_2$.
\end{cor}

\begin{proof}
	Suppose that $(g,u,[l]_u)\in \mathfrak{P}^{\alpha_1}_{sk}(M)$ and $(g,u,[l]_u)\in \mathfrak{P}^{\alpha_2}_{sk}(M)$. Then, there exist $(u,v,l,n) \in \mathbb{E}^{-1}(g,u,[l]_u)$ satisfying equations \eqref{eq:duvtorsion} and \eqref{eq:dlntorsion} with respect to both torsions $\alpha_1$ and $\alpha_2$ and relative to pairs of one-forms $(\kappa_1,\rho_1)$ and $(\kappa_2,\rho_2)$, respectively. Combining these systems of equations, we obtain
	\begin{eqnarray*}
		&    u \wedge (\alpha_1 - \alpha_2) = 0  \, , \qquad   - (\kappa_1 - \kappa_2)\wedge l - (\rho_1 - \rho_2)\wedge n + \ast_g (v\wedge (\alpha_1 - \alpha_2)) = 0  \\
		&  (\kappa_1 - \kappa_2)\wedge u + \ast_g(l \wedge ( \alpha_1 - \alpha_2)) = 0  \, , \qquad   (\rho_1 -\rho_2)\wedge u + \ast_g(n \wedge (\alpha_1 - \alpha_2)) = 0  
	\end{eqnarray*}
	
	\noindent
	Hence $\alpha_1 = \alpha_2 + f_1 u$ for a function $f_1\in C^{\infty}(M)$. Plugging this equation into the remaining equations above, we obtain:
	\begin{eqnarray*}
		&       (\kappa_1 - \kappa_2)\wedge l + (\rho_1 - \rho_2)\wedge n = f_1 l\wedge n  \\
		&  (\kappa_1 - \kappa_2)\wedge u + f_1 n\wedge u = 0  \, , \qquad   (\rho_1 -\rho_2)\wedge u  - f_1 l\wedge u = 0  
	\end{eqnarray*}
	
	\noindent
	Therefore, from the second line we obtain:
	\begin{eqnarray*}
		\kappa_1 = \kappa_2 - f_1 n + f_2 u\, , \quad \rho_1 = \rho_2 + f_1 l + f_3 u
	\end{eqnarray*}
	
	\noindent
	in terms of functions $f_2, f_3\in C^{\infty}(M)$. Plugging these expressions into the first line above, we obtain:
	\begin{equation*}
		(f_2 u - f_1\, n)\wedge l + ( f_1 l + f_3 u)\wedge n = f_1 l\wedge n 
	\end{equation*}
	
	\noindent
	Hence $f_1 = f_2 = f_3 = 0$ and thus $\alpha_1 = \alpha_2$. 
\end{proof}

\noindent
Hence, a skew-torsion parabolic structure, equivalentely, a skew-torsion isotropic parallelism on a given four-manifold can only be so with respect to a unique completely skew-torsion tensor $H=\ast_g \alpha\in \Omega^1(M)$. In particular, we can refer to skew-torsion parabolic pairs without explicit mention of their torsion. Hence, we obtain a disjoint partition into full subcategories:
\begin{eqnarray*}
	\mathfrak{P}_{sk}(M) = \bigcup_{\alpha\in\Omega^1(M)} \mathfrak{P}^{\alpha}_{sk}(M) \, , \qquad \mathfrak{F}_{sk}(M) = \bigcup_{\alpha\in\Omega^1(M)} \mathfrak{F}^{\alpha}_{sk}(M)
\end{eqnarray*}

\noindent
and a natural map:
\begin{equation*}
	\mathfrak{P}_{sk}(M) \to \Omega^1(M)\, , \qquad (u,[l]_u) \mapsto \alpha
\end{equation*}

\noindent
and similarly for $\mathfrak{F}_{sk}(M)$. This maps descends to the moduli space $\mathfrak{M}_{sk}(M)$:
\begin{equation*}
	\mathfrak{M}_{sk}(M) \to \frac{\Omega^1(M)}{\Diff(M)}\, , \qquad [(u,[l]_u)] \mapsto [\alpha]
\end{equation*}

\noindent
It would be interesting to further elucidate the basic properties of this map, which gives a fibration onto its image with fiber given by the set of skew-torsion parabolic pairs with fixed torsion. We end this section with a nice corollary that follows immediately from Remark \ref{remark:Liebrackets} or, alternatively, Proposition \ref{prop:existencenullcoframeII}. As explained in Chapter \ref{chapter:IrreducibleSpinors4d}, every skew-torsion parabolic structure $(g,u,[l]_u)$ defines a screen bundle $\mathfrak{G}_u$ associated to its Dirac current $u\in\Omega^1(M)$. A choice of conjugate vector field to $u$ defines a rank-two distribution in $TM$ that is isomorphic to $\mathfrak{G}_u$. It is then natural to ask if there exists a choice of conjugate vector field to $u$ such that the rank-two distribution that it defines is integrable. This leads us to the following result.

\begin{cor}
Let $[u,v,l,n]$ be a skew-torsion isotropic parallelism on $M$. There exists an integrable realization of the screen bundle determined by $u$ if and only if there exists a null coframe $(u,v,l,n)\in [u,l,v,n]$ satisfying:
\begin{eqnarray*}
\alpha(u) = 0\, , \qquad    \alpha(v)  =    \kappa(n)   - \rho(l)
\end{eqnarray*}
	
\noindent
In particular $\Ker(u)\subset TM$ is integrable and $\dd l = \dd n = 0$.
\end{cor}

\noindent
The foliation determined by an integrable screen bundle corresponds to the \emph{wave front} of $(M,g)$ when the latter is interpreted as a \emph{gravitational wave}. Hence, and interestingly enough, the previous corollary implies that the gravitational waves defined by skew-torsion parallel spinors have flat wave fronts when the latter assemble into a codimension-two foliation.


\subsection{Natural reduction to a Riemann surface}


In this section we will use the theory of isotropic parallelisms to study a class of skew-torsion parallel spinors that reflects naturally the local structure determined by the existence of an isotropic isometry, which by Proposition \ref{prop:existencenullcoframe} always exists on Lorentzian four-manifolds equipped with a skew-torsion parallel spinor and is given by its Dirac current $u^{\sharp_g}\in \mathfrak{X}(M)$. Throughout this section $X$ will denote a two-dimensional oriented manifold.  

\begin{definition}
An isotropic parallelism $[u,v,l,n]$ defined on $M = \mathbb{R}^2 \times X$ is \emph{adapted} if $u^{\sharp_g} = \partial_{x_v}$, where $(x_u,x_v)$ are the Cartesian coordinates of the $\mathbb{R}^2$ factor.  
\end{definition} 

\begin{definition}
A nowhere vanishing spinor $\varepsilon$ on $(\mathbb{R}^2\times X, g)$ is adapted if its associated isotropic parallelism is adapted. 
\end{definition}

\noindent
By construction, every skew-torsion Lorentzian four-manifold is locally isomorphic to an adapted skew-torsion parabolic pair on $M = \mathbb{R}^2\times X$.  

\begin{lemma}
\label{lemma:naturalrepresentative}
Let $[u,v,l,n]$ be an adapted isotropic parallelism. Then, there exists a representative $(u,v,l,n)\in [u,v,l,n]$ given by:
\begin{equation*}
(u,v,l,n) = (e^{\cF} \dd x_u + \theta , \dd x_v +  e^{\cK}  \dd x_u + \omega, l^{\perp} , n^{\perp})
\end{equation*}
	
\noindent
where $\cF, \cK \in C^{\infty}(\mathbb{R}^2\times \Sigma)$ are functions on $\mathbb{R}^2\times X$ and $\theta , \omega , l^{\perp}$ and $n^{\perp}$ are sections of the pull-back of $T^{\ast}X$ along the canonical projection $\mathbb{R}^2 \times \Sigma \to \Sigma$. 
\end{lemma}

\begin{proof}
Given any representative $(u,v,l,n)\in [u,v,l,n]$, we expand each of its elements in a basis given by $\dd x_u$, $\dd x_v$ and its projection to $T^{\ast} X$ using the fact that $T^{\ast} M = T^{\ast}\mathbb{R}^2 \oplus T^{\ast}X$. Since $u^{\sharp_g} = \partial_{x_v}$ and $u(u^{\sharp_g}) = u(\partial_{x_v}) = 0$, the expansion of $u$ in the aforementioned basis cannot have any $\dd x_v$ term, and similarly for $l$ and $n$. On the other hand, since $v(u^{\sharp_g}) = v(\partial_{x_v}) = 1$ the coefficient of the $\dd x_v$ term in the expansion of $v$ needs to be 1. Then, performing a \emph{gauge transformation} generated by $w = l(\partial_{x_u}) l + n(\partial_{x_u}) n$ we obtain the expression in the statement of the lemma.
\end{proof}

\noindent
We will refer to such a $(u,v,l,n) \in [u,v,l,n]$ as the \emph{natural} representative of the adapted isotropic parallelism $[u,v,l,n]$. Given a natural representative $(u,v,l,n) \in [u,v,l,n]$ as in Lemma \ref{lemma:naturalrepresentative}, we denote by:
\begin{equation*}
q = l^{\perp} \otimes l^{\perp} + n^{\perp}\otimes n^{\perp}
\end{equation*} 

\noindent
the induced Riemannian metric on $X$. Note  that this is well-defined since $l  , n  \in \Omega^1(X)$ are linearly independent for every $x_u \in \mathbb{R}$ and $x_v \in \mathbb{R}$. Using this induced metric on $X$, we can neatly express the metric dual $(u^{\sharp_g},v^{\sharp_g},l^{\sharp_g},n^{\sharp_g})$ of a natural representative $(u,v,l,n)$ in terms of the metric dual with respect to $q_{x_u}$ on $X$, which will be useful in later computations. We obtain:
\begin{eqnarray*}
	& u^{\sharp_g} = \partial_{x_v}\, , \qquad v^{\sharp_g} = e^{-\cF} \partial_{x_u}-  e^{\cK-\cF}\partial_{x_v}\\ 
	& l^{\sharp_g} = - e^{-\cF} \theta(l^{\perp}) \partial_{x_u} + (e^{\cK-\cF} \theta(l^{\perp}) - \omega(l^{\perp})) \partial_{x_v} + (l^{\perp})^{\sharp_{q}}\\
	&  n^{\sharp_g} = - e^{-\cF} \theta(n^{\perp}) \partial_{x_u} + (e^{\cK - \cF} \theta(n^{\perp}) - \omega(n^{\perp})) \partial_{x_v} + (n^{\perp})^{\sharp_{q}}
\end{eqnarray*}

\noindent
In particular, note that a priori $[l^{\sharp_g},n^{\sharp_g}] \neq 0$ and the span of $l^{\sharp_g},n^{\sharp_g} \subset \mathfrak{X}(\mathbb{R}^2\times \Sigma)$ is a rank-two distribution in $TM$ that in general has non-trivial intersection with both factors $T\mathbb{R}^2$ and $TX$ in the splitting  $TM = T\mathbb{R}^2\oplus TX$. On the other hand, the span of $u^{\sharp_g},v^{\sharp_g} \subset \mathfrak{X}(\mathbb{R}^2\times \Sigma)$ is precisely the tangent bundle factor $T\mathbb{R}^2$, although again in general we will have $[u^{\sharp_g},v^{\sharp_g}] \neq 0$. Evaluated at an adapted isotropic coframe, the differential system \eqref{eq:duvtorsion} and \eqref{eq:dlntorsion} can be equivalently written as follows:
\begin{eqnarray}
& \dd u =    \alpha(u) \nu_{q} - u \wedge\ast_{q} \alpha^{\perp}  \, , \,\, \dd v =  - \alpha(v) \nu_{q} + v \wedge\ast_{q} \alpha^{\perp} - \kappa\wedge l - \rho\wedge n  \label{eq:torsion2d} \\
& \dd l = \kappa\wedge u - (\alpha(u) v - \alpha(v) u) \wedge n - \alpha(n) u\wedge v \, , \,\, \dd n = \rho\wedge u + (\alpha(u) v - \alpha(v) u) \wedge l + \alpha(l) u\wedge v \nonumber
\end{eqnarray}

\noindent
where $\alpha^{\perp}\in\Omega^1(X)$ is the natural projection of $\alpha$ to $T^{\ast}X$, $\nu_{q} = l^{\perp}\wedge n^{\perp}$ is the induced Riemannian volume form and $\ast_{q} \colon \Omega^{\bullet}(X) \to \Omega^{\bullet}(X)$ is the Hodge operator defined by $q$ on $X$. On the other hand, the exterior derivative of an adapted isotropic coframe is readily found to be:
\begin{eqnarray}
& \dd u  = e^{\cF} \dd_X \cF \wedge \dd x_u + e^{\cF} \partial_{x_v}\cF\, \dd x_v \wedge \dd x_u + \dd_X \theta + \dd x_u \wedge \partial_{x_u}\theta + \dd x_v \wedge \partial_{x_v}\theta \nonumber \\
& \dd v =  e^{\cK} \dd_X \cK \wedge \dd x_u + e^{\cK} \partial_{x_v} \cK\,  \dd x_v \wedge \dd x_u + \dd_X \omega + \dd x_u \wedge \partial_{x_u}\omega + \dd x_v \wedge \partial_{x_v}\omega \label{eq:exteriorparallelism2d}\\
& \dd l  =   \dd_X l^{\perp}  + \dd x_u \wedge \partial_{x_u} l^{\perp} + \dd x_v \wedge \partial_{x_v} l^{\perp} \, ,\quad \dd n =    \dd_X n^{\perp}  + \dd x_u \wedge \partial_{x_u}n^{\perp}  + \dd x_v \wedge \partial_{x_v}n^{\perp} \nonumber
\end{eqnarray}

\noindent
where $\dd_X\colon \Omega^{\bullet}(X) \to \Omega^{\bullet}(X)$ denotes the exterior derivative on $X$.

\begin{lemma}
\label{lemma:zeroconditions}
Let $(u,v,l,n)$ be adapted isotropic coframe with skew-torsion relative to one-forms $\kappa, \rho \in \Omega^1(\mathbb{R}^2\times X)$. Then:
\begin{eqnarray*}
& \partial_{x_v} \cF = \partial_{x_v} \cK = 0\, , \qquad \partial_{x_v} \theta = \partial_{x_v} \omega = 0\, , \qquad \kappa_v = - \alpha^{\perp}(n^{\perp})\, , \qquad \rho_v =  \alpha^{\perp}(l^{\perp})\\
& \partial_{x_v} l^{\perp} =  - \alpha_v n^{\perp}  \, , \qquad \partial_{x_v} n^{\perp} =   \alpha_v l^{\perp}
\end{eqnarray*}
	
\noindent
where we have set $\alpha = \alpha_u \dd x_u + \alpha_v \dd x_v + \alpha^{\perp}$, and similarly for $\kappa$ and $\rho$ in the splitting given by $T^{\ast}M = T^{\ast}\mathbb{R}^2 \oplus T^{\ast}X$.
\end{lemma}

\begin{proof}
Plugging equations \eqref{eq:exteriorparallelism2d} into the differential system \eqref{eq:torsion2d} and evaluating it on the vector field $\partial_{x_v}$ we obtain:
\begin{eqnarray*}
\partial_{x_v}e^{\cF} \dd x_u  + \partial_{x_v}\theta = 0\, , \qquad  \partial_{x_v}e^{\cK} \dd x_u + \partial_{x_v}\omega = \ast_{q} \alpha^{\perp}- \kappa_v l^{\perp} - \rho_v n^{\perp}\\
\partial_{x_v} l^{\perp} = \kappa_v u - \alpha_v n^{\perp} + \alpha(n) u\, , \qquad \partial_{x_v} n^{\perp} = \rho_v u + \alpha_v l^{\perp} - \alpha(l) u
\end{eqnarray*}

\noindent
Hence, we readily conclude that:
\begin{equation*}
\partial_{x_v} \cF = 0\, , \quad \partial_{x_v} \cK = 0\, , \quad \partial_{x_v} \theta = 0\, , \quad \kappa_v = - \alpha^{\perp}(n^{\perp})\, , \quad \rho_v =  \alpha^{\perp}(l^{\perp})
\end{equation*}

\noindent
which in turn implies the following conditions:
\begin{equation*}
\partial_{x_v}\omega = \ast_{q} \alpha^{\perp}- \kappa_v l^{\perp} - \rho_v n^{\perp} = 0\, , \qquad \partial_{x_v} l^{\perp} =  - \alpha_v n^{\perp}  \, , \qquad \partial_{x_v} n^{\perp} =   \alpha_v l^{\perp}
\end{equation*}

\noindent
and hence we conclude.
\end{proof}

\noindent
By the previous lemma it follows that an adapted isotropic coframe $(u,v,l,n)$ only depends on the Cartesian coordinate $x_v$ through $(l^{\perp},n^{\perp})$. Nonetheless, the \emph{transverse metric} $q$ constructed in terms of $(l^{\perp},n^{\perp})$ does not depend on $x_v$, that is:
\begin{eqnarray*}
\cL_{\partial_{x_v}}q = \partial_{x_v} q = l^{\perp}\odot \partial_{x_v}l^{\perp}  + n^{\perp}\odot \partial_{x_v}n^{\perp} = - \alpha_v  n^{\perp}  \odot \partial_{x_v}l^{\perp} + \alpha_v n^{\perp}\odot l^{\perp} = 0
\end{eqnarray*}

\noindent
which immediately implies that $\partial_{x_v}$ is an isometry of $g = u\odot v + q$, as expected. Hence, the \emph{evolution} of $(l^{\perp},n^{\perp})$ in $x_v$ preserves the isometry type of $(X,q)$. For the remainder of this subsection we will assume that $\partial_{x_v}l^{\perp} = \partial_{x_v}n^{\perp} = 0$ whence the adapted isotropic coframe $(u,v,l,n)$ is independent of $x_v$. We will refer to such adapted isotropic coframes as \emph{invariant}. By Lemma \ref{lemma:zeroconditions} this occurs if and only if $\alpha_v = 0$, which we will assume from now on. We will encounter this condition again in Section \ref{sec:susyconf} when we consider supersymmetric configurations and solutions. Below we write $\alpha^{\perp} = \alpha^{\perp}_l l + \alpha^{\perp}_n n$ for ease of notation.

\begin{lemma}
\label{lemma:reducedinterequations}
An invariant adapted isotropic coframe $(u,v,l,n)$ is skew-torsion if and only if the following differential system is satisfied:
\begin{eqnarray*}
&  \partial_{x_u} \theta = \dd_X e^{\cF} - e^{\cF} \ast_{q} \alpha^{\perp} \, , \quad    \partial_{x_u}\omega =  \dd_X e^{\cK} + e^{\cK} \ast_{q} \alpha^{\perp} -\kappa_u l^{\perp} - \rho_u n^{\perp}\\
& \partial_{x_u} l =   \alpha_u  n^{\perp} - e^{\cF}\kappa^{\perp} + (\kappa_u + \alpha^{\perp}_n e^{\cK}) \theta - \alpha^{\perp}_n e^{\cF} \omega \\
& \partial_{x_u} n =  - \alpha_u  l^{\perp} + \alpha^{\perp}_l e^{\cF} \omega - e^{\cF}\rho^{\perp} + (\rho_u - \alpha^{\perp}_l e^{\cK}) \theta\\
& \dd_X \theta =  -  \theta\wedge \ast_{q}\alpha^{\perp} \, , \quad  \dd_X\omega  =  - \alpha_u   e^{-\cF}   \nu_{q} + \omega \wedge \ast_{q} \alpha^{\perp} - \kappa^{\perp} \wedge l^{\perp} - \rho^{\perp} \wedge n^{\perp}    \\
& \dd_X l  =  \alpha_u    e^{-\cF } \theta \wedge n  - \alpha^{\perp}_n\theta \wedge \omega  - \theta  \wedge \kappa^{\perp}   \\
& \dd_X n =  - \alpha_u  e^{-\cF } \theta \wedge l  + \alpha^{\perp}_l\theta \wedge \omega  - \theta  \wedge \rho^{\perp}  
\end{eqnarray*}
	
\noindent
together with the conditions contained in Lemma \ref{lemma:zeroconditions}.
\end{lemma}

\begin{proof}
Follows by expanding and combining equations \eqref{eq:torsion2d} and \eqref{eq:exteriorparallelism2d}.
\end{proof}

\noindent
The way to proceed with the seemingly daunting differential system given in the previous lemma is to first identically solve some of the equations by isolating for $\kappa$ and $\rho$, which are not variables of the system, as \emph{composite} in terms of the underlying adapted isotropic coframe and torsion.

\begin{prop}
Let $(u,v,l,n)$ be an invariant adapted isotropic coframe $(u,v,l,n)$. Then:
\begin{eqnarray*}
& \kappa_u   =  (\dd_X e^{\cK})(l^{\perp}) - \alpha^{\perp}_n e^{\cK} - (\partial_{x_u}\omega)(l^{\perp})\, , \quad \kappa_v = - \alpha^{\perp}_n\\
& \rho_u   = (\dd_X e^{\cK})(n^{\perp})  + \alpha^{\perp}_l e^{\cK} - (\partial_{x_u}\omega)(n^{\perp})\, , \quad \rho_v = \alpha^{\perp}_l \\
& e^{\cF} \kappa^{\perp} = - \alpha^{\perp}_n e^{\cF} \omega + \alpha_u   n^{\perp}  +  ((\dd_X e^{\cK}) (l^{\perp})   -  (\partial_{x_u}\omega)(l^{\perp}))  \theta -  \partial_{x_u} l\\
& e^{\cF} \rho^{\perp} = \alpha^{\perp}_l e^{\cF} \omega  - \alpha_u l^{\perp} +   ((\dd_X e^{\cK}) (n^{\perp})  - (\partial_{x_u}\omega)(n^{\perp}) ) \theta - \partial_{x_u} n
\end{eqnarray*}
	
\noindent
where $\kappa , \rho \in \Omega^1(M)$ are the one-forms relative to which $(u,v,l,n)$ satisfies the differential system \eqref{eq:duvtorsion} and \eqref{eq:dlntorsion}.
\end{prop}

\begin{proof}
We begin with the second equation in the first line of Lemma \ref{lemma:reducedinterequations}, which we immediately solve by isolating for $\kappa_u$ and $\rho_u$ via evaluation on $l^{\perp}$ and $n^{\perp}$. This gives: 
\begin{eqnarray*}
& \kappa_u   = (\dd_X e^{\cK})(l^{\perp})  - \alpha^{\perp}_n e^{\cK} - (\partial_{x_u}\omega)(l^{\perp})\\
& \rho_u   =  (\dd_X e^{\cK})(n^{\perp})  + \alpha^{\perp}_l e^{\cK} - (\partial_{x_u}\omega)(n^{\perp}) 
\end{eqnarray*}
	
\noindent
which correspond to the equations for $\kappa_u$ and $\rho_u$ in the statement of the lemma. Plugging these expressions into the second and third lines of Lemma \ref{lemma:reducedinterequations} and solving for $\kappa^{\perp}$ and $\rho^{\perp}$ we obtain the third and fourth equations in the statement of the lemma. This, together with Lemma \ref{lemma:zeroconditions} gives the desired result.
\end{proof}

\noindent
Once we have identically solved the equations that determine $\kappa$ and $\rho$ for a given adapted isotropic coframe $(u,v,l,n)$, we identify the reduced system in which both $\kappa$ and $\rho$ have been substituted in terms of $(u,v,l,n)$ and the underlying torsion $\alpha$.
\begin{prop}
\label{prop:reduceddifferential}
An invariant adapted isotropic coframe $(u,v,l,n)$ is skew-torsion if and only if the following differential system is satisfied:
\begin{eqnarray*}
&  \partial_{x_u} \theta  = \dd_X e^{ \cF} - e^{ \cF} \ast_{q} \alpha^{\perp}  \, , \qquad \dd_X \theta =   \ast_{q}\alpha^{\perp} \wedge \theta\\
& e^{\cF} \dd_X\omega  =  \alpha_u    \nu_{q} +  (\dd_X e^{\cK}  - \partial_{x_u} \omega) \wedge \theta - l \wedge \partial_{x_u} l  - n \wedge \partial_{x_u} n \\
& e^{\cF }\dd_X l = \theta \wedge \partial_{x_u} l   \, , \qquad e^{\cF }\dd_X n = \theta \wedge \partial_{x_u} n 
\end{eqnarray*}
\end{prop}

\begin{proof}
Follows by plugging Proposition \ref{prop:reduceddifferential} into Lemma \ref{lemma:reducedinterequations}.
\end{proof} 
 
\noindent
We arrive at a differential system, equivalent to the original differential system in \eqref{eq:duvtorsion} and \eqref{eq:dlntorsion} under the given assumptions, that becomes a system of constrained evolution equations in the Cartesian coordinate $x_u$ on a two-dimensional oriented manifold $X$. The first line in the differential system of Proposition \ref{prop:reduceddifferential}, namely:
\begin{eqnarray}
\label{eq:constrainproblem}
\partial_{x_u} \theta  = \dd_X e^{ \cF} - e^{ \cF} \ast_{q} \alpha^{\perp}  \, , \qquad \dd_X \theta =   \ast_{q}\alpha^{\perp} \wedge \theta
\end{eqnarray}

\noindent
can be considered as a \emph{constrained} evolution problem for $\theta$, that determines the latter in terms of \emph{given data} $( \cF , \alpha^{\perp}, q)$. Once $\theta$ has been determined by this constrained evolution problem, it can be plugged into the remaining equations, which now conform a separate independent system. 

\begin{prop}
	\label{prop:reduceddifferential2}
A solution $\theta$ to the constrained evolution problem \eqref{eq:constrainproblem} with initial data $\theta_o\in \Omega^1(X)$ exists if and only if:
\begin{equation}
\label{eq:conditions2d}
\dd_X \theta_o =  \ast_{q_o}\alpha^{\perp}_o \wedge \theta_o\, , \qquad e^{\cF}\dd_X \ast_{q} \alpha^{\perp} = \theta \wedge \partial_{x_u} \ast_{q} \alpha^{\perp} 
\end{equation}
	
\noindent
in which case, $\theta$ is uniquely given by:
\begin{equation}
	\label{eq:integratedtheta}
\theta = \int_0^{x_u} (\dd_X e^{\cF} - e^{\cF} \ast_{q}\alpha^{\perp}) \, \dd x_u + \theta_o
\end{equation}
	
\noindent
where $\alpha_o := \alpha\vert_{x_u = 0} $ and $q_o = q\vert_{x_u = 0} $.
\end{prop}

\begin{proof}
Suppose that $\theta$ is a solution of \eqref{eq:constrainproblem} for given data $(\cF,\alpha^{\perp},q)$. Integrating the first equation in \eqref{eq:constrainproblem} with initial data $\theta_o\in \Omega^1(X)$ we immediately obtain \eqref{eq:integratedtheta}. The first equation in \eqref{eq:conditions2d} is recovered by evaluating the second equation in \eqref{eq:constrainproblem} at $x_u = 0$. On the other hand, taking the derivative of the second equation in \eqref{eq:constrainproblem} with respect to $x_u$ we obtain:
\begin{eqnarray*}
&\partial_{x_u}\dd_X \theta = \dd_X \partial_{x_u}\theta = - \dd_X (e^{\cF}\ast_q \alpha^{\perp}) =  \partial_{x_u}(\ast_{q}\alpha^{\perp}) \wedge \theta +  \ast_{q}\alpha^{\perp} \wedge \partial_{x_u}\theta =   \\
& = \partial_{x_u}(\ast_{q}\alpha^{\perp}) \wedge \theta  + \ast_q \alpha^{\perp}\wedge \dd_X e^{ \cF} 
\end{eqnarray*}

\noindent
Simplifying this expression we obtain the second equation in \eqref{eq:conditions2d}. For the converse, assume that $(\theta,\cF,\alpha^{\perp},q)$ satisfies Equation \eqref{eq:conditions2d} and take $\theta$ to be given as in \eqref{eq:integratedtheta}. Since the first equation in \eqref{eq:constrainproblem} is clearly satisfied, we only need to prove that the second equation in \eqref{eq:constrainproblem} holds. Taking the exterior derivative of $\theta$ we obtain: 
\begin{eqnarray*}
	& \dd_X \theta = -\int_0^{x_u} \dd_X ( e^{\cF} \ast_{q}\alpha^{\perp}) \dd x_u + \dd_X \theta_o= -\int_0^{x_u} \partial_{x_u} \theta \wedge \ast_{q}\alpha^{\perp} \dd x_u -\int_0^{x_u}  e^{\cF} \dd_X\ast_{q}\alpha^{\perp}  \dd x_u  + \dd_X \theta_o	\\
	& = \int_0^{x_u}  (\theta \wedge \partial_{x_u}  \ast_{q}\alpha^{\perp} -  e^{\cF} \dd_X\ast_{q}\alpha^{\perp} ) \dd x_u - \theta \wedge \ast_{q}\alpha^{\perp} +\theta_o \wedge \ast_{q}\alpha^{\perp}_o + \dd_X \theta_o
\end{eqnarray*}

\noindent
and therefore by Equation \eqref{eq:conditions2d} this gives the second equation in \eqref{eq:constrainproblem}.
\end{proof}

\noindent
By the previous proposition, after substituting for $\theta$ the condition satisfied by $\alpha$ becomes an integro-differential equation given by:
\begin{equation*}
	e^{\cF}\dd_X \ast_{q} \alpha^{\perp} = (\int_0^{x_u} (\dd_X e^{\cF} - e^{\cF} \ast_{q }\alpha^{\perp} ) \, \dd x_u + \theta_o)\wedge \partial_{x_u} \ast_{q} \alpha^{\perp} 
\end{equation*}

\noindent
Once we solve \eqref{eq:constrainproblem} for $\theta$ and $\alpha$ we can consider them as given data for the remaining equations, which need to be solved for $\omega$, $l$ and $n$. Note that within this approach, $\cF$ is considered as given data, otherwise we should take into account that in general $\theta$ depends on $\cF$ as stated in formula \eqref{eq:integratedtheta}. Combining propositions \ref{prop:reduceddifferential} and \ref{prop:reduceddifferential2} we obtain the following reduced characterization of invariant adapted isotropic coframes. 

\begin{cor}
\label{cor:2devol}
An invariant adapted isotropic coframe $(u,v,l,n)$ is skew-torsion if and only if $(\omega,l,n)$ satisfies the following differential system:
\begin{eqnarray*} 
& e^{\cF}\dd_X l = \theta \wedge \partial_{x_u} l\, , \qquad  e^{\cF}\dd_X n = \theta \wedge \partial_{x_u} n \, , \qquad e^{\cF}\dd_X \ast_{q} \alpha^{\perp}  = \theta \wedge \partial_{x_u} \ast_{q } \alpha^{\perp} \\
& e^{\cF} \dd_X\omega  =  \alpha_u    \nu_{q} +  (\dd_X e^{\cK}  - \partial_{x_u} \omega) \wedge \theta - l \wedge \partial_{x_u} l  - n \wedge \partial_{x_u} n 
\end{eqnarray*}
	
\noindent
where:
\begin{equation*}
\theta  = \int_0^{x_u} (\dd_X e^{\cF } - e^{\cF } \ast_{q }\alpha^{\perp} ) \, \dd x_u + \theta_o
\end{equation*}
	
\noindent
with $\dd_X \theta_o = \ast_{q_o}\alpha^{\perp}_o \wedge \theta_o $.
\end{cor}

\noindent
Note that by a standard application of the existence and uniqueness of solutions for first-order symmetric hyperbolic systems it follows that all equations in the previous corollary except for the equation for $\alpha^{\perp}$ can be solved \emph{locally} in $x_u$. Hence, the existence of adapted isotropic coframes reduces to proving existence of a solution $\alpha^{\perp}$ to the following integro-differential equation:
\begin{eqnarray*}
e^{\cF}\dd_X \ast_{q} \alpha^{\perp}  = ( \int_0^{x_u} (\dd_X e^{\cF } - e^{\cF } \ast_{q }\alpha^{\perp} ) \, \dd x_u + \theta_o) \wedge \partial_{x_u} \ast_{q } \alpha^{\perp} 
\end{eqnarray*}

\noindent
for given $(l,n,\cF)$ satisfying $\dd_X \theta_o = \ast_{q_o}\alpha^{\perp}_o \wedge \theta_o $ as initial data condition.

 
\section{Supersymmetric configurations}
\label{sec:susyconf}

 
In this section we study the \emph{supersymmetric configurations} or \emph{BPS states} of four-dimensional NS-NS supergravity. This is precisely the supergravity theory whose supersymmetry conditions involve skew-torsion parallel spinors. 
 
\begin{lemma}
\label{lemma:dilatinoeq}
A NS-NS configuration $(g,b,\phi)\in \Conf(\cC,\cX)$ satisfies the dilatino equation in \eqref{eq:KSE4d} with respect to $\varepsilon\in\Gamma(S)$ if and only if:
\begin{eqnarray*}
\varphi_{\phi}\wedge u = \ast_g(\alpha_b\wedge u)\, , \quad \varphi_{\phi} \wedge l\wedge u = \alpha_b (l) \ast_g u\, , \quad  \varphi_{\phi}(l)\, u =  \ast_g (u\wedge l\wedge \alpha_b)\, , \quad  \varphi_{\phi}(u) = \alpha_b(u) = 0 
\end{eqnarray*}

\noindent
where $(u,[l]_u)$ is the parabolic pair associated to $\varepsilon$.  
\end{lemma}
 
\begin{proof}
By Proposition \ref{prop:constraintendopoly}, the dilatino equation holds if and only if:
\begin{equation}
\label{eq:dilatino}
\varphi_{\phi}\diamond_g \alpha_{\varepsilon} = H_b\diamond_g \alpha_{\varepsilon}
\end{equation}

\noindent
where $\alpha_{\varepsilon} = u + u\wedge l$ is the square of $\epsilon$, see Definition \ref{def:squarespinor}, that determines its associated parabolic pair $(u,[l]_u)$. We compute:
\begin{eqnarray*}
& \varphi_{\phi}\diamond_g \alpha_{\varepsilon} = \varphi_{\phi}(u) + \varphi_{\phi}\wedge u + \varphi_{\phi}\wedge u \wedge l +  \varphi_{\phi}(u)\, l -  \varphi_{\phi}(l)\, u \\
& H_b \diamond_g \alpha_{\varepsilon} = \alpha_b\diamond_g\nu_g \diamond_g \alpha_{\varepsilon} = \ast_g (\alpha_b\wedge u) - \alpha(u) \nu_h  - \ast_g(u \wedge l \wedge \alpha_b) + \alpha_b(u)\, \ast_g l - \alpha_b(l)\, \ast_g u
\end{eqnarray*}

\noindent
where in the second equation we have used Equation \eqref{eq:nuaction} to rewrite the geometric product in terms of $\alpha_b = \ast_g H_b$. Plugging these expressions in \eqref{eq:dilatino} and separating by degree in \eqref{eq:dilatino} gives the conclusion.
\end{proof}

\noindent
Recall that given a parabolic pair $(u,[l]_u)$, we denote by $u^{\perp_g}\subset T^{\ast}M$ the distribution spanned by the orthogonal of $u\in \Omega^1(M)$ in $T^{\ast}M$, and we denote by:
\begin{equation*}
\mathfrak{G}_u = \frac{u^{\perp_g}}{u}
\end{equation*}

\noindent
the corresponding dual screen bundle. In particular, any section $w \in \Gamma(u^{\perp_g})$ defines a section of $\mathfrak{G}_u$ which we denote by $[w] \in \Gamma(\mathfrak{G}_u)$

\begin{prop}
\label{prop:susydilatino}
A NS-NS configuration $(g,b,\phi)\in \Conf(\cC,\cX)$ satisfies the dilatino equation in \eqref{eq:KSE4d} with respect to $\varepsilon\in\Gamma(S)$ with associated parabolic pair $(u,[l]_u)$ if and only if $\varphi_{\phi} , \alpha_b \in \Gamma(u^{\perp_g})$ and there exists a section $[\mathfrak{m}] \in \Gamma(\mathfrak{G}_u)$ such that:
\begin{equation*}
[\varphi_{\phi}] = \ast_q [\mathfrak{m}]  \, , \qquad [\alpha_b] =  [\mathfrak{m}] 
\end{equation*} 

\noindent
where $q$ is the Riemannian metric induced by $g$ on $\mathfrak{G}_u$ and $\ast_q\colon \wedge \frG_u \to \wedge \frG_u$ is the associated Hodge dual. 
\end{prop}

\begin{proof}
Conditions $\varphi_{\phi}(u) = \alpha_b(u) = 0$ in Lemma \ref{lemma:dilatinoeq} are equivalent to $\phi_{\phi} , \alpha_b \in \Gamma(u^{\perp_g})$. Furthermore:
\begin{equation*}
u\wedge \varphi_{\phi} = u\wedge [\varphi_{\phi}]\, , \qquad u\wedge \alpha_{b} = u\wedge [\alpha_{b}]
\end{equation*}

\noindent
which, plugged into the first equation in Lemma \ref{lemma:dilatinoeq}, gives:
\begin{equation*}
u\wedge [\varphi_{\phi}] = \ast_g(u\wedge \alpha_b ) = \ast_g(u\wedge  [ \alpha_b] ) = u\wedge \ast_q [ \alpha_b] 
\end{equation*}

\noindent
and therefore we obtain $[\varphi_{\phi}] = \ast_q [\alpha_b]$. With this condition the remaining equations in Lemma \ref{lemma:dilatinoeq} are automatically satisfied and thus we conclude.
\end{proof}

\noindent
By the previous proposition, if  $(g,b,\phi)\in \Conf(\cC,\cX)$ satisfies the dilatino equation in \eqref{eq:KSE4d} with respect to an isotropic parallelism $[u,v,l,n]$, then for every representative $(u,v,l,n)$ there exists a unique section $\frm \in \Gamma(\langle\mathbb{R}u\rangle \oplus \langle\mathbb{R}v\rangle)^{\perp_g}$ such that:
\begin{equation*}
\varphi_{\phi} = \frc_{\phi} u + \ast_{q} \frm\, , \qquad \alpha_b = \frf_b u +   \frm
\end{equation*}

\noindent
for functions $\frc_{\phi} , \frf_b \in C^{\infty}(M)$ depending not only on $\phi$ and $b$ respectively, but also on the choice of representative $(u,v,l,n)\in [u,v,l,n]$. If $(u,v^{\prime},l^{\prime},n^{\prime})\in [u,v,l,n]$ is any other representative, then with respect to this choice of representative we have:
\begin{equation*}
\varphi_{\phi} = \frc^{\prime}_{\phi} u + \ast_{q^{\prime}} \frm^{\prime}_b\, , \qquad \alpha^{\prime}_b = \frc^{\prime}_b u +   \frm^{\prime}_b
\end{equation*}

\noindent
where:
\begin{equation*}
	\frc_{\phi}^{\prime} = \frc_{\phi} + \varphi_{\phi}(\frw) \, , \qquad 	\frf_b^{\prime} = \frf_b + \alpha_b(\frw) 
\end{equation*}

\noindent
in terms of the unique $\frw \in \Gamma(\langle\mathbb{R}u\rangle \oplus \langle\mathbb{R}v\rangle)^{\perp_g}$ for which:
\begin{equation*}
(u,v^{\prime},l^{\prime},n^{\prime}) = \frw\cdot (u,v,l,n) = (u , v - \frac{1}{2} \vert\frw\vert_g^2 u + \frw , l - \frw(l) u , n - \frw(n) u)
\end{equation*}

\noindent
Based on the previous formulae, it may seem possible to choose the representative $(u,v,l,n)\in [u,v,l,n]$ wisely so as to have $\frc_\phi = \frf_b = 0$. However, this may not be possible in general since we cannot guarantee that the supports of $\frf_b , \frc_{\phi} \in C^{\infty}(M)$ are contained in the supports of $\alpha(\frw) , \varphi_{\phi}(\frw) \in C^{\infty}(M)$, respectively. Proposition \ref{prop:susydilatino} together with the characterization of skew-torsion parallel spinors in terms of isotropic parallelisms given in Proposition \ref{prop:existencenullcoframeII} gives the following characterization of supersymmetric configurations.  

\begin{prop}
\label{prop:susyconfiguration}
A configuration $(g,b,\phi)\in \Conf(\cC,\cX)$ is supersymmetric if and only if there exists a skew-torsion isotropic parallelism $[u,v,l,n]$ such that:
\begin{equation*}
\varphi_{\phi}(u) = 0 \, , \qquad \alpha_b(u) = 0\, , \qquad [\varphi_{\phi}]  = \ast_q [\alpha_b] 
\end{equation*}

\noindent
where $[\varphi_{\phi}] , [\alpha_b]\in \Gamma(\mathfrak{G}^{\ast}_u )$ are sections of the screen bundle associated to $u\in \Omega^1(M)$.
 \end{prop}

\noindent
As an immediate consequence of the previous proposition we obtain the following result, which justifies the title of this chapter.
\begin{cor}
	\label{cor:susykundt}
Let $(g,b,\phi,\varepsilon)\in \Conf(\cC,\cX)$ be a supersymmetric configuration. Then, $(M,g)$ is a Kundt Lorentzian four-manifold.
\end{cor}

\noindent
Furthermore, again as a consequence of Proposition \ref{prop:susyconfiguration}, we obtain the following causal character for the dilaton and torsion of a supersymmetric configuration.
\begin{cor}
Let $(g,b,\phi,\varepsilon)\in \Conf_s(\cC,\cX)$ be a supersymmetric configuration. Then, $\varphi_{\phi}$, respectively $\alpha_b$, is nowhere time-like and is isotropic at a point $m\in M$ if and only if $[\varphi_{\phi}]\vert_m = 0$, respectively $[\alpha_b]\vert_m = 0$.
\end{cor}

\noindent
Given a supersymmetric configuration $(g,b,\phi,\varepsilon)\in \Conf_s(\cC , \cX)$, we canonically obtain an skew-torsion isotropic parallelism $[u,l,v,n]$ on $M$ given by the unique isotropic parallelism associated to $\varepsilon$, which by Corollary \ref{cor:equaltorsion} determines the torsion uniquely. This defines a natural functor:
\begin{equation*}
\Conf_s(\cC , \cX) \to  \mathfrak{F}_{sk}(M)
\end{equation*}   
 
\noindent
This functor is not essentially surjective in general, since by Lemma \ref{lemma:dilatinoeq} the isotropic parallelism associated to a supersymmetric configuration necessarily satisfies $\alpha(u) = 0$, which may not be the case for general skew-torsion isotropic parallelisms. Instead, we have the following result.
 
\begin{prop}
\label{prop:iffsusyconf}
A skew-torsion isotropic parallelism $[u,v,l,n] \in \mathfrak{F}_{sk}(M)$ defines a supersymmetric configuration if and only if:
\begin{equation*}
\alpha(u) = 0\, , \qquad \frac{1}{2\pi} [\ast_g \alpha] \in H^3(M,\mathbb{Z})
\end{equation*}

\noindent
and in addition there exists a function $\mathfrak{c}\in C^{\infty}(M)$ such that:
\begin{equation*}
\frac{1}{2\pi}[\mathfrak{c}\, u + \ast_q(\alpha - \alpha(v))]\in H^1(M,\mathbb{Z})
\end{equation*}
 	
\noindent
for any, and hence all, representatives $(u,v,l,n) \in [u,v,l,n]$.
\end{prop}

\begin{remark}
In the previous proposition we are implicitly using that by Corollary \ref{cor:equaltorsion} a skew-torsion isotropic parallelism determines uniquely its torsion.
\end{remark}

\begin{proof}
The \emph{only if} direction follows directly from Proposition \ref{prop:susyconfiguration}. To prove the converse, we observe that if:
\begin{equation*}
\frac{1}{2\pi} [\ast_g \alpha] \in H^3(M,\mathbb{Z})
\end{equation*}

\noindent
then there exists a gerbe with connective structure and curving $b\in \Omega^2(Y)$ such that $\alpha = \alpha_b$. Furthermore
if $[\mathfrak{c}  u + \ast_q(\alpha - \alpha(v))]\in 2\pi\, H^1(M,\mathbb{Z})$, then there exists a $\mathbb{Z}$-covering and a function $\phi$ on the total space of this covering such that:
\begin{equation*}
\varphi_{\phi} := \mathfrak{c} u + \ast_q(\alpha - \alpha(v)) \in \Omega^1(M)
\end{equation*}

\noindent
In particular, $\varphi_{\phi}$ is closed and satisfies $\varphi_{\phi}(u) = 0$ and thus by Proposition \ref{prop:susyconfiguration} we conclude.
\end{proof}
 
\noindent
If, in the situation of the previous proposition, we write:
\begin{equation*}
\alpha = \frf\, u +  \frm
\end{equation*}

\noindent
in terms of a function $\frf\in C^{\infty}(M)$, then $ \mathfrak{c} u + \ast_q(\alpha - \alpha(v)) = \frc u +  \ast_q \frm$ as required by the characterization of supersymmetric configurations given in Proposition \ref{prop:susyconfiguration}. 

\begin{remark}
The function $\frc \in C^{\infty}(M)$ occurring in the expression $[\mathfrak{c} u + \ast_q(\alpha - \alpha(v))]\in H^1(M,\mathbb{Z})$ definitely depends on the choice of representative $(u,v,l,n) \in [u,v,l,n]$. If $(u,v^{\prime},l^{\prime},n^{\prime}) \in [u,v,l,n]$ is another representative, then a calculation gives:
\begin{equation*}
c^{\prime} = c + \ast_q (\frw\wedge (\alpha-\alpha(v)))
\end{equation*}

\noindent
in terms of the unique $\frw \in \Gamma(\langle\mathbb{R}u\rangle \oplus \langle\mathbb{R}v\rangle)^{\perp_g}$ for which $(u,v^{\prime},l^{\prime},n^{\prime}) = \frw\cdot (u,v,l,n)$.
\end{remark} 

\noindent
\noindent
By Proposition  \ref{prop:susyconfiguration}, the differential system satisfied by the skew-torsion isotropic parallelism $[u,v,l,n]$ of a supersymmetric configuration reduces to:
\begin{eqnarray}
& \dd u = u\wedge (\alpha_n l - \alpha_l n) \, , \qquad \dd v = -\kappa\wedge l-\rho\wedge n + v\wedge (\alpha_l n - \alpha_n l) - \frf_b l\wedge n \label{eq:reducedintegrability1}\\
& \dd l = u \wedge (\frf_b n - \alpha_n v - \kappa)  \, , \qquad \dd n = u \wedge (\alpha_l v - \frf_b l - \rho)  \label{eq:reducedintegrability2}  
\end{eqnarray}

\noindent
Consequently, the \emph{integrability conditions} satisfied by the skew-torsion isotropic parallelism $[u,v,l,n]$ of a supersymmetric configuration also simplify. 

\begin{lemma}
\label{lemma:skewtorsionintegrabilityI}
Let $(g,b,\phi,\varepsilon)\in \Conf_s(\cC,\cX)$ be a supersymmetric configuration with associated isotropic parallelism $[u,v,l,n]\in \mathfrak{F}_{sk}(M)$. Then, the following formulas hold:
\begin{eqnarray*}
& u^{\sharp_g}(\alpha_l) = u^{\sharp_g}(\alpha_n) = 0\, , \quad l^{\sharp_g}(\alpha_l) + n^{\sharp_g}(\alpha_n) = 0\\
& u\wedge (\dd\kappa + \dd\alpha_n \wedge v - \dd\frf_b \wedge n + n \wedge (\alpha_l \kappa + \alpha_n \rho)) = 0\\
& u\wedge (\dd\rho + \dd\frf_b \wedge l - \dd\alpha_l \wedge v - l \wedge (\alpha_l \kappa + \alpha_n \rho)) = 0\\
& \dd\kappa\wedge l + \dd\rho \wedge n +  (\alpha_l \kappa + \alpha_n \rho) \wedge l\wedge n + v \wedge (\dd\alpha_l \wedge n - \dd \alpha_n\wedge l) + \dd \frf_b \wedge l \wedge n = 0
\end{eqnarray*}

\noindent
In particular, $\dd\alpha (u,l) = 0$ and $\dd \alpha (u,n) = 0$
\end{lemma}

\begin{proof}
The first line in the statement of lemma is equivalent to the exterior derivative of the first equation in \eqref{eq:reducedintegrability1}. The second line is equivalent to the exterior derivative of the first equation in \eqref{eq:reducedintegrability2}, whereas the third line is equivalent to the exterior derivative of the second equation in \eqref{eq:reducedintegrability2}. Finally, the fourth line in the statement of the lemma is equivalent to the exterior derivative of the second equation in \eqref{eq:reducedintegrability1}.
\end{proof}

\begin{prop}
\label{prop:skewtorsionintegrabilityI}
Let $(g,b,\phi,\varepsilon)\in \Conf_s(\cC,\cX)$ be a supersymmetric configuration with associated isotropic parallelism $[u,v,l,n]\in \mathfrak{F}_{sk}(M)$. Then, the following formulas hold:
\begin{eqnarray*}
& \dd\kappa(l,n) = l(\frf_b)  + \alpha_l \kappa(l) + \alpha_n \rho(l)\, , \quad \dd\kappa(u,n) = n(\alpha_n) \, , \quad \dd\kappa(u,l) = l(\alpha_n)  \\
& \dd\rho(l,n) = n(\frf_b)  + \alpha_l \kappa(n) + \alpha_n \rho(n)\, , \quad \dd\rho(u,n) = - n(\alpha_l) \, , \quad \dd\rho(u,l) = - l(\alpha_l)  \\
& \alpha_l \kappa(u) + \alpha_n \rho(u) + u(\frf_b) = 0\, , \quad \dd\kappa(u,v) = v(\alpha_n)\, , \quad \dd\rho(u,v) = - v(\alpha_l)\\
& \dd \rho (v,l) - \dd \kappa(v,n)  + v(\frf_b) + \alpha_l \kappa(v) + \alpha_n \rho(v) = 0
\end{eqnarray*}
	
\noindent
Together with conditions $u^{\sharp_g}(\alpha_l) = u^{\sharp_g}(\alpha_n) = 0$ and $l^{\sharp_g}(\alpha_l) + n^{\sharp_g}(\alpha_n) = 0$.
\end{prop}

\begin{proof}
The result follows from Lemma \ref{lemma:skewtorsionintegrabilityI} after expanding all equations in the given isotropic coframe $(u,v,l,n)\in [u,v,l,n]$ and combining them appropriately.
\end{proof}

\noindent
By the previous proposition, if $(u,v,l,n) \in [u,v,l,n]$ is a skew-torsion isotropic coframe associated to a supersymmetric configuration, then have:
\begin{eqnarray*}
& \dd\kappa = \dd\kappa (v,l) u\wedge l + \dd\kappa(v,n) u\wedge n + v(\alpha_n) v\wedge u \\
& + l(\alpha_n) v\wedge l + n(\alpha_n) v\wedge n + (l(\frf_b)  + \alpha_l \kappa(l) + \alpha_n \rho(l)) l\wedge n \\
& \dd \rho = \dd\rho (v,l) u\wedge l + \dd\rho(v,n) u\wedge n - v(\alpha_l) v\wedge u \\
& - l(\alpha_l) v\wedge l  - n(\alpha_l) v\wedge n + ( n(\frf_b)  + \alpha_l \kappa(n) + \alpha_n \rho(n)) l\wedge n 
\end{eqnarray*}

\noindent
for the characteristic one-forms $\kappa$ and $\rho$ of $(u,v,l,n)$.
 
\begin{lemma}
\label{lemma:dvarphi}
Let $(g,b,\phi,\varepsilon)\in \Conf_s(\cC,\cX)$ be a supersymmetric configuration with isotropic parallelism $[u,v,l,n]\in \mathfrak{F}_{sk}(M)$. Then:
\begin{eqnarray}
&  u(\frc_{\phi}) = \alpha_n \kappa(u) - \alpha_l \rho(u)  + \alpha_l^2 + \alpha_n^2 \nonumber \\
&  l(\frc_{\phi}) =\frc_{\phi}\alpha_n - v(\alpha_n) - \frf_{b} \alpha_l - \alpha_l \rho(l) + \alpha_n \kappa(l) \label{eq:dvarphi}  \\
&   n(\frc_{\phi}) = v(\alpha_l) - \frf_{b} \alpha_n - \frc_{\phi}\alpha_l  + \alpha_n \kappa(n) - \alpha_l \rho(n) \nonumber
\end{eqnarray}

\noindent
where $\varphi_{\phi} = \mathfrak{c}_{\phi} u + \ast_q \frm = \mathfrak{c}_{\phi} u + \alpha_l n - \alpha_n l$.
\end{lemma}

\begin{remark}
The previous system of equations can be understood as a differential system for $\frc_{\phi}\in C^{\infty}(M)$, which is a function not occurring as a variable in the differential system satisfied by the isotropic parallelism associated to a supersymmetric configuration. 
\end{remark}

\begin{proof}
By Proposition \ref{prop:susyconfiguration} we have $\varphi_{\phi} = \mathfrak{c}_{\phi} u + \ast_q \frw$ for a unique $\frw \in \Gamma(\langle\mathbb{R}u\rangle \oplus \langle\mathbb{R}v\rangle)^{\perp_g}$ in any given conjugate parallelism $(u,v,l,n)\in [u,v,l,n]$. Since $\varphi_{\phi}$ is necessarily closed, we must have:
\begin{eqnarray*}
& 0 = \dd\varphi_{\phi} = \dd\mathfrak{c}_{\phi}\wedge u + \mathfrak{c}_{\phi} \dd u +  \dd\ast_q\frw = \dd\mathfrak{c}_{\phi}\wedge u + \mathfrak{c}_{\phi} u\wedge (\alpha_n l - \alpha_l n)+  \dd (\alpha_l n - \alpha_n l) = \dd\mathfrak{c}_{\phi}\wedge u \\
&  + \mathfrak{c}_{\phi} u\wedge (\alpha_n l - \alpha_l n) -  \dd\alpha_n \wedge l + \dd\alpha_l \wedge n - \alpha_n u \wedge (\frf\, n - \alpha_n v - \kappa) + \alpha_l u \wedge (\alpha_l v - \frf\, l - \rho) 
\end{eqnarray*}

\noindent
Further expanding the previous equation and using the integrability conditions of Lemma \ref{lemma:skewtorsionintegrabilityI}, we obtain the relations given in \eqref{eq:dvarphi}.
\end{proof}

\noindent
We compute below the covariant derivative of both the curvature $\varphi_{\phi}$ of the dilaton $\phi$ and (the dual of) the curvature $\alpha_b$ of the $b$-field of a supersymmetric configuration. This formulae will be used extensively in the following.

\begin{lemma}
\label{lemma:nablavarphi}
Let $(g,b,\phi,\varepsilon)\in \Conf_s(\cC,\cX)$ be a supersymmetric configuration with isotropic parallelism $[u,v,l,n]\in \mathfrak{F}_{sk}(M)$. Then:
\begin{eqnarray*}
& \nabla^g\alpha_{b} =  (v(\frf_b) + \alpha_l \kappa(v) + \alpha_n \rho(v)) u\otimes u    + l(\alpha_l) l\otimes l + n(\alpha_l) n\otimes n + v(\alpha_l) u\otimes l   \\
& + \dd\kappa (l,n) l\otimes u + v(\alpha_n)  u\otimes n + \dd\rho(l,n) n\otimes u + l(\alpha_n) l\otimes n + n(\alpha_l) n\otimes l\\
& \nabla^g\varphi_{\phi} = (v(\frc_{\phi}) - \alpha_n \kappa(v) + \alpha_l \rho(v)) u\otimes u - l(\alpha_n) l\otimes l + n(\alpha_l) n\otimes n + \frac{1}{2} (\alpha_l^2 + \alpha_n^2) u\odot v \\
& + ( \frac{1}{2} \alpha_n \frc_{\phi}  - v(\alpha_n)  - \frac{1}{2} \alpha_l \frf_b) u\odot l + (v(\alpha_l) - \frac{1}{2}\alpha_n \frf_{b} - \frac{1}{2} \alpha_l \frc_{\phi} ) u\odot n + l(\alpha_l) l\otimes n - n(\alpha_n) n\otimes l  \label{eq:nablavarphi}
\end{eqnarray*}
 	
\noindent
where $\alpha_{b} = \mathfrak{f}_b u + \frm_b$ and $\varphi_{\phi} = \mathfrak{c}_b u + \ast_q\frm_b$. 
\end{lemma}
 
\begin{proof}
By Proposition \ref{prop:existencenullcoframe} adapted to the case of the skew-torsion isotropic parallelism of a supersymmetric configurations, we have: 
\begin{eqnarray*}
& \nabla^g u =  \frac{1}{2} u\wedge (\alpha_n l - \alpha_l n)  \, , \qquad \nabla^g v = - \kappa\otimes l - \rho\otimes n + \frac{1}{2}v\wedge (\alpha_l n - \alpha_n l) - \frac{1}{2} \frf\, l\wedge n  \\
& \nabla^g l = \kappa\otimes u + \frac{1}{2}  u \wedge (\frf\, n - \alpha_n v)   \, , \qquad \nabla^g n = \rho\otimes u + \frac{1}{2} u \wedge (\alpha_l v - \frf\, l)  
\end{eqnarray*}

\noindent
Using this formulae together with Proposition \ref{prop:skewtorsionintegrabilityI}, we compute:
\begin{eqnarray*}
& \nabla^g_u \alpha_b = 0\, , \quad \nabla^g_v \alpha_b =(v(\frf_b) + \alpha_l \kappa(v) + \alpha_n \rho(v)) u + v(\alpha_l) l + v(\alpha_n) n \\
& \nabla^g_l \alpha_b =  l(\alpha_l) l + l(\alpha_n) n + \dd\kappa(l,n) u \, , \quad \nabla^g_n \alpha_b = n(\alpha_l) l + n(\alpha_n) n + \dd\rho(l,n) u
\end{eqnarray*}

\noindent
Plugging these relations into the following expansion for $\nabla^g\alpha_b$:
\begin{eqnarray*}
\nabla^g \alpha_b =  u\otimes\nabla^g_v \alpha_b + l\otimes\nabla^g_l \alpha_b + n\otimes\nabla^g_n \alpha_b  
\end{eqnarray*}

\noindent
we obtain the desired result. To obtain the expression for $\nabla^g\varphi_{\phi}$ we proceed analogously by using that $\varphi_{\phi} = \frc_{\phi} u + \ast_q \frm = \frc_{\phi} u + \alpha_l n - \alpha_n l$:
\begin{eqnarray*}
& \nabla^g_u \varphi_{\phi} = (u(\frc_{\phi}) - \frac{1}{2} \alpha_l^2 - \frac{1}{2} \alpha_n^2 - \alpha_n \kappa(u) + \alpha_l \rho(u)) u = \frac{1}{2} (\alpha_l^2 + \alpha_n^2) u\\
&\nabla^g_v \varphi_{\phi} = (v(\frc_{\phi}) - \alpha_n \kappa(v) + \alpha_l \rho(v))u + \frac{1}{2} (\alpha^2_l + \alpha^2_n) v + ( \frac{1}{2} \alpha_n \frc_{\phi}  - v(\alpha_n)  - \frac{1}{2} \alpha_l \frf_b)l \\
&   + (v(\alpha_l) - \frac{1}{2}\alpha_n \frf_{b} - \frac{1}{2} \alpha_l \frc_{\phi} ) n \\
& \nabla^g_l \varphi_{\phi} = (l(\frc_{\phi}) + \alpha_l \rho(l) - \alpha_n \kappa(l) - \frac{1}{2}\frc_{\phi}\alpha_n + \frac{1}{2} \frf_b \alpha_l) u + l(\alpha_l) n - l(\alpha_n) l\\
& = (\frac{1}{2}\frc_{\phi}\alpha_n - \frac{1}{2}\frf_b \alpha_l  - v(\alpha_n)) u + l(\alpha_l) n - l(\alpha_n) l\\
& \nabla^g_n \varphi_{\phi} = (n(\frc_{\phi}) + \alpha_l \rho(n) - \alpha_n \kappa(n) + \frac{1}{2} \frc_{\phi}\alpha_l + \frac{1}{2} \frf_b \alpha_n) u + n(\alpha_l) n - n(\alpha_n) l \\
& = (v(\alpha_l) - \frac{1}{2}\frf_b \alpha_n - \frac{1}{2}\frc_{\phi}\alpha_l) u + n(\alpha_l) n - n(\alpha_n) l 
\end{eqnarray*}

\noindent
where we have used the differential system \eqref{eq:dvarphi}. Plugging these relations into the following expansion for $\nabla^g\varphi_{\phi}$:
\begin{eqnarray*}
\nabla^g \varphi_{\phi} = v\otimes\nabla^g_u \varphi_{\phi} +  u\otimes\nabla^g_v \varphi_{\phi}  + l\otimes\nabla^g_l \varphi_{\phi}  + n\otimes\nabla^g_n \varphi_{\phi}  
\end{eqnarray*}

\noindent
we obtain the desired result and hence we conclude.
\end{proof}

\noindent
Note that the expression given for $\nabla^g\varphi_{\phi}$ in the previous lemma is manifestly symmetric after using the identity $l(\alpha_l) + n(\alpha_n) = 0$ obtained in Lemma \ref{lemma:skewtorsionintegrabilityI}.
\begin{prop}
\label{prop:furtheridentities}
Let $(g,b,\phi,\varepsilon)\in \Conf_s(\cC,\cX)$ be a supersymmetric configuration with associated isotropic parallelism $[u,v,l,n]\in \mathfrak{F}_{sk}(M)$. Then:
\begin{equation}
\label{eq:furtheridentities}
\nabla^{g\ast}\varphi_{\phi} = l(\alpha_n) - n(\alpha_l)  - \alpha_l^2 - \alpha_n^2   \, , \quad \nabla^{g\ast}\alpha_b = 0 \, , \quad  \cL_u \varphi_{\phi} = 0\, , \quad  \cL_u \alpha_{b} = 0
\end{equation}
	
\noindent
where $\cL_u$ denotes the Lie derivative with respect to $u^{\sharp_g} \in \mathfrak{X}(M)$.
\end{prop}

\begin{proof}
Recall that, by definition:
\begin{eqnarray*}
\nabla^{g\ast}\varphi_{\phi} = - \mathrm{Tr}_g(\nabla^{g}\varphi_{\phi})\, , \qquad \nabla^{g\ast}\alpha_{b} = - \mathrm{Tr}_g(\nabla^{g}\alpha_{b})
\end{eqnarray*}

\noindent
and thus computing the trace of the expressions given for $\nabla^g\varphi_{\phi}$ and $\nabla^g\alpha_{b}$ in Lemma \ref{lemma:nablavarphi} we obtain the first two equations in \eqref{eq:furtheridentities}. Equation $\cL_u \varphi_{\phi} = 0$ follows directly from $\varphi_{\phi}(u) = 0$ by Proposition \ref{prop:susyconfiguration} together with the fact that $\varphi_{\phi}$ is a closed one-form. To compute the Lie derivative of $\alpha_b$ we obtain first its exterior derivative by skew-symmetrization of the $\nabla^g\alpha_b$ in Lemma \ref{lemma:nablavarphi}. We obtain: 
\begin{equation*}
\dd\alpha = (l(\alpha_n) - n(\alpha_l)) l\wedge n + (v(\alpha_l) - \dd\kappa(l,n)) u\wedge l + (v(\alpha_n)-\dd\rho(l,n)) u\wedge n   
\end{equation*}

\noindent
and thus $u\lrcorner_g\dd\alpha_b = 0$. We conclude since by Proposition \ref{prop:susyconfiguration} we have $\alpha(u) = 0$. 
\end{proof}

\noindent 
Note that by the previous proposition we also have $\cL_u \ast_g \alpha_b = \cL_u H_b = 0$, since by Proposition \ref{prop:existencenullcoframe} the vector field $u^{\sharp_g} \in \mathfrak{X}(M)$ is Killing. Hence, we have:
\begin{equation*}
\cL_u g = 0 \, , \qquad  \cL_u \varphi_{\phi} = 0\, , \qquad  \cL_u \alpha_{b} = 0
\end{equation*}

\noindent
for every supersymmetric configuration $(g,b,\phi,\varepsilon) \in \Conf_s(\cC,\cX)$. However, we cannot yet say that the supersymmetric configuration is \emph{itself} invariant under the action of $\mathbb{R}$ induced by $u^{\sharp_g}\in\mathfrak{X}(M)$, since this would require \emph{lifting} this action to an action on curvings on $\cC$ and connections on $\cX$. We will deal with this type of questions in the future when we study the \emph{moduli space} of supersymmetric configurations and solutions in NS-NS supergravity.

In Section \ref{sec:generalcurvature} we computed the curvature tensors of the Lorentzian metric associated to a skew-torsion isotropic parallelism, which are all highly constrained. In the case of the skew-torsion parallelism of a supersymmetric configuration these curvature tensors are further constrained. We compute them in the following to end this section. First, as a direct application of Corollary \ref{cor:curvaturetorsion}, for every $w, w_1, w_2 \in \mathfrak{X}(M)$ we obtain:
\begin{eqnarray*}
& \cR^{g,b}_{w_1 w_2}   =   - \dd \kappa(w_1,w_2) \,u\wedge l - \dd\rho(w_1,w_2)\, u\wedge n \\
&\Ric^{g,b}(w)   =   (\dd\kappa(w,l) + \dd\rho(w,n)) u - \dd\kappa(w,u) l - \dd\rho(w,u) n \\
& s^{g,b} = 2( \dd\kappa(u,l) + \dd\rho(u,n))  
\end{eqnarray*}

\noindent
for the curvature tensors of the connection with skew-torsion $H_b = \ast_g\alpha_b$. Expanding as prescribed in Remark \ref{remark:curvatureskewtorsion} and using Lemma \eqref{lemma:nablavarphi}, we obtain the following explicit expressions for the curvature of the Lorentzian metric of a supersymmetric configuration. Alternatively, they can be directly obtained by computing using the differential system given in \eqref{eq:reducedintegrability1} and \eqref{eq:reducedintegrability2}.
\begin{prop}
\label{prop:curvaturessusy}
Let $(g,b,\phi,\varepsilon)\in \Conf_s(\cC,\cX)$ be a supersymmetric configuration with associated isotropic parallelism $[u,v,l,n]\in \mathfrak{F}_{sk}(M)$. The Ricci and scalar curvatures of $g$ are given by:
\begin{eqnarray*}
& 2\Ric^g(u) = (l(\alpha_n) - n(\alpha_l) - \alpha^2_l - \alpha_n^2) u\\
& 2\Ric^g(v) = (l(\alpha_n) - n(\alpha_l) - \alpha_l^2 - \alpha_n^2) v + (\dd\rho(l,n)  + \alpha_l \frf_b + v(\alpha_n) ) l\\
& + (\alpha_n \frf_b - v(\alpha_l) - \dd\kappa(l,n)) n + (2 \dd\kappa(v,l) + 2 \dd\rho(v,n) + \frf_b^2) u\\
& \Ric^g(l) = \frac{1}{2} (v(\alpha_n) + \dd\rho(l,n) + \alpha_l \frf_b)u + (l(\alpha_n) - \frac{1}{2}\alpha_n^2) l + (\frac{1}{2} \alpha_l\alpha_n - l(\alpha_l)) n \\
& \Ric^g(n) =  \frac{1}{2} (\alpha_n \frf_b - v(\alpha_l) - \dd\kappa(l,n) ) u + (n(\alpha_n) + \frac{1}{2} \alpha_l \alpha_n) l - (n(\alpha_l) + \frac{1}{2} \alpha_l^2) n 
\end{eqnarray*}
	
\noindent
where $(u,v,l,n)\in [u,l,v,n]$ is any representative.
\end{prop}
 
\noindent
We will apply this proposition extensively in the following section. We remind the reader that for simplicity in the exposition we are using the notation $\alpha = \frf_b u + \frm_b$ with $\frm_b = \alpha_l l + \alpha_n n$, where $\frm \in \Gamma(\langle\mathbb{R}u\rangle \oplus \langle\mathbb{R}v\rangle)^{\perp_g}$.

 
\section{Supersymmetric NS-NS solutions}
\label{sec:susysolutions}

 
In this section we finally arrive to the study of the supersymmetric solutions of the four-dimensional NS-NS system on a bundle gerbe $\cC = (\cP,Y,\cA)$ and a principal $\mathbb{Z}$-bundle $\cX$ using the results obtained in previous sections of this chapter.  

\begin{definition}
A \emph{supersymmetric NS-NS solution} is a supersymmetric NS-NS configuration:
\begin{equation*}
(g,b,\phi,\varepsilon)\in \Conf_s(\cC,\cX)	
\end{equation*}

\noindent 
such that $(g,b,\phi)$ satisfies the NS-NS system \eqref{eq:NSNSsystem} on $(\cC,\cX)$.
\end{definition}
 
\noindent
By Corollary \ref{cor:susykundt}, supersymmetric NS-NS solutions define a special class of four-dimensional Kundt manifolds. We denote by $\Sol_s(\cC,\cX)$ the category of supersymmetric solutions on $(\cC,\cX)$. The starting point in our study of supersymmetric solutions is the characterization of supersymmetric configurations obtained in the previous section, see Proposition \ref{prop:iffsusyconf} and Proposition \ref{prop:furtheridentities} as well as Lemma  \ref{lemma:dvarphi} and  Lemma \ref{lemma:nablavarphi}. Building on this characterization of supersymmetric NS-NS configurations we have to impose the equations of the NS-NS system \eqref{eq:NSNSsystem}, namely the Einstein, Maxwell, and dilaton equations in \eqref{eq:NSNSsystem}. We begin with the dilaton equation.

\begin{lemma}
Let $(g,b,\phi,\varepsilon)\in \Conf_s(\cC,\cX)$ be a supersymmetric configuration with associated isotropic parallelism $[u,v,l,n]\in \mathfrak{F}_{sk}(M)$. Then, $(g,b,\phi)$ satisfies the dilaton equation in \eqref{eq:NSNSsystem} if and only if:
\begin{eqnarray}
\label{lemma:dilatoncondition}
l(\alpha_n) - n(\alpha_l) + \alpha_l^2 + \alpha^2_n = 0
\end{eqnarray}
	
\noindent
where $(u,v,l,n)\in [u,l,v,n]$ is any representative.
\end{lemma}

\begin{proof}
Using Proposition \ref{prop:furtheridentities} we compute:
\begin{equation*}
\nabla^{g\ast}\varphi_{\phi} + \vert\varphi_{\phi}\vert_g^2 + \vert\alpha_b\vert_g^2 = l(\alpha_n) - n(\alpha_l)  + \alpha_l^2 + \alpha_n^2 = 0
\end{equation*}

\noindent
since $\vert\varphi_{\phi}\vert_g^2 = \vert\alpha_b\vert_g^2 = \alpha_l^2 + \alpha_n^2$.
\end{proof} 

\noindent
Equation \eqref{lemma:dilatoncondition} can be understood as a differential equation for the \emph{b-field}, although it originates from the dilaton equation. This is because by Proposition \ref{prop:susydilatino}, the \emph{space-like} part of $\varphi_{\phi}$ is determined entirely by $b$ through $\alpha_b$.
\begin{lemma}
\label{lemma:Maxwellcondition}
Let $(g,b,\phi,\varepsilon)\in \Conf_s(\cC,\cX)$ be a supersymmetric configuration with associated isotropic parallelism $[u,v,l,n]\in \mathfrak{F}_{sk}(M)$. Then, $(g,b,\phi)$ satisfies the Maxwell equation in \eqref{eq:NSNSsystem} if and only if:
\begin{eqnarray*}
& l(\alpha_n) - n(\alpha_l) + \alpha_l^2 + \alpha^2_n = 0\\
& v(\alpha_l)   =  l(\frf_b)  + \alpha_l \kappa(l) + \alpha_n \rho(l) + c_{\phi} \alpha_l + \frf_{b} \alpha_n\\
& v(\alpha_n)  = n(\frf_b)  + \alpha_l \kappa(n) + \alpha_n \rho(n) + c_{\phi} \alpha_n - \frf_{b} \alpha_l
\end{eqnarray*}
	
\noindent
where $(u,v,l,n)\in [u,l,v,n]$ is any representative. In particular, if $(g,b,\phi,\varepsilon)$ satisfies the Maxwell equation then it necessarily satisfies the dilaton equation.
\end{lemma}

\begin{proof}
We impose the Maxwell equation in \eqref{eq:NSNSsystem} using the expression for the exterior derivative of $\alpha_b$ obtained in the proof of Proposition \ref{prop:furtheridentities}. We have:
\begin{eqnarray*}
& \dd\alpha = (l(\alpha_n) - n(\alpha_l)) l\wedge n + (v(\alpha_l) - \dd\kappa(l,n)) u\wedge l + (v(\alpha_n)-\dd\rho(l,n)) u\wedge n  \\
& = (\frc_{\phi} u + \alpha_l n - \alpha_n l)\wedge (\frf_b u + \alpha_l l + \alpha_n n)
\end{eqnarray*}

\noindent
Expanding this equation and isolating terms by type we obtain:
\begin{eqnarray*}
& l(\alpha_n) - n(\alpha_l) + \alpha_l^2 + \alpha^2_n = 0\\
& v(\alpha_l) - \dd\kappa(l,n) = c_{\phi} \alpha_l + \frf_{b} \alpha_n\, , \quad v(\alpha_n) - \dd\rho(l,n) = c_{\phi} \alpha_n - \frf_{b} \alpha_l
\end{eqnarray*}

\noindent
On the other hand, by Proposition \ref{prop:skewtorsionintegrabilityI}, we know that:
\begin{equation*}
\dd\kappa(l,n) = l(\frf_b)  + \alpha_l \kappa(l) + \alpha_n \rho(l)\, ,  \qquad \dd\rho(l,n) = n(\frf_b)  + \alpha_l \kappa(n) + \alpha_n \rho(n)
\end{equation*} 

\noindent
which plugged back into the previous equation gives the desired result. 
\end{proof}

\begin{lemma}
\label{lemma:Einsteincondition}
Let $(g,b,\phi,\varepsilon)\in \Conf_s(\cC,\cX)$ be a supersymmetric configuration with associated isotropic parallelism $[u,v,l,n]\in \mathfrak{F}_{sk}(M)$. Then, $(g,b,\phi)$ satisfies the Einstein equation in \eqref{eq:NSNSsystem} if and only if:
\begin{eqnarray*}
& l(\alpha_n) - n(\alpha_l) + \alpha_l^2 + \alpha^2_n = 0\\
& \dd\kappa(v,l) + \dd\rho(v,n) + v(\frc_{\phi}) - \alpha_n \kappa(v) + \alpha_l \rho(v) = 0\\
& v(\alpha_l)   =  l(\frf_b)  + \alpha_l \kappa(l) + \alpha_n \rho(l) + c_{\phi} \alpha_l + \frf_{b} \alpha_n\\
& v(\alpha_n)  = n(\frf_b)  + \alpha_l \kappa(n) + \alpha_n \rho(n) + c_{\phi} \alpha_n - \frf_{b} \alpha_l
\end{eqnarray*}
	
\noindent
where $(u,v,l,n)\in [u,l,v,n]$ is any representative. In particular, if $(g,b,\phi,\varepsilon)$ satisfies the Maxwell equation then it necessarily satisfies the dilaton equation.
\end{lemma}

\begin{proof}
We expand the Einstein equation in \eqref{eq:NSNSsystem} into four separate equations by evaluating it on a given istropic coframe $(u,l,v,n)\in [u,l,v,n]$. This gives:
\begin{eqnarray*}
& \mathrm{Ric}^g(u) + \nabla^g_u \varphi_{\phi} + \frac{1}{2} \vert\alpha_b\vert_g^2 u =0	\\
& \mathrm{Ric}^g(v) + \nabla^g_v \varphi_{\phi} - \frac{1}{2}  \frf_b\alpha_b  + \frac{1}{2} \vert\alpha_b\vert_g^2 v =0	\\
& \mathrm{Ric}^g(l) + \nabla^g_l \varphi_{\phi} - \frac{1}{2}  \alpha_l \alpha_b  + \frac{1}{2} \vert\alpha_b\vert_g^2 l =0	\\
& \mathrm{Ric}^g(n) + \nabla^g_n \varphi_{\phi} - \frac{1}{2}  \alpha_n \alpha_b  + \frac{1}{2} \vert\alpha_b\vert_g^2 n =0	
\end{eqnarray*}

\noindent
To proceed further, we write $\varphi_{\phi} = \frc_{\phi} + \alpha_l n - \alpha_n l$, $\varphi_{\phi} = \frc_{\phi} + \alpha n - \alpha_n l$ and we use Lemma \ref{lemma:nablavarphi} together with Proposition \ref{prop:curvaturessusy} to analyze each of the previous for equations separately. For the first equation, we simply have:
\begin{equation*}
\mathrm{Ric}^g(u) + \nabla^g_u \varphi_{\phi} + \frac{1}{2} \vert\alpha_b\vert_g^2 u  = \frac{1}{2} (l(\alpha_n) - n(\alpha_l) + \alpha^2_l + \alpha_n^2) u
\end{equation*}

\noindent
which gives the first equation of the lemma. For the second equation, we further split it in four components by evaluating on $(u,v,l,n)$. We have:
\begin{eqnarray*}
& \mathrm{Ric}^g(v,u) + (\nabla^g_v \varphi_{\phi})(u)  + \frac{1}{2} \vert\alpha_b\vert_g^2 = \frac{1}{2} (l(\alpha_n) - n(\alpha_l) + \alpha_l^2 + \alpha_n^2)\\
& \mathrm{Ric}^g(v,v) + (\nabla^g_v \varphi_{\phi})(v) - \frac{1}{2} \frf_b^2  = \dd\kappa(v,l) + \dd\rho(v,n) + v(\frc_{\phi}) - \alpha_n \kappa(v) + \alpha_l \rho(v)\\
& \mathrm{Ric}^g(v,l) + (\nabla^g_v \varphi_{\phi})(l) -\frac{1}{2} \frf_b \alpha_l  = \frac{1}{2} (\dd\rho(l,n) - v(\alpha_n) - \alpha_l \frf_b + \alpha_n \frc_{\phi})\\
& \mathrm{Ric}^g(v,n) + (\nabla^g_v \varphi_{\phi})(n) -\frac{1}{2} \frf_b \alpha_n  =  \frac{1}{2} (v(\alpha_l) - \dd\kappa(l,n)   - \alpha_n\frf_b - \alpha_l \frc_{\phi})
\end{eqnarray*} 

\noindent
The remaining components of the Einstein equation do not yield any new conditions, since:
\begin{eqnarray*}
& \mathrm{Ric}^g(l) + \nabla^g_l \varphi_{\phi} - \frac{1}{2}  \alpha_l \alpha_b  + \frac{1}{2} \vert\alpha_b\vert_g^2 l =  \frac{1}{2} (\dd\rho(l,n) - v(\alpha_n) - \alpha_l \frf_b + \alpha_n \frc_{\phi}) u\\	 
& \mathrm{Ric}^g(n) + \nabla^g_n \varphi_{\phi} - \frac{1}{2}  \alpha_n \alpha_b  + \frac{1}{2} \vert\alpha_b\vert_g^2 n =  \frac{1}{2} (v(\alpha_l) - \dd\kappa(l,n)   - \alpha_n\frf_b - \alpha_l \frc_{\phi}) u
\end{eqnarray*}

\noindent
Substituting now:
\begin{equation*}
\dd\kappa(l,n) = l(\frf_b)  + \alpha_l \kappa(l) + \alpha_n \rho(l)\, ,  \qquad \dd\rho(l,n) = n(\frf_b)  + \alpha_l \kappa(n) + \alpha_n \rho(n)
\end{equation*} 

\noindent
into the previous expressions we conclude.
\end{proof}

\noindent
From the previous lemmata, it easily follows that the Einstein, Maxwell, and dilaton equations of the NS-NS system \eqref{eq:NSNSsystem} are not fully independent when evaluated on supersymmetric configurations. More precisely, we have the following result.

\begin{cor}
If a supersymmetric NS-NS configuration satisfies the Einstein equation, then it also satisfies the Maxwell and dilaton equations. If a supersymmetric NS-NS configuration satisfies the Maxwell equation then it satisfies the dilaton equation.  
\end{cor}

\noindent
Schematically:
\begin{center}
Einstein equation $\,\, \Rightarrow\,\,$ Maxwell equation $\,\,\Rightarrow\,\,$ Dilaton equation
\end{center}

\noindent
Lemma \ref{lemma:Einsteincondition} contains all the equations that a supersymmetric NS-NS configuration needs to satisfy in order to be a supersymmetric NS-NS solution. The first equation in Lemma \ref{lemma:Einsteincondition} can be interpreted as a second-order equation for the $b$-field. The second equation in Lemma \ref{lemma:Einsteincondition} is the genuine differential equation that comes from imposing the Einstein equation of the NS-NS system (the remaining equations are equivalent to the Maxwell and dilaton equations). 

\begin{remark}
We can conclude that the \emph{genuine} Einstein equation of the NS-NS system evaluated on a supersymmetric configuration is the second equation in Lemma \ref{lemma:Einsteincondition}, since the other equations are implied by the Maxwell and dilaton equations of the NS-NS system. 
\end{remark}

\noindent
Using Lemma \ref{lemma:dvarphi}, the third and fourth equations in Lemma \ref{lemma:Einsteincondition} can be better analyzed by writing them in the following, equivalent, matrix form:

\begin{equation*}
 \left( \begin{array}{ccc} \alpha_l & \alpha_n \\ 
		\alpha_n & -\alpha_l \end{array} \right) \left( \begin{array}{ccc} \kappa(l) - \rho(n)  \\ 
		\kappa(n) + \rho(l)   \end{array} \right)
	= \left( \begin{array}{ccc} n(\frc_{\phi}) - l(\frf_b)  \\ 
		n(\frf_b) + l(\frc_{\phi}) \end{array} \right) 
\end{equation*}	

\noindent
where we have set $\alpha_l := \alpha(l)$ and $\alpha_n := \alpha(n)$ for ease of notation. The differential system of Lemma \ref{lemma:Einsteincondition} needs to be solved for appropriately chosen $\kappa ,\rho \in \Omega^1(M)$. Hence, at every point $m\in M$ such that:
\begin{equation*}
	\det \left( \begin{array}{ccc} \alpha_l & \alpha_n \\ 
		\alpha_n & -\alpha_l \end{array} \right) \left( \begin{array}{ccc} \kappa(l) - \rho(n)  \\ 
		\kappa(n) + \rho(l)   \end{array} \right)\vert_m = - (\alpha^2_l + \alpha_n^2)\vert_m \neq 0
\end{equation*}	

\noindent
the previous equations is solved \emph{algebraically} by:
\begin{equation*}
 \left( \begin{array}{ccc} \kappa(l) - \rho(n)  \\ 
		\kappa(n) + \rho(l)   \end{array} \right)\vert_m
	=  \frac{1}{\alpha_l^2 + \alpha_n^2}	\left( \begin{array}{ccc} \alpha_l & \alpha_n \\ 
		\alpha_n & -\alpha_l \end{array} \right) \left( \begin{array}{ccc} n(\frc_{\phi}) - l(\frf_b)  \\ 
		n(\frf_b) + l(\frc_{\phi}) \end{array} \right) \vert_m
\end{equation*}

\noindent
This condition turns out to be crucial for the NS-NS system, as the following result shows. 

\begin{lemma}
\label{lemma:constrainkapparho}
Let $(g,b,\phi,\varepsilon)\in \Conf_s(\cC,\cX)$ be a supersymmetric configuration with associated isotropic parallelism $[u,v,l,n]\in \mathfrak{F}_{sk}(M)$ relative to $\kappa , \rho \in \Omega^1(M)$. Then:
\begin{eqnarray*}
& (\alpha_l^2 + \alpha_n^2) (\kappa(u) + \alpha_n) + \alpha_l u(\frf_b) - \alpha_n u(\frc_{\phi})  = 0\\
& (\alpha_l^2 + \alpha_n^2) (\rho(u) - \alpha_l) + \alpha_n u(\frf_b) + \alpha_l u(\frc_{\phi})  = 0\\
& (\alpha_l^2 + \alpha_n^2)( \kappa(l) + \frc_{\phi}) - \frac{1}{2} v (\alpha_l^2 + \alpha_n^2) + \alpha_l l(\frf_b) - \alpha_n l(\frc_{\phi})  = 0\\
& (\alpha_l^2 + \alpha_n^2) (\rho(l) + \frf_b) + \alpha_l v(\alpha_n) - \alpha_n v(\alpha_l) + \alpha_n l(\frf_b) + \alpha_l l(\frc_{\phi})    = 0\\
& (\alpha_l^2 + \alpha_n^2)( \kappa(n) -  \frf_b)- \alpha_l v(\alpha_n) + \alpha_n v(\alpha_l) + \alpha_l n(\frf_b) - \alpha_n n(\frc_{\phi})   = 0\\
& (\alpha_l^2 + \alpha_n^2)( \rho(n) + \frc_{\phi}) - \frac{1}{2}v (\alpha_l^2 + \alpha_n^2) + \alpha_n n(\frf_b) + \alpha_l n(\frc_{\phi})   = 0
\end{eqnarray*}

\noindent
where we have set $\alpha_l := \alpha(l)$ and $\alpha_n := \alpha(n)$.
\end{lemma}

\begin{proof}
The result follows by appropriately combining Proposition \ref{prop:skewtorsionintegrabilityI} together with Lemmas \ref{lemma:dvarphi} and \ref{lemma:Einsteincondition}. More precisely, multiplying equation:
\begin{equation*}
\alpha_l \kappa(u) + \alpha_n \rho(u) + u(\frf_b) = 0 
\end{equation*}

\noindent
of Proposition \ref{prop:skewtorsionintegrabilityI} by $\alpha_l$ and adding the result to the multiplication of the first equation in \eqref{eq:dvarphi} by $\alpha_n$ we obtain the first equation in the statement. Multiplying the previous equation by $\alpha_n$ and combining the result with the multiplication of the first equation in \eqref{eq:dvarphi} by $\alpha_l$ we obtain the second equation in the statement. Multiplying the second equation in \eqref{eq:dvarphi} by $\alpha_n$ and combining the result with the the third equation in Lemma \ref{lemma:Einsteincondition} multiplied by $\alpha_l$ we obtain the third equation in the statement. Similarly, multiplying the second equation in \eqref{eq:dvarphi} by $\alpha_l$ and combining the result with the the third equation in Lemma \ref{lemma:Einsteincondition} multiplied by $\alpha_n$ we obtain the fourth equation in the statement. The remaining equations are obtained analogously and hence we conclude. 
\end{proof}

\noindent
At every point $m\in M$ such that $\vert\frm_b\vert_q^2\vert_m \neq 0$, the relations of the previous lemma conform an algebraic system for all the components of $\kappa$ and $\rho$ except for $\kappa(v)$ and $\rho(v)$, which is solved by isolating for them. The key question now is what is the structure of the set of points $m\in M$ such that $\vert\frm_b\vert_q^2\vert_m = 0$ or, equivalently, the set of points $m\in M$ such that $\alpha_l\vert_m = 0$ and $\alpha_n\vert_m = 0$. Define:
\begin{equation*}
N = \left\{ m\in M \,\, \vert\,\, \alpha_l\vert_m = 0\,\, \& \,\, \alpha_n\vert_m = 0\right\} \subset M
\end{equation*}

\noindent
Alternatively, we have $N:= \frm^{-1}_b(0)$. To proceed further we need to study more closely the Einstein equation on $M\backslash N$. It is convenient to define, associated to every supersymmetric solution $(g,b,\phi,\varepsilon)\in \Sol_s(\cC,\cX)$, the following smooth complex function:
\begin{equation*}
\frF \colon M \to \mathbb{C}\, , \qquad m \mapsto  \alpha_l(m) + i \alpha_n(m) 
\end{equation*}

\noindent
where $i\in \mathbb{C}$ denotes the imaginary unit. Note that $\frF\colon M \to \mathbb{C}$ does not depend on the representative chosen $(u,v,l,n)\in [u,v,l,n]$ in the isotropic parallelism determined by $\varepsilon$. The map $\frF \colon M \to \mathbb{C}$ is fundamental to understand the geometric and topological structure of supersymmetric solutions and is reminiscent of the \emph{triholomorphic moment map} underlying a complete hyper-Kähler four-manifold equipped with a local tri-hamiltonian $S^1$-action \cite{Bielawski}. This analogy becomes more transparent after noticing that $\frF \colon M \to \mathbb{C}$ descends to the quotient of $M$ by the action $\mathbb{R}\times M \to M$ generated by $u^{\sharp_g}\in\mathfrak{X}(M)$:
\begin{equation*}
\frF_u \colon M/\mathbb{R} \to \mathbb{C}
\end{equation*}

\noindent
where we are assuming that $u^{\sharp_g}\in\mathfrak{X}(M)$ is complete. This defines a continuous function on the topological space $M/\mathbb{R}$ equipped with the quotient topology. Further assumptions on this quotient give rise to a casuistry that can be used to develop a partial classification of supersymmetric solutions. We plan to explore this perspective in the future. Similarly, it is also convenient to introduce the complex function:
\begin{equation*}
\frG   \colon M \to \mathbb{C}\, , \qquad m \mapsto  \frc_{\phi}(m) + i \frf_b(m) 
\end{equation*}

\noindent
In contrast to $\frF \in C^{\infty}(M,\mathbb{C})$, as defined above $\frG \in C^{\infty}(M,\mathbb{C})$ does depend on the representative $(u,v,l,n)\in [u,v,l,n]$ chosen in the isotropic parallelism $[u,v,l,n]$ of the corresponding supersymmetric solution. To deal with this dependence, consider the following equivalence relation in the complex functions $C^{\infty}(M,\mathbb{C})$:
\begin{equation*}
f_1 \sim f_2 \quad \mathrm{iff}\quad f_1 = f_2  - \bar{\frF} F
\end{equation*}

\noindent
for a complex function $F \in C^{\infty}(M,\mathbb{C})$. We denote by:
\begin{equation*}
 C^{\infty}_{\frF}(M,\mathbb{C}) :=  C^{\infty}(M,\mathbb{C})/\frF
\end{equation*}

\noindent
the corresponding quotient, which inherits the structure of a ring and a multiplicative module over  $C^{\infty}(M,\mathbb{C})$. Then, a supersymmetric solution determines a class $[\frG]_{\frF} \in C^{\infty}_{\frF}(M,\mathbb{C})$.  Evidently, if $\frF$ is nowhere vanishing, then $ C^{\infty}_{\frF}(M,\mathbb{C}) = \left\{ 0\right\}$ is the trivial ring. In practice, this means that we can choose a representative $(u,v,l,n)\in [u,v,l,n]$ in terms of which $\frf_b =  \frc_{\phi} = 0$. Whereas this is always possible locally around every point $m\in M$ at which $\vert\alpha_b\vert^2_g\vert_m \neq 0$, it might be globally obstructed. This obstruction is measured by both $[\frG]_{\frF}\in C^{\infty}_{\frF}(M,\mathbb{C})$ and $C^{\infty}_{\frF}(M,\mathbb{C})$, which becomes an \emph{invariant} of the given supersymmetric solution.  

\begin{remark}
By the previous discussion, we obtain a new \emph{invariant} associated to a supersymmetric NS-NS solution, namely the ring $C^{\infty}_{\frF}(M,\mathbb{C})$. This complements the rank-one invariants defined by $\kappa$ and $\rho$ as described in Section \ref{sec:generalcurvature} in Chapter \ref{chapter:parallelspinorstorsion}.
\end{remark}

\noindent
Before studying the geometry and topology of general supersymmetric NS-NS solutions, we first need to consider separately the case $\alpha_l = \alpha_n = 0$, to which we will refer as \emph{null}.

\begin{definition}
A supersymmetric solution $(g,b,\phi,\varepsilon)\in \Sol_s(\cC,\cX)$ is \emph{null} if both $\varphi_{\phi}$ and $\alpha_{b}$ are collinear with the Dirac current $u\in \Omega^1(M)$ associated to $\varepsilon$.
\end{definition}

\noindent
We will refer to a supersymmetric solution $(g,b,\phi,\varepsilon)\in \Sol_s(\cC,\cX)$ as having \emph{everywhere space-like flux} if $\vert \alpha_b\vert^2_g\in C^{\infty}(M)$ is a nowhere vanishing function on $M$, which is then necessarily positive. The null case is therefore \emph{complementary} to the nowhere vanishing flux case, and both constitute the building blocks of NS-NS supersymmetric solutions. Note that if $(g,b,\phi,\varepsilon)\in \Sol_s(\cC,\cX)$ is a null solution with associated isotropic parallelism $[u,v,l,n]$, then Proposition \ref{prop:existencenullcoframeII} implies that $(M,g,u)$ is a Brinkmann four-manifold, that is, $(M,g)$ is a Lorentzian four-manifold and $u\in \Omega^1(M)$ is a parallel with respect to the Levi-Civita connection on $(M,g)$.  

\begin{thm}
\label{thm:susynullsolutions}
An oriented and strongly spin four-manifold $M$ admits a null supersymmetric NS-NS solution for a certian pair $(\cC,\cX)$ if and only if it admits a global coframe $(u,v,l,n)$ satisfying the following differential system:    
\begin{eqnarray}
\label{eq:nulleqs} 
\dd u = 0 \, , \quad \dd v = -\kappa\wedge l-\rho\wedge n - \frf\, l\wedge n\, , \quad \dd l = u \wedge (\frf\, n  - \kappa)  \, , \quad \dd n = - u \wedge (\frf\, l + \rho)   
\end{eqnarray}
	
\noindent
for a pair of closed one-forms $\kappa,\rho \in \Omega^1(M)$ and constants $\frf , \frc \in \mathbb{R}$ satisfying the following cohomological conditions:
\begin{equation*}
\frac{1}{2\pi}  [\frc\, u]\in H^1(M,\mathbb{Z})\, , \qquad \frac{1}{2\pi} [\frf \ast_g u] \in H^3(M,\mathbb{Z})
\end{equation*}
	
\noindent
In particular, $(M,g)$ is a Brinkmann four-manifold with isotropic Ricci curvature. 
\end{thm} 

\begin{proof}
By Proposition \ref{prop:iffsusyconf} and Lemma \ref{lemma:Einsteincondition}, $M$ admits a null supersymmetric solution for some pair $(\cC,\cX)$ consisting of a bundle gerbe $\cC$ and a principal $\mathbb{Z}$-bundle $\cX$ if and only if there exists an isotropic parallelism $[u,v,l,n]$, a pair of functions $\frc , \frf \in C^{\infty}(M)$, and a pair of one-forms $\kappa , \rho \in \Omega^1(M)$ satisfying the differential system \eqref{eq:nulleqs}:
\begin{equation*}
\dd u = 0 \, , \quad \dd v = -\kappa\wedge l-\rho\wedge n - \frf\, l\wedge n\, , \quad \dd l = u \wedge (\frf\, n  - \kappa)  \, , \quad \dd n = - u \wedge (\frf\, l + \rho)  
\end{equation*}

\noindent
which corresponds to the condition that $[u,v,l,n]$ is skew-torsion with torsion $H = \frf \ast_g u$, together with the differential system:
\begin{eqnarray*}
\dd\kappa(v,l) + \dd\rho(v,n) + v(\frc_{\phi})   = 0\, , \qquad   l(\frf_b) = 0\, , \qquad n(\frf_b)  = 0
\end{eqnarray*}

\noindent
and the condition that the following define integral classes:
\begin{equation*}
\frac{1}{2\pi}[\frc\, u] \in H^1(M,\mathbb{Z})\, , \qquad \frac{1}{2\pi}[\frf \ast_g u] \in H^3(M,\mathbb{Z})
\end{equation*}

\noindent
This gives the \emph{if} condition. For the \emph{only if} condition, suppose that $(g,b,\phi,\varepsilon) \in \Sol_s(\cC,\cX)$ is a null supersymmetric solution with associated skew-torsion isotropic parallelism $[u,v,l,n]$. Then, by Proposition \ref{prop:furtheridentities}, Lemma \ref{lemma:dvarphi} and Lemma \ref{lemma:Maxwellcondition} we have:
\begin{equation*}
u(\frc) = 0\, , \qquad l(\frc) = 0 \, , \qquad n(\frc) = 0\, , \qquad u(\frf) = 0\, , \qquad l(\frf) = 0 \, , \qquad n(\frf) = 0
\end{equation*}
 
\noindent
so the non-vanishing derivatives of $\frc$ and $\frf$ can only happen along $v$. The previous relations together with Lemma \ref{lemma:skewtorsionintegrabilityI} give the following \emph{integrability conditions}:
\begin{eqnarray}
u\wedge (\dd\kappa - \dd\frf_b \wedge n) = 0\, , \quad u\wedge (\dd\rho + \dd\frf_b \wedge l ) = 0\, , \quad \dd\kappa\wedge l + \dd\rho \wedge n  + \dd \frf_b \wedge l \wedge n = 0\label{eq:auxnulli}
\end{eqnarray}
 
\noindent
which in turn imply:
\begin{eqnarray}
& \dd\kappa(l,n) = 0 \, , \quad \dd\kappa(u,n) = 0 \, , \quad \dd\kappa(u,l) = 0\, , \quad \dd\kappa(u,v) = 0 \nonumber \\
& \dd\rho(l,n) = 0  \, , \quad \dd\rho(u,n) = 0 \, , \quad \dd\rho(u,l) = 0\, , \quad \dd\rho(u,v) = 0  \label{eq:auxnullii}\\
& \dd \rho (v,l) - \dd \kappa(v,n)  + v(\frf_b)  = 0 \nonumber
\end{eqnarray}

\noindent
From Equation \eqref{eq:auxnulli} we obtain that:
\begin{eqnarray*}
\dd\kappa = \dd\frf_b \wedge n + u \wedge \omega_{\kappa}\, , \qquad \dd\rho = - \dd\frf_b \wedge l + u \wedge \omega_{\rho}
\end{eqnarray*}

\noindent
for certain one-forms $\omega_{\kappa} , \omega_{\rho} \in \Omega^1(M)$. Imposing the first two lines in \eqref{eq:auxnullii} we conclude that:
\begin{eqnarray*}
 \omega_{\kappa} = \omega_{\kappa}(v) u - v(\frf)\, n \, , \qquad  \omega_{\rho} = \omega_{\kappa}(v) u + v(\frf)\, l
\end{eqnarray*}

\noindent
and therefore $\dd\kappa = 0$ and $\dd\rho = 0$. This, together with the third line in \eqref{eq:auxnullii} implies that $v(\frf) = 0$ and thus $\frf$ is constant. On the other hand, Lemma \ref{lemma:Einsteincondition} implies:
\begin{eqnarray*}
\dd\kappa(v,l) + \dd\rho(v,n) + v(\frc) =  v(\frc) = 0
\end{eqnarray*}

\noindent
and consequently $\frc$ is also constant. Since satisfying all equations in Lemma \ref{lemma:Einsteincondition} is equivalent to satisfying the Einstein equation on an supersymmetric configuration we conclude that $\Ric^g + \frac{\frf^2}{2} u\otimes u = 0$ and therefore $\Ric^g$ is isotropic, namely it satisfies $\vert \Ric^g\vert^2_g = 0$.
\end{proof}

\noindent
The previous theorem together with Proposition \ref{prop:curvaturetorsion} implies the following corollary, which can be understood as a characterization of the universal cover of supersymmetric null solutions.

\begin{cor}
A simply connected Lorentzian four-manifold $(M,g)$ admits a compatible null supersymmetric solution $(g,b,\phi)$ if and only if it admits a parallel isotropic one-form $u\in \Omega^1(M)$ and there exists a constant $\frf \in \mathbb{R}$ such that $\nabla^{g,\alpha}$ is flat with $\alpha = \frf\, u$. 
\end{cor}

\noindent
Similarly, we have the following result as an immediate consequence of Theorem \ref{thm:susynullsolutions}. 

\begin{cor}
A Brinkmann four-manifold $(M,g,u)$ locally admits a compatible null supersymmetric solution if and only if there exists a constant $\frf \in \mathbb{R}$ such that the metric connection $\nabla^{g,\alpha}$ with torsion $\alpha = \frf\, u$ is flat.
\end{cor}

\begin{remark}
By $(M,g,u)$ \emph{locally admitting} a compatible null supersymmetric solution we mean that around every point in $M$ there exists an open set admitting a null supersymmetric solution with metric $g$ and Dirac current $u$. 
\end{remark}

\noindent
These corollaries motivate introducing the following notion of Brinkmann four-manifold. 

\begin{definition}
A \emph{torsion-flat} Brinkmann four-manifold is a Brinkmann manifold $(M,g,u)$ for which there exists a constant $\frf$ such that  $\nabla^{g,H}$ is flat with $H = \frf\, \ast_g u$.
\end{definition}

\noindent
This seems like an intrinsically interesting class of Brinkmann four-manifolds to study, with the added bonus that they are locally supersymmetric. To the best of our knowledge, these have not been systematically considered yet in the literature. We will consider the local structure of null supersymmetric solutions in Section \ref{sec:dilatonflux}. Having dealt with the null case, we consider now the flux supersymmetric solutions with everywhere space-like flux.  

\begin{thm}
\label{thm:susyfluxsolutions}
An oriented and strongly spin four-manifold $M$ admits a supersymmetric NS-NS solution $(g,b,\phi,\varepsilon)\in \Sol_s(\cC,\cX)$ with everywhere space-like flux on a certain pair $(\cC,\cX)$ if and only if it admits a global coframe $(u,v,l,n)$ and a pair of complex functions $\frF, \frK \in C^{\infty}(M,\mathbb{C})$ satisfying the following differential system:    
 \begin{eqnarray}
& \dd u =   u\wedge \Im(\frF \bar{\Psi}) \, , \quad \dd v +   u\wedge \Re(\frK \bar{\Psi}) + 2 \frF \bar{\frF}  \Im (v(\bar{\frF}) \frF) \nu_q = 0\, , \quad \bar{\frF} \dd \Psi + v(\bar{\frF}) u \wedge \Psi = 0 \quad \label{eq:spaceflux1} \\
&   \bar{\Psi}(\frF) + \frF\bar{\frF} = 0\, , \qquad 2 v((\bar{\frF}\frF)^{-1}\Re(v(\bar{\frF})\frF)) = \Re(\Psi(\bar{\frK})) + \Im(\frF \bar{\frK}) \label{eq:spaceflux2}
\end{eqnarray}

\noindent
together with the following cohomological conditions:
\begin{equation*}
\frac{1}{2\pi}[\Im(\bar{\frF} \Psi)]\in H^1(M,\mathbb{Z})\, , \qquad \frac{1}{2\pi} [\ast_g \Re(\bar{\frF}\Psi)] \in H^3(M,\mathbb{Z})
\end{equation*}

\noindent
where $g= u\odot v + l\otimes l + n\otimes n$ is the Lorentzian metric associated to $(u,v,l,n)$ and $\Psi := l +i n\in \Omega^1(M,\mathbb{C})$. If that is the case, the supersymmetric solutions $(g,b,\phi)$ associated to such $(u,v,l,n)$ and $(\frF,\frK)$ are determined by:
\begin{eqnarray*}
g = u\odot v + l\otimes l + n\otimes n\, , \qquad \alpha_b = \ast_g \Re(\bar{\frF}\Psi)\, , \qquad \varphi_{\phi} = \Im(\bar{\frF} \Psi) 
\end{eqnarray*}

\noindent
where $\varphi_{\phi}$ is the curvature of $\phi$ and $\alpha_b$ is the Hodge dual of the curvature of $\phi$.
\end{thm}

\begin{proof}
Let $(g,b,\phi,\varepsilon)\in \Sol_s(\cC,\cX)$ be a NS-NS supersymmetric solution with nowhere vanishing flux and let $[u,v,l,n]$ be its associated skew-torsion isotropic parallelism with torsion $\alpha_b \in \Omega^1(M)$, which by Proposition \ref{prop:iffsusyconf} satisfies the differential system \eqref{eq:reducedintegrability1} and \eqref{eq:reducedintegrability2} for certain characteristic one-forms $\kappa , \rho \in \Omega^1(M)$. Since by assumption $\vert \frm_b \vert^2_q\in C^{\infty}(M)$ is nowhere vanishing, there exists a unique isotropic coframe $(u,v,l,n) \in [u,v,l,n]$ such that:
\begin{equation*}
\alpha_b = \alpha_l l + \alpha_n n\, , \qquad  \varphi_{\phi} = \alpha_l n - \alpha_n l
\end{equation*}

\noindent
that is, such that $\frc_{\phi} = \frf_b = 0$. By Lemma \ref{lemma:constrainkapparho} the characteristic one-forms $\kappa , \rho \in \Omega^1(M)$ satisfy the following relations:
\begin{eqnarray*}
& \kappa(u) = - \alpha_n   \, , \qquad     \rho(u) = \alpha_l  \\
&  \kappa(l)  =   \frac{1}{2} \vert \frm_b \vert^{-2}_q   v (\vert \frm_b \vert^2_q) \, , \qquad  \rho(l) = \vert \frm_b \vert^{-2}_q  (\alpha_n v(\alpha_l) - \alpha_l v(\alpha_n) )   \\
&  \kappa(n)   = \vert \frm_b \vert^{-2}_q  (\alpha_l v(\alpha_n) - \alpha_n v(\alpha_l)) \, , \qquad\rho(n)   =  \frac{1}{2} \vert \frm_b \vert^{-2}_q  v (\vert \frm_b \vert^2_q) 
\end{eqnarray*}
	
\noindent
where we have used that $\vert \frm_b \vert^2_q\in C^{\infty}(M)$ is nowhere vanishing. Using the previous relations we can \emph{solve} for all components of $\kappa$ and $\rho$ except for $\kappa(v)$ and $\rho(v)$. Plugging these relations into the differential system given in equations \eqref{eq:reducedintegrability1} and \eqref{eq:reducedintegrability2}, we obtain: 
\begin{eqnarray*}
	& \dd u = u\wedge (\alpha_n l - \alpha_l n) \\
	& \dd v = -\kappa_v u\wedge l-\rho_v u\wedge n + 2 (\alpha_l^2 + \alpha_n^2) (\alpha_l v(\alpha_n) - \alpha_n v(\alpha_l)) l\wedge n\\
	& (\alpha_l^2 + \alpha_n^2) \dd l = -  u \wedge (\frac{1}{2} v(\alpha^2_l+\alpha_n^2) l + (\alpha_l v(\alpha_n) - \alpha_n v(\alpha_l)) n)  \\
	& (\alpha_l^2 + \alpha_n^2) \dd n = -  u \wedge ( (\alpha_n v(\alpha_l) - \alpha_l v(\alpha_n)) l + \frac{1}{2} v(\alpha^2_l+\alpha_n^2) n)
\end{eqnarray*}

\noindent
Using the complex function $\frF = \alpha_l + i \alpha_n$ introduced earlier, the previous system can be written as follows:
\begin{eqnarray*}
	& \dd u = u\wedge (\alpha_n l - \alpha_l n) \\
	& \dd v = -\kappa_v u\wedge l-\rho_v u\wedge n - 2 \frF \bar{\frF}  \Im (v(\bar{\frF}) \frF) l\wedge n\\
	&  \frF\, \bar{\frF}  \dd l = -  u \wedge (\Re (v(\bar{\frF}) \frF) l - \Im (v(\bar{\frF}) \frF) n)  \\
	& \frF\, \bar{\frF} \dd n = -  u \wedge ( \Im (v(\bar{\frF}) \frF) l + \Re (v(\bar{\frF}) \frF) n)
\end{eqnarray*} 

\noindent
or, equivalently using the complex tetrad $(u,v,\Psi = l + in,\bar{\Psi} = l - in)$ associated to $(u,v,l,n)$, as follows:
\begin{eqnarray*}
\dd u =   u\wedge \Im(\frF \Psi^{\ast}) \, , \quad \dd v +   u\wedge \Re(\frK \Psi^{\ast}) + 2 \frF \bar{\frF}  \Im (v(\bar{\frF}) \frF) \nu_q = 0\, , \quad \bar{\frF}  \dd \Psi + v(\bar{\frF}) u \wedge \Psi = 0   
\end{eqnarray*} 

\noindent
where we have introduced the complex function:
\begin{equation*}
\frK := \kappa(v) + i\rho(v)
\end{equation*}

\noindent
and $\nu_q = l\wedge n$ denotes the \emph{transverse} Riemannian volume form associated to $(u,v,l,n)$. Since $(g,b,\phi,\varepsilon)\in \Sol_s(\cC,\cX)$ is supersymmetric, it satisfies the NS-NS system, which by Lemma \ref{lemma:Einsteincondition} in our current situation is equivalent to:
 \begin{eqnarray}
& l(\alpha_n) - n(\alpha_l) + \alpha_l^2 + \alpha^2_n = 0 \nonumber \\
& \dd\kappa(v,l) + \dd\rho(v,n)   - \alpha_n \kappa(v) + \alpha_l \rho(v) = 0 \label{eq:auxiliarmotion}\\
& v(\alpha_l)   =   \alpha_l \kappa(l) + \alpha_n \rho(l) \, , \quad v(\alpha_n)  =   \alpha_l \kappa(n) + \alpha_n \rho(n)  \nonumber
\end{eqnarray}

\noindent
The first line in \eqref{eq:auxiliarmotion} can be written as follows:
\begin{equation*}
l(\alpha_n) - n(\alpha_l) + \alpha_l^2 + \alpha^2_n  = \Psi^{\ast}(\frF) + \frF\bar{\frF} = 0
\end{equation*}

\noindent
where we are implicitly using that $l(\alpha_l) + n(\alpha_n) = 0$ by Lemma \ref{lemma:skewtorsionintegrabilityI}. Furthermore, the third line in \eqref{eq:auxiliarmotion} can be written in the following form:
\begin{eqnarray*}
	& v(\alpha_l)   =   \alpha_l \kappa(l) + \alpha_n \rho(l) =   (\frF \bar{\frF})^{-1}( \alpha_l  \Re(v(\bar{\frF})\frF) + \alpha_n \Im(v(\bar{\frF})\frF) )  \\
	&  v(\alpha_n)  =   \alpha_l \kappa(n) + \alpha_n \rho(n) = (\frF \bar{\frF})^{-1}( \alpha_n  \Re(v(\bar{\frF})\frF) - \alpha_l \Im(v(\bar{\frF})\frF) )
\end{eqnarray*}

\noindent
from which we obtain:
\begin{eqnarray*}
& \frF \bar{\frF}\, v(\frF) = \frF \bar{\frF}\, v(\alpha_l + i\alpha_n) =  \alpha_l  \Re(v(\bar{\frF})\frF) + \alpha_n \Im(v(\bar{\frF})\frF)   + i( \alpha_n  \Re(v(\bar{\frF})\frF) - \alpha_l \Im(v(\bar{\frF})\frF) )\\
& = \alpha_l (\Re(v(\bar{\frF})\frF) - i \Im(v(\bar{\frF})\frF)) + \alpha_n (i\Re(v(\bar{\frF})\frF) +  \Im(v(\bar{\frF})\frF)) \\
& = \alpha_l (\Re(v(\bar{\frF})\frF) - i \Im(v(\bar{\frF})\frF)) + i \alpha_n ( \Re(v(\bar{\frF})\frF) - i \Im(v(\bar{\frF})\frF)) = v(\frF )\bar{\frF} \frF
\end{eqnarray*}

\noindent
and therefore we conclude that the third line in \eqref{eq:auxiliarmotion} holds automatically if the the first line is satisfied. It remains to consider the second line in Equation \eqref{eq:auxiliarmotion}. We have: 
\begin{eqnarray*}
& \dd\kappa(v,l) + \dd\rho(v,n) + v(\frc_{\phi}) - \alpha_n \kappa(v) + \alpha_l \rho(v) \\
& =  v (\kappa(l)) - l(\kappa(v)) - \kappa([v,l]) + v (\rho(l)) - l(\rho(v)) - \rho([v,l]) + v(\frc_{\phi}) - \alpha_n \kappa(v) + \alpha_l \rho(v) = 0
\end{eqnarray*}

\noindent
which implies:
\begin{eqnarray*}
& 0 = \dd\kappa(v,l) + \dd\rho(v,n)   - \alpha_n \kappa(v) + \alpha_l \rho(v) = \\
& = v(\kappa(l)) - l(\Re(\frK)) + v(\rho(n)) - n(\Im(\frK)) - \alpha_n \Re(\frK) + \alpha_l\Im(\frK)\\
& = 2 v((\bar{\frF}\frF)^{-1}\Re(v(\bar{\frF})\frF)) - \Re(\Psi(\bar{\frK})) - \Im(\frF \bar{\frK})
\end{eqnarray*}

\noindent
Here we have used that the Lie brackets of $(u,v,l,n)\in [u,v,l,n]$ are given by:
\begin{eqnarray*}
	& \left[ u^{\sharp_g} , v^{\sharp_g} \right]= - (\alpha_n + \kappa(u))\, l^{\sharp_g} + (\alpha_l - \rho(u))\, n^{\sharp_g} = 0\\ 
	& \left[ u^{\sharp_g} , l^{\sharp_g}\right] = (\kappa(u) + \alpha_n)\, \mu^{\sharp_g} = 0 \, , \quad \left[ u^{\sharp_g} , n^{\sharp_g}\right] = (\rho(u) - \alpha_l)\, \mu^{\sharp_g} = 0  \\ 
	& \left[ v^{\sharp_g} , l^{\sharp_g}\right] = - \alpha_n v + \kappa(v) u + \kappa(l) l + \rho(l) n\\
	& \left[ v^{\sharp_g} , n^{\sharp_g}\right] =  \alpha_l v + \rho(v) u + \kappa(n) l + \rho(n) n\\
	& \left[ l^{\sharp_g} , n^{\sharp_g}\right] =   (\rho(l) - \kappa(n)) u 
\end{eqnarray*}

\noindent
These relations follow from Remark \ref{remark:Liebrackets} after setting $\alpha(u) = 0$. This gives all the equations and conditions in the stating of the theorem. For the converse, suppose that we have a solution $(u,v,l,n,\frF,\frK)$ to all the equations in the statement of the theorem. Then, $(u,v,l,n)$ satisfies equations \eqref{eq:reducedintegrability1} and \eqref{eq:reducedintegrability2} relative to the one-forms:
\begin{eqnarray*}
\kappa + i \rho = \frK\, u +  \bar{\frF}\, v +  \bar{\frF}^{-1}  v(\bar{\frF}) \Psi 
\end{eqnarray*}

\noindent
Hence, $(u,v,l,n)$ defines a skew-torsion isotropic parallelism $[u,v,l,n]$ with torsion $\alpha = \Re(\bar{\frF} \Psi)$, and therefore the underlying Lorentzian metric $g = u\odot v + l\otimes l + n\otimes n$ admits a skew-torsion parallel spinor with torsion $\alpha = \Re(\bar{\frF} \Psi)$. Since the normalized Hodge dual of the latter is by assumption an integral closed three-form it follows that there exists a bundle gerbe $\cC = (\cP , Y , \cA)$ and a curving $b$ on $\cC$ such that $\alpha = \alpha_b$. Similarly, since by assumption $(2\pi)^{-1} \Im (\bar{\frF} \Psi)\in \Omega^1(M)$ is closed integral one-form on $M$, there exists a principal $\mathbb{Z}$-bundle $\cX$ on $M$ and an equivariant real function $\phi \in C^{\infty}(M)$ such that $\varphi_{\phi}  = \Im(\bar{\frF} \Psi)$. With this choices, all the conditions stated in Proposition \ref{prop:susyconfiguration} are automatically satisfied and therefore $(g,b,\phi)$ defines a supersymmetric configuration on $(\cC,\cX)$ with respect to the the spinor determined (modulo a global sign) by $[u,v,l,n]$. Finally, the following equations:
\begin{eqnarray*}
\bar{\Psi}(\frF) + \frF\bar{\frF} = 0\, , \qquad 2 v((\bar{\frF}\frF)^{-1}\Re(v(\bar{\frF})\frF)) = \Re(\Psi(\bar{\frK})) + \Im(\frF \bar{\frK})
\end{eqnarray*}

\noindent
that hold by assumption are equivalent to the NS-NS system evaluated on the given supersymmetric configuration, and therefore we conclude.
\end{proof}
 
\begin{remark}
In a supergravity context, the \emph{supersymmetry parameter} $\varepsilon$, that is, the underlying skew-torsion parallel spinor for the case of supersymmetric NS-NS solutions, is usually irrelevant in itself; only the geometric and topological consequences of its existence being important. Within this mindset, the previous corollary captures the if and only if conditions for a supersymmetry solution to exist with no mention of the underlying Lorentzian metric or supersymmetry parameter. This seems to realize, at least for the class of supersymmetric solutions considered here, the motivation and ideology explained in \cite{Tomasiello:2011eb}, where the \emph{complete} Type-II theory, of which NS-NS supergravity is a subsector, is considered in ten dimensions.
\end{remark}  

\noindent
By Proposition \ref{prop:susyconfiguration}, for every supersymmetric configuration $(g,b,\phi,\varepsilon)$ we have $\alpha_b(u) = 0$ and therefore by the first equation in \eqref{eq:reducedintegrability1} the Dirac current $u\in \Omega^1(M)$ satisfies the Frobenius integrability and the kernel $\Ker(u)\subset TM$ integrates to a codimension-one transversely orientable foliation $\cF_u \subset M$. The Godbillon-Vey class of this foliation is given by:
\begin{equation*}
	\sigma_u := [(\alpha_n l - \alpha_l n)\wedge \dd (\alpha_n l - \alpha_l n)] \in H^3(M,\mathbb{R})
\end{equation*}

\noindent
Hence, as an immediate consequence of the cohomological conditions contained in Theorem \ref{thm:susyfluxsolutions} we obtain the following result.

\begin{cor}
The Godbillon-Vey invariant of a supersymmetric NS-NS solution with everywhere space-like flux vanishes.
\end{cor}

\noindent


\section{The dilaton-flux foliation}
\label{sec:dilatonflux}


Every NS-NS solution $(g,b,\phi)\in \Sol(\cC,\cX)$ defines three natural distributions in $TM$, generally singular and possibly trivial, namely:

\begin{itemize}
	\item The \emph{dilaton distribution} $\Ker(\varphi_{\phi})\subset TM$, given by the kernel the one-form $\varphi_{\phi} \in \Omega^1(M)$. This distribution is singular at the zeroes of $\varphi_{\phi}$ and is trivial for constant-dilaton NS-NS solutions.
	
	\item The \emph{flux distribution} $\Ker(\alpha_b)\subset TM$, given by the kernel the one-form $\alpha_{b} \in \Omega^1(M)$. This distribution is singular at the zeroes of $\alpha_{b}$ and it is trivial for constant-dilaton NS-NS solutions.
	
	\item The \emph{dilaton-flux distribution} $\Ker(\varphi_{\phi})\cap\Ker(\alpha_b)$, given by the intersection of $\Ker(\varphi_{\phi})$ and $\Ker(\alpha_b)$. This is therefore a codimension-two distribution at those points, if any, where both $\varphi_{\phi}$ and $\alpha_b$ are non-zero and linearly independent.  
\end{itemize}

\noindent
The main interesting point about these singular distributions is that all of them are integrable, since $\varphi_{\phi}$ is closed and $\alpha_b$ satisfies the \emph{Cartan criteria} $\dd \alpha_b = \varphi_{\phi} \wedge \alpha_b$ precisely with respect to $\varphi_{\phi}$. In particular, assuming that the flux distribution is regular, then its \emph{Godbillon-Vey class} vanishes since $\varphi_{\phi}$ is closed. For nontrivial null supersymmetric solutions the dilaton-flux distribution degenerates into the flux or dilaton distributions, or both if  both $\alpha_b$ and $\varphi_{\phi}$ are nowhere vanishing, and therefore becomes the standard distribution and associated foliation that every Brinkmann manifold carries. For supersymmetric NS-NS solutions with everywhere space-like flux the dilaton-flux distribution is regular and therefore defines a regular codimension-two transversely orientable foliation. We consider both these cases in this section. 


\subsection{Everywhere isotropic flux}


Let $(g,b,\phi)$ be a null supersymmetric NS-NS solution with associated skew-torsion isotropic parallelism $[u,v,l,n]$. Then, every representative $(u,v,l,n) \in [u,v,l,n]$ satisfies the exterior differential system \eqref{eq:nulleqs} for a non-zero real constant $\frf \in \mathbb{R}^{\ast}$ and a pair of closed one-forms $\kappa , \rho \in \Omega^1(M)$. Around every point in $M$, there exists an open set isomorphic to $U\times X$, where $U\subset \mathbb{R}^2$ with Cartesian coordinates $(x_u,x_v)$, and $X$ is an oriented two-dimensional manifold and through this isomorphism there exists a representative of $(u,v,l,n) \in [u,v,l,n]$ of the form:
\begin{equation}
\label{eq:adaptednullflux}
(u,v,l,n) =   (\dd x_u   , \dd x_v +  \cH_{x_u}  \dd x_u + \omega_{x_u}, l_{x_u} , n_{x_u})
\end{equation}
 
\noindent
to which we will refer as being \emph{adapted}, similarly to the notion of adapted representative introduced in  Section \ref{sec:sktorsionparallelspinors}. Here $\cH_{x_u}$ is a family of functions on $X$ parametrized by $x_u$ and independent of $x_v$ whereas $\omega_{x_u}, l_{x_u} , n_{x_u} \in \Omega^1(X)$ are families of one-forms on $X$ parametrized by $x_u$ and again independent of $x_v$. Note that $u^{\sharp_g} = \partial_{x_v}$. Plugging \eqref{eq:adaptednullflux} into  \eqref{eq:nulleqs} we obtain:
\begin{eqnarray*}
& \dd_X \cH_{x_u} \wedge \dd x_u + \dd x_u \wedge \partial_{x_u}\omega_{x_u} + \dd_X \omega_{x_u} = -\kappa\wedge l_{x_u}-\rho\wedge n_{x_u} - \frf\, l_{x_u}\wedge n_{x_u}\\
& \dd x_u \wedge \partial_{x_u} l_{x_u} + \dd_X l_{x_u} = \dd x_u \wedge (\frf\, n_{x_u}  - \kappa^{\perp})  \\
& \dd x_u \wedge  \partial_{x_u} n_{x_u} + \dd n_{x_u} = -  \dd x_u \wedge (\frf\, l_{x_u} + \rho^{\perp})   
\end{eqnarray*}

\noindent
as well as $\kappa_v = \rho_v = 0$. Isolating by tensor type, the previous system is equivalent to:
\begin{eqnarray*}
& \kappa^{\perp} = \frf\, n_{x_u} - \partial_{x_u} l_{x_u}\, , \qquad \rho^{\perp} = - \frf\, l_{x_u} - \partial_{x_u} n_{x_u}\, , \qquad \dd_X l_{x_u} = \dd_x n_{x_u} = 0\\
&    \partial_{x_u}\omega_{x_u} - \dd_X \cH_{x_u} =  \partial_{x_u} l_{x_u}\wedge l_{x_u} + \partial_{x_u} n_{x_u}\wedge n_{x_u} \, , \qquad \dd_X \omega_{x_u} = \frf\, l_{x_u}\wedge n_{x_u}
\end{eqnarray*}

\noindent
Assuming for simplicity that $\partial_{x_u} l_{x_u} = c^l_{x_u} l_{x_u}$ and $\partial_{x_u} n_{x_u} = c^n_{x_u} n_{x_u}$ for families of functions $c^l_{x_u}$ and $c^n_{x_u}$ on $X$, the previous equations are equivalent to:
\begin{eqnarray*}
\omega_{x_u} = \int_0^{x_u} \dd_X\cH_{z} \dd z + \omega\, , \quad \dd_X \int_0^{x_u} \dd_X\cH_{z} \dd z + \dd_X \omega = \frf\, l_{x_u}\wedge n_{x_u}
\end{eqnarray*}

\noindent
where $\omega$ is a one-form on $X$ and $(l_{x_u} , n_{x_u})$ is a family of closed one-forms, everywhere linearly independent.


\subsection{Everywhere space-like flux}


Let $(g,b,\phi)$ be a supersymmetric NS-NS solutions with everywhere space-like flux and associated skew-torsion isotropic parallelism $[u,v,l,n]$. We choose the representive $(u,v,l,n) \in [u,v,l,n]$ such that:
\begin{eqnarray*}
\varphi_{\phi} = \ast_q \frm_b \, , \qquad \alpha_b = \frm_b
\end{eqnarray*}

\noindent
for a unique $\frm_b \in \Gamma(\langle\mathbb{R}u\rangle \oplus \langle\mathbb{R}v\rangle)^{\perp_g}$. We have:
\begin{eqnarray*}
\dd \frm_b  = \ast_q\frm_b \wedge \frm_b\, , \qquad \dd \ast_q\frm_b = 0
\end{eqnarray*}

\noindent
Therefore, for supersymmetric NS-NS solutions with everywhere space-like flux the flux-dilaton distribution can be identified with:
\begin{eqnarray*}
\Ker(\ast_q\frm_b)\cap \Ker(\frm_b) \subset TM
\end{eqnarray*}

\noindent
and since $\frm_b \in \Gamma(\langle\mathbb{R}u\rangle \oplus \langle\mathbb{R}v\rangle)^{\perp_g}$ is nowhere vanishing we conclude that it is regular and given by the span of $u^{\sharp_g}$ and $v^{\sharp_g}$. Furthermore, we have the following vanishing Lie brackets:
\begin{eqnarray*}
\left[ u^{\sharp_g} , v^{\sharp_g} \right] = 0\, , \quad  \left[ u^{\sharp_g} , l^{\sharp_g}\right]   = 0 \, , \quad \left[ u^{\sharp_g} , n^{\sharp_g}\right]  = 0  
\end{eqnarray*} 

\noindent
From this together with the first equation in \eqref{eq:spaceflux1} it follows that there exists local coordinates $(x_u,x_v,x_1,x_2)$ such that:
\begin{eqnarray*}
 e^{-x_2} u  = \dd x_u \, , \qquad e^{-x_2} \frm_b  = \dd x_1 \, , \qquad  \ast_q \frm_b = \dd x_2
\end{eqnarray*}

\noindent
On the other hand, by construction we have:
\begin{equation*}
\left( \begin{array}{ccc} \alpha_b  \\ 
\varphi_{\phi} \end{array} \right)  = \left( \begin{array}{ccc} \frm_b  \\ 
\ast_q \frm_b \end{array} \right)  = 	\left( \begin{array}{ccc} \alpha_l & \alpha_n \\ 
- \alpha_n & \alpha_l \end{array} \right) \left( \begin{array}{ccc} l  \\ 
n  \end{array} \right)
\end{equation*}	

\noindent
and thus:
\begin{equation*}
\left( \begin{array}{ccc} l  \\ 
	n  \end{array} \right) =\frac{1}{\alpha_l^2 + \alpha_n^2}	\left( \begin{array}{ccc} \alpha_l & - \alpha_n \\ 
	\alpha_n & \alpha_l \end{array} \right)   \left( \begin{array}{ccc} \frm_b  \\ 
		\ast_q \frm_b \end{array} \right)     
\end{equation*}	

\noindent
From this we obtain:
\begin{eqnarray*}
l = \frac{1}{\alpha_l^2 + \alpha_n^2} (\alpha_l e^{x_2} \dd x_1 - \alpha_n \dd x_2)\, , \qquad n = \frac{1}{\alpha_l^2 + \alpha_n^2} (\alpha_n e^{x_2} \dd x_1 + \alpha_l \dd x_2)
\end{eqnarray*}

\noindent
Furthermore, the coordinate $x_v$ can be rearranged so as to $u^{\sharp_g} = \partial_{x_v}$, implying that all the local coefficients of the metric $g = u\odot v + l\otimes l + n \otimes n$ in the coordinates $(x_u,x_v,x_1,x_2)$ are independent of $x_v$. A quick computation gives:
\begin{eqnarray*}
g = e^{x_2}\dd x_u \odot (\cH_{x_u} \dd x_u +   \dd x_v + \omega_{x_u}) + \frac{1}{\alpha_l^2 + \alpha_n^2} (e^{2 x_2} \dd x_1\otimes \dd x_1 + \dd x_2 \otimes \dd x_2)
\end{eqnarray*}

\noindent
where we have written $v = \cH_{x_u}  \dd x_u +  \dd x_v + \omega_{x_u}$ in terms a local family of functions $\cH_{x_u}$ parametrized by $x_u$ and depending on $(x_1,x_2)$ and a local family of  one-forms $\omega_{x_u}$ along the coordinates $(x_1 , x_2)$ and parametrized by $x_u$. Imposing Equation \eqref{eq:spaceflux1} on the representation $(u,v,l,n)$ that we have constructed we obtain:
\begin{eqnarray}
& e^{-x_2}\dd (e^{x_2} v) = -\kappa\wedge (\bar{\alpha}_l e^{x_2} \dd x_1 - \bar{\alpha}_n \dd x_2) -\rho\wedge  (\bar{\alpha}_n e^{x_2} \dd x_1 + \bar{\alpha}_l \dd x_2)   \label{eq:dilatonfluxaux1} \\
& \dd l = - e^{x_2} \dd x_u \wedge (\alpha_n   \omega_{x_u}  + \kappa^{\perp})  \, , \qquad \dd n = e^{x_2} \dd x_u\wedge (\alpha_l   \omega_{x_u}  - \rho^{\perp})  \label{eq:dilatonfluxaux2}
\end{eqnarray}

\noindent
where for simplicity we have set:
\begin{eqnarray*}
\bar{\alpha}_l = \frac{\alpha_l}{\alpha_l^2 + \alpha_n^2}\, , \qquad \bar{\alpha}_n = \frac{\alpha_n}{\alpha_l^2 + \alpha_n^2}
\end{eqnarray*}
Equations \eqref{eq:dilatonfluxaux2} are equivalent to:
\begin{eqnarray*}
& \kappa^{\perp} = - e^{- x_2}\partial_{x_u} l - \alpha_n \omega_{x_u}\, , \qquad \rho^{\perp} = - e^{- x_2} \partial_{x_u} n + \alpha_l \omega_{x_u}\\
& \partial_{x_1} \bar{\alpha}_n +\partial_{x_2} (e^{x_2} \bar{\alpha}_l) = 0\, , \qquad \partial_{x_1} \bar{\alpha}_l - \partial_{x_2} (e^{x_2} \bar{\alpha}_n) = 0
\end{eqnarray*}

\noindent
In particular, the Riemannian metric $q_{x_u} = l\otimes l + n\otimes n$ is a family of flat two-dimensional metrics parametrized by $x_u$. Plugging these relations into Equation \eqref{eq:dilatonfluxaux1}, we obtain:
\begin{eqnarray*}
& \dd x_u \wedge (\partial_{x_u}\omega_{x_u} - e^{-x_2} \dd_{X}(e^{x_2} \cH_{x_u})) + \dd x_2 \wedge \dd x_v + e^{-x_2} \dd_{X}(e^{x_2}\omega_{x_u}) \\
& = - \dd x_v \wedge (\kappa_v l + \rho_v n) + \dd x_u \wedge (\alpha_n l - \alpha_l n)  + (e^{-x_2}\partial_{x_u}l + \alpha_n \omega_{x_u})\wedge l + (e^{-x_2}\partial_{x_u} n - \alpha_l \omega_{x_u}) \wedge  n
\end{eqnarray*}

\noindent
which is equivalent to:
\begin{eqnarray*}
& \partial_{x_u}\omega_{x_u} - e^{-x_2} \dd_{X}(e^{x_2} \cH_{x_u}) = \alpha_n l - \alpha_l n\, , \qquad \dd x_2 = \kappa_v l + \rho_v n\\
& e^{-x_2}  \dd_{X}(e^{x_2}\omega_{x_u}) = (e^{-x_2}\partial_{x_u}l + \alpha_n \omega_{x_u})\wedge l + (e^{-x_2}\partial_{x_u} n - \alpha_l \omega_{x_u}) \wedge  n
\end{eqnarray*}

\noindent
Hence:
\begin{eqnarray*}
&  \partial_{x_u}\omega_{x_u} - e^{-x_2} \dd_{X}(e^{x_2} \cH_{x_u}) + \dd x_2 = 0\, , \qquad \kappa_v \alpha_l + \rho_v \alpha_n = 0\, , \qquad \rho_v \bar{\alpha}_l - \kappa_v \bar{\alpha}_n = 1\\
& \dd_{X}\omega_{x_u} = 2\, (\bar{\alpha}_l \partial_{x_u}\bar{\alpha}_n  - \bar{\alpha}_n \partial_{x_u}\bar{\alpha}_l) \dd x_1 \wedge \dd x_2  
\end{eqnarray*}

\noindent
These equations are solved by:
\begin{eqnarray*}
\omega_{x_u} = e^{-x_2} \dd_{X} (e^{x_2} \int_0^{x_u} (\cH_{z} - 1)\dd z) + \omega\, , \quad \kappa_v = - \alpha_n\, , \quad \rho_v = \alpha_l\\
\int_0^{x_u} \dd_{X} \cH_{z} \dd z \wedge \dd x_2 + \dd_X \omega = 2\, (\bar{\alpha}_l \partial_{x_u}\bar{\alpha}_n  - \bar{\alpha}_n \partial_{x_u}\bar{\alpha}_l) \dd x_1 \wedge \dd x_2  
\end{eqnarray*}

\noindent
where $\omega$ is a one-form along $(x_1,x_2)$.


\renewcommand{\chaptername}{Chapter}

\renewcommand{\leftmark}{Chapter \thechapter. Globally hyperbolic supersymmetric configurations}

\chapter{The NS-NS supergravity evolution flow}
\label{chapter:Globallyhyperbolicsusy}
 

In this chapter we investigate the supersymmetric solutions of four-dimensional NS-NS supergravity whose associated Lorentzian metric is globally hyperbolic, to which we will refer simply as \emph{globally hyperbolic supersymmetric NS-NS solutions}. When restricted to configurations with globally hyperbolic metric, the NS-NS system defines a second-order Riemannian evolution flow on any given Cauchy hypersurface, called the \emph{NS-NS evolution flow}, whereas the supersymmetry conditions of NS-NS supergravity define a first-order evolution flow, called the \emph{supersymmetric NS-NS evolution flow}. The main purpose of this chapter is to initiate the study of the interaction between these two flows, which can be expected to be closely related, and set up a geometric framework to investigate in the future the moduli of initial data sets of both of these flows.   
 

\section{Globally hyperbolic isotropic parallelisms}
\label{sec:ReductionGloballyHyperbolic}


Let $M$ be an oriented four-manifold and let $\mathfrak{F}(M)$ denote the category of isotropic parallelisms on $M$. As explained in Chapter \ref{chapter:IrreducibleSpinors4d}, there is a natural functor:
\begin{equation*}
\mathfrak{F}(M) \to \mathrm{Lor}(M)\, , \qquad [u,v,l,n] \mapsto g = u\odot v + l\otimes l + n\otimes n
\end{equation*}

\noindent
where $\mathrm{Lor}(M)$ is considered as a the category of Lorentzian metrics and isometries on $M$.  

\begin{definition}
\label{def:globallyhyperbolic}
An isotropic parallelism $[u,v,l,n]\in \mathfrak{F}(M)$ is \emph{globally hyperbolic} if its associated metric:
\begin{eqnarray*}
g = u\odot v + l\otimes l + n\otimes n
\end{eqnarray*}

\noindent
is a globally hyperbolic metric on $\cI\times \Sigma$ with Cauchy hypersurface $\Sigma$.
\end{definition}

\noindent
We denote by $\mathfrak{F}_o(M)$ the full subcategory of globally hyperbolic isotropic parallelisms $\mathfrak{F}(M)$. Let $[u,v,l,n]$ be a globally hyperbolic isotropic parallelism on $M$ with associated globally hyperbolic metric $g$. Then, a celebrated theorem of Bernal and S\'anchez \cite{Bernal:2003jb,Bernal:2004gm} states that in this case $(M,g)$ has the following isometry type:
\begin{equation}
\label{eq:globahyp}
(M,g) = (\cI\times \Sigma, -\lambda^2_t \dd t\otimes \dd t + h_t) 
\end{equation}

\noindent
where $t$ is the canonical coordinate of the interval $\cI \subset \mathbb{R}$ containing the origin, $\lambda_t$ is a smooth family of nowhere vanishing functions on $\Sigma$ parametrized by $t\in\cI$ and $h_t$ is a family of complete Riemannian metrics on $\Sigma$, again parametrized by $t\in\cI$. From now on we consider the identification \eqref{eq:globahyp} to be fixed for each $[u,v,l,n] \in \mathfrak{F}_o(M)$. We set:
\begin{equation*}
\Sigma_t := \left\{ t\right\}\times \Sigma \hookrightarrow M\, , \qquad \Sigma := \left\{ 0\right\}\times \Sigma \hookrightarrow M 
\end{equation*}

\noindent
and define:
\begin{equation*}
\mathfrak{t}_t = \lambda_t\, \dd t 
\end{equation*}

\noindent
to be the outward-pointing unit time-like one-form orthogonal to $T^{\ast}\Sigma_t$ for every $t\in \cI$. We will consider $\Sigma\hookrightarrow M$, endowed with the induced Riemannian metric:
\begin{equation*}
h := h_0\vert_{T\Sigma\times T\Sigma} 
\end{equation*}

\noindent
to be the Cauchy hypersurface of $(M,g)$. The \emph{shape operator} or scalar second fundamental form $\Theta_t \in \Gamma (T^{\ast}\Sigma_t \odot T^{\ast}\Sigma_t)$ of the embedded manifold $\Sigma_t\hookrightarrow M$ is defined in the usual way as follows:
\begin{equation}
\label{eq:secondfundamentalform}
\Theta_t  := \nabla^g \mathfrak{t}_t\vert_{T\Sigma_t\times T\Sigma_t} 
\end{equation}

\noindent
Using standard theory of immersed Riemannian manifolds, this definition can be shown to be equivalent to:
\begin{equation*}
\Theta_t = - \frac{1}{2\lambda_t} \partial_t h_t \in \Gamma(T^{\ast}\Sigma_t\odot T^{\ast}\Sigma_t)
\end{equation*}

\noindent
Moreover, it can be seen that:
\begin{equation*}
\nabla^g \alpha \vert_{T\Sigma_t\times TM} = \nabla^{h_t} \alpha + \Theta_t(\alpha)\otimes \mathfrak{t}_t\, , \qquad \forall\,\,\alpha\in \Omega^1(\Sigma_t) 
\end{equation*}

\noindent
where $\nabla^{h_t}$ denotes the Levi-Civita connection on $(\Sigma_t,h_t)$ and $\Theta_t(\alpha) := \Theta_t(\alpha^{\sharp_{h_t}})$ is by definition the evaluation of $\Theta_t$ on the metric dual of $\alpha$. For every representative  $(u,v,l,n)\in [u,v,l,n]\in \mathfrak{F}_o(M)$ we write:
\begin{equation}
\label{eq:hyperbolicsplittingframe}
u = u^o_t\, \mathfrak{t}_t + u^{\perp}_t\, , \qquad v = v^o_t\, \mathfrak{t}_t + v^{\perp}_t  \, , \qquad l = l^o_t\, \mathfrak{t}_t + l^{\perp}_t \, , \qquad n = n^o_t\, \mathfrak{t}_t + n^{\perp}_t 
\end{equation}

\noindent
where the superscript $\perp$ denotes orthogonal projection to $T^{\ast}\Sigma_t$ with respect to $h_t$ and where we have defined:
\begin{equation*}
u^o_t = - g(u,\mathfrak{t}_t)\, , \qquad v^o_t = - g(v,\mathfrak{t}_t)\, , \qquad l^o_t = - g(l,\mathfrak{t}_t)\, ,\qquad n^o_t = - g(n,\mathfrak{t}_t) 
\end{equation*}

\noindent
This splitting is not preserved within the equivalence class $[u,v,l,n]$. Instead, we have:
\begin{eqnarray*}
& \bar{v}^o_t = v^0_t -\frac{1}{2} \vert \frw\vert^2_g u^o_t + \frw^o_t\, , \quad \bar{v}^{\perp}_t = v^{\perp}_t - \frac{1}{2} \vert \frw\vert^2_g u^{\perp}_t + \frw^{\perp}_t\\
& \bar{l}_t^{o} = l^o_t - \frw(l) u^o_t\, , \quad  \bar{l}_t^{\perp} = l^{\perp}_t - \frw (l) u^{\perp}_t\, , \quad \bar{n}_t^{o} = n^o_t - \frw(n) u^o_t\, , \quad  \bar{n}_t^{\perp} = n^{\perp}_t - \frw (n) u^{\perp}_t
\end{eqnarray*}

\noindent
where we have set:
\begin{equation*}
(u,\bar{v},\bar{l},\bar{n}) = \mathfrak{w}\cdot (u,v,l,n)  = (u,v-\frac{1}{2} \vert \frw\vert^2_g u + \frw,l - \frw (l) u,n - \frw (n) u)\, , \quad \frw\in \Gamma (u,v)^{\perp_g}
\end{equation*}
 
\noindent
and where we have split $(u,\bar{v},\bar{l},\bar{n})$ analogously to the splitting of $(u,v,l,n)$ in Equation \eqref{eq:hyperbolicsplittingframe}. Hence, every isotropic parallelism $[u,v,l,n]$ defines an equivalence class of functions and frames on $\Sigma$ through its globally hyperbolic reduction. Rather than considering this equivalence class in the following, we are interested in obtaining the privileged representative that is most convenient for the study of globally hyperbolic supersymmetric solutions and more generally differential spinors on globally hyperbolic Lorentzian four-manifolds. 

\begin{lemma}
\label{lemma:ghcoframe}
Let $[u,v,l,n]$ be a globally hyperbolic isotropic parallelism on $M$. Then, there exists a unique representative $(u,v,l,n)\in[u,v,l,n]$ such that:
\begin{equation}
\label{eq:adaptedcoframe}
u = u^o_t\, \mathfrak{t}_t + u^{\perp}_t\, , \qquad v = \frac{1}{2(u^o_t)^2}\,( -u^o_t \mathfrak{t}_t + u^{\perp}_t)  \, , \qquad l =   l^{\perp}_t \, , \qquad n = n^{\perp}_t 
\end{equation}

\noindent
in the splitting \eqref{eq:hyperbolicsplittingframe}.
\end{lemma}

\begin{proof}
We choose $\frw\in\Gamma (u,v)^{\perp_g}$ such that:
\begin{equation*}
0 = \bar{l}_t^{o} = l^o_t - \frw(l) u^o_t\, , \qquad 0 = \bar{n}_t^{o} = n^o_t - \frw(n) u^o_t
\end{equation*}

\noindent
This already determines $\frw$ uniquely  as follows:
\begin{equation*}
\frw = \frac{1}{u^o_t} (l^o_t l + n^o_t n)
\end{equation*}

\noindent
Hence $(u^o_t)^2 \vert \frw \vert_g^2 =  (l^o_t)^2   + (n^o_t)^2$. Plugging these equations into the previous expression for $\bar{v}$ in terms of $\frw$ and $u$, we obtain:
\begin{eqnarray}
\label{eq:splittingchange}
\bar{v}^o_t = v^o_t + \frac{1}{2  u^o_t} ((l^o_t)^2   + (n^o_t)^2)  \, , \qquad \bar{v}^{\perp}_t = v^{\perp}_t - \frac{1}{2 (u^o_t)^2}((l^o_t)^2   + (n^o_t)^2) u^{\perp}_t + \frac{1}{u^o_t} (l^o_t l^{\perp}_t + n^o_t n^{\perp}_t)
\end{eqnarray}

\noindent
On the other hand, by assumption we must have:
\begin{eqnarray*}
& g = u\odot \bar{v}  + \bar{l}\otimes \bar{l} + \bar{n}\otimes \bar{n}  = (2 u^o_t v^o_t + (l^o_t)^2 + (n^o_t)^2) \frt_t\otimes \frt_t + \frt_t \odot (v^o_t u^{\perp}_t + u^o_t v^{\perp}_t + l^o_t l^{\perp}_t + n^o_t n^{\perp}_t)\\
& + u^{\perp}_t \odot v^{\perp}_t + l^{\perp}_t \otimes l^{\perp}_t + n^{\perp}_t \otimes n^{\perp}_t = \frt_t\otimes \frt_t + h_t
\end{eqnarray*}

\noindent
and thus:
\begin{equation*}
2 u^o_t v^o_t + (l^o_t)^2 + (n^o_t)^2 = -1\, , \qquad v^o_t u^{\perp}_t + u^o_t v^{\perp}_t + l^o_t l^{\perp}_t + n^o_t n^{\perp}_t = 0
\end{equation*}

\noindent
Plugging this equation back into \eqref{eq:splittingchange}, we obtain:
\begin{equation*}
\bar{v}^o_t = -  \frac{1}{2  u^o_t} \, , \qquad \bar{v}^{\perp}_t =  \frac{1}{2 (u^o_t)^2}  u^{\perp}_t  
\end{equation*}

\noindent
and thus we conclude. 
\end{proof}

\noindent
We will refer to the representative $(u,v,l,n)\in [u,v,l,n]$ of the previous lemma as the \emph{adapted} representative in $[u,v,l,n]$. Hence, every globally hyperbolic isotropic parallelism $[u,v,l,n]$ naturally defines a family of tuples:
\begin{equation*}
(\lambda_t , u^o_t , u^{\perp}_t , l^{\perp}_t , n^{\perp}_t ) 
\end{equation*}

\noindent
as explained in the previous discussion. Note that metric coefficient $\lambda_t$, or alternatively the unit time-like one-form $\frt_t$, needs to be considered as part of the globally hyperbolic reduction of $[u,v,l,n]$ in order to be able to recover the globally hyperbolic metric associated to $[u,v,l,n]$ from its globally hyperbolic reduction. By direct substitution, we obtain:
\begin{eqnarray*}
g = u\odot v + l \otimes l + n\otimes n = -\frt_t\otimes \frt_t + \frac{1}{(u^o_t)^2} u^{\perp}_t \otimes  u^{\perp}_t +  l^{\perp}_t\otimes  l^{\perp}_t +  n^{\perp}_t\otimes  n^{\perp}_t
\end{eqnarray*}

\noindent
Note that $u^{\perp}_{t}$ is not necessarily of unit norm, whereas $(u^o_t)^{-1} u^{\perp}_{t}$ is because of $u$ being isotropic. Since $l^{\perp}_t$ and $n^{\perp}_t$ are already unit-norm one-forms on $(\Sigma,h_t)$, it is convenient to introduce the following notation:
\begin{equation*}
e^t_u := \frac{1}{u^o_t} u^{\perp}_{t} \, , \qquad e_l := l^{\perp}_{t} \, , \qquad e_n := n^{\perp}_{t} 
\end{equation*}

\noindent
in terms of which the metric $g$ associated to $[u,v,l,n]$ adopts the standard form:
\begin{equation*}
g =  -\frt_t\otimes \frt_t + e^t_u \otimes  e^t_u + e^t_l \otimes  e^t_l + e^t_n \otimes  e^t_n
\end{equation*}

\noindent
In particular, $(e^t_u,e^t_l,e^t_n)$ is a global orthonormal coframe on $\Sigma$. For ease of notation, in the following we will set $\fre_t = (e^t_u,e^t_l,e^t_n)$. Furthermore, it will be convenient for future computations to introduce a family of functions $\fra_t$ on $\Sigma$ as follows:
\begin{equation*}
u_t^o = e^{\fra_t}\, , \qquad \forall\,\, t\in \cI
\end{equation*}

\noindent
Note that $u_t^o$ is nowhere vanishing for every $t\in\cI$ and furthermore without lost of generality we can take it to be strictly positive.
 
\begin{definition}
The tuple $(\lambda_t,\fra_t,\fre_t)$ is the \emph{globally hyperbolic reduction} of $[u,v,l,n]$ on the three-manifold $\Sigma$.
\end{definition}

\noindent
This globally hyperbolic reduction of a globally hyperbolic isotropic parallelisms $[u,v,l,n]$ will play a fundamental role in our study of globally hyperbolic solutions of NS-NS supergravity.


\section{The evolution problem of skew-torsion parallel spinors}
\label{sec:evolutionskewtorsion}


In this section we consider the evolution problem defined by a globally hyperbolic skew-torsion parallel spinor, namely a skew-torsion parallel spinor whose associated isotropic parallelism is globally hyperbolic in the sense of Definition \ref{def:globallyhyperbolic}. This evolution flow yields the general framework to study the evolution problem defined by globally hyperbolic supersymmetric NS-NS solutions, which is our main object of study and involves, as explained in Chapter \ref{chapter:susyKundt4d}, a particular class of skew-torsion parallel spinors as their main ingredient. By the results of Section \ref{sec:sktorsionparallelspinors}, a globally hyperbolic Lorentzian four-manifold admits a skew-torsion parallel spinor if and only if its associated isotropic parallelism is skew-torsion. Hence, we can study globally hyperbolic Lorentzian four-manifolds equipped with skew-torsion parallel spinors equivalently in terms of globally hyperbolic skew-torsion isotropic parallelisms. We remind the reader that, by virtue of Theorem \ref{thm:existencenullcoframeII}, a four-manifold $M$ admits a skew-torsion isotropic parallelism with torsion $H = \ast_g \alpha \in \Omega^3(M)$ if and only if it admits an isotropic parallelism satisfying the following differential system:
\begin{eqnarray}
& \dd u = \ast_g(\alpha\wedge u)  \, , \qquad \dd v = \ast_g(\alpha\wedge v)  - \kappa\wedge l - \rho\wedge n   \label{eq:duvtorsiongeneralII}\\
& \dd l = \ast_g(\alpha\wedge l) + \kappa\wedge u \, , \qquad \dd n = \ast_g(\alpha\wedge n) + \rho\wedge u   \label{eq:dlntorsiongeneralII}
\end{eqnarray}

\noindent
where $\alpha = \ast_g H$. Our goal now is to \emph{reduce} the differential system \eqref{eq:duvtorsiongeneralII} and \eqref{eq:dlntorsiongeneralII} when evaluated on the adapted representative $(u,v,l,n)\in [u,v,l,n]$ of a globally hyperbolic isotropic parallelism $[u,v,l,n]$. Recall that by Lemma \ref{lemma:ghcoframe} we can write:
\begin{equation*}
u = e^{\fra_t} (\mathfrak{t}_t + e_u^t)\, , \qquad v = \frac{e^{-\fra_t}}{2}\,( -  \mathfrak{t}_t + e_u^t)  \, , \qquad l =   e_l^t \, , \qquad n = e_n^t 
\end{equation*}

\noindent
Using these equations, we compute now the exterior derivative of adapted representative $(u,v,l,n) \in [u,v,l,n]$ of a skew-torsion globally hyperbolic isotropic parallelism $[u,v,l,n]$ in terms of its globally hyperbolic reduction, namely in terms of objects and differentials defined on the three-manifold $\Sigma$.

\begin{lemma}
\label{lemma:derivativeadaptedcoframe}
Let $(u,v,l,n)\in [u,v,l,n]$ be the adapted representative of a globally hyperbolic parallelism. The following formulas hold:
\begin{eqnarray*}
& \dd u = e^{\fra_t} ( \frt_t \wedge (\lambda_t^{-1}\partial_t \fra_t e^t_u - \dd_{\Sigma} \fra_t - \lambda_t^{-1} \dd_{\Sigma} \lambda_t + \lambda_t^{-1} \partial_t e^t_u) + \dd_{\Sigma} e^t_u + \dd_{\Sigma} \fra_t \wedge e^t_u)\\
& \dd v = \frac{1}{2} e^{-\fra_t} ( \frt_t \wedge (\lambda_t^{-1} \partial_t e^t_u - \lambda_t^{-1}\partial_t \fra_t e^t_u - \dd_{\Sigma} \fra_t + \lambda_t^{-1} \dd_{\Sigma} \lambda_t) + \dd_{\Sigma} e^t_u - \dd_{\Sigma} \fra_t \wedge e^t_u)\\
& \dd e^t_l = \frt_t \wedge \lambda_t^{-1} \partial_t e^t_l + \dd_{\Sigma} e^t_l \, , \qquad \dd e^t_n = \frt_t \wedge \lambda_t^{-1} \partial_t e^t_n + \dd_{\Sigma} e^t_n 
\end{eqnarray*}
	
\noindent
where $(\lambda_t,\fra_t,\fre_t)$ is the globally hyperbolic reduction of $[u,v,l,n]$.
\end{lemma}

\begin{proof}
We compute:
\begin{eqnarray*}
& \dd u = e^{\fra_t} \dd \fra_t \wedge (\mathfrak{t}_t + e_u^t) + e^{\fra_t} (\dd\mathfrak{t}_t + \dd e_u^t) = e^{\fra_t} (\dd_{\Sigma} \fra_t \wedge \mathfrak{t}_t + \partial_t \fra_t \dd t\wedge e^t_u + \dd_{\Sigma} a_t \wedge e^t_u )\\
& + e^{\fra_t} (\dd_{\Sigma} \lambda_t \wedge \dd t + \dd_{\Sigma} e_u^t + \dd t \wedge \partial_t e^t_u)
\end{eqnarray*}

\noindent
After isolating for type of tensor we obtain the first equation in the lemma. Similarly, we have:
\begin{eqnarray*}
& \dd v = \frac{1}{2}e^{-\fra_t}  \dd \fra_t \wedge (\mathfrak{t}_t - e_u^t) +  \frac{e^{-\fra_t}}{2} (- \dd\mathfrak{t}_t + \dd e_u^t) =  \frac{e^{-\fra_t}}{2} ( \dd_{\Sigma} \fra_t \wedge \mathfrak{t}_t - \dd_{\Sigma} \fra_t\wedge e_u^t - \partial_t \fra_t \dd t \wedge e^t_u)\\
& +  \frac{e^{-\fra_t}}{2} (\dd_{\Sigma} e_u^t + \dd t \wedge \partial_t e^t_u - \lambda_t^{-1}\dd_{\Sigma}\lambda_t \wedge \mathfrak{t}_t)
\end{eqnarray*}

\noindent
which gives the second equation in the lemma. The third and fourth equations are proven similarly and thus we conclude. 
\end{proof}

\noindent
Let $(\lambda_t,a_t,\fre_t)$ be the globally hyperbolic reduction of a skew-torsion globally hyperbolic parallelism $[u,v,l,n]$ with torsion $H = \ast_g \alpha$. We split $\alpha$ as expected:
\begin{equation}
\label{eq:splittingxialpha}
\alpha = \alpha_t^o \frt_t +  \alpha^{\perp}_t
\end{equation}

\noindent
in terms of uniquely determined families of functions $\alpha_t^o$ and one-forms $\alpha^{\perp}_t$ on $\Sigma$. For further reference in upcoming calculations, we recall the following identities:
\begin{equation}
\label{eq:hodgehyperbolic}
\ast_g(\frt_t \wedge \beta) = -\ast_{h_t}\beta\, , \qquad \ast_{g}\beta = (-1)^{\vert\beta \vert} \frt_t \wedge \ast_{h_t} \beta\, , \qquad \beta \in \Gamma(\wedge M) 
\end{equation}

\noindent
relating the Hodge dual $\ast_g$ associated with $g$ to the Hodge dual $\ast_{h_t}$ associated with $h_t$. Using these relations the following lemma follows by a direct computation.

\begin{lemma}
\label{lemma:dualphasplitting}
Let $(u,v,l,n)$ the globally hyperbolic reduction of  a globally hyperbolic isotropic parallelism $[u,v,l,n]$ on $\Sigma$. Then:
\begin{eqnarray*}
& \ast_g(\alpha\wedge u) = e^{\fra_t}  (\frt_t \wedge \ast_{h_t}(\alpha^{\perp}_t \wedge e^t_u)  - \alpha^o_t \ast_{h_t} e^t_u + \ast_{h_t}\alpha^{\perp}_t) \\
& \ast_g(\alpha\wedge v) = \frac{1}{2} e^{-\fra_t} (\frt_t \wedge \ast_{h_t}(\alpha^{\perp}_t \wedge e^t_u)  - \alpha^o_t \ast_{h_t} e^t_u - \ast_{h_t}\alpha^{\perp}_t) \\
& \ast_g(\alpha\wedge e^t_l) =  \frt_t \wedge \ast_{h_t}(\alpha^{\perp}_t \wedge e^t_l)  - \alpha^o_t \ast_{h_t} e^t_l \\
& \ast_g(\alpha\wedge e^t_n) =  \frt_t \wedge \ast_{h_t}(\alpha^{\perp}_t \wedge e^t_n)  - \alpha^o_t \ast_{h_t} e^t_n 
\end{eqnarray*}

\noindent
where $\alpha = \alpha^o_t \frt_t + \alpha^{\perp}_t \in \Omega^1(\cI\times \Sigma)$.
\end{lemma}

\noindent
The following lemma is the first step in our characterization of skew torsion globally hyperbolic parallelisms.

\begin{lemma}
\label{lemma:skewghevolution1}
A globally hyperbolic parallelism $[u,v,l,n]$ on $\cI\times \Sigma$ is skew-torsion with torsion $H =  \ast_g\alpha$ if and only if there exists a family of functions $(\kappa^o_t,\rho^o_t)$ and a family of one-forms $(\kappa^{\perp}_t , \rho^{\perp}_t)$ on $\Sigma$ such that the globally hyperbolic reduction $(\lambda_t,\fra_t,\fre_t)$ of $[u,v,l,n]$ satisfies the following differential equations on $\Sigma$:
\begin{eqnarray*}
& \lambda_t^{-1} \partial_t e^t_u = \ast_{h_t} (\alpha^{\perp}_t \wedge e^t_u) - \lambda_t^{-1}\partial_t \fra_t e^t_u + \dd_{\Sigma} \fra_t + \lambda_t^{-1} \dd_{\Sigma} \lambda_t\\
& \lambda_t^{-1} \partial_t e^t_l = \ast_{h_t} (\alpha^{\perp}_t \wedge e^t_l) +  e^{\fra_t} (\kappa^o_t e^t_u - \kappa^{\perp}_t) \, , \quad  \lambda_t^{-1} \partial_t e^t_n = \ast_{h_t} (\alpha^{\perp}_t\wedge e^t_n)  + e^{\fra_t} (\rho^o_t e^t_u - \rho^{\perp}_t) \\
& \dd_{\Sigma}  e^t_u  = e^t_u\wedge \dd_{\Sigma} \fra_t + \ast_{h_t}(\alpha^{\perp}_t  - \alpha^o_t e^t_u) \\
& \dd_{\Sigma} e^t_l = e^{\fra_t} \kappa^{\perp}_t\wedge e^t_u - \alpha^o_t \ast_{h_t} e^t_l\, , \quad  \dd_{\Sigma} e^t_n = e^{\fra_t} \rho^{\perp}_t\wedge e^t_u  - \alpha^o_t \ast_{h_t} e^t_n 
\end{eqnarray*}
	
\noindent
together with the equations:
\begin{eqnarray*}
 e^{\fra_t}(\kappa^o_t e^t_l +   \rho^o_t e^t_n )= \lambda_t^{-1}\partial_t \fra_t e^t_u  - \lambda_t^{-1} \dd_{\Sigma} \lambda_t   \, , \quad    \dd_{\Sigma} a_t \wedge e^t_u =  \ast_{h_t}\alpha^{\perp}_t     +    e^{\fra_t} (\kappa^{\perp}_t \wedge e^t_l + \rho^{\perp}_t \wedge e^t_n  )
\end{eqnarray*}

\noindent
further restricting the evolution of $(\lambda_t,\fra_t,\fre_t)$ and $(\kappa^o_t,\rho^o_t, \kappa^{\perp}_t , \rho^{\perp}_t)$.
\end{lemma}
 
\begin{proof}
Let $[u,v,l,n]$ be a globally hyperbolic parallelism on $\cI\times \Sigma$ relative to $\kappa , \rho \in \Omega^1(\cI\times\Sigma)$ and with torsion $H = \ast_{g} \alpha$. We split $\alpha$ as prescribed in Equation \eqref{eq:splittingxialpha}. We split $\kappa$ and $\rho$ similarly to $\alpha$, namely, we set:
\begin{equation*}
\kappa = \kappa^o_t \frt_t + \kappa^{\perp}_t\, ,\qquad \rho = \rho^o_t \frt_t + \rho^{\perp}_t
\end{equation*}
	
\noindent
where $(\kappa^o_t , \rho^o_t)$ are families of functions on $\Sigma$ and $(\kappa^{\perp}_t , \rho^{\perp}_t)$ are families of one-forms on $\sigma$. By Theorem \ref{thm:existencenullcoframeII} it follows that $[u,v,l,n]$ is a globally hyperbolic skew-torsion parallelism if and only if its adapted representative $(u,v,l,n)\in [u,v,l,n]$ on $\cI\times \Sigma$ satisfies \eqref{eq:duvtorsiongeneralII} and \eqref{eq:dlntorsiongeneralII}. Using Lemmas \ref{lemma:derivativeadaptedcoframe} and \ref{lemma:dualphasplitting} we obtain that the first equation in \eqref{eq:duvtorsiongeneralII} is equivalent to:
\begin{equation*}
\frt_t \wedge (\lambda_t^{-1}\partial_t \fra_t e^t_u - \dd_{\Sigma} \fra_t - \lambda_t^{-1} \dd_{\Sigma} \lambda_t + \lambda_t^{-1} \partial_t e^t_u) + \dd_{\Sigma} e^t_u + \dd_{\Sigma} \fra_t \wedge e^t_u =  \frt_t  \wedge  \ast_{h_t} (\alpha^{\perp}_t \wedge e^t_u)  +  \ast_{h_t}(\alpha^{\perp}_t   - \alpha^o_t e^t_u) 
\end{equation*}
	
\noindent
Isolating by tensor type, this equation becomes equivalent to:
\begin{eqnarray}
\label{eq:partialeuproof}
\lambda_t^{-1} \partial_t e^t_u = \ast_{h_t} (\alpha^{\perp}_t \wedge e^t_u) - \lambda_t^{-1}\partial_t \fra_t e^t_u + \dd_{\Sigma} \fra_t + \lambda_t^{-1} \dd_{\Sigma} \lambda_t\, , \,\,  e^{-\fra_t} \dd_{\Sigma} (e^{\fra_t} e^t_u ) = \ast_{h_t}(\alpha^{\perp}_t  - \alpha^o_t e^t_u) 
\end{eqnarray}
	
\noindent
Proceeding similarly for both equations in \eqref{eq:dlntorsiongeneralII}, we obtain that they are equivalent to the following equations:
\begin{eqnarray*}
& \frt_t \wedge \lambda_t^{-1} \partial_t e^t_l + \dd_{\Sigma} e^t_l =  \frt_t  \wedge (\ast_{h_t} (\alpha^{\perp}_t  \wedge e^t_l) + e^{\fra_t} (\kappa^o_t e^t_u - \kappa^{\perp}_t)) + e^{\fra_t} \kappa^{\perp}_t\wedge e^t_u - \alpha^o_t \ast_{h_t} e^t_l \\
& \frt_t \wedge \lambda_t^{-1} \partial_t e^t_n + \dd_{\Sigma} e^t_n = \frt_t  \wedge (\ast_{h_t} (\alpha^{\perp}_t \wedge e^t_n) +   e^{\fra_t} (\rho^o_t e^t_u - \rho^{\perp}_t)) +  e^{\fra_t} \rho^{\perp}_t\wedge e^t_u  - \alpha^o_t \ast_{h_t} e^t_n  
\end{eqnarray*}
	
\noindent
Isolating by tensor type, we obtain the following equivalent system:
\begin{eqnarray*}
& \lambda_t^{-1} \partial_t e^t_l = \ast_{h_t} (\alpha^{\perp}_t \wedge e^t_l) +  e^{\fra_t} (\kappa^o_t e^t_u - \kappa^{\perp}_t) \, , \qquad  \lambda_t^{-1} \partial_t e^t_n = \ast_{h_t} (\alpha^{\perp}_t\wedge e^t_n)  + e^{\fra_t} (\rho^o_t e^t_u - \rho^{\perp}_t) \\
& \dd_{\Sigma} e^t_l = e^{\fra_t} \kappa^{\perp}_t\wedge e^t_u - \alpha^o_t \ast_{h_t} e^t_l\, , \qquad \dd_{\Sigma} e^t_n = e^{\fra_t} \rho^{\perp}_t\wedge e^t_u  - \alpha^o_t \ast_{h_t} e^t_n   
\end{eqnarray*}
	
\noindent
We consider now the second equation in  \eqref{eq:duvtorsiongeneralII}, which is the last that remains to be considered. Using Lemmas \ref{lemma:derivativeadaptedcoframe} and \ref{lemma:dualphasplitting} it follows that it is equivalent to:
\begin{eqnarray*}
& \frac{1}{2} e^{-\fra_t}  \frt_t \wedge (\lambda_t^{-1} \partial_t e^t_u - \lambda_t^{-1}\partial_t \fra_t e^t_u - \dd_{\Sigma} \fra_t + \lambda_t^{-1} \dd_{\Sigma} \lambda_t) + \frac{1}{2} e^{-\fra_t}  \dd_{\Sigma} e^t_u - \frac{1}{2} e^{-\fra_t}  \dd_{\Sigma} \fra_t \wedge e^t_u \\
& =  \frt_t  \wedge ( \frac{1}{2} e^{-\fra_t} \ast_{h_t} (\alpha^{\perp}_t \wedge e^t_u) - \kappa^o_t e^t_l - \rho^o_t e^t_n) -   \frac{1}{2} e^{-\fra_t} \ast_{h_t}(\alpha^{\perp}_t + \alpha^o_t e^t_u)  - \kappa^{\perp}_t \wedge e^t_l - \rho^{\perp}_t \wedge e^t_n
\end{eqnarray*}
	
\noindent
Isolating by tensor type, we obtain:
\begin{eqnarray}
&  \lambda_t^{-1} \partial_t e^t_u    = \lambda_t^{-1}\partial_t \fra_t e^t_u + \dd_{\Sigma} \fra_t - \lambda_t^{-1} \dd_{\Sigma} \lambda_t +   \ast_{h_t} (\alpha^{\perp}_t \wedge e^t_u) -2 e^{\fra_t}\kappa^o_t e^t_l - 2 e^{\fra_t} \rho^o_t e^t_n  \\
&  \dd_{\Sigma} e^t_u  = \dd_{\Sigma} \fra_t \wedge e^t_u - \ast_{h_t}(\alpha^{\perp}_t +  \alpha^o_t e^t_u)   +  2 e^{\fra_t}( e^t_l \wedge \kappa^{\perp}_t +   e^t_n \wedge \rho^{\perp}_t)   \label{eq:auxiliaryeuII}
\end{eqnarray}
	
\noindent
The first equation prescribes the \emph{time-derivative} of $e^t_u$, and therefore needs to be compared with the first equation in \eqref{eq:partialeuproof}. Similarly, the second equation prescribes the exterior derivative of $e^t_u$ and therefore needs to be compared with the second equation in \eqref{eq:partialeuproof}. This results, together with the previous equations, in the \emph{constraints} appearing in the statement of the lemma. The converse follows by tracing back the previous steps and hence we conclude. 
\end{proof}

\noindent
The previous lemma can be conveniently refined as follows.
\begin{prop}
\label{prop:evolutionsystem}
A globally hyperbolic parallelism $[u,v,l,n]$ on $\cI\times \Sigma$ is skew-torsion with torsion $H =  \ast_g\alpha$ if and only if there exists a family of one-forms $(\kappa^{\perp}_t , \rho^{\perp}_t)$ on $\Sigma$ such that the globally hyperbolic reduction $(\lambda_t,\fra_t,\fre_t)$ of $[u,v,l,n]$ satisfies the following system of evolution equations on $\Sigma$:
\begin{eqnarray}
& \partial_t \fra_t =  \dd_{\Sigma} \lambda_t( e^t_u) \nonumber \\
&  \lambda_t^{-1} \partial_t e^t_u = \ast_{h_t} (\alpha^{\perp}_t \wedge e^t_u) + \dd_{\Sigma} \fra_t + \lambda_t^{-1} \dd_{\Sigma} \lambda_t (e_l^t) e^t_l + \lambda_t^{-1} \dd_{\Sigma} \lambda_t (e_n^t) e^t_n \nonumber \\
& \lambda_t^{-1} \partial_t e^t_l = \ast_{h_t} (\alpha^{\perp}_t \wedge e^t_l) - (\lambda_t^{-1} \dd_{\Sigma} \lambda_t (e^t_l) e^t_u + e^{\fra_t}\kappa^{\perp}_t) \label{eq:evolutionequationsI} \\
& \lambda_t^{-1} \partial_t e^t_n = \ast_{h_t} (\alpha^{\perp}_t\wedge e^t_n)  - (\lambda_t^{-1} \dd_{\Sigma} \lambda_t (e^t_n) e^t_u + e^{\fra_t}\rho^{\perp}_t) \nonumber
\end{eqnarray}

\noindent
together with the following time-dependent constraint equations:
\begin{eqnarray}
&  \dd_{\Sigma} \fra_t \wedge e^t_u =  \ast_{h_t}\alpha^{\perp}_t     +    e^{\fra_t} \kappa^{\perp}_t \wedge e^t_l +   e^{\fra_t} \rho^{\perp}_t \wedge e^t_n   \nonumber \\
& \dd_{\Sigma}  e^t_u  =  - e^{\fra_t} \kappa^{\perp}_t \wedge e^t_l - e^{\fra_t} \rho^{\perp}_t \wedge e^t_n -   \alpha^o_t \ast_{h_t} e^t_u \nonumber \\
& \dd_{\Sigma} e^t_l = e^{\fra_t} \kappa^{\perp}_t\wedge e^t_u - \alpha^o_t \ast_{h_t} e^t_l \label{eq:constraintequationsI}\\
& \dd_{\Sigma} e^t_n = e^{\fra_t} \rho^{\perp}_t\wedge e^t_u  - \alpha^o_t \ast_{h_t} e^t_n \nonumber
\end{eqnarray}

\noindent
where $\alpha = \alpha^o_t \frt_t + \alpha^{\perp}_t$. 
\end{prop}

\begin{remark}
From the differential system \eqref{eq:evolutionequationsI} and \eqref{eq:constraintequationsI} it can be easily recovered that the exterior derivative of $e^t_u$ can be equivalently written as follows:
\begin{equation*}
\dd_{\Sigma}  e^t_u  = e^t_u\wedge \dd_{\Sigma} \fra_t + \ast_{h_t}(\alpha^{\perp}_t  - \alpha^o_t e^t_u)
\end{equation*}

\noindent
as it can be verified via a simple combination with the equation for exterior derivative of $\fra_t$. 
\end{remark}

\begin{proof}
We consider the following equation:
\begin{equation*}
 e^{\fra_t}\kappa^o_t e^t_l + e^{\fra_t} \rho^o_t e^t_n = \lambda_t^{-1}\partial_t \fra_t e^t_u  - \lambda_t^{-1} \dd_{\Sigma} \lambda_t   
\end{equation*}

\noindent
in Lemma \ref{lemma:skewghevolution1}. Evaluating this equation on the family of vector fields $(e^t_u)^{\sharp_{h_t}}$ gives equation $\partial_t \fra_t =  \dd_{\Sigma} \lambda_t( e^t_u)$, whereas evaluating it on $(e^t_l)^{\sharp_{h_t}}$ and $(e^t_n)^{\sharp_{h_t}}$ and solving for $\kappa^o_t$ and $\rho^o_t$ we obtain the remaining equations after an educated manipulation.
\end{proof}

\noindent
As already indicated in the statement of the previous proposition, the differential system characterizing skew-torsion parallel spinors splits, when transformed into the equivalent system for skew-torsion globally hyperbolic parallelisms, into two natural blocks, namely block \eqref{eq:evolutionequationsI} consisting on the evolution equations for $(\fra_t,\fre_t)$ together with block \eqref{eq:constraintequationsI} which contains the natural \emph{time-dependent constraint equations} for the aforementioned evolution equations. This is therefore an evolution problem for $(\fra_t,\fre_t)$ for which we can consider $\lambda_t$ as given data. We elaborate on the previous lemma in the following in order to obtain a geometric interpretation of the globally hyperbolic skew-torsion isotropic parallelisms in terms of the second fundamental form of the underlying Cauchy hypersurface. 

\begin{lemma}
\label{lemma:ThetatI}
Let $(\lambda_t,\fra_t,\fre_t)$ be the globally hyperbolic reduction of a skew-torsion globally hyperbolic parallelism $[u,v,l,n]$ relative to $(\kappa^{\perp}_t , \rho^{\perp}_t)$. Then, the second fundamental form $\Theta_t \in \Gamma(T^{\ast}\Sigma \odot T^{\ast}\Sigma)$ of $(\Sigma,h_t)$ is given by:
\begin{eqnarray}
\label{eq:ThetatI}
-2 \Theta_t = \dd_{\Sigma} \fra_t \odot e^t_u - e^{\fra_t} (\kappa^{\perp}_t \odot e^t_l + \rho^{\perp}_t \odot e^t_n)
\end{eqnarray}

\noindent
for every $t\in \cI$.
\end{lemma}

\begin{proof}
Using equation \eqref{eq:evolutionequationsI} we compute:
\begin{eqnarray*}
& -2\Theta_t  = \frac{1}{\lambda_t} (e^t_u \odot \partial_t e^t_u  + e^t_l\odot \partial_t e^t_l + e^t_n \odot\partial_t e^t_n) \\
& = (\ast_{h_t} (\alpha^{\perp}_t \wedge e^t_u)  + \dd_{\Sigma} \fra_t  + \lambda_t^{-1} \dd_{\Sigma} \lambda_t (e_l^t) e^t_l  + \lambda_t^{-1} \dd_{\Sigma} \lambda_t (e_n^t) e^t_n) \odot e^t_u \\
& (\ast_{h_t} (\alpha^{\perp}_t \wedge e^t_l) - (\lambda_t^{-1} \dd_{\Sigma} \lambda_t (e^t_l) e^t_u + e^{\fra_t}\kappa^{\perp}_t))\odot e^t_l \\
& + (\ast_{h_t} (\alpha^{\perp}_t\wedge e^t_n)  - (\lambda_t^{-1} \dd_{\Sigma} \lambda_t (e^t_n) e^t_u + e^{\fra_t}\rho^{\perp}_t)) \odot e^t_n\\
& = \dd_{\Sigma} \fra_t \odot e^t_u - e^{\fra_t} (\kappa^{\perp}_t \odot e^t_l + \rho^{\perp}_t \odot e^t_n)
\end{eqnarray*}

\noindent
and thus we conclude. 
\end{proof}

\begin{prop}
Let $(\lambda_t,\fra_t,\fre_t)$ be the globally hyperbolic reduction of a skew-torsion globally hyperbolic parallelism $[u,v,l,n]$ relative to $(\kappa^{\perp}_t , \rho^{\perp}_t)$. Then:
\begin{eqnarray}
& \Theta_t (e^t_u) = -\dd_{\Sigma} \fra_t - \frac{1}{2}\ast_{h_t} (\alpha^{\perp}_t \wedge e^t_u)\nonumber\\
& \Theta_t (e^t_l) = e^{\fra_t} \kappa^{\perp}_t - \frac{1}{2} \ast_{h_t} (\alpha^{\perp}_t \wedge e^t_l) \label{eq:ThetaexpressionII}\\
& \Theta_t (e^t_n) = e^{\fra_t} \rho^{\perp}_t - \frac{1}{2} \ast_{h_t} (\alpha^{\perp}_t \wedge e^t_n)\nonumber
\end{eqnarray}

\noindent
for every $t\in \cI$, where $\Theta_t$ is the second fundamental form of $(\Sigma,h_t)$.
\end{prop}

\begin{proof}
Evaluating Equation \eqref{eq:ThetatI} on $(e^t_u)^{\sharp_{h_t}}$ we obtain:
\begin{eqnarray*}
\Theta_t(e^t_u) = \dd_{\Sigma} \fra_t   + \dd_{\Sigma} \fra_t (e^t_u) e^t_u - e^{\fra_t} (\kappa^{\perp}_t (e^t_u) e^t_l + \rho^{\perp}_t(e^t_u) e^t_n)
\end{eqnarray*}

\noindent
Combining this equation with the evaluation of the first equation in \eqref{eq:constraintequationsI} on $(e^t_u)^{\sharp_{h_t}}$, we obtain the first equation in \eqref{eq:ThetaexpressionII}. Evaluating Equation \eqref{eq:ThetatI} on $(e^t_l)^{\sharp_{h_t}}$ we obtain:
\begin{equation*}
\Theta_t (e^t_l)= \dd_{\Sigma} \fra_t(e^t_l)  e^t_u - e^{\fra_t} (\kappa^{\perp}_t (e^t_l)   e^t_l + \rho^{\perp}_t(e^t_l)  e^t_n + \kappa^{\perp}_t)
\end{equation*}

\noindent
Combining this equation with the evaluation of the first equation in \eqref{eq:constraintequationsI} on $(e^t_l)^{\sharp_{h_t}}$, we obtain the second equation in \eqref{eq:ThetaexpressionII}. The third equation in \eqref{eq:ThetaexpressionII} is obtained similarly by evaluating Equation \eqref{eq:ThetatI} on $(e^t_n)^{\sharp_{h_t}}$ and then combining the result with the evaluation of the first equation in \eqref{eq:constraintequationsI} on $(e^t_n)^{\sharp_{h_t}}$.
\end{proof}

\noindent
The previous proposition allows to obtain an equivalent expression for the evolution problem defined by the differential system \eqref{eq:evolutionequationsI} and \eqref{eq:constraintequationsI} which is obtained by isolating $\kappa^{\perp}_t$ and $\rho^{\perp}_t$ in \eqref{eq:ThetaexpressionII} and plugging the result back into \eqref{eq:evolutionequationsI} and \eqref{eq:constraintequationsI}.

\begin{cor}
\label{cor:Thetaevolutionsystem}
A globally hyperbolic parallelism $[u,v,l,n]$ on $\cI\times \Sigma$ is skew-torsion with torsion $H =  \ast_g\alpha$ if and only if its globally hyperbolic reduction $(\lambda_t,\fra_t,\fre_t)$ satisfies the following system of evolution equations on $\Sigma$:
\begin{eqnarray}
& \partial_t \fra_t =  \dd_{\Sigma} \lambda_t( e^t_u) \nonumber \\
&  \lambda_t^{-1} \partial_t e^t_u + \Theta_t(e^t_u) = \frac{1}{2}\ast_{h_t} (\alpha^{\perp}_t \wedge e^t_u) +   \lambda_t^{-1} \dd_{\Sigma} \lambda_t (e_l^t) e^t_l + \lambda_t^{-1} \dd_{\Sigma} \lambda_t (e_n^t) e^t_n \nonumber \\
& \lambda_t^{-1} \partial_t e^t_l + \Theta_t(e^t_l)  = \frac{1}{2} \ast_{h_t} (\alpha^{\perp}_t \wedge e^t_l) -  \lambda_t^{-1} \dd_{\Sigma} \lambda_t (e^t_l) e^t_u   \label{eq:evolutionequationsII} \\
& \lambda_t^{-1} \partial_t e^t_n + \Theta_t(e^t_n) = \frac{1}{2} \ast_{h_t} (\alpha^{\perp}_t\wedge e^t_n)  -  \lambda_t^{-1} \dd_{\Sigma} \lambda_t (e^t_n) e^t_u   \nonumber
\end{eqnarray}

\noindent
together with the following time-dependent constraint equations:
\begin{eqnarray}
&  \dd_{\Sigma} \fra_t \wedge e^t_u =  - \Theta_t (e^t_u) \wedge e^t_u + \frac{1}{2}e^t_u \wedge \ast_{h_t} (\alpha^{\perp}_t \wedge e^t_u)  \nonumber \\
& \dd_{\Sigma} e^t_u = \Theta_t(e^t_u)\wedge e^t_u + \ast_{h_t} (\alpha^{\perp}_t   - \alpha^o_t e^t_u) - \frac{1}{2} e^t_u \wedge \ast_{h_t} (\alpha^{\perp}_t \wedge e^t_u) \nonumber \\
& \dd_{\Sigma} e^t_l = \Theta_t(e^t_l)\wedge e^t_u +  (\frac{1}{2}\alpha^{\perp}_t (e^t_u) - \alpha^o_t) \ast_{h_t} e^t_l \label{eq:constraintequationsII} \\
& \dd_{\Sigma} e^t_n = \Theta_t(e^t_n)\wedge e^t_u +  (\frac{1}{2} \alpha^{\perp}_t (e^t_u) - \alpha^o_t) \ast_{h_t} e^t_n \nonumber
\end{eqnarray}

\noindent
where $\Theta_t$ is the second fundamental form of $(\Sigma,h_t)$ and $\alpha = \alpha^o_t \frt_t + \alpha^{\perp}_t$. 
\end{cor}

\begin{remark}
The first equation in \eqref{eq:constraintequationsII} can be equivalently written as follows:
\begin{equation*}
e^{-\fra_t}\dd_{\Sigma} (e^{\fra_t} e^t_u) =  \ast_{h_t} (\alpha^{\perp}_t   - \alpha^o_t e^t_u)
\end{equation*}

\noindent
which will be useful in the following. 
\end{remark}

\noindent
The previous reformulation of the evolution problem posed by a globally hyperbolic skew-torsion spinor is particularly convenient to define the associated initial data and constraint equations. Still, it can be refined as to eliminate the family of functions $\fra_t$ in terms of a cohomological condition in de Rahm cohomology. 

\begin{prop}
\label{prop:Thetaevolutionsystem}
A globally hyperbolic parallelism $[u,v,l,n]$ on $\cI\times \Sigma$ is skew-torsion with torsion $H =  \ast_g\alpha$ if and only if its globally hyperbolic reduction $(\lambda_t,\fra_t,\fre_t)$ satisfies the following system of evolution equations on $\Sigma$:
\begin{eqnarray}
	&  \lambda_t^{-1} \partial_t e^t_u + \Theta_t(e^t_u) = \frac{1}{2}\ast_{h_t} (\alpha^{\perp}_t \wedge e^t_u) +   \lambda_t^{-1} \dd_{\Sigma} \lambda_t (e_l^t) e^t_l + \lambda_t^{-1} \dd_{\Sigma} \lambda_t (e_n^t) e^t_n \nonumber \\
	& \lambda_t^{-1} \partial_t e^t_l + \Theta_t(e^t_l)  = \frac{1}{2} \ast_{h_t} (\alpha^{\perp}_t \wedge e^t_l) -  \lambda_t^{-1} \dd_{\Sigma} \lambda_t (e^t_l) e^t_u   \label{eq:evolutionequationsIII} \\
	& \lambda_t^{-1} \partial_t e^t_n + \Theta_t(e^t_n) = \frac{1}{2} \ast_{h_t} (\alpha^{\perp}_t\wedge e^t_n)  -  \lambda_t^{-1} \dd_{\Sigma} \lambda_t (e^t_n) e^t_u   \nonumber
\end{eqnarray}

\noindent
together with the following time-dependent constraint equations:
\begin{eqnarray}
& \dd_{\Sigma} e^t_u = \Theta_t(e^t_u)\wedge e^t_u + \ast_{h_t} (\alpha^{\perp}_t   - \alpha^o_t e^t_u) - \frac{1}{2} e^t_u \wedge \ast_{h_t} (\alpha^{\perp}_t \wedge e^t_u) \nonumber \\
& \dd_{\Sigma} e^t_l = \Theta_t(e^t_l)\wedge e^t_u +  (\frac{1}{2} \alpha^{\perp}_t (e^t_u) - \alpha^o_t) \ast_{h_t} e^t_l \label{eq:constraintequationsIII} \\
& \dd_{\Sigma} e^t_n = \Theta_t(e^t_n)\wedge e^t_u +  (\frac{1}{2} \alpha^{\perp}_t (e^t_u) - \alpha^o_t) \ast_{h_t} e^t_n \nonumber 
\end{eqnarray}

\noindent
and the following conditions:
\begin{equation}
\label{eq:cohomologicalccondition}
0 = [\Theta_t (e^t_u) + \frac{1}{2}\ast_{h_t} (\alpha^{\perp}_t \wedge e^t_u)]\in H^1(\Sigma,\mathbb{R})\, , \quad \dd_{\Sigma}(\dd_{\Sigma}\lambda_t (e^t_u)) = \partial_t (\Theta_t (e^t_u) + \frac{1}{2}\ast_{h_t} (\alpha^{\perp}_t \wedge e^t_u))
\end{equation}

\noindent
where $\Theta_t$ is the second fundamental form of $(\Sigma,h_t)$ and $\alpha = \alpha^o_t \frt_t + \alpha^{\perp}_t$. 
\end{prop}

\begin{remark}
The first equation in \eqref{eq:constraintequationsIII} can be equivalently written as follows:
\begin{equation*}
\dd_{\Sigma} e^t_u = \Theta_t(e^t_u)\wedge e^t_u +  (\alpha^{\perp}_t (e^t_u)   - \alpha^o_t  ) \ast_{h_t} e^t_u  +\frac{1}{2} (\alpha^{\perp}_t(e^t_l) \ast_{h_t} e^t_l + \alpha^{\perp}_t(e^t_n) \ast_{h_t} e^t_n) 
\end{equation*}

\noindent
which is sometimes useful in computations. 
\end{remark}

\begin{proof}
Let $(\lambda_t,\fra_t,\fre_t)$ be a solution to the differential system \eqref{eq:evolutionequationsII} and \eqref{eq:constraintequationsII}. Then, $(\lambda_t , \fre_t)$ automatically satisfies equations \eqref{eq:evolutionequationsIII} and \eqref{eq:constraintequationsIII} and thus we only need to deal with the conditions contained in Equation \eqref{eq:cohomologicalccondition}. By the first equation in \eqref{eq:ThetaexpressionII} the first equation in \eqref{eq:cohomologicalccondition} follows. Taking the exterior derivative of the first equation in \eqref{eq:ThetaexpressionII} and combining the result with the first equation in \eqref{eq:evolutionequationsII}  we obtain the second equation in \eqref{eq:cohomologicalccondition}. For the converse, suppose that $(\lambda_t,\fre_t)$ satisfies equations \eqref{eq:evolutionequationsIII}, \eqref{eq:constraintequationsIII} and \eqref{eq:cohomologicalccondition}. We only need to prove that the first equation \eqref{eq:evolutionequationsII} and the first equation in \eqref{eq:constraintequationsII} both hold. By the first equation in \eqref{eq:cohomologicalccondition}, there exist a family of functions $\bar{a}_t$ such that:
\begin{equation*}
\dd_{\Sigma} \bar{a}_t + \Theta_t (e^t_u) + \frac{1}{2}\ast_{h_t} (\alpha^{\perp}_t \wedge e^t_u) = 0
\end{equation*}

\noindent
Wedging this equation with $e^t_u$ yields the first equation in \eqref{eq:constraintequationsII}. On the other hand, taking its time derivative and using the second equation in \eqref{eq:cohomologicalccondition} we obtain:
\begin{equation*}
\dd_{\Sigma} \partial_t \bar{a}_t = \dd_{\Sigma}(\dd_{\Sigma}\lambda_t (e^t_u))
\end{equation*}

\noindent
and thus:
\begin{equation*}
\partial_t \bar{a}_t = \dd_{\Sigma}\lambda_t (e^t_u) + c_t
\end{equation*}

\noindent
where $c_t$ is a family of constants on $\Sigma$. Defining $\fra_t = \bar{a}_t - \int c_t \,\dd t$ we recover the first equations in \eqref{eq:evolutionequationsII} and \eqref{eq:constraintequationsII} and thus we conclude. 
\end{proof}

\noindent
Therefore, we can consider the evolution problem defined by a globally hyperbolic isotropic parallelism as being given either by the differential system  \eqref{eq:evolutionequationsI} and \eqref{eq:constraintequationsI}, the differential system \eqref{eq:evolutionequationsII} and \eqref{eq:constraintequationsII} or the differential system given by \eqref{eq:evolutionequationsIII}, \eqref{eq:constraintequationsIII} and \eqref{eq:cohomologicalccondition}. In the first two cases the data that is being \emph{evolved} are pairs of the form $(\fra_t,\fre_t)$ for fixed data $(\lambda_t,\alpha^o_t,\alpha^{\perp}_t)$, whereas in the latter case are families $(\fre_t)$ again for fixed data $(\lambda_t,\alpha^o_t,\alpha^{\perp}_t)$. The evolution problem as given in \eqref{eq:evolutionequationsI} and \eqref{eq:constraintequationsI} or \eqref{eq:evolutionequationsII} and \eqref{eq:constraintequationsII} is more convenient from the analytic point of view, whereas the equivalent evolution problem given in \eqref{eq:evolutionequationsIII}, \eqref{eq:constraintequationsIII} and \eqref{eq:cohomologicalccondition} is arguably more convenient from a geometric point of view, especially to study the corresponding constraint equations.

\begin{definition}
Given data $(\lambda_t,\alpha^o_t,\alpha^{\perp}_t)$, the differential system \eqref{eq:evolutionequationsI} and \eqref{eq:constraintequationsI}, equivalently, the differential system \eqref{eq:evolutionequationsII}  and \eqref{eq:constraintequationsII}, for pairs $(\fra_t, \fre_t)$ is the \emph{skew-torsion Cauchy flow}.
\end{definition}

\noindent
We will also refer to the evolution problem defined by equations \eqref{eq:evolutionequationsIII}, \eqref{eq:constraintequationsIII} and \eqref{eq:cohomologicalccondition} for $(\fre_t)$ as the skew-torsion Cauchy flow. In Proposition \ref{prop:evolutionsystem} we have expressed the evolution problem associated to  skew-torsion globally hyperbolic parallelisms in terms of the existence of families of one-forms $(\kappa^{\perp}_t,\rho^{\perp}_t)$ for which the evolution and constraint equations respectively given in \eqref{eq:evolutionequationsI} and \eqref{eq:constraintequationsI} are satisfied. This family of one-forms can be considered as \emph{auxiliary}, in the sense that they can be completely eliminated from \eqref{eq:evolutionequationsI} and \eqref{eq:constraintequationsI}, as the following result shows, without involving the second fundamental form $\Theta_t$ as we did in Corollary \ref{cor:Thetaevolutionsystem} and Proposition \ref{prop:Thetaevolutionsystem}. 

\begin{prop}
\label{prop:simplifiedevolution}
A globally hyperbolic parallelism $[u,v,l,n]$ on $\cI\times \Sigma$ is skew-torsion with torsion $H =  \ast_g\alpha$ if and only if its globally hyperbolic reduction $(\lambda_t,\fra_t,\fre_t)$ satisfies the following evolution equations on $\Sigma$:
\begin{eqnarray*}
& \lambda_t^{-1} \partial_t e^t_u = \ast_{h_t} (\bar{\alpha}_t \wedge e^t_u) + \dd_{\Sigma} \fra_t + \lambda_t^{-1}( \dd_{\Sigma} \lambda_t (e_l^t) e^t_l + \dd_{\Sigma} \lambda_t (e_n^t) e^t_n)\, , \quad  \dd_{\Sigma}e^t_u = e^t_u\wedge \dd_{\Sigma} \fra_t +  \ast_{h_t}\bar{\alpha}_t \\
& \dd_{\Sigma} e^t_l = e_u^t\wedge (\lambda_t^{-1} \partial_t e^t_l - \bar{\alpha}_t (e^t_u) e^t_n)\, , \quad \dd_{\Sigma} e^t_n = e_u^t\wedge (\lambda_t^{-1} \partial_t e^t_n + \bar{\alpha}_t (e^t_u)  e^t_l)\\
& \partial_t \fra_t = (\dd_{\Sigma}\lambda_t )(e^t_u) \, , \quad \lambda_t^{-1} (e^t_l\wedge \partial_t e^t_l + e^t_n\wedge \partial_t e^t_n)  = \dd_{\Sigma} \fra_t \wedge e^t_u  +  \alpha^{\perp}_t(e^t_u) \ast_{h_t} e^t_u - \lambda_t^{-1} \dd_{\Sigma}\lambda_t \wedge e_u^t  
\end{eqnarray*}

\noindent
where we have set $\bar{\alpha}_t := \alpha^{\perp}_t  - \alpha^o_t e^t_u$ for every $t\in \cI$.
\end{prop}

\begin{proof}
A globally hyperbolic parallelism $[u,v,l,n]$ on $\cI\times \Sigma$ is skew-torsion with torsion $H =  \ast_g\alpha$ if and only if its globally hyperbolic reduction $(\lambda_t,\fra_t,\fre_t)$ satisfies \eqref{eq:evolutionequationsI} and \eqref{eq:constraintequationsI}. Isolating for $\kappa^{\perp}_t$ and $\rho^{\perp}_t$ in \eqref{eq:constraintequationsI}, we obtain:
\begin{eqnarray*}
& e^{\fra_t}\kappa^{\perp}_t  =  \ast_{h_t} (\alpha^{\perp}_t \wedge e^t_l) -  \lambda_t^{-1} \partial_t e^t_l - \lambda_t^{-1} \dd_{\Sigma}\lambda_t (e^t_l) e^t_u\\
& e^{\fra_t}\rho^{\perp}_t =   \ast_{h_t} (\alpha^{\perp}_t \wedge e^t_n) -    \lambda_t^{-1} \partial_t e^t_n - \lambda_t^{-1} \dd_{\Sigma} \lambda_t (e^t_n) e^t_u 
\end{eqnarray*}

\noindent
which \emph{solve} the third and fourth equations in \eqref{eq:evolutionequationsI}. This implies:
\begin{eqnarray*}
& e^{\fra_t} e^t_l \wedge \kappa^{\perp}_t + e^{\fra_t} e^t_n \wedge \rho^{\perp}_t  = \alpha^{\perp}_t(e^t_u) \ast_{h_t} e^t_u + \ast_{h_t} \alpha^{\perp}_t   -  \lambda_t^{-1} (e^t_l\wedge \partial_t e^t_l + e^t_n\wedge \partial_t e^t_n)  - \lambda_t^{-1} \dd_{\Sigma}\lambda_t \wedge e_u^t
\end{eqnarray*}

\noindent
which plugged into the remaining equations in \eqref{eq:evolutionequationsI} and \eqref{eq:constraintequationsI} together with the explicit expression for $\kappa^{\perp}_t$ and $\rho^{\perp}_t$ above, gives the equations in the statement of the proposition. 
\end{proof}

\noindent
As a corollary of the proof of the previous proposition we obtain the following explicit expressions for the one forms $\kappa , \rho \in \Omega^1(\cI\times\Sigma)$ relative to which a globally hyperbolic parallelism is isotropic.

\begin{cor}
Let $[u,v,l,n]$ be a skew-torsion globally hyperbolic parallelism with hyperbolic reduction $(\lambda_t,\fra_t,\fre_t)$ relative to $\kappa,\rho \in \Omega^1(\cI\times \Sigma)$. Then:
\begin{eqnarray*}
& \kappa^o_t = - e^{-\fra_t} \lambda_t^{-1} \dd \lambda_t(e^t_l)\, , \quad\kappa^{\perp}_t = e^{-\fra_t} ( \ast_{h_t} (\alpha^{\perp}_t \wedge e^t_l) -  \lambda_t^{-1} \partial_t e^t_l - \lambda_t^{-1} \dd_{\Sigma}\lambda_t (e^t_l) e^t_u) \\
& \rho^o_t = - e^{-\fra_t} \lambda_t^{-1} \dd \lambda_t(e^t_n)\, , \quad \rho^{\perp}_t = e^{-\fra_t}(\ast_{h_t} (\alpha^{\perp}_t \wedge e^t_n) -    \lambda_t^{-1} \partial_t e^t_n - \lambda_t^{-1} \dd_{\Sigma} \lambda_t (e^t_n) e^t_u )
\end{eqnarray*}

\noindent
where $\kappa = \kappa^o_t \frt_t + \kappa^{\perp}_t$ and $\rho = \rho^o_t \frt_t + \rho^{\perp}_t$.
\end{cor}

\noindent
On a tubular neighborhood around the Cauchy hypersurface $\Sigma$ we can choose coordinates such that $\lambda_t = 1$, which is a condition that can be assumed without loss of generality to study the skew-torsion Cauchy flow locally \emph{in time}. With this assumption the previous corollary gives:
\begin{equation}
\label{eq:reducedkapparho}
\kappa = \kappa^{\perp}_t = e^{-\fra_t} ( \ast_{h_t} (\alpha^{\perp}_t \wedge e^t_l) -\partial_t e^t_l)\, , \qquad \rho = \rho^{\perp}_t = e^{-\fra_t}(\ast_{h_t} (\alpha^{\perp}_t \wedge e^t_n) -  \partial_t e^t_n )
\end{equation}

\noindent
Similarly, the skew-torsion Cauchy flow reduces to the following evolution equations:
\begin{eqnarray}
\label{eq:evolutionequationsInolambda}
\partial_t e^t_u = \ast_{h_t} (\alpha^{\perp}_t \wedge e^t_u) + \dd_{\Sigma} a \, , \quad \partial_t e^t_l = \ast_{h_t} (\alpha^{\perp}_t \wedge e^t_l) -  e^{a}\kappa^{\perp}_t \, ,\quad  \partial_t e^t_n = \ast_{h_t} (\alpha^{\perp}_t\wedge e^t_n)  - e^{a}\rho^{\perp}_t 
\end{eqnarray}

\noindent
together with the following time-dependent constraints:
\begin{eqnarray}
&  \dd_{\Sigma} \fra_t \wedge e^t_u =  \ast_{h_t}\alpha^{\perp}_t     +    e^{\fra_t} \kappa^{\perp}_t \wedge e^t_l +   e^{\fra_t} \rho^{\perp}_t \wedge e^t_n   \nonumber \\
& \dd_{\Sigma}  e^t_u  =  - e^{\fra_t} \kappa^{\perp}_t \wedge e^t_l - e^{\fra_t} \rho^{\perp}_t \wedge e^t_n -   \alpha^o_t \ast_{h_t} e^t_u \nonumber \\
& \dd_{\Sigma} e^t_l = e^{\fra_t} \kappa^{\perp}_t\wedge e^t_u - \alpha^o_t \ast_{h_t} e^t_l \label{eq:constraintequationsInolambda}\\
& \dd_{\Sigma} e^t_n = e^{\fra_t} \rho^{\perp}_t\wedge e^t_u  - \alpha^o_t \ast_{h_t} e^t_n \nonumber
\end{eqnarray}

\noindent
or, equivalently, to the following evolution equations:
\begin{eqnarray}
\label{eq:evolutionequationsIInolambda}
\partial_t e^t_a + \Theta_t(e^t_a)  = \frac{1}{2} \ast_{h_t} (\alpha^{\perp}_t \wedge e^t_a) \, , \qquad a = u,l,n
\end{eqnarray}

\noindent
together with the following time-dependent constraints:
\begin{eqnarray}
& \dd_{\Sigma} e^t_u = \Theta_t(e^t_u)\wedge e^t_u + \ast_{h_t} (\alpha^{\perp}_t   - \alpha^o_t e^t_u) - e^t_u \wedge \frac{1}{2}\ast_{h_t} (\alpha^{\perp}_t \wedge e^t_u) \nonumber \\
& \dd_{\Sigma} e^t_l = \Theta_t(e^t_l)\wedge e^t_u +  (\frac{1}{2} \alpha^{\perp}_t (e^t_u) - \alpha^o_t) \ast_{h_t} e^t_l \label{eq:constraintequationsIInolambda} \\
& \dd_{\Sigma} e^t_n = \Theta_t(e^t_n)\wedge e^t_u +  (\frac{1}{2} \alpha^{\perp}_t (e^t_u) - \alpha^o_t) \ast_{h_t} e^t_n \nonumber \\
& 0 = [\Theta_t (e^t_u) + \frac{1}{2}\ast_{h_t} (\alpha^{\perp}_t \wedge e^t_u)]\in H^1(\Sigma,\mathbb{R})\, , \quad \ \partial_t (\Theta_t (e^t_u) + \frac{1}{2}\ast_{h_t} (\alpha^{\perp}_t \wedge e^t_u)) = 0\nonumber 
\end{eqnarray}

\noindent
We will refer to this evolution flow as the \emph{normal} skew-torsion Cauchy flow. As expected, the constraint equations are not affected by the choice $\lambda_t = 1$ since they do not depend on the latter. 

\begin{thm}
\label{thm:Cauchytorsion}
An oriented three-manifold $\Sigma$ admits an embedding as a Cauchy hypersurface into a skew-torsion flat globally hyperbolic Lorentzian four-manifold equipped with a skew-torsion flat parallel spinor if and only if $\Sigma$ admits a normal skew-torsion Cauchy flow $\fre_t$ with torsion $(\alpha^o_t, \alpha^{\perp}_t)$ satisfying: 
\begin{eqnarray*}
 \beta_l  = \ast_{h_t} (\alpha^{\perp}_t \wedge e^t_l) -\partial_t e^t_l \, , \quad  \beta_n  = \ast_{h_t} (\alpha^{\perp}_t \wedge e^t_n) -\partial_t e^t_n
\end{eqnarray*}
 
\noindent
for a pair of closed one-forms $\beta_l , \beta_n \in \Omega^1(\Sigma)$.
\end{thm}

\begin{proof}
Follows from Equation \eqref{eq:reducedkapparho} together with Corollary \ref{cor:torsionflat}.
\end{proof}

\noindent
Interestingly enough, the skew-torsion flat condition on the globally hyperbolic metric $g$, which is a second order partial differential equation, is guaranteed in the previous theorem by a first-order \emph{time} condition in terms of representatives of the rank-one cohomological invariants of the underlying skew-torsion parallel spinor.

\begin{cor}
Let $(M,g)$ be a skew-torsion globally hyperbolic Lorentzian four-manifold whose rank-one invariants do not admit any time-independent representative on a normal tubular neighborhood of $\Sigma\subset M$. Then $(M,g)$ is not skew-torsion flat.
\end{cor}

\noindent
We end this section with a brief discussion of the constraint equations of the skew-torsion Cauchy flow. Assume momentarily that $\Sigma$ is compact to avoid any analytic issues in the discussion below. Due to the fact that the skew-torsion Cauchy flow \emph{evolves} pairs of the form $(\fra_t,\fre_t)$, the associated Cauchy problem can be expected to require prescribing both $(\fra_t,\fre_t)$ and its derivative at $t=0$. However, since the time-dependent constraint equations \eqref{eq:constraintequationsII} do not contain any term involving $\partial_t \fra_t$, we do not need to specify its value at $t = 0$. Instead, it is prescribed by $\partial_t e^t_u$ and the given data $\lambda_t$ at $t=0$ through the expression:
\begin{equation*}
	\partial_t a_t \vert_{t=0} = \dd_{\Sigma}\lambda_0(e^0_u)
\end{equation*}

\noindent
which corresponds to the evaluation of the first equation in \eqref{eq:evolutionequationsII} at $t=0$. Therefore, we only need to worry about the time derivative of $\fre_t$ at $t=0$, which we denote by $\fre$. Let $F(\Sigma)$ be the bundle of oriented coframes on $M$, which is a trivializable principal bundle over $M$ with structure group given by the identity component $\Gl_o(3,\mathbb{R})$ of the general linear group $\Gl(3,\mathbb{R})$ in three dimensions. We denote by $\Gamma(F(\Sigma))$ the space of sections of $F(\Sigma)$. Families $(\fre_t)$ define curves $(\fre_t)  \colon \cI\to \Gamma(F(\Sigma))$ which are smooth in the sense that they are defined by smooth sections of the pull-back bundle $\fre\in \Gamma(\pr^{\ast}F(\Sigma))$, where $\pr\colon \cI \times \Sigma  \to \Sigma$ is the canonical projection on $\Sigma$. Every element $\fre \in \Gamma(F(\Sigma))$ determines a canonical trivialization:
\begin{eqnarray*}
	F(\Sigma) =\Sigma\times \Gl_o(3,\mathbb{R})
\end{eqnarray*}

\noindent
which extends to a natural identification at the level of sections:
\begin{equation*}
	\Gamma(F(\Sigma)) = C^{\infty}(\Sigma,\Gl_o(3,\mathbb{R}))
\end{equation*}

\noindent
where $C^{\infty}(\Sigma,\Gl_o(3,\mathbb{R}))$ denotes the space of smooth maps from $\Sigma$ to $\Gl_o(3,\mathbb{R})$. In particular, if $\fre \in \Gamma(F(\Sigma))$ is used to determine the previous identifications, then it is mapped to the constant map from $\Sigma$ to the identity in $\Gl_o(3,\mathbb{R})$. If $(\fre_t)\colon \cI \to \Gamma(F(\Sigma))$ is a smooth curve such that $\fre := \fre_0$, then the tangent space of $\Gamma(F(\Sigma))$ at $\fre$ can be identified as follows:
\begin{equation*}
	T_{\fre}\Gamma(F(\Sigma)) = C^{\infty}(\Sigma,\frg(3,\mathbb{R})) = \Gamma(T^{\ast}\Sigma\otimes \mathbb{R}^3)
\end{equation*}

\noindent
where $\frg(3,\mathbb{R})$ denotes the Lie algebra of $\Gl(3,\mathbb{R})$. Hence, the tangent space of $\Gamma(F(\Sigma))$ at a given point is canonically identified with triples of one-forms on $M$, which we will denote generically by $\frv\in \Gamma(T^{\ast}M\otimes \mathbb{R}^3)$. The components of such a triplet of one-forms will be denoted by $\frv = (\frv_u,\frv_l,\frv_n)$. Alternatively, if we consider $(\fre_t)\colon \cI \to \Gamma(F(\Sigma))$ as a section of $\pr^{\ast}F(\Sigma)$, then $\partial_t \fre_t = \cL_{\partial_t} \fre_t$ and its restriction to $t=0$ clearly define a triplet of one-forms on $M$ given by the time Lie derivatives of each of one-forms that constitute the coframe $\fre$ on $\cI\time \Sigma$. Given an element $(\fre,\frv)\in TF(\Sigma)$, we obtain a naturally symmetric two-form on $\Sigma$ given by:
\begin{equation*}
	\Theta_{\fre\frv} := \frv_u\odot \fre_u + \frv_l\odot \fre_l + \frv_n \odot \fre_n
\end{equation*}

\noindent
In this way we obtain a natural map $\Theta\colon TF(M) \to \Gamma(T^{\ast}\Sigma\odot T^{\ast}\Sigma)$ given by $(\fre,\frv) \mapsto \Theta_{\fre\frv}$. Fixing $\fre\in \Gamma(F(M))$, we obtain by restriction:
\begin{equation*}
	\Theta_{\fre} \colon T_{\fre}\Gamma(F(M)) \to \Gamma(T^{\ast}\Sigma\odot T^{\ast}\Sigma)\, , \qquad \frv \mapsto \Theta_{\fre\frv}
\end{equation*}

\noindent
Hence, and by the previous discussion, evaluating at $t = 0\in \cI$ the time-dependent constraints given in \eqref{eq:constraintequationsII} we obtain the \emph{constraint equations} of the skew-torsion spinor Cauchy flow, which are explicitly given by:
\begin{eqnarray}
	& \dd_{\Sigma} e_u = \Theta_{\fre\frv} (e_u)\wedge e_u + \ast_{h} (\alpha^{\perp}   - \alpha^o e_u) +  \frac{1}{2}\ast_{h} (\alpha^{\perp} \wedge e_u)\wedge e_u \nonumber \\
	& \dd_{\Sigma} e_i = \Theta_{\fre\frv} (e_i)\wedge e_u +  (\frac{1}{2} \alpha^{\perp}_t (e_u) - \alpha^o) \ast_{h} e_i\, , \quad i = l, n \label{eq:constraintsnotime} \\
	& 0 = [\Theta_{\fre\frv} (e_u) + \frac{1}{2}\ast_{h} (\alpha^{\perp} \wedge e_u)]\in H^1(\Sigma,\mathbb{R})\nonumber 
\end{eqnarray}

\noindent
for elements $(\fre,\frv) \in  TF(\Sigma)$. Hence, $TF(\Sigma)$ is the \emph{configuration space} for the initial data of the skew-torsion spinor Cauchy flow.

\section{Left-invariant parallel spinor flows}
\label{sec:leftinvariant}


The general theory of differential spinors and globally hyperbolic istropic parallelisms that we have developed in previous sections is especially well-adapted to study spinors parallel under very general connections. However, it already offers a clear and transparent framework in the simplest case of spinors parallel under the Levi-Civita connection on a Lorentzian four-manifold. The study of such irreducible isotropic parallel spinors is classical in the literature \cite{EhlersKundt}. In particular, for parallel isotropic spinors on globally hyperbolic Lorentzian manifolds, the seminal work of Baum, Leistner and Lischewski \cite{BaumLeistnerLischewski,LeistnerLischewski,Lischewski} proved the well-posedness of the corresponding Cauchy problem via a careful analysis of suitable hyperbolic evolution equations given in terms of the Ricci tensor and other geometric objects. Despite these fundamental results, producing explicit parallel spinor flows, namely the evolution flows determined by isotropic parallel spinors on globally hyperbolic Lorentzian four-manifolds, is a difficult task. The framework of globally hyperbolic isotropic parallelisms that we have developed in the previous subsections allows to attack this problem directly. We will do so in this section in the left-invariant case, classifying as a result all left-invariant parallel spinor flows on a Cauchy hypersurface given by a fixed simply connected three-dimensional Lie group. The simply-connected condition is not arbitrary since Proposition \ref{prop:Thetaevolutionsystem}, which characterizes parallel spinor flows by taking $\alpha^{\perp}_t = \alpha^o_t = 0$, involves a cohomological condition that needs to be preserved when taking discrete quotients.  


\subsection{Left-invariant parallel Cauchy pairs}


Let $\Sigma = \G$ be a simply connected three-dimensional Lie group. Denote by $\Conf(\G)$ the set of pairs $(\fre,\Theta)$ consisting of a coframe $\fre$ and a symmetric $(2,0)$ tensor $\Theta$ on $\G$. A pair $(\fre,\Theta) \in \Conf(\G)$ is said to be \emph{left-invariant} if both $\fre$ and $\Theta$ are left-invariant. We will refer to pairs $(\fre,\Theta)$ satisfying the constraint equations \eqref{eq:constraintequationsIInolambda} with vanishing torsion as \emph{left-invariant Cauchy pairs}. Given a left-invariant pair $(\fre,\Theta) \in \Conf(\G)$ we write:
\begin{equation*}
	\Theta = \Theta_{ab}\, e_a\otimes e_b\, , \quad a, b = u, l, n\, ,
\end{equation*}

\noindent
where summation over repeated indices is understood. For further reference we introduce the following notation:
\begin{equation*}
	\lambda := \sqrt{\Theta_{ul}^2 + \Theta_{un}^2}\, , \quad \theta := \begin{pmatrix}
		\tll & \tln \\
		\tln & \tnn
	\end{pmatrix} \, , \quad T := \mathrm{Tr}(\theta)\, ,  \quad \Delta := \mathrm{Det}(\theta) = \Theta_{ll}\Theta_{nn} - \Theta_{ln}^2\, .
\end{equation*}

\noindent
These will play an important role in the classification of left-invariant Cauchy pairs, which was completed in \cite{Murcia:2020zig} and which we proceed to summarize.

\begin{thm}\cite[Theorem 4.9]{Murcia:2020zig}
	\label{thm:allcauchygroups}
	A connected and simply-connected Lie group $\G$ admits left-invariant parallel Cauchy pairs (respectively constrained Ricci flat parallel Cauchy pairs) if and only if $\G$ is isomorphic to one of the Lie groups listed in the table below. If that is the case, a left-invariant shape operator $\Theta$ belongs to a Cauchy pair $(\fre,\Theta)$ for certain left-invariant coframe $\fre$ if and only if $\Theta$ is of the form listed below when written in terms of $\fre = (e_u,e_l,e_n)$:
	\renewcommand{\arraystretch}{1.5}
	\begin{center}
		\begin{tabular}{|  p{1cm}| p{8.8cm} | p{3.8cm} |}
			\hline
			$\mathrm{G}$ & \emph{Cauchy parallel pair} & \emph{Constrained Ricci flat}  \\ \hline
			$\mathbb{R}^3$ & $\Theta=\Theta_{uu} e_u \otimes e_u$ &  $\Theta=\Theta_{uu} e_u \otimes e_u$   \\ \hline
			\multirow{2}*{$ \mathrm{E}(1,1)$} & $\Theta=\Theta_{uu} e_u \otimes e_u+ \Theta_{i j} e_i \otimes e_j$ & \multirow{2}*{\emph{Not allowed}} \\ & $i,j=l,n,\, \quad \Theta_{ll}=-\Theta_{nn}$ &  \\ 
			\hline \multirow{11}*{$\tau_2 \oplus \mathbb{R}$} & $\Theta=(\tul e_l +\tun e_n) \odot e_u$  & \multirow{2}*{\emph{Not allowed}} \\ & $\tul^2+\tun^2 \neq 0$ &   \\ \cline{2-3} & $\Theta=\Theta_{uu} e_u \otimes e_u+ \Theta_{i j} e_i \otimes e_j$ &  $\Theta=T e_u \otimes e_u+ \Theta_{i j} e_i \otimes e_j$  \\ & $\begin{aligned} &i,j=l,n,\, \\ &T \neq 0\, , \Delta=0 \end{aligned}$ & $\begin{aligned} &i,j=l,n,\, \\ &T \neq 0\, , \Delta=0 \end{aligned}$\\ \cline{2-3} & $\Theta=-T e_u \otimes e_u+\tul e_u \odot e_l+\tll e_l \otimes e_l\, , \quad \tul, \tll \neq 0$ & \emph{Not allowed}  \\ \cline{2-3} & $\Theta=-T e_u \otimes e_u+\tun e_u \odot e_n+\tnn e_n \otimes e_n\, , \quad \tun, \tnn \neq 0$ & \emph{Not allowed}  \\\cline{2-3} & $\Theta=-T e_u\otimes e_u + \tul e_u \odot e_l+ \tun e_u \odot e_n+\Theta_{i j} e_i \otimes e_j\,$ & \multirow{3}*{\emph{Not allowed}} \\  & $\begin{aligned}&i,j=l,n,\, \, \tln\tul\tun \neq 0\, , \\&\tnn=\frac{\tun}{\tul} \tln\, , \tll=\frac{\tul}{\tun} \tln \end{aligned}$ &  \\ & &   \\ \hline\multirow{2}*{$ \tau_{3,\mu}$} & $\Theta=\Theta_{uu} e_u \otimes e_u+ \Theta_{i j} e_i \otimes e_j$ & $\Theta=\left (\frac{T^2-2\Delta}{T} \right ) e_u \otimes e_u+ \Theta_{i j} e_i \otimes e_j$   \\ & $i,j=l,n,\, \quad T, \Delta \neq 0$& $i,j=l,n,\, \quad T, \Delta \neq 0$  \\ \hline
		\end{tabular}		
	\end{center}
	Regarding the case $\G =\tau_{3,\mu}$:
	\begin{itemize}[leftmargin=*]
		\item If $\tln \neq 0$, then
		\begin{equation*}
			\mu = \frac{T-\text{\emph{sign}}(T)\sqrt{T^2-4\Delta}}{T+\text{\emph{sign}}(T)\sqrt{T^2-4\Delta}}\, .
		\end{equation*}
		\item If $\tln=0$ and $\vert\tll\vert \geq \vert \tnn \vert $, then
		\begin{equation*}
			\mu=\frac{\tnn}{\tll}\,.
		\end{equation*}
		\item If $\tln=0$ and $\vert \tnn \vert \geq \vert \tll \vert $, then
		\begin{equation*}
			\mu=\frac{\tll}{\tnn}\,.
		\end{equation*}
	\end{itemize}
	\renewcommand{\arraystretch}{1}
\end{thm}

\noindent
The previous theorem will be used extensively in the next section. We have the following corollary.

\begin{cor}
	\label{cor:isog}
	Let $\G$ be a connected and simply connected Lie group equipped with a left-invariant Cauchy pair. Then the isomorphism type of $\G$ is prescribed by $T$, $\Delta$ and $\lambda$ as follows:
	\begin{itemize}
		\item If $T=\Delta=\lambda=0$, then $\G \simeq \mathbb{R}^3$.
		\item If $T=\lambda=0$ but $\Delta \neq 0$, then $\G \simeq \mathrm{E}(1,1)$.
		\item If  $\Delta = 0$ but $\lambda^2 + T^2 \neq 0$, then $\G \simeq \tau_2\oplus \mathbb{R}$.
		\item If $T, \Delta \neq 0$ and $\lambda=0$, then $\G \simeq \tau_{3,\mu}$.
	\end{itemize}
	
	\noindent
	Observe that the case $\lambda \neq 0$ and $\Delta \neq 0$ is not allowed. 
\end{cor}

\noindent
We are using standard notation for the groups $\G$ as explained for example in \cite[Appendix A]{Freibert}.


\subsection{Left-invariant parallel spinor flows}
\label{sec:leftinvariantspinorflow}


Let $\G$ be a simply connected three-dimensional Lie-group. We say that a parallel spinor flow $\left\{ \lambda_t ,\fre^t \right\}_{t\in \cI}$ defined on $\G$ is \emph{left-invariant} if both $\lambda_t$ and $\fre^t$ are left-invariant for every $t\in \cI$. The latter condition immediately implies that $h_{\fre^t}$ is a left-invariant Riemannian metric and $\lambda_t$ is constant for every $t\in \cI$. Let $\left\{ \fre^t\right\}_{t\in\cI}$ be a family of left-invariant coframes on $\G$. Any square matrix $\cA\in \mathrm{Mat}(3,\mathbb{R})$ acts naturally on $\left\{ \fre^t\right\}_{t\in \cI}$ as follows: 
\begin{equation*}
	\cA(\fre^t) := 
	\begin{pmatrix}
		\sum_b \cA_{ub} e^t_b  \\
		\sum_b \cA_{lb} e^t_b  \\
		\sum_b \cA_{nb} e^t_b
	\end{pmatrix}
\end{equation*}

\noindent
where we label the entries $\cA_{ab}$ of $\cA$ by the indices $a,b = u, l, n$. As a direct consequence of Proposition \ref{prop:Thetaevolutionsystem} we have the following result.

\begin{prop}
	A simply connected three-dimensional Lie group  $\G$ admits a left-invariant parallel spinor flow if and only if there exists a smooth family of non-zero constants $\left\{ \lambda_t\right\}_{t\in \cI}$ and a family $\left\{ \fre^t \right\}_{t\in \cI}$ of left-invariant coframes on $\G$ satisfying the following differential system:
	\begin{eqnarray}
		\label{eq:leftinv}
		\partial_t \fre^t  +  \lambda_t\Theta_t(\fre^t) = 0\, , \quad \dd \fre^t  = \Theta_t(\fre^t) \wedge e^t_u\, , \quad \partial_t(\Theta_t(e^t_u )) = 0\, , \quad \dd\Theta_t(e^t_u ) = 0\, ,
	\end{eqnarray}
	
	\noindent
	to which we will refer as the left-invariant (real) parallel spinor flow equations.
\end{prop}

\noindent
We will refer to solutions $\left\{ \lambda_t, \fre^t \right\}_{t\in \cI}$ of the left-invariant parallel spinor flow equations as \emph{left-invariant parallel spinor flows}. Given a parallel spinor flow $\left\{ \lambda_t, \fre^t \right\}_{t\in \cI}$, we write:
\begin{equation*}
	\Theta^t = \sum_{a,b} \Theta^t_{ab} e^t_a\otimes e^t_n\, , \qquad a, b = u, l , n\, ,
\end{equation*}

\noindent
in terms of uniquely defined functions $(\Theta^t_{ab})$ on $\cI$.

\begin{lemma}
	\label{lemma:intecondli}
	Let $\{\lambda_t, \fre^t\}_{t \in \mathcal{I}}$ be a left-invariant parallel spinor flow. The following equations hold:
	\begin{equation*}
		\begin{split}
			\partial_t \tuu^t=\lambda_t ((\tuu^t)^2 + (\tul^t )^2+(\tun^t)^2)\,,& \quad \partial_t \tul^t=\partial_t \tun^t=0\, , \quad \partial_t \tll^t=\lambda_t\tll^t \tuu^t-\lambda_t (\tul^t)^2\, \\ 
			\partial_t \tln^t= \lambda_t  \tln^t\tuu^t-\lambda_t \tun^t \tul^t\,, & \quad \partial_t \tnn^t=\lambda_t\tnn^t \tuu^t-\lambda_t
			(\tun^t)^2\, , \\
			\tln^t \tul^t=\tll^t \tun^t\, ,& \quad \tln^t \tun^t=\tnn^t \tul^t\, , \\ 
			\tll^t \tul^t+\tln^t\tun^t+\tul^t\tuu^t=0\,, &\quad \tln^t \tul^t+\tnn^t \tun^t+\tun^t \tuu^t=0\, .
		\end{split}
	\end{equation*}
	
	\noindent
	In particular, $\tul^t=\tul$ and $\tun^t=\tun$ for some constants $\tul,\tun \in \mathbb{R}$.
\end{lemma}

\begin{proof}
	A direct computation shows that equation $\partial_t(\Theta_t(e^t_u )) = 0$ is equivalent to:
	\begin{equation*}
		\partial_t \Theta^t_{ub} = \lambda_t \Theta^t_{ua} \Theta^t_{ab}\, .
	\end{equation*}
	
	\noindent
	On the other hand, equation $\dd \Theta_t(e^t_u) = 0$ is equivalent to:
	\begin{equation*}
		\Theta^t_{ua} \Theta^t_{al} = 0\, , \qquad \Theta^t_{ua} \Theta^t_{an} = 0\, .
	\end{equation*}
	
	\noindent
	The previous equations can be combined into the following equivalent conditions:
	\begin{eqnarray*}
		& \partial_t \tuu^t=\lambda_t ((\tuu^t)^2 + (\tul^t )^2+(\tun^t)^2)\, , \quad \partial_t \tul^t=\partial_t \tun^t=0\, ,\\
		& \tll^t \tul^t+\tln^t\tun^t+\tul^t\tuu^t=0\, , \quad \tln^t \tul^t+\tnn^t \tun^t+\tun^t \tuu^t=0\, ,
	\end{eqnarray*}
	
	\noindent
	which recover five of the equations in the statement. Similarly, equation $\dd (\Theta_t(\fre^t) \wedge e^t_u) = 0$ is equivalent to:
	\begin{equation*}
		\tln^t \tul^t=\tll^t \tun^t\, ,\qquad \tln^t \tun^t=\tnn^t \tul^t\, ,
	\end{equation*}
	
	\noindent
	which yields the third line of equations in the statement. We take now the exterior derivative of the first equation in \eqref{eq:leftinv} and combine the result with the second equation in \eqref{eq:leftinv}:
	\begin{equation*}
		\dd (\partial_t e^t_a  +  \lambda_t\Theta_t(e^t_a) ) =  \partial_t (\Theta^t_{ab} e^t_b\wedge e^t_u) + \lambda_t \Theta^t_{ab} \Theta^t_{bc} e^t_c\wedge e^t_u  =  (\partial_t \Theta^t_{ab} \delta_{uc} - \lambda_t \Theta^t_{ab} \Theta^t_{uc}) e^t_b\wedge e^t_c = 0\, .
	\end{equation*}
	\noindent
	Expanding the previous equation we obtain the remaining three equations in the statement and we conclude.
\end{proof}

\begin{remark}
	We will refer to the equations of Lemma \ref{lemma:intecondli} as the \emph{integrability conditions} of the left-invariant parallel spinor flow.
\end{remark}

\noindent
The following observation is crucial in order to \emph{decouple} the left-invariant parallel spinor flow equations.

\begin{lemma}
	\label{lemma:cKTheta}
	A pair $\left\{ \lambda_t, \fre^t \right\}_{t\in \cI}$ is a left-invariant parallel spinor flow if and only if there exists  a family of left-invariant two-tensors $\left\{ \cK_t \right\}_{t\in \cI}$ such that the following equations are satisfied:
	\begin{eqnarray*}
		\partial_t \fre^t  +  \lambda_t\cK_t(\fre^t) = 0\, , \quad \dd \fre^t  = \cK_t(\fre^t) \wedge e^t_u\, , \quad \partial_t(\cK_t(e^t_u )) = 0\, , \quad \dd(\cK_t(e^t_u)) = 0\, .
	\end{eqnarray*}
\end{lemma}

\begin{proof}
	The \emph{only if} direction follows immediately from the definition of left-invariant parallel Cauchy pair by taking $\{\cK_t\}_{t \in \cI}=\{\Theta_t\}_{t \in \cI}$. For the \emph{if} direction we simply compute:
	\begin{equation*}
		\Theta_t = - \frac{1}{2\lambda_t}\partial_t h_{\fre^t} = - \frac{1}{2\lambda_t} (( \partial_t e^t_a)\otimes e^t_a + e^t_a \otimes ( \partial_t e^t_a)) = \cK_t\, ,
	\end{equation*}
	
	\noindent
	hence equations \eqref{eq:leftinv} are satisfied and $\left\{ \lambda_t, \fre^t \right\}_{t\in \cI}$ is a left-invariant parallel spinor flow.
\end{proof}

\noindent
By the previous Lemma we promote the components of $\left\{\Theta_t \right\}_{t\in \cI}$ with respect to the basis $\left\{\fre^t \right\}_{t\in \cI}$ to be independent variables of the left-invariant parallel spinor flow equations \eqref{eq:leftinv}. Within this interpretation, the variables of left-invariant parallel spinor flow equations consist of triples $\left\{\lambda_t,\fre^t,\Theta^t_{ab} \right\}_{t\in \cI}$, where $\left\{\Theta^t_{ab} \right\}_{t\in \cI}$ is a family of symmetric matrices. On the other hand, the integrability conditions of Lemma \ref{lemma:intecondli} are interpreted as a system of equations for a pair $\left\{\lambda_t,\Theta^t_{ab} \right\}_{t\in \cI}$. In particular, the first equation in \eqref{eq:leftinv} is linear in the variable $\fre^t$ and can be conveniently rewritten as follows. For any family of coframes $\left\{\fre^t\right\}_{t\in\cI}$, set $\fre = \fre^0$ and consider the unique smooth path:
\begin{equation*}
	\U^t \colon \cI \to \Gl_{+}(3,\mathbb{R})\, , \qquad t\mapsto \U^t\, ,
\end{equation*}

\noindent
such that $\fre^t = \U^t(\fre)$, where  $\Gl_{+}(3,\mathbb{R})$ denotes the identity component in the general linear group $\Gl(3,\mathbb{R})$. More explicitly:
\begin{equation*}
	\fre^t_a = \sum_b \U^t_{ab} \fre_b\, , \qquad a, b = u, l ,n\, ,
\end{equation*}

\noindent
where $\U^t_{ab} \in C^{\infty}(\G)$ are the components of $\U^t$. Plugging $\fre^t = \U^t(\fre)$ in the first equation in \eqref{eq:leftinv} we obtain the following equivalent equation:
\begin{equation}
	\label{eq:ptu}
	\partial_t \U^t_{ac} + \lambda_t \Theta^t_{ab} \U^t_{bc} = 0\, ,  \quad a, b, c = u, l, n\, ,
\end{equation}

\noindent
with initial condition $\U^0 = \mathrm{Id}$. 

\noindent 
A necessary condition for a solution $\left\{\lambda_t,\Theta^t_{ab} \right\}_{t\in \cI}$ of the integrability conditions to arise from an honest left-invariant parallel spinor pair is the existence of a left-invariant coframe $\fre$ on $\Sigma$ such that $(\fre,\Theta)$ is a Cauchy pair, where $\Theta = \Theta^0_{ab} e_a\otimes e_b$. Consequently we define the set $\mathbb{I}(\Sigma)$ of \emph{admissible} solutions to the integrability equations as the set of pairs $(\left\{\lambda_t,\Theta^t_{ab} \right\}_{t \in \cI} , \fre )$ such that $\left\{\lambda_t,\Theta^t_{ab} \right\}_{t\in \cI}$ is a solution to the integrability equations and $(\fre,\Theta)$ is a left-invariant parallel Cauchy pair. 

\begin{prop}
	\label{prop:bijectionsolutions}
	There exists a natural bijection $\varphi\colon \mathbb{I}(\Sigma) \to \cP(\Sigma)$ which maps every pair:
	\begin{equation*}
		(\left\{\lambda_t,\Theta^t_{ab} \right\}_{t\in\cI} , \fre )\in \mathbb{I}(\Sigma)\, ,
	\end{equation*}
	
	\noindent
	to the pair $\left\{ \lambda_t,\fre^t = \U^t(\fre)\right\}_{t\in\cI}\in \cP(\Sigma)$, where $\left\{\U^t\right\}_{t\in\cI}$ is the unique solution of \eqref{eq:ptu} with initial condition $\U^0 = \mathrm{Id}$. 
\end{prop}

\begin{remark}
	\label{remark:inverse}
	The inverse of $\varphi$  maps every left-invariant parallel spinor flow $\left\{ \lambda_t,\fre^t\right\}_{t\in\cI}$ to the pair $(\left\{\lambda_t,\Theta^t_{ab} \right\} , \fre )$, where $\Theta^t_{ab}$ are the components of the shape operator associated to $\left\{ \lambda_t,\fre^t\right\}_{t\in\cI}$ in the basis $\left\{\fre^t\right\}_{t\in\cI}$ and $\fre = \fre^0$.
\end{remark}

\begin{proof}
	Let $(\left\{\lambda_t,\Theta^t_{ab} \right\} , \fre )\in \mathbb{I}(\Sigma)$ and let  $\left\{\U^t\right\}_{t\in\cI}$ be the solution of \eqref{eq:ptu} with initial condition $\U^0 = \mathrm{Id}$, which exists and is unique on $\cI$ by standard ODE theory \cite[Theorem 5.2]{CoddingtonLevinson}. We need to prove that $\left\{ \lambda_t,\fre^t = \U^t(\fre)\right\}_{t\in\cI}$ is a left-invariant parallel spinor flow. Since $\left\{\U^t\right\}_{t\in\cI}$  satisfies \eqref{eq:ptu} for the given $\left\{\lambda_t,\Theta^t_{ab} \right\}$, it follows that $\Theta^t = \Theta^t_{ab} e^t_a\otimes e^t_b$ is the shape operator associated to $\left\{ \lambda_t,\fre^t \right\}_{t\in\cI}$ whence the first equation in \eqref{eq:leftinv} is satisfied. On the other hand, the third and fourth equations in \eqref{eq:leftinv} are immediately implied by the integrability conditions satisfied by  $\left\{\lambda_t,\Theta^t_{ab} \right\}$. Regarding the second equation in \eqref{eq:leftinv}, we observe that the integrability conditions contain the equation $\dd (\Theta_t(\fre^t) \wedge e_u^t)=0$ and thus:
	\begin{equation}
		\label{eq:constraintintegrada}
		\dd \fre^t= \Theta^t(\fre^t) \wedge e_u^t+ \mathfrak{w}^t\, ,
	\end{equation}
	
	\noindent
	where $\{ \mathfrak{w}^t\}_{t \in \cI}$ is a family of triplets of closed two-forms on $\Sigma$. Taking the time derivative of the previous equations, plugging the exterior derivative of the first equation in \eqref{eq:leftinv} and using again the integrability conditions, we obtain that $\mathfrak{w}^t$ satisfies the following differential equation:
	\begin{equation}
		\label{eq:frwode}
		\partial_t \mathfrak{w}^t_a = -\lambda_t \Theta^t_{ad} \mathfrak{w}^t_d\, ,
	\end{equation}
	
	\noindent
	with initial condition $\mathfrak{w}^0 = \mathfrak{w}$. Restricting equation \eqref{eq:constraintintegrada} to $t=0$ it follows that $\mathfrak{w}$ satisfies:
	
	\begin{equation*}
		\dd \fre= \Theta(\fre) \wedge e_u +  \mathfrak{w}\, ,
	\end{equation*}
	
	\noindent
	Since by assumption $(\fre,\Theta)$ is left-invariant Cauchy pair, the previous equation is satisfied if and only if $\mathfrak{w} =0$ whence $\mathfrak{w}^t = 0$ by uniqueness of solutions of the linear differential equation \eqref{eq:frwode}. Therefore, the second equation in \eqref{eq:leftinv}  follows and $\varphi$ is well-defined. The fact that $\varphi$ is in addition a bijection follows directly by Remark \ref{remark:inverse} and hence we conclude. 
\end{proof}

\begin{cor}
	\label{cor:solutionsiff}
	A pair $\left\{\lambda_t,\fre^t \right\}_{t\in \cI}$ is a parallel spinor flow if and only if $(\left\{\lambda_t,\Theta^t_{ab} \right\} , \fre )$ is an admissible solution to the integrability equations.
\end{cor}

\noindent
Therefore, solving the left-invariant parallel spinor flow is equivalent to solving the integrability conditions with initial condition $\Theta_{ab}$ being part of a left-invariant parallel Cauchy pair $(\fre,\Theta)$. We remark that $\left\{\lambda_t \right\}_{t\in \cI}$ is of no relevance locally since it can be eliminated through a reparametrization of time after possibly shrinking $\cI$. However, regarding the long time existence of the flow as well as for applications to the construction of four-dimensional Lorentzian metrics it is convenient to keep track of $\cI$, whence we maintain $\left\{\lambda_t \right\}_{t\in \cI}$ in the equations.

For further reference we define a \emph{quasi-diagonal}  left-invariant parallel spinor flow as one for which $\lambda=\sqrt{\tul^2+\tun^2} =0$. Since the function $t\to \int_{0}^t \lambda_\tau \dd \tau$ is going to be a common occurrence in the following, we define:
\begin{equation*}
	\cB_{t} := \int_{0}^t \lambda_\tau \dd \tau\, .
\end{equation*}

\noindent
We distinguish now between the cases $\lambda = 0$ and $\lambda \neq 0$.

\begin{lemma}
	\label{lemma:Thetaqdiagonal}
	Let $\{\lambda_t, \fre^t\}_{t \in \cI}$ be a quasi-diagonal left-invariant parallel spinor flow. Then, the only non-zero components of $\Theta^t$ are:
	\begin{eqnarray*}
		\tuu^t =\frac{\tuu}{1- \tuu \cB_t}\, , \quad \tll^t=\frac{\tll}{1- \tuu \cB_t}\, , \quad\tln^t =\frac{\tln}{1- \tuu \cB_t}\, , \quad \tnn^t=\frac{\tnn}{1- \tuu \cB_t}\, ,
	\end{eqnarray*}
	
	\noindent
	where $\Theta^t$ is the shape operator associated to $\{\lambda_t, \fre^t\}_{t \in \cI}$  and $\Theta = \Theta^0$. Furthermore, every such $\Theta^t$ satisfies the integrability equations with quasi-diagonal initial data.
\end{lemma}

\begin{proof}
	Setting $\tul =\tun =0$ in the integrability conditions we obtain the following equations:
	\begin{equation*}
		\partial_t \tuu^t=\lambda_t (\tuu^t)^2\, , \quad \partial_t \tll^t=\lambda_t \tll^t \tuu^t\, , \quad \partial_t \tln^t=\lambda_t \tln^t\tuu^t\, , \quad \partial_t \tnn^t=\lambda_t \tnn^t\tuu^t\, .
	\end{equation*}
	
	\noindent
	whose general solution is given in the statement of the lemma.
\end{proof}

\begin{remark}
	\label{rem:maxintquasi}
	Let $\Theta_{uu}\neq 0$ and define $t_0$ to be the real number (in case it exists) with the smallest absolute value such that:
	\begin{equation*}
		\int_{0}^{t_0} \lambda_\tau \dd \tau = \Theta_{uu}^{-1}\, .
	\end{equation*}
	
	\noindent
	Then the maximal interval on which $\Theta^t$ is defined is $\cI= (-\infty, t_0)$ if $\Theta_{uu} > 0$ and $\cI= (t_0,\infty)$ if $\Theta_{uu} < 0$. This is also the maximal interval on which the left-invariant parallel spinor flow in the quasi-diagonal case can be defined. If such $t_0$ does not exist, then $\cI=\mathbb{R}$.
\end{remark}

\noindent
We consider now the non-quasi-diagonal case $\lambda  \neq 0$. Given a pair $\left\{\lambda_t,\Theta^t_{ab} \right\}_{t\in \cI}$, we introduce for convenience the following function: 
\begin{equation*}
	\cI \ni t\mapsto y_t = \lambda \, \cB_{t} + \mathrm{Arctan} \left[ \frac{\tuu}{\lambda}\right]\, ,
\end{equation*}

\noindent
where $\Theta_{ab}$ are the components of $\Theta$ in the basis $\fre$.

\begin{lemma}
	\label{lemma:Thetageneral}
	A pair $\left\{\lambda_t,\Theta^t_{ab} \right\}_{t\in \cI}$ satisfies the integrability equations with non-quasi-diagonal initial value $\Theta_{ab}$ if and only if:
	\begin{eqnarray*}
		& \tuu^t = \lambda\,\mathrm{Tan}\left[y_t \right]\, , \quad \tul^t = \tul\, , \quad \tun^t = \tun\, , \quad \tll^t = c_{ll}\, \mathrm{Sec} \left[ y_t \right] - \frac{\tul^2}{\lambda} \mathrm{Tan}\left[ y_t\right]\, ,\\
		& \tnn^t = c_{nn}\, \mathrm{Sec} \left[y_t\right] - \frac{\tun^2}{\lambda} \mathrm{Tan}\left[y_t\right]\, , \quad \tln^t = c_{ln}\, \mathrm{Sec} \left[y_t \right] - \frac{\tul \,\tun}{\lambda} \mathrm{Tan}\left[y_t\right]\, ,
	\end{eqnarray*}
	
	\noindent
	where $c_{ll}, c_{nn}, c_{ln}\in\mathbb{R}$ are real constants given by:
	\begin{equation*}
		c_{ll} =  \frac{\tll \lambda^2 + \tul^2 \tuu}{\lambda\,\sqrt{\lambda^2 +   \tuu^2}}\, , \quad c_{nn} =   \frac{\tnn \lambda^2 + \tun^2 \tuu}{\lambda\,\sqrt{\lambda^2 +   \tuu^2}}\, , \quad c_{ln} =  \frac{\tln \lambda^2 + \tul \tun \tuu}{\lambda\,\sqrt{\lambda^2 +   \tuu^2}}\, ,
	\end{equation*}
	
	\noindent
	such that the following algebraic equations are satisfied:
	\begin{eqnarray}
		\label{eq:algebraicTheta}
		& \tln \tul = \tll \tun\, , \quad  \tnn \tul = \tln \tun\, , \nonumber \\
		&\tln \tun + \tul (\tll + \tuu) = 0\, , \quad \tln \tul + \tun (\tnn + \tuu) = 0\, ,
	\end{eqnarray}
	
	\noindent
	where $\Theta_{ab}$, $a,b = u,l,n$, denote the entries of $\Theta^t_{ab}$ at $t=0$.
\end{lemma}

\begin{remark}
	Note that equations \eqref{eq:algebraicTheta} form an algebraic system for the entries of the initial condition $\Theta$, therefore restricting the allowed initial data that can be used to solve the integrability conditions. This is a manifestation of the fact that the initial data of the parallel spinor flow is constrained by the parallel Cauchy equations. The latter were solved in the left-invariant case in \cite{Murcia:2020zig}, as summarized in Theorem \ref{thm:allcauchygroups}, and its solutions can be easily verified to satisfy equations \eqref{eq:algebraicTheta} automatically.
\end{remark}

\begin{proof}
	By Lemma \ref{lemma:intecondli} we have $\partial_t\Theta^t_{ul} = \partial_t\Theta^t_{un} = 0$ whence $\tul^t = \tul, \tun^t = \tun$ for some real constants $\tul, \tun \in \mathbb{R}$. Plugging these constants into the first equation of Lemma \ref{lemma:intecondli} it becomes immediately integrable with solution:
	\begin{equation*}
		\Theta^t_{uu} =   \lambda \mathrm{Tan}\left[ \lambda (\cB_t + k_1)\right]\, ,
	\end{equation*} 
	
	\noindent
	for a certain constant $k_1\in \mathbb{R}$. Imposing $\Theta_{uu}^0 = \Theta_{uu}$ we obtain:
	\begin{equation*}
		k_1 = \frac{1}{\lambda}(\mathrm{Arctan} \left[ \frac{\tuu}{\lambda}\right] + n\pi)\, , \qquad n\in \mathbb{Z}\, ,
	\end{equation*}
	
	\noindent
	and the expression for $\Theta^t_{uu}$ follows. Plugging now $\tuu^t = \lambda\,\mathrm{Tan}\left[y_t \right]$ in the remaining differential equations of Lemma \ref{lemma:intecondli} they can be directly integrated, yielding the expressions in the statement after imposing $\Theta^0_{ab} = \Theta_{ab}$. Plugging the explicit expressions for $\Theta^t_{ab}$ in the algebraic equations of Lemma \ref{lemma:intecondli}, these can be equivalently reformulated as the algebraic system \eqref{eq:algebraicTheta} for $\Theta^t_{ab}$ at $t=0$ and we conclude. 
\end{proof}

\begin{remark}
	\label{rem:maxintnonquasi}
	Let $t_-<0$ denote the largest value for which $\lambda \cB_{t_-}+ \mathrm{Arctan}\left[ \frac{\tuu}{\lambda} \right]=-\frac{\pi}{2}$ and let $t_+>0$ denote the smallest value for which $\lambda \cB_{t_+}+ \mathrm{Arctan}\left[ \frac{\tuu}{\lambda} \right]=\frac{\pi}{2}$ (if $t_-$, $t_+$ or both do not exist, we take by convention $t_{\pm} =\pm \infty$). Then, the maximal interval of definition on which $\Theta^t$ is defined is $\mathcal{I}=(t_-,t_+)$. 
\end{remark}

\subsection{Classification of left-invariant spinor flows}


Proposition \ref{lemma:intecondli} states that $\tul^t=\tul$ and $\tun^t=\tun$ for constants $\tul, \tun \in \mathbb{R}$. Therefore, we proceed to classify left-invariant parallel spinor flows in terms of the possible values of $\tul$ and $\tun$. We begin with the classification of quasi-diagonal left-invariant parallel spinor flows, defined by the condition $\tul=\tun=0$, that is, $\lambda = 0$.  

\begin{prop}
	\label{prop:liqd} 
	Let $\{\lambda_t, \fre^t\}_{t \in \cI}$ be a quasi-diagonal left-invariant parallel spinor flow with initial data $(\fre,\Theta)$ satisfying $\Theta_{uu} \neq 0$. Define $Q$ to be the orthogonal two by two matrix diagonalizing $\theta/\Theta_{uu}$ as follows:
	\begin{equation*}
		\frac{\theta}{\Theta_{uu}} =   Q	\begin{pmatrix}
			\rho_{+} & 0 \\
			0 & \rho_{-}
		\end{pmatrix}  
		Q^{\ast}
	\end{equation*}
	
	\noindent
	with eigenvalues $\rho_{+}$ and $\rho_{-}$ and where $Q^\ast$ denotes the matrix transpose of $Q$. Then:
	\begin{equation*}
		e^t_u = (1-\tuu \cB_t) e_u\, , \quad 	\begin{pmatrix}
			e^t_l \\
			e^t_n  
		\end{pmatrix} = Q 
		\begin{pmatrix}
			\left[1- \Theta_{uu}\cB_t \right]^{\rho_{+}} & 0 \\
			0 &  \left[1- \Theta_{uu}\cB_t \right]^{\rho_{-}}
		\end{pmatrix}  Q^{\ast} \begin{pmatrix}
			e_l \\
			e_n  
		\end{pmatrix}
	\end{equation*}

	\noindent
	Conversely, for every family of functions $\left\{ \lambda_{t}\right\}_{t\in\cI}$ the previous expression defines a parallel spinor flow on $\G$. The case $\tuu=0$ is recovered by taking the formal limit $\tuu \rightarrow 0$. 
\end{prop}

\begin{proof}
	Define the function $\cI\ni t\to x_t := \mathrm{Log} \left[1-\tuu \cB_t \right]$. By Proposition \ref{prop:bijectionsolutions} and Corollary \ref{cor:solutionsiff} it suffices to use the explicit expression for $\Theta^t$ obtained in Lemma \ref{lemma:Thetaqdiagonal} to solve Equation \eqref{eq:ptu} with initial condition $\U^0 = \mathrm{Id}$ on a simply connected Lie group admitting quasi-diagonal parallel Cauchy pairs. Plugging the explicit expression of $\Theta^t$ in \eqref{eq:ptu} we obtain:
	\begin{equation}
		\label{eq:ptu2}
		\begin{pmatrix}
			\partial_t \U_{uu}^t & \partial_t \U_{ul}^t & \partial_t \U_{un}^t \\
			\partial_t \U_{lu}^t & \partial_t \U_{ll}^t & \partial_t \U_{ln}^t \\
			\partial_t \U_{nu}^t & \partial_t \U_{nl}^t & \partial_t \U_{nn}^t 
		\end{pmatrix} = \frac{\partial_t x_t}{\tuu} 
		\begin{pmatrix}
			\U_{uu}^t \Theta_{uu} & \U_{ul}^t \Theta_{uu} & \U_{un}^t \Theta_{uu}\\
			\U_{lu}^t \Theta_{ll} + \U_{nu}^t \Theta_{ln} & \U_{ll}^t \Theta_{ll} + \U_{nl}^t \Theta_{ln} & \U_{ln}^t \Theta_{ll} + \U_{nn}^t \Theta_{ln} \\
			\U_{lu}^t \Theta_{ln} + \U_{nu}^t \Theta_{nn} & \U_{ll}^t \Theta_{ln} + \U_{nl}^t \Theta_{nn} & \U_{ln}^t \Theta_{ln} + \U_{nn}^t \Theta_{nn}
		\end{pmatrix}  
	\end{equation}
	
	\noindent
	in terms of the initial data  $\Theta_{ab}$. The previous differential system can be equivalently written as follows:
	\begin{equation*}
		\partial_t\, \U^t_{uc} = \partial_t x_t \, \U^t_{uc}\, ,  \qquad \begin{pmatrix}
			\partial_t \U_{lc}^t \\
			\partial_t \U_{nc}^t   
		\end{pmatrix} = \frac{\partial_t x_t}{\tuu} 
		\begin{pmatrix}
			\Theta_{ll}   &   \Theta_{ln} \\
			\Theta_{ln}   &  \Theta_{nn}  
		\end{pmatrix}  \begin{pmatrix}
			\U_{lc}^t \\
			\U_{nc}^t   
		\end{pmatrix} \, , \quad c=u,l,n\, .
	\end{equation*}
	
	\noindent
	The general solution to the equations for $\U^t_{uc}$ with initial condition $\U^0 = \mathrm{Id}$ is given by:
	\begin{equation*}
		\U^t_{uu}=1-\tuu \cB_t\, , \quad \U^t_{ul}=\U^t_{un} =0\, .
	\end{equation*} 
	
	\noindent
	Consider now the diagonalization of the constant matrix occurring in the differential equations for $\U^t_{ic}$:
	\begin{equation*}
		\frac{1}{\tuu} 
		\begin{pmatrix}
			\Theta_{ll}   &   \Theta_{ln} \\
			\Theta_{ln}   &  \Theta_{nn}  
		\end{pmatrix}   = Q\begin{pmatrix}
			\rho_{+}   &   0 \\
			0   &  \rho_{-}  
		\end{pmatrix} Q^{\ast} \, ,
	\end{equation*}
	
	\noindent
	where $Q$ is a two by two orthogonal matrix and $Q^{\ast}$ is its transpose. The eigenvalues are explicitly given by:
	\begin{equation*}
		\rho_{\pm} = \frac{T \pm \sqrt{T^2- 4\Delta}}{2\Theta_{uu}} \, .  
	\end{equation*}
	
	\noindent
	We obtain:
	\begin{equation*}
		Q^{\ast} \begin{pmatrix}
			\partial_t \U_{lc}^t \\
			\partial_t \U_{nc}^t   
		\end{pmatrix} = \partial_t x_t 
		\begin{pmatrix}
			\rho_{+}   &   0 \\
			0   &  \rho_{-}  
		\end{pmatrix} Q^{\ast}  \begin{pmatrix}
			\U_{lc}^t \\
			\U_{nc}^t   
		\end{pmatrix} \, , \quad c=u,l,n\, , 
	\end{equation*}
	
	\noindent
	whose general solution is given by:
	\begin{equation*}
		\begin{pmatrix}
			\U_{lc}^t \\
			\U_{nc}^t   
		\end{pmatrix} = Q \begin{pmatrix}
			k^{+}_c e^{\rho_{+} x_t} \\
			k^{-}_c e^{\rho_{-} x_t}   
		\end{pmatrix} = Q \begin{pmatrix}
			k^{+}_c \left[1- \Theta_{uu}\cB_t \right]^{\rho_{+}} \\
			k^{-}_c \left[1- \Theta_{uu}\cB_t \right]^{\rho_{-}}  
		\end{pmatrix} \, , \quad c=u,l,n\, , 
	\end{equation*}
	
	\noindent
	for constants $k^{+}_c, k^{-}_c\in \mathbb{R}$. Imposing the initial condition $\U^0 =\mathrm{Id}$ we obtain the following expression for $k^{+}_c$ and $k^{-}_c$:
	\begin{equation*}
		\begin{pmatrix}
			k^{+}_c   \\
			k^{-}_c     
		\end{pmatrix} = Q^{\ast}\begin{pmatrix}
			\delta_{lc} \\
			\delta_{nc}   
		\end{pmatrix}  \, , \quad c=u,l,n\, , 
	\end{equation*}
	
	\noindent
	whence:
	\begin{equation*}
		\begin{pmatrix}
			k^{+}_u   \\
			k^{-}_u     
		\end{pmatrix} =0 \, , \quad \begin{pmatrix}
			k^{+}_l   \\
			k^{-}_l     
		\end{pmatrix} = Q^{\ast}\begin{pmatrix}
			1 \\
			0    
		\end{pmatrix} = \begin{pmatrix}
			Q^{\ast}_{ll} \\
			Q^{\ast}_{nl}  
		\end{pmatrix} \, , \quad \begin{pmatrix}
			k^{+}_n   \\
			k^{-}_n     
		\end{pmatrix} = Q^{\ast}\begin{pmatrix}
			0 \\
			1   
		\end{pmatrix} = \begin{pmatrix}
			Q^{\ast}_{ln} \\
			Q^{\ast}_{nn}  
		\end{pmatrix}  \, . 
	\end{equation*}
	
	\noindent
	We conclude that:
	\begin{equation*}
		\begin{pmatrix}
			\U_{ll}^t &  \U_{ln}^t \\
			\U_{nl}^t &   \U_{nn}^t 
		\end{pmatrix} =  Q 
		\begin{pmatrix}
			e^{\rho_{+} x_t} & 0 \\
			0 &  e^{\rho_{-} x_t}
		\end{pmatrix}  Q^{\ast} = Q 
		\begin{pmatrix}
			\left[1- \Theta_{uu}\cB_t \right]^{\rho_{+}} & 0 \\
			0 &  \left[1- \Theta_{uu}\cB_t \right]^{\rho_{-}}
		\end{pmatrix}  Q^{\ast}
	\end{equation*}
	
	\noindent
	and the statement is proven. The converse follows by construction upon use of Lemma \ref{lemma:cKTheta} and Proposition \ref{prop:bijectionsolutions}. It can be easily seen that the case $\Theta_{uu} =0$ is obtained by taking the formal limit $\tuu \rightarrow 0$ and we conclude.
\end{proof}

\begin{remark}
	The Ricci tensor of the family of Riemannian metrics $\left\{ h_{\fre^t}\right\}_{t\in \cI}$ associated to a left-invariant quasi-diagonal parallel spinor flow $\left\{\lambda_t,\fre^t\right\}_{t\in\cI}$ is given by:
	\begin{equation}
		\ric^{h_{\fre^t}}=-T^t\Theta_t+\frac{\mathcal{H}_t}{2} e_u^t \otimes e_u^t\, ,
	\end{equation}
	
	\noindent
	where $T^t = \Theta^t_{ll} + \Theta^t_{nn}$. If $\cH_t = 0$ for every $t\in \cI$, that is, if the parallel Cauchy pair defined by $\left\{\lambda_t,\fre^t\right\}_{t\in\cI}$ is constrained Ricci flat, then:
	\begin{equation}
		\ric^{h_{\fre^t}}=\frac{T^t}{2\lambda_t}\partial_t h_{\fre^t}\, ,
	\end{equation}
	
	\noindent
	which, after a reparametrization of the time coordinate can be brought into the form $\ric^{h_{\fre^{\tau}}} = -2\partial_{\tau} h_{\fre^{\tau}}$ after possibly shrinking $\cI$. Hence, this gives a particular example of a left-invariant Ricci flow on $\G$.
\end{remark}

\noindent
We consider now $\tul \tun=0$ but $\tul^2+\tun^2 \neq 0$. This case necessarily corresponds to $\G = \tau_2 \oplus \mathbb{R}$.

\begin{prop}
	\label{prop:linocero}
	Let $\{\lambda_t, \fre^t\}_{t \in \cI}$ be a left-invariant parallel spinor flow  with initial parallel Cauchy pair $(\fre,\Theta)$ satisfying $\tul \tun=0$ and $\lambda \neq 0$. Then:
	\begin{itemize}[leftmargin=*]
		\item If $\tul=0 $ the following holds:
		\begin{eqnarray*}
			& e_u^t =\left (1-\tuu \cB_t \right )\, e_u- \tun\cB_t\, \, e_n\, , \quad e_l^t = e_l \, , \\
			& e_n^t = \left ( \frac{\tuu}{\tun}-\frac{\lambda}{\tun}(1-\tuu\cB_t  )\mathrm{Tan}\left[y_t \right]\right) e_u +\left (1+\lambda \cB_t \, \mathrm{Tan}\left[y_t\right]\right) e_n\, .
		\end{eqnarray*}
		
		\item  If $\tun=0 $ the following holds:
		\begin{eqnarray*}
			& e_u^t =\left (1-\tuu \cB_t \right )\, e_u- \tul \cB_t\, \, e_l\, , \quad e_n^t = e_n \, , \\
			& e_l^t = \left ( \frac{\tuu}{\tul}-\frac{\lambda}{\tul}(1-\tuu\cB_t  )\mathrm{Tan}\left[y_t \right]\right) e_u +\left (1+\lambda \cB_t \, \mathrm{Tan}\left[y_t\right]\right) e_l\, .
		\end{eqnarray*}
	\end{itemize}
	
	\noindent
	Conversely, every such family $\{\lambda_t, \fre^t\}_{t \in \cI}$ is a left-invariant parallel spinor flow for every $\{\lambda_t\}_{t \in \cI}$. 
\end{prop}

\begin{proof}
	We prove the case $\tul=0$ and $\tun \neq 0$ since the case $\tun=0$ and $\tul \neq 0$ follows similarly. Setting $\tul=0$ and assuming $\tun \neq 0$ in Lemma \ref{lemma:Thetageneral} we immediately obtain:
	\begin{equation*}
		\tuu^t=-\tnn^t=\lambda\mathrm{Tan} \left[ y_t\right]\, , \qquad \Theta^t_{ll} = \Theta^t_{ln} = 0\, .
	\end{equation*}
	
	\noindent
	where we have also used that, in this case, $\Theta_{ln} = \Theta_{ll} = 0$ and $\tuu=-\tnn$ as summarized in Theorem \ref{thm:allcauchygroups}. Hence: 
	\begin{equation*}
		\Theta^t = \Theta_{un}
		\begin{pmatrix}
			0 & 0 & 1\\
			0 & 0 & 0\\
			1 & 0 & 0
		\end{pmatrix}  +  \lambda\mathrm{Tan} \left[ y_t\right]
		\begin{pmatrix}
			1 & 0 & 0\\
			0 & 0 & 0\\
			0 & 0 & -1
		\end{pmatrix} \, ,
	\end{equation*}
	
	\noindent
	and Equation \eqref{eq:ptu} reduces to:
	\begin{equation*}
		\partial_t \U^t + \lambda_t \Theta_{un} \begin{pmatrix}
			\U_{nu}^t & \U_{nl}^t & \U_{nn}^t\\
			0 & 0 & 0\\
			\U_{uu}^t & \U_{ul}^t & \U_{un}^t 
		\end{pmatrix} + \lambda \lambda_t \mathrm{Tan} \left[ y_t\right]
		\begin{pmatrix}
			\U_{uu}^t & \U_{ul}^t & \U_{un}^t\\
			0 & 0 & 0\\
			-\U_{nu}^t & -\U_{nl}^t & -\U_{nn}^t
		\end{pmatrix} = 0 \, ,
	\end{equation*}
	
	\noindent
	or, equivalently:
	\begin{eqnarray*}
		\partial_t \U^t_{uc} + \lambda_t \lambda\mathrm{Tan} \left[ y_t\right] \U^t_{uc}+\lambda_t \tun \U^t_{nc} = 0\, , \quad \partial_t\, \U^t_{lc} = 0\, , \quad
		\partial_t\, \U^t_{nc} -\lambda_t \lambda\mathrm{Tan} \left[ y_t\right] \U_{nc}+\lambda_t \tun \U^t_{uc} =0\, . 
	\end{eqnarray*}
	
	\noindent
	The general solution to this system with initial condition $\U^{0}=\mathrm{Id}$ is given by:
	\begin{eqnarray*}
		& \U^t_{uu} = 1-\tuu \cB_t\, , \quad \U^t_{un}=-\tun \cB_t\, , \quad \U^t_{ul}=\U^t_{lu}=\U^t_{ln}=\U^t_{nl}=0\, , \quad \U^t_{ll}=1\,, \\ 
		& \U^t_{nn} = 1+ \lambda \cB_t  \mathrm{Tan}\left[ y_t \right]\, , \quad \U^t_{nu} = \frac{\tuu}{\tun}-\frac{\lambda}{\tun}(1-\tuu\cB_t  ) \mathrm{Tan}\left[ y_t \right] \, ,
	\end{eqnarray*}
	
	\noindent
	which implies the statement. The converse follows by construction upon use of Lemma \ref{lemma:cKTheta} and Proposition \ref{prop:bijectionsolutions}.
\end{proof}

\begin{remark}
	The Ricci tensor of the family of metrics $\left\{ h_{\fre^t}\right\}_{t\in\cI}$ associated to a left-invariant parallel spinor with  if $\tul=0$ but $\tun \neq 0$ reads:
	\begin{equation*}
		\ric^{h_{\fre^t}}=-\Theta_t \circ \Theta_t=\frac{\mathcal{H}_t}{4}  (h_{\fre^t}-e_n^t \otimes e_n^t)\, .
	\end{equation*}
	
	\noindent
	Recall that $\nabla^{h_{\fre^t}} e_n^t=0$ and thus $\left\{h_{\fre^t}, e_n^t \right\}_{t\in\cI}$ defines a family of $\eta\,$-Einstein cosymplectic structures \cite{Capelletti,Olszak}. On the other hand, if  $\tun=0$ but $\tul \neq 0$ the curvature of $\left\{ h_{\fre^t}\right\}_{t\in\cI}$ is given by: 
	\begin{equation*}
		\ric^{h_{\fre^t}}=-\Theta_t \circ \Theta_t=\frac{\mathcal{H}_t}{4}  (h_{\fre^t}-e_l^t \otimes e_l^t)\, .
	\end{equation*}
	
	\noindent
	whence $\{ h_{\fre^t}, e_l^t\}_{t\in\cI}$ defines as well a family of $\eta\,$-Einstein cosymplectic structures on $\G$.
\end{remark}

\noindent
Finally we consider $\tul \tun \neq 0$, a case that again corresponds to $\G = \tau_{2} \oplus \mathbb{R}$.
\begin{prop}
	\label{prop:tultunneq0}
	Let $\{\lambda_t, \fre^t\}_{t \in \cI}$ be a left-invariant parallel spinor flow  with initial parallel Cauchy pair $(\fre,\Theta)$ satisfying $\tul \tun \neq 0$. Then:
	\begin{eqnarray*}
		& e_u^t = e_u+\Bt (T e_u-\tul e_l-\tun e_n)\, ,\\
		& e_l^t = -\frac{\tul}{\lambda} \left( \frac{T}{\lambda} + \left (1+ T \Bt\right ) \mathrm{Tan} [y_t] \right)\, e_u +\left ( 1+\frac{\tul^2 \mathcal{B}_t}{\lambda} \mathrm{Tan} [y_t]  \right) e_l+\frac{\tul \tun \Bt }{\lambda} \mathrm{Tan} [y_t] \, \, e_n\, ,\\ 
		& e_n^t = -\frac{\tun}{\lambda} \left (\frac{T}{\lambda} + \left ( 1+ T \Bt \right) \mathrm{Tan} [y_t]\right) \, e_u  + \frac{\tul \tun \Bt}{\lambda}  \mathrm{Tan} [y_t] \, e_l+\left (1+\frac{\tun^2 \Bt}{ \lambda} \ \mathrm{Tan} [y_t]  \right )\, e_n\, ,  
	\end{eqnarray*}
	
	\noindent
	Conversely, every such family $\{\lambda_t, \fre^t\}_{t \in \cI}$ is a left-invariant parallel spinor flow for every $\{\lambda_t\}_{t \in \cI}$. 
\end{prop}

\begin{proof}
	Assuming $\tul,\tun\neq 0$ in Lemma \ref{lemma:Thetageneral} we obtain:
	\begin{eqnarray*}
		\tuu^t = \lambda \mathrm{Tan}\left[ y_t \right]\, , \quad \tll^t =\frac{\tul}{\tun}\tln^t\, , \quad \tnn^t=\frac{\tun}{\tul} \tln^t \, , \quad \tln^t=-\frac{\tul \tun}{\tul^2+\tun^2} \tuu^t\, .
	\end{eqnarray*}
	
	\noindent
	Note that $\tuu^t=-\tll^t-\tnn^t$. Hence:
	\begin{equation*}
		\Theta^t = 
		\begin{pmatrix}
			0 & \Theta_{ul} & \Theta_{un}\\
			\Theta_{ul} & 0 & 0\\
			\Theta_{un} & 0 & 0
		\end{pmatrix}  -  \frac{\mathrm{Tan} \left[ y_t\right]}{\lambda}
		\begin{pmatrix}
			-\lambda^2  & 0 & 0\\
			0 & \Theta_{ul}^2 & \Theta_{ul}\Theta_{un}\\
			0 & \Theta_{ul}\Theta_{un} & \Theta_{un}^2
		\end{pmatrix} \, ,
	\end{equation*}
	
	\noindent
	and Equation \eqref{eq:ptu} reduces to:
	\begin{eqnarray*}
		& \frac{1}{\lambda_t}\partial_t \U^t_{uc} + \U^t_{lc} \tul + \U^t_{nc} \tun + \U^t_{uc} \tuu^t  =0\, ,\\
		& \frac{1}{\lambda_t}\partial_t \U^t_{lc} + \tul \left (\U^t_{uc}- \lambda^{-1} (\tul \U^t_{lc}+\tun \U_{nc}^t) \mathrm{Tan} \left[ y_t\right] \right ) =0\, ,\\
		& \frac{1}{\lambda_t}\partial_t \U^t_{nc} + \tun \left (\U^t_{uc}- \lambda^{-1} (\tul \U^t_{lc}+\tun \U_{nc}^t) \mathrm{Tan} \left[ y_t\right] \right ) =0\, .
	\end{eqnarray*}
	
	\noindent
	The general solution to this system with initial condition $\U^{0}=\mathrm{Id}$ is given by:
	\begin{eqnarray*}
		& \U^t_{uu} =1 - \Theta_{uu} \Bt\, , \quad \U^t_{ul}=-\tul \Bt\, , \quad \U^t_{un}=-\tun \Bt\, ,\\ 
		& \U^t_{lu} =  \frac{\tul}{\lambda} \left(\frac{\Theta_{uu}}{\lambda} - \left (1 - \Theta_{uu}\Bt\right ) \mathrm{Tan}[y_t] \right)\, ,\quad \U^t_{ll} = 1+\frac{\tul^2 \mathcal{B}_t}{\lambda} \mathrm{Tan}[y_t] \, , \quad \U^t_{ln}=\frac{\tul \tun \Bt }{\lambda} \mathrm{Tan}[y_t]\, , \\ 
		& \U^t_{nu} =\frac{\tun}{\lambda} \left(\frac{\Theta_{uu}}{\lambda} - \left (1 - \Theta_{uu}\Bt\right ) \mathrm{Tan}[y_t] \right)\, , \quad  \U^t_{nl} =\frac{\tul \tun \Bt}{\lambda} \mathrm{Tan}[y_t]\, , \quad 
		\U^t_{nn}=1+\frac{\tun^2 \Bt}{\lambda}  \mathrm{Tan}[y_t] \, , 
	\end{eqnarray*}
	
	\noindent
	and we conclude.
\end{proof}

\begin{remark}
	The three-dimensional Ricci tensor of the family of Riemannian metrics $\left\{h_{\fre^t}\right\}_{t\in \cI}$ associated to a left invariant parallel spinor flow with $\Theta_{ul}\Theta_{un} \neq 0$ reads:
	\begin{equation*}
		\ric^{h_{\fre^t}}=-\Theta_t \circ \Theta_t =\frac{\mathcal{H}_t}{4} (h_{\fre^t}-\eta_t \otimes \eta_t)\, , \quad \eta_t=\frac{1}{\sqrt{\tul^2+\tun^2}}( \tun e_l^t-\tul e_n^t)\, ,
	\end{equation*}
	
	\noindent
	Note that $\nabla^{h_{\fre^t}} \eta_t=0$, so $\{ h_{\fre^t},\eta_t\}_{t \in \cI}$ defines a family of $\eta\,$-Einstein cosymplectic Riemannian structures on $\G$.
\end{remark}

\noindent
As a corollary to the classification of left-invariant parallel spinor flows presented in Propositions \ref{prop:liqd}, \ref{prop:linocero} and \ref{prop:tultunneq0} we can explicitly obtain the evolution of the Hamiltonian constraint in each case. 

\begin{cor}
	Let $\left\{\lambda_t,\fre^t\right\}_{t\in\cI}$ be a left-invariant parallel spinor in $(M,g)$. 
	\begin{itemize}[leftmargin=*]
		\item If $\tul=\tun=0$, then $\mathcal{H}_t=\frac{\mathcal{H}_{0}}{(1-\tuu \cB_t)^2}$.
		
		\item If $\tul=0$ but $\tun \neq 0$ then $\mathcal{H}_t=\frac{\tun^2 \mathcal{H}_{0}}{\tuu^2+\tun^2} \sec^2 \left[    \lambda  \cB_t  + \mathrm{Arctan}\left[ \frac{\tuu}{\lambda} \right] \right]$.

		\item If $\tun=0$ but $\tul \neq 0$ then $\mathcal{H}_t=\frac{\tul^2 \mathcal{H}_{0}}{\tuu^2+\tul^2} \sec^2 \left[     \lambda \cB_t  + \mathrm{Arctan}\left[ \frac{\tuu}{\lambda}\right] \right]$.

		\item If $\tul,\tun \neq 0$ then $\mathcal{H}_t=\frac{\lambda^2 \mathcal{H}_{0}}{\lambda^2 + \tuu^2} \sec^2 \left[ \lambda \cB_t + \mathrm{Arctan}\left[ \frac{\tuu}{\lambda}\right] \right] $.
	\end{itemize}
	
	\noindent
	where $\cH_0$ is the Hamiltonian constraint at time $t=0$. 
\end{cor}

\noindent
Since the secant function has no zeroes, the Hamiltonian constraint vanishes for a given $t\in \cI$, and hence for every $t\in \cI$, if and only if it vanish at $t=0$. Theorem \ref{thm:allcauchygroups} implies that only quasi-diagonal left-invariant parallel spinor flows admit constrained Ricci flat initial data. Therefore the Hamiltonian constraint of left-invariant parallel spinor flows with $\lambda \neq 0$ is non-vanishing for every $t\in\cI$ and such left-invariant parallel spinor flows cannot produce four-dimensional Ricci flat Lorentzian metrics.


\section{Globally hyperbolic reduction of abelian bundle gerbes}
\label{sec:reductionbundlegerbe}


In Section \ref{sec:evolutionskewtorsion} we have obtained the globally hyperbolic reduction of a skew-torsion isotropic parallelism, since it plays a fundamental role in the reduction of globally hyperbolic supersymmetric NS-NS solutions. In this section we take care of the globally hyperbolic reduction of another fundamental ingredient of NS-NS supergravity and its supersymmetric solutions and configurations, namely an abelian bundle gerbe equipped with a connection and a curving \cite{Murray}. The reader is referred to Appendix \ref{chapter:BundleGerbes} for the basic details on bundle gerbes. Let $[u,v,l,n]$ be a globally hyperbolic parallelism on $M = \cI \times \Sigma$ with associated metric:
\begin{equation*}
g = u \odot v + l\otimes + n\otimes n = -\lambda^2_t \dd t\otimes \dd t + h_t\, .
\end{equation*}

\noindent
Let $\bar{\cC} := (\bar{\cP} , \bar{\cA}, \bar{Y})$ be a bundle gerbe on $M = \cI\times \Sigma$. Since the goal is to apply this reduction to the evolution problem defined by NS-NS supergravity, we will assume that the \emph{topological data} contained in $\cC$, namely $(\bar{\cP},\bar{Y})$, is given by the pull-back of a bundle gerbe $(\cP,Y)$ defined over $\Sigma$, while the connective structure $\bar{\cA}$ is allowed to be genuinely time-dependent, not necessarily given by the pull-back of a fixed connective structure on $(\cP,Y)$. We will refer to such bundle gerbes as \emph{reducible} bundle gerbes on $\cI\times \Sigma$. This set-up captures the time-dependence of the evolution problem on a bundle gerbe in its full generality. Since the submersion $\bar{Y} \to M$ of a reducible bundle gerbe on $\cI\times\Sigma$ is by assumption a pull-back submersion it follows that it is isomorphic to:
\begin{equation*}
\bar{Y} = \cI \times Y \to \cI \times \Sigma	
\end{equation*}

\noindent
where the projection on the first factor is the identity map and $Y\to \Sigma$ is the submersion underlying the bundle gerbe $(\cP , Y)$. Similarly:
\begin{equation*}
\bar{Y}\times_M  \bar{Y} = \cI\times Y \times_\Sigma Y \to \cI \times \Sigma
\end{equation*}

\noindent
with the identity map $\cI\to \cI$ acting on the first factor. Since $\bar{\cP}$ is the pull-back of $\cP\to Y\times_{\Sigma} Y$ to $\cI \times Y\times_{\Sigma} Y$, we have:
\begin{eqnarray*}
\bar{\cP}= \cI \times \cP \to \cI \times Y\times_{\Sigma} Y
\end{eqnarray*}

\noindent
again with the identity map $\cI\to \cI$ acting on the first factor. The principal $\U(1)$ action of $\bar{\cP} = \cI \times \cP$ is the one induced by the principal $\U(1)$ action of $\cP$ and the trivial action on the $\cI$ factor. Let $\bar{\cA} \in \Omega^1(\mathbb{R}\times \cP , \mathbb{R})$ be a connective structure on $\bar{\cP}$. Then we can write:
\begin{equation*}
\bar{\cA} = \Psi_t \dd t + \cA_t
\end{equation*}

\noindent
where $\Psi_t$ is a family of functions on $\cP$ and $\cA_t$ is a family of connections on $\cP$. Recall that $t$ denotes the Cartesian coordinate on $\cI$. Since $\bar{\cA}$ is a connection on $\bar{\cP}$, it is in particular invariant under the $\U(1)$ action of the principal bundle $\bar{\cP}$, and since this action is ineffective on $\cI$ we conclude that $\left\{ \Psi_t\right\}$ is in fact a family of invariant functions on $\cP$ whence it descends to a family of invariant functions on $Y\times_{\Sigma} Y$ that we denote by the same symbol for ease of notation. Since $(\bar{\cP} , \bar{\cA}, \bar{Y})$ is a bundle gerbe, it comes equipped with an isomorphism:
\begin{equation*}
 (\bar{\pi}^{\ast}_{12}\bar{\cP}\otimes\bar{\pi}^{\ast}_{23}\bar{\cP} , \bar{\pi}^{\ast}_{12}\bar{\cA}\otimes \bar{\pi}^{\ast}_{23}\bar{\cA})  \xrightarrow{\simeq} (\bar{\pi}^{\ast}_{13} \bar{\cP} , \bar{\pi}^{\ast}_{13}\bar{\cA})
\end{equation*}

\noindent
where $\bar{\pi}\colon \bar{Y} \to M$ and the notation $\pi_{ij} \colon \bar{Y}\times_M \bar{Y}\times_M \bar{Y} \to \bar{Y}\times_M \bar{Y}$ forgets the entry not labelled neither by $i$ nor $j$. This implies $\bar{\delta} \Psi_t = 0$, where, as explained in Appendix \ref{chapter:BundleGerbes}, $\bar{\delta}$ is the differential of the simplicial manifold determined by $\bar{\pi}\colon \bar{Y} \to M$. Since this differential is exact, it follows that $\Psi_t$ descends through $\bar{\delta}$ to a family of functions on $Y$ that we denote again by $\Psi_t$. The preceding discussion thus leads to the following correspondence. 

\begin{prop}
There is a natural one-to-one correspondence between reducible bundle gerbes $(\bar{\cP},\bar{Y},\bar{\cA})$ on $M = \cI\times \Sigma$ and tuples $(\cP,Y,\cA_t , \Psi_t)$, where  $(\cP,Y)$ is a bundle gerbe on $\Sigma$, $\Psi_t$ is a family of functions on $Y$, and $\cA_t$ is a family of connective structures on $(\cP,Y)$.
\end{prop}

\noindent
We will refer to $(\cP,Y, \cA_t , \Psi_t )$ as the \emph{reduction} of $(\bar{\cP},\bar{Y},\bar{\cA})$. A direct computation gives the following result.

\begin{lemma}
\label{lemma:gerbecurvaturedecomposition}
The curvature $\cF_{\bar{\cA}} \in \Omega^2(\mathbb{R}\times Y\times_{\Sigma} Y)$ of $\bar{\cA} = \Psi_t \dd t + \cA_t$ satisfies the following equation:
\begin{eqnarray*}
\cF_{\bar{\cA}} = \dd_{Y^{[2]}} \Psi_t \wedge \dd t + \dd t \wedge \partial_t \cA_t + \cF_{\cA_t}
\end{eqnarray*}

\noindent
where $\dd_{Y^{[2]}}\colon \Omega^{\bullet}(Y\times_{\Sigma} Y) \to \Omega^{\bullet}(Y\times_{\Sigma} Y) $ is the exterior derivative on $\Omega^{\bullet}(Y\times_{\Sigma} Y)$ and $\cF_{\cA_t}$ is the family of curvatures of $\cA_t$.
\end{lemma}

\noindent
A curving on $(\bar{\cP},\bar{Y},\bar{\cA})$ is by definition a two-form $\bar{b}\in \Omega^2(\bar{Y}\times_M \bar{Y})$ satisfying:
\begin{equation*}
\bar{\delta} \bar{b} = F_{\bar{\cA}}
\end{equation*}

\noindent 
Since $\bar{Y}\times_M \bar{Y} = \cI \times Y \times_{\Sigma} Y$, we can write every curving $\bar{b}$ on $(\bar{\cP},\bar{Y},\bar{\cA})$ as follows:
\begin{equation}
\label{eq:decompositionbarb}
\bar{b} = \dd t \wedge a_t + b_t
\end{equation}

\noindent
for uniquely determined families of one forms $a_t$ and two-forms $b_t$ on $Y$.
\begin{lemma}
\label{lemma:reduciblegerbe}
Let $(\bar{\cP},\bar{Y},\bar{\cA})$ be a reducible bundle gerbe with connective structure on $M = \cI \times \Sigma$ and let  $(\cP,Y, \cA_t, \Psi_t)$ be its reduction on $\Sigma$. A two-form $\bar{b}\in \Omega^2(\bar{Y}\times_M \bar{Y})$ is a curving on $(\bar{\cP},\bar{Y},\bar{\cA})$ if and only if, in the decomposition given in \eqref{eq:decompositionbarb}, $\left\{ b_t\right\}$ is a family of curvings on $(\cP,Y,\cA)$ and $a_t$ is a family of one-forms on $Y$ satisfying:
\begin{equation}
\label{eq:conditionat}
\delta a_t = \partial_t \cA_t - \dd_{Y^{[2]}} \Psi_t
\end{equation}

\noindent
for every $t\in \cI$, where $\delta$ is the simplicial differential of the simplicial manifold determined by submersion $\pi \colon Y \to M$.
\end{lemma}

\begin{proof}
Plugging equation \eqref{eq:decompositionbarb} in $\bar{\delta} \bar{b} = F_{\bar{\cA}}$ we have:
\begin{equation*}
\bar{\delta}\bar{b} = \bar{\delta}(\dd t \wedge a_t + b_t) = \dd t \wedge \delta a_t + \delta b_t = \cF_{\bar{\cA}} = \dd_{Y^{[2]}} \Psi_t \wedge \dd t + \dd t \wedge \partial_t \cA_t + \cF_{\cA_t}
\end{equation*}

\noindent
where we have used Lemma \ref{lemma:gerbecurvaturedecomposition}. Isolating by tensor type this gives:
\begin{equation*}
\delta b_t = F_{\cA_t}\, , \qquad \delta a_t = \partial_t \cA_t - \dd_{Y^{[2]}} \Psi_t
\end{equation*}

\noindent
and thus we conclude. 
\end{proof}

\begin{cor}
There is a natural one-to-one correspondence between curvings on a reducible bundle gerbe $(\bar{\cP},\bar{Y},\bar{\cA})$ with reduction $(\cP,Y,\cA_t,\Psi_t)$ and pairs $(a_t,b_t)$ consisting of a family $b_t$ of curvings on $(\cP,Y,\cA_t)$ and a family of one-forms $a_t$ on $Y$ satisfying Equation \eqref{eq:conditionat}.
\end{cor}

\noindent
Let $\bar{b} = \dd t \wedge a_t + b_t$ be a curving on a reducible bundle gerbe over $\cI\times \Sigma$. Then:
\begin{equation*}
\dd_{\bar{Y}}\bar{b} = \dd t \wedge (\partial_t b_t - \dd_Y a_t) + \dd_Y b_t
\end{equation*}

\noindent
where $\dd_{\bar{Y}}$ is the exterior derivative on $\bar{Y}$ and $\dd_Y$ is the exterior derivative on $Y$. 

\begin{lemma}
The family of two-forms $(\partial_t b_t - \dd_Y a_t) \in \Omega^2(Y)$ satisfies $\delta (\partial_t b_t - \dd_Y a_t) = 0$ for every $t\in \cI$, where $\delta$ denotes the simplicial differential of the simplicial manifold determined by $\pi \colon Y\to M$. 
\end{lemma}

\begin{proof}
We compute:
\begin{equation*}
\delta (\partial_t b_t - \dd_Y a_t) = \partial_t \delta b_t - \dd_{Y^{[2]}} \delta a_t =  \partial_t \cF_{\cA_t}- \dd_{Y^{[2]}} (\partial_t \cA_t - \dd_{Y^{[2]}} \Psi_t) = 0
\end{equation*}

\noindent
where we have used Lemma \ref{lemma:reduciblegerbe}.
\end{proof}

\noindent
Since $\delta$ is exact, by the previous lemma it follows that there exists a unique family of two-forms $c_t \in \Omega^2(\Sigma)$ on $\Sigma$ such that:
\begin{equation*}
\bar{\pi}^{\ast} c_t = \partial_t b_t - \dd_Y a_t
\end{equation*}

\noindent
Since $\dd_{\bar{Y}} \bar{b} = \bar{\pi}^{\ast} H_{\bar{b}}$ for a uniquely determined three-form $\bar{H}_{\bar{b}}\in \Omega^3(M)$, it follows that:
\begin{equation*}
H_{\bar{b}} = \dd t \wedge c_t + H_{b_t}
\end{equation*}

\noindent
where $H_{b_t}$ is the curvature of the family of curvings $b_t$. Therefore, every curving $\bar{b}$ on a reducible bundle gerbe is equivalent to a pair $(a_t , b_t)$ as introduced in Lemma \ref{lemma:reduciblegerbe}, and this pair defines in turn a family of two-forms $c_t$ on $\Sigma$ as described above. 
\begin{definition}
The family of two-forms $c_t$ is the \emph{derived} family of two-forms of the curving $\bar{b}$ on the reducible bundle gerbe $(\bar{\cP},\bar{\cA},\bar{Y})$.
\end{definition}

\noindent
The fact that $H_{\bar{b}}\in \Omega^3(\cI\times \Sigma)$ is closed translates into:
\begin{equation*}
\partial_t H_{b_t} = \dd_{\Sigma} c_t\, , \qquad \dd_{\Sigma} H_{b_t} = 0
\end{equation*}

\noindent
Since $b_t$ is a family of curvings on $(\cP,Y,\cA)$, the second condition holds automatically.  We have:
\begin{equation*}
H_{b_t}  = \int_0^t \dd_{\Sigma} c_{\tau} \dd \tau  + \omega
\end{equation*}

\noindent
for a time-independent integral closed three-form $\omega/2\pi \in\Omega^3(M)$. For dimensional reasons it is convenient to introduce a family of functions $f_{b_t}$ and a family of one-forms $\theta_t$ associated to every family of curvings $b_t$ and metrics $h_t$ through:
\begin{eqnarray*}
f_{b_t} =  \ast_{h_t} H_{b_t}\, , \qquad \theta_t = \ast_{h_t} \int_0^t \dd_{\Sigma} c_{\tau} \dd \tau
\end{eqnarray*}

\noindent
Hence:
\begin{equation*}
f_{b_t} = \ast_{h_t} \dd_{\Sigma} \ast_{h_t} \theta_t + \ast_{h_t}\omega = - \dd_{\Sigma}^{h_t\ast} \theta_t + \ast_{h_t}\omega
\end{equation*}
where $\dd_{\Sigma}^{h_t\ast} \colon \Omega^1(\Sigma) \to C^{\infty}(\Sigma)$ is the formal adjoint of the exterior derivative on $\Sigma$ with respect to $h_t$. We will also refer to $\theta_t$ as the \emph{derived family} of one-forms associated to $\bar{b}$ and the given family of metrics $h_t$.





\section{The supersymmetric NS-NS evolution flow on a bundle gerbe}
\label{sec:NSNSevolutionflow}


In this section we consider the evolution problem posed by the globally hyperbolic supersymmetric solutions of four-dimensional NS-NS supergravity on a fixed pair $(\bar{\cC},\bar{\cX})$. For this, we fix an oriented three-dimensional manifold $\Sigma$ and we take $\bar{\cC} = (\bar{\cP}, \bar{Y}, \bar{\cA})$ to be a reducible bundle gerbe with connective structure on $\cI\times \Sigma$, whose reduction is denoted by $(\cP, Y , \cA_t ,  \Psi_t )$. Similarly, we consider $\bar{\cX}$ to be the pull-black of a principal $\mathbb{Z}$-bundle defined on $\Sigma$ to $\cI\times \Sigma$ via the canonical projection $\cI\times \Sigma\to \Sigma$.


\subsection{Globally hyperbolic NS-NS solutions}


Our first task is to characterize the globally hyperbolic NS-NS solutions on $(\bar{\cC},\bar{\cX})$ in terms of families of solutions on $(\cP, Y , \cA_t  , \Psi_t ,\cX)$.

\begin{definition}
A NS-NS configuration $(g ,\bar{b},\bar{\phi})\in \Conf(\bar{\cC},\bar{\cX})$ is \emph{globally hyperbolic} if $g$ is a globally hyperbolic metric on $\cI\times \Sigma$ of the form $g = - \lambda_t^2 \dd t\otimes \dd t + h_t$ for a family of functions $\lambda_t$ and a family of Riemannian metrics $h_t$ on $\Sigma$.   
\end{definition}

\begin{remark}
By the seminal work of Bernal and S\'anchez \cite{Bernal:2003jb,Bernal:2004gm}, every globally hyperbolic Lorentzian manifold is \emph{globally} isometric to a Lorentzian manifold of the form $(\mathbb{R}\times \Sigma, - \lambda_t^2 \dd t\otimes \dd t + h_t)$. Hence, there is no loss of generality in our previous definition. 
\end{remark}

\noindent
Since $\bar{\cX}$ is the pull-back of a $\mathbb{Z}$-bundle $\cX$ on $\Sigma$, it follows that $\bar{\cX} = \cI \times \cC$, which immediately implies that every equivariant function $\bar{\phi} \colon \bar{\cC}\to \mathbb{R}$ is equivalent to a family of equivariant functions $\phi_t$ on $\cX$. Therefore, globally hyperbolic configurations $(g , \bar{b}, \bar{\phi}) \in \Conf(\bar{\cC},\bar{\cX})$ are in one-to-one correspondence with tuples $( h_t , a_t , b_t , \phi_t , \lambda_t)$ as described above. We will refer to the latter as the \emph{reduction} of the former. 
\begin{definition}
A NS-NS configuration $(g ,\bar{b},\bar{\phi})\in \Conf(\bar{\cC},\bar{\cX})$ is \emph{normal globally hyperbolic} if it is globally hyperbolic with Lorentzian metric $g$ of the form $g = -  \dd t\otimes \dd t + h_t$ for a family of Riemannian metrics $h_t$ on $\Sigma$.   
\end{definition}

\noindent
Every reduced globally hyperbolic NS-NS configuratioon is equivalent to a normal one on a tubular neighborhood of the Cauchy hypersurface in $M = \cI\times \Sigma$. Since for this dissertation we are interested in the \emph{local} behaviour of reduced globally hyperbolic configurations, in the following we will only consider normal ones. A direct calculation proves the following lemma.

\begin{lemma}
\label{lemma:hyperbolicderivatives}
Let $g = -  \dd t \otimes \dd t + h_t$ be a globally hyperbolic metric on $M = \cI\times \Sigma$. Then, for every family of one-forms $\beta_t \in \Omega^1(\Sigma)$ we have:
\begin{eqnarray*}
\nabla^g \dd t = \Theta_t\, , \quad  \nabla^g\beta_t =  \dd t \otimes \partial_{t} \beta_t + \dd t \odot \Theta_t (\beta^{\sharp_{h_t}}_t) + \nabla^{h_t}\beta_t\\
\dd_{\nabla^g} \Theta_t =  \dd_{\nabla^{h_t}} \Theta_t  + \dd t  \wedge \partial_t\Theta_t + \dd t\wedge \Theta_t\circ_{h_t} \Theta_t 
\end{eqnarray*}

\noindent
where $\Theta_t , \Theta_t\circ_{h_t}\Theta_t \in \Gamma(T^{\ast}\Sigma \odot T^{\ast}\Sigma)$ are considered as one-forms taking values on $T^{\ast}\Sigma$ and we have set:
\begin{eqnarray*}
\Theta_t\circ_{h_t}\Theta_t (v_1,v_2) = \langle \Theta_t(v_1) , \Theta_t(v_2)\rangle_{h_t}
\end{eqnarray*} 

\noindent
for every $v_1,v_2 \in T\Sigma$.
\end{lemma}

\noindent
Using the previous lemma we compute the Ricci and scalar curvatures of the \emph{evolving} metric $h_t$ on $\Sigma$ in terms of families of objects on $\Sigma$.
\begin{lemma}
\label{lemma:hyperboliccurvatures}
Let $g = -  \dd t \otimes \dd t + h_t$ be a globally hyperbolic metric on $M = \cI\times \Sigma$. The Ricci and scalar curvatures of $g$ are given by:
\begin{eqnarray*}
& \Ric^g = (\Tr_{h_t}(\partial_t\Theta_t) + \vert \Theta_t \vert_{h_t}^2) \dd t \otimes \dd t  + \Ric^{h_t} + \Tr_{h_t}(\Theta_t) \Theta_t - 2 \Theta_t\circ_{h_t} \Theta_t\\
&  -  \partial_t \Theta_t + \dd t \odot (\dd_{\Sigma} \Tr_{h_t}(\Theta_t) + \nabla^{h_t\ast}\Theta_t)\\
& s^g = s^{h_t} + \Tr_{h_t}(\Theta_t)^2 - 2 \Tr_{h_t}(\partial_t\Theta_t) - 3 \vert \Theta_t\vert_{h_t}^2 
\end{eqnarray*}

\noindent
where $\Theta_t\in \Gamma(T^{\ast}\Sigma\odot T^{\ast}\Sigma)$ is understood as a one-form taking values on $T^{\ast}\Sigma$.
\end{lemma}

\noindent
We decompose now each of the equations of the NS-NS system \eqref{eq:NSNSsystem} evaluated on globally hyperbolic configurations in terms of \emph{flow equations} on $\Sigma$.

\begin{lemma}
Let $(g,\bar{b},\bar{\phi})\in \Conf(\bar{\cC},\bar{\cX})$ be a normal globally hyperbolic NS-NS configuration with reduction $( h_t , a_t , b_t , \phi_t)$. Then:
\begin{eqnarray*}
\varphi_{\bar{\phi}} = \partial_t \phi_t \dd t + \varphi_{\phi_t}
\end{eqnarray*}
	
\noindent
where $\varphi_{\phi_t}$ is the curvature of $\phi_t$. Furthermore:
\begin{equation*}
\varphi_{\phi_t} = \dd_{\Sigma}\phi_t + \rho
\end{equation*}

\noindent
in terms of a time-independent integral closed one-form $\rho \in \Omega^1 (\Sigma)$.  
\end{lemma}

\begin{proof}
The exterior derivative $\dd_{\cC}\bar{\phi}$ of $\bar{\phi}$ as a function on $\bar{\cC}$ expands as follows:
\begin{equation*}
\dd_{\bar{\cC}}\bar{\phi} = \partial_t \phi_t \dd t + \dd_{\cC}\phi_t \in \Omega^1(\cI\times \cC)
\end{equation*}

\noindent
This implies that:
\begin{equation*}
	\varphi_{\bar{\phi}} =  \partial_t \phi_t \dd t + \varphi_{\phi_t}
\end{equation*}

\noindent
where $\varphi_{\phi_t}$ is the family of curvatures of $\phi_t$. Furthermore, since $\varphi_{\bar{\phi}}$ is closed, we have:
\begin{equation*}
0 =	\dd \varphi_{\bar{\phi}} =  \partial_t \dd_{\Sigma}\phi_t \wedge \dd t + \dd t \wedge \partial_t\varphi_{\phi_t} + \dd_{\Sigma}\varphi_{\phi_t}
\end{equation*}

\noindent
and hence we conclude. 
\end{proof}

\begin{lemma}
\label{lemma:Einsteingh}
A normal globally hyperbolic NS-NS configuration $(g,\bar{b},\bar{\phi})\in \Conf(\bar{\cC},\bar{\cX})$ satisfies the Einstein equation of the NS-NS system \eqref{eq:NSNSsystem} if and only if its reduction $( h_t , a_t , b_t , \phi_t)$ satisfies the following differential system:
\begin{eqnarray*}
& \partial_t^2 \phi_t + \Tr_{h_t}(\partial_t\Theta_t) + \vert \Theta_t \vert_{h_t}^2 - \frac{1}{2} \vert\theta_t\vert_{h_t}^2 = 0\\
& \dd_{\Sigma} \Tr_{h_t}(\Theta_t) + \nabla^{h_t\ast}\Theta_t + \partial_t \varphi_{\phi_t} +\Theta_t (\varphi_{\phi_t}^{\sharp_{h_t}}) - \frac{1}{2}  f_{b_t} \theta_t  = 0\\
& \Ric^{h_t} + (\Tr_{h_t}(\Theta_t) + \partial_t \phi_t) \Theta_t - 2 \Theta_t\circ_{h_t} \Theta_t -  \partial_t \Theta_t + \nabla^{h_t} \varphi_{\phi_t}  - \frac{1}{2} \theta_t \otimes \theta_t + \frac{1}{2} (\vert \theta_t\vert_{h_t}^2 - f_{b_t}^2) h_t  = 0
\end{eqnarray*}

\noindent
where $\Theta_t$ is the fundamental form of $(\left\{ t\right\} \times \Sigma, h_t)$, $f_{b_t}\in C^{\infty}(\Sigma)$ is the Hodge dual of $H_{b_t}$ and $\theta_t$ is the derived family of one-forms of $( h_t , a_t , b_t , \phi_t)$. 
\end{lemma}

\begin{remark}
Note that equations:
\begin{equation*}
\varphi_{\phi_t} = \dd_{\Sigma}\phi_t + \rho\, , \qquad f_{b_t} = - \dd_{\Sigma}^{h_t\ast} \theta_t + \ast_{h_t}\omega
\end{equation*}
need not be included in the previous lemma since they are a consequence of $( h_t , a_t , b_t , \phi_t)$ being the \emph{reduction} of $(g,\bar{b},\bar{\phi})\in \Conf(\bar{\cC},\bar{\cX})$. 
\end{remark}

\begin{proof}
Let $(g,\bar{b},\bar{\phi})$ be a globally hyperbolic configuration with reduction $( h_t , a_t , b_t , \phi_t)$. We evaluate the Einstein equation of the NS-NS system, namely the first equation in \eqref{eq:NSNSsystem}, on the following configuration:
\begin{eqnarray*}
& g = - \dd t\otimes \dd t + h_t \in \Gamma(T^{\ast}M\odot T^{\ast}M) \, , \quad \bar{b} = \dd t\wedge a_t + b_t \in \Omega^2(\cI\times Y)\\
& \bar{\phi}(t,p) = \phi_t(p)\in C^{\infty}(\Sigma)
\end{eqnarray*}

\noindent
for every $(t,p)\in \cI\times \Sigma$. By Lemma \ref{lemma:hyperboliccurvatures}, the Ricci tensor of $g$ is given by:
\begin{eqnarray*}
& \Ric^{g} (\partial_t , \partial_t) =  \Tr_{h_t}(\partial_t\Theta_t) + \vert \Theta_t\vert^2_{h_t} \,  ,\quad \Ric^g(\partial_t)\vert_{T\Sigma} = \dd_{\Sigma} \Tr_{h_t}(\Theta_t) + \nabla^{h_t\ast}\Theta_t\\
& \Ric^{g}\vert_{T\Sigma\times T\Sigma} = \Ric^{h_t} + \Tr_{h_t}(\Theta_t) \Theta_t - 2 \Theta_t\circ_{h_t} \Theta_t -  \partial_t \Theta_t 
\end{eqnarray*}

\noindent
where $\Theta_t$ is the second fundamental form of the embedding $\left\{ t\right\}\times \Sigma \hookrightarrow \cI \times \Sigma$ as defined in Equation \eqref{eq:secondfundamentalform} with respect to the metric $g$. By the discussion of Section \ref{sec:reductionbundlegerbe} we have:
\begin{equation*}
H_{\bar{b}} = \dd t \wedge c_t + H_{b_t} = \dd t \wedge \ast_{h_t} \theta_t + H_{b_t}
\end{equation*}

\noindent
from which we compute:
\begin{eqnarray*}
& (H_{\bar{b}}\circ_{g} H_{\bar{b}}) (\partial_t , \partial_t) =  \vert \theta_t \vert_{g}^2 = \vert \theta_t \vert_{h_t}^2\\
& (H_{\bar{b}}\circ_{g} H_{\bar{b}}) (\partial_t , v) = \langle \ast_{h_t}\theta_t, - \dd t\wedge \ast_{h_t}(\theta_t \wedge v) + f_{b_t} \ast_{h_t} v  \rangle_{g} =  f_{b_t} \theta_t(v)\, , \quad v\in T\Sigma
\end{eqnarray*}

\noindent
where we have defined the function $f_{b_t} \in C^{\infty}(\Sigma)$ associated to every family of curvings $b_t$ and Riemannian metrics $h_t$ by $H_{b_t} := f_{b_t} \nu_{h_t}$. Furthermore:
\begin{eqnarray*}
& (H_{\bar{b}}\circ_{g} H_{\bar{b}})(v_1,v_2) = \langle -\dd t \wedge v_1\lrcorner_{h_t}\ast_{h_t} \theta + f_{b_t} \ast_{h_t} v_1 , -\dd t \wedge v_2\lrcorner_{h_t}\ast_{h_t} \theta + f_{b_t} \ast_{h_t} v_2 \rangle_{h_t} \\
& = -  \langle \ast_{h_t}(\theta_t \wedge v_1) , \ast_{h_t}(\theta_t \wedge v_2) \rangle_{h_t} + f_{b_t}^2 h_t(v_1,v_2) = - \vert \theta_t\vert_{h_t}^2 h_t(v_1,v_2) + \theta_t(v_1) \theta_t(v_2) + f_{b_t}^2 h_t(v_1,v_2)
\end{eqnarray*}

\noindent
for every $v_1 ,v_2 \in T\Sigma$. On the other hand, using Lemma \ref{lemma:hyperbolicderivatives} we obtain:
\begin{eqnarray*}
& \nabla^{g}_{\partial_t}\varphi_{\bar{\phi}} =  \partial^2_t \phi_t  \dd t + \partial_t \varphi_{\phi_t} +\Theta_t (\varphi_{\phi_t}^{\sharp_{h_t}})\\
& \nabla^{g}_v \varphi_{\bar{\phi}} = (v(\partial_t\varphi_{\phi_t})   + \Theta_t (\varphi^{\sharp_{h_t}}_{\phi_t},v) ) \dd t + \Theta_t (v) \partial_t \phi_t + \nabla^{h_t}_v \varphi_{\phi_t}
\end{eqnarray*}

\noindent
for every $v \in T\Sigma$. With these provisos in mind, we compute:
\begin{eqnarray*}
\Ric^g(\partial_t , \partial_t) + (\nabla^g_{\partial_t}\varphi_{\bar{\phi}})(\partial_t) - \frac{1}{2} (H_{\bar{b}}\circ_{g} H_{\bar{b}}) (\partial_t , \partial_t) = \partial_t^2 \phi_t + \Tr_{h_t}(\partial_t\Theta_t) + \vert \Theta_t \vert_{h_t}^2 - \frac{1}{2} \vert\theta_t\vert_{h_t}^2 = 0
\end{eqnarray*}

\noindent
for the \emph{time-time} component of the Einstein equation in \eqref{eq:NSNSsystem}. For the mixed \emph{time-space} components we have:
\begin{eqnarray*}
& (\Ric^g(\partial_t) + (\nabla^g_{\partial_t}\varphi_{\bar{\phi}})  - \frac{1}{2} (H_{\bar{b}}\circ_{g} H_{\bar{b}}) (\partial_t ))\vert_{T\Sigma} = \\
& = \dd_{\Sigma} \Tr_{h_t}(\Theta_t) + \nabla^{h_t\ast}\Theta_t + \partial_t \varphi_{\phi_t} +\Theta_t (\varphi_{\phi_t}^{\sharp_{h_t}}) - \frac{1}{2}  f_{b_t} \theta_t  = 0
\end{eqnarray*}

\noindent
Finally, for the \emph{space-space} components of the Einstein equation we have:
\begin{eqnarray*}
& 0 = (\Ric^g  + (\nabla^g \varphi_{\bar{\phi}})  - \frac{1}{2} H_{\bar{b}}\circ_{g} H_{\bar{b}})\vert_{T\Sigma\times T\Sigma} = \\
& = \Ric^{h_t} + (\Tr_{h_t}(\Theta_t) + \partial_t \phi_t) \Theta_t - 2 \Theta_t\circ_{h_t} \Theta_t -  \partial_t \Theta_t + \nabla^{h_t} \varphi_{\phi_t}  - \frac{1}{2} \theta_t \otimes \theta_t + \frac{1}{2} (\vert \theta_t\vert_{h_t}^2 - f_{b_t}^2) h_t 
\end{eqnarray*}

\noindent
and hence we conclude. 
\end{proof}

\begin{lemma}
\label{lemma:Maxwellgh}
A normal globally hyperbolic NS-NS configuration $(g,\bar{b},\bar{\phi})\in \Conf(\bar{\cC},\bar{\cX})$ satisfies the Maxwell equation of the NS-NS system \eqref{eq:NSNSsystem} if and only if its reduction $( h_t , a_t , b_t , \phi_t)$ satisfies the following differential system:
\begin{equation*}
\dd_{\Sigma} f_{b_t} - \partial_t \theta_t + \partial_t \phi_t \theta_t - f_{b_t} \varphi_{\phi_t} = 0 \, , \qquad \dd_{\Sigma} \theta_t = \varphi_{\phi_t} \wedge \theta_t
\end{equation*}
	
\noindent
where $\frf_{b_t}\in C^{\infty}(\Sigma)$ is the Hodge dual of $H_{b_t}$ and $\theta_t$ is the derived family of one-forms of $( h_t , a_t , b_t , \phi_t)$. 
\end{lemma}

\begin{proof}
We consider the Maxwell equation as written in \eqref{eq:NSNSsystem4d}, namely:
\begin{equation*}
\dd \alpha_{\bar{b}} = \varphi_{\bar{\phi}} \wedge \alpha_{\bar{b}}
\end{equation*}

\noindent
We have:
\begin{equation*}
\alpha_{\bar{b}} = \ast_g (\dd t \wedge \ast_{h_t}\theta_t + f_{b_t} \nu_{h_t}) = - (\theta_t + f_{b_t} \dd t)
\end{equation*}

\noindent
where we have used equations \eqref{eq:hodgehyperbolic}. Hence, we have:
\begin{eqnarray*}
& \dd \alpha_{\bar{b}} + \alpha_{\bar{b}}\wedge \varphi_{\bar{\phi}} = -\dd (\theta_t + f_{b_t} \dd t) - (\theta_t + f_{b_t} \dd t) \wedge (\partial_t \phi_t \dd t + \varphi_{\phi_t})\\
& = \dd t\wedge (\dd_{\Sigma} f_{b_t} - \partial_t \theta_t) - \dd_{\Sigma} \theta_t + \dd t \wedge ( \partial_t \phi_t \theta_t - f_{b_t} \varphi_{\phi_t})   - \theta_t \wedge \varphi_{\phi_t} = 0
\end{eqnarray*}

\noindent
Isolating by tensor type we conclude.
\end{proof}

\begin{lemma}
\label{lemma:Dilatongh}
A normal globally hyperbolic NS-NS configuration $(g,\bar{b},\bar{\phi})\in \Conf(\bar{\cC},\bar{\cX})$ satisfies the dilaton equation of the NS-NS system \eqref{eq:NSNSsystem} if and only if its reduction $( h_t , a_t , b_t , \phi_t)$ satisfies the following differential system:
\begin{equation*}
\partial_t^2\phi_t + \nabla^{h_t\ast}\varphi_{\phi_t}  - (\Tr_{h_t}(\Theta_t) + \partial_t \phi_t) \partial_t \phi_t - f_{b_t}^2 + \vert \theta_t \vert_{h_t}^2 + \vert \varphi_{\phi_t} \vert_{h_t}^2 = 0
\end{equation*}
 
\noindent
where $\Theta_t$ is the fundamental form of $(\left\{ t\right\} \times \Sigma, h_t)$, $\frf_{b_t}\in C^{\infty}(\Sigma)$ is the Hodge dual of $H_{b_t}$ and $\theta_t$ is the derived family of one-forms of $( h_t , a_t , b_t , \phi_t)$. 
\end{lemma}

\begin{proof}
We compute:
\begin{eqnarray*}
\vert \varphi_{\phi}\vert_g^2 + \vert\alpha_b\vert^2_g =  \vert \theta_t \vert_{h_t}^2 + \vert \varphi_{\phi_t} \vert_{h_t}^2 - (\partial_t \phi_t)^2 - f_{b_t}^2 
\end{eqnarray*}

\noindent
On the other hand we have:
\begin{equation*}
\nabla^{g\ast}\varphi_{\bar{\phi}}  = \nabla^{g\ast}(\partial_t \phi_t + \varphi_{\phi_t}) =  \partial_t^2\phi_t - \Tr_{h_t}(\Theta_t) \partial_t \phi_t + \nabla^{h_t\ast}\varphi_{\phi_t}
\end{equation*}

\noindent
and hence  we conclude.
\end{proof}

\noindent
The previous lemmata allows to characterize globally hyperbolic NS-NS solutions on a reducible pair $(\bar{\cP},\bar{\cA},\bar{Y})$ in terms of constrained evolution flows for tuples $( h_t , a_t , b_t , \phi_t)$ on the reduction $(\cP,Y,\cA_t,\Psi_t,\cX)$ of $(\bar{\cP},\bar{\cA},\bar{Y})$.
\begin{prop}
\label{prop:normalhyperbolicflow}
A normal globally hyperbolic NS-NS configuration $(g,\bar{b},\bar{\phi})\in \Conf(\bar{\cC},\bar{\cX})$ satisfies the NS-NS system \eqref{eq:NSNSsystem} if and only if its reduction $( h_t , a_t , b_t , \phi_t)$ satisfies the following system of evolution equations:
\begin{eqnarray*}
& \partial_t \Theta_t = \Ric^{h_t} + (\Tr_{h_t}(\Theta_t) + \partial_t \phi_t) \Theta_t - 2 \Theta_t\circ_{h_t} \Theta_t   + \nabla^{h_t} \varphi_{\phi_t}  - \frac{1}{2} \theta_t \otimes \theta_t + \frac{1}{2} (\vert \theta_t\vert_{h_t}^2 - f_{b_t}^2) h_t \nonumber \\
& \partial_t \theta_t = \dd_{\Sigma} f_{b_t}  + \partial_t \phi_t \theta_t - f_{b_t} \varphi_{\phi_t} \label{eq:normalevolution} \\
& \partial_t^2\phi_t + \nabla^{h_t\ast}\varphi_{\phi_t}  = (\Tr_{h_t}(\Theta_t) + \partial_t \phi_t) \partial_t \phi_t + f_{b_t}^2 - \vert \theta_t \vert_{h_t}^2 - \vert \varphi_{\phi_t} \vert_{h_t}^2  \nonumber
\end{eqnarray*}
together with the following system of time-dependent constraints:
\begin{eqnarray*}
& s^{h_t}  + (\Tr_{h_t}(\Theta_t) +  \partial_t \phi_t)^2 -   \vert\Theta_t\vert^2_{h_t} - \vert \varphi_{\phi_t} \vert_{h_t}^2  - \frac{1}{2}( f_{b_t}^2 + \vert \theta_t \vert^2_{h_t} )  - 2 \nabla^{h_t\ast} \varphi_{\phi_t}  = 0\\
& \dd_{\Sigma} \theta_t = \varphi_{\phi_t} \wedge \theta_t \, , \quad \dd_{\Sigma} \Tr_{h_t}(\Theta_t) + \nabla^{h_t\ast}\Theta_t + \partial_t \varphi_{\phi_t} +\Theta_t (\varphi_{\phi_t}^{\sharp_{h_t}}) - \frac{1}{2}  \frf_{b_t} \theta_t  = 0
\end{eqnarray*}

\noindent
where $\Theta_t$ is the fundamental form of $(\left\{ t\right\} \times \Sigma, h_t)$, $f_{b_t}\in C^{\infty}(\Sigma)$ is the Hodge dual of $H_{b_t}$ and $\theta_t$ is the derived family of one-forms of $( h_t , a_t , b_t , \phi_t)$. 
\end{prop}

\begin{proof}
Lemmas \eqref{lemma:Einsteingh}, \eqref{lemma:Maxwellgh} and \eqref{lemma:Dilatongh} immediately imply that a normal globally hyperbolic NS-NS configuration $(g,\bar{b},\bar{\phi})\in \Conf(\bar{\cC},\bar{\cX})$ satisfies the NS-NS system \eqref{eq:NSNSsystem} if and only if its reduction $( h_t , a_t , b_t , \phi_t)$ on $(\cP,Y,\cA_t,\Psi_t,\cX)$ satisfies the following system of evolution equations:
\begin{eqnarray*}
& \partial_t \Theta_t = \Ric^{h_t} + (\Tr_{h_t}(\Theta_t) + \partial_t \phi_t) \Theta_t - 2 \Theta_t\circ_{h_t} \Theta_t   + \nabla^{h_t} \varphi_{\phi_t}  - \frac{1}{2} \theta_t \otimes \theta_t + \frac{1}{2} (\vert \theta_t\vert_{h_t}^2 - f_{b_t}^2) h_t  \\
& \partial_t \theta_t = \dd_{\Sigma} f_{b_t}  + \partial_t \phi_t \theta_t - f_{b_t} \varphi_{\phi_t} \\
& \partial_t^2\phi_t + \nabla^{h_t\ast}\varphi_{\phi_t}  = (\Tr_{h_t}(\Theta_t) + \partial_t \phi_t) \partial_t \phi_t + f_{b_t}^2 - \vert \theta_t \vert_{h_t}^2 - \vert \varphi_{\phi_t} \vert_{h_t}^2  
\end{eqnarray*}
together with the following system of time-dependent constraints:
\begin{eqnarray*}
& \s^{h_t}  + (\Tr_{h_t}(\Theta_t)^2 +  \partial_t \phi_t)^2 -   \vert\Theta_t\vert^2_{h_t} - \vert \varphi_{\phi_t} \vert_{h_t}^2  - \frac{1}{2}( f_{b_t}^2 + \vert \theta_t \vert^2_{h_t} )  - 2 \nabla^{h_t\ast} \varphi_{\phi_t}  = 0\\
& \dd_{\Sigma} \Tr_{h_t}(\Theta_t) + \nabla^{h_t\ast}\Theta_t + \partial_t \varphi_{\phi_t} +\Theta_t (\varphi_{\phi_t}^{\sharp_{h_t}}) - \frac{1}{2}  f_{b_t} \theta_t  = 0 
\end{eqnarray*}

\noindent
Taking the trace of the first evolution equation above and combining it with the third evolution equation we obtain:
\begin{eqnarray*}
\partial_t^2 \phi_t + \Tr_{h_t}(\partial_t\Theta_t) = s^{h_t} - 2 \nabla^{h_t\ast} \varphi_{\phi_t}  + (\Tr_{h_t}(\Theta_t)^2 +  \partial_t \phi_t)^2 - 2 \vert\Theta_t\vert^2_{h_t} - \vert \varphi_{\phi_t} \vert_{h_t}^2  - \frac{1}{2} f_{b_t}^2 
\end{eqnarray*}

\noindent
Together with the first time-dependent constraint above, this implies:
\begin{eqnarray*}
s^{h_t}  + (\Tr_{h_t}(\Theta_t)^2 +  \partial_t \phi_t)^2 -   \vert\Theta_t\vert^2_{h_t} - \vert \varphi_{\phi_t} \vert_{h_t}^2  - \frac{1}{2}( f_{b_t}^2 + \vert \theta_t \vert^2_{h_t} )  - 2 \nabla^{h_t\ast} \varphi_{\phi_t}  = 0
\end{eqnarray*} 

\noindent
and hence we conclude. 
\end{proof}

\noindent
The previous proposition presents the evolution and constraint equations of the NS-NS system in the string frame and form the starting point for the study of the evolution problem posed by NS-NS supergravity. To the best of our knowledge we have not even seen this constrained evolution system written explicitly in the literature, especially not in the global geometric context we are considering. By the seminal results of Oskar Schiller, to be presented in his upcoming doctoral dissertation at Hamburg University, the Cauchy problem for the NS-NS system is well-posed.


\subsection{The supersymmetric NS-NS evolution flow}


In this section we consider the evolution flow defined by the globally hyperbolic supersymmetric configurations on a reducible triple $(\bar{\cP},\bar{\cC},\bar{\cA})$ with reduction $\cC = (\cP, Y , \cA_t  , \Psi_t )$. Our first result builds on the theory developed so far for skew-torsion parallelisms and supersymmetric NS-NS solutions and characterizes normal globally hyperbolic supersymmetric solutions in terms of solutions of a constrained evolution system, as explained in the following.

\begin{prop}
\label{prop:susyevol}
A normal globally hyperbolic NS-NS configuration $(g,\bar{b},\bar{\phi})\in \Conf(\bar{\cC},\bar{\cX})$ is supersymmetric if and only if there exists a family of functions $\fra_t$ and a family of coframes $\fre_t$ on $\Sigma$ such that the following differential system:
\begin{eqnarray}
& \partial_t \fra_t = 0 \, , \quad \partial_t e^t_i + \Theta_t(e^t_i)  + \frac{1}{2} c_t (e^t_i)= 0 \, , \quad i = u, l , n\nonumber \\
&   \dd_{\Sigma} e^t_i = \Theta_t(e^t_i)\wedge e^t_u + \frac{1}{2} f_{b_t}  \ast_{h_t} e^t_i - \frac{1}{2} \delta_{ui}c_t\, , \quad i = u, l , n \label{eq:susyevol} \\
& e^{-\fra_t} \dd_{\Sigma} (e^{\fra_t} e^t_u) =  f_{b_t} \ast_{h_t} e^t_u - c_t\, , \quad \varphi_{\phi_t} = \partial_t\phi_t e^t_u  + c_t(e^t_u)\, , \quad   f_{b_t} = \ast_{h_t} (c_t\wedge e^t_u)\nonumber
\end{eqnarray}

\noindent
is satisfied, where $( h_t , a_t , b_t , \phi_t)$  is the reduction of $(g,\bar{b},\bar{\phi})$.
\end{prop}

\begin{remark}
In the previous proposition we are using the Kronecker delta $\delta_{ui}$ for the indices $i = u, l, n$. Hence, $\delta_{uu} = 1$ and $\delta_{ul} = \delta_{un} = 0$. 
\end{remark}

\begin{proof}
By Proposition \ref{prop:Thetaevolutionsystem} together with the discussion in Section \ref{sec:ReductionGloballyHyperbolic}, a given globally hyperbolic configuration $(g,\bar{b},\bar{\phi}) \in \Conf(\bar{\cC},\bar{\cX})$ is supersymmetric if and only if there exists a globally hyperbolic skew-torsion parallelism $[u,v,l,n]$ such that Lemma \ref{lemma:dilatinoeq} holds. By Proposition \ref{prop:Thetaevolutionsystem}, $[u,v,l,n]$ is globally hyperbolic and skew-torsion if and only if its globally hyperbolic reduction $(\fra_t,\fre_t)$ satisfies equations \eqref{eq:evolutionequationsIInolambda} and \eqref{eq:constraintequationsIInolambda} with:
\begin{eqnarray*}
\varphi_{\bar{\phi}} = \partial_t\phi_t \dd t + \varphi_{\phi_t} \, , \qquad \alpha_{\bar{b}} = - (\theta_t + f_{b_t} \dd t)
\end{eqnarray*}

\noindent
By Lemma \ref{lemma:dilatinoeq}, we impose:
\begin{eqnarray*}
\varphi_{\bar{\phi}}(u) = e^{\fra_t}( \partial_t \phi_t \dd t + \varphi_{\phi_t})(-\partial_t + (e^t_u)^{\sharp_{h_t}})= e^{\fra_t}( \varphi_{\phi_t}(e^t_u) -\partial_t \phi_t )  = 0\\
\alpha_{\bar{b}}(u) = - e^{\fra_t}(f_{b_t} \dd t + \theta_t)(-\partial_t + (e^t_u)^{\sharp_{h_t}})= e^{\fra_t}(f_{b_t} -\theta_t(e^t_u))  = 0
\end{eqnarray*}

\noindent
which is equivalent to:
\begin{equation*}
\partial_t \phi_t = \varphi_{\phi_t}(e^t_u) \, , \qquad f_{b_t} =\theta_t(e^t_u)
\end{equation*}

\noindent
Furthermore:
\begin{eqnarray*}
& \varphi_{\bar{\phi}} \wedge u -  \ast_g (\alpha_{\bar{b}} \wedge u) = e^{\fra_t}(\partial_t\phi_t \dd t + \varphi_{\phi_t}) \wedge (\dd t + e^t_u) + e^{\fra_t}\ast_g((f_{b_t} \dd t + \theta_t) \wedge (\dd t + e^t_u))\\
& = e^{\fra_t}(\dd t \wedge (\partial_t\phi_t e^t_u -\varphi_{\phi_t} + \ast_{h_t} (\theta_t\wedge e^t_u)) + \varphi_{\phi_t}\wedge e^t_u - f_{b_t} \ast_{h_t} e^t_u + \ast_{h_t} \theta_t)
\end{eqnarray*}

\noindent
This implies:
\begin{equation*}
\partial_t\phi_t e^t_u -\varphi_{\phi_t} + \ast_{h_t} (\theta_t\wedge e^t_u) = 0\, , \qquad 	\varphi_{\phi_t}\wedge e^t_u - f_{b_t} \ast_{h_t} e^t_u + \ast_{h_t} \theta_t = 0
\end{equation*}

\noindent
or, equivalently:
\begin{equation*}
\varphi_{\phi_t} = \partial_t\phi_t e^t_u  + \ast_{h_t} (\theta_t\wedge e^t_u)\, , \qquad \theta_t = f_{b_t} e^t_u + \theta_t(e^t_l) e^t_l  + \theta_t(e^t_n) e^t_n 
\end{equation*}

\noindent
Hence, there exists a family of one-forms $w_{b_t} \in \Gamma(e^t_u)^{\perp_{h_t}} $ orthogonal to $e^t_u$ such that:
\begin{equation}
\label{eq:wbt}
\varphi_{\phi_t} = \partial_t \phi_t e^t_u + \ast_{q_t} w_{b_t} \, , \qquad \theta_t = f_{b_t} e^t_u + w_{b_t}
\end{equation}

\noindent
where $\ast_{q_t}$ is the Hodge dual determined on the distribution $(e^t_u)^{\perp_{h_t}}\subset T^{\ast}\Sigma$ by the metric $q_t = e^t_l \otimes e^t_l + e^t_n \otimes e^t_n$. Similarly, we compute:
\begin{eqnarray*}
& \alpha_{\bar{b}}(l) \ast_g u - \varphi_{\bar{\phi}}\wedge l\wedge u =  e^{\fra_t} (\theta_t(e^t_l) \nu_{h_t} - \varphi_{\phi_t}\wedge e^t_l\wedge e^t_u) \\
& + e^{\fra_t} \dd t \wedge (\theta_t(e^t_l) e^t_l \wedge e^t_n  - \partial_t \phi_t e^t_l \wedge e^t_u - \varphi_{\phi_t} \wedge e^t_l) = 0 
\end{eqnarray*}

\noindent
which is thus automatically satisfied by \eqref{eq:wbt}. Analogously, equation $\varphi_{\phi}(l) = \ast_g (u\wedge l \wedge \alpha_{\bar{b}})$ is automatically solved by \eqref{eq:wbt}. Furthermore, we notice that:
\begin{eqnarray*}
& \ast_{h_t} (\alpha^{\perp}_t   - \alpha^o_t e^t_u) - \frac{1}{2} e^t_u \wedge \ast_{h_t} (\alpha^{\perp}_t \wedge e^t_u) = f_{b_t} \ast_{h_t} e^t_u - c_t + \frac{1}{2} e^t_u \wedge \ast_{h_t} (\theta_t \wedge e^t_u)\\
& = (f_{b_t} - \theta_t(e^t_u)) \ast_{h_t} e^t_u - \theta_t(e^t_l) \ast_{h_t} e^t_l- \theta_t(e^t_n) \ast_{h_t} e^t_n + \frac{1}{2} e^t_u \wedge \ast_{h_t} (\theta_t \wedge e^t_u) \\
& = -\frac{1}{2} (\theta_t(e^t_l) \ast_{h_t} e^t_l + \theta_t(e^t_n) \ast_{h_t} e^t_n) = -\frac{1}{2} e^t_u \wedge c_t (e^t_u) = \frac{1}{2}   (f_{b_t} \ast_{h_t}e^t_u - c_t) 
\end{eqnarray*}

\noindent
Substituting the previous expressions and relations into \eqref{eq:evolutionequationsIInolambda} and \eqref{eq:constraintequationsIInolambda} and massaging the equations we obtain the differential system \eqref{eq:susyevol} in the statement of the lemma. These equations guarantee that the globally hyperbolic reduction of $[u,v,l,n]$ is skew-torsion with torsion given by the curvature of a curving on a reducible bundle gerbe. The converse follows by tracing back the previous steps and hence we conclude. 
\end{proof}

\noindent
By the previous proposition, the family of functions $\fra_t \in C^{\infty}(\Sigma)$ that conforms the reduced isotropic parallelism, not to be confused with the family of one-forms $a_t \in \Omega^1(Y)$ associated to the globally hyperbolic reduction of a curving, is actually constant \emph{in time}, namely $\partial_t a_t = 0$, so it reduces to a function $\fra \in C^{\infty}(\Sigma)$ on $\Sigma$. This function can be decoupled from the evolution system as the following result shows.

\begin{cor}
\label{cor:susyevolcor}
A normal globally hyperbolic NS-NS configuration $(g,\bar{b},\bar{\phi})\in \Conf(\bar{\cC},\bar{\cX})$ is supersymmetric if and only if there exists a family of coframes $\fre_t$ on $\Sigma$ such that the following differential system:
\begin{eqnarray}
& \partial_t \fra_t = 0 \, , \qquad \partial_t e^t_i + \Theta_t(e^t_i)  + \frac{1}{2} c_t (e^t_i)= 0 \, , \quad i = u, l , n\nonumber \\
&   \dd_{\Sigma} e^t_i = \Theta_t(e^t_i)\wedge e^t_u + \frac{1}{2} f_{b_t}  \ast_{h_t} e^t_i - \frac{1}{2} \delta_{ui}c_t\, , \quad i = u, l , n\label{eq:susyevolcor} \\
& \varphi_{\phi_t} = \partial_t\phi_t e^t_u  + c_t(e^t_u)\, , \quad   f_{b_t} = \ast_{h_t} (c_t\wedge e^t_u)\nonumber
\end{eqnarray}

\noindent
together with the following conditions:
\begin{eqnarray*}
0 = [\Theta_t (e^t_u) - \frac{1}{2} c_t(e^t_u)]\in H^1(\Sigma,\mathbb{R})\, , \qquad  \partial_t (\Theta_t (e^t_u)) = \frac{1}{2} \partial_t( c_t(e^t_u))
\end{eqnarray*}

\noindent
are satisfied, where $( h_t , a_t , b_t , \phi_t)$  is the reduction of $(g,\bar{b},\bar{\phi})$.
\end{cor}

\noindent
In practical situations, we are interested in the existence of normal globally hyperbolic configurations without necessarily fixing the underlying reducible bundle gerbe and principal $\mathbb{Z}$ bundle. 

\begin{prop}
\label{prop:susyevolII}
A four-manifold of the form $M = \cI \times \Sigma$ admits a supersymmetric normal globally hyperbolic NS-NS configuration with Cauchy hypersurface $\Sigma$ if and only if there exists a tuple $(\fre_t,c_t, f_t,\phi_t)$ on $\Sigma$ satisfying the following differential system:
\begin{eqnarray*}
& \partial_t( \Theta_t (e^t_u)) = \frac{1}{2} \partial_t (c_t(e^t_u))\, , \qquad  \partial_t e^t_i + \Theta_t(e^t_i)  + \frac{1}{2} c_t (e^t_i)= 0 \, ,  \qquad i = u, l , n\nonumber \\
& f_{t} = \ast_{h_t}(c_t\wedge e^t_u)\, , \quad \dd_{\Sigma} e^t_i = \Theta_t(e^t_i)\wedge e^t_u + \frac{1}{2} f_{b_t}  \ast_{h_t}  e^t_i - \frac{1}{2} \delta_{ui} c_t \, ,  \quad i = u, l , n   
\end{eqnarray*}
	
\noindent
together with the following cohomological conditions:
\begin{eqnarray*}
[\partial_t\phi_t e^t_u  + c_t(e^t_u)]\in H^1(\Sigma,\mathbb{Z}) \, , \quad[\frac{f_t}{2\pi} \nu_{h_t}]\in H^3(\Sigma,\mathbb{Z})\, , \quad 0 = [\Theta_t (e^t_u) - \frac{1}{2} c_t(e^t_u)]\in H^1(\Sigma,\mathbb{R})   
\end{eqnarray*}

\noindent
for every $t\in \cI$. 
\end{prop}

\noindent
We will refer to the differential system \eqref{eq:susyevol} as occurring in Proposition \ref{prop:susyevol}, or equivalently in Corollary \ref{cor:susyevolcor} or Proposition \ref{prop:susyevolII}, as the \emph{normal NS-NS supersymmetry flow}, which is an evolution flow for families of the form $(\fre_t,b_t,\phi_t)$. Given such a family, we obtain a canonical candidate of NS-NS evolution flow given simply by $(h_{\fre_t},b_t,\phi_t)$, where $h_{\fre_t}$ is the family of Riemannian metrics on $\Sigma$ given by $h_{\fre_t} = e^t_u\otimes e^t_u + e^t_l\otimes e^t_l  + e^t_n\otimes e^t_n$. 
 
\begin{remark}
\label{remark:covariantevolution}
From equations \eqref{eq:susyevol} or \eqref{eq:susyevolcor} it immediately follows that the covariant derivative of the family of coframes $\fre_t$ belonging to a NS-NS supersymmetry flow $(\fre_t,b_t,c_t)$ is given by:
\begin{eqnarray*}
\nabla^{h_t} e^t_i + \delta_{ui}( \Theta_t + \frac{1}{4} c_t)= \Theta_t(e^t_i)\otimes e^t_u + \frac{1}{4} f_{b_t} \ast_{h_t} e^t_i  - \frac{1}{4} e^t_u\odot c_t(e^t_i)
\end{eqnarray*}

\noindent
for $i = u, l, n$. The previous covariant derivatives can be equivalently written as follows: 
\begin{eqnarray*}
& \nabla^{h_t} e^t_u + \Theta_t = \Theta_t(e^t_u)\otimes e^t_u - \frac{1}{2} e^t_u\otimes c_t(e^t_u) \\
& \nabla^{h_t} e^t_l = \Theta_t(e^t_l)\otimes e^t_u - \frac{1}{2} f_{b_t} e_u^t \otimes e^t_n + \frac{1}{2} c_t (e^t_u , e^t_l) e^t_u \otimes e^t_u\\
& \nabla^{h_t} e^t_n  = \Theta_t(e^t_n)\otimes e^t_u + \frac{1}{2} f_{b_t} e_u^t \otimes e^t_l + \frac{1}{2} c_t (e^t_u , e^t_n) e^t_u \otimes e^t_u
\end{eqnarray*}
	
\noindent
From this we immediately obtain:
\begin{eqnarray*}
\nabla^{h_t\ast} e^t_i = - \Tr_{h_t}(\nabla^{h_t} e^t_i) = \delta_{ui} \Tr_{h_t}(\Theta_t) -  \Theta_t(e^t_u,e^t_i) - \frac{1}{2} c_t(e^t_u,e^t_i)
\end{eqnarray*}

\noindent
or, equivalently:
\begin{eqnarray*}
& \nabla^{h_t \ast } e^t_u   = \Tr_{h_t}(\Theta_t) - \Theta_t(e^t_u,e^t_u)\\ 
& \nabla^{h_t \ast} e^t_l = - \Theta_t(e^t_l , e^t_u) -   \frac{1}{2} c_t (e^t_u , e^t_l)\, , \qquad  \nabla^{h_t \ast} e^t_n  = -\Theta_t(e^t_n, e^t_u) -  \frac{1}{2} c_t (e^t_u , e^t_n) 
\end{eqnarray*}
	
\noindent
These equations will be useful in the following.
\end{remark}

\noindent
It is natural to compare the evolution flow prescribed by the NS-NS supersymmetry conditions for $(\fre_t,b_t,\phi_t)$ with the NS-NS evolution flow for the associated family $(h_{\fre_t},b_t,\phi_t)$. For this, we first need to compute the scalar curvature of $h_{\fre_t}$, which we denote by $s^{\fre_t}$.

\begin{lemma}
\label{lemma:curvaturesusyflow}
Let $(\fre_t,b_t,\phi_t)$ be a normal NS-NS supersymmetry flow. Then:
\begin{equation*}
s^{\fre_t} =  \vert \Theta_t  \vert^2_{h_t} - \Tr_{h_t}(\Theta_t)^2 + 2 (\nabla^{h_t\ast}\Theta_t)(e^t_u) + 2 \Tr_{h_t}(\nabla^{h_t}_{e^t_u}\Theta_t)  -   (\nabla^{h_t\ast} c_t) (e^t_u)  - \frac{1}{2} \vert c_t (e^t_u)\vert_{h_t}^2
\end{equation*}

\noindent
where $s^{\fre_t}$ denotes the scalar curvature of $h_{\fre_t}$ at a given time $t\in \cI$. 
\end{lemma}

\begin{proof}
Using equations \eqref{eq:susyevol} together with Remark \ref{remark:covariantevolution}, we compute:
\begin{eqnarray*}
&\nabla^{h_t}_{e^t_i}\nabla^{h_t}_{e^t_j} e^t_u = \nabla^{h_t}_{e^t_i}(- \Theta_t(e^t_j) + \Theta_t(e^t_u,e^t_j) e^t_u - \frac{1}{2} \delta_{uj} c_t(e^t_u))\\
& = - (\nabla^{h_t}_{e^t_i}\Theta_t )(e^t_j) - \Theta_t(\nabla^{h_t}_{e^t_i} e^t_j) + ((\nabla^{h_t}_{e^t_i}\Theta_t)(e^t_u,e^t_j) + \Theta_t(\nabla^{h_t}_{e^t_i} e^t_u,e^t_j) +  \Theta_t(e^t_u, \nabla^{h_t}_{e^t_i} e^t_j)) e^t_u\\
& +  \Theta_t(e^t_u,  e^t_j) \nabla^{h_t}_{e^t_i} e^t_u - \frac{1}{2} \delta_{uj} ((\nabla^{h_t}_{e^t_i} c_t)(e^t_u) + c_t (\nabla^{h_t}_{e^t_i} e^t_u))
\end{eqnarray*}

\noindent
On the other hand:
\begin{eqnarray*}
\nabla^{h_t}_{[e_i,e_j]} e^t_u = - \Theta_t([e_i,e_j]) + \Theta_t(e^t_u,[e_i,e_j]) e^t_u - \frac{1}{2} e^t_u([e_i,e_j])  c_t(e^t_u)
\end{eqnarray*}

\noindent
From this formulae we obtain, after a tedious calculation, the following expression:
\begin{eqnarray*}
&   (\cR^{\fre_t}_{e_i e_u} e_u)(e_i) = \vert \Theta_t(e^t_u) \vert^2_{h_t} - \Theta_t(e^t_u,e^t_u) \Tr_{h_t}(\Theta_t) + (\nabla^{h_t\ast}\Theta_t)(e^t_u) + \Tr_{h_t}(\nabla^{h_t}_{e^t_u}\Theta_t)  \\ 
& - \frac{1}{2} (\nabla^{h_t\ast} c_t) (e^t_u) + \frac{1}{2} c_t (\Theta_t (e^t_l) , e^t_l) + \frac{1}{2} c_t (\Theta_t (e^t_n) , e^t_n) - \frac{1}{2} \langle \Theta_t(e^t_u) , c_t(e^t_u)\rangle_{h_t} - \frac{1}{4} \vert c_t (e^t_u)\vert_{h_t}^2
\end{eqnarray*}

\noindent
where summation over repeated indices is assumed. Regarding $e^t_l$ we have:
\begin{eqnarray*}
& \nabla^{h_t}_{e^t_i} \nabla^{h_t}_{e^t_j} e^t_l = \nabla^{h_t}_{e^t_i} ( \Theta_t(e^t_l,e^t_j) e^t_u - \frac{1}{2} f_{b_t} \delta_{uj} e^t_n + \frac{1}{2} \delta_{uj} c_t (e^t_u , e^t_l) e^t_u) \\
& = ( (\nabla^{h_t}_{e^t_i}\Theta_t)(e^t_l,e^t_j) +  \Theta_t(\nabla^{h_t}_{e^t_i}e^t_l,e^t_j) +  \Theta_t(e^t_l, \nabla^{h_t}_{e^t_i}e^t_j) ) e^t_u +   \Theta_t(e^t_l,e^t_j) \nabla^{h_t}_{e^t_i} e^t_u\\
& + \frac{1}{2} \delta_{uj} ((\nabla^{h_t}_{e^t_i}c_t) (e^t_u , e^t_l) + c_t (\nabla^{h_t}_{e^t_i} e^t_u , e^t_l) + c_t (e^t_u ,\nabla^{h_t}_{e^t_i} e^t_l))e^t_u\\
& + \frac{1}{2} \delta_{uj} (c_t ( e^t_u , e^t_l)\nabla^{h_t}_{e^t_i} e^t_u - \dd f_{b_t}(e^t_i) e^t_n - f_{b_t} \nabla^{h_t}_{e^t_i} e^t_n)
\end{eqnarray*}
 
\noindent
as well as:
\begin{eqnarray*}
\nabla^{h_t}_{[e^t_i, e^t_j]} e^t_l = \Theta_t(e^t_l,[e^t_i, e^t_j])  e^t_u - \frac{1}{2} f_{b_t} e_u^t([e^t_i, e^t_j]) e^t_n + \frac{1}{2} c_t (e^t_u , e^t_l) e^t_u([e^t_i, e^t_j]) e^t_u
\end{eqnarray*}

\noindent
from which we obtain:
\begin{eqnarray*}
& (\cR^{\fre_t}_{e_i e_l} e_l)(e_i) = (\nabla_{e^t_u} \Theta_t) (e^t_l,e^t_l) - (\nabla_{e^t_l} \Theta_t) (e^t_l,e^t_u)   - \Theta_t (e^t_l, e^t_l) \Tr_{h_t}\Theta_t + \vert \Theta_t(e^t_l)\vert^2_{h_t} \\
&   - \frac{1}{2} ( c_t(e^t_u , e^t_l) \Theta_t (e^t_u, e^t_l) + (\nabla^{h_t}_{e^t_l} c_t) (e^t_u , e^t_l) + \langle \Theta_t(e^t_l) , c_t(e^t_l)\rangle_{h_t}) - \frac{1}{4} c_t(e^t_u,e^t_l)^2
\end{eqnarray*}

\noindent
We compute similarly for $e^t_n$. We have:
\begin{eqnarray*}
& \nabla^{h_t}_{e^t_i} \nabla^{h_t}_{e^t_j} e^t_n  = \nabla^{h_t}_{e^t_i}( \Theta_t(e^t_n, e^t_j) e^t_u + \frac{1}{2} \delta_{uj} (f_{b_t} e^t_l + \frac{1}{2} c_t (e^t_u , e^t_n)  e^t_u) ) = \\
& = ( (\nabla^{h_t}_{e^t_i}\Theta_t)(e^t_n,e^t_j) +  \Theta_t(\nabla^{h_t}_{e^t_i}e^t_n,e^t_j) +  \Theta_t(e^t_n, \nabla^{h_t}_{e^t_i}e^t_j) ) e^t_u +   \Theta_t(e^t_n,e^t_j) \nabla^{h_t}_{e^t_i} e^t_u\\
& + \frac{1}{2} \delta_{uj} ((\nabla^{h_t}_{e^t_i}c_t) (e^t_u , e^t_n) + c_t (\nabla^{h_t}_{e^t_i} e^t_u , e^t_n) + c_t (e^t_u ,\nabla^{h_t}_{e^t_i} e^t_n) ) e^t_u \\
& + \frac{1}{2} \delta_{uj} ( c_t ( e^t_u , e^t_n)\nabla^{h_t}_{e^t_i} e^t_u + \dd f_{b_t}(e^t_i) e^t_l + f_{b_t} \nabla^{h_t}_{e^t_i} e^t_l)
\end{eqnarray*}

\noindent
as well as:
\begin{eqnarray*}
\nabla^{h_t}_{[e^t_i, e^t_j]} e^t_n = \Theta_t(e^t_n,[e^t_i, e^t_j])  e^t_u + \frac{1}{2} f_{b_t} e_u^t([e^t_i, e^t_j]) e^t_l + \frac{1}{2} c_t (e^t_u , e^t_n) e^t_u([e^t_i, e^t_j]) e^t_u
\end{eqnarray*}

\noindent
from which we obtain:
\begin{eqnarray*}
	& (\cR^{\fre_t}_{e_i e_n} e_n)(e_i) = (\nabla_{e^t_u} \Theta_t) (e^t_n,e^t_n) - (\nabla_{e^t_n} \Theta_t) (e^t_n,e^t_u)   - \Theta_t (e^t_n, e^t_n) \Tr_{h_t}\Theta_t + \vert \Theta_t(e^t_n)\vert^2_{h_t} \\
	&   - \frac{1}{2} ( c_t(e^t_u , e^t_n) \Theta_t (e^t_u, e^t_n) + (\nabla^{h_t}_{e^t_l} c_t) (e^t_u , e^t_n) + \langle \Theta_t(e^t_n) , c_t(e^t_n)\rangle_{h_t}) - \frac{1}{4} c_t(e^t_u,e^t_n)^2
\end{eqnarray*}

\noindent
The previous formulae yields the scalar curvature $s^{h_t}$ of the family of metrics $h_t$ as follows:
\begin{eqnarray*}
s^{\fre_t} = (\cR^{\fre_t}_{e_i e_u} e_u)(e_i)  + (\cR^{\fre_t}_{e_i e_l} e_l)(e_i)  + (\cR^{\fre_t}_{e_i e_n} e_n)(e_i) 
\end{eqnarray*}

\noindent
which after some manipulations gives the expression in the statement of the lemma.
\end{proof}

\noindent
We end this subsection by proving what to the best of our knowledge seems to be the first \emph{compatibility result} between the NS-NS evolution flow and its supersymmetric counterpart.  

\begin{thm}
\label{thm:compatibilityflows}
A triple $(h_{\fre_t},b_t,\phi_t)$ associated to a NS-NS supersymmetric flow $(\fre_t,b_t,\phi_t)$ preserves the constraint equations of the NS-NS system if and only if:
\begin{eqnarray*}
& e^t_u (\partial_t \phi_t) + \Theta_t(e^t_u,e^t_u) \partial_t \phi_t - \frac{1}{2} f_{b_t} = 0\\
& (\nabla^{g\ast}c_t)(e^t_u)  = 0\, ,\qquad c_t(e^t_u) = 0\, , \qquad e_u\wedge( \dd_{\Sigma} f_{b_t} - f_{b_t}\Theta_t(e^t_u)) = 0 \\
&   (\nabla^{h_t}\Theta_t)(e^t_u , e^t_u) - (\nabla^{h_t}_{e^t_u}\Theta_t)(e^t_u) + (\partial_t^2 \phi_t + \Tr_{h_t}(\nabla^{h_t}_{e^t_u}\Theta_t) + (\nabla^{h_t\ast}\Theta_t)(e^t_u)) e^t_u\\
& + \frac{1}{2} \ast_{h_t} (e^t_u \wedge (\dd_{\Sigma} f_{b_t} + \Theta_t(e^t_u))  - \frf_{b_t}  c_t ) = 0 
\end{eqnarray*}
\end{thm}

\begin{proof}
Let $(\fre_t,b_t,\phi_t)$ be a normal NS-NS supersymmetric flow and let $(h_{\fre_t},b_t,\phi_t)$ be its associated NS-NS candidate flow $h_{\fre_t} = e^t_u \otimes e^t_u + e^t_l\otimes e^t_l + e^t_n \otimes e^t_n$. By Proposition \ref{prop:susyevol}, $(\fre_t,b_t,\phi_t)$ is a normal supersymmetric flow if and only if it satisfies all equations in \eqref{eq:susyevol}. On the other hand, by Proposition \ref{prop:normalhyperbolicflow}, the time-dependent constraints of the NS-NS evolution flow are given by:
\begin{eqnarray}
& s^{h_t}  + (\Tr_{h_t}(\Theta_t) +  \partial_t \phi_t)^2 -   \vert\Theta_t\vert^2_{h_t} - \vert \varphi_{\phi_t} \vert_{h_t}^2  - \frac{1}{2}( f_{b_t}^2 + \vert \theta_t \vert^2_{h_t} )  - 2 \nabla^{h_t\ast} \varphi_{\phi_t}  = 0 \label{eq:normalconstraintI} \\
& \dd_{\Sigma} \theta_t = \varphi_{\phi_t} \wedge \theta_t \, , \quad \dd_{\Sigma} \Tr_{h_t}(\Theta_t) + \nabla^{h_t\ast}\Theta_t + \partial_t \varphi_{\phi_t} +\Theta_t (\varphi_{\phi_t}^{\sharp_{h_t}}) - \frac{1}{2}  \frf_{b_t} \theta_t  = 0 \label{eq:normalconstraintII} 
\end{eqnarray}

\noindent
By the last row in \eqref{prop:susyevol} we have:
\begin{eqnarray*}
\varphi_{\phi_t} = \partial_t\phi_t e^t_u  + c_t(e^t_u)\, , \qquad   f_{b_t} = \ast_{h_t} (c_t\wedge e^t_u)
\end{eqnarray*}

\noindent
and therefore:
\begin{eqnarray*}
c_t =  c_t (e^t_u)\wedge e^t_u -   H_{b_t}(e^t_u) =  \varphi_{\phi_t} \wedge e^t_u -   f_{b_t} \ast_{h_t} e^t_u
\end{eqnarray*}

\noindent
Applying $\nabla^{h_t\ast}$ to these expressions, we compute:
\begin{eqnarray*}
& \nabla^{g\ast}\varphi_{\phi_t} = - e_u(\partial_t \phi_t) + (\Tr_{h_t}(\Theta_t) - \Theta_t(e^t_u,e^t_u)) \partial_t \phi_t \\ 
& - (\nabla^{h_t\ast} c_t)(e^t_u) - \langle \Theta_t(e^t_u) , c_t(e^t_u)\rangle_{h_t} - \frac{1}{2} \vert c_t(e^t_u) \vert^2_{h_t}\\
& \nabla^{h_t\ast}c = \nabla^{h_t\ast}\varphi_{\phi_t} e^t_u + \nabla^{h_t}_{e^t_u}\varphi_{\phi_t} -  \nabla^{h_t}_{\varphi_{\phi_t}} e^t_u - \nabla^{h_t\ast}e^t_u \varphi_{\phi_t} + \ast_{h_t} (e^t_u \wedge \dd_{\Sigma} f_{b_t})- f_{b_t} \ast_{h_t}\dd_{\Sigma} e^t_u
\end{eqnarray*}

\noindent
The second equation implies:
\begin{eqnarray*}
(\nabla^{h_t\ast}c)(e^t_u) = \nabla^{h_t\ast}\varphi_{\phi_t}  + (\nabla^{h_t}_{e^t_u}\varphi_{\phi_t})(e^t_u) -  (\nabla^{h_t}_{\varphi_{\phi_t}} e^t_u)(e^t_u) - \nabla^{h_t\ast}e^t_u \partial_t\phi_t 
\end{eqnarray*}

\noindent
which substituted back into the previous expression for $\nabla^{g\ast}\varphi_{\phi_t}$ gives:
\begin{eqnarray*}
& 2 \nabla^{g\ast}\varphi_{\phi_t} = - e^t_u(\partial_t \phi_t) + 2 \nabla^{h_t\ast}e^t_u \partial_t\phi_t - \langle \Theta_t(e^t_u) , \Theta_t(e^t_u)\rangle_{h_t} - \frac{1}{2} \vert c_t (e^t_u)\vert^2_{h_t} \\
&- (\nabla^{h_t}_{e^t_u}\varphi_{\phi_t})(e^t_u) +  (\nabla^{h_t}_{\varphi_{\phi_t}} e^t_u)(e^t_u)
\end{eqnarray*}

\noindent
Substituting this equation back into the previous expression for $(\nabla^{h_t\ast}c)(e^t_u)$, we obtain:
\begin{eqnarray}
& 2 (\nabla^{h_t\ast}c)(e^t_u) =  (\nabla^{h_t}_{e^t_u}\varphi_{\phi_t})(e^t_u) -   (\nabla^{h_t}_{\varphi_{\phi_t}} e^t_u)(e^t_u) -   e^t_u(\partial_t \phi_t) -  \langle \Theta_t(e^t_u) , c_t(e^t_u)\rangle_{h_t} - \frac{1}{2} \vert c_t(e^t_u) \vert^2_{h_t}\nonumber\\
&  = \Theta_t(\varphi_{\phi_t},e^t_u) - \Theta_t(e^t_u,e^t_u) \partial_t\phi_t + \frac{1}{2} c_t(e^t_u , \varphi_{\phi_t}) -  \langle \Theta_t(e^t_u) , c_t(e^t_u)\rangle_{h_t} - \frac{1}{2} \vert c_t(e^t_u) \vert^2_{h_t} = 0 \label{eq:cu0}
\end{eqnarray}

\noindent
If the triple $(h_{\fre_t},b_t,\phi_t)$ associated to $(\fre_t,b_t,\phi_t)$  is a NS-NS evolution flow then it needs to satisfy the first equation in \eqref{eq:normalconstraintII}, which can be equivalently written as follows:
\begin{equation*}
\nabla^{g\ast}c_t + c_t(\varphi_{\phi_t}) = \nabla^{g\ast}c_t + \partial_t \phi_t c_t(e^t_u) + c_t(c_t(e^t_u)) =  0
\end{equation*} 
 
\noindent
This implies in turn:
\begin{equation*}
(\nabla^{g\ast}c_t)(e^t_u) + c_t(\varphi_{\phi_t}, e^t_u) = (\nabla^{g\ast}c_t)(e^t_u)   + c_t(c_t(e^t_u),e^t_u) = (\nabla^{g\ast}c_t)(e^t_u)   - \vert c_t(e^t_u)\vert^2_{h_t} =  0
\end{equation*} 

\noindent
Since by Equation \eqref{eq:cu0} we have $(\nabla^{g\ast}c_t)(e^t_u)  = 0$, we conclude that $(\fre_t,b_t,\phi_t)$ is a supersymmetric NS-NS evolution flow that preserves the constraint equations \eqref{eq:normalconstraintI} and \eqref{eq:normalconstraintII} of the NS-NS evolution flow only if:
\begin{eqnarray*}
(\nabla^{g\ast}c_t)(e^t_u)  = 0\, ,\qquad c_t(e^t_u) = 0
\end{eqnarray*}

\noindent
Hence, we will assume both conditions in the following. Elaborating on the expression obtained above for $\nabla^{h_t\ast}c_t$, we obtain:
\begin{eqnarray*}
\nabla^{h_t\ast} c_t + \ast_{h_t}(e_u\wedge \dd_{\Sigma} f_{b_t}) + f_{b_t} \ast_{h_t} \dd_{\Sigma} e^t_u = 0 
\end{eqnarray*}

\noindent
and hence we conclude that the triple $(h_{\fre_t},b_t,\phi_t)$ associated to $(\fre_t,b_t,\phi_t)$  satisfies the first equation in \eqref{eq:normalconstraintII} if and only if:
\begin{equation*}
(\nabla^{g\ast}c_t)(e^t_u)  = 0\, ,\qquad c_t(e^t_u) = 0\, , \qquad e_u\wedge( \dd_{\Sigma} f_{b_t} - f_{b_t}\Theta_t(e^t_u)) = 0 
\end{equation*}

\noindent
Note that here we have substituted $\dd_{\Sigma} e^t_u$ by its expression as given in Proposition \ref{prop:susyevolII}. On the other hand, by the first equation in Remark \ref{remark:covariantevolution} we have:
\begin{eqnarray*}
\nabla^{h_t\ast}\nabla^{h_t} e^t_u  = - \nabla^{h_t\ast}\Theta_t + \nabla^{h_t\ast}(\Theta_t(e^t_u)\otimes e^t_u)  	
\end{eqnarray*}

\noindent
We compute:
\begin{eqnarray*}
 \nabla^{h_t\ast}(\Theta_t(e^t_u)\otimes e^t_u) = ((\nabla^{h_t\ast}\Theta_t)(e^t_u) + \vert \Theta_t\vert^2_{h_t} - 2 \vert \Theta_t(e^t_u) \vert_{h_t}^2)e^t_u  + \Theta_t(\Theta_t(e^t_u))
\end{eqnarray*}

\noindent
Hence:
\begin{equation*}
\nabla^{h_t\ast}\nabla^{h_t} e^t_u  = - \nabla^{h_t\ast}\Theta_t + ((\nabla^{h_t\ast}\Theta_t)(e^t_u) + \vert \Theta_t\vert^2_{h_t} - 2 \vert \Theta_t(e^t_u) \vert_{h_t}^2)e^t_u  + \Theta_t(\Theta_t(e^t_u))
\end{equation*}

\noindent
The previous equation gives the \emph{rough Laplacian} of $e^t_u$, which by virtue of the Weitzenb\"ock formula we are going to compare with the Hodge Laplacian applied to $e^t_u$. We compute:
\begin{eqnarray*}
& \dd_{\Sigma} \nabla^{h_t \ast} e^t_u  = \dd_{\Sigma} (\Tr_{h_t}(\Theta_t) - \Theta_t(e^t_u,e^t_u)) = \\
& = \dd_{\Sigma} \Tr_{h_t}(\Theta_t) - (\nabla^{h_t}\Theta_t)(e^t_u , e^t_u) + 2 \Theta_t(\Theta_t(e^t_u)) -2 \Theta_t(e^t_u,e^t_u) \Theta_t(e^t_u)\\
& \dd_{\Sigma}^{h_t \ast} \dd_{\Sigma} e^t_u  = \nabla^{h_t \ast} \dd_{\Sigma} e^t_u = \nabla^{h_t \ast} (\Theta_t(e^t_u)\wedge e^t_u) =  ((\nabla^{h_t\ast}\Theta_t)(e^t_u) + \vert \Theta_t\vert^2_{h_t} - 2 \vert \Theta_t(e^t_u) \vert_{h_t}^2)e^t_u   \\
& + (\nabla^{h_t}_{e^t_u}\Theta_t)(e^t_u) + (2 \Theta_t(e^t_u,e^t_u) - \Tr_{h_t}(\Theta_t)) \Theta_t(e^t_u)
\end{eqnarray*}

\noindent
and thus:
\begin{eqnarray*}
& \Delta_{h_t} e^t_u = (\dd_{\Sigma}^{h_t \ast} \dd_{\Sigma} + \dd_{\Sigma} \dd_{\Sigma}^{h_t \ast}) e^t_u = ((\nabla^{h_t\ast}\Theta_t)(e^t_u) + \vert \Theta_t\vert^2_{h_t} - 2 \vert \Theta_t(e^t_u) \vert_{h_t}^2)e^t_u  + (\nabla^{h_t}_{e^t_u}\Theta_t)(e^t_u)  \\
& - \Tr_{h_t}(\Theta_t)\Theta_t(e^t_u) + \dd_{\Sigma} \Tr_{h_t}(\Theta_t) - (\nabla^{h_t}\Theta_t)(e^t_u , e^t_u) + 2 \Theta_t(\Theta_t(e^t_u)) 
\end{eqnarray*}

\noindent
The last ingredient to apply Weitzenb\"ock formula is the evaluation of $e_u^t$ in the Ricci tensor $\Ric^{h_t}$ of $h_t$, namely:
\begin{eqnarray*}
\Ric^{\fre_t}(e^t_u) = \cR_{e^t_u e^t_l} e^t_l + \cR_{e^t_u e^t_n} e^t_n
\end{eqnarray*}

\noindent
We compute:
\begin{eqnarray*}
& \nabla^{h_t}_{e^t_u}\nabla^{h_t}_{e^t_l} e^t_l = \nabla^{h_t}_{e^t_u} (\Theta_t(e^t_l , e^t_l) e^t_u) = (\nabla^{h_t}_{e^t_u}\Theta_t)(e^t_l , e^t_l) e^t_u + 2\Theta_t(\nabla^{h_t}_{e^t_u}e^t_l , e^t_l) e^t_u + \Theta_t(e^t_l , e^t_l) \nabla^{h_t}_{e^t_u} e^t_u\\
& = ((\nabla^{h_t}_{e^t_u}\Theta_t)(e^t_l,e^t_l) + 2\Theta_t(e^t_u,e^t_l)^2 - f_{b_t} \Theta_t(e^t_l,e^t_n) + \Theta_t(e^t_u , e^t_u) \Theta_t(e^t_l , e^t_l))e^t_u - \Theta_t(e^t_l , e^t_l) \Theta_t(e^t_u)   \\
& \nabla^{h_t}_{e^t_l}\nabla^{h_t}_{e^t_u} e^t_l = ((\nabla^{h_t}_{e^t_l}\Theta_t)(e^t_u,e^t_l) - \vert \Theta_t(e^t_l)\vert^2_{h_t} + \Theta_t(e^t_u,e^t_l)^2 + \Theta_t(e^t_u,e^t_u)\Theta_t(e^t_l,e^t_l)) e^t_u\\
& - \Theta_t(e^t_u , e^t_l) ( \Theta_t(e^t_l , e^t_l) e^t_l +  \Theta_t(e^t_l , e^t_n) e^t_n)  - \frac{1}{2} (e^t_l(f_{b_t}) e^t_n + f_{b_t} \Theta_t(e^t_l,e^t_n) e^t_u)\\
&  = ((\nabla^{h_t}_{e^t_l}\Theta_t)(e^t_u,e^t_l) - \vert \Theta_t(e^t_l)\vert^2_{h_t} + 2 \Theta_t(e^t_u,e^t_l)^2 + \Theta_t(e^t_u,e^t_u)\Theta_t(e^t_l,e^t_l)) e^t_u\\
& - \Theta_t(e^t_u , e^t_l) \Theta_t(e^t_l)  - \frac{1}{2} (e^t_l(f_{b_t}) e^t_n + f_{b_t} \Theta_t(e^t_l,e^t_n) e^t_u)
\end{eqnarray*}

\noindent
together with:
\begin{eqnarray*}
& \nabla^{h_t}_{[e^t_u,e^t_l]}e^t_l = \Theta_t(e^t_l, [e^t_u,e^t_l]) e^t_u - \frac{1}{2} f_{b_t} e_u^t([e^t_u,e^t_l]) e^t_n \\
& = (\vert \Theta(e^t_l)\vert^2_{h_t}  - \frac{1}{2} f_{b_t} \Theta_t(e^t_l,e^t_n)) e^t_u - \frac{1}{2} f_{b_t}  \Theta_t(e^t_u,e^t_l) e^t_n  
\end{eqnarray*}

\noindent
Here we have used that:
\begin{eqnarray*}
[e^t_u,e^t_l] = \nabla^{h_t}_{e^t_u} e^t_l - \nabla^{h_t}_{e^t_l} e^t_u = \Theta_t(e^t_l) - \frac{1}{2} f_{b_t} e^t_n 
\end{eqnarray*}

\noindent
Similarly:
\begin{eqnarray*}
& \nabla^{h_t}_{e^t_u}\nabla^{h_t}_{e^t_n} e^t_n = \nabla^{h_t}_{e^t_u} (\Theta_t(e^t_n , e^t_n) e^t_u) = (\nabla^{h_t}_{e^t_u}\Theta_t)(e^t_n , e^t_n) e^t_u + 2\Theta_t(\nabla^{h_t}_{e^t_u}e^t_n , e^t_n) e^t_u + \Theta_t(e^t_n , e^t_n) \nabla^{h_t}_{e^t_u} e^t_u\\
& = ((\nabla^{h_t}_{e^t_u}\Theta_t)(e^t_n,e^t_n) + 2\Theta_t(e^t_u,e^t_n)^2 + f_{b_t} \Theta_t(e^t_l,e^t_n) + \Theta_t(e^t_u , e^t_u) \Theta_t(e^t_n , e^t_n))e^t_u - \Theta_t(e^t_n , e^t_n) \Theta_t(e^t_u)   \\
& \nabla^{h_t}_{e^t_n}\nabla^{h_t}_{e^t_u} e^t_n = ((\nabla^{h_t}_{e^t_n}\Theta_t)(e^t_u,e^t_n) - \vert \Theta_t(e^t_n)\vert^2_{h_t} + \Theta_t(e^t_u,e^t_n)^2 + \Theta_t(e^t_u,e^t_u)\Theta_t(e^t_n,e^t_n)) e^t_u\\
& - \Theta_t(e^t_u , e^t_n) ( \Theta_t(e^t_n , e^t_n) e^t_n +  \Theta_t(e^t_l , e^t_n) e^t_l)  + \frac{1}{2} (e^t_n(f_{b_t}) e^t_l + f_{b_t} \Theta_t(e^t_l,e^t_n) e^t_u)\\
&  = ((\nabla^{h_t}_{e^t_n}\Theta_t)(e^t_u,e^t_n) - \vert \Theta_t(e^t_n)\vert^2_{h_t} + 2 \Theta_t(e^t_u,e^t_n)^2 + \Theta_t(e^t_u,e^t_u)\Theta_t(e^t_n,e^t_n)) e^t_u\\
& - \Theta_t(e^t_u , e^t_n) \Theta_t(e^t_n)  + \frac{1}{2} (e^t_n(f_{b_t}) e^t_l + f_{b_t} \Theta_t(e^t_l,e^t_n) e^t_u)
\end{eqnarray*}

\noindent
together with:
\begin{eqnarray*}
& \nabla^{h_t}_{[e^t_u,e^t_n]}e^t_n = \Theta_t(e^t_n, [e^t_u,e^t_n]) e^t_u + \frac{1}{2} f_{b_t} e_u^t([e^t_u,e^t_n]) e^t_l \\
& = (\vert \Theta(e^t_n)\vert^2_{h_t}  + \frac{1}{2} f_{b_t} \Theta_t(e^t_l,e^t_n)) e^t_u + \frac{1}{2} f_{b_t}  \Theta_t(e^t_u,e^t_n) e^t_l  
\end{eqnarray*}

\noindent
Here we have used that:
\begin{eqnarray*}
[e^t_u,e^t_n] = \nabla^{h_t}_{e^t_u} e^t_n - \nabla^{h_t}_{e^t_n} e^t_u = \Theta_t(e^t_n) + \frac{1}{2} f_{b_t} e^t_l
\end{eqnarray*}

\noindent
Therefore:
\begin{eqnarray*}
& \Ric^{\fre_t}(e^t_u) = (\Tr_{h_t}(\nabla_{e^t_u}\Theta_t) + (\nabla^{h_t\ast}\Theta_t)(e^t_u)) e^t_u + \Theta_t(\Theta_t(e^t_u)) - \Tr_{h_t}(\Theta_t) \Theta_t(e^t_u)\\
& \frac{1}{2} \ast_{h_t} (e^t_u\wedge (\Theta_t(e^t_u) + \dd_{\Sigma} f_{b_t}))
\end{eqnarray*}

\noindent
and the Weitzenböck formula $\nabla^{h_t\ast} \nabla^{h_t} e^t_u = \Delta_{h_t} e^t_u + \Ric(e^t_u)$ applied to $e^t_u$ gives:
\begin{eqnarray*}
& \dd_{\Sigma} \Tr_{h_t}(\Theta_t) + \nabla^{h_t \ast} \Theta_t = (\nabla^{h_t}\Theta_t)(e^t_u , e^t_u) - (\nabla^{h_t}_{e^t_u}\Theta_t)(e^t_u) + (\Tr_{h_t}(\nabla^{h_t}_{e^t_u}\Theta_t) + (\nabla^{h_t\ast}\Theta_t)(e^t_u)) e^t_u\\
& + \frac{1}{2} \ast_{h_t} (e^t_u \wedge (\dd_{\Sigma} f_{b_t} + \Theta_t(e^t_u)))
\end{eqnarray*}
 
\noindent
Hence we conclude that the triple $(h_{\fre_t},b_t,\phi_t)$ associated to $(\fre_t,b_t,\phi_t)$  satisfies both equations in \eqref{eq:normalconstraintII} if and only if:
\begin{eqnarray*}
& (\nabla^{g\ast}c_t)(e^t_u)  = 0\, ,\qquad c_t(e^t_u) = 0\, , \qquad e^t_u\wedge( \dd_{\Sigma} f_{b_t} - f_{b_t}\Theta_t(e^t_u)) = 0 \\
& \partial_t^2 \phi_t\, e^t_u  - \frac{1}{2}  \frf_{b_t} \ast_{h_t} c_t + (\nabla^{h_t}\Theta_t)(e^t_u , e^t_u) - (\nabla^{h_t}_{e^t_u}\Theta_t)(e^t_u) + (\Tr_{h_t}(\nabla^{h_t}_{e^t_u}\Theta_t) + (\nabla^{h_t\ast}\Theta_t)(e^t_u)) e^t_u\\
& + \frac{1}{2} \ast_{h_t} (e^t_u \wedge (\dd_{\Sigma} f_{b_t} + \Theta_t(e^t_u))) = 0
\end{eqnarray*}

\noindent
Comparing now Equation \eqref{eq:normalconstraintI} with Lemma \ref{lemma:curvaturesusyflow} using the expression for $\nabla^{g\ast}\varphi_{\phi_t}$ obtained above we conclude. 
\end{proof}

\noindent
This is the starting point for a deeper analysis of the rich interaction between the NS-NS and the supersymmetric NS-NS evolution flows, which we plan to develop in the future. In the same sense that supersymmetry provides first order \emph{partial} integrability of the NS-NS system, we expect that the supersymmetric NS-NS evolution flow can provide first order partial integrability of the full NS-NS evolution flow. Our main conjecture, based on prior results for parallel spinors on globally hyperbolic Lorentzian four-manifolds \cite{Murcia:2020zig,Murcia:2021dur}, is that evolving initial data admissible to \emph{both} flows by the supersymmetric NS-NS flow will produce a NS-NS evolution flow. An interesting consequence of Theorem \ref{thm:compatibilityflows} that reflects the \emph{rigidity} of satisfying both the NS-NS evolution flow and its supersymmetric counterpart is given in the following corollary.

\begin{cor}
A triple $(h_{\fre_t},b_t,\phi_t)$ associated to a NS-NS supersymmetric flow $(\fre_t,b_t,\phi_t)$ preserves the constraint equations of the NS-NS system only if $(h_{\fre_t},\Theta_t)$ satisfies the Hamiltonian constraint equation of a globally hyperbolic Lorentzian four-manifold equipped with a spinor parallel with respect to the Levi-Civita connection. 
\end{cor}


\subsection{The NS-NS constraint equations}

 
Using the explicit form of the normal supersymmetry NS-NS evolution system as given in \eqref{eq:susyevol}, or equivalently in Corollary \eqref{cor:susyevolcor} or Proposition \ref{prop:susyevolII}, we can extract the \emph{constraint equations} of the evolution system. These can be considered as the general constraint equations of  the evolution flow determined by the NS-NS supersymmetry conditions as they do not depend on the choice $\lambda_t = 1$. 

\begin{definition}
Let $(\cC,\cX)$ respectively be a bundle gerbe $\cC$ and a principal $\mathbb{Z}$ bundle $\cX$ on $\Sigma$. The \emph{NS-NS sypersymmmetry constraint equations} on $(\cC,\cX)$ is given by:
\begin{eqnarray}
& f_{b} = \ast_{h}(c\wedge e^t_u) \, , \quad \dd_{\Sigma} e_u = \Theta_{\fre\frv}(e_u)\wedge e_u -\frac{1}{2} e_u \wedge c (e_u) \nonumber \\
& \dd_{\Sigma} e
_l = \Theta_{\fre\frv}(e_l)\wedge e^t_u + \frac{1}{2} f_{b}  \ast_{h} e_l\, , \qquad  \dd_{\Sigma} e_n = \Theta_{\fre\frv}(e_n)\wedge e^t_u + \frac{1}{2} f_{b}  \ast_{h} e_n  \label{eq:susyNSNSconstraints}\\
& \varphi_{\phi} = \psi e_u + c(e_u)\, , \quad  0 = [\Theta_{\fre\frv} (e_u) - \frac{1}{2} c(e_u)]\in H^1(\Sigma,\mathbb{R}) \nonumber
\end{eqnarray}
	
\noindent
for tuples $(\fre,\frv,b,c,\phi,\psi)$ consisting on a global coframe $\fre$ on $\Sigma$, a triplet of one-forms $\frv \in \Omega^1(\Sigma,\mathbb{R}^3)$, a curving $b$ on $\cC$, a two-form $c\in \Omega^2(\Sigma)$ and  a pair of functions $\phi , \psi\in C^{\infty}(\Sigma)$.
\end{definition}

\noindent
As explained in Section \eqref{sec:evolutionskewtorsion} here we have set:
\begin{equation*}
\Theta_{\fre\frv} = \fre_u \odot \frv_u + \fre_l\odot \frv_l + \fre_n\odot \frv_n
\end{equation*}

\noindent
and furthermore we can consider pairs $(\fre,\frv)$ as sections of the tangent bundle to the bundle of oriented frames $F(\Sigma)$ of $\Sigma$, or equivalently, the tangent bundle of the space of sections of $F(\Sigma)$, which we will denote by $TF(\Sigma)$. Similarly, a pair $(b,c)$ can be considered as an element of the tangent bundle of the set of curvings $B(\cC)$ on $\cC$, which is an affine space modeled on the two-forms $\Omega^2(\Sigma)$, and the pair of functions $(\phi,\psi)\in TC^{\infty}(\Sigma)$ can be understood as an element in the tangent bundle of the space of smooth functions on $\Sigma$. Hence, the \emph{configuration space} of the NS-NS sypersymmmetry constraint equations on $\cC$ is the direct product space:
\begin{eqnarray*}
\mConf(\cC,\cX) = TF(\Sigma)\times TB(\cC) \times T C^{\infty}(\Sigma) = T(F(\Sigma)\times B(\cC)\times C^{\infty}(\Sigma))
\end{eqnarray*}

\noindent
which is closely related to the configuration space of the constraint equations of the NS-NS evolution flow. Hence, the NS-NS supersymmetry constraint equations should be understood as the conditions imposed by \emph{supersymmetry} on the initial data of the NS-NS evolution flow, rather than being understood as the constraint equations of the NS-NS supersymmetry flow itself, which is of first order in all of its variables. Consequently, it is natural to study the  NS-NS supersymmetry constraint equations in combination with NS-NS constraint equations.  

\begin{definition}
A \emph{Cauchy NS-NS tuple} $(\fre,\frv,b,c,\phi,\psi)\in \mConf(\cC,\cX)$ is an element of the configuration space of the NS-NS supersymmetry constraint equations. A \emph{supersymmetric NS-NS Cauchy tuple} is a Cauchy tuple that satisfies the NS-NS supersymmetry constraint equations given in \eqref{eq:susyNSNSconstraints}.
\end{definition}

\noindent
We denote by:
\begin{equation*}
\mSol(\cC,\cX)\subset \mConf(\cC,\cX) \, , 
\end{equation*}

\noindent
the subspace of solutions of the NS-NS supersymmetry constraint equations. We have a canonical map:
\begin{equation*}
\mathbb{Q}\colon \mConf(\cC,\cX) \to \Conf_{\frc}(\cC,\cX)\, , \qquad (\fre,\frv,b,c,\phi,\psi) \mapsto (h_{\fre},\Theta_{\fre\frv},b,c,\phi,\psi)\, ,
\end{equation*}

\noindent
where as usual:
\begin{equation*}
h_{\fre} = e_u\otimes e_u + e_l\otimes e_l + e_n \otimes e_n\, , \qquad \Theta_{\fre\frv} = e_u\odot \frv_u +  e_l\odot \frv_l +  e_n \odot \frv_n
\end{equation*}

\noindent
We will denote by $\mathbb{Q}_s\colon \mSol(\cC,\cX) \to \Conf_{\frc}(\cC,\cX)$ the restriction of $\mathbb{Q}$ to the solutions of the NS-NS supersymmetry constraint equations. These maps relate the initial data of the NS-NS equations and the NS-NS supersymmetry conditions and lead to a rich interplay which, to the best of our knowledge, has not been studied systematically in the literature.

\begin{definition}
An element $(h,\Theta,b,c,\phi,\psi)\in \Conf_{\frc} (\cC,\cX)$ is supersymmetric if it belongs to the image of $\mathbb{Q}_s\colon \mSol(\cC,\cX) \to \Conf_{\frc}(\cC,\cX)$. If in addition $(h,\Theta,b,c,\phi,\psi)\in \Sol_{\frc} (\cC,\cX)$, namely it is also a solution of the NS-NS constraint equations then $(h,\Theta,b,c,\phi,\psi)$ is a \emph{NS-NS supersymmetric initial data}.
\end{definition}
 
\noindent
As a direct consequence of Theorem \ref{thm:compatibilityflows}, we obtain the following result. 

\begin{cor}
A three-manifold $\Sigma$ is the Cauchy hypersurface of a NS-NS supersymmetric solution only if it admits a pair $(\cC,\cX)$ and a Cauchy NS-NS tuple $(\fre,\frv,b,c,\phi,\psi)\in \mConf(\cC,\cX)$ satisfying the following differential system:
\begin{eqnarray*}
& e_u (\psi) + \Theta_{\fre\frv}(e_u,e_u) \psi - \frac{1}{2} f_{b} = 0\\
& (\nabla^{h\ast}c)(e_u)  = 0\, ,\qquad c(e_u) = 0\, , \qquad e_u\wedge( \dd_{\Sigma} f_{b} - f_{b}\Theta_{\fre\frv}(e_u)) = 0 
\end{eqnarray*}

\noindent
where $\fre = (e_u,e_l,e_n)$.
\end{cor}

\noindent
This corollary is only the starting point of the study of the initial conditions to both the NS-NS evolution flow and its supersymmetric counterpart, which we hope can lead to exciting new mathematical results and applications to the differential geometry and topology of three-dimensional Riemannian manifolds. 




\appendix




\renewcommand{\leftmark}{Chapter \thechapter. Bundle Gerbes}

\chapter{Abelian bundle gerbes}
\label{chapter:BundleGerbes}

 



Let $M$ be an oriented manifold. The notion of \emph{abelian $\U(1)$ gerbe} \cite{Murray}, or abelian gerbe for short, on a manifold $M$ appears as one of the cases in the \emph{hierarchy} of \emph{geometric realizations} of the singular cohomology groups of $M$. A geometric realization of a singular cohomology group $H^k(M,\mathbb{Z})$ is, roughly speaking, a possibly higher category whose objects, modulo \emph{one-isomorphisms}, are in natural bijection with $H^k(M,\mathbb{Z})$. The interest in finding these geometric realizations steams, among other reasons, from the fact that a geometric realization of a cohomology class in $H^k(M,\mathbb{Z})$ generally has a rich higher groupoid of symmetries which is lost when descending to cohomology. This higher groupoid of symmetries plays a fundamental role in mathematical gauge theory. The first steps in the geometric realization of singular cohomology are well-known:
\begin{itemize}
\item $H^0(M,\mathbb{Z})$ is in bijection with the $\mathbb{Z}$-valued continuous functions on $M$. Indeed, on every connected component of $M$, the set of continuous functions into $\mathbb{Z}$ is in bijection with $\mathbb{Z}$ and hence we obtain as many copies of $\mathbb{Z}$ as connected components of $M$, as expected.
	
\item $H^1(M,\mathbb{Z})$ is in bijection with the $\U(1)$-valued smooth functions on $M$ modulo homotopy. Indeed, by Brown's representability theorem, we have $H^1(M,\mathbb{Z}) = [M , K(1,\mathbb{Z})]$, where $[M , K(1,\mathbb{Z})]$ are the the homotopy classes of maps from $M$ to the Eilenberg-McLane space $K(1,\mathbb{Z}) = S^1$.
	
\item  $H^2(M,\mathbb{Z})$ is in bijection with the set of $\U(1)$-bundles on $M$ modulo based isomorphism classes of principal $\U(1)$-bundles on $M$. The bijection is established by the first Chern class, which assigns to each bundle its characteristic class in $H^2(M,\mathbb{Z})$. This is a celebrated result that is extensively used throughout mathematics and theoretical physics.
\end{itemize}

\noindent
The next step in this hierarchy, namely the geometric realization of $H^3(M,\mathbb{Z})$, is where we encounter the notion of Abelian gerbe and the first instance of what is usually referred to as \emph{higher geometry} \cite{BaezSchreiber} in the literature. There are various variations on the notion of abelian gerbe \cite{Brylinski,Giraud}. In this dissertation, we consider abelian gerbes as introduced by Murray in \cite{Murray}. These are commonly called \emph{bundle gerbes} in the literature. Before introducing them, let's recall some basic notions of the theory of simplicial manifolds and simplicial line bundles, which can be consider as precursor notions of that of a bundle gerbe. First, remember that the \emph{simplex category} $\Delta$ is the category whose objects are finite ordinal sets, usually denoted by $[0] = \left\{ 0\right\}$, $[1] = \left\{ 0,1\right\}$, $[2] = \left\{ 0,1,2\right\}$ and so on, and whose morphisms are order-preserving maps.

\begin{definition}
	Let $C$ be a category. A \emph{simplicial object} in $C$ is a contravariant functor from $\Delta$ to $C$. 
\end{definition}

\noindent
A morphism between two simplicial objects is a natural transformation between their corresponding functors. Let Man denote the category of smooth manifolds and smooth maps. 

\begin{definition}
	A \emph{simplicial manifold} is a simplicial object on Man.
\end{definition}

\noindent
Unraveling the definition of simplicial manifold as contravariant functor $X\colon \Delta \to \mathrm{Man}$, it follows that a simplicial manifold $X\colon \Delta \to \mathrm{Man}$ can be equivalently defined as a sequence of manifolds $X_{\bullet} := \left\{ X_i \right\}_{i\in \mathbb{N}}$ given by:
\begin{equation*}
	X_i := X([i])\, , \quad [i]\in \Ob(\Delta)\, , \quad i\in\mathbb{N}
\end{equation*}

\noindent
together with smooth maps $\Delta (o) \colon X_j \to X_i$ for every arrow $o \colon [i]\to [j]$ in $\Delta$ satisfying the compatibility condition:
\begin{eqnarray*}
	\Delta(o_2)\Delta(o_1) = \Delta(o_1\circ o_2)
\end{eqnarray*}

\noindent
for every pair of composable arrows $o_1$ and $o_2$ in $\Delta$. It is well-known that every such smooth map $\Delta (o) \colon X_j \to X_i$ can  be written as a composition of collections of smooth maps:
\begin{eqnarray*}
	d_i \colon X_n \to X_{n-1}\, , \qquad s_i \colon X_n \to X_{n+1}\, , \qquad n\in \mathbb{N}\, , \qquad i = 0, \hdots, n
\end{eqnarray*}

\noindent
satisfying, for each $0\leq i \leq n$, the following identities:
\begin{eqnarray*}
	& d_i d_j = d_{j-1} d_i\, , \qquad d_i s_j = s_{j-1} d_i \, , \qquad i< j\\
	& d_j s_j = d_{j+1} s_j = \Id\\
	& d_i s_j = s_j d_{j-1} \, , \qquad i>j+1\\
	& s_i s_j = s_{j+1} s_i\, , \qquad i \leq j
\end{eqnarray*}

\noindent
The functions $d_i \colon X_n \to X_{n-1}$ are commonly called the \emph{face maps} whereas the functions $s_i \colon X_n \to X_{n+1}$ are commonly referred to as the \emph{degeneracy maps}. In this context a morphism of simplicial manifolds $X_{\bullet} \to Y_{\bullet}$ consists of a sequence of maps $X_j \to Y_j$ that commute with the face and degeneracy maps. Canonically associated to any simplicial manifold $X_{\bullet}$, there is a bicomplex $\left\{ \Omega^i(X_j) \right\}_{i,j\in\mathbb{N}}$ with differentials:
\begin{eqnarray*}
	\frD_{i,j} \colon \Omega^i(X_j) \to \Omega^{i+1}(X_j)\oplus  \Omega^i(X_{j+1}) \, , \qquad \alpha \mapsto (-1)^j \dd\alpha \oplus \delta \alpha
\end{eqnarray*}

\noindent
where $\dd \colon \Omega^i(X_j) \to \Omega^{i+1}(X_j)$ is the standard exterior derivative on $X_j$ and:
\begin{equation*}
	\delta\colon \Omega^i(X_j) \to   \Omega^i(X_{j+1}) \, , \qquad \alpha \mapsto \delta \alpha = \sum_{k=0}^{j+1} (-1)^k d^{\ast}_k\alpha
\end{equation*}

\noindent
These differentials combine into a total differential:
\begin{eqnarray*}
	\frD \colon \bigoplus_{p+q =r} \Omega^p(X_q) \to \bigoplus_{p+q =r +1}\Omega^p(X_q) 
\end{eqnarray*}

\noindent
The cohomology of this complex is by definition the simplicial de Rham cohomology of $X_{\bullet}$. Let $\left\{ X_i \right\}_{i\in \mathbb{N}}$ be a simplicial manifold and let $\cP\to X_j$ be a principal $\U(1)$ bundle for some fixed $j\geq 0$. We define a new principal $\U(1)$ bundle $\delta(\cP) \to X_{j+1}$ on $X_{j+1}$ by:
\begin{eqnarray*}
	\delta(\cP) = \bigotimes_{k=0}^j d_k^{\ast}(\cP)^{(-1)^k}
\end{eqnarray*}

\noindent
It follows that $\delta^2(\cP)$ is canonically trivial for every $\U(1)$ bundle $\cP$. We will denote the canonical trivializing section of $\delta^2(\cP)$ by $1\in \Gamma(\delta^2(\cP))$.

\begin{definition}
	A simplicial line bundle over $X_{\bullet}$ is a pair $(\cP,\sigma)$ consisting of a $\U(1)$ bundle $\cP\to X_1$ and a section $\sigma \in \Gamma(\delta \cP)$ such that $\delta(\sigma) = 1 \in \Gamma(\delta^2(\cP))$.
\end{definition}

\begin{definition}
	A bundle gerbe on $M$ is a pair $(\cP,Y)$ consisting of a submersion $Y\to M$ and a simplicial line bundle $\cP$ on the shifted simplicial manifold $Y^{[\bullet+1]}$ generated by $Y\to M$. 
\end{definition}

\noindent
Here $Y^{[\bullet +1]}$ is the simplicial manifold given by the sequence of manifolds:

\begin{equation*}
	Y_i = Y^{[i+1]}\, , \qquad i\in \mathbb{N}
\end{equation*}

\noindent
where $Y^{[i+1]}$ is the $(i+1)$-fold fibered product of $Y\to M$ with itself. Hence, if $(\cP,Y)$ is a bundle gerbe then $\cP\to Y^{[2]}$ is a $\U(1)$ bundle on $Y^{[2]}$. By definition, we have:
\begin{equation*}
	\sigma(y_1,y_2,y_3) \in \cP_{(y_2,y_3)}\otimes \cP_{(y_1,y_3)}^{\ast} \otimes \cP_{(y_1,y_2)}
\end{equation*}

\noindent
for every $y_1 , y_2 , y_3 \in Y$, where the subscript in $\cP$ denotes restriction to the given point in $Y^[2]$. Since $\cP$ is a $\U(1)$ bundle, it follows that there exists a bundle morphism:
\begin{eqnarray*}
	\mu \colon \cP_{(y_1,y_2)}\otimes \cP_{(y_2,y_3)}\to \cP_{(y_1,y_3)}
\end{eqnarray*}

\noindent
such that:
\begin{equation*}
	\sigma(y_1,y_2,y_3) = z_1 \otimes \mu(z_1 , z_2)\otimes z_2 \in \cP_{(y_2,y_3)}\otimes \cP_{(y_1,y_3)}^{\ast} \otimes \cP_{(y_1,y_2)}	
\end{equation*}

\noindent
for every $y_1 , y_2 , y_3 \in Y$. The bundle morphism $\mu\colon \cP_{(y_1,y_2)}\otimes \cP_{(y_2,y_3)}\to \cP_{(y_1,y_3)}$ is the \emph{gerbe multiplication}, which satisfies an associativity condition over $Y^{[4]}$ as a consequence of $\delta(\sigma) = 1$. These remarks show that bundle gerbes can be equivalently defined as follows. 

\begin{definition}
\label{def:bundlegerbe}
An \emph{abelian bundle gerbe} on $M$ is a triple $(\cP,Y,\mu)$ consisting of:
\begin{itemize}
\item A surjective submersion $Y\to M$.
		
\item A $\U(1)$ bundle $\cP\to Y^{[2]}$ defined over the fibered product $Y^{[2]} = Y\times_M Y$ of $Y\to M$ with itself. 
		
\item An isomorphism $\mu\colon \pi^{\ast}_{12} \cP \otimes \pi^{\ast}_{13} \cP \to \pi^{\ast}_{13} \cP$ of principal $\U(1)$ bundles over $Y^{[3]} = Y\times_M Y \times_M Y$ such that the following diagram of isomorphisms:
\[ \begin{tikzcd}
\pi^{\ast}_{12}\cP \otimes \pi^{\ast}_{23}\cP \otimes \pi^{\ast}_{34}\cP  \arrow{r}{\pi^{\ast}_{123}\mu\otimes \Id} \arrow[swap]{d}{\Id\otimes \pi^{\ast}_{234}\mu} & \pi^{\ast}_{13}\cP \otimes \pi^{\ast}_{34}\cP \arrow{d}{\pi^{\ast}_{134}\mu} \\%
\pi^{\ast}_{12}\cP \otimes \pi^{\ast}_{24}\cP \arrow{r}{\pi^{\ast}_{124}\mu}& \pi^{\ast}_{14}\cP
\end{tikzcd}
\]
		
\noindent
is commutative over $Y^{[4]} = Y\times_M Y \times_M Y \times_M Y$. Here, the map:
		\begin{equation*}
			\pi_{ij} \colon Y^{[4]} \to Y^{[2]}
		\end{equation*}
		
		\noindent
		with $1 \leq i,j \leq 5$ is the map that omits the factor in $Y^{[4]}$ different from $i$ and $j$.
	\end{itemize}
\end{definition}

\noindent
The last condition in the previous definition ensures the \emph{associativity} of the isomorphism $\mu$. In the main text of the dissertation the symbol $\mu$ in the triple $(\cP,Y,\mu)$ is omitted for ease of notation. 

\begin{definition}
A \emph{connective structure} on an abelian bundle gerbe $(\cP,Y,\mu)$ is a connection $\cA$ on $\cP$ preserved by the isomorphism $\mu\colon \pi^{\ast}_{12} \cP \otimes \pi^{\ast}_{23} \cP \to \pi^{\ast}_{13} \cP$.
\end{definition}

\noindent
That is, a connective structure $\cA$ on $(\cP,Y,\mu)$ induces connections on $\pi^{\ast}_{12} \cP$, $\pi^{\ast}_{23} \cP$ and $\pi^{\ast}_{13} \cP$ that are preserved by $\mu$.

\begin{definition}
	A \emph{curving} on an abelian bundle gerbe with connective structure $(\cP,Y,\mu,\cA)$ is a two-form $b\in \Omega^2(Y)$ such that:
	\begin{equation*}
		F_{\cA} = \pi^{\ast}_2 b - \pi^{\ast}_1 b
	\end{equation*}
	
	\noindent
	as an equation of two-forms defined on $Y^{[2]}$, where $F_{\cA}$ is the curvature of $\cA$.
\end{definition}

\noindent
Since $F_{\cA}$ is a closed two-form on $Y^{[2]}$, applying the exterior derivative to the previous equation we obtain:
\begin{eqnarray*}
	\pi^{\ast}_2 \dd b = \pi^{\ast}_1 \dd b
\end{eqnarray*}

\noindent
This implies that there exists a unique three-form $H\in \Omega^3(M)$ on $M$ satisfying:
\begin{equation*}
	\dd b = \pi^{\ast} H
\end{equation*}

\noindent
This three-form is necessarily closed and is called the \emph{curvature} of $b$.

As mentioned above, the definition of a bundle gerbe $(\mathcal{P}, Y, \mu)$ given in \ref{def:bundlegerbe}, and originally discovered by Murray \cite{Murray}, provides a concrete presentation of a more abstract and fundamental object: a \emph{gerbe} in the sense of Giraud's non-abelian cohomology \cite{Giraud}. Roughly speaking, according to Giraud, see also \cite{Brylinski}, a gerbe on $M$ with band $\mathcal{A}$, where $\mathcal{A}$ denotes a sheaf of abelian groups on $M$, is a particular type of \emph{stack in groupoids} over $M$. For our purposes, that is, for the purpose of obtaining a mathematical object equivalent to an abelian bundle gerbe as defined above,  the band of the gerbe must be the sheaf of smooth $\U(1)$-valued functions on $M$, denoted $\U(1)_M$. To assimilate the definition proposed by Giraud, we must first recall why the notion of a sheaf, which occurs naturally in the theory of principal and vector bundles, needs to be elevated to that of a stack. A \emph{sheaf} on a manifold $M$ assigns a set of data, for example an abelian group of smooth functions, to every open set $U \subseteq M$, with the crucial property that compatible local data can be glued together uniquely to yield global data. This works wonderfully for objects like functions or sections of bundles and fibrations. However, it is insufficient for geometric objects that possess non-trivial automorphisms, such as principal bundles. More precisely, two principal bundles can be locally isomorphic in multiple ways, and a sheaf of sets cannot capture this extra layer of information. A \emph{stack} remedies this by assigning not just a set but a \emph{category}, specifically, a groupoid $\mathcal{G}(U)$, to each open set $U\subset M$. It demands that not only objects (like principal bundles) can be glued, but that their \emph{morphisms} (bundle isomorphisms) can also be glued. More formally, a stack in groupoids is a category-valued presheaf $\mathcal{G}$ that satisfies a descent condition: for any open cover $\{U_i\}$ of an open set $U$, the category $\mathcal{G}(U)$ is equivalent to the category of "descent data". This data consists of:
\begin{itemize}
\item A collection of objects $a_i \in \operatorname{Ob}(\mathcal{G}(U_i))$ for each $i$.
\item A collection of morphisms $\phi_{ij} \colon a_i|_{U_{ij}} \to a_j|_{U_{ij}}$ for each pair $(i,j)$.
\item These morphisms must satisfy a cocycle condition on triple overlaps $U_{ijk}$: $\phi_{ik} = \phi_{jk} \circ \phi_{ij}$.
\end{itemize}

\noindent
This structure precisely captures the data needed to glue objects that have a notion of isomorphism between them. With these provisos in mind, a gerbe can be defined now a stack $\mathcal{G}$ that satisfies two further conditions:
\begin{enumerate}
\item \emph{Local Non-Emptiness:} The stack is locally non-empty. This means that for any point in $M$, there is a neighborhood $U$ over which the category $\mathcal{G}(U)$ contains at least one object. Intuitively, this means the geometric object that the gerbe represents can always be localized.

\item \emph{Transitivity:} Any two objects in $\mathcal{G}(U)$ are locally isomorphic. This means that for any two objects $a, b \in \operatorname{Ob}(\mathcal{G}(U))$, there exists an open cover $\{U_j\}$ of $U$ such that the restrictions $a|_{U_j}$ and $b|_{U_j}$ are isomorphic in the category $\mathcal{G}(U_j)$. This is the essential \emph{gerby} property: it asserts that while there may be many local descriptions of our geometric object, they are all locally equivalent. The only variation is in \emph{how} they are equivalent.
\end{enumerate}

\noindent
A bundle gerbe $(\mathcal{P}, Y \xrightarrow{\pi} M, \mu)$ allows us to explicitly construct such a stack, which we will also call $\mathcal{G}$. For any open set $U \subseteq M$, we define the category $\mathcal{G}(U)$ as follows:

\begin{enumerate}
\item \emph{Objects of the Category $\mathcal{G}(U)$}: An \emph{object} in the category $\mathcal{G}(U)$ is a smooth section of the submersion $Y \to M$ over $U$, that is:
\begin{equation*}
\operatorname{Ob}(\mathcal{G}(U)) := \{ s \colon U \to Y \mid \pi \circ s = \mathrm{id}_U \}	 
\end{equation*}
 
\item \emph{Morphisms of the Category $\mathcal{G}(U)$}: Given two objects $s_1, s_2 \in \operatorname{Ob}(\mathcal{G}(U))$, a \emph{morphism} from $s_1$ to $s_2$ is a smooth section $\alpha$ of the pullback line bundle $(s_1, s_2)^*\mathcal{P}$ over $U$.
\begin{equation*}
\operatorname{Hom}_{\mathcal{G}(U)}(s_1, s_2) := \Gamma(U, (s_1, s_2)^*\mathcal{P})	 
\end{equation*}

\noindent
Here, $(s_1, s_2)$ is the map $U \to Y^{[2]}$ given by $x \mapsto (s_1(x), s_2(x))$.
	
\item \emph{Composition of Morphisms:} To define composition, consider three objects $s_1, s_2, s_3$ and two morphisms $\alpha \in \operatorname{Hom}(s_1, s_2)$ and $\beta \in \operatorname{Hom}(s_2, s_3)$. We need to define their composite, $\beta \circ \alpha \in \operatorname{Hom}(s_1, s_3)$. The gerbe multiplication $\mu$ provides an isomorphism $\pi_{12}^*\mathcal{P} \otimes \pi_{23}^*\mathcal{P} \xrightarrow{\sim} \pi_{13}^*\mathcal{P}$ over $Y^{[3]}$. We can pull this back along the map $(s_1, s_2, s_3)\colon U \to Y^{[3]}$. This gives an isomorphism on $U$:
\begin{equation*}
\mu_{(s_1,s_2,s_3)}\colon (s_1,s_2)^*\mathcal{P} \otimes (s_2,s_3)^*\mathcal{P} \xrightarrow{\sim} (s_1,s_3)^*\mathcal{P}	 
\end{equation*}

\noindent
The composition of morphisms is then defined as:
\begin{equation*}
\beta \circ \alpha := \mu_{(s_1,s_2,s_3)}(\alpha \otimes \beta)
\end{equation*}

\noindent
The associativity of this composition is guaranteed by the associativity condition that $\mu$ satisfies over $Y^{[4]}$. The identity morphism for an object $s \colon U \to Y$ is given by a canonical trivializing section of $(s,s)^*\mathcal{P}$, which is assumed to exist for any bundle gerbe.
\end{enumerate}

\begin{remark}
A section $s$ can be thought of as a local "choice of trivialization" for the gerbe. Since $\pi$ is a submersion, such sections are guaranteed to exist locally, but there is no canonical choice. The collection of all possible local choices constitutes the objects of our category.
\end{remark}

\begin{remark}
If objects are local trivializations, then morphisms are the "gauge transformations" between them. The line bundle $\mathcal{P} \to Y^{[2]}$ can be thought of as the "space of transformations" between any two points $(y_1, y_2)$ in the same fiber of $Y \to M$. By pulling $\mathcal{P}$ back via the sections $s_1$ and $s_2$, we obtain a line bundle over $U$ whose sections are precisely the smooth family of transformations relating $s_1$ to $s_2$.
\end{remark}

\begin{remark}
The gerbe multiplication $\mu$ acts as the \emph{associator} for our structure. The fact that $\mu$ itself satisfies a higher associativity condition (represented by a commutative pentagon diagram over $Y^{[4]}$) is precisely what ensures that our definition of composition is associative. The identity morphism for an object $s$ is a canonical non-vanishing section of $(s,s)^*\mathcal{P}$, whose existence is part of the bundle gerbe data.
\end{remark}

\noindent
The structure we have just defined, $\mathcal{G} = \{\mathcal{G}(U)\}_{U \subseteq M}$, forms a stack in groupoids. We can check Giraud's axioms schematically:
\begin{itemize}
\item \emph{It is a stack:} The definitions of objects and morphisms are local in nature, and gluing conditions for sections ensure that $\mathcal{G}$ is a stack.
\item \emph{It is a groupoid:} Every morphism $\alpha \in \operatorname{Hom}(s_1, s_2)$ is invertible because the fibers of the line bundle $\mathcal{P}$ are one-dimensional vector spaces, so any non-zero section $\alpha(x)$ is invertible.
\item \emph{Non-empty:} Since the map $\pi\colon Y \to M$ is a surjective submersion, it admits local sections. Therefore, for any point in $M$, there is a neighborhood $U$ for which $\mathcal{G}(U)$ is non-empty.
\item \emph{Transitive:} Given any two objects $s_1, s_2 \in \operatorname{Ob}(\mathcal{G}(U))$, the space of morphisms between them, $\Gamma(U, (s_1, s_2)^*\mathcal{P})$, is the space of sections of a line bundle. Line bundles are always locally trivial, meaning we can find a cover $\{U_j\}$ of $U$ and a non-vanishing (hence invertible) section of $(s_1, s_2)^*\mathcal{P}$ over each $U_j$. This section is precisely a local isomorphism between $s_1|_{U_j}$ and $s_2|_{U_j}$.
\end{itemize}

\noindent
The \emph{band} of this gerbe is the sheaf of automorphisms of any of its objects. Let's compute the automorphism group for an object $s \in \operatorname{Ob}(\mathcal{G}(U))$:
\begin{equation*}
\operatorname{Aut}_{\mathcal{G}(U)}(s) = \operatorname{Hom}_{\mathcal{G}(U)}(s, s) = \Gamma(U, (s,s)^*\mathcal{P})
\end{equation*}

\noindent
Since $\mathcal{P} \to Y^{[2]}$ is a $\U(1)$-bundle, the pullback bundle $(s,s)^*\mathcal{P}$ is a principal $\U(1)$-bundle over $U$. However, the bundle gerbe structure includes a canonical trivialization (identity multiplication) of this bundle. Therefore, its sections are simply maps from $U$ to $\U(1)$:
\begin{equation*}
\operatorname{Aut}_{\mathcal{G}(U)}(s) \cong C^\infty(U, \U(1))
\end{equation*}

\noindent
The sheaf of these automorphism groups is the sheaf $\U(1)_M$ of smooth $\U(1)$-valued functions on $M$. This is precisely the \emph{band} of the gerbe. In conclusion, a bundle gerbe $(\mathcal{P}, Y \to M, \mu)$ as defined previously is a concrete presentation of a Giraud gerbe with band $\U(1)_M$. This construction provides the essential bridge between the convenient, hands-on calculus of bundle gerbes and the powerful, abstract theory of stacks and higher-order cohomology.








\bibliographystyle{ThesisStyle}

\end{document}